\documentclass[12pt,a4paper]{article}

\usepackage{stmaryrd}

\newcommand{\podot}{\varolessthan}

\newcommand{\bodot}{\varogreaterthan}

\newcommand{\botimes}{\varogreaterthan}

\usepackage{hyperref}

\newcommand{\pageformat}[6]{\setlength{\hoffset}{-1in}
                  \setlength{\voffset}{-1in}
                  \addtolength{\hoffset}{#5}
                            \addtolength{\voffset}{#6}
                            \setlength{\oddsidemargin}{#1}
                            \setlength{\evensidemargin}{#2}
                            \setlength{\textwidth}{\paperwidth}
                  \addtolength{\textwidth}{-\oddsidemargin}
                  \addtolength{\textwidth}{-\evensidemargin}
                  \addtolength{\textwidth}{-\marginparsep}
                  \addtolength{\textwidth}{-\marginparwidth}
                            \setlength{\topmargin}{#3}
                            \setlength{\textheight}{\paperheight}
                  \addtolength{\textheight}{-\topmargin}
                  \addtolength{\textheight}{-\headheight}
                  \addtolength{\textheight}{-\headsep}
                  \addtolength{\textheight}{-\footskip}
                  \addtolength{\textheight}{-#4}}
\pageformat{3cm}{2cm}{25mm}{25mm}{1pt}{0pt}

\usepackage{ifthen}
\newboolean{article}
    \setboolean{article}{true}
\newboolean{report}
\newboolean{book}
\newboolean{letter}
\newboolean{german}
\newboolean{italian}
\newboolean{nobaselinestretch}
\newboolean{nosectionappendix}
\newboolean{oldtoc}
\newboolean{nosectionequn}
\newboolean{notheorem}
\ifthenelse{\boolean{german}}{
    \usepackage{german}}{}

\usepackage[latin1]{inputenc}
\usepackage{arttoc}

\usepackage{amsmath}
\usepackage{amssymb}
\usepackage[mathscr]{eucal}

\ifthenelse{\boolean{notheorem}}{}{
    \usepackage{theorem}}

\ifthenelse{\boolean{nobaselinestretch}}{}{
    \renewcommand{\baselinestretch}{1.25}}

\newenvironment{env}[2]{\begin{#1}#2\end{#1}}{}
    \newcommand{\beq}[1]{\begin{env}{equation}{#1}}
    \newcommand{\beqn}[1]{\begin{env}{equation*}{#1}}
    \newcommand{\bal}[1]{\begin{env}{align}{#1}}
    \newcommand{\baln}[1]{\begin{env}{align*}{#1}}
    \newcommand{\bga}[1]{\begin{env}{gather}{#1}}
    \newcommand{\bgan}[1]{\begin{env}{gather*}{#1}}
    \newcommand{\bflal}[1]{\begin{env}{flalign}{#1}}
    \newcommand{\bflaln}[1]{\begin{env}{flalign*}{#1}}
    \newcommand{\bmu}[1]{\begin{env}{multline}{#1}}
    \newcommand{\bmun}[1]{\begin{env}{multline*}{#1}}
    \newcommand{\bsp}[1]{\begin{env}{split}{#1}}

    \newcommand{\eeq}{\end{env}}
    \newcommand{\eeqn}{\end{env}}
    \newcommand{\eal}{\end{env}}
    \newcommand{\ealn}{\end{env}}
    \newcommand{\ega}{\end{env}}
    \newcommand{\egan}{\end{env}}
    \newcommand{\eflal}{\end{env}}
    \newcommand{\eflaln}{\end{env}}
    \newcommand{\emu}{\end{env}}
    \newcommand{\emun}{\end{env}}
    \newcommand{\esp}{\end{env}}

\newcommand{\lf}{\vspace{2ex}}
\newcommand{\bulletline}[1][]{\lf\noindent~\hfill$\bullet\bullet\bullet$\hfill~

\lf\noindent\bf{#1}}

\renewcommand{\bf}[1]{\textbf{#1}}
\renewcommand{\it}[1]{\textit{#1}}
\renewcommand{\sc}[1]{\textsc{#1}}
\renewcommand{\sf}[1]{\textsf{#1}}

\renewcommand{\tt}[1]{\texttt{#1}}
\newcommand{\hl}[1]{\bf{\it{#1}}}
\newcommand{\mrm}[1]{\mathrm{#1}}
\newcommand{\mbf}[1]{\mathbf{#1}}
\newcommand{\msf}[1]{\text{\small $\sf{#1}$}}

\newcommand{\cmc}[1]{\mathcal{#1}}
\newcommand{\eus}[1]{\mathscr{#1}}
\newcommand{\euf}[1]{\mathfrak{#1}}
\newcommand{\bb}[1]{\mathbb{#1}}
\newcommand{\msmall}[1]{{\setlength{\arraycolsep}{.6ex}\text{\small$#1$}}}

\newcommand{\mscriptsize}[1]{{\setlength{\arraycolsep}{.3ex}\text{\scriptsize$#1$}}}
\newcommand{\mtiny}[1]{{\setlength{\arraycolsep}{.3ex}\text{\tiny$#1$}}}
\newcommand{\nbd}[1]{$#1$\nobreakdash--}
\newcommand{\ol}[1]{\overline{#1}}
\newcommand{\ul}[1]{\underline{#1}}
\newcommand{\wt}[1]{\widetilde{#1}}
\newcommand{\wh}[1]{\widehat{#1}}
\newcommand{\ve}{\varepsilon}
\newcommand{\vt}{\vartheta}
\newcommand{\vk}{\varkappa}
\newcommand{\vp}{\varphi}
\newcommand{\om}{\omega}
\newcommand{\Om}{\Omega}
\newcommand{\abs}[1]{\left\lvert#1\right\rvert}
\newcommand{\norm}[1]{\left\lVert#1\right\rVert}

\newcommand{\bnorm}[1]{\bigl\lVert#1\bigr\rVert}

\newcommand{\snorm}[1]{\norm{\smash{#1}}}

\newcommand{\family}[1]{\left(#1\right)}

\newcommand{\bfam}[1]{\bigl(#1\bigr)}
\newcommand{\Bfam}[1]{\Bigl(#1\Bigr)}
\newcommand{\AB}[1]{\langle#1\rangle}
\newcommand{\bAB}[1]{\bigl\langle#1\bigr\rangle}
\newcommand{\BAB}[1]{\Bigl\langle#1\Bigr\rangle}
\newcommand{\CB}[1]{\{#1\}}
\newcommand{\bCB}[1]{\bigl\{#1\bigr\}}
\newcommand{\BCB}[1]{\Bigl\{#1\Bigr\}}
\newcommand{\SB}[1]{[#1]}
\newcommand{\bSB}[1]{\bigl[#1\bigr]}
\newcommand{\BSB}[1]{\Bigl[#1\Bigr]}

\newcommand{\RO}[1]{[#1)}

\newcommand{\Matrix}[1]{\begin{pmatrix}#1\end{pmatrix}}
\newcommand{\SMatrix}[1]{\msmall{\Matrix{#1}}}

\newcommand{\sMatrix}[1]{\mscriptsize{\Matrix{#1}}}
\newcommand{\tMatrix}[1]{\mtiny{\Matrix{#1}}}

\newcommand{\rtMatrix}[1]{\raisebox{.3ex}{\tMatrix{#1}}}

\newcommand{\sbars}[1]{\:\bar{#1}^s\:}
\newcommand{\sodots}{\sbars{\odot}}
\newcommand{\set}[2][]{
    \ifthenelse{\equal{#1}{}}{
        \CB{#2}}{
        \CB{#1~|~#2}}}
\newcommand{\bset}[2][]{
    \ifthenelse{\equal{#1}{}}{
        \bCB{#2}}{
        \bCB{#1~|~#2}}}
\newcommand{\Bset}[2][]{
    \ifthenelse{\equal{#1}{}}{
        \BCB{#2}}{
        \BCB{#1~\big|~#2}}}
\newcommand{\zero}{\CB{0}}

\DeclareMathOperator{\ls}{\normalfont\msf{span}}

\DeclareMathOperator{\ssls}{\scriptstyle\sf{span}}
\DeclareMathOperator{\cls}{\ol{\ls}}

\DeclareMathOperator{\sscls}{\ol{\ssls}}
\DeclareMathOperator*{\limind}{lim\,ind}

\DeclareMathOperator{\id}{\normalfont\msf{id}}
\DeclareMathOperator{\sid}{\mtiny{\sf{id}}}
\DeclareMathOperator{\ssid}{\scriptstyle\sf{id}}

\DeclareMathOperator{\alg}{\normalfont\msf{alg}}
\renewcommand{\ker}{\operatorname{\msf{ker}}}
\renewcommand{\dim}{\operatorname{\msf{dim}}}

\renewcommand{\Im}{\operatorname{\msf{Im}}}

\newcommand{\C}{\bb{C}}

\newcommand{\E}{\bb{E}}

\newcommand{\bJ}{\bb{J}}
\newcommand{\K}{\bb{K}}

\newcommand{\N}{\bb{N}}

\newcommand{\Q}{\bb{Q}}
\newcommand{\R}{\bb{R}}
\newcommand{\bS}{\bb{S}}
\newcommand{\T}{\bb{T}}
\newcommand{\Z}{\bb{Z}}
\newcommand{\cA}{\cmc{A}}
\newcommand{\cB}{\cmc{B}}
\newcommand{\cC}{\cmc{C}}
\newcommand{\cD}{\cmc{D}}
\newcommand{\cF}{\cmc{F}}
\newcommand{\cG}{\cmc{G}}
\newcommand{\cI}{\cmc{I}}
\newcommand{\cJ}{\cmc{J}}

\newcommand{\cO}{\cmc{O}}

\newcommand{\sB}{\eus{B}}

\newcommand{\sE}{\eus{E}}
\newcommand{\sF}{\eus{F}}

\newcommand{\sK}{\eus{K}}
\newcommand{\sL}{\eus{L}}

\newcommand{\sN}{\eus{N}}

\newcommand{\sS}{\eus{S}}
\newcommand{\sT}{\eus{T}}
\newcommand{\sU}{\eus{U}}

\newcommand{\ee}{\euf{e}}

\newcommand{\en}{\euf{n}}

\newcommand{\ep}{\euf{p}}

\newcommand{\eG}{\euf{G}}
\newcommand{\eH}{\euf{H}}

\newcommand{\eL}{\euf{L}}

\newcommand{\U}{\mbf{1}}
\newcommand{\F}{{\mrm{F}}}
\newcommand{\G}{\Gamma}

\newcommand{\DG}{{\mrm{I}\hspace{-0.3ex}\G}}

\newcommand{\I}{{I\!\!\!\;I}}
\newcommand{\f}{\text{\scriptsize$\sF$}}
\newcommand{\s}{\text{\scriptsize$\sS$}}

\ifthenelse{\boolean{nosectionequn}}{}{
    \numberwithin{equation}{section}
    }

\ifthenelse{\boolean{article}\or\boolean{letter}\or\boolean{nosectionequn}}{
    \setboolean{nosectionappendix}{true}}{}
\ifthenelse{\boolean{nosectionappendix}}{}{
    \renewcommand{\appendix}{
        \chapter*{\appendixname}
        \addcontentsline{toc}{chapter}{\appendixname}
        \renewcommand{\thesection}{\Alph{section}}
        \setcounter{section}{0}}}
   
\ifthenelse{\boolean{report}\or\boolean{book}}{
    }{}

\ifthenelse{\boolean{notheorem}}{}{
        \newcommand{\mnname}{Mathematical note.}
        \newcommand{\enname}{End of the note.}
        \newcommand{\definame}{Definition.}
        \newcommand{\propname}{Proposition.}
        \newcommand{\lemname}{Lemma.}
        \newcommand{\exname}{Example.}
        \newcommand{\exername}{Exercise.}
        \newcommand{\remname}{Remark.}
        \newcommand{\obname}{Observation.}
        \newcommand{\thmname}{Theorem.}
        \newcommand{\corname}{Corollary.}
        \newcommand{\proofname}{Proof.}
    \ifthenelse{\boolean{german}}{
        \renewcommand{\mnname}{Mathematische Notiz.}
        \renewcommand{\enname}{Ende der Notiz.}
        \renewcommand{\exname}{Beispiel.}
        \renewcommand{\exername}{ï¿½bung.}
        \renewcommand{\remname}{Bemerkung.}
        \renewcommand{\obname}{Beobachtung.}
        \renewcommand{\thmname}{Satz.}
        \renewcommand{\corname}{Korollar.}
        \renewcommand{\proofname}{Beweis.}}{}
    \ifthenelse{\boolean{italian}}{
        \renewcommand{\mnname}{Nota matematica.}
        \renewcommand{\enname}{Fina della nota.}
        \renewcommand{\definame}{Definizione.}
        \renewcommand{\propname}{Proposizione.}
        \renewcommand{\exname}{Esempio.}
        \renewcommand{\exername}{Esercizio.}
        \renewcommand{\remname}{Nota.}
        \renewcommand{\obname}{Osservazione.}
        \renewcommand{\thmname}{Teorema.}
        \renewcommand{\corname}{Corollario.}
        \renewcommand{\proofname}{Dimostrazione.}
       
       \renewcommand{\refname}{Riferimenti bibliografici}

       \renewcommand{\appendixname}{Appendice}

       }{}
    \theoremheaderfont{\normalfont\bfseries}
    \theoremstyle{change}
        \theorembodyfont{\rmfamily}
            \newtheorem{emp}{}[section]
                \newcommand{\bemp}[1][]{
                    \begin{emp}\hskip-\labelsep\bf{#1}\hskip\labelsep}
                \newcommand{\eemp}{\end{emp}}
\newtheorem{itemp}[emp]{}
                \newcommand{\bitemp}[1][]{
                    \begin{itemp}\hskip-\labelsep\bf{#1}\hskip\labelsep\normalfont\itshape}
                \newcommand{\eitemp}{\end{itemp}}
            \newtheorem{mn}[emp]{\mnname}
                \newcommand{\bnm}{\begin{mn}~\begin{quotation}\renewcommand{\baselinestretch}{1}\small\noindent\ignorespaces}
                \newcommand{\enm}{\end{quotation}\hfill\bf{\enname}\end{mn}}
            \newtheorem{ex}[emp]{\exname}
                \newcommand{\bex}{\begin{ex}}
                \newcommand{\eex}{\end{ex}}
            \newtheorem{exer}[emp]{\exername}
                \newcommand{\bexer}{\begin{exer}}
                \newcommand{\eexer}{\end{exer}}
            \newtheorem{defi}[emp]{\definame}
                \newcommand{\bdefi}{\begin{defi}}
                \newcommand{\edefi}{\end{defi}}
            \newtheorem{rem}[emp]{\remname}
                \newcommand{\brem}{\begin{rem}}
                \newcommand{\erem}{\end{rem}}
            \newtheorem{ob}[emp]{\obname}
                \newcommand{\bob}{\begin{ob}}
                \newcommand{\eob}{\end{ob}}
        \theorembodyfont{\normalfont\itshape}
            \newtheorem{thm}[emp]{\thmname}
                \newcommand{\bthm}{\begin{thm}}
                \newcommand{\ethm}{\end{thm}}
            \newtheorem{prop}[emp]{\propname}
                \newcommand{\bprop}{\begin{prop}}
                \newcommand{\eprop}{\end{prop}}
            \newtheorem{cor}[emp]{\corname}
                \newcommand{\bcor}{\begin{cor}}
                \newcommand{\ecor}{\end{cor}}
            \newtheorem{lem}[emp]{\lemname}
                \newcommand{\blem}{\begin{lem}}
                \newcommand{\elem}{\end{lem}}
\newenvironment{empn}[1]{\lf\noindent\bf{#1}\ignorespaces\hskip\labelsep}{\lf}
		\newcommand{\bempn}[1]{\begin{empn}{#1}}
		\newcommand{\eempn}{\end{empn}}
		\newcommand{\bitempn}[1]{\begin{empn}{#1}\normalfont\itshape}
		\newcommand{\eitempn}{\end{empn}}
                \newcommand{\bnmn}{\begin{empn}{\mnname}~\begin{quotation}\renewcommand{\baselinestretch}{1}\small\noindent\ignorespaces}
                \newcommand{\enmn}{\end{quotation}\hfill\bf{\enname}\end{empn}}
		\newcommand{\bexn}{\begin{empn}{\exname}}
		\newcommand{\eexn}{\end{empn}}
		\newcommand{\bexern}{\begin{empn}{\exername}}
		\newcommand{\eexern}{\end{empn}}
		\newcommand{\bdefin}{\begin{empn}{\definame}}
		\newcommand{\edefin}{\end{empn}}
		\newcommand{\bremn}{\begin{empn}{\remname}}
		\newcommand{\eremn}{\end{empn}}
		\newcommand{\bobn}{\begin{empn}{\obname}}
		\newcommand{\eobn}{\end{empn}}
		\newcommand{\bthmn}{\bitempn{\thmname}}
		\newcommand{\ethmn}{\eitempn}

\newcommand{\qedsymbol}{~\rule[-0.35mm]{2mm}{2mm}}
    \newcounter{proof}[emp]
    \newenvironment{Proof}[1]{
        \vspace{1ex}
        \renewcommand{\item}[1][\stepcounter{proof}(\roman{proof})]%
            {##1\hskip\labelsep}
        \noindent\textsc{#1\hskip\labelsep}}{
        \nolinebreak\qedsymbol}
    \newcommand{\proof}[1][\proofname]{
        \begin{Proof}{#1}\ignorespaces}
    \newcommand{\qed}{\end{Proof}}
    \newcommand{\noqed}{
        \renewcommand{\qedsymbol}{}
        \end{Proof}}}
    \ifthenelse{\boolean{italian}}{
        \renewcommand{\proofname}{Dimostrazione.}}{}

\usepackage[varg]{txfonts}
\usepackage{bbm}
\newcommand{\bbm}[1]{\mathbbm{#1}}
\newcommand{\be}{\bbm{e}}
\newcommand{\bbf}{\bbm{f}}

\newcommand{\bn}{\bbm{n}}
\newcommand{\bm}{\bbm{m}}
\newcommand{\br}{\bbm{r}}
\newcommand{\bs}{\bbm{s}}
\newcommand{\bt}{\bbm{t}}

\DeclareMathOperator{\inv}{\msf{inv}}

\newcommand{\lcom}{\,\overset{<}{,}\,}

\usepackage[matrix,arrow,curve]{xy}

\renewcommand{\thefootnote}{[\alph{footnote}]}

\usepackage{color}

\makeindex
\usepackage{makeidx}

\begin{document}

\bibliographystyle{amsalpha}

\title{CP-Semigroups and Dilations,\\Subproduct Systems and Superproduct Systems:\\The Multi-Parameter Case and Beyond}

%%%% BO postarxiv
%\author{Orr Moshe Shalit and Michael Skeide}
\author{Orr Moshe Shalit and Michael Skeide{\renewcommand{\thefootnote}{}}\setcounter{footnote}{0}}
%%%% EO postarxiv

\date{March 2020; this revision May 2022\\}

\maketitle 

{\setlength{\baselineskip}{1ex}

% \vspace{-3ex}
\begin{abstract}
\noindent
%%%% BO postarxiv
%
% 
% MOVED FOOTNOTE LINK TO OTHER PLACE TO MOVE THE RED MARKER, I THINK IT LOOKS BETTER THIS WAY
% 
%{\renewcommand{\thefootnote}{}\footnote{MSC 2010: 46L55; 46L07; 46L53 . Keywords: Quantum dynamics; quantum probability; product systems; $E_0$-semigroups and $E$-semigroups; Markov semigroups and CP-semigroups; dilations.}}\setcounter{footnote}{0}
%%%% EO postarxiv
These notes are the output of a decade of research on how the results about dilations of one-parameter CP-semigroups with the help of product systems, can be put forward to \nbd{d}para\-meter semigroups -- and beyond. While existing work on the two- and \nbd{d}parameter case is based on the approach via the \it{Arveson-Stinespring correspondence} of a CP-map by Muhly and Solel \cite{MuSo02} (and limited to von Neumann algebras), here we explore consequently the approach via Paschke's \it{GNS-correspondence} of a CP-map \cite{Pas73} by Bhat and Skeide \cite{BhSk00}. (A comparison is postponed to Appendix \ref{vNAPP}\ref{vNcomm}.)

The generalizations are multi-fold, the difficulties often enormous. In fact, our only true if-and-only-if theorem, is the following: A \it{Markov semigroup} over (the opposite of) an Ore monoid admits a \it{full} (strict or normal) \it{dilation} if and only if its \it{GNS-subproduct system} embeds into a product system. Already earlier, it has been observed that the GNS- (respectively, the Arveson-Stinespring) correspondences form a \it{subproduct system}, and that the main difficulty is to embed that into a product system. Here we add, that every dilation comes along with a \it{superproduct system} (a product system if the dilation is full). The latter may or may not contain the GNS-subproduct system; it does, if the dilation is \it{strong} -- but not only. 

Apart from the many positive results pushing forward the theory to large extent, we provide plenty of counter examples for almost every desirable statement we could not prove. Still, a small number of open problems remains. The most prominent: Does there exist a CP-semigroup that admits a dilation, but no strong dilation? Another one: Does there exist a Markov semigroup that admits a (necessarily strong) dilation, but no full dilation?
\end{abstract}

\vfill
\newpage
}

\tableofcontents

% }
\newpage

\section{Introduction}\label{intro}

As stated in the abstract, we present here the output of a decade of efforts to push forward the results about dilations of CP-semigroup with the help of product systems, from the one-parameter case to the multi-parameter case -- and beyond. Only rarely will we be able to derive results as ``round and nice'' as we know them from the one-parameter case. Many theorems have only a forward implication, but do not allow, when put together, to close the circle. Still, they frequently deliver powerful methods for constructing dilations or powerful criteria for establishing non-existence of such. And for most circles we could not close, we provide counter examples that prove they really do not close. (We even shed more light on the one-parameter case, illustrating that also for that case the situation is far from being as finalized as the literature makes us believe.) So, the readers may not expect a round and closed treatment bringing the theory to an end. But, they may expect a comprehensive up-to-date toolbox for tackling the dilation problem even under the most general circumstances.

\lf
We start this introduction with a brief account about those known results in the one-pa\-rameter case that motivate our approach.

Let $\cB$ denote a unital \nbd{C^*}algebra. If $T=\bfam{T_t}_{t\in\R_+}$ is a \hl{CP-semigroup}\phantomsection\index{CP-semigroup} on $\cB$ (that is, a semigroup of completely positive maps $T_t$ on $\cB$), then we may associate with each $T_t$ Paschke's \cite{Pas73} \hl{GNS-construction}\index{GNS-construction}\index{CP-map!GNS-construction} $(\sE_t,\xi_t)$. That is, $\sE_t$ is a correspondence over $\cB$ that is generated by a single element $\xi_t$ fulfilling $T_t(b)=\AB{\xi_t,b\xi_t}$. We easily verify that $\xi_{s+t}\mapsto\xi_s\odot\xi_t$ extends as an isometric bimodule map $w_{s,t}\colon\sE_{s+t}\rightarrow\sE_s\odot\sE_t$, that these \it{coproducts} iterate associatively, that $(\sE_0,\xi_0)=(\cB,\U)$, and that the marginal maps $w_{0,t}$ and $w_{t,0}$ are just the canonical identifications $x_t\mapsto\U\odot x_t$ and $x_t\mapsto x_t\odot\U$. In other words, the $w_{s,t}$ turn the family $\bfam{\sE_t}_{t\in\R_+}$ into a \hl{subproduct system}\index{subproduct system}. By definition, the $\xi_t$ fulfill $w_{s,t}\xi_{s+t}=\xi_s\odot\xi_t$, that is, they form a \hl{unit}\index{unit} $\xi^\odot=\bfam{\xi_t}_{t\in\R_+}$.

Every subproduct system $\bfam{\sE_t}_{t\in\R_+}$ embeds into a \hl{product system}\index{product system} $E^\odot$. By this, we mean that $E^\odot$ is a subproduct system that has unitary coproduct maps $u_{s,t}^*$, and that $\sE_t\subset E_t$ such that $u_{s,t}^*$ restricted to $\sE_{s+t}$ is $w_{s,t}$. In particular, it is clear that the $\xi_t\in\sE_t\subset E_t$ form a unit also for $E^\odot$.

Given a product system $E^\odot$ and a \hl{contractive} unit $\xi^\odot$ (that is, $\norm{\xi_t}\le1$ for all $t$), it is possible to construct a Hilbert \nbd{\cB}module $E$, a \hl{unit vector}\index{unit vector} $\xi\in E$ (that is, $\AB{\xi,\xi}=\U$), and an \hl{\nbd{E}semigroup}\index{E-semigroup@\nbd{E}semigroup} (that is, a semigroup of endomorphisms) $\vt=\bfam{\vt_t}_{t\in\R_+}$ on $\sB^a(E)$ such that
\beqn{
\AB{\xi,\vt_t(\xi b\xi^*)\xi}
~=~
\AB{\xi_t,b\xi_t}.
}\eeqn
In other words, with the embedding $i\colon b\mapsto\xi b\xi^*$ of $\cB$ into $\sB^a(E)$ and the \it{expectation} $\ep\colon a\mapsto\AB{\xi,a\xi}$, we get a \phantomsection\hl{dilation} of the contractive CP-semigroup $T_t:=\AB{\xi_t,b\xi_t}$, that is, the diagram
\beq{ \label{1stdil}
\parbox{10cm}{\xymatrix{
\cB	\ar[d]_i	\ar[rr]^{T_t}	&&	\cB
\\
\sB^a(E)		\ar[rr]_{\vt_t}	&&	\sB^a(E)	\ar[u]_\ep
}}
}\eeq
commutes for all $t$.

In a minute, we shall review, very briefly, the basic citations for the aforementioned results; more details in much more general circumstances shall follow, later on throughout these notes. For this introduction, let us just mention that the preceding construction of a dilation of a contractive \hl{one-parameter} CP-semigroup $T$ depends in different stages and in different ways on the order structure of the monoid $\R_+$. For instance, the construction of the product system from the subproduct system of GNS-correspondences depends on the fact that $\R_+$ is totally ordered; the construction of the \nbd{E}semigroup, instead, only depends on the fact that $\R_+$ is directed. It is the scope of these notes to find out what we can say about CP-semigroups and their dilations, when the indexing semigroup is $\R_+^d$ or $\N_0^d$ (\phantomsection\index{multi-parameter case!discrete}\index{multi-parameter case!continuous time}\index{discrete case}\index{continuous time case}the continuous time and the discrete \nbd{d}parameter semigroups) or even a general monoid. It turns out that, like in the one-parameter case, product systems of correspondences play a crucial role in constructing dilations and in understanding them. While, starting with the GNS-correspondences, the construction of a suitable product system has to be replaced by different constructions (which often only work under additional conditions), the construction of an \nbd{E}semigroup works (at least in the unital case) as soon as the monoid is suitably directed, namely, \it{left reversible}. We spend a considerable amount of energy to find out when existence of a dilation guarantees that a CP-semigroup has a product system, and we give several constructions for multi-parameter product systems.

The approach in the beginning of this introduction is from Bhat and Skeide \cite{BhSk00}, where product systems of correspondences occur for the first time; it has not yet been applied directly to the multi-parameter case. There are, however, a number of results by Solel \cite{Sol06} in the discrete two- and \nbd{d}parameter case and by Shalit \cite{Sha08,Sha09} in the continuous two-parameter case, applying the methods from Muhly and Solel \cite{MuSo02}. The approach in \cite{MuSo02} starts from the so-called \it{Arveson-Stinespring correspondences} of the CP-maps $T_t$, which are correspondences over the commutant of $\cB$, $\cB'$, and only works for von Neumann algebras. The relation between the constructions has been made precise in Skeide \cite{Ske03c,Ske08,Ske09} (in terms of the \it{commutant of von Neumann correspondences}) and in Muhly and Solel \cite{MuSo07} (in terms of \it{\nbd{\sigma}duals of \nbd{W^*}correspondences}). Multi-parameter product systems have first been considered by Fowler \cite{Fow02}. \it{Subproduct systems} have been around as the input for several inductive limit constructions; see Sch\"urmann \cite{MSchue93}, Arveson \cite{Arv97a}, Bhat and Skeide \cite{BhSk00}, Muhly and Solel \cite{MuSo02}, and Skeide \cite{Ske06d,Ske03c}. Only rather recently Shalit and Solel \cite{ShaSo09} gave a formal definition and started investigating them and their application to dilations systematically. Almost simultaneously, Bhat and Mukherjee \cite{BhMu10} introduced one-parameter subproduct systems of Hilbert spaces under the name of \it{inclusion systems}, and proved that every such inclusion system is contained in a unique minimal product system.

As opposed with the known observation that CP-semigroups come shipped with subproduct systems, it is surely a key insight of these notes that existence of dilations necessarily leads to \index{superproduct system}\it{superproduct systems}, opening up the way to find necessary criteria for existence of dilations in terms of superproduct systems.

\lf
\lf
The readers -- and we --  have to digest generalizations not just in one, but in several directions; we try our best to make sure that these generalizations do not occur at once but in appetizing portions. (Dilations that are not to $\sB^a(E)$ (not full); product systems, but over more general monoids; subproduct systems and superproduct systems, but isolated from their occurrence from CP-semigroups and dilations; \nbd{E_{(0)}}semigroups that come only with superproduct systems; a notion of minimality that splits into several.)

Sections with general theory take turns with example sections (the latter clearly marked as such). The example sections either illustrate applications of the positive results, or provide counter examples for what can go wrong. In particular, the ``multi-examples'' in the  Example Sections \ref{EXBexSEC} and \ref{EXN02SEC}, which may be considered a sort of culmination point of these notes, unite both aspects. The examples (or, in Section \ref{EXN02SEC}, classes of such) for dilations obtained by applying the the constructive parts of our results, exhibit almost all sorts of bad behaviour that we could not exclude before by theorems.

\lf\noindent
\bf{Overview.~}
In Section \ref{monoSEC}, we put out the general concept of dilation of CP-semigroups. As compared with the situation in the beginning of this introduction, captured in Diagram \eqref{1stdil}, the dilating endomorphism semigroup acts on a general unital \nbd{C^*}algebra $\cA$ and the indexing monoid $\R_+$ is replaced by a general monoid $\bS$. Quickly, we turn to dilations (so-called \it{weak} dilations) where $\cB$ sits as a corner $p\cA p$ in $\cA$ for some projection $p\in\cA$, so that $i$ is just the natural injection and $\ep$ is the compression map $p\bullet p$. In quantum dynamics, where \it{dilation} is thought of as a model for understanding the irreversible evolution (=Markov semigroup=unital CP-semigroup) of a \it{small} system $\cB$ as a projection $\ep$ from a reversible evolution of a \it{big} system $\cA$ into which the small one $\cB$ is immersed, one would rather expect unital embeddings. However, it is an \it{empiric fact} (for instance, if the Markov semigroup is \it{spatial}; see Skeide \cite[Chapters 10+13]{Ske16}; see also Remark \ref{dilhistrem}) that all known unital dilations do ``contain'' a weak dilation sitting inside and that frequently weak dilations may be promoted to unital dilations: Weak dilations are the building blocks of more general dilations; in these notes we restrict, like in large parts of the literature, our attention exclusively to weak dilations.

Dilations of Markov semigroups fulfill an extra property: They are \it{strong} dilations in the sense that $p\vt_t(a)p=T_t(pap)$ for all $a\in\cA$. Large parts of the literature, with Bhat's work probably the only exception, considers only weak dilations that are strong (or \it{regular} in Bhat's terminology). It is a feature of these notes that we do not pose this restriction. (Simply put, we do not see a motivation to justify this restriction; it is also not motivated from classical dilation theory for contractions on a Hilbert space, as most of the (co)isometric dilations of contractions would lead to dilations of CP-semigroups that are not strong; see Section \ref{EXwnsSEC}. Our life would have been much easier, if we put that restriction; but why should we exclude one of the most inspiring sources, classical dilations, from our mind?) Another extra property a dilation may fulfill or not, is being a \it{full} dilation, that is, $\cB$ sitting as a (strictly) full corner in $\cA$. Equivalently, $\cA=\sB^a(E)$ for the Hilbert \nbd{\cB}module $E:=\cA p$. (Being full, is part of practically all notions of minimality in the literature; we see, it is not so strange to expect Hilbert modules playing a role in dilation theory.) We conclude Section \ref{monoSEC}, by discussing the powerful tool of \it{unitalization}. By unitalization, we transform a (contractive) CP-semigroup into a Markov semigroup. At least for strong dilations, this allows to apply large parts of the stronger results about the dilation of Markov semigroups (so-called \it{Markov dilations}) to general CP-semigroups. (It is noteworthy that restricting to $\sB(H)$, frequently done in the literature, excludes such a powerful method, because the unitalization of $\sB(H)$ is not another $\sB(H)$.)

In the Example Section \ref{EXwnsSEC}, we discuss the connection between dilation of CP-semigroups and classical dilation theory in terms of \it{elementary CP-semigroups} and their \it{solidly elementary dilations}. This section also serves the purpose to put at ease both readers who know only the classical side of the theory and readers who only know the quantum dynamical side.

In Section \ref{PSmonoSEC}, we discuss the notion of product system -- but over monoids -- and how units for product systems relate to CP-semigroups. We illustrate that product systems (and the like) and the related semigroups are indexed by opposite monoids, and we explain the choice we fix throughout these notes, whenever (that is, almost always) product systems and semigroups occur in the same context: Product systems (and the like) are indexed by the monoid $\bS$; CP-semigroups and their dilations are indexed by the opposite monoid $\bS^{op}$. Most constructive results in these notes depend on special properties the monoid $\bS$ has to satisfy. In the last part of Section \ref{PSmonoSEC}, we present the necessary theory of algebraic semigroups.

In Section \ref{SPSpbSEC}, we define superproduct systems and subproduct systems -- for formal reasons (superproduct systems are formally nearer to product systems than subproduct systems) in the historically wrong order -- and explore their basic properties, independently of their relation to semigroups or dilations. Much space is reserved for the notion of subsystems of such systems and what it means to embed one into another. This was easy, if we stood in the same category. But, later on, it turns out that we have to be interested, in particular, in how a subproduct systems sits as a subsystem of a superproduct system --  and this (together with the fact that, in applications, the structure maps of a subproduct system need not be adjointable) makes it tricky. We also provide some results about the generation of subsystems.

A single correspondence $E$ generates a discrete one-parameter product system $\bfam{E^{\odot n}}_{n\in\N_0}$. The \it{time ordered} product system, may be considered as what one obtains by \it{exponentiating} this discrete product system to a continuous time one-parameter product system. (The time ordered Fock module is the module analogue of the symmetric Fock space; see Appendix \ref{FockAPP}.) In the Example Section \ref{EXexpSEC}, we apply this insight to discrete \nbd{d}parameter superproduct systems and subproduct systems, \it{exponentiating} them to continuous time \nbd{d}parameter superproduct systems and subproduct systems. We point out that exponentiation respects inclusions and embeddings. This machinery allows to promote many discrete (counter) examples to continuous time (counter) examples. Of course, this also opens up to the study of exponential (super)(sub)product systems, question we do not address in these notes.

In Section \ref{CPspsSEC}, we discuss the \it{GNS-subproduct systems} (see the beginning of this introduction) of CP-semigroups over arbitrary monoids. A major result is, certainly, Theorem \ref{adSPS-CPthm}, which asserts that every adjointable subproduct system is (strictly) Morita equivalent to one of a (strict) CP-semigroup on some $\sB^a(E)$. This is crucial to establish, in the Example Section \ref{EXsubnsupSEC}, existence of CP-semigroups with no strong dilations and Markov semigroups with no dilations at all. Examples of CP-semigroups whose subproduct system is not adjointable, can be found in the Example Section \ref{EXnonadSEC}.

CP-semigroups give rise to subproduct systems. In Section \ref{DilSPSpSEC}, imitating the construction of the product system of an endomorphism semigroup on $\sB^a(E)$ in Skeide \cite{Ske02} we, finally, show that dilations give rise to superproduct systems. The insight that CP-semigroups and other irreversible dynamics lead, in the one or the other way, to subproduct systems, has been around now for quite a while. (The idea to start with the GNS-subproduct system, is new, though, in these notes.) The observation that dilations lead to superproduct systems, is entirely new. It turns out that the superproduct system of a dilation may or may not contain the GNS-subproduct system of the dilated CP-semigroup; a somewhat surprising observation. (Bhat's Example\index{Bhat's example} \ref{Bex} illustrates a failure.) This is, why we introduce the notion of \it{good dilation}: A dilation whose superproduct system does contain the GNS-subproduct system. Strong dilations are good, but not \it{vice versa}; in fact, every non-strong solidly elementary dilation (that is, dilations coming from the classical dilation problem as discussed in Section \ref{EXwnsSEC}) is still a good dilation. By a minimalization procedure, at least in the von Neumann case, every good dilation gives rise to a strong dilation (Observation \ref{algminvNob}). An (adjointable) subproduct system that does not embed into a superproduct system (discussed just before), gives rise to a CP-semigroup that does not only not admit a strong dilation; it also does not admit a good dilation. While the unitalization of this CP-semigroup is a Markov semigroup with no dilation (Markov dilations are strong), we do neither have an example for a non-Markov CP-semigroup that admits no dilation, nor an example for a (necessarily non-Markov) CP-semigroup that admits dilations but only such that are not good.

These and more questions that arise after having established the occurrence of superproduct systems from dilations and the relation with the GNS-subproduct system of the dilated CP-semigroup, we formulate in Section \ref{QSEC}. We formulate the answers we get throughout these notes, but we also state clearly where we do not have answers. Section \ref{QSEC} may be considered an introduction to the second half of our notes; but, of course, it cannot be understood without having first appreciated at least the theoretic part of the first half.

The following sections explore superproduct systems of dilations and what we can do with them. Section \ref{unisupSEC} addresses how superproduct systems behave under unitalization of dilations. Section \ref{leftdilSEC} addresses the question of so-called left dilations for superproduct systems. (Left dilations of a product system furnish a unital endomorphism semigroup for that product system. They would do the same for superproduct systems -- if we had not shown that they can exist only if the superproduct system is a product system. But, if we have a product system and a unital unit (so that the product system contains the GNS-subproduct system of the Markov semigroup determined by the unit), and if the monoid is sufficiently``nice'' (an Ore monoid), then we may construct a left dilation and the endomorphism semigroup it determines, is a dilation of the Markov semigroup; Theorem \ref{Oreindthm}.) In the Example Section \ref{EXpropsupSEC} we furnish dilations with proper superproduct systems.

Appealing to Theorem \ref{Oreindthm}, which grants existence of Markov dilations provided we manage to embed the GNS-subproduct system into a product system, Sections \ref{compSEC}--\ref{EXsubnsupSEC} address the problem to construct such product systems -- basically in the \nbd{d}parameter case(s), but not only. The basic observation is that the \phantomsection\index{one-parameter case}\index{multi-parameter case}\nbd{d}parameter monoids $\N_0^d$ or $\R_+^d$ are products of the one-parameter monoids $\N_0$ or $\R_+$; and for the one-parameter monoids we know how to construct product systems out of subproduct systems. So, assuming we have product systems ${E^k}^\odot$ over $\bS^k$ ($k=1,\ldots,d$), can we put these \it{marginal} product systems together, to turn the family
\beqn{
\bfam{E^1_{t_1}\odot\ldots\odot E^d_{t_d}}_{(t_1,\ldots,t_d)\in\bS}
}\eeqn
into a product system over $\bS=\bS^1\times\ldots\times\bS^d$? Section \ref{compSEC}, gives a complete answer to this question in Theorem \ref{prodthm}. The proof depends on the analysis of structures regarding the permutation groups, which are outsourced to Appendix \ref{popAPP}. A special case, \it{strongly commuting} CP-semigroups, is discussed in Section \ref{strcomSEC}. Not only do we have to push forward the notion of \it{strongly commuting} appropriately; we think that the discussion here might also illuminate a bit more the ``secrets'' of this somewhat mysterious notion. In the Example Section \ref{qconvSEC}, we discuss so-called \it{quantized convolutions semigroups} (introduced by Arveson and studied by Markiewicz). We show that they are strongly commuting in an even stronger sense, and applying our methods we construct dilations for \nbd{d}parameter quantized convolution semigroups. The theory from Section \ref{compSEC} simplifies enormously in the discrete \nbd{d}parameter case and justifies a separate treatment in Section \ref{N0dSEC}. This section is the basis for the Example Section \ref{EXN02SEC}, where we construct for each
%%%% BO new
discrete 
%%%% EO 
two-parameter CP-semigroup on a von Neumann algebra a dilation. Understanding the precise form of three-parameter product systems is also the basis for understanding why the three-parameter subproduct system in the Example Section \ref{EXsubnsupSEC} does not embed into a product system (giving rise to a three-parameter Markov semigroup with no dilation).

Up to that point, the discussion could be kept at a level that is surprisingly algebraic. Every now and then, the condition that certain homomorphisms or CP-maps be strict, pops up. (See the discussion in the second half of Section \ref{CPspsSEC}, starting after Remark \ref{CPASrem} and leading, especially, to Corollary \ref{stri2cor} and to Example \ref{EPSex}. Strictness, there, appears dressed as a simple nondegeneracy condition that makes things work.)
%(Strictness comes, usually, packed as a nondegeneracy condition. For instance, a unital endomorphism $\vt$ of $\sB^a(E)$ is strict if and only if $\vt(\sK(E))$ acts nondegenerately on $E$. The fact that the superproduct system of a full dilation is a product system, depends on strictness of the dilation.)
But with the end of Section \ref{EXsubnsupSEC}, we also have reached the end of attempts to keep topological questions ``behind the scenes''. In particular, we are meeting more and more statements that can no longer be proved in the \nbd{C^*}setting (usually, generalizing easily to the von Neumann setting), but only for von Neumann algebras, modules, and correspondences. Section \ref{topSEC} addresses these questions -- and prepares right away for Section \ref{minSEC} about minimality, where these things will be applied. A thorough introduction to the necessary facts about von Neumann modules (which we will apply exclusively, as opposed with \nbd{W^*}modules) is outsourced to Appendix \ref{vNAPP}. (This appendix also contains an approximately detailed comparison of our approach here based on \cite{BhSk00} and the approach based on \cite{MuSo02}. This comparison cannot be done without discussing the \it{commutant} of von Neumann correspondences -- a concept that we do not need anywhere else in these notes, and which, therefore, is banished to the appendix.) Let us emphasize that we do not at all address questions of ``time''-continuity (continuity with $t$ in a topological monoid). With one exception, all our constructions will deliver continuous things if we start with continuous things -- and work also without continuity. The exception is the construction of a CP-semigroup for a subproduct system in Theorem \ref{adSPS-CPthm}. Here, in order to get something continuous, we would have to replace the direct sum over the members of the subproduct system by a direct integral -- which, of course, requires to have at hand a continuous structure of that field. (One could say measurable field. But continuous fields have always shown to be enough to produce a satisfactory theory, while requiring just measurability definitely causes continuity problem for the semigroups at $t=0$.) Section \ref{minSEC} about minimality, a topic (only) apparently well-understood in the one-parameter case, is probably the toughest of the whole work, with all the ramifications caused by different notions and with an enormous lot of results that depend on hypotheses that cannot always be achieved at the same time; we refer the reader to the section introduction.

After these preparations, our notes culminate, as mentioned, in the last two Example Sections \ref{EXBexSEC} and \ref{EXN02SEC}. Apart from the promised examples of bad behaviour, in the one-parameter non-Markov case in Section \ref{EXBexSEC} and in the two-parameter Markov case in Section \ref{EXN02SEC}, the former presents a thorough treatment of the relation between dilations of CP-semigroups on $\sB(H)$ and the (still classical) theory of dilation of row contractions, while the latter contains the (already mentioned) existence result for dilations of arbitrary (not necessarily normal) discrete two-parameter CP-semigroups on von Neumann algebras.

\lf\lf\lf\noindent
\bf{Conventions, notation, and other preliminaries.~} \phantomsection\label{IC}\index{conventions}
The letter $\bS$ stands for a monoid whose neutral element we always denote by $0$. If we assume that $\bS$ is abelian, we shall write its operation additively, $(s,t)\mapsto s+t$. Otherwise, we write it multiplicatively, $(s,t)\mapsto s\cdot t=:st$, but we continue denoting the neutral element by $0$. Consequently, a \hl{semigroup over $\bS$ on a set $\cB$}\index{semigroup!over a monoid} is a family $T=\bfam{T_t}_{t\in\bS}$ of maps $T_t$ on $\cB$ fulfilling $T_s\circ T_t=T_{st}$ and the monoid condition $T_0=\id_\cB$.

If $A,B,C$ are sets and $(a,b)\mapsto ab$ is a map from $A\times B$ into $C$, by $AB$ we mean the set $\CB{ab\colon a\in A,b\in B}$. Even if $A,B,C$ are (topological) vector spaces, we \bf{do not} adopt any convention where $AB$ would mean the (closed) linear span of all $ab$.\index{product notation!for sets}

Regarding the choice between (abstract) \nbd{W^*}algebras and (concrete) von Neumann algebras, we opt to work with von Neumann algebras. Consequently, we will work with von Neumann modules, not \nbd{W^*}modules. Letters $\cA,\cB,\ldots$ stand for (usually unital) \nbd{C^*}algebras or for von Neumann algebras, letters $E,F,\ldots$ stand for Hilbert (or \nbd{C^*})modules or for von Neumann modules. Recall that \hl{Hilbert module}\index{module!Hilbert} means Hilbert \bf{right} module. The space of \hl{adjointable}\phantomsection\index{adjointable!operator}\index{operator!adjointable} (hence, bounded) operators on $E$ (from $E$ to $F$) we denote by $\sB^a(E)$ ($\sB^a(E,F)$). (In a few occasions, we will write $\sL^a(E,F)$ for the adjointable maps between pre-Hilbert modules; if one of them is Hilbert, then $\sL^a(E,F)=\sB^a(E,F)$.) By $x^*\in E^*$ we denote the map $x^*\colon y\mapsto\AB{x,y}$ with adjoint $x\colon b\mapsto xb$. The \nbd{C^*}algebra of \hl{compact operators}\index{operator!compact (in $\sK(F,E)$)} on $E$ is $\cls\CB{xy^*\colon x,y\in E}$, the closed linear span in $\sB^a(E)$ of the \hl{rank-one operators}\index{operator!rank-one ($xy^*$)} $xy^*$. Their linear span, the so-called \hl{finite-rank operators}\index{operator!finite rank (in $\sF(F,E)$)}, is denoted by $\sF(E)$. The compact operators $\sK(E,F)$ and the finite-rank operators $\sF(E,F)$ from $E$ to $F$ are defined in a similar manner. A \hl{correspondence}\index{correspondence} from $\cA$ to $\cB$ (or Hilbert \nbd{\cA}\nbd{\cB}bimodule) is a Hilbert \nbd{\cB}module with a nondegenerate({\bf !}) left action of $\cA$ by adjointable operators such that the adjoint of the action of $a\in\cA$ is the action of $a^*$. A correspondence is \hl{faithful}\index{correspondence!faithful} if it has a faithful left action. Recall that the \hl{internal tensor product}\index{tensor product}\index{correspondence!(internal) tensor product of} of a correspondence $E$ from $\cA$ to $\cB$ and a correspondence $F$ from $\cB$ to $\cC$ is that unique correspondence $E\odot F$ from $\cA$ to $\cC$ which is generated by elements $x\odot y$ subject to the inner product $\AB{x\odot y,x'\odot y'}=\AB{y,\AB{x,x'}y'}$ and the relation $a(x\odot y)=(ax)\odot y$.

%%%% BO 
% I think we need also to recall multiplier algebras and the strict topology. In particular, I never did like that sentence "A linear map... we call strict, if it is strictly continuous...". Before such a sentence it is required to remind something. 

Recall that for every \nbd{C^*}algebra $\cA$, the \hl{multiplier algebra}\index{multiplier algebra} of $\cA$ is the (up to \nbd{\cA}intertwining isomorphism) unique maximal \nbd{C^*}algebra $M(\cA)$ that contains $\cA$ as an essential ideal. ($M(\cA)$ may be realized as the well-known double centralizers; we discuss this, even for 
%%%% BO new
pre-\nbd{C^*}al\-ge\-bras,
%pre-\nbd{C^*}al\-ge\-bras)
%%%% EO 
in Proposition \ref{convprop}.) The \phantomsection\hl{strict topology}\index{strict!topology} on $M(\cA)$ is the topology induced by the family of seminorms given by $\norm{\bullet a}$ and $\norm{a \bullet}$ ($a \in \cA$); also the strict completion of $\cA$ is a realization of $M(\cA)$. We say a linear map $T\colon M(\cA)\rightarrow M(\cB)$ between two multiplier \nbd{C^*}algebras is \hl{strict}\index{strict!linear map}, if it is strictly continuous on bounded subsets of $M(\cA)$ (by a straightforward application of the \it{closed graph theorem} necessarily into bounded subsets of $M(\cB)$). 
%%%% EO

If $\cA = \sK(E)$, then $M(\cA) = \sB^a(E)$; Kasparov \cite{Kas80}. (See our proof in Proposition \ref{convprop}.) In this case, the strict topology and the \nbd{*}strong (operator) topology coincide on bounded subsets of $\sB^a(E)$; \cite[Proposition 8.1]{Lan95}.  A homomorphism $\vt\colon\sB^a(E)\rightarrow\sB^a(F)$ is strict if and only if the set $\vt(\sK(E))F$ is total in $\vt(\id_E)F$. (This is more or less \cite[Proposition 2.5]{Lan95}. We prove a more general statement in Proposition \ref{striGNSprop}.)

Recall that a linear map $T\colon\cA\rightarrow\cB$ is \phantomsection\hl{completely positive}\index{completely positive|see {CP-map}}\index{CP-map} (we also say $T$ is a \hl{CP-map}) if
\beq{\label{CPeq}
\sum_{i,j}b_i^*T(a_i^*a_j)b_j
~\ge~
0
}\eeqn
for all finite choices of $a_i\in\cA,b_i\in\cB$. If $\cA$ is unital, then there exists a pair $(E,\xi)$, the so-called \hl{GNS-construction}\index{CP-map!GNS-construction|bf}\index{GNS-construction|bf} of $T$ (Paschke \cite{Pas73}), consisting of a \hl{GNS-correspondence} $E$ from $\cA$ to $\cB$ and a \hl{cyclic vector}\index{GNS-construction!cyclic vector}\index{cyclic vector|see {GNS-construction}} $\xi\in E$ such that
\bal{\label{GNSeq}
T(a)
&
~=~
\AB{\xi,a\xi},
&
\cls\cA\xi\cB
&
~=~
E.
}\eal
We denote the situation in \eqref{GNSeq} as $(E,\xi)=$ GNS-$T$.
%%%% BO
The GNS-correspondence of $T$ can be constructed explicitly by considering the $\cA$-$\cB$ bimodule $\cA \otimes \cB$, endowing it with the semi-inner product
\beqn{
\bAB{a \otimes b, a' \otimes b'} := b^* T(a^*a') b' , 
}\eeqn
quotienting out the kernel, and completing to obtain $E$.  The image of $\U_\cA \otimes \U_\cB$ in $E$ is, then, the required cyclic vector $\xi$. 
%%%% EO 
(If $\cA$ is nonunital, then either we have to put stronger conditions on $T$ (for instance, $T$ extends to a strict CP-map $M(\cA)\rightarrow M(\cB)$), or we get weaker statements (for instance, $\xi\notin\cls\cA\xi\cB$); see, for instance, Skeide \cite[Section 4.1]{Ske01}.) Conversely, if $F$ is a correspondence from $\cA$ to $\cB$ and if $\zeta$ is in $F$, then $S=\AB{\zeta,\bullet\zeta}$ defines a CP-map $S$. If $S=T$, \phantomsection\index{minimal!GNS-construction}then $\xi\mapsto\zeta$ extends as a (unique) bilinear \hl{isometry}\index{isometry}\index{operator!isometry} (that is, an inner product preserving map) from $E$ to $F$. In particular, if also $(F,\zeta)$ is a GNS-construction for $T$, then there is a unique cyclic vector intertwining isomorphism from $(E,\xi)$ to $(F,\zeta)$. In other words, the GNS-construction is unique up to such an isomorphism.%
\footnote{ \label{StineFN}
Recall that every Hilbert \nbd{\cB}module $E$ (even a pre-Hilbert module over a pre-\nbd{C^*}algebra) may be transformed into a module of operators\phantomsection\index{operator module}\index{module!operator} in the following way: Identify $\cB$ as a concrete operator algebra $\cB\subset\sB(G)$ of operators on a Hilbert space $G$; form the Hilbert space $H:=E\odot G$; for $x\in E$ define the operator 
%%%% BO 
%  $L_x\odot\id_G\colon g\mapsto x\odot g$ in $\sB(G,H)$. 
$L_x = x\odot\id_G\colon g\mapsto x\odot g$ in $\sB(G,H)$. 
%%%% EO
Then the subspace $L_E\subset\sB(G,H)$ satisfies $L_E\cB\subset L_E$, $L_E^*L_E\subset\cB$, and $\sscls L_EG=H$, that is, $L_E$ is a concrete (pre-)Hilbert \nbd{\cB}module\index{concrete!Hilbert module@(pre-)Hilbert module}. Moreover, $L_{xb}=L_xb$ and $L_x^*L_y=x^*y=\AB{x,y}$ so that $x\mapsto L_x$ is an isomorphism of (pre-)Hilbert \nbd{\cB}modules. Moreover, if $L'\colon x\mapsto L'_x\in\sB(G,H')$ is another isomorphism such that $L'^*_xL'_y=\AB{x,y}$ and $\sscls L'_EG=H'$, then $x\odot g\mapsto L'_xg$ defines a unitary $u\colon H\rightarrow H'$ such that $uL_x=L'_x$. Therefore, as soon as the identification $\cB\subset\sB(G)$ is chosen (for instance, if $\cB$ is a von Neumann algebra), there is nothing arbitrary in the identification of $E=L_E\subset\sB(G,H)$.

Note that $\sB^a(E)=\sB^a(E)\odot\id_G\subset\sB(H)$. Therefore, if $E$ is an \nbd{\cA}\nbd{\cB}correspondence, then $H$ inherits a representation $\rho$ of $\cA$. If $\xi$ is a vector in $E$ 
%%%% BO 
% (for instance, if  $(E,\xi)=$ GNS-$T$), then $L_\xi^*\rho(a)L_\xi=T(a)$, that is, we recover the \it{Stinespring construction} \cite{Sti55}.  
then $L_\xi^*\rho(a)L_\xi=\AB{\xi, a \xi}$. In particular if  $(E,\xi)=$ GNS-$T$, then $L_\xi^*\rho(a)L_\xi=T(a)$, that is, we recover the \phantomsection\hl{Stinespring construction}\index{Stinespring construction}\index{minimal!Stinespring construction} \cite{Sti55}.  
%%%% EO
However, in the context of CP-semigroups, the Stinespring construction is not even approximately as useful as the GNS-construction. Barreto, Bhat, Liebscher, and Skeide \cite[2.1.7. Functoriality]{BBLS04}: ``A [...] Hilbert \nbd{\cA}\nbd{\cB}module $E$ is a \phantomsection\hl{functor}\index{correspondence!as functor} sending (non-degenerate) representations of $\cB$ on $F$  to (non-degenerate) representations of $\cA$ on $E\odot F$,, and the composition of two such functors is the tensor product. The Stinespring construction is a dead end for this functoriality.''

In the von Neumann case, the full information about a (GNS or not) correspondence $E$ can, however, be reconstructed from the \it{Stinespring representation} $\rho$, if we add as second input a representation of $\cB'$, the so-called  \it{commutant lifting} $\rho'$ defined by $\rho'(b')=\id_E\odot b'$. We will explain this in Appendix \ref{vNAPP}.
}

% \newpage
\lf\noindent
\bf{Unusual conventions and ways of writing.~} \phantomsection\label{IUC}\index{conventions!unusual}
Here we collect some conventions, in the wide sense, that will be applied without mention throughout these notes, and that might puzzle the reader when ignoring them.
\begin{itemize}
\item
\hl{Observations and remarks:}\index{observations and remarks!how to read them} The statements made in an \hl{observation} have the rank of a proposition, but in an observation both statement and proof are incorporated in a single bit of text. Like propositions, observations may be important in their own right (beyond these notes) or may communicate facts that are necessary where they stand or elsewhere in these notes.

A \hl{remark}  furnishes additional information that is \bf{not} logically needed where it stands. A remark may try to guide reader's intuition by adding mathematical or historical information, in particular, information about things that might be known to the reader but do not occur here. It is up to the reader to decide whether it helps, but logically a remark can be ignored where it stands. It is, however, possible that a remark might be referred to in other places of these notes.

\item
By a \hl{(co)restriction}\index{corestriction@(co)restriction} of a map $f\colon A\rightarrow B$ we mean a map $A'\rightarrow B'$ for subsets $A'\subset A$ and $B'\subset B$ satisfying $f(A')\subset B'$, obtained by $A'\ni a'\mapsto f(a')\in B'$. This should not be confused with formulations containing several parenthetical insertions that either are all there or are all not there to give two different (both -- hopefully -- correct) meanings to a sentence. In that case, we would write \it{(co-)restriction}. (Example: The adjointable contraction $a\colon E\rightarrow F$ is a (co-)isometry if $a^*$ possesses a (co-)restriction to a unitary; see also Appendix \ref{coisoAPP}. Here, we do not mean that both domain and codomain can be made smaller (for instance, $\zero$ for both is always possible), but that in the first case (without (co-)) the domain $F$ of $a^*$ has to be made smaller (namely, $aE$), while in the second case (with (co-)) the codomain $E$ of $a^*$ has to be made smaller (namely, $(\ker a)^\perp$). See also the use of (super)(sub)product subsystems in Definition \ref{subdefi}; and see also the following item.)

\item
For us, an \hl{approximate unit}\index{approximate unit!in a pre-\nbd{C^*}algebra} in a pre-\nbd{C^*}algebra is a net $\bfam{u_\lambda}_{\lambda\in\Lambda}$ such that $\lim_\lambda u_\lambda a=a=\lim_\lambda au_\lambda$ for all $a\in\cA$. An approximate unit may be \hl{self-adjoint}, that is, $u_\lambda=u_\lambda^*$ (in which case it is sufficient to check only one of the two limits for all $a$); it may be \hl{positive}, that is, $u_\lambda\ge0$; it may be \hl{bounded} ((\hl{strictly}) \hl{contractive}), that is, $\norm{u_\lambda}\le M$ for some $M>0$ ($\norm{u_\lambda}(<)\le 1$).

\item
Instead of \phantomsection\it{\nbd{W^*}modules}\index{von Neumann!W@\nbd{W^*}} and \it{\nbd{W^*}correspondences}, we consequently use \hl{von Neumann modules}\index{von Neumann!module} and \hl{von Neumann correspondences}\index{von Neumann!correspondence}. Despite the latter appear technically simpler -- or, better, the algebraic advantages of the notion lead (more) automatically to the solutions of technical problems --, the former appear still much more common. Therefore, we give a brief introduction to von Neumann modules and correspondences in Appendix \ref{vNAPP}.

\item
Last but not least, the terms defined in definitions are typeset in \hl{boldface italics}. A word typeset in \bf{boldface}, means just to emphasize that word.
\end{itemize}

\lf\lf
\noindent
In this revised version we include and index and, as Appendix \ref{rmAPP}, a \it{Road Map} through these notes.

\lf\lf\noindent
\bf{Acknowledgments.~}
After some loose e-mail discussions starting with the Fields Workshop ``Advances in Quantum Dynamics'' 2007, we started working seriously on this pro\-ject when we met, again in Canada, in November 2009 in Waterloo, where OS was spending his first year as postdoc, and in January 2010 in Kingston, where MS was spending the first four months of his sabbatical. We are deeply grateful to our hosts, Ken Davidson in Waterloo and Roland Speicher in Kingston, who have made these encounters possible. MS is grateful to the math departments of Ben-Gurion University and the Technion for having him hosted numerous times, as well as for (partial) economic support from his department of economy. MS is also grateful to OS for having made this possible, and for having taken on once the adventure of a visit in Campobasso. OS is grateful to MS for the wonderful hospitality granted during that adventure. OS is also grateful for partial support granted by ISF Grants no. 195/16 and 431/20.

Last but not least, both authors wish to express their gratitude to the two anonymous referees, for kindly taking the time to give this manuscript a careful reading.

\newpage

\section{CP-Semigroups over monoids}\label{monoSEC}

In these notes we are interested in dilations of CP-semigroups $T=\bfam{T_t}_{t\in\bS}$ on unital \nbd{C^*}algebras $\cB$ that are indexed by a monoid $\bS$. While our primary interest is to construct such dilations, it is also indispensable to find out which necessary conditions on the CP-semigroup arise from the assumption that $T$ \it{does} admit a dilation. For instance, when $T_t=\ep\circ\vt_t\circ i$ (as in Diagram \eqref{1stdil}), it is forced that all $T_t$ are contractions.

CP-Semigroups come shipped with \it{subproduct systems}; \nbd{E}semigroups on $\sB^a(E)$ come shipped with \it{product systems}; in the construction of a dilating \nbd{E}semigroup the problem to \it{transform} the subproduct system into a product system, plays a crucial role. We shall investigate these topics and how they have to be modified in this generality in the subsequent sections. In the present section, fixing also some notations, we discuss the concept of \it{dilation} briefly in general. Soon, we restrict our attention to the so-called \it{weak} dilations. Among the weak dilations we put emphasis on the so-called \it{strong} ones. Examples for non-strong weak dilations will haunt us throughout these notes. Only in Section \ref{DilSPSpSEC}, it turns out that strong dilations come along with a \it{superproduct system} containing the subproduct system of the CP-semigroup, and only then we define \it{good} dilations as dilations sharing that property; good dilations sit strictly in between strong and weak ones. The superproduct system of a dilation need not be a product system; the question whether or not it embeds into a product system is key.

\lf
Let us fix a monoid $\bS$. Recall (see the conventions in the introduction) that a \phantomsection\hl{semigroup over $\bS$ on a set $\cB$}\index{semigroup!over a monoid} is a \hl{monoid map} $t\mapsto T_t$ from $\bS$ into the maps on $\cB$, that is, $T_{st}=T_s\circ T_t$ and $T_0=\id_\cB$. A semigroup over $\bS$ on a \nbd{C^*}algebra $\cB$ is a \hl{CP-semigroup}\index{semigroup!CP-|bf}\index{CP-semigroup!|bf} if all $T_t$ are completely positive (CP) maps.

\bemp[Convention.~] \label{ucCONV}\index{conventions!numbered!contractive CP-semigroups@\ref{ucCONV} (contractive CP-semigroups)}\index{CP-semigroup!contractive}
We shall assume that all our CP-semigroups act on \hl{unital} \nbd{C^*}algebras. We shall assume that all our CP-semigroups are \hl{contractive}.
\eemp

The latter is referred to by Bhat \cite{Bha96} as a \it{quantum dynamical semigroup}\index{semigroup!quantum dynamical}. We said already that CP-semigroups possessing dilations, have to be contractive. Also the following classes are contractive automatically.

\bdefi \label{spCPdefi}
A CP-semigroup $T=\bfam{T_t}_{t\in\bS}$ is
\begin{enumerate}
\item
a \hl{Markov semigroup}\index{semigroup!Markov}\index{Markov semigroup} if $T$ is \hl{unital}, that is, if $T_t(\U)=\U$ for all $t\in\bS$;

\item
an \hl{\nbd{E}semigroup}\index{semigroup!E@\nbd{E}}\index{E-semigroup@\nbd{E}semigroup} if all $T_t$ are endomorphisms;

\item
an \hl{\nbd{E_0}semigroup}\index{semigroup!E0@\nbd{E_0}}\index{E0-semigroup@\nbd{E_0}semigroup} if it is a Markov semigroup and an \nbd{E}semigroup. 
\end{enumerate}
\edefi
Even if $\bS$ is a topological monoid, we shall \bf{not} assume any continuity condition with $t\in\bS$, unless explicitly mentioned. By homomorphism, representation, and so forth, of \nbd{C^*}algebras, we \bf{always} mean \nbd{*}homomorphism, \nbd{*}representation, and so forth; if not, we shall say algebra homomorphism, algebra representation, and so forth.

\bdefi \label{dildef}
Let $T$ be a CP-semigroup over $\bS$ on $\cB$. By a \hl{dilation}\index{dilation, general} of $T$ we understand a quadruple $(\cA,\theta,i,\ep)$ consisting of a unital \nbd{C^*}algebra $\cA$, an \nbd{E}semigroup $\theta$ over $\bS$ on $\cA$, an embedding $i\colon\cB\rightarrow\cA$, and an \hl{expectation}\index{expectation} $\ep\colon\cA\rightarrow\cB$ (that is, the map $i\circ\ep$ is a conditional expectation onto $i(\cB)$), such that the diagram
\beq{ \label{basdil}
\parbox{10cm}{\xymatrix{
\cB	\ar[d]_i	\ar[rr]^{T_t}	&&	\cB
\\
\cA				\ar[rr]_{\theta_t}	&&	\cA	\ar[u]_\ep
}}
}\eeq
commutes for all $t$.

More generally, we shall say $(\cA,\theta,i,\ep)$ is a \hl{dilation}, whenever the maps $T_t:=\ep\circ\theta_t\circ i$ define a semigroup (necessarily CP) on $\cB$. (It is very convenient to be able to say `dilation' without having to indicate in advance the dilated semigroup.) We speak of a \hl{Markov dilation}\index{dilation, general!Markov}\index{dilation, weak!Markov} if the $T_t$ form a Markov semigroup.%
\footnote{
Be aware that several authors, following K\"{u}mmerer \cite{Kuem85}, mean by \it{Markov dilation} something considerably more restrictive.
\vspace{1ex}
}

A dilation is \hl{unital}\index{dilation, general!unital} if $i$ is unital.

A dilation is \hl{weak}\index{dilation, weak}\index{weak dilation} if $i(\cB)$ is a \index{corner}\hl{corner} in $\cA$ (that is, $i(\cB)=p\cA p$, where $p:=i(\U_\cB)$, so that $pap=i\circ\ep(pap)=p(i\circ\ep(a))p=i\circ\ep(a)$).

A weak dilation is \hl{full}\index{dilation, weak!full}\index{full dilation}, if the corner $i(\cB)$ is \hl{strictly full}, that is, if $\cA$ is (isomorphic to) the multiplier algebra of the ideal \,$\ls\cA p\cA$. (By this, we mean that the canonical homomorphism $\cA\rightarrow M(\ls\cA p\cA)$ is an isomorphism. See the proof of Proposition \ref{convprop} for an explanation.)

A dilation is \hl{reversible}\index{dilation, general!reversible} if $\theta$ is an automorphism semigroup.

A dilation is \hl{semireversible}\index{dilation, general!semireversible} if $\theta$ is an injective \nbd{E_0}semigroup.%
\footnote{Frequently, a semireversible dilation may be extended to a reversible dilation on a containing $\sB(H)\supset\cA$; see Remark \ref{dilhistrem}. Surely, for a reversible extension to exist, being semireversible is a necessary condition.}
\edefi

\bemp[Convention.~]\label{subalgconv}\index{conventions!numbered!weak dilations@\ref{subalgconv} (weak dilations)}
When we consider dilations to general $\cA$, then we usually shall assume that $\cB\subset\cA$ and that $i$ is the canonical embedding, so that $\ep$ is a usual conditional expectation onto $\cB$.
\eemp

After Convention \ref{subalgconv}, we shall denote a weak dilation as $(\cA,\theta,p)$ so that $\cB:=p\cA p\subset\cA$ and $i$ is the canonical injection.

\vspace{-.2ex}
\bemp[Convention.~]\label{Ba(E)conv}\index{conventions!numbered!module dilations@\ref{Ba(E)conv} (module dilations)}
When we consider dilations to $\cA=\sB^a(E)$ where $E$ is a Hilbert  \nbd{\cB}module, then we will usually have a \hl{unit vector}\index{unit vector|bf} $\xi\in E$ (that is, $\AB{\xi,\xi}=\U_\cB$) such that $\ep=\AB{\xi,\bullet\xi}$ and such that $i=\xi\bullet\xi^*\colon b\mapsto\xi b\xi^*$. It follows that $p=i(\U_\cB)=\xi\xi^*$. In this situation, we also shall use the form $\vt$ of the letter theta for the \nbd{E}semigroup, instead of $\theta$.
\eemp

In the situation of Convention \ref{Ba(E)conv}, we shall write the dilation as $(E,\vt,\xi)$ and refer to it as \hl{module dilation}\index{dilation, weak!module}\index{module dilation}.

It is important to observe that \index{dilation, weak!module!$\cong$ full dilation}\index{module dilation!$\cong$ full dilation}\it{module dilation} and \it{full dilation} are essentially the same thing. Of course, a module dilation is full. ($\sB^a(E)p\sB^a(E)=EE^*$ and the multiplier algebra of $\ls EE^*=\sF(E)$ is $\sB^a(E)$.) But also the converse is true in the following sense:

\bprop\label{convprop}{~}

\begin{enumerate}
\item \label{cp1}
If $p$ is a projection in a \nbd{C^*}algebra $\cA$, then $E:=\cA p$ with inner product $\AB{ap,a'p}:=pa^*a'p$ is a Hilbert module over $\cB:=p\cA p$. Moreover, $\xi:=p$ is a unit vector in $E$ and the action $a\colon a'p\mapsto aa'p$ defines a homomorphism $\cA\rightarrow\sB^a(E)$, the \hl{canonical homomorphism}.

\item \label{cp2}
If $\cC$ is a pre-\nbd{C^*}algebra, then the \hl{multiplier algebra}\index{multiplier algebra!of a pre-\nbd{C^*}algebra}
\beqn{
M(\cC)
~:=~
\bset[(L,R)]{L,R\in\sL(\cC);aL(b)=R(a)b}
}\eeqn
with `the usual operations', coincides with $\sL^a(\cC)$ (considering $\cC$ as a pre-Hilbert module over itself).

\item \label{cp3}
If $E$ is a Hilbert \nbd{\cB}module, then $M(\sF(E))=M(\sK(E))=\sB^a(E)$.
\end{enumerate}
\noindent
In conclusion, if $(\cA,\theta,p)$ is a full weak dilation, then it is conjugate to the module dilation $(E,\vt,\xi)$ under the canonical isomorphism $\cA\rightarrow\sB^a(E)$.
\eprop

\proof
\eqref{cp1}
The only question with not entirely obvious answer is why $E$ is complete. Since $E\subset\cA$ and $\cA$ is complete, we get $\ol{E}\subset\cA$, so that $E=\cA p\supset\ol{E}p=\ol{E}\supset E$, so, $\ol{E}=E$.

\eqref{cp2}
The double centralizers $(L,R)$ in $M(\cC)$ are equipped with the usual operations known from the \nbd{C^*}case $\cC=\ol{\cC}$, just forgetting everything about `norm' and `bounded'. Then the `forgetful map' $(L,R)\mapsto L$ defines a homomorphism into $\sL^a(\cC)$. Likewise, the map $a\mapsto (L_a,R_a)$ where $L_a:=a\bullet$ and $R_a:=\bullet a$, defines a homomorphism from $\sL^a(\cC)$ into $M(\cC)$. One easily checks that the two are inverses of each other. (See Skeide \cite[Lemma 1.7.10]{Ske01}, which is promoted here from bounded double centralizers to arbitrary double centralizers.)

\eqref{cp3}
By \eqref{cp2}, we have $M(\sF(E))=\sL^a(\sF(E))$. Since $E$ is assumed complete, we have $E= EE^*E=\ls EE^*E=\sF(E)\,\ul{\odot}\,E$ and, of course, $\sF(E)=E\,\ul{\odot}\,E^*$. 
%%%% BO 
(Here and elsewhere we write $\ul{\odot}$ for the algebraic version of the internal tensor, where no completion is carried out; see Skeide \cite[Appendix]{Ske12p}) 
%%%% EO
We get $\sL^a(\sF(E))=\sL^a(E)$ via $a\mapsto a\odot\id_E$. So, $M(\sF(E))=\sL^a(E)$. Again, since $E$ is complete, $\sL^a(E)=\sB^a(E)$. On the other hand, the isomorphism $a\mapsto a\odot\id_{E^*}$ from $\sL^a(E)=\sB^a(E)$ onto $M(\sF(E))$, clearly, maps into the bounded double centralizers, which, therefore, extend to elements in $M(\sK(E))$. (This proof of Kasparov's result $M(\sK(E))=\sB^a(E)$ \cite{Kas80} follows very much the proof of \cite[Corollary 1.7.14]{Ske01}, but is, maybe, yet another bit more stringent.)\qed

\lf
The property of a weak dilation being full is an intrinsic way to find out if a weak dilation is a module dilation. For a Markov semigroup, existence of module dilations is almost equivalent to existence of a product system containing its subproduct system; see Example \ref{EPSex}. It depends, however, on a topological question, namely the question whether the dilating \nbd{E}semigroup consists of strict endomorphisms. Both strictness of endomorphisms of $\cA$ as multiplier algebra and strictness of endomorphisms of $\sB^a(E)$ make sense. It is the nice property of multiplier algebras that apart from being strict completions (of $\ls\cA p\cA$ and of $\sF(E)$, respectively), they also may be captured purely algebraically in terms of double centralizers. This is why we can allow ourselves to discuss the notion of full dilation already here, while for the time being we still ignore topological questions. They will enter and be attacked starting from Sections \ref{topSEC}. As in these notes we intend to study dilations in terms of product systems, it will be a vital question to find out whether existence of a dilation guarantees existence of a module (that is, of a full) dilation.

While in these notes we shall restrict our attention to weak dilations, the following historical remark tries to explain, that weak dilations also are important to get other types of dilations, and to understand the latter.

\brem \label{dilhistrem}
The setting of unital reversible dilations of a (necessarily!) Markov one-para\-meter semigroup is what one wishes in quantum probability. Many dilations of such type have been obtained with the help of \it{quantum stochastic calculus} on Fock type objects (see, for instance, Hudson and Parthasarathy \cite{HuPa84}, K\"{u}mmerer and Speicher \cite{KueSp92}, Goswami and Sinha \cite{GoSi99}, Skeide \cite{Ske00}) or in more algebraic situations (see, for instance, Accardi, Fagnola, and Quae\-gebeur \cite{AFQ92}, Hellmich, K\"{o}stler, and K\"{u}mmerer \cite{HKK98,Koes00}). These dilations are all obtained as a perturbation of a so-called \it{noise} (a dilation of the trivial semigroup $T_t=\id_\cB$, as defined and examined in Skeide \cite{Ske06d}) by a \it{unitary cocycle}. It is well known that all semireversible one-parameter \nbd{E_0}semigroups can be extended to automorphism semigroups; see Arveson and Kishimoto \cite{ArKi92} for von Neumann algebras (or a new proof in Skeide \cite[Theorem B.36]{Ske16}), and Skeide \cite[Corollary 3.2]{Ske11a} for (not necessarily unital) \nbd{C^*}algebras (where \nbd{E_0}semigroup means that $\theta_t(\cA)$ acts nondegenerately on $\cA$). In fact, all known unital reversible dilations can be thought of as obtained in that way.

Weak dilation is a sort of minimum requirement a dilation must fulfill. We do not know of a dilation constructed as a cocycle perturbation of a noise that did not have sitting inside also a weak dilation. In fact, Skeide \cite{Ske16} has shown that a Markov one-parameter semigroup admits a unital reversible dilation that is the cocycle perturbation of a noise if and only if that Markov semigroup is \it{spatial}. (Spatiality is a property of the order structure of the set of CP-semigroups \it{dominated} by the Markov semigroup. It has been defined by Arveson \cite{Arv97a} for Markov semigroups on $\sB(H)$ and it has been generalized in Bhat, Liebscher, and  Skeide \cite{BLS10}, while the definition by Powers \cite{Pow04} for $\sB(H)$ is considerably more restrictive.) The construction of a weak dilation plays a crucial role in the proof of that result.

Non-Markov CP-semigroups and their dilations do occur in problems of single and multivariate operator theory. See Sections \ref{EXwnsSEC} and \ref{EXBexSEC} (in particular, \ref{EXBexSEC}\ref{KrowSSEC}).
\erem

\lf
Let us return to CP-semigroups and their weak dilations. Markov semigroups form a subclass with many particularly nice properties. Some results we have for Markov semigroups are not true for CP-semigroups, or at least have proofs that run considerably less smoothly.

A powerful tool to reduce problems about general CP-semigroups to results about Markov semigroups, is the so-called \it{unitalization} of a CP-semigroup, to be discussed in the last part of this section after Proposition \ref{strequivprop}. But first we discuss a nice extra property of Markov dilations (Equation \eqref{strdileq}, below) -- an extra property that turns out to be crucial in order to be able to understand dilations of non-Markov semigroups via unitalization.

Recall that, by Convention \ref{subalgconv}, we assume the situation where $\cB\subset\cA$. Speaking about weak dilations, we shall emphasize the projection $p\in\cA$, putting $\cB:=p\cA p$ and $\ep:=p\bullet p$, and write dilations as triples $(\cA,\theta,p)$.

Recall that a projection $p\in\cA$ is \hl{increasing}\index{projection!increasing} for $\theta$, if $\theta_t(p)\ge p$ for all $t\in\bS$. The following result is \it{folklore}, but we include a proof for convenience -- rather to illustrate how our terminology and conventions are applied, than to do something new.

\bprop\label{Mcharprop}
Let $(\cA,\theta,p)$ be a unital \nbd{C^*}algebra $\cA$, an \nbd{E}semigroup $\theta$ on $\cA$, and a projection $p\in\cA$, and put $T_t(pap):=p\theta_t(pap)p$. Then the following are equivalent:
\begin{enumerate}
\item\label{id1}
$p$ is increasing.

\item\label{id2}
The maps $T_t$ are unital.

\item\label{id3}
$(\cA,\theta,p)$ is a weak Markov dilation.
\end{enumerate}
Moreover, under any of the conditions we have
\beq{\label{strdileq}
p\theta_t(a)p
~=~
T_t(pap)
}\eeq
for all $a\in\cA$.
\eprop

\proof
Inserting $p=\U_\cB$ into $T_t$, we see that \eqref{id1} and \eqref{id2} are equivalent.

Clearly, if \eqref{id1} is true, then \eqref{strdileq} holds. By iterating \eqref{strdileq}, we see that the $T_t$ form a semigroup, which, by \eqref{id2}, is Markov, that is, we have \eqref{id3}. And, by definition, \eqref{id3} implies \eqref{id2}.\qed

\lf
The bonus property of weak Markov dilations manifested in \eqref{strdileq}, is illustrated in the diagram
\beq{ \label{strodil}
\parbox{10cm}{\xymatrix{
\cB	\ar[rr]^{T_t}							&&\cB
\\
\cA	\ar[rr]_{\theta_t}		\ar[u]^{\ep}	&&\cA	\ar[u]_{\ep}
}}
}\eeq
which might be called a \it{coextension} of $T$. (An \it{extension} would by the ``dual'' diagram with two times $i$ instead of two times $\ep$; of course, an extension of a non-\nbd{E}semigroup can never be an \nbd{E}semigroup.) The property is, clearly, stronger than the dilation diagram in \eqref{basdil}, so we call a dilation fulfilling the diagram in \eqref{strodil} a \hl{strong dilation}\index{strong dilation}\index{dilation, weak!strong}\index{dilation, general!strong}.

\bob \label{semgenob}
Suppose we have  a quadruple $(\cA,\theta,i,\ep)$ such that the strongness condition in Diagram \eqref{strodil} holds for all $t$ from a subset $G\subset\bS$ that generates $\bS$ as a monoid (for instance, $G=\bS$). Then $(\cA,\theta,i,\ep)$ is a strong dilation. (Indeed, $\ep\circ\theta_t=\ep\circ\theta_{g_1}\circ\ldots\circ\theta_{g_n}=T_{g_1}\circ\ldots\circ T_{g_n}\circ\ep=T_t\circ\ep$.)
\eob

Note that this neither assumes that $T$ is Markov nor that the dilation is weak.  However, since in these notes we are exclusively interested in weak dilations, and in order to avoid constructs such as \it{strong weak dilations}, as a \bf{convention}\index{conventions!strong dilations are assumed weak}\index{strong dilation!is weak (convention)}, we always shall understand (in these notes) by a strong dilation a strong dilation that is also a weak dilation.

\brem \label{regdilrem}
Strong dilations in that sense (being also weak) are those that occurred in larger parts of the literature, especially, in multivariate operator theory where non-Markov semigroups do play a role. They also have occurred under the name of \it{\nbd{E}dilations}; we avoid that name.

Markov dilations are strong. But since we are interested also in applications to multivariate operator theory, we have to distinguish carefully between strong and weak dilations. In Section \ref{EXwnsSEC} (and, further, in Subsection \ref{EXBexSEC}\ref{KrowSSEC}), we discuss examples of weak dilations that are not strong, coming from classical dilation theory of operators. Example \ref{E0weakex} below, is an important appetizer.
\erem

\bob \label{E0strongob}
If $(\cA,\theta,p)$ is a strong dilation and $\theta$ an \nbd{E_0}semigroup, then $p$ is increasing. We shall call a dilation where $\theta$ is an \nbd{E_0}semigroup, an \hl{\nbd{E_0}dilation}\index{E0dilation@\nbd{E_0}dilation}\index{dilation, weak!E0@\nbd{E_0}}\index{dilation, general!E0@\nbd{E_0}}. So, CP-semigroups admitting a strong \nbd{E_0}dilation are necessarily Markov.
\eob

The following example illustrates that the condition to be strong cannot be dropped.

\bex \label{E0weakex}
Let $u_t$ denote the unitary right shift group on $H:=L^2(\R)$ defined by setting $\SB{u_tf}(x)=f(x-t)$, and define the \nbd{E_0}semigroup $\vt$ on $\sB(H)$ by setting $\vt_t:=u_t\bullet u_t^*$. Define $p\in\sB(H)$ to be the projection onto $G:=L^2(\R_+)$. Then $v_t:=pu_tp=u_tp$ is the isometric right shift semigroup on $G$. We have
\beqn{
p\vt_t(p\bullet p)p
~:=~
(pu_tp)\bullet(pu_tp)^*
~=~
(u_tp)\bullet(u_tp)^*
~=~
\vt_t(p\bullet p).
}\eeqn
It follows that the maps $T_t:=p\vt_t(p\bullet p)p=v_t\bullet v_t^*$ are homomorphisms and form an \nbd{E}semigroup $T$ on $p\sB(H)p=\sB(G)$ which is not an \nbd{E_0}semigroup; in particular, the CP-semigroup $T$ dilated by the \nbd{E_0}semigroup $\vt$, is not a Markov semigroup. By Observation \ref{E0strongob}, this \nbd{E_0}dilation is not strong.

Suppose, on the one hand, that (as discussed in the Example-Section \ref{EXwnsSEC}) the $w_t:=v_t^*$ form a coisometric dilation of a (proper) contraction semigroup $c_t$, so that there is a projection $q\in\sB(G)\subset\sB(H)$ such that $qv_t^*q=c_t$. (For instance, take $q=ff^*$, where the unit vector $f$ is the function $x\mapsto e^{-x}$, so that $qv_t^*q=qe^{-t}$.) Then $(\sB(H),\vt,q)$ is an \nbd{E_0}dilation of the (non-Markov!) CP-semigroup $S_t:=c_t^*\bullet c_t$ (by Observation \ref{E0strongob}, necessarily a weak dilation) which \it{compresses} (by $p$) to the (strong!) non-$E_0$ dilation $(\sB(G),T,q)$ of $S$. (See Subsection \ref{EXBexSEC}\ref{BexSSEC} and the part about incompressible dilations in Section \ref{minSEC}.)

Adding, on the other hand, $\C$ to $H$ and defining $\wh{q}$ to be the projection onto that $\C$, we get by $(\sB(H\oplus\C),((u_t+\wh{q})\bullet(u_t+\wh{q})^*),\wh{q})$ an \nbd{E_0}dilation of the Markov semigroup $\wh{S}_t=\id_\C$ that compresses via $p+\wh{q}$ to a dilation $(\sB(G\oplus\C),((v_t+\wh{q})\bullet(v_t+\wh{q})^*),\wh{q})$ that is not an \nbd{E_0}dilation.
\eex

\lf
We may ask, if there is a similar characterization of weak and of strong dilations of general CP-semigroups in terms of the projection $p$ as for Markov semigroups in Proposition \ref{Mcharprop}. Well, for weak dilations this is hoping for too much. Bhat has a description of general weak dilations in terms of a pair of decreasing projections in the case of one-parameter semigroups and $\cA=\sB(H)$; see \cite[Corollary 2.3]{Bha02}. (The proof of that result depends on existence of projections onto closed subspaces and does not fit the representation free framework we discuss here. At best, it might be provable for von Neumann algebras.) However, at least for strong dilations we can say a bit more.

Recall that Bhat \cite{Bha03} calls (in the one-parameter case and for $\cA=\sB(H)$) a weak dilation $(\cA,\theta,p)$ \hl{regular}\index{dilation, weak!regular} if $\theta_t(\U-p)\le\U-p$ for all $t$. Arveson \cite{Arv03} called such a projection \it{coinvariant}; we prefer to call it \hl{coincreasing}\index{projection!coincreasing (coinvariant)}. (We would expect that also the following proposition is \it{folklore} and might be found somewhere in the works of Bhat and of Arveson.)

\bprop \label{strequivprop}
Under the same hypotheses as in Proposition \ref{Mcharprop}. Equivalent are:
\begin{enumerate}
\item \label{sp1}
$\theta_t(\U-p)p=0$ ~for all~ $t\in\bS$.

\item \label{sp2}
$(\cA,\theta,p)$ is a strong dilation.

\item \label{sp3}
$(\cA,\theta,p)$ is a regular dilation.

\end{enumerate}
\eprop

\proof
Like in the proof of Proposition \ref{Mcharprop}, a key problem is to show that $\theta_t(\U-p)p=0$. Indeed, that latter equation implies both that $p\theta_t(a)p=p\theta_t(pap)p$ and that $(\U-p)\theta_t(\U-p)(\U-p)=\theta_t(\U-p)-0-0+0=\theta_t(\U-p)$ so that $\U-p$ is decreasing. Moreover, from \eqref{strdileq}, like before, it follows that the $T_t$ form a semigroup.

On the other hand, if $(\cA,\theta,p)$ is a strong dilation, then $p\theta_t(\U-p)p=p\theta_t(p(\U-p)p)p=0$ implies $\theta_t(\U-p)p=0$.

Likewise, $\theta_t(\U-p)\le\U-p$ implies $0\le p\theta_t(\U-p)p\le p(\U-p)p=0$ so that, again, $\theta_t(\U-p)p=0$.\qed

\lf
We now discuss \phantomsection\label{ulize}how unitalization helps to reduce problems regarding strong dilations of CP-semigroups to those regarding Markov semigroups (where we know all dilations are strong).

If $\cB$ is a unital \nbd{C^*}algebra, then we define its \phantomsection\hl{unitalization}\index{unitalization!of a (unital) \nbd{C^*}algebra} $\wt{\cB}=\cB\oplus\C\wt{\U}$; that is, we add a new unit $\wt{\U}$. The old unit $\U$ is, of course, different from the new one; in fact, excluding the trivial case $\U=0$, both $\U$ and $\wt{\U}-\U$ are nontrivial central projections. Algebraically, $\wt{\cB}$ is isomorphic to the \nbd{*}algebraic direct sum $\C\oplus\cB$ via $b+\lambda\wt{\U}\mapsto(\lambda,b+\lambda\U)$, respectively, $(\lambda,b)\mapsto b+\lambda(\wt{\U}-\U)$. From $\C\oplus\cB$ we see how to norm $\wt{\cB}$ as a \nbd{C^*}algebra, while from $\wt{\cB}=\cB\oplus\C\wt{\U}$ it is easy to define the unitalization of maps and to understand their properties.

In fact, if $T$ is a linear map on $\cB$, we define its \hl{unitalization}\index{unitalization!of a map}\index{CP-map!unitalization of} as $\wt{T}(b+\lambda\wt{\U}):=T(b)+\lambda\wt{\U}$. From this definition it is immediate that $\wt{\big.S\circ T\big.}=\wt{S}\circ\wt{T}$. In particular, if $T=\bfam{T_t}_{t\in\bS}$ is a semigroup on $\cB$, then $\wt{T}=\bfam{\wt{T}_t}_{t\in\bS}$ is a semigroup\index{unitalization!of a semigroup}\index{semigroup!CP-!unitalization of}\index{CP-semigroup!unitalization of} of unital maps on $\wt{\cB}$. In the picture, $\C\oplus\cB$ this looks
\vspace{-1.5ex}
\bmun{
~~~~~~~~~
\wt{T}_t(\lambda,b)
~=~
\wt{T}_t(b+\lambda(\wt{\U}-\U))
~=~
T_t(b-\lambda\U)+\lambda\wt{\U}
\\[-.5ex]
~=~
\bfam{\lambda,T_t(b-\lambda\U)+\lambda\U}
~=~
\bfam{\lambda,T_t(b)+\lambda(\U-T_t(\U))},
~~~~~~~~~
}\emun

\vspace{-2ex}\noindent
and it is not really fun to check directly that these form a semigroup on $\C\oplus\cB$. Also, if $T$ is a homomorphism, then clearly so is $\wt{T}$. If $T$ is a contractive(!) CP-map, then $\wt{T}$ is a CP-map, too. The easiest way to see this (see Bhat and Skeide \cite[Section 8]{BhSk00} or \cite{Ske08}), is using the GNS-construction for $T$ to indicate explicitly the GNS-construction for $\wt{T}$:

\bemp \label{GNSuni}
Let $T\colon\cA\rightarrow\cB$ be a contractive CP-map and let $(E,\xi):=$ GNS-$T$. (See around Equation \eqref{GNSeq} in the introduction.) Now $\norm{T(\U)}\le1$, so $\wt{\U}-T(\U)\ge\U-T(\U)\ge0$. We put $\wh{\xi}=\sqrt{\rule[1.3ex]{0pt}{1ex}\,\smash{\wt{\U}-T(\U)}}\in\wt{\cB}$ and define the Hilbert \nbd{\wt{\cB}}module $\wh{E}:=\ol{\wh{\xi}\wt{\cB}}$. We turn $\wh{E}$ into a correspondence by defining a left action of $\wt{\cA}$ as the unique unital extension of the \nbd{0}representation of $\cA$. Likewise, we view $E$ as a correspondence from from $\wt{\cA}$ to $\wt{\cB}$ extending in the only possible way the actions of $\cA$ and $\cB$. We put $\wt{E}:=\wh{E}\oplus E$ and $\wt{\xi}:=\wh{\xi}\oplus\xi$. Then
\vspace{-1ex}
\baln{
\wt{T}
&
~=~
\AB{\wt{\xi},\bullet\wt{\xi}},
&
\cls\wt{\cA}\wt{\xi}\wt{\cB}
&
~=~
\wt{E}.
}\ealn
{~}\vspace{-4.5ex}

\noindent
This both establishes $\wt{T}$ as a CP-map and identifies $(\wt{E},\wt{\xi})$ as GNS-$\wt{T}$.
\eemp
\vspace{-2ex}

% \bthm\label{wstrunithm}
% Let $(\cA,\theta,p)$ be a unital \nbd{C^*}algebra $\cA$, an \nbd{E}semigroup $\theta$ on $\cA$, and a projection $p\in\cA$. Then the following are equivalent:
% \vspace{-1ex}
% \begin{enumerate}
% \item
% $(\cA,\theta,p)$ is a strong dilation.
% \vspace{-1ex}

% \item
% $(\wt{\cA},\wt{\theta},\wt{p})$ with $\wt{p}:=(1,p)=p+\wt{\U}-\U$ is an  \nbd{E_0}dilation.
% \end{enumerate}
% \vspace{-1ex}
% Moreover, in either case, $(\wt{\cA},\wt{\theta},\wt{p})$ is a dilation of the unitalization of the semigroup dilated by $(\cA,\theta,p)$.
% \ethm

% \proof
% Suppose $(\cA,\theta,p)$ is a strong dilation, and denote by $T_t(pap):=p\theta_t(a)p$ the dilated semigroup. Then
% \beqn{
% \wt{p}\wt{\theta}_t(a+\lambda\wt{\U})\wt{p}
% ~=~
% \wt{p}(\theta_t(a)+\lambda\wt{\U})\wt{p}
% ~=~
% p\theta_t(a)p+\lambda\wt{p}
% ~=~
% T_t(pap)+\lambda\U_{\wt{\cB}}
% ~=~
% \wt{T}_t(\wt{p}(a+\lambda\wt{\U})\wt{p}).
% }\eeqn
% So, $(\wt{\cA},\wt{\theta},\wt{p})$ is an \nbd{E_0}dilation of $\wt{T}$.

% Conversely, if $(\wt{\cA},\wt{\theta},\wt{p})$ is an \nbd{E_0}dilation, then $\wt{p}$ is increasing. Consequently,
% \beqn{
% \theta_t(p-\U)+\wt{\U}
% ~=~
% \wt{\theta}_t(\wt{p})
% ~\ge~
% \wt{p}
% ~=~
% p+\wt{\U}-\U,
% }\eeqn
% that is, $\theta_t(\U-p)\le\U-p$, so, $(\cA,\theta,p)$ is a strong dilation. Again, by the first part, the semigroup dilated by $(\wt{\cA},\wt{\theta},\wt{p})$ is unitalization of the semigroup dilated by $(\cA,\theta,p)$.\qed

\bthm\label{wstrunithm}
Let $(\cA,\theta,p)$ be a unital \nbd{C^*}algebra $\cA$, an \nbd{E}semigroup $\theta$ on $\cA$, and a projection $p\in\cA$. \index{unitalization!of a strong dilation}\index{dilation, weak!strong!unitalization of}Then the following are equivalent:
\vspace{-1ex}
\begin{enumerate}
\item
$(\cA,\theta,p)$ is a strong dilation.
\vspace{-1ex}

\item
$(\wt{\cA},\wt{\theta},\wt{p})$ with $\wt{p}:=(1,p)=p+\wt{\U}-\U$ is a strong dilation.
\end{enumerate}
\vspace{-1ex}
Moreover, in either case, $(\wt{\cA},\wt{\theta},\wt{p})$ is a dilation of the unitalization of the semigroup dilated by $(\cA,\theta,p)$.
\ethm

\proof
Recall from Observation \ref{semgenob} that just the strongness condition $p\theta_t(a)p=p\theta_t(pap)p=:T_t(pap)$ is enough to know that $(\cA,\theta,p)$ is a(n of course, strong) dilation, and likewise for $(\wt{\cA},\wt{\theta},\wt{p})$.

We have 
\beqn{
\wt{p}\wt{\theta}_t(a+\lambda\wt{\U})\wt{p}
~=~
\wt{p}(\theta_t(a)+\lambda\wt{\U})\wt{p}
~=~
p\theta_t(a)p+\lambda\wt{p}
}\eeqn
and
\beqn{
\wt{T}_t(\wt{p}(a+\lambda\wt{\U})\wt{p})
~=~
T_t(pap)+\lambda\U_{\wt{\cB}}
~=~
T_t(pap)+\lambda\wt{p}.
}\eeqn
So, $(\cA,\theta,p)$ is a strong dilation (of $T$) if and only if $(\wt{\cA},\wt{\theta},\wt{p})$ is a strong dilation (of $\wt{T}$).\qed

\bob \label{preurecob}
By construction, $\U\wt{\theta}_t(\U\bullet\U)\U=\theta_t$, so that $(\wt{\cA},\wt{\theta},\U)$ is a dilation of $\theta$. Since $\wt{\theta}_t(\wt{\U}-\U)\U=\U-\theta_t(\U)$, this dilations is strong if and only if $\theta$ is an \nbd{E_0}semigroup. By strongness of $(\cA,\theta,p)$, this would make $T$ a Markov semigroup.
\eob

\bthm\label{uninonunithm}
A CP-semigroup $T$ admits a strong dilation if and only if \,$\wt{T}$ admits a weak (and, therefore, strong) dilation.
\ethm

\proof
The \it{only if}-part is settled by the preceding theorem. For the \it{if}-part let $(\wh{\cA},\wh{\theta},\wh{p})$ be a weak (and, therefore, strong) dilation of $\wt{T}$. Denote by $\wh{\U}$ the unit of $\wh{\cA}$. Recall that $\wh{p}=\U_{\wt{\cB}}$ is the unit of $\wt{\cB}$, and define the projections $p:=\U_\cB\in\cB\subset\wt{\cB}=\wh{p}\wh{\cA}\wh{p}\subset\wh{\cA}$ and $q:=\wh{p}-p\in\wt{\cB}\subset\wh{\cA}$. Put $\cA:=(\wh{\U}-q)\wh{\cA}(\wh{\U}-q)$, so that $\U:=\wh{\U}-q$ is the unit of $\cA$.

Since $p=\wh{p}-q\le(\wh{\U}-\wh{p})+(\wh{p}-q)=\wh{\U}-q=\U$, we have $\cB=p\wh{\cA}p=p\cA p\subset\cA$. We shall show that $\wh{\theta}_t(\cA)\subset\cA$ and that $p\wh{\theta}_t(a)p=T_t(pap)$ $(a\in\cA)$. In other words, if we define $\theta$ as the (co)restriction of $\wh{\theta}$ to $\cA$, then $(\cA,\theta,p)$ is a strong dilation of $T$.

We, first, show that $\wh{\theta}_t(\U)\le\U$, so that
\beqn{
\wh{\theta}_t(\cA)
~=~
\wh{\theta}_t(\U)\wh{\theta}_t(\cA)\wh{\theta}_t(\U)
~=~
\U\wh{\theta}_t(\U)\wh{\theta}_t(\cA)\wh{\theta}_t(\U)\U
~\subset~
\U\wh{\cA}\U
~=~
\cA.
}\eeqn
%%%% BO 
% I thought to add a couple of details
Indeed, $q\wh{\theta}_t(q)q=q\wh{p}\wh{\theta}_t(q)\wh{p}q=q\wt{T}_t(q)q=q(\wh{p}-T_t(p))q=q$, thus, $\wh{\theta}_t(q)\ge q$. 
%%%% EO 
So,
\beqn{
\U+q
~=~
\wh{\U}
~\ge~
\wh{\theta}_t(\wh{\U})
~=~
\wh{\theta}_t(\U)+\wh{\theta}_t(q)
~\ge~
\wh{\theta}_t(\U)+q,
}\eeqn
that is, $\wh{\theta}_t(\U)\le\U$.

Now for $a=\U a\U\in\cA$ we have $\theta_t(a)=\wh{\theta}_t(a)\in\cA$, so that $p\theta_t(a)p=\wh{p}\theta_t(a)\wh{p}$. Since $\U\wh{p}=p$, we find
\beqn{
p\theta_t(a)p
~=~
\wh{p}\wh{\theta}_t(\U a\U)\wh{p}
~=~
\wt{T}_t(\wh{p}\U a\U\wh{p})
~=~
\wt{T}_t(pap)
~=~
T_t(pap).\qedsymbol
}\eeqn
\noqed

%%%% BO
% \vspace{-3ex}
% \brem
% In \cite[Theorem 2.2]{ShaSk11} a somewhat different unitalization was considered --- CP-semigroups on $\sB(H)$ were unitized to maps on $\sB(H \oplus \C)$ --- and the ``if" part of the above theorem was obtained. 
% \erem
%%%% EO

\vspace{-1ex}
If $(\wh{\cA},\wh{\theta},\wh{p})=(\wt{\cA},\wt{\theta},\wt{p})$ for some strong dilation $(\cA,\theta,p)$ of $T$ as in Theorem \ref{wstrunithm} (so that $\wh{\U}=\wt{\U}$ is the unit of $\wt{\cA}$ and $\U$ is the unit of $\cA$), then the procedure in the proof gives back $(\cA,\theta,p)$. (Indeed, $q=\wt{p}-p=(1,0)=\wt{\U}-\U$. Therefore, the algebra $\U\wt{\cA}\U$ constructed in the proof is $\cA$, and $\wt{\theta}$ (co)restricted to $\cA$ is $\theta$.)

But $\wh{\cA}$ need not be the unitalization of some unital \nbd{C^*}algebra $\cA$, and $\wh{\theta}$ need not be the unitalization of an \nbd{E}semigroup $\theta$. In fact, the algebra $q\wh{\cA}q$ need not be one-dimensional, and the compression of $\wh{\theta}$ to that corner can be a quite general Markov semigroup. Also, there is no reason, why $q$ should be central in $\wh{\cA}$. However:

\bcor \label{qcentcor}
If $(\wh{\cA},\wh{\theta},\wh{p})$ is full and if $q$ is central, then $(\cA,\theta,p)$ is full, too.
\ecor

\proof
We have $\ls\wh{\cA}\wh{p}\wh{\cA}=\ls\wh{\cA}p\wh{\cA}+\wh{\cA}q\wh{\cA}$, hence, $\ls\cA p\cA=\ls\U\wh{\cA}\wh{p}\wh{\cA}\U$. And $qM(\cC)=M(q\cC)$ for any central $q\in M(\cC)$ and any pre-\nbd{C^*}algebra $\cC$.\qed

\lf
Here are some more simple consequences of the two theorems.\vspace{-1ex}

\bcor \label{unicor1}
If $(\wh{\cA},\wh{\theta},\wh{p})$ is a weak dilation of the Markov semigroup $\wt{T}$, then there exists an \nbd{E_0}dilation of $\wt{T}$ of the form $(\wt{\cA},\wt{\theta},\wt{p})$ for some strong dilation $(\cA,\theta,p)$ of $T$.
\ecor

\bcor \label{unicor2}
If $(\wh{\cA},\wh{\theta},\wh{p})$ is a weak \nbd{E_0}dilation of the Markov semigroup $\wt{T}$, then the dilation $(\wt{\cA},\wt{\theta},\wt{p})$ of $\wt{T}$ constructed from the strong dilation $(\cA,\theta,p)$ of $T$ in the proof of Theorem \ref{uninonunithm}, ``sits inside'' $(\wh{\cA},\wh{\theta},\wh{p})$. More precisely, if we put $\wt{p}:=\wh{p}$ and $\wt{\cA}:=\cA+\C q=\cA+\C\wh{\U}$ (clearly, isomorphic to $\C\oplus\cA$), then  $\wh{\theta}$ (co)restricts to $\wt{\theta}$ on $\wt{\cA}$ and does the job.
\ecor

If $(\wh{\cA},\wh{\theta},\wh{p})$ is not an \nbd{E_0}dilation, there is no reason why $\wh{\theta}_t(\wh{\U})$ should be an element of $\wt{\cA}$.

\bob
In the following section, we will see that coisometric dilations of contraction semigroups in $\cB$ give rise to dilations of the corresponding \it{elementary} CP-semigroups by elementary \nbd{E}semigroups; they may be strong or not. The more important it is, to see that \it{being elementary} is something that gets lost under unitalization. More precisely:

The unitalization $\wt{T}$ of an \hl{elementary CP-map}\index{CP-map!elementary}\index{elementary!CP-map} $T$ (that is, $T=c^*\bullet c$ for some contraction $c\in\cB$) is elementary if and only $T$ is unital. (Indeed, the elementary map $\rtMatrix{\lambda\\c}^*\bullet\rtMatrix{\lambda\\c}=\rtMatrix{\bar{\lambda}\bullet\lambda\\c^*\bullet c}$ on $\wt{\cB}$ is unital if and only if $\abs{\lambda}=1$ and $c^*c=\U$. But then, it is nothing but the unitalization of the (unital) elementary CP-map $c^*\bullet c$.)
\eob

Apart from \index{Bhat's example}Bhat's amazing case study \cite{Bha03} (see Section \ref{EXBexSEC}) and his related work, these notes may be the only place where in the studying weak dilations is not limited to the strong ones, only. In fact, (in particular, in Sections \ref{DilSPSpSEC}, \ref{QSEC}, \ref{EXsubnsupSEC}, \ref{topSEC}, \ref{EXBexSEC}, and \ref{minSEC}) we invest quite a bit of effort in illustrating what works for weak dilations that are not strong, and what does not work. In particular, we point out that our criteria for non-existence of a strong dilation, do not resolve the question whether there is a CP-semigroup (necessarily non-Markov) that has no weak dilation. In fact, this is one of the open problems we leave. We only show that the dilations whose existence we cannot exclude by our criteria, are \it{really bad}: They cannot be \it{good} in the sense of the definition following Example \ref{pPSelemex}.

% \newpage

\index{CP-map!elementary}\index{elementary!CP-map}

\newpage

\section[\sc{Examples:} Weak dilations that are not strong for elementary CP-semigroups]{Examples: Weak dilations that are not strong and classical dilation theory}\label{EXwnsSEC}

This section merges coisometric dilation of contraction semigroups in $\cB=\sB(G)$ with elementary dilations of CP-semigroups and serves to make the reader who knows only classical dilation theory feel comfortable in dilations of  CP-semigroups -- and \it{vice versa}. It can be skipped, until it is referenced to.

A (contraction) semigroup $c=\bfam{c_t}_{t\in\bS}$ over a monoid $\bS$ of elements $c_t$ in a (necessarily unital) \nbd{C^*}algebra $\cB$ gives rise to a CP-semigroup $T=\bfam{T_t}_{t\in\bS^{op}}$ over the opposite monoid of $\bS$, $\bS^{op}$, via $T_t:=c_t^*\bullet c_t$.%
\footnote{ \label{FNelem}
The reason for our unusual choice, writing $c_t^*\bullet c_t$ instead of the more common $c_t\bullet c_t^*$ leading to the fact that the two semigroups are indexed by opposite monoids\phantomsection\index{opposite monoid!semigroups and their (super)(sub)product systems are over opposite monoids}, will be addressed immediately in Section \ref{PSmonoSEC}. For instance, a different choice would cause unpleasant anti-linear relation in formulae like Equation \eqref{unitcomp}, when comparing semigroups $c$ with so-called \it{units} of \it{product systems}.
\vspace{1ex}
}
We call a CP-semigroup $T$ that arises in that way from a contraction semigroup an \phantomsection\hl{elementary} CP-semigroup\index{semigroup!elementary CP-}\index{elementary!CP-semigroup}. An elementary CP-semigroup $T$ is Markov if and only if the semigroup $c$ is isometric, and $T$ is an \nbd{E}semigroup if and only if $c$ is coisometric. Consequently, elementary \nbd{E_0}semigroups are automorphism semigroups.

Suppose $\theta$ is an elementary \nbd{E}semigroup on $\cA$ implemented by a coisometric semigroup $w$ in $\cA$ as $\theta_t=w^*_t\bullet w_t$. Choose a projection $p\in\cA$. Then $(\cA,\theta,p)$ is dilation if and only if the elements $c_t:=pw_tp$ satisfy $c_t^*c_s^*\bullet c_sc_t=c_{st}^*\bullet c_{st}$. It does not follow that the $c_t$ form a semigroup in their own right%
\footnote{
In the case of normal automorphism semigroups $\alpha$ on $\sB(G)$, it is easy to see that each $\alpha_t$ is implemented by a unitary $u_t$ as $\alpha_t=u^*_t\bullet u_t$. But the problem to find these unitaries such that they form a semigroup $u$
% and the question to what extent this semigroup is unique, are 
is nontrivial. Apart from the (obviously affirmative) answer in the discrete one-parameter case, already for the continuous time one-parameter case the answer depends on technical conditions: If $\alpha$ is strongly continuous, then Wigner's theorem \cite{Wig39} states that the answer is affirmative; however, there are non-measurable examples that violate the statement. (This is, essentially, the question, whether or not a one-dimensional algebraic (=without measurability requirements) Arveson system has to be isomorphic to the trivial one or not. See Liebscher \cite[Section 7.2]{Lie09}, in particular, \cite[Example 7.17]{Lie09}.)
\vspace{1ex}
}%
; see Example \ref{nonsolex}. But if the $c_t$ form a semigroup, then also the dilated CP-semigroup $T_t:=c_t^*\bullet c_t$ is elementary. We call such a dilation \hl{solidly elementary}\index{dilation!solidly elementary}\index{elementary!dilation, solidly}.%
\footnote{ \label{elemFN}
Since Skeide \cite{Ske11a}, \it{elementary dilation} is occupied by a different notion; see Appendix \ref{dilcomm}.
}

We may ask, when a solidly elementary dilation is strong. Actually, we can say even a bit more.

\bprop \label{strsolprop}
Let $\bfam{w_t}$ be a semigroup of coisometries in $\cA$ and let $p\in\cA$ be a projection. Then $(\cA,\bfam{w_t^*\bullet w_t},p)$ is a strong dilation if and only if $pw_tp=w_tp$ for all $t$. Moreover, such a strong dilation is solidly elementary, too.
\eprop

\proof
$\theta_t(\U-p)p=0$ ~$\Longleftrightarrow$~ $p\theta_t(\U-p)p=0$ ~$\Longleftrightarrow$~ $\abs{(\U-p)w_tp}^2=0$ ~$\Longleftrightarrow$~ $pw_tp=w_tp$.\qed

\lf
In Example \ref{E0weakex} we have seen a solidly elementary dilation -- even an \nbd{E_0}dilation -- of a CP-semigroup -- even an \nbd{E}semigroup --  on $\sB(G)$ that, by Observation \ref{E0strongob}, is not strong. More will follow soon.
%%%% BO
% I erased the incomplete sentence on Examples \ref{dwnsex} and \ref{cwnsex} 
% I dit not guess what you meant, and they will appear shortly, anyway. 
%%%% EO

\lf
Putting emphasis on the semigroup(s) $w$ (and $c$) rather than the elementary semigroup(s) $\theta$ (and $T$), we say $(\cA,w,p)$ is a \phantomsection\hl{coisometric dilation}\index{dilation!coisometric, of a contraction semigroup} (of the contraction semigroup $c$) if the $pw_tp$ form a semigroup (respectively, if $pw_t p=c_t$).
%%%% BO 
% if $c_t = pw_t p$ for all $t$. 
%%%% EO
The coisometric dilation is \hl{strong}\index{dilation!coisometric, of a contraction semigroup!strong} if $pw_tp=w_tp$ for all $t$. (For instance, the unitary semigroup $\bfam{u_t}_{t\in\R_+}$ in Example \ref{E0weakex} is a strong coisometric dilation of the isometric semigroup $\bfam{v_t}_{t\in\R_+}$. The unitary semigroup $\bfam{u^*_t}_{t\in\R_+}$, instead, is a coisometric dilation of the coisometric semigroup $\bfam{v^*_t}_{t\in\R_+}$ which is not strong.)
%%%% BO
% (When $\cA = \sB(H)$, this is the case if and only if $w_t$ is an extension of $c_t$ for all $t$.) 
%%%% EO 
Every coisometric dilation gives rise to a solidly elementary dilation $(\cA,\theta,p)$; every solidly elementary dilation arises, by definition, from a (though, not unique) coisometric dilation. 

One can say in general that any (not necessarily coisometric) semigroup $\bfam{w_t}_{t\in\bS}$ in $\cA$  is a  \hl{dilation}  of a semigroup $\bfam{c_t}_{t\in\bS}$ in $\cB = p\cA p$ if $pw_t p = c_t$ for all $t \in \bS$.  (Care: While taking adjoints transforms a dilation into a dilation, this need not be so with strong dilations.)

%%%% BO
\brem \label{FSzNrem}
The case where $\bS = \N_0$, $\cB = \sB(G)$ and $\cA = \sB(K)$ has been in investigated in depth since the 1950s, and is the subject of \it{classical dilation theory}; see the monograph Sz.-Nagy and Foias \cite{SzNF2010}. (The discrete multi-parameter case and the continuous one-parameter case have also been addressed; see Sections I.6, I.7, I.9 and Sections I.8, III.8 in \cite{SzNF2010}, respectively.) 
In the theory of Sz.-Nagy and Foias, emphasis in the analysis of a contraction is put on its (minimal) isometric and unitary dilations. (Of course, since the adjoint of a contraction is a contraction, isometric dilations and coisometric dilations just translate into each other under adjoint. And it is comparably easy to promote a (co-)isometric dilation to a unitary dilation. See again Example \ref{E0weakex}, and interpret appropriately the transition $\bfam{u_t}_{t\in\R_+}\leftrightarrow\bfam{u^*_t}_{t\in\R_+}$. And recall also Footnote \ref{FNelem}. In order to not create an unnecessary source of confusion, we reserve \it{strong dilation} to the coisometric ones that interest us.)
By Sz.-Nagy's isometric (unitary) dilation theorem (see Theorems I.4.1 and I.4.2 in \cite{SzNF2010}), for every contraction $c\in\sB(G)$, there exists a Hilbert space $H\supset G$ and an isometry (a unitary) $w\in\sB(H)$ such that $(\sB(H),\bfam{w^n},p)$ ($p$ the projection onto $G$) is an isometric (a unitary) dilation of $\bfam{c^n}$. 
The dilation can also be chosen minimal in an appropriate sense, and the minimal isometric (unitary) dilation is determined uniquely up to unitary equivalence. (The adjoint of the unique minimal isometric dilation is strong in our sense.) We shall not require the deep ramifications of Sz.-Nagy and Foias's theory, but we shall use it as a source of inspiration, intuition (oftentimes misleading) and examples. 

We can recover existence of coisometric dilations of a contraction semigroup in $\sB(G)$, even the minimal one, (and, \it{a fortiori}, existence of a solidly elementary dilation of any elementary CP-semigroup on $\sB(G)$) as the special case $d=1$ from existence of coisometric dilation of (the adjoint of) a (\it{row contractive}) \nbd{d}tuple in Subsection \ref{EXBexSEC}\ref{KrowSSEC}. But since it is so easy, we repeat here the direct classical proof. Given a contraction $c \in \sB(G)$, we put where $\delta := \sqrt{\id_G - cc^*}$, we form the infinite direct sum $H := G \oplus G \oplus \ldots\,$, and we put
\vspace{-1ex}
\beq{ \label{classmin}
w 
~:=~ 
\SMatrix{
c\,\,&\delta & & & 
\\ 
& 0\,\, &\id_G & & 
\\
 & & 0\,\, &\id_G & 
\\ 
& & & \ddots & \ddots}.
\vspace{-1ex}
}\eeq
One checks readily that $(\sB(H),\bfam{w^n},p)$ is a coisometric dilation of $\bfam{c^n}$. (Note, however, that this dilation is, in general, not the minimal one; but, it is strong and can, therefore, easily be minimalized, what we do not discuss.)
\erem

Let us briefly recall that a (not necessarily coisometric) dilation $\bfam{w^n}_{n\in\N_0}$ ($w\in\sB(H)$) of a semigroup $\bfam{c^n}_{n\in\N_0}$ ($c\in\sB(G)$) has (up to unitary equivalence) the general form given by
%%%% EO
\beq{ \label{coisodil}
w
~=~
\SMatrix{\alpha&\beta&\gamma\\&c&\delta\\&&\ve}
~\in~
\sB\rtMatrix{F\\G\\K},
\vspace{-2ex}
}\eeq
where $H=\rtMatrix{F\\G\\K}\supset G$ and $p=\rtMatrix{0&&\\&\sid_G&\\&&0}$.%
\footnote{
The result is due to Sarason \cite{Sar65}. For being self-contained, sufficiency of the form in \eqref{coisodil} being obvious, we sketch a proof of necessity. 
%%%% BO 
Define the subspace $L:=\sscls w^{\N_0}G\supset G$ of $H$, the smallest one containing $G$ and being invariant under $w$. Put $K=L^\perp$. Then $w$ decomposes as $w=\rtMatrix{w_{11}&w_{12}\\&w_{22}}\in\sB\rtMatrix{L\\K}$ and $w^n=\rtMatrix{w_{11}^n&W_n\\&w_{22}^n}$ for some $W_n\in\sB(K,L)$. Since $G\subset L$, we have $p\in\sB(L)\subset\sB(H)$, so $pw_{11}^np=c^n$. Now, put $F:=G^\perp\cap L$ (the complement of $G$ in $L$). The elements $(\ssid_H-p)w^ng$ are total in $F$. Since $L\ni pw_{11}(\ssid_H-p)w^np=pw(\ssid_H-p)w^np=pw^{n+1}p-pwpw^np=0$, both $w_{11}$ and $w$ map $F$ into $F$.
%%%% EO
Therefore, $w_{11}=\rtMatrix{\alpha&\beta\\&c}\in\sB\rtMatrix{F\\G}$.
}
%%%% BO
It may occur that $F$ or $K$ (or both, when $w=c$) are $\zero$. 
%%%% EO

After having secured (by Section \ref{EXBexSEC}\ref{KrowSSEC} or by the construction in Remark \ref{FSzNrem}) existence of a coisometric dilation, we are able to give an example of a dilation to an elementary \nbd{E}semigroup that is not solidly elementary.

\bex \label{nonsolex}
Define the rotation matrix $M:=\rtMatrix{x&-\sqrt{1-x^2}\\\sqrt{1-x^2}&x}\in M_2$ with $x:=\frac{1}{\sqrt{3}}$. For the projection $Q=e_1e_1^*$ onto the first coordinate, we find $QM^2Q=(2x^2-1)Q=-\frac{1}{3}Q$ and $QMQMQ=x^2Q=\frac{1}{3}Q$. So, with $C:=QMQ$, we have ${C^*}^2\bullet C^2=(QM^2Q)^*\bullet(QM^2Q)$, but $C^2=-QM^2Q\ne QM^2Q$. Let $N:=\rtMatrix{0&1&0\\0&0&1\\0&0&0}\in M_3$, and denote by $(\sB(H),\bfam{W^n},P)$ some coisometric dilation of $\bfam{N^n}$. Then putting
\baln{
c
&
~:=~
C\otimes N,
&
w
&
~:=~
M\otimes W,
&
p
&
~:=~
Q\otimes P,
}\ealn
we find that $(M_2\otimes\sB(H),\bfam{{w^*}^n\bullet w^n},p)$ is a dilation of $\bfam{{c^*}^n\bullet c^n}$ satisfying $pw^2p\ne pwpwp$. By Proposition \ref{strsolprop}, this dilation is also not strong.
\eex

We now examine when a coisometric dilation is strong and give some more concrete examples where this is not the case.

\bex \label{dwnsex}
For that a dilation $\bfam{w^n}_{n\in\N_0}$ , given by $w$ as in \eqref{coisodil}, is strong, it is necessary and sufficient that
\vspace{-1ex}
\beqn{
\theta_1(\U-p)p
~=~
\SMatrix{\alpha^*&&\\\beta^*&c^*&\\\gamma^*&\delta^*&\ve^*}
\SMatrix{\ssid_F&&\\&0&\\&&\ssid_K}
\SMatrix{\alpha&\beta&\gamma\\&c&\delta\\&&\ve}
\SMatrix{0&&\\&\ssid_G&\\&&0}
~=~
\SMatrix{\alpha^*&&\\\beta^*&0&\\\gamma^*&0&\ve^*}
\SMatrix{0&\beta&0\\&c&0\\&&0}
~=~
\SMatrix{0&\alpha^*\beta&0\\0&\beta^*\beta&0\\0&\gamma^*\beta&0}
}\eeqn
is $0$. This happens if and only if $\beta$ is $0$. 

First, let us look at the strong case, so $\beta=0$. In this case, conjugation with the canonical isomorphism between $\rtMatrix{F\\G\\K}$ and $\rtMatrix{G\\F\\K}$ transforms $w$ into $\rtMatrix{c&0&\delta\\&\alpha&\gamma\\&&\ve}=:\rtMatrix{c&\delta'\\&\ve'}\in\sB\rtMatrix{G\\K'}$ with $K':=\rtMatrix{F\\K}$. (Also the matrix $w$ in Remark \ref{FSzNrem} has this block-form.) The general form of a strong dilation is, therefore, (up to unitary equivalence) given by $\rtMatrix{c&\delta\\&\ve}$.

Now let us find some concrete $c$ and $w$ with $\beta\ne0$. Recall that $w$ has to be a coisometry, so we have the necessary and sufficient condition
\beqn{
ww^*
~=~
\SMatrix{\alpha&\beta&\gamma\\&c&\delta\\&&\ve}
\SMatrix{\alpha^*&&\\\beta^*&c^*&\\\gamma^*&\delta^*&\ve^*}
~=~
\SMatrix{
\alpha\alpha^*+\beta\beta^*+\gamma\gamma^*~&~~\beta c^*+\gamma\delta^*~~&~\gamma\ve^*
\\
\hfill c\beta^*+\delta\gamma^*~&~~cc^*+\delta\delta^*~~&~\delta\ve^*
\\
\hfill\ve\gamma^*~&~~\hfill\ve\delta^*~~&~\ve\ve^*
}
~=~
\SMatrix{\ssid_F&&\\&\ssid_G&\\&&\ssid_K}.
}\eeqn
$\ve$ must be a coisometry. If $\ve=0$ (meaning $K=\zero$), then $\delta=0$, so $c$ is a coisometry. Dilating a coisometry $c$ to a coisometry $w$, is not really what one typically wants. But, formally, there is no problem in doing also that. (Any $u_t^*$ for fixed $t>0$ in Example \ref{E0weakex}, is an example. And in general, since also $\gamma=0$, we are left with the form $w=\rtMatrix{\alpha&\beta\\&c}$ and the conditions $\beta c^*=0$ and $\alpha\alpha^*+\beta\beta^*=\id_F$. They can be satisfied with any contraction $\beta^*\in\sB(F,G)$ mapping into the orthogonal complement of $c^*G$, meaning $c$ is a proper coisometry if $\beta$ should be nonzero, and $\alpha^*=\sqrt{\id_F-\beta\beta^*}$.) Generally, if $\ve$ is a unitary $u$, then $\gamma$ and $\delta$ have to be $0$. So, we get just a direct sum of the preceding case $\ve=0$ and the unitary $u$ on the third summand $K$.

For getting something less trivial, $\ve$ has to be a proper coisometry, making $K$ necessarily infinite-dimensional. If $F=\zero$, we are back in the strong form $w=\rtMatrix{c&\delta\\&\ve}$, we discussed already.
%%%% BO
% I fixed something but this sentence is still not clear. 
% $\delta^*\in\sB(G,K) $ has to map into the complement of $\ve^*K$ and has to satisfy $\delta\delta^*=\id_G-cc^*$, so, depending on $c$, the dimension of $(cG)^\perp$ may not be too small.)
%%%% EO

Let us try one-dimensional $F=\C$. Then $\beta^*$ and $\gamma^*$ are just vectors, \,$\beta^*\in G$ and $\gamma^*\in(\ve^* K)^\perp\subset K$ (to satisfy $\gamma\ve^*=0$), acting as $\lambda\mapsto\beta^*\lambda$ and $\lambda\mapsto\gamma^*\lambda$, respectively. Additionally, let us also assume that $c$ is the simplest nontrivial (that is, nonzero and nonunitary) contraction possible, namely, $G=\C$ and $c\in\C=\sB(G)$ with $0<\abs{c}<1$. Then also $\delta^*$ is simply a vector in $(\ve^* K)^\perp$ (to satisfy $\delta^*\ve=0$) with length $\norm{\delta^*}=\sqrt{1-\abs{c}^2}$. We have $c\ne0$; we fixed a proper coisometry $\ve$; we chose appropriate $\delta^*$ (with conditions depending only on $c$ and $\ve$); we wish $\beta\ne0$. Of course, $\abs{\alpha}^2+\abs{\beta}^2+\norm{\gamma^*}^2=1$, so all of the three summands are bounded by $1$. Among all $\gamma^*\in(\ve K)^\perp$ satisfying $\beta c^*+\gamma\delta^*=0$, there is a unique one, $\gamma^*_0$, of minimal length $\snorm{\gamma^*_0}=\frac{\abs{\beta}\,\abs{c}}{\norm{\delta^*}}=\frac{\abs{\beta}\,\abs{c}}{\sqrt{1-\abs{c}^2}}$. (All other possible $\gamma^*$ differ from $\gamma^*_0$ by an element in $(\ve^*K)^\perp$ perpendicular also to $\gamma^*_0$. Therefore, if $\ve^*$ is the ``smallest'' proper isometry possible, the \it{one-sided shift}, then $\gamma^*$ is unique.) We must have
\beqn{
\textstyle
1
~\ge~
\abs{\beta}^2+\norm{\gamma^*}^2
~=~
\abs{\beta}^2\bfam{1+\frac{\abs{c}^2}{1-\abs{c}^2}}
~=~
\frac{\abs{\beta}^2}{1-\abs{c}^2},
}\eeqn
that is, $\abs{\beta}\le\norm{\delta^*}$. In particular, we may choose $\beta\ne0$. Choosing also $\alpha$ accordingly, all conditions to make $w$ a coisometry are satisfied. (Note that also $\alpha=0$ is possible.)

To have something fundamentally different, let us now assume that $G=\C^2$ and let $c$ be the projection onto the first basis vector $e_1\in\C^2$ (and still $F=\C$). Then still $\delta$ is characterized by a single vector, namely, $\delta^*_2=\delta^*(e_2)\in(\ve^*K)^\perp$ ($\delta^*(e_1)$ forced to being $0$). Note that $\delta^*_2$ has to be a unit vector, now. The most distinguishing difference, we obtain if we choose for $\beta^*$ the vector $e_2$ so that, now, $\beta c^*=0$. This also forced $\gamma\delta^*$ to being $0$. So, $\gamma^*$ has to be a vector in $(\ve^*K)^\perp$ perpendicular to the unit vector $\delta^*_2$. If $\gamma=0$, then we may stay with the \it{one-sided shift} for $\ve^*$. If we want $\gamma\ne0$, then we need at least two copies of the shift. Anyway, any choice of such $\gamma$ is possible, as long as $\norm{\gamma^*}^2\le1-\norm{\beta^*}^2$ (and $\alpha$ appropriately).
\eex

\bex \label{cwnsex}
In Example \ref{cnonadex} we explain the machinery from Bhat and Skeide \cite{BhSk15} that allows to interpolate a discrete semigroup of operators on  a Hilbert space after having it amplified to the tensor product with $L^2\SB{0,1}$, obtaining a (strongly continuous) semigroup over $\R_+$. Doing this to both the contraction semigroup $c^n$ and its coisometric dilation $w^n$, we obtain a weak continuous time one-parameter dilation that is not strong. Needless to say that every homomorphism $\vp\colon\bS\rightarrow\N_0$ or $\R_+$ turns a one-parameter example into an example over $\bS$. Thus, taking for $\vp$ \nbd{d}fold addition provides \nbd{d}parameter examples from one-parameter examples both discrete and continuous time; see also Remark \ref{homomrem}.
\eex

% \lf
Elementary CP-semigroups and their dilations (elementary and not) will follow us throughout these notes as a source of examples: For an \nbd{E_0}semigroup dilating (necessarily weakly!) a non-Markov \nbd{E}semigroup (Example \ref{E0weakex}); for \it{subproduct systems} (Observation \ref{elemob} and Example \ref{nonadex}); for \it{product systems}, when the dilation is elementary (Example \ref{pPSelemex}); for dilations that are not \it{good} (see after Example \ref{nonadex}) and either not \it{algebraically minimal} (Subsection \ref{minSEC}\ref{algminSSEC}) but with product system, or algebraically minimal but with proper \it{superproduct system} (Section \ref{DilSPSpSEC}), when the dilation is not elementary in Section \ref{EXBexSEC}. The latter section deals with Bhat's example\index{Bhat's example} \cite{Bha03}: A dilation of a scalar (hence, elementary) CP-semigroup that is is not elementary, and exhibits all sorts of bad behaviour that a discrete one-parameter semigroup can possibly exhibit.

In Subsection \ref{EXBexSEC}\ref{KrowSSEC}, we explain the relation between dilations of (normal) discrete one-parameter CP-semigroups on $\sB(G)$ and the \it{coisometric dilations} of (the adjoint of) \it{row contractions} -- a ramification of classical dilation theory. \it{En passant}, we recover with our methods (this time without any countability hypothesis) the result that every row contraction (of arbitrary dimension) admits such a dilation.

Another ramification is the passage from the one-parameter case to the \nbd{d}parameter case. For instance, in Example \ref{minnuniex}, we have an elementary discrete two-parameter CP-semigroup and two (solidly elementary) dilations that are both minimal in the best sense possible, but not conjugate. Elementary discrete three-parameter CP-semigroups give rise to examples that do not admit any strong nor any solidly elementary dilation. We conclude this section by discussing elementary discrete $d$-parameter semigroups and Parrot's classical example \cite{Par70}.

%%%% BO
Given $d$ commuting contractions $c_i$, we obtain a contraction semigroup $c=\bfam{c_\bn}_{\bn\in\N_0^d}$ given by $c_\bn:=c_1^{n_1}\ldots c_d^{n_d}$. Conversely, every contraction semigroup over $\N_0^d$ arises this way. By And\^{o}'s dilation theorem \cite{And63}, every pair of commuting contractions has a coisometric extension, that is, every discrete two-parameter semigroup of contractions has a strong coisometric dilation. \it{A fortiori}, every elementary CP-semigroup  over $\N_0^2$ admits a solidly elementary strong dilation. (While by the note following Theorem \ref{N02Mdilthm} we \bf{do} recover that every such CP-semigroup does have some (strong module) dilation, this time we did not yet find out whether our methods allowed to recover And\^{o}'s result.)

And\^{o}'s dilation theorem raised the question whether every $d$-tuple (also $d\ge3$) of commuting contractions has a commuting (co-)isometric or unitary dilation. The answer is negative, as the following example due to Parrot shows. We adapt the short argument Halmos adopted \cite[p.909]{Hal70}, which shows somewhat more than we need. 

\bex \label{Parex}
Let $F$ be a Hilbert space, put $G = F \oplus F$, and define
\vspace{-.5ex}
\beqn{
c_1 
~:=~ 
\SMatrix{0&v_1\\0&0} 
\,\, , \,\, 
c_2
~:=~ 
\SMatrix{0&v_2\\0&0} 
\,\, , \,\, 
c_3 
~:=~ 
\SMatrix{0&v_3\\0&0} 
\vspace{-.5ex}
}\eeqn
(so that $c_ic_j=0=c_jc_i$) where $v_1 = \id_F$ and $v_2,v_3$ are two coisometries on $F$. Suppose there exists a Hilbert space $H=\rtMatrix{F\\F\\K}$ and commuting coisometries $w_1, w_2, w_3 \in \sB(H)$ such that $pw_ip = c_i$, $i=1,2,3$, where $p$ is the orthogonal projection of $H$ onto $G=F\oplus F$. Let us write
\vspace{-2ex}
\beqn{
w_i
~=~ 
\SMatrix{0&v_i&d_i\\0&0&e_i\\ \ast& \ast & \ast}.
}\eeqn
Then both $v_iv_i^*+d_id_i^*=\id_F+d_id_i^*$ and $e_ie_i^*$ have to be $\id_F$, so that $d_i = 0$ and $e_i$ is a coisometry. Thus, the \nbd{13}element of the product $w_i w_j$ is equal to $v_i e_j$ for all $i,j$. If $w_i w_j = w_j w_i$ for all $i,j$, recalling also that $v_1=\id_F$, we obtain 
\vspace{-1ex}
\baln{
e_2
&
~=~ 
v_2 e_1,
&
e_3
&
~=~ 
v_3 e_1,
&
v_2 e_3 
&
~=~ 
v_3 e_2.
}\ealn

\vspace{-1ex}\noindent
Therefore, $v_2 v_3 e_1 = v_2 e_3 = v_3 e_2 = v_3 v_2 e_1$. Since $e_1$ is a coisometry, multiplying with $e_1^*$ from the right we obtain $v_2 v_3 = v_3 v_2$. Therefore, if we choose $v_2$ and $v_3$ noncommuting (so that $\dim F \geq 2$), then no such $w_i$ can exist.

As a consequence, the maps $T_i := c_i^*\bullet c_i$ generate an elementary CP-semigroup $T$ over $\N_0^3$ that has no solidly elementary dilation.  However, we do not know whether $T$ has a dilation that is not solidly elementary. (Bhat's example\index{Bhat's example} (see Subsection \ref{EXBexSEC}\ref{BexSSEC}) shows that an elementary CP-semigroup might have an \it{incompressible} (weak) dilation that is not elementary.)

% \OW[RETHINK:]{If we replace $c_i$ by $1\oplus c_i$ acting on $\C\oplus G$, then again we have three commuting contractions that have no commuting coisometric dilation; thus the resulting elementary CP-semigroup $T_\bn := c_\bn^* \bullet c_\bn$ has no solidly elementary dilation, as above. In \cite[Theorem 5.14]{ShaSo09} it is proved that in fact $T$ cannot have any strong dilation acting on some $\sB(H)$, and thus we obtain an example of an elementary CP-semigroup with no full strong dilation. We do not know whether or not it has some dilation, perhaps weak, or acting on a different kind of algebra.}
\eex

\newpage

\section{An intermezzo on product systems over monoids, units, and CP-semigroups}\label{PSmonoSEC}

Product systems are crucial almost everywhere in the subject of these notes. In this section, we do not much more than setting up some notation that will be used over and over again. Additionally, we point out that semigroups and their product systems are indexed by \it{opposite} monoids, and we motivate our choice to concentrate on the monoid that indexes the product systems rather than its opposite that indexes the semigroups. We also fix some notations regarding the pre-order structure of a monoid.

\bdefi
\label{PSdefi}
Let $\bS$ be a monoid. A \hl{product system}\index{product system|bf}\index{systems!product} over $\bS$ is a family $E^\odot=\bfam{E_t}_{t\in\bS}$ of correspondences $E_t$ over a \nbd{C^*}algebra $\cB$ with bilinear unitaries $u_{s,t}\colon E_s\odot E_t\rightarrow E_{st}$ such that the \hl{product}\index{product notation!for product systems} $(x_s,y_t)\mapsto x_sy_t:=u_{s,t}(x_s\odot y_t)$ is associative and fulfills the \hl{marginal conditions}\index{marginal conditions} $E_0=\cB$ (the trivial \nbd{\cB}correspondence) with $u_{0,t}$ and $u_{t,0}$ being left and right action of $\cB=E_0$ on $E_t$.
\edefi

We obtain that the diagrams
\beqn{
\xymatrix{
(E_r\odot E_s)\odot E_t	\ar@{=}[r]	\ar[d]_{u_{r,s}\odot\,\sid_t}
							&E_r\odot E_s\odot E_t		&E_r\odot(E_s\odot E_t)
												\ar@{=}[l]	\ar[d]^{\sid_r\odot\,u_{s,t}}	
\\
E_{rs}\odot E_t	\ar[dr]_{u_{rs,t}}	&					&E_r\odot E_{st}	\ar[dl]^{u_{r,st}}
\\
							&E_{rst}				&
}
}\eeqn
and
\beqn{
\xymatrix{
E_0\odot E_t	\ar@{=}[r]	\ar[dr]_{u_{0,t}}		&E_t	\ar[d]^{\sid_t}	&E_t\odot E_0	\ar@{=}[l]	\ar[dl]^{u_{t,0}}	
\\
									&E_t				&
}
}\eeqn
commute for all $r,s,t\in\bS$. (The horizontal lines with $=\!=$ indicate that we are working in a tensor category, where all possible bracketings are identified by (unique) isomorphisms with an object without brackets, and that tensor products of an object with the \it{neutral element $E_0=\cB$} are identified with that object. Thinking in terms of tensor category, frequently helps to switch between \nbd{C^*}versions and von Neumann-versions, but we do not insist in their formal use. In the tensor category of sets with the set product, it is crystal how the identifications have to be chosen. Therefore, it is good to always think of a tensor product as being constructed in the end always starting from a product of sets and quotienting out some equivalence relation. The identifications we indicate by  $=\!=$ will, then, be always obtained by doing to simple tensors $x_1\odot\ldots\odot x_n$ what the formulae suggest for the element $(x_1,\ldots,x_n)$ in the set product.%
\footnote{
\bf{Not} canonical and, therefore, very dangerous, are the identifications $E_s\odot E_t$``$=$''$E_{st}$ via $x_s\odot y_t$``$=$''$x_sy_t$; in fact, on the same family of correspondence (Hilbert spaces) there are non-isomorphic product system structures. See also Observation \ref{noncob} in Appendix \ref{vNAPP}.)
}

\brem
Families of Hilbert spaces that factor as a tensor product, have been known since the sixties of the last century. However, the first formal definition of a one-parameter product system of Hilbert spaces (\hl{Arveson system}\index{Arveson system}, for short) is due to Arveson \cite{Arv89}. In the context of Arveson's theory, Arveson systems arise from (normal) \nbd{E_0}semigroups on $\sB(H)$. The first product systems of general correspondences occurred in Bhat and Skeide \cite{BhSk00} in a construction that associates with each CP-semigroup on a \nbd{C^*}algebra $\cB$ a product system of correspondences over $\cB$. (See the beginning of the introduction.) Subsequently, the connection between \nbd{E_0}semigroups and Arveson systems has been generalized to a construction of product systems of correspondences over $\cB$ from strict \nbd{E_0}semigroups on $\sB^a(E)$ (or normal in the von Neumann-case) where $E$ is a Hilbert \nbd{\cB}module in Skeide \cite{Ske02,Ske03c} and \cite{Ske09} (preprint 2004), and from normal \nbd{E}semigroups by Bhat and Lindsay \cite{Bha96,BhLi05}. Since Skeide \cite{Ske16}, the theory of the relation between one-parameter \nbd{E_0}semigroups and one-parameter product systems has reached a more or less final status. Meanwhile, many more papers about product systems came up; see, for instance, Muhly and Solel \cite{MuSo02}, Hirshberg \cite{Hir04}, Alevras \cite{Ale04}. All these papers refer to the one-parameter cases $\bS=\R_+$ or $\bS=\N_0$. 
%%%% BO 
The first one considering product systems over general semigroups was probably Fowler \cite{Fow02}, who constructed generalized Cuntz-Pimsner algebras from such product systems. Product systems over general semigroups -- in particular, $\N_0^d$ -- continued to be studied in several works on \nbd{C^*} and nonselfadjoint operator algebras. In the context of CP- and \nbd{E}semigroups, product systems over $\N_0^2$ were used by Solel \cite{Sol06} to dilate two commuting CP maps (see also \cite{Sol08}), product systems over $\R_+^2$ were studied by Shalit \cite{Sha08,Sha08b} and \cite{Sha11} to dilate continuous time two-parameter CP-semigroups (see also \cite{Sha10}), and then product systems over subsemigroups of $\R_+^d$ appeared in Shalit and Solel \cite{ShaSo09} in order to tackle dilations of multi-parameter CP-semigroups. Product systems over certain quotients of the free  semigroup were used by Vernik \cite{Ver16} in order to construct dilations for CP maps commuting according to a graph. Recently, product systems over cones in $\R_+^d$ and their connection to multi-parameter \nbd{E_0}semigroups were studied by Murugan and Sundar \cite{MurSu17p,MurSu18,MurSu19}. 
%%%% EO 
\erem

\lf
There is no problem to consider product systems where $\cB$ is nonunital. At first sight, one also may replace $\cB$ with any bigger algebra containing $\cB$ as an ideal. However, as soon as the product system has to satisfy continuity conditions at $t=0$, this will determine $\cB=E_0$ uniquely. (\it{Continuous} product systems are spanned by their \it{continuous} sections. Therefore, in a continuous product system, $E_0$ must be the closure of $\bigcup_{t\ne0}\cB_{E_t}$, where the topology in which we close is that in which the product system should be continuous. On the other hand, since the left action of a correspondence is required nondegenerate, we also cannot make $\cB$ smaller arbitrarily.) So with this in mind, we see that limiting ourselves to unital \nbd{C^*}algebras, as in most of our applications, is a substantial restriction. The following concept of a unit for a product system is, for instance, undefined for nonunital $\cB$; see, however, Bhat, Liebscher, and Skeide \cite[Remark 5.3]{BLS10}.

\bdefi
\label{unitdef}
A \hl{unit}\index{unit|bf}\index{unit!for a product system} for a product system $E^\odot$ is a family $\xi^\odot=\bfam{\xi_t}_{t\in\bS}$ of elements $\xi_t\in E_t$ such that $\xi_s\xi_t=\xi_{st}$ and $\xi_0=\U\in\cB=E_0$.

A unit $\xi^\odot$ is \hl{unital}\index{unit!unital} (\hl{contractive}\index{unit!contractive}, etc.), if $\AB{\xi_t,\xi_t}=\U$ ($\le\U$, etc.) for all $t$.
\edefi

Units are intimately related with CP-semigroups\index{unit!units and CP-semigroups}. Indeed, if $\xi^\odot$ is a unit for a product system $E^\odot$, then
\beq{ \label{unitcomp}
\AB{\xi_s\odot\xi_t,b\xi_s\odot\xi_t}
~=~
\AB{\xi_t,\AB{\xi_s,b\xi_s}\xi_t}.
}\eeq
%%%% BO
Now, letting $T_t:=\AB{\xi_t,\bullet\xi_t}$, we have that the right hand side it equal to $T_t(T_s(b))$, while the left hand side is equal to $\AB{\xi_s\xi_t,b \xi_s \xi_t} = \AB{\xi_{st},b\xi_{st}} = T_{st}(b)$. In other words, the CP-maps $T_t$ form a semigroup $T$. 
%%%% EO
However, they form a semigroup over the \hl{opposite monoid}\index{opposite monoid!semigroups and their (super)(sub)product systems are over opposite monoids|bf} of $\bS$, $\bS^{op}$, and not over $\bS$, because $T_{st}=\AB{\xi_{st},\bullet\xi_{st}}=T_t\circ T_s$. The CP-semigroup $T$ is Markov (contractive) if and only if the unit $\xi^\odot$ is unital (contractive).

Effectively, it was a major task in \cite{BhSk00} to construct not only a product system for a one-parameter CP-semigroup on $\cB$, but, for unital $\cB$ also a unit for that product system giving the CP-semigroup back in the described way; \cite[Theorem 4.8]{BhSk00}. Doing something similar for CP-semigroups over more general monoids, is one of the problems we have to face these notes.

\lf
Let us rest here for a moment to make some considerations on the occurrence of the opposite monoid\phantomsection\index{opposite monoid!semigroups and their (super)(sub)product systems are over opposite monoids|bf} $\bS^{op}$. The fact that product system and semigroup are indexed by opposite monoids, also occurs in the relation between product systems and \nbd{E_0}semigroups. (If an \nbd{E_0}semigroup dilates a Markov semigroup, then the two semigroups are indexed by the same monoid.) The fact that we see it only now, is caused by the fact that this is the first paper where the connection between semigroups and product systems is considered for possibly non-abelian monoids. We also should not hide that the product systems constructed in Arveson \cite{Arv89} from \nbd{E_0}semigroups on $\sB(H)$, and in Muhly and Solel \cite{MuSo02} (and, likewise, in Shalit and Solel \cite{ShaSo09}) from CP-semigroups on von Neumann algebras inspired by a construction from Arveson \cite{Arv97} for CP-maps (paralleling \cite{Arv89}), would be indexed by the \bf{same} monoid as the semigroup, not by the opposite. But these methods work only for von Neumann algebras and von Neumann correspondences, not for \nbd{C^*}algebras and correspondences. Moreover, these product systems do not have units that would give back the original CP-semigroup. The relation between the two product systems has been clarified in Skeide \cite{Ske03c,Ske08}. One product system is the \hl{commutant}\phantomsection\index{opposite monoid!semigroups and their (super)(sub)product systems are over opposite monoids!and commutant} of the other (also called \hl{\nbd{\sigma}dual}\index{opposite monoid!semigroups and their (super)(sub)product systems are over opposite monoids!and \nbd{\sigma}dual} in Muhly and Solel \cite{MuSo04}, where $\sigma$ is a faithful representation of the \nbd{W^*}algebra that must be chosen); units for one product system correspond to \it{covariant representations} of the other. The commutant is anti-multiplicative for the tensor product, and this  explains why if a product system is indexed by $\bS$ its commutant is indexed by $\bS^{op}$. See Appendix \ref{vNAPP}\ref{vNcomm}.

Still, we have to decide if we wish product systems indexed by $\bS$ and semigroups indexed by $\bS^{op}$, or \it{vice versa}. We give two reasons why we opt for the first possibility. The first reason is completely subjective: We prefer to have product systems indexed by $\bS$, because in the end this is a work based on product systems and in this way we reduce to a minimum the attention that has to be paid to notation. The second reason has to do with the way how, usually, semigroups are directed in literature; see the monograph by Clifford and Preston \cite{ClPr61}: When $\bS$ is directed, as required in some of our theorems, then we wish that it be directed following the conventions in literature. We explain briefly, why this is compatible with indexing semigroups by opposites of directed monoids.

If $\bS$ is an abelian semigroup, then one usually says that $s\le t$ if $t$ can be written as $r+s$ for some $r\in\bS$. In order that $\le$ be a direction, a partial order, a total order, and so forth, $\bS$ has to fulfill extra conditions. If $\bS$ is non-abelian, we have two possibilities to compose $t$ out of $s$ and $r$, namely $t=rs$ and $t=sr$. It is the first one,
\beqn{
s \le t
~:\Longleftrightarrow~
t \in \bS s,
}\eeqn
that is related to the following important property used in literature; see \cite[p.\ 34]{ClPr61}.

\bdefi \label{rrefdefi}
A semigroup $\bS$ is \hl{right-reversible}\index{semigroup (algebraic)!right-reversible} if
$\bS s \cap\bS t \neq \emptyset$ for all $s,t \in \bS$.
\edefi

\bprop
\label{r-i-dir-prop}
Let $\bS$ be a semigroup.
\begin{enumerate}
\item
The relation $\le$ is transitive.

\item
If $\bS$ is a monoid, then the relation $\le$ is also reflexive, and $0\le t$ for all $t\in\bS$.

\item
If $\bS$ is a monoid, then the relation $\le$ turns $\bS$ into a directed set if and only if $\bS$ is right-reversible\index{semigroup (algebraic)!right-reversible!direction for monoids}.
\end{enumerate}
\eprop

\proof
Let $r,s,t\in\bS$ such that $r\le s$ and $s\le t$. This means there exists $r'\in\bS$ such that $s=r'r$, and there exists $s'\in\bS$ such that $t=s's$. So, $t=s's=s'r'r$, that is, $r\le t$. In other words, $\le$ is transitive.

If $E$ has a neutral element $0$, then $t=0t\in\bS t$, so $t\le t$ (that is, $\le$ is reflexive), and $t=t0\in\bS 0$, so $0\le t$.

Let $r,s\in\bS$. In order that the monoid $\bS$ be directed by $\le$, we have to find $t\in\bS$ such that $r\le t$ and $s\le t$. In other words, we have to find $r',s'\in\bS$ such that $r'r=s's(=:t)$. But this is possible for all $r,s$ if and only if $\bS$ is right-reversible.\qed
\index{semigroup (algebraic)!right-reversible}
\bdefi \label{Oredefi}
An \hl{Ore semigroup}\index{semigroup (algebraic)!Ore} is a right-reversible semigroup $\bS$ that is also cancellative.
\edefi

Ore semigroups are precisely those which may be embedded into a group $\cG$, the \hl{universal covering group}\index{semigroup (algebraic)!Ore!universal covering group of}, in such a way that $\bS^{-1}\bS=\cG$. (This is particularly interesting, when we wish to embed an \nbd{E_0}semigroup into an automorphism group, as discussed in Laca \cite{Lac00}.) If $\bS$ is an Ore semigroup, then $\le$ is a partial order if and only if $\bS$ has no invertible elements except, possibly, a neutral element $0$.

\bex
The abelian semigroups $\N_0^d$ and $\R_+^d$ are Ore semigroups with no invertible elements but $0$, and their universal covering groups are $\Z^d$ and $\R^d$, respectively.
\eex

\bex
\label{expl:Ore}
For $k\in\N$, let $\bS$ be the universal semigroup generated by two elements $a,b$, such that $a^kb = ba$. Then $\bS$ is an Ore semigroup with no invertible elements; see \cite[Page 36]{ClPr61}. Since $\bS$ is universal and since for $k\ge2$ there exist noncommuting operators $a$ and $b$ satisfying the relations, $\bS$ is nonabelian for $k\ge2$.
\eex

We see that the part of the theory we develop that is based on Ore monoids, covers \nbd{d}para\-meter semigroups of CP-maps, which commute. But it also covers, for example, the monoid generated by two CP maps $T$ and $S$, satisfying $S\circ  T^k = T \circ S$. Setting $T:=a^*\bullet a$ and $S:=b^*\bullet b$ ($a$ and $b$ from the preceding example), gives an example.

Of course, there is an equivalent setting of \hl{left-reversible}\index{semigroup (algebraic)!right-reversible!or left-reversible} semigroup with a direction given by
\beqn{
s \curlyeqprec t
~:\Longleftrightarrow
~t \in s\bS.
}\eeqn
$\bS$ is right-reversible if and only if $\bS^{op}$ is left-reversible. Literature concentrates on the right-reversible case and the corresponding direction in Proposition \ref{r-i-dir-prop}, and so do we. The following example sorts out the question, whether we should consider semigroups over $\bS$ or over $\bS^{op}$. (Note that we only need the relation $\le$, but not that it be a direction. Consequently, $\bS$ need not be right-reversible.)

\bex \label{incrpex}
Recall (Proposition \ref{Mcharprop}) that $(\cA,\theta,p)$ is the dilation of a Markov semigroup if and only if the projection $p\in\cA$ is increasing for $\theta$, that is, if and only if $p\le\theta_t(p)$ for all $t\in\bS$. The idea is, of course, that this should imply $\theta_s(p)\le\theta_t(p)$, whenever $s\le t$. Let us apply $\theta_s$ to $p\le\theta_r(p)$. We get $\theta_s(p)\le\theta_s(\theta_r(p))=\theta_{sr}(p)$ for all $r,s\in\bS$. In other words, $\theta_s(p)\le\theta_t(p)$, whenever $t=sr$ for some $r\in\bS$. As compared with the relation in Proposition \ref{r-i-dir-prop}, this is the wrong order.

Consequently, if $p$ is increasing for a semigroup $\theta$ over $\bS^{op}$, then $s\le t$
%%%% BO 
% I thought this needs emphasizing: 
in $\bS$
%%%% EO 
implies that $\theta_s(p)\le\theta_t(p)$.
\eex
\index{semigroup (algebraic)!right-reversible!or left-reversible}
\lf
Being directed is a crucial property when we will construct dilations from product systems with unital units in Theorem \ref{Oreindthm}. Alone that is enough to understand, why it is crucial to know if we are able to construct a product system with unit for a given CP-semigroup. The discrete one-parameter case $\N_0$ is easy. The proof in \cite{BhSk00} for both one-parameter cases, the continuous time case $\R_+$ and the discrete case, relies on an extra property of the relation $\le$ for $\R_+$ and $\N_0$, which is particularly simple for $\N_0$. (Roughly, the sets $\bJ_t$, our viewpoint of the interval partition to be introduced below, are directed for $\R_+$ and $\N_0$, and for $\N_0$ they even have unique maximal elements so that the inductive limit in Theorem \ref{indlimthm} may be avoided by just writing down things for the unique maximal elements.) We discuss a sufficient property of $\le$ for general $\bS$. 

\bdefi
\label{totdirdefi}
A semigroup $\bS$ is \hl{totally directed}\index{semigroup (algebraic)!totally directed} if for all $s,t\in\bS$ we have $s\le t$ or $t\le s$.
\edefi

A totally directed monoid $\bS$ is right-reversible. (If $s\le t$, then $\bS s\ni t\in\bS t$, and if $t\le s$, then $\bS t\ni s\in\bS s$.) But even then, $\le$ need not be a total order. (For instance, if $\bS$ is a group, then $s\le t$ \bf{and} $t\le s$ for all $s,t\in\bS$.) If $\bS$ is an Ore semigroup, then the set $\set[r\in\bS]{\exists s,t\in\bS\colon rs=t,t\le s}$ is a subgroup of $\bS$. So, if $\bS$ is a totally directed monoid with no invertible elements but $0$, then it is totally ordered. (Indeed, if $s\le t$ and $t\le s$ with $s\ne t$, then $r$ such that $rs=t$ is a nonneutral element of that subgroup.)

\blem\label{totdirlem}
A right-reversible semigroup $\bS$ is totally directed if and only if for each chain $s_1\le\ldots\le s_n=t$ and $s\le t$, the element $s$ can be ``inserted'' into that chain, that is, there is a chain $s'_1\le\ldots\le s'_{n+1}:=t$ and $1\le k\le n$ such that $s'_i=s_i$ for $i<k$, $s'_k=s$, and $s'_i=s_{i-1}$ for $i>k$.
\elem

\proof
If $\bS$ is totally directed, then for each $i$ we have $s_i\le s$ or $s\le s_i$. Choose for $k$ the smallest number such that $s_k\not\le s$, for then $s\le s_k$, so $s_{k-1}\le s\le s_k$.

Conversely, suppose $\bS$ fulfills the second property and choose $r,s\in\bS$. Since $\bS$ is right-reversible, there are $r',s'\in\bS$ such that $r'r=s's=:t$. So, $r\le t$ and $s\le t$. By hypothesis, $s$ can be inserted into the chain $0\le r\le t$, so that $s\le r$ or $r\le s$.\qed

\lf
The following is a refinement of the notation for tuples \phantomsection\index{tuple notation!$\bJ_t$} of elements in $\R_+$, invented in \cite{BhSk00}. Let $\bS$ be a semigroup. For $t\in\bS$ we define
\beqn{
\bJ_t
~:=~
\bCB{\bt=(t_n,\ldots,t_1)\colon n\in\N,t_n\ldots t_1=t}.
}\eeqn
For $\bs=(s_m,\ldots,s_1)\in\bJ_s$ and $\bt=(t_n,\ldots,t_1)\in\bJ_t$, we define their \hl{join}\index{tuple notation! join of tuples} $\bs\smallsmile\bt$ in $\bJ_{st}$ as
\beqn{
\bs\smallsmile\bt
~:=~
(s_m,\ldots,s_1,t_n,\ldots,t_1).
}\eeqn
For $\bs=(s_m,\ldots,s_1)\in\bJ_t$ we say $\bs\le\bt\in\bJ_t$ if there exist tuples $\bs_i\in\bJ_{s_i}$ such that
\beqn{
\bt
~=~
\bs_m\smallsmile\ldots\smallsmile\bs_1.
}\eeqn
We also denote by $()$ the \hl{empty tuple}\index{tuple notation! empty tuple}. We put $()\smallsmile\bt:=\bt=:\bt\smallsmile()$. Note, however, that $()\notin\bJ_t$!

\bprop \label{Jpoprop}
$\le$ is a partial order on $\bJ_t$, and $(t)$ is the unique minimal element.
\eprop

\proof
Obviously, the relation $\le$ on $\bJ_t$ is reflexive and transitive. It is also anti-symmetric. (Indeed, if $\bt=(t_n,\ldots,t_1)=\bs_m\smallsmile\ldots\smallsmile\bs_1$, then $n\ge m$. If also $\bt\le\bs$, then $n\le m$, too. So, the only way to write $\bt$ in the described way, is with $\bs_i=(t_i)$, whence, $t_i=s_i$.) So, $\le$ is a partial order. Of course, $(t)$ is strictly smaller than any other element in $\bJ_t$.\qed

\lf
Much of the construction of a product system for a CP-semigroup in \cite{BhSk00}, depends on that in the case $\bS=\R_+$ the partial order on $\bJ_t$ is a direction. In fact, in that case $\bJ_t$ is a lattice, a lattice isomorphic to the lattice of interval partitions (with double points) of $\SB{0,t}$ which takes its operations from ``union'' and ``intersection'' of chains $s_1\le\ldots\le s_{n-1}\le s_n=t$. In general, if $s_1\le\ldots\le s_{n-1}\le s_n=t$ is a chain in a semigroup $\bS$, then there are elements $t_i\in\bS$ such that $t_1=s_1$ and $s_i=t_is_{i-1}$ ($1<i\le n)$. In particular, we have $t_i\ldots t_1=s_i$ and $\bt=(t_n,\ldots,t_1)\in\bJ_t$. However, this decomposition is not necessarily unique. It is unique, if $\bS$ is cancellative. Also, given two chains the union of whose elements can be arranged to form a chain, the order in which the elements appear in the united chain is not necessarily unique. It is unique, if the relation $\le$ on $\bS$ is a partial order. Even worse: The union of two chains need not even allow to be arranged into a chain. In that case, $\bJ_t$ has no chance to be directed.

For these reasons, several constructions of product systems, like that in \cite{BhSk00} but also others to be discussed later on in these notes, are possible only under the following, quite restrictive, condition.

\bthm\label{totdirthm}
Let $\bS$ be a totally directed, cancellative semigroup. Then $\bJ_t$ with the partial order $\le$ is directed.
\ethm

\proof
We have to show that for each $\bt=(t_n,\ldots,t_1)$ and $\bt'=(t'_{n'},\ldots,t'_1)$ in $\bJ_t$ there exists $\br\in\bJ_t$ such that $\bt\le\br$ and $\bt'\le\br$. If one of the two tuples $\bt$ or $\bt'$ is minimal, then there is nothing to show. So let us assume that $n>1$ and $n'>1$. Let us start with the chains $s_1\le\ldots\le s_{n-1}\le t$ and $s'_1\le\ldots\le s'_{n'-1}\le t$, where we put $s_i:=t_i\ldots t_1$ and $s'_i:=t'_i\ldots t'_1$. Extending the procedure from the proof of Lemma \ref{totdirlem}, we start constructing a new chain in which each $s_i$ and each $s'_i$ occurs and in which the primed and unprimed members stay in the order they had originally (a \it{shuffle product}) in the following way. Write down the first members $s_1,\ldots,s_{k_1-1}$ from the unprimed chain that fulfill $s_i\le s'_1$. (If $s_1\not\le s'_1$, then write nothing.) Then take the first members $s'_1,\ldots,s'_{k'_1-1})$ of the primed chain that fulfill $s'_i\le s_{k_1}$. (This partial chain contains at least $s'_1$.) If we are not yet done, take the next members $s_{k_1},\ldots,s_{k_2-1}$ from the unprimed chain that satisfy $s_i\le s'_{k'_1}$. If we are not yet done, take the next members $s'_{k'_1},\ldots,s'_{k'_2-1}$ from the primed chain that satisfy $s'_i\le s_{k_2}$. Continue until there are no elements left (so that we have a chain of $n+n'-2$ elements) and add the \nbd{(n+n'-1)}st element $s=t$. We get a chain of the following form
\bmun{
\bfam{\sigma_1\le\ldots\le\sigma_{n+n'-1}=t}
~=~
\\
\bfam{s_1\le\ldots\le s_{k_1-1}\le s'_1\le\ldots\le s'_{k'_1-1}\le\ldots\le s_{k_{\ell-1}}\le\ldots\le s_{k_\ell-1}\le s'_{k'_{\ell-1}}\le\ldots\le s'_{k'_\ell-1}\le s=t}
}\emun
where $k_\ell=n$ and $k'_\ell=n'$. It is possible that $k_1=1$, so that the first block of unprimed $s_i$ is empty, and it is possible that $k'_{\ell-1}=k'_\ell$, so that the last block of primed $s'_i$ is empty. (This chain is not unique. For instance, we could also have started with the primed quantities. If $s_1\le s'_1$, $s'_1\le s_1$, and $s_1\ne s'_1$ the two procedures would have resulted in two different chains, one starting with $s_1$ the other with $s'_1$. But the procedure ``minimizes'' in a sense the number $\approx 2\ell-1$ of switches between primed and unprimed parts.) So far, this works for all totally directed semigroups.

To get $\br$, we do to that chain $\sigma_1\le\ldots\le\sigma_{n+n'-1}=t$ what we described in front of the theorem, getting $r_i$ such that $r_i\ldots r_1=\sigma_i$. Now, since $\bS$ is cancellative, all partial products $r_i\ldots r_j$ $(i\ge j)$ are uniquely fixed by the property that $r_i\ldots r_{j+1}\sigma_j=\sigma_i$. Now let $\vk_i$ denote  the position in the chain where $s_i$ appears. Then $\vk_{i+1}>\vk_i$ and
\beqn{
r_{\vk_{i+1}}\ldots r_{\vk_i+1}s_i
~=~
s_{i+1},
}\eeqn
so that $r_{\vk_{i+1}}\ldots r_{\vk_i+1}=t_{i+1}$. If we put $\vk_1:=0$, then that formula remains true also for $i=0$. So, $\br_i:=(r_{\vk_i},\ldots,r_{\vk_{i-1}})\in\bJ_{t_i}$ such that $\br=\br_n\smallsmile\ldots\smallsmile\br_1$, that is, $\bt\le\br$. The same argument for primed quantities (with $\vk'_i$ giving tuples $\br'_i\in\bJ_{t'_i}$) shows that $\bt'\le\br$.\qed

\lf
We do not address the following natural questions or similar ones:
\begin{itemize}
\item
How far can an abelian totally directed cancellative monoid be away from the positive elements of an ordered field or its singly generated submonoids?

\item
Are there nonabelian examples? This question has been answered by Levy \cite{Lev20p} in the affirmative sense. (In \cite[Section 2.2]{Lev20p} he exhibits an example of a nonabelian monoid which is even  totally ordered. That paper also contains a discussion on various constraints that a total order imposes on the structure of a monoid.)

\item
Are there monoids $\bS$ such that the $\bJ_t$ are directed which are not totally directed and cancellative? (If the answer was no, this would simplify quite a bit the hypothesis of quite a number of theorems, later on.)

\end{itemize}

% \OW[OPEN/DECIDE WHAT TO INCLUDE]{QUESTIONS:

% How far is an abelian directed monoid from a submonoid of $\R_+$?

% (We have the positive numbers in an ordered field, should we mention?)

% Are there non-abelian directed monoids without invertible elements?

% (We have and example, but it is not Ore)

% }

\newpage

\section{Super-  and subproduct systems}\label{SPSpbSEC}

Let us return to 
%%%% BO 
the notion of 
%%%% EO
product systems and some of its ramifications. A product system is, loosely speaking, a family of correspondences with associative identifications $E_s\odot E_t=E_{st}$. What we get from CP-semigroups, in a first step, is a family with identifications $E_s\odot E_t\supset E_{st}$, a \it{subproduct system}; see Section \ref{CPspsSEC}. What we get from a strong dilation without further constraints is a family with identifications $E_s\odot E_t\subset E_{st}$, a \it{superproduct system}; see Section \ref{DilSPSpSEC}. Motivated by the recurrence of subproduct systems in the construction of L\'{e}vy processes (Sch\"{u}rmann \cite{MSchue93}), of dilations (Bhat and Skeide \cite{BhSk00}, Muhly and Solel \cite{MuSo02}), of new product systems from given ones (Skeide \cite{Ske06d}, Bhat and Mukherjee \cite{BhMu10}), and from multivariate operator theory (Shalit and Solel \cite{ShaSo09}, Davidson, Ramsey, and Shalit \cite{DRS11}), Shalit and Solel \cite{ShaSo09} have given a formal definition of \it{subproduct systems} of \nbd{W^*}correspondences, followed independently by \it{inclusion systems} of Hilbert spaces by Bhat and Mukherjee \cite{BhMu10}. Instances of \it{superproduct systems} have occurred in Hellmich, K\"{o}stler, and K\"{u}mmerer \cite[Section 4.5]{HKK04p}, in Skeide \cite[Sections 9 and 10.2]{Ske06d} (made explicit in Skeide \cite{Ske09r2}), in Margetts and Srinivasan \cite{MaSr13}, and in \cite{Ske08p2}. Margetts and Srinivasan \cite{MaSr13} gave the first formal definition of superproduct system of Hilbert spaces. It should be noted that their definition of product system of von Neumann correspondences in \cite{MaSr14p} is not compatible with their definition of superproduct system in \cite{MaSr13,MaSr14p}. (An Arveson system, according to them, would not be a product system of \nbd{\C}correspondences, unless one-dimensional; but it would be a superproduct system.)

There are several pitfalls about terminology (in particular, regarding units and morphisms) around, which we -- hopefully -- avoided bothering the reader with in this section and the main body of the text. Only in the Example-Section \ref{EXexpSEC} we faced the problems head on, culminating in the extended Remark \ref{copurem} about consequences of terminology.

\bdefi \label{SPSUdef}
A \hl{superproduct system}\index{superproduct system|bf}\index{systems!superproduct} $E^\podot=\bfam{E_t}_{t\in\bS}$ over $\bS$ is like a product system over $\bS$ in Definition \ref{PSdefi} (including $E_0=\cB$) with structure maps $v_{s,t}$, just that the structure maps are only required to be isometries. (The marginal ones, $v_{0,t}$ and $v_{t,0}$, remain unitaries, automatically. We continue with the product notation\index{product notation!for superproduct systems} $x_sy_t:=v_{s,t}(x_s\odot y_t)$.)

A \hl{subproduct system}\index{subproduct system|bf}\index{systems!subproduct} over $\bS$ is a family $E^\bodot=\bfam{E_t}_{t\in\bS}$ (including $E_0=\cB$) with isometric structure maps $w_{s,t}\colon E_{st}\rightarrow E_s\odot E_t$ satisfying the opposite of the diagrams following Definition \ref{PSdefi}.

A \phantomsection\hl{unit}\index{unit!for a superproduct system}\index{unit!for a subproduct system} for a super(sub)product system is a family $\xi^\odot=\bfam{\xi_t}_{t\in\bS}$ of elements $\xi_t\in E_t$ fulfilling 
%%%% BO 
$\xi_s\xi_t=\xi_{st}$ 
%%%% EO
$(w_{s,t}\xi_{st}=\xi_s\odot\xi_t)$ and $\xi_0=\U$.

A super(sub)product system with adjointable structure maps is called an \phantomsection\hl{adjointable}\index{adjointable!(super)(sub)product system}\index{systems!adjointable}\index{superproduct system!adjointable}\index{subproduct system!adjointable} super(sub)product system.
\edefi

Note that a product system is an adjointable superproduct system, and that a product system may be turned into an adjointable subproduct system by defining the structure maps $w_{s,t}:=u_{s,t}^*$. A super(sub)product system is called a \phantomsection\hl{proper}\index{superproduct system!proper}\index{subproduct system!proper} super(sub)product system if it is not a product system.

We mentioned already that there are plenty of subproduct systems (and units) arising from CP-semigroups, and examples of superproduct systems arising from their dilations. (Concrete examples of a dilations with proper superproduct systems are in Section \ref{EXpropsupSEC}.) Simple examples of proper subproduct systems can be obtained by cutting a product system to be $\zero$ starting from $t\ge t_0>0$. Note that these examples are adjointable. Simple examples of superproduct systems can be obtained by replacing the members $E_t$ of a product system for $t\le t_0>0$ with $F_t$ from a product subsystem (for instance, $E_t=\cB$ in the trivial product system and $F_t=\cI$, an ideal). While it is not so easy to obtain non-adjointable subproduct systems (see Section \ref{EXnonadSEC}), we easily obtain non-adjointable superproduct systems (for instance, by choosing the ideal $\cI$ to be non-complemented in $\cB$). See Remark \ref{adrem} below and the forthcoming results for more about adjointability. See also Section \ref{EXexpSEC} for how to obtain continuous time multi-parameter examples by \it{exponentiating} discrete multi-parameter examples.

Among sub- and superproduct systems, although much more recent, the notion of superproduct system is formally much closer to the notion of product system. In fact, if a subproduct system is not adjointable, then there is no \it{product} there, but only a \it{coproduct}. That makes some of the following definitions for subproduct systems look different from the corresponding definitions for product systems, while the versions for superproduct systems are perfect analogues of those for product systems. Therefore, we continue giving, like in the preceding definition, the definitions for superproduct systems first.

\bdefi \label{subdefi}
Let $E^\podot=\bfam{E_t}_{t\in\bS}$ be a superproduct system over a monoid $\bS$. A family $\bfam{F_t}_{t\in\bS}$ of subcorrespondences $F_t\subset E_t$ with $F_0=E_0$ is a \phantomsection\hl{(super)(sub)product subsystem}\index{systems!supersubproduct@(super)(sub)product subsystems}\index{superproduct system!supersubproduct@(super)(sub)product subsystems of}\index{subproduct system!supersubproduct@(super)(sub)product subsystems of} of the superproduct system $E^\podot$ if
%%%% BO
\beqn{
v_{s,t}(F_s\odot F_t)
~~~(\subset)~
~(\supset)
~=~~
F_{st}
}\eeqn
%%%% EO
for all $s,t\in\bS$. (This definition applies, in particular, to product systems $E^\odot$ which are particular superproduct systems.)

Let $E^\bodot=\bfam{E_t}_{t\in\bS}$ be a subproduct system over a monoid $\bS$. A family $\bfam{F_t}_{t\in\bS}$ of subcorrespondences $F_t\subset E_t$ with $F_0=E_0$ is a \phantomsection\hl{(super)(sub)product subsystem}\index{subproduct system!supersubproduct@(super)(sub)product subsystems of} of the subproduct system $E^\bodot$ if
\beqn{
F_s\odot F_t
~~~(\subset)~
~(\supset)
~=~~
w_{s,t}F_{st}
}\eeqn
for all $s,t\in\bS$. (This definition applies, in particular, to product systems $E^\odot$ which are particular subproduct systems with structure maps $w_{s,t}:=u_{s,t}^*$, and in that case it coincides with the preceding definition.)
\edefi

In the sequel, in a superproduct system,  we frequently shall write $v_{s,t}(F_s\odot F_t)$ as ~$\cls F_sF_t$.

It is clear that a super(sub)prod\-uct subsystem of a super(sub)product system is a super(sub) product system with the (co)re\-stricted structure maps. It requires a moments thought, to see that a subproduct subsystem of a superproduct system is, indeed, a subproduct system: 

\blem \label{isubcoilem}
Suppose we have an isometry $v\colon\sE\rightarrow\sF$ and submodules $E\subset\sE$ and $F\subset\sF$ such that $vE\supset F$. Then there exists a (unique) isometry $w\colon F\rightarrow E$ such that $vwy=y$ for all $y\in F$.
\elem

\proof
$w$ is simply the (co)restriction of the inverse of the map $v$ considered onto $v\sE$ to the claimed domain and codomain. We present this in more detail.

$v$ corestricts to a unitary $u\colon\sE\rightarrow v\sE$ and $v\sE\supset vE\supset F$.  Of course, $u^*\colon v\sE\rightarrow\sE$ satisfies $vu^*y=y$ for all $y\in v\sE$. The restriction $w$ of $u^*$ to $F$ maps $F$ (isometrically) into $E$. (Indeed, suppose $u^*y\notin E$. Then $y=vu^*y\notin v E\supset F$.) Again $vwy=vu^*y=y$. Since $v$ is isometric, $w$ is determined uniquely.\qed

\lf
Of course, if all modules are bimodules and $v$ is (apart from being automatically right linear) also left linear, then, being the (co)restriction of the inverse of a left linear map, also $w$ is left linear. Applying all this to ~~$v:=v_{s,t}$, ~~~$\sE:=E_s\odot E_t\supset F_s\odot F_t=:E$, ~and ~~$\sF:=E_{st}\supset F_{st}=:F$, we get left linear isometries $w_{s,t}\colon F_{st}\rightarrow F_s\odot F_t$. Coassociativity follows easily from associativity of the $v_{s,t}$. The proof that superproduct subsystems of subproduct systems form superproduct systems, is analogue.

It is noteworthy that the definition of subproduct subsystem is compatible with \cite[Definition 5.1]{ShaSo09}.

\bex \label{unitgenex}
If $\sU$ is a set of units for a superproduct system $\sE^\podot$, then the correspondences $S\!_t:=\cls\bigcup_{\xi^\odot\in\sU}\cB\xi_t\cB$ form a subproduct subsystem $S^\bodot$ of $\sE^\podot$.\index{superproduct system!subproduct subsystem spanned by a set of units} (Indeed,
\beqn{
\cls S\!_sS\!_t
~=~
\cls\bigcup_{\xi^\odot,\xi'^\odot\in\sU}\cB\xi_s\cB\xi'_t\cB
~\supset~
\cls\bigcup_{\xi^\odot\in\sU}\cB\xi_s\xi_t\cB
~=~
\cls\bigcup_{\xi^\odot\in\sU}\cB\xi_{st}\cB
~=~
S\!_{st},
}\eeqn
that is, $\cls S\!_sS\!_t\supset S\!_{st}$.) Every unit in $\sU$ is a unit for $S^\bodot$, too. (If $w_{s,t}$ denotes the coproduct of $S^\bodot$ inherited from $\sE^\podot$, then $w_{s,t}\xi_{st}$ is that unique element $\Xi_{s,t}\in S\!_s\odot S\!_t\subset\sE_s\odot\sE_t$ such that $v_{s,t}\Xi_{s,t}=\xi_{st}$. Since $v_{s,t}(\xi_s\odot\xi_t)=\xi_s\xi_t=\xi_{st}$, by uniqueness of $\Xi_{s,t}$ we get $w_{s,t}\xi_{st}=\xi_s\odot\xi_t$.) Obviously, $w_{s,t}$ is the unique subproduct system structure on the family $S\!_t$ making every $\xi^\odot\in\sU$ a unit.
\eex

\lf
We illustrate that the structures recovered  from the containing structure, work together nicely.

\bprop \label{bicontprop}
Choose $\#_1,\#_2,\#_3$ from $\CB{\podot,\bodot}$ and suppose that $E^{\#_1}\subset F^{\#_2}\subset G^{\#_3}$, in the sense of Definition \ref{subdefi}, so that also $E^{\#_1}\subset G^{\#_3}$.

Then the structure of $E^{\#_1}$ being recovered from $E^{\#_1}\subset F^{\#_2}$ with the structure of $F^{\#_2}$ being recovered from $F^{\#_2}\subset G^{\#_3}$, coincides with the structure of $E^{\#_1}$ being recovered from $E^{\#_1}\subset G^{\#_3}$.
\eprop

\bemp[Convention.] \label{subconv}\index{conventions!numbered!subsystems@\ref{subconv} (subsystems)}
After this proposition, we include in the notation $E^{\#_1}\subset F^{\#_2}$ to mean that also the structure of  $E^{\#_1}$ is that inherited from $F^{\#_2}$.
\eemp

\proof[Proof of Proposition \ref{bicontprop}.~]
The cases of three equal $\#_i$ being trivial and the cases where for at least one $i$ we have $\#_i=\#_{i+1}$ being easy once the remaining cases are understood, we treat only the most tricky cases $\#_1\ne\#_2\ne\#_3$ (so that $\#_1=\#_3$). Actually, we treat only one of these two cases, because the other one is perfectly analogue.

So let $E^\podot\subset F^\bodot\subset G^\podot$ and denote by $v_{s,t}$ the product of the superproduct system $E^\podot$ inherited from the containing subproduct system $F^\bodot$, denote by $w_{s,t}$ the coproduct of the subproduct system $F^\bodot$ inherited from the containing superproduct system $G^\podot$, and denote by $V_{s,t}$ the product of the superproduct system $G^\podot$. Choose $x_s\in E_s$ and $y_t\in E_t$. Then, $v_{s,t}$ being induced by $w_{s,t}$,
\beqn{
z_{s,t}
~:=~
v_{s,t}(x_s\odot y_t)
}\eeqn
is the unique element in $E_{st}\subset F_{st}$ such that $w_{s,t}z_{s,t}=x_s\odot y_t$. Further, $w_{s,t}$ being induced by $V_{s,t}$,
\beqn{
Z_{s,t}
~:=~
w_{s,t}z_{s,t}
}\eeqn
is the unique element in $F_s\odot F_t\subset G_s\odot G_t$ such that $V_{s,t}Z_{s,t}=z_{s,t}$. Putting everything together we get that, by construction,
\beqn{
v_{s,t}(x_s\odot y_t)
~=~
z_{s,t}
~=~
V_{s,t}Z_{s,t}
~=~
V_{s,t}(w_{s,t}z_{s,t})
~=~
V_{s,t}(x_s\odot y_t),
}\eeqn
so that the product of $E^\podot$ coincides with that inherited by (co)restriction from the product of $G^\podot$.\qed

\bprop \label{PSsubprop}
If $E^\#$, respectively, $F^{\#'}$ is a product system, then validity of $E^\#\subset F^{\#'}$ does not depend on the choice whether we consider $E^\#$, respectively, $F^{\#'}$ as superproduct system or as a subproduct system.
\eprop

\proof
\item
$E^\odot\subset F^\podot$, that is, $\cls E_sE_t=v_{s,t}(E_s\odot E_t)=E_{st}$. The product $u_{s,t}$ of $E^\odot$ inherited from $F^\podot$, it the (co)restriction of the product $v_{s,t}$. Its adjoint coproduct $u_{s,t}^*$ sends $x_{st}:=y_sz_t$ to $y_s\odot z_t$, the unique element $X_{s,t}$ in $E_s\odot E_t$ such that $v_{s,t}X_{s,t}=v_{s,t}(y_s\odot z_z)=y_sz_t=x_{s,t}$. So, $u_{s,t}^*$ coincides with the coproduct inherited from $F^\podot$.

\item
$E^\odot\subset F^\bodot$, that is, $E_s\odot E_t=w_{s,t}E_{st}$. The product $u_{s,t}$ of $E^\odot$ inherited from $F^\bodot$ sends $x_s\odot y_t$ to that unique element $Z_{s,t}$ in $E_{st}$ such that $w_{s,t}Z_{s,t}=x_s\odot y_t$. By surjectivity of $w_{s,t}$ and uniqueness of $Z_{s,t}$, the $u_{s,t}$ have to be surjective. Obviously, $u_{s,t}$ is the adjoint of the (co)restriction of the coproduct $w_{s,t}$, that is, the coproduct inherited from $F^\bodot$.

\item
$E^\podot\subset F^\odot$. If we consider $F^\odot$ as a superproduct system, this means $\cls E_sE_t=u_{s,t}(E_s\odot E_t)\subset E_{st}$ with the (co)restriction $v_{s,t}$ of the product $u_{s,t}$ as product. If we consider $F^\odot$ as a subproduct system, this means $E_s\odot E_t\subset u_{s,t}^*E_{st}$. This is, obviously, equivalent to the containment $u_{s,t}(E_s\odot E_t)\subset E_{st}$. Moreover, the product inherited from the coproduct $u_{s,t}^*$ send $x_s\odot y_t$ to the unique $Z_{s,t}\in E_{st}$ such that $u_{s,t}^*Z_{s,t}=x_s\odot y_t$, that is, $Z_{s,t}=u_{s,t}(x_s\odot y_t)=x_sy_t$ coincides with the product inherited from the superproduct system $F^\odot$.

\item
$E^\bodot\subset F^\odot$.  If we consider $F^\odot$ as a superproduct system, this means $\cls E_sE_t=u_{s,t}(E_s\odot E_t)\supset E_{st}$. The coproduct $w_{s,t}$ inherited from the product $u_{s,t}$ sends $x_{st}$ to that unique $X_{s,t}\in E_s\odot E_t$ such that $u_{s,t}X_{s,t}=x_{s,t}$, that is, such that $X_{s,t}=u_{s,t}^*x_{s,t}$. It follows that this coproduct coincides with the coproduct inherited from the subproduct system $F^\odot$ by (co)restriction.\qed

\lf\lf
No discussion of a structure 
%%%% BO 
without the discussion of the notion of \it{morphism}; no discussion of substructure without the discussion of the notion of \it{embedding}. 
%%%% EO
It is quite clear, what morphisms and embeddings are among superproduct systems, and what morphisms and embeddings are among subproduct systems. But in order to deal, later on, with our concrete embedding problems (basically, embedding subproduct systems into superproduct systems and, further, into product systems), we need these notions in the same generality as in Proposition \ref{bicontprop}.

A morphism between product systems $E^\odot$ and $E'^\odot$ over $\bS$ is a family of maps $a_t\in\sB^{bil}(E_t,E'_t)$ that respects the products, that is, $a_{st}(x_sy_t)=(a_sx_s)(a_ty_t)$, respectively, $a_{st}u_{s,t}=u'_{s,t}(a_s\odot a_t)$. We require by hand that $a_0=\id_\cB$. (Note that this is not automatic. It is forced that $a_0$ is an algebra homomorphism, but more we cannot deduce; for instance, without any further conditions, $a_0$ can very well be $0$. (There are some considerations relevant for the one-parameter case in Skeide \cite[Remark 3.8]{Ske06c}, Abbaspour and Skeide \cite[Section 4]{AbSk07}, and Dor-On and Markiewicz \cite{D-OMa14}.) Note, too, that in many other papers, it would be required that the $a_t$ be adjointable; we do not do that, as our particular interest lies in embeddings that are non-adjointable.) As usual, this notion generalizes immediately to morphism between superproduct systems: $a_{st}v_{s,t}=v'_{s,t}(a_s\odot a_t)$.
\begin{subequations} \label{diam}
\beq{ \label{m1}
\parbox{10cm}{\xymatrix{
E_s\odot E_t		\ar[d]_{a_s\odot\,a_t}		\ar[rr]^{v_{s,t}}			&&		E_{st}		\ar[d]^{a_{st}}
\\
E'_s\odot E'_t														\ar[rr]_{v'_{s,t}}			&&		E'_{st}
}}
}\eeq
While for subproduct systems, we have to take the condition that the coproduct is preserved: $w'_{s,t}a_{st}=(a_s\odot a_t)w_{s,t}$.
\beq{ \label{m2}
\parbox{10cm}{\xymatrix{
E_s\odot E_t		\ar[d]_{a_s\odot\,a_t}												&&		E_{st}		\ar[d]^{a_{st}}	\ar[ll]_{w_{s,t}}			
\\
E'_s\odot E'_t																								&&		E'_{st}									\ar[ll]^{w'_{s,t}}							
}}
}\eeq
Of course, in either case, isomorphisms are morphisms that consist of unitaries, and embeddings are morphisms that consist of isometries.

We see that in the diagram with the same vertices and the same vertical arrows, the horizontal arrows are simply the structure maps with arrow pointing into the corresponding (and, in the non-adjointable case, only possible) direction. We, therefore, require from a morphism $E^\podot\rightarrow E'^\bodot$ that $w'_{s,t}a_{st}v_{s,t}=a_s\odot a_t$:
\beq{ \label{m3}
\parbox{10cm}{\xymatrix{
E_s\odot E_t		\ar[d]_{a_s\odot\,a_t}	\ar[rr]^{v_{s,t}}				&&		E_{st}		\ar[d]^{a_{st}}
\\
E'_s\odot E'_t																								&&		E'_{st}									\ar[ll]^{w'_{s,t}}							
}}
}\eeq
And for a morphism $E^\bodot\rightarrow E'^\podot$ we require that $a_{st}=v'_{s,t}(a_s\odot a_t)w_{s,t}$:
\beq{ \label{m4}
\parbox{10cm}{\xymatrix{	
E_s\odot E_t		\ar[d]_{a_s\odot\,a_t}												&&		E_{st}		\ar[d]^{a_{st}}		\ar[ll]_{w_{s,t}}
\\
E'_s\odot 														E'_t	\ar[rr]_{v'_{s,t}}		&&		E'_{st}
}}
}\eeq
\end{subequations}

\bdefi \label{moemisdef}
A \hl{morphism}\index{systems!morphisms}\index{morphisms|see {systems}} from $E^{\#}$ to $E'^{\#'}$ is a family $a^\#=\bfam{a_t}_{t\in\bS}$ of maps $a_t\in\sB^{bil}(E_t,E'_t)$ fulfilling $a_0=\id_\cB$ and the one relevant diagram among the Diagrams \eqref{m1}, \eqref{m2}, \eqref{m3}, \eqref{m4}.

An \hl{embedding}\index{systems!embeddings} is an isometric morphism.

An \hl{isomorphism}\index{systems!isomorphisms} is a surjective embedding.
\edefi

\noindent
Why did we indicate a morphism from $E^{\#}$ to $E'^{\#'}$ only as $a^\#$ with the operation $\#$ from the domain system, ignoring what the operation $\#'$ is? The answer is that $a^\#$ is a morphism for the operation $\#$ in the following sense:

\bprop \label{submoprop}
Suppose for $E^\#$ and $E'^{\#'}$ we have a family of maps $a_t\in\sB^{bil}(E_t,E'_t)$ fulfilling $a_0=\id_\cB$. Then $\bfam{a_t}_{t\in\bS}$ is a morphism from $E^{\#}$ to $E'^{\#'}$ (denoted $a^\#$) if and only if:
\begin{enumerate}
\item \label{submo1}
For the family $\bfam{\ol{a_tE_t}}_{t\in\bS}$ of subcorrespondences of $E'_t$ we have $\bfam{\ol{a_tE_t}}_{t\in\bS}=(\ol{aE})^\#\subset E'^{\#'}$ in the sense of Convention \ref{subconv}.

\item \label{submo2}
$a^\#$ is a morphism $E^{\#}$ to $(\ol{aE})^\#$.
\end{enumerate}
\eprop

\lf\noindent
Condition \ref{submo1} means that the family  $\bfam{\ol{a_tE_t}}_{t\in\bS}$ is a \nbd{\#}subsystem of $E'^{\#'}$. Condition \ref{submo2} means that $a^\#$ is a morphism for the operation $\#$ on $\bfam{\ol{a_tE_t}}_{t\in\bS}$ and  inherited from being a \nbd{\#}subsystem of $E'^{\#'}$, that is, $a^\#$ among the Diagrams \ref{diam} satisfies the relevant for the operations $\#$ to $\#$, not $\#$ to $\#'$.

\lf
\proof[Proof of Proposition \ref{submoprop}.~]
The two cases $\#'=\#$ being obvious by standard considerations, among the remaining two cases $\#'\ne\#$ we deal only with $\#=\bodot$ and $\#'=\podot$, the other case being analogue.

If $\bfam{\ol{a_tE_t}}_{t\in\bS}$ is a subproduct subsystem of superproduct system $E'^\podot$, then its inherited coproduct $w'_{s,t}$ sends $x'_{st}\in\ol{a_{st}E_{st}}$ to that unique $X'_{s,t}\in\ol{a_sE_s}\odot\ol{a_tE_t}$ fulfilling $v'_{s,t}X'_{s,t}=x'_{st}$. Applying this to $x'_{st}:=a_{st}x_{st}$, we get
\beqn{
(a_s\odot a_t)w_{s,t}x_{st}
~=~
w'_{s,t}a_{st}x_{st}
~~~\Longleftrightarrow~~~
v'_{s,t}(a_s\odot a_t)w_{s,t}x_{st}
~=~
a_{st}x_{st},
}\eeqn
that is, Diagrams \eqref{m2} and \eqref{m4} are equivalent. So, the only remaining thing to be shown, is that validity of Diagram \eqref{m4} implies that $\bfam{\ol{a_tE_t}}_{t\in\bS}$ is a subproduct subsystem of $E'^\podot$. But this follows from
\beqn{
v'_{s,t}(\ol{a_sE_s}\odot\ol{a_tE_t})
~=~
v'_{s,t}\ol{(a_s\odot a_t)(E_s\odot E_t)}
~\supset~
v'_{s,t}\ol{(a_s\odot a_t)w_{s,t}E_{st}}
~=~
\ol{a_{st}E_{st}}.\qedsymbol
}\eeqn
\noqed

It is noteworthy that when $E^\#$ or $E'^{\#'}$ is a product system, then the set of morphisms does not depend on whether we consider that product system as a superproduct system or as a subproduct system. (This follows, essentially, from Proposition \ref{PSsubprop}.) We, therefore, see that \it{isomorphism} in the case $\#'\ne\#$ is somewhat `boring', because in this case $E'^{\#'}$ and, therefore, also $E^\#$ both are necessarily product systems.

\brem \label{adrem}
An isometry between Hilbert modules need not be adjointable. (It is adjointable if and only if its range is complemented.) In fact, in Example \ref{nonadex} we shall see a subproduct system (coming from a CP-semigroup) with nonadjointable structure maps. The definition in \cite{ShaSo09} is for subproduct systems of \nbd{W^*}correspondences, where all maps have adjoints automatically. The definition in Viselter \cite{Vis10} is for \nbd{C^*}correspondences, and requires adjointable structure maps, explicitly. However, our subproduct systems do come from CP-semigroups, so we prefer not to include adjointability in the definition.

One might wonder what happens if we define subproduct systems not via isometric coproduct maps, but via \it{coisometric} product maps. While it is clear how isometries have to be defined -- and it turns out they need not be adjointable --, it is quite unclear, how to define \it{coisometry} without requiring explicitly that they be adjointable. We propose the following `\it{co-pair}' of definitions: An \hl{isometry} is, then, a contractive right linear map $E\rightarrow F$ between Hilbert modules for which there exists a submodule $F'$ of $F$ such that the map corestricts to a unitary $E\rightarrow F'$; a \phantomsection\hl{coisometry}\index{operator!isometry!isometric \it{versus} coisometric} is a contractive right linear map $E\rightarrow F$ between Hilbert modules for which there exists a submodule $E'$ of $E$ such that the map restricts to a unitary $E'\rightarrow F$. Obviously, the new definition of isometry is equivalent to the definition to be an inner product preserving map, and, therefore, much less than required in the new definition has actually to be checked. (A map that corestricts to a unitary, is right linear and contractive automatically; also the submodule $F'$, obviously is bound to be the image of $E$ under this map.) In Appendix \ref{coisoAPP} we show that a coisometry according to the proposed definition \bf{is} adjointable. (Here, the hypotheses right linear and contractive are indispensable, but we get automatic adjointability; but, like for isometries, the submodule to which we restrict is unique.) Therefore, describing the structure of a subproduct system (a superproduct system) in terms of coisometric products (coproducts), would be a possibility to incorporate automatic adjointability.

We conclude this remark with: A map $v\colon E\to F$ is a \phantomsection\hl{partial isometry}\index{operator!isometry!partial} if it corestricts to a coisometry $E\rightarrow vE$. It is adjointable if and only if $vE$ is complemented in $F$. (A definition as unitary (co)restriction to a map $E'\rightarrow F'$ is not feasible, because we may choose $E'$ and $F'$ jointly too small (for instance, $\zero$). We do not know how to avoid the asymmetry in this definition. But, after all the different automatic properties of isometries and coisometries show already that there is no symmetry to be expected right from the beginning. Note however: If the map $v$ is adjointable, then $v$ is a partial isometry if and only if $v^*$ is a partial isometry.)
\erem

Some results on super- and subproduct systems which we discuss in the following sections, depend on adjointability of the structure maps.

\blem
Suppose we have an adjointable isometry $v\colon\sE\rightarrow\sF$ and submodules $E\subset\sE$ and $F\subset\sF$ such that $vE\subset F$. If $E$ is complemented in $\sE$, then the (co)restriction $w\colon E\rightarrow F$ is adjointable, too.
\elem

\proof
Denote by $p\in\sB^a(\sE)$ the projection onto $E$. Then $pv^*\upharpoonright F$ is an adjoint of $w$. (Indeed, $\AB{pv^*y,x}=\AB{y,vpx}=\AB{y,wx}$ for all $y\in F$ and all $x=px\in E$.)\qed

\lf
Note that the condition that $E$ is complemented in $\sE$ is indispensable. Indeed, suppose $\sE=\sF=F$ with $E$ not complemented in $\sE$ and let $v=\id_\sE$. Then $w$ is not adjointable. This may also be used to construct a counter example when the hypothesis of being complemented is dropped from the following result.

\bcor \label{asubscor}
A super(sub)product subsystem sitting complementedly in an adjointable super(sub)product system, is adjointable, too.
\ecor

% We briefly discuss two constructions of product systems, one from superproduct systems, the other from subproduct systems, under the hypothesis of Theorem \ref{totdirthm} where $\bJ_t$ is directed.

We now switch our attention to the generation of super- and subproduct systems and, closely related, to the intersection of subsystems. It is well-known that for subspaces $V_i$ and $W_i$ of vector spaces $V$ and  $W$, respectively, the tensor products fulfill \phantomsection\index{tensor product!intersection property}$(V_1\cap V_2)\otimes(W_1\cap W_2)=(V_1\otimes W_1)\cap(V_2\otimes W_2)$. But already for Hilbert spaces one has to work. In Appendix \ref{laAPP} we discuss this and illustrate that this \it{intersection property} remains true even for families of von Neumann subcorrespondences. We do not know if it holds for \nbd{C^*}correspondences. However, for subcorrespondences $E_i$ of $E$ and subcorrespondences $F_i$ of $F$, it is obvious that $(\bigcap_i E_i)\odot(\bigcap_i F_i)$ is contained in $\bigcap_i(E_i\odot F_i)$. Therefore:

\bprop
The intersection of a family of superproduct subsystems of a given superproduct system is again a superproduct subsystem.
\eprop

The analogue statement for subproduct subsystems of a given subproduct system is true for von Neumann correspondences, but presently we do not know if it is true for \nbd{C^*}cor\-res\-pond\-ences. In particular, we do not know if the intersection of product subsystems of \nbd{C^*}correspondences is again a product system.

\brem
It should be noted that the intersection of nonzero Arveson systems (which are continuous) can be the zerosystem (which is not continuous at $t=0$); just take two different one-dimensional subsystems. Therefore, the intersection of continuous product systems can be critical. (For Arveson systems, Liebscher \cite[Theorem 5.7]{Lie09} has shown that the intersection is at least measurable.)
\erem

By the usual abstract intersection, we obtain:

\bcor
For every family of subsets $S_t\subset \sE_t$ of a superproduct system $\sE^\podot$, there is a unique smallest superproduct subsystem containing them.
\ecor

We can specify this better:

\bthm \label{supintthm}
Let $\sE^\podot=\bfam{\sE_t}_{t\in\bS}$ be a superproduct system over a monoid $\bS$ with subsets $S_t\subset\sE_t$. Then the subcorrespondences
\beqn{
E_t
~:=~
\cls\bigcup_{\bt\in\bJ_t}\cB S_{t_n}\ldots\cB S_{t_1}\cB
}\eeqn
of $\sE_t$ for $t\ne0$ and $E_0:=\cB$ form a superproduct subsystem $E^\podot$ of $\sE^\podot$, the superproduct subsystem \phantomsection\hl{generated}\index{superproduct system!superproduct subsystem generated by a subfamily} by $\bfam{S_t}_{t\in\bS}$.
\ethm

\proof
Clearly, the $E_t$ are subcorrespondences of $\sE_t$. Now
\bmu{\label{SupPSincl}
\Bfam{\bigcup_{\bs\in\bJ_s}\cB S_{s_m}\ldots\cB S_{s_1}\cB}\Bfam{\bigcup_{\bt\in\bJ_t}\cB S_{t_n}\ldots\cB S_{t_1}\cB}
~=~
\bigcup_{\substack{\br\in\bJ_{st}\\(s,t)\le\br}}\cB S_{r_k}\ldots\cB S_{r_1}\cB
\\
~\subset~
\bigcup_{\br\in\bJ_{st}}\cB S_{r_k}\ldots\cB S_{r_1}\cB.
~~~~~~~~~~~~
}\emu
Taking the closed linear span of both sides, it follows that $E_sE_t\subset\cls E_sE_t\subset E_{st}$.\qed

\lf
The following general result is interesting in its own right. Its full power comes out, however, only in the totally directed case in Theorem \ref{supsubthm}.

\bprop \label{PsubSubgenprop}
Let $\bfam{{F^i}^\odot}_{i\in I}$ be a family of product subsystems of the superproduct system $\sE^\podot$. Then the $S_t:=\cls\bigcup_{i\in I}F^i_t$ form a subproduct subsystem of $\sE^\podot$.\index{superproduct system!subproduct subsystem spanned by a subfamily of product subsystems}
\eprop

\proof
Clearly $S\!_t$ is a subcorrespondence of $\sE_t$. Since ${F^i}^\odot$ is a product system, we have $F^i\!\!_{st}=\cls F^i\!\!_sF^i_t$, so
\beqn{
\bigcup_{i\in I}F^i\!\!_{st}
~=~
\bigcup_{i\in I}\cls F^i\!\!_sF^i_t
\subset
\cls\Bfam{\bigcup_{i\in I}F^i\!\!_s}\Bfam{\bigcup_{i\in I}F^i_t}
~=~
\cls S\!_sS\!_t,
}\eeqn
hence, $S\!_{st}=\cls\bigcup_{i\in I}F^i\!\!_{st}\subset\cls S\!_sS\!_t$. In other words, the $S\!_t$ form a subproduct subsystem.\qed

\lf
Some results follow along the lines of \cite{BhSk00}, when $\bS$ is such that each $\bJ_t$ is directed.

\bthm\label{supsubthm}
In the situation of Theorem \ref{supintthm}, let $\bS$ be totally directed and cancellative (or, more generally, let $\bJ_t$ be directed for all $t$). Then:

If $S^\bodot=\bfam{S_t}_{t\in\bS}$ is a subproduct subsystem of $\sE^\podot$, then $E^\podot$ is a product subsystem.
\ethm
              
\proof
If the $S\!_t$ form a subproduct system, so that $S\!_t\subset\cls S\!_{t_n}\ldots S\!_{t_1}=\cls\cB S\!_{t_n}\ldots\cB S\!_{t_1}\cB$, then $\cB S\!_{s_m}\ldots\cB S\!_{s_1}\cB\subset\cls\cB S\!_{t_n}\ldots\cB S\!_{t_1}\cB$ whenever $\bs\le\bt\in\bJ_t$. If $\bJ_{st}$ is directed (for instance, if $\bS$ is a totally directed, cancellative monoid), then $\CB{\br\in\bJ_{st}\colon(s,t)\le\br}$ is a cofinal subset of $\br\in\bJ_{st}$. So, in the limit, the inclusion in the closed linear span of \eqref{SupPSincl} becomes equality.\qed

\bcor\label{totdircor}
Let $\bS$ be totally directed and cancellative (or, more generally, let $\bJ_t$ be directed for all $t$). Then, by Proposition \ref{PsubSubgenprop} and Theorem \ref{supsubthm}, every superproduct system $\sE^\podot=\bfam{\sE_t}_{t\in\bS}$ over $\bS$ has a unique maximal product subsystem, namely, the superproduct subsystem $E^\odot$ generated by the family of all product subsystems. (In case this family is empty, we put $E_t:=\zero$ for $t\ne0$.)
\ecor

\bcor\label{unitgencor}
Let $\sU$ be a set of units for a superproduct system $\sE^\podot$, and denote by $S\!_t:=\cls\bigcup_{\xi^\odot\in\sU}\cB\xi_t\cB$ the correspondences that form the subproduct subsystem $S^\bodot$ of $\sE^\podot$ generated by $\sU$ as in Example \ref{unitgenex}. If $\bS$ is totally directed and cancellative (or, more generally, if $\bJ_t$ is directed for all $t$), then the superproduct subsystem generated by $S^\bodot$ is a product system, the product system generated by the units in $\sU$.
\ecor

\lf
In Example \ref{unitgenex}, we have seen that units of a superproduct system (\nbd{\cB})span a subproduct subsystem. By Corollary \ref{unitgencor}, they generate a product subsystem containing that subproduct system, if all $\bJ_t$ are directed. Since a superproduct system with a unit is what we get from a strong dilation in Theorem \ref{sdilunithm}, and since, for Ore monoids, unit plus product(\bf{!}) system would allow to construct a module dilation of that CP-semigroup (Theorem \ref{Oreindthm}), it would be great if we could get the statement of Corollary \ref{unitgencor} for Ore monoids. However, in the Example-Section \ref{EXN02SEC} (Theorem \ref{facpsgenthm}) we see a superproduct system over the Ore monoid $\N_0^2$ with a unit that has no product subsystem containing the unit.

On the other hand, once the $\bJ_t$ are directed, the statement that a subproduct system generates a product system, unlike in Theorem \ref{supsubthm}, does actually not depend on having this subproduct system sitting in a superproduct system; we can construct the containing product system directly from the subproduct system. (That is, in the totally directed case, it is not necessary to hypothesize existence of a strong dilation to show that the subproduct system of a CP-semigroup embeds into a product system.) The following theorem is in several ways an adaptation of Bhat and Skeide \cite[Section 4]{BhSk00}.

\bthm\label{indlimthm}
Let $\bS$ be a totally directed and cancellative monoid (or, more generally, let $\bJ_t$ be directed for all $t$) and let $\sE^\bodot=\bfam{\sE_t}_{t\in\bS}$ be a subproduct system over $\bS$. \index{product system!generated by a subproduct system over a totally directed monoid}Then there is a (unique) minimal product system $E^\odot=\bfam{E_t}_{t\in\bS}$ containing $\sE^\bodot$ and all its units, the product system \hl{generated} by its subproduct subsystem $\sE^\bodot$.
\ethm

We could have simply imitated the proof of \cite[Theorem 5]{BhMu10} in the one-parameter Hilbert space case, defining the maps
\beqn{
\beta_{\bs_m\smallsmile\ldots\smallsmile\bs_1,\bs}
\colon
\sE_{s_m}\odot\ldots\odot\sE_{s_1}
~\longrightarrow~
(\sE_{s^m_{m_m}}\odot\ldots\odot\sE_{s^m_1})\odot\ldots\odot(\sE_{s^1_{m_1}}\odot\ldots\odot\sE_{s^1_1})
}\eeqn
by iterating the subproduct embeddings $\sE_{s_i}\rightarrow\sE_{s^i_{m_i}}\odot\ldots\odot\sE_{s^i_1}$, showing they form an inductive system, and observing that the $E_t:=\limind_{\bt\in\bJ_t}\sE_{t_n}\odot\ldots\odot\sE_{t_1}$ form the desired (and unique) minimal product system $E^\odot=\bfam{E_t}_{t\in\bS}$. But, verifying the necessary associativity conditions of the $\beta_{\bs_m\smallsmile\ldots\smallsmile\bs_1,\bs}$ is painful (though not really difficult), and had been dealt with in \cite{BhMu10} only quite stepmotherly. Therefore, we hesitate to reference to that proof. On the other hand, we also hesitate to repeat the same type of arguments over and over again (like \cite{BhMu10} repeats \cite{BhSk00}, \cite{BBLS04}, and \cite{Ske06d}). Despite we trust (like \cite{BhMu10}) that the interested reader will be able to verify the necessary associativity conditions and follow through the proof from \cite{BhMu10}, we prefer to be formally complete. So, here is a proof by reduction:

\proof[Proof of Theorem \ref{indlimthm}.~]
We reduce the statement to the statement from \cite{BhSk00} (suitably generalized to monoids as in our hypotheses) that the GNS-subproduct system of any CP-semigroup (on a unital \nbd{C^*}algebra) embeds into a product system. Then, we appeal to Theorem \ref{adSPS-CPthm} saying that every adjointable subproduct system is the subproduct system of a strict CP-semigroup on a suitable $\sB^a(E)$. Finally, we boil down the general case (non-adjointable subproduct systems) to the adjointable case by passing to von Neumann correspondences. (The last two parts require knowledge form sections still to come.) 

\item
The proof in \cite{BhSk00} is not suffering the defect of unclear associativity, because it starts from the GNS-subproduct system of a CP-semigroup which is spanned by a single unit. Here associativity of the $\beta$ is crystal from the more explicit character of the construction, and adapting the proof from \cite{BhSk00} to more general monoids is practically word by word. We get: Every GNS-subproduct system (over $\bS$ as in the hypotheses) embeds into a product system.

\item
 By Theorem \ref{adSPS-CPthm}, every adjointable subproduct system arises by applying the Morita equivalence operation $E^*\odot\bullet\odot E$ to the GNS-subproduct system of a suitable strict CP-semigroup on a suitable $\sB^a(E)$. Applying the same operation to the product system containing the GNS-subproduct system (which exists by the first part), we get a product system of \nbd{\cB}cor\-re\-spond\-ences containing our original superproduct system. 
%%%% BO 
% Hello Michael. I think this statement is somewhat "unfair" (or maybe speculative would be a better way to say it), since the preprints appeared essentially at the same time, and independently. How could they reduce to a result they did not know? So I commented out the following line: (I suggest we delete it). 
% (By the way, already the proof in \cite{BhMu10} could have been reduced in that way to \cite[Corollary 2.10]{ShaSo09}.)
%%%% EO

\item
By embedding a non-adjointable subproduct system into a subproduct system of von Neumann correspondences (automatically, adjointable), by the first two parts, we get an embedding into a product system of von Neumann correspondences. Now, looking at how the product system in the first part has been obtained as an inductive limit, taking also note of the fact that the inductive system of von Neumann \nbd{\sB^a(E)}correspondences gives rise (by the Morita equivalence operations) to an inductive system of von Neumann \nbd{\cB}correspondences and that, further, the $\beta$ of these von Neumann correspondences restrict precisely to the $\beta$ on the \nbd{C^*}correspondences as stated in front of this proof, we get an increasing family of \nbd{C^*}correspondences whose limit, necessarily, form a product system, the one we seek. In this way, also the uniqueness statement is made evident.\qed

\bob
The proof depends on an inductive limit over the directed set $\bJ_t$. So, there is no hope to provide a construction (possibly, leading to a superproduct system containing the given subproduct system and being generated by it as a superproduct system) if the $\bJ_t$ are not directed. In the case of CP-semigroups, the superproduct system containing the subproduct system derives from the assumption that there exists a strong dilation -- and it need not contain a product system also containing the subproduct system.
\eob

\brem
It is noteworthy that, following the construction indicated before the proof of Theorem \ref{indlimthm}, it is not even necessary to worry to prove that the inductive limit is a product system. It is, clearly, a superproduct system; it contains the subproduct system; and it is, clearly, generated by the latter. So, by Theorem \ref{supsubthm}, it is a product system. The given proof cannot take advantage of this observation, though.
\erem

% \brem\label{SPSCPrem}
%In the end we are mainly interested in subproduct systems coming from CP-semigroups. We will show (Theorem \ref{adSPS-CPthm}) that every adjointable subproduct system of \nbd{\cB}cor\-respondences arises in a certain way from a (strict) CP-semigroup on some $\sB^a(E)$. If the subproduct system of the latter embeds into a product system, then so does the former. Therefore, whenever we are able to show that the subproduct system of a CP-semigroup over a certain monoid embeds into a product system, then this is true for all adjointable subproduct systems over that monoid. This argument includes, in particular, the case of one-parameter subproduct systems of Hilbert spaces dealt with in \cite{BhMu10}, which are adjointable, and reduces the proofs there to the case of one-parameter CP-semigroups dealt with in \cite{BhSk00}. Taking into account also Shalit and Solel \cite[Theorem ???]{ShaSo09}, which asserts that every subproduct system of \nbd{W^*}correspondences is the subproduct system of Arveson-Stinespring correspondences of a CP-semigroup, taking into account that this subproduct system is isomorphic to the \it{commutant system} in the se sense of Skeide \cite{Ske03c,Ske08} of the CP-semigroup, and taking also into account the result from \cite{BhSk00} that the latter embeds into a product system, \cite{BhMu10} could have been proved by reduction to these results.
% \erem

The following discussion, which concludes this section, is towards the comparison of the superproduct systems of dilations of $T$ and of $\wt{T}$ in Section \ref{unisupSEC}. Recall that $\wt{\cB}$ contains $\cB\ni\U$ as direct summand. But we prefer to do some of the discussion, where $\cI$ is just any closed ideal of $\cB$.

Let $\cI$ be a closed ideal in the (for the time being not necessarily unital) \nbd{C^*}algebra $\cB$, and let $\sE^\podot$ be a superproduct system of \nbd{\cB}correspondences $\sE_t$ over the monoid $\bS$. We, then, may define the \nbd{\cI}correspondences $\sF_t:=\cls\cI\sE_t\cI$. Since
\beqn{
\sF_s\sF_t
~=~
\cls\cI\sE_s\cI\sE_t\cI
~\subset~
\cls\cI\sE_s\sE_t\cI
~\subset~
\cls\cI\sE_{st}\cI
~=~
\sF_{st},
}\eeqn
the $\sF_t$ with the (co)restrictions of the product of $\sE^\podot$ form a superproduct system $\sF^\podot$, too.

\blem
Let $E$ be a Hilbert \nbd{\cB}module and let $F$ be a correspondence from $\cB$ to $\cC$. Then $\cls E\odot\cI F=E\odot F$ if and only if $\cls\cB_E\cI F=\cls\cB_E F$.
\elem

\proof
This follows from $\cls E\cB_E=E$ and $\cls(E^*\odot\id_F)(E\odot F)=\cls\cB_EF$.\qed

\bcor \label{PSIPScor}
Suppose that the \nbd{\cI}correspondences $\sF_t$ defined above satisfy $\sF_t=\cls\cI\sE_t$ for all $t\in\bS$. Then $\sF^\podot$ is a product system, if $\sE^\podot$ is a product system.
\ecor

The opposite statement need not be true. (Indeed, suppose $\cB=\cJ\oplus\cI$ and take the \it{external direct sum} of a superproduct system of \nbd{\cJ}correspondences that is not a product system, and of a product system of \nbd{\cI}correspondences.)

\brem \label{PSIPSrem}
In Skeide \cite{Ske08p}, we call a product system $\sE^\odot$ fulfilling (as in the corollary) $\cls\cI\sE_t=\cls\cI\sE_t\cI$ (upper) \hl{triangular}\index{product system!triangular}, and we call an external direct sum \hl{diagonal}\index{product system!diagonal}. The terminology is clear in the case when $\cB=\cJ\oplus\cI$, because, then, a \nbd{\cB}correspondence $\sE$ decomposes into a matrix $\cls\rtMatrix{\cJ\sE\cJ&\cJ\sE\cI\\\cI\sE\cJ&\cI\sE\cI}$, and the condition in the corollary just means that the matrix is (upper) triangular.
\erem

Recall that $\sF_t$ is not only an \nbd{\cI}correspondence but that it is also a \nbd{\cB}subcorrespondence of $\sE_t$. So, (modulo the replacement of $\sF_0=\cI$ with $\cB$) we may apply Corollary \ref{asubscor}. We collect what these results mean for unitalizations.

\bcor
For a \nbd{C^*}algebra $\cB\ni\U$ let $\wt{\sE}^\podot$ be a superproduct system of \nbd{\wt{\cB}}cor\-re\-spon\-dences, and put $\sE_t:=\U\wt{\sE}_t\U$.
\begin{enumerate}
\item
The $\sE_t$ form a superproduct system $\sE^\podot$ of \nbd{\cB}correspondences.

\item
If $\wt{\sE}^\podot$ is adjointable, then so is $\sE^\podot$.

\item
If 
%%%% BO 
$\wt{\sE}^\podot$ is triangular (that is, if $\U\wt{x}_t=\U\wt{x}_t\U$, for all $\wt{x}_t \in \wt{\sE}_t$) and a product system, 
%%%% EO 
then $\sE^\podot$ is a product system.
\end{enumerate}
\ecor

\brem
It is easy to promote a product system $\sE^\odot$ of \nbd{\cB}correspondences to a diagonal product system of \nbd{\wt{\cB}}correspondences; simply put the trivial Arveson system in the other place of the diagonal. This will, however, not allow to promote contractive units to unital ones. In Section \ref{DilSPSpSEC} (the discussion on unitalizations culminating in Lemma \ref{pswtpslem}) and Section \ref{EXpropsupSEC}, we will see in how far that this \bf{can} be done by embedding into a triangular product system, when unit and product system come from a strong dilation.
\erem

\newpage

\section[\sc{Examples:} Exponentiating (super)(sub)product systems]{Examples: Exponentiating (super)(sub)product systems}\label{EXexpSEC}

To be precise, we wish to \it{exponentiate} discrete \nbd{d}parameter super- or subproduct systems (that is, over $\N_0^d$) to obtain continuous time \nbd{d}parameter super- or subproduct systems (that is, over $\R_+^d$). Many properties of the discrete systems are preserved under \it{exponentiation}, so we get a tool transforming discrete (counter)examples into continuous time (counter)examples.

We first discuss the case of one-parameter product systems, which is well known. We illustrate, without going too much into details, how a slight change of point of view allows us to get immediately the exponential of discrete one-parameter super- and subproduct systems. Then we discuss thoroughly the \nbd{d}parameter case.

The term \it{exponential product system} has been used right from the beginning by Arveson \cite{Arv89} for Arveson systems consisting of symmetric Fock spaces. It generalizes to modules as \it{time ordered product systems}, which are formed by \it{time ordered Fock modules} as introduced by Bhat and Skeide \cite{BhSk00}. We could define and discuss them right here. But as we wish, also in view of later sections, to provide a thorough discussion, we take this out to Appendix \ref{FockAPP}. (Also the ``function'' spaces $L^2(S,F)$, together with an explanation why we'd better put quotation marks and also what the word ``pointwise'' in the sequel means, are discussed there in detail.)

Let us just recall from Appendix \ref{FockAPP} that for a correspondence $F$ over $\cB$ the \phantomsection\hl{time ordered product system}\index{product system!time orered} $\DG^\odot(F)=\bfam{\DG_t(F)}_{t\in\R_+}$ consists of the correspondences
\beqn{
\DG_t(F)
~:=~
\om_t\cB\oplus\bigoplus_{n\in\N}\Delta_nL^2(\RO{0,t}^n,F^{\odot n})
}\eeqn
where $\om_t=\U\in\cB$ (so that the \nbd{0}particle sector is just isomorphic to $\cB$) and where $\Delta_n$ is the indicator function of $\CB{\bbm{a}=(\alpha_n,\ldots,\alpha_1)\colon\alpha_n\ge\ldots\ge\alpha_1\ge0}$ (acting by multiplication as a projection) and with product system structure given by
\beqn{
u_{s,t}(G_m\odot H_n)
~:=~
\BSB{
(\beta_m,\ldots,\beta_1,\alpha_n,\ldots,\alpha_1)
~\longmapsto~
G_m(\beta_m-t,\ldots,\beta_1-t)\odot H_n(\alpha_n,\ldots,\alpha_1)
}
}\eeqn
for $G_m$ in the \nbd{m}particle sector of $\DG_s(F)$ and $H_n$ in the \nbd{n}particle sector of $\DG_t(F)$.

Now we put emphasis on the fact that the pointwise tensor product of functions with values in $F^{\odot m}$ and $F^{\odot n}$, respectively, is just the pointwise product within the discrete product system $\bfam{F^{\odot n}}_{n\in\N_0}$. Effectively:

\bob \label{d1pob}
Every discrete one-parameter product system $F^\odot=\bfam{F_n}_{n\in\N_0}$ is isomorphic to the product system $\bfam{F_1^{\odot n}}_{n\in\N_0}$. To see this concretely, it is convenient to introduce for every product system $E^\odot$ (with product maps $u_{s,t}$) over any monoid the notation $u_{t_n,\ldots,t_1}$ for the $n$th iterated product $E_{t_n}\odot\ldots\odot E_{t_1}\rightarrow E_{t_n\ldots t_1}$. By multiple associativity, this does not depend on how we iterate. For the discrete product system $F^\odot$ we get maps $u_{1,\ldots,1}\colon F_1^{\odot n}\rightarrow F_n$. Again by multiple associativity, these bilinear unitaries form an isomorphism $\bfam{F_1^{\odot n}}_{n\in\N_0}\rightarrow F^\odot$.
\eob

With this observation in mind, we may form for any discrete one-parameter product system $F^\odot$ its \hl{exponential} continuous time one-parameter product system $\DG^\odot(F^\odot)=\bfam{\DG_t(F^\odot)}_{t\in\R_+}$ with
\beq{ \label{GFo}
\DG_t(F^\odot)
~:=~
\om_t\cB\oplus\bigoplus_{n\in\N}\Delta_nL^2(\RO{0,t}^n,F_n)
}\eeq
and product
\beq{ \label{TOproddef}
G_mH_n
~:=~
\BSB{
(\beta_m,\ldots,\beta_1,\alpha_n,\ldots,\alpha_1)
~\longmapsto~
G_m(\beta_m-t,\ldots,\beta_1-t)H_n(\alpha_n,\ldots,\alpha_1)
},
}\eeq
the pointwise product of elements in $F_m$ and $F_n$. Obviously, $\DG^\odot(F^\odot)$ is isomorphic to $\DG^\odot(F_1)$ via the (pointwise) isomorphisms $u_{1,\ldots,1}^*$.

It is the point of view taken in $\DG^\odot(F^\odot)$ that generalizes -- quite obviously -- to discrete one-parameter systems of more general type, and with more work and some extra idea to the \nbd{d}para\-meter case.

Associativity of the product of $\DG^\odot(F^\odot)$ in \eqref{TOproddef} follows simply from associativity of the product of $F^\odot$. (Note, for later generalization to the \nbd{d}parameter case, how easy it is to convince yourself of associativity of the operation on time arguments in this one-parameter case.)  Of course, by the discussion in Appendix \ref{FockAPP}, the product is also unitary because the product of $F^\odot$ is unitary; but for associativity we don't need that. In fact, if the $F_n$ form only a superproduct system, then we are concerned with isometric product maps. 

\brem
Even if the $F_n$ formed an adjointable subproduct system, we would get coisometric product maps, whose adjoints turn the exponential family into a(n adjointable) subproduct system. In fact, the only property of the product maps of $F^\odot$ which we need to define associative product maps for the family with members \eqref{GFo}, is contractivity. We will address some (useful) aspects of greater generality and notation in Remark \ref{copurem}, which concludes this section. But, since we wish to exponentiate also non-adjointable subproduct systems, we need in any case a separate discussion in terms of their coproduct maps.
\erem

After the preceding discussion, the following things, which essentially sum up the one-para\-meter case, are clear:

% \newpage

\bob \label{1pexpob}
{~}

\begin{enumerate}
\item \label{1pexp1}
Let $F^\podot=\bfam{F_n}_{n\in\N_0}$ be a discrete one-parameter superproduct system.  Then the \hl{exponential superproduct system} $\DG^\podot(F^\podot)=\bfam{\DG_t(F^\podot)}_{t\in\R_+}$ with $\DG_t(F^\podot)$ defined in exactly the same way as $\DG_t(F^\odot)$ in \eqref{GFo} and with product defined by \eqref{TOproddef}, is a superproduct system.

\item \label{1pexp2}
$\DG^\podot(F^\podot)$ is a(n adjointable) (proper) superproduct system if and only if $F^\podot$ is a(n adjointable) (proper) superproduct system.

We shall see later in Theorem \ref{disexpthm} that $F^\podot$ is, in some sense, contained in $\DG^\podot(F^\podot)$. From this, the only-if parts are clear. Likewise, the if part for adjointable is clear. Surjectivity in the case of product systems is already known.
\end{enumerate}
\eob

\noindent
We mentioned already that the preceding discussion can be generalized to coisometric products, which would cover adjointable subproduct systems. Of course, also for a not necessarily adjointable discrete one-parameter subproduct system $F^\bodot$, we wish to turn the $\DG_t(F^\bodot)$ defined as  in \eqref{GFo} into a subproduct system, by defining a coproduct in a way that is ``adjoint'' to the product in \eqref{TOproddef}. The point is, that it is quite tedious (as usual with coproducts, which always require the invention of tricky notation) to get hold of what the adjoint of such a product actually does to a typical element. For instance a general ``function'' in $L^2(\RO{0,t+s}^n,F_n)$ will not go into a single $L^2(\RO{0,s}^k,F_k)\odot L^2(\RO{0,t}^{n-k},F_{n-k})$ for a fixed $k$, but we will have to take the direct sum over $0\le k\le n$. In order to avoid that, we we can look at functions of the form
\begin{subequations} \label{TOcodef}
\beq{
(\gamma_n,\ldots,\gamma_1)
~\longmapsto~
\I_{\RO{r'_n,r_n}}(\gamma_n)\ldots\I_{\RO{r'_1,r_1}}(\gamma_1)G_n
}\eeq
with $G_n\in F_n$ and $s+t\ge r_n>r'_n\ge\ldots\ge r_{k+1}>r'_{k+1}\ge t\ge r_k>r'_k\ge\ldots\ge r_1>r'_1=0$, which form a total subset. We send such a function to the function
\bmu{
\Bfam{(s_{n-k},\ldots,s_1),(t_k,\ldots,t_1)}
\\
~\longmapsto~
\I_{\RO{r'_n,r_n}}(s_{n-k}-t)\ldots\I_{\RO{r'_{k+1},r_{k+1}}}(s_1-t)\I_{\RO{r'_k,r_k}}(t_k)\ldots\I_{\RO{r'_1,r_1}}(t_1)w_{n-k,k}G_n
~~~~~~
}\emu
\end{subequations}
in $L^2(\RO{0,s}^k,F_k)\odot L^2(\RO{0,t}^{n-k},F_{n-k})$. This, indeed, extends to an isometry $\DG_{s+t}(F^\bodot)\rightarrow\DG_s(F^\bodot)\odot\DG_t(F^\bodot)$. It is also sufficiently clear that the coassociativity condition is fulfilled.

We do not provide details, but content ourselves with that Equations \eqref{TOcodef} contain a hint how to write things down in a way that is useful to deal with both (super and sub) \nbd{d}parameter cases and contain a hint how to recover a discrete system from its exponential system in Theorem \ref{disexpthm}. Of course, the discussion of \nbd{d}parameter subproduct systems includes a proof of the following analogue of Observation \ref{1pexpob}.

\bob \label{1pexbob}
{~}

\begin{enumerate}
\item \label{1'exp1}
Let $F^\bodot=\bfam{F_n}_{n\in\N_0}$ be a discrete one-parameter subproduct system.  Then the \hl{exponential subproduct system} $\DG^\bodot(F^\bodot)=\bfam{\DG_t(F^\bodot)}_{t\in\R_+}$ with $\DG_t(F^\bodot)$ defined in exactly the same way as $\DG_t(F^\odot)$ in \eqref{GFo} and with coproduct defined by \eqref{TOcodef}, is a subproduct system.

\item \label{1'pexp2}
$\DG^\bodot(F^\bodot)$ is a(n adjointable) (proper) subproduct system if and only if $F^\bodot$ is a(n adjointable) (proper) subproduct system.
 \end{enumerate}
\eob

\lf\noindent
After this sketchy discussion of the one-parameter case, which only serves as motivation, we now come to the \nbd{d}parameter case:

For whatever discrete one-parameter family $F^\#=\bfam{F_n}_{n\in\N_0}$ ($\#=\podot,\bodot, ...$), the \nbd{n}particle sector of $\DG_t(F^\#)$ was defined as $\Delta_nL^2(\RO{0,t}^n,F_n)$. In analogy, for whatever discrete \nbd{d}parameter family $F^\#=\bfam{F_\bn}_{\bn\in\N_0^d}$,  we wish to define the \nbd{\bn=(n_1,\ldots,n_d)}particle sector of $\DG_\bt(F^\#)$ (with $\bt=(t_1,\ldots,t_d)$) as a suitable subspace of
\beqn{
L^2(\RO{0,t_1}^{n_1}\times\ldots\times\RO{0,t_d}^{n_d},F_\bn)
~=~
L^2\RO{0,t_1}^{n_1}\otimes\ldots\otimes L^2\RO{0,t_d}^{n_d}\otimes F_\bn.
}\eeqn
The advantage of factoring everything out into tensor products, is two-fold. Firstly, it allows us to get easily the subspaces of this $L^2$ we really want; secondly, it provides us with the basic tool to reduce (co)associativity of the (co)product we are going to define to (co)associativity of our \nbd{d}parameter system $F^\#$ and (co)associativity of the (co)product in the scalar time ordered system $\DG^\otimes(\C)$. We define
\beqn{
\DG_\bt(F^\#)
~:=~
\bigoplus_{\bn\in\N_0^d}
\Delta_{n_1}L^2\RO{0,t_1}^{n_1}\otimes\ldots\otimes\Delta_{n_d}L^2\RO{0,t_d}^{n_d}\otimes F_\bn,
}\eeqn
with the usual conventions that $L^2\RO{0,t_i}^0=\Om_{t_i}\C$ and $L^2\RO{0,t_1}^0\otimes\ldots\otimes L^2\RO{0,t_d}^0=\Om_\bt\C$ (so that $\Om_{t_1}\otimes\ldots\otimes\Om_{t_d}=\Om_\bt$), and we consider each of these \nbd{\bn}particle sectors as a subset of
\beqn{
\DG_{t_1}(\C)\otimes\ldots\otimes\DG_{t_d}(\C)\otimes F_\bn,
}\eeqn
based on the inclusions $\Delta_{n_k}L^2\RO{0,t_k}^{n_k}\subset\DG_{t_k}(\C)$ in each factor $k=1,\ldots,d$. In the one-parameter case, this would read $\Delta_nL^2\RO{0,t}^n\otimes F_n\subset\DG_t(\C)\otimes F_n$, and the product in \eqref{TOproddef} amounts to
\beq{ \label{1pTpdef}
(g_s\otimes G_m)\odot(h_t\otimes H_m)
~\longmapsto~
g_sh_t\otimes G_mH_n,
}\eeq
where, really, $g_s\otimes h_t\mapsto g_sh_t$ is the product $u_{s,t}$ of the scalar time-ordered product system $\DG^\otimes(\C)$ restricted to elements $g_s$ and $h_t$ in the \nbd{m} and \nbd{n}particle sector of $\DG_s(\C)$ and $\DG_t(\C)$, respectively, and where $G_m\odot H_n\mapsto G_mH_n$ is the product $v_{m,n}$ of elements $G_m$ and $H_n$ in the discrete one-parameter superproduct system in question.

We wish to do the same on the \nbd{d}parameter case. To that goal, we note that for the trivial product system $\C^{d\otimes}:=\bfam{\C}_{\bn\in\N_0^d}$ we have
\beqn{
\DG_\bt(\C^{d\otimes})
~=~
\DG_{t_1}(\C)\otimes\ldots\otimes\DG_{t_d}(\C).
}\eeqn
We observe that the family $\DG^{d\otimes}(\C):=\bfam{\DG_\bt(\C)}_{\bt\in\R_+^d}$ with product defined by
\beqn{
u^d_{\bs,\bt}
\colon
(g_{s_1}\otimes\ldots\otimes g_{s_d})\otimes(h_{t_1}\otimes\ldots\otimes h_{t_d})
~\longmapsto~
g_{s_1}h_{t_1}\otimes\ldots\otimes  g_{s_d}h_{t_d},
}\eeqn
is a product system over $\R_+^d$. (This is nothing but the product of $d$ product systems $\DG^\otimes(\C)$ of Hilbert spaces as explained in the beginning of Section \ref{compSEC} for $d=2$.) With this in mind, returning to the \nbd{d}parameter case,
on $\DG^\podot(F^\podot)=\bfam{\DG_\bt(F^\podot)}_{\bt\in\R_+^d}$ we define the product
\beq{ \label{dTOproddef}
(g_\bs\otimes G_\bm)\odot(h_\bt\otimes H_\bn)
~\longmapsto~
g_\bs h_\bt\otimes G_\bm H_\bn,
}\eeq
where $g_\bs\otimes h_\bt\mapsto g_\bs h_\bt$ is the product $u^d_{\bs,\bt}$ (restricted to the \nbd{\bm} and \nbd{\bn}particle sector of $\DG_\bs(\C^{d\otimes})$ and $\DG_\bt(\C^{d\otimes})$, respectively) and where $G_\bm\otimes H_\bn\mapsto G_\bm H_\bn$ is the product $v_{\bm,\bn}$ of the discrete \nbd{d}parameter superproduct system $F^\podot$.

It is plain that this product is associative. It is also plain that this product sends elements from different pairs $(\bm,\bn)$ (which, therefore, are orthogonal) into different direct summands on the right-hand side. It is also plain, that adjointability is preserved. Finally, it is plain (by the same discussion as for the one-parameter case) that the product is surjective if (and, of course, only if) the product of $F^\podot$ is surjective. We, therefore, proved:

\bthm \label{dpexpthm}
Let $F^\podot=\bfam{F_\bn}_{\bn\in\N_0^d}$ be a  discrete \nbd{d}parameter superproduct system.  Then $\DG^\podot(F^\podot)$ with the product maps defined by \eqref{dTOproddef} is a continuous time \nbd{d}parameter superproduct system, the \hl{exponential superproduct system}\index{superproduct system!exponential (of a superproduct system over $\N_0^d$)}\index{systems!exponential!superproduct system} of $F^\podot$.

$\DG^\podot(F^\podot)$ is a(n adjointable) (proper) superproduct system if and only if $F^\podot$ is a(n adjointable) (proper) superproduct system.
\ethm

We now wish to promote that theorem to subproduct systems. There are several problems that make this case more tricky.

Firstly, the product maps $u_{s,t}$ (and the like) in the case of superproduct systems have the big notational advantage that we can really work with the notation of product $g_sh_t=u_{s,t}(g_s\otimes h_t)$. It is this product notation that allowed us in \eqref{1pTpdef} to capture the product of $g_s\otimes G_m$ and $h_t\otimes H_n$ as $g_sh_t\otimes G_mH_n$. But actually, when trying to write down the product map with the help of $u_{s,t}\otimes v_{n,m}$, we do have to involve the canonical flip that identifies $(\DG_s(\C)\otimes F_m)\odot(\DG_t(\C)\otimes F_n)$ with $(\DG_s(\C)\otimes\DG_t(\C))\otimes(F_m\odot F_n)$ by exchanging the 2nd and the 3rd site. For the coproduct in \eqref{TOcodef} we ``made disappear'' the problem by means of the function picture. ($L^2(\Om,F)$ does not ``see'' if we think of is as $L^2(\Om)\otimes F$ or as $F\otimes L^2(\Om)$.) Now since we want to exploit the tensor product picture of \nbd{L^2}spaces, we can no longer really avoid the flip.

Secondly, the coproduct does not generally map an element to a simple tensor. Fortunately, at least for the part $u_{s,t}^*$ of the coproduct, we can neutralize that problem by evaluating that part on products $g_sh_t$. The indicator functions in \eqref{TOcodef} exploit this.

Thirdly, and most severely, we are not actually looking at the product $u_{s,t}$ on its full domain $\DG_s(\C)\otimes\DG_t(\C)$ but only on the \nbd{n} and \nbd{m}particle sector, and which $n$ and which $m$ depends on the second factor $F_n$ and $F_m$, respectively. Even if we write down elements $g_s$ and $h_t$ such that the product is in a fixed \nbd{n}particle sector of $\DG_{s+t}(\C)$, it is not at all so that necessarily the factors $g_s$ and $h_t$ are in fixed number of particles sector. Again, the indicator functions  in \eqref{TOcodef} do fulfill this requirement.

After explaining, basically in the one-parameter case, the problems in the passage from product to coproduct and how they have been resolved by \eqref{TOcodef}, we turn to a discrete \nbd{d}para\-meter subproduct system $F^\bodot$. For fixed $\bs,\bt\in\R_+^d$ and $\bbm{k}\in\N_0^d$, observe that the elements
\beqn{
g_\bs h_\bt\otimes G_{\bbm{k}}
}\eeqn
with $\bn+\bm=\bbm{k}$, with $g_\bs$ in the \nbd{\bm}particle sector of $\DG_\bs(F^\bodot)$, with $h_\bt$ in the \nbd{\bn}particle sector of $\DG_\bt(F^\bodot)$, and with $G_{\bbm{k}}\in F_{\bbm{k}}$, are total in the \nbd{\bbm{k}}particle sector of $\DG_{\bs+\bt}(F^\bodot)$. Moreover, for different choices of $\bm+\bn=\bbm{k}$ they are orthogonal. So varying $\bm$ and $\bn$ freely (with $\bbm{k}:=\bm+\bn$) we get a total subset of $\DG_{\bs+\bt}(F^\bodot)$. The map
\beq{ \label{dTOcodef}
g_\bs h_\bt\otimes G_{\bm+\bn}
~\longmapsto~
\f_{2,3}(g_\bs\otimes h_\bt\otimes w_{\bm,\bn}G_{\bm+\bn}),
}\eeq
where $\f_{2,3}$ is the flip of the 2nd and 3rd tensor site (out of four, since $w_{\bm,\bn}G_{\bm+\bn}\in F_\bm\odot F_\bn$), is our definition of the coproduct. After having been very explicit about  what is sorted into which direct summand, associativity and the other statements in the following theorem are clear.

\bthm \label{dbexpthm}
Let $F^\bodot=\bfam{F_\bn}_{\bn\in\N_0^d}$ be a  discrete \nbd{d}parameter subproduct system.  Then $\DG^\bodot(F^\bodot)$ with the coproduct maps defined by \eqref{dTOcodef} is a continuous time \nbd{d}parameter subproduct system, the \hl{exponential subproduct system}\index{subproduct system!exponential (of a subproduct system over $\N_0^d$)}\index{systems!exponential!subproduct system} of $F^\bodot$.

$\DG^\bodot(F^\bodot)$ is a(n adjointable) (proper) superproduct system if and only if $F^\bodot$ is a(n adjointable) (proper) superproduct system.
\ethm

We invented the machinery of exponentiating discrete systems, to be able to transform discrete (counter) examples into continuous time (counter) examples. Two of the most important questions we meet in these notes, are:  Does a given subproduct system embed into a superproduct system? Does a given superproduct system embed into a product system? (Recall the discussions around  Definitions \ref{subdefi} and \ref{moemisdef}, where subsystems and embeddings are defined.) We will now show that the answers to both questions coincide for a discrete \nbd{d}parameter system and for its exponential system. For one direction we will show that embeddings lift to embeddings. For the other direction we will recover the original system ``sitting inside'' appropriately in its exponential system.

Forgetting about whether super or sub (or other structures; see Remark \ref{copurem}), for every discrete \nbd{d}parameter family $F^\circ=\bfam{F_\bn}_{\bn\in\N_0^d}$ we may define its exponential continuous time  \nbd{d}parameter family $\DG^\circ(F^\circ)=\bfam{\DG_\bt(F^\circ)}_{\bt\in\R_+^d}$, precisely as we did for super- and subproduct systems. For two such families $F^\circ$ and $F'^\circ$ and a family $a^\circ=\bfam{a_\bn}_{\bn\in\N_0^d}$ of maps $a_\bn\in\sB^{bil}(F_\bn,F'_\bn)$, we wish to define its \phantomsection\hl{second quantized family}\index{morphisms!second quantized (on an exponential system)} $\DG^\circ(a^\circ)=\bfam{\DG_\bt(a^\circ)}_{\bt\in\R_+^d}$ by letting act $\DG_\bt(a^\circ)\in\sB^{bil}(\DG_\bt(F^\circ),\DG_\bt(F'^\circ))$ as pointwise multiplication by $a_\bn$ on the \nbd{\bn}particle sector of $\DG_\bt(F^\circ)$. This is possible  if and only if the norms $\snorm{a_\bn}$ are bounded uniformly in $\bn$. Note that second quantization preserves properties like adjointability, contractivity, isometricity, coisometricity, and, therefore, unitarity. (We have already argued several times, why in the latter two cases surjectivity is preserved.) It is routine to convince ourselves that when $\#,\#'\in\CB{\podot,\bodot}$ and $F^\circ=F^\#$, $F'^\circ=F'^{\#'}$, the family $\DG^\circ(a^\circ)$ is a morphism if (and only if) $a^\circ$ is a morphism.
% \OW[MICHAEL]{Better check once more the four cases. (I checked! --Orr)}
So, an embedding of discrete systems (no matter whether that of a subproduct system into a superproduct system, or that of a superproduct system into a product system, or any other sort of embedding) gives rise to an embedding of its exponentials.

Now let us try to recover $F^\circ$ from $\DG^\circ(F^\circ)$. Let us start with two observations: Firstly, every continuous time \nbd{d}parameter system $E^\circ=\bfam{E_\bt}_{\bt\in\R_+^d}$ (of a type $\podot,\bodot,\ldots$) gives rise to a discrete \nbd{d}parameter system $\bfam{E_\bt}_{\bt\in\N_0^d}$ (of the same type) simply by restricting the indices (and operations) to integer times.

\brem \label{homomrem}
This is a special case of the observation that for any monoid morphism $\vp\colon\bS'\rightarrow\bS$ we can turn a system over $\bS$ into a system over $\bS'$ in the obvious way. Here the homomorphism is the canonical embedding of $\N_0^d$ into $\R_+^d$; this was a crucial observation in the construction of an \nbd{E_0}semigroup for every Arveson system in Skeide \cite{Ske06} in the one-parameter case. A non-injective, but surjective example is the \nbd{d}fold addition map $\N_0^d\rightarrow\N_0$ which turns the discrete one-parameter system $E^{\odot n}$ with just tensor product as multiplication, into a discrete \nbd{d}parameter product system $E_\bn=E^{\odot(n_1+\ldots+n_d)}$ with just tensor product as multiplication; in the Example-Section \ref{EXsubnsupSEC}, this corresponds to the trivial (but untypical!) case, where all flips $\sF_{j,i}$ are the identity.
\erem

Secondly, every $\DG_\bt(F^\circ)$ contains many copies of each $F_\bn$. More precisely, for every unit vector $\xi_\bn$ in the \nbd{\bn}particle sector of $\DG_\bt(\C^{d\otimes})$, the map is $x_\bn\mapsto\xi_\bn\otimes x_\bn$ is an isometry onto a (complemented!) copy of of $F_\bn$. We wish to choose the $\xi_\bn$ in such a way that a product/coproduct of the discrete exponential system $\bfam{\DG_\bt(F^\circ)}_{\bt\in\N_0^d}$ restricted to the copy $\xi_\bn\otimes F_\bn$ gives back the product/coproduct of $F^\circ$.

The best way (and, thinking carefully about it, the only way) to achieve this, is to seek discrete units that sit in the \nbd{\bn}particle sectors of the scalar discrete exponential product system $\bfam{\DG_\bt(\C^{d\otimes})}_{\bt\in\N_0^d}$. (The scalar exponential product system has lots of (continuous time) units. But these units are (multiples of) exponential vectors and, except for the vacuum, exponential vectors never sit in a single \nbd{\bn}particle sector. The fact that we have to seek discrete units that do not extend to continuous time units, is, again, very much reminding of \cite{Ske06}.) Once we have such a discrete (unital) unit $\xi_\bn$, so that $u_{\bm,\bn}(\xi_\bm\otimes\xi_\bn)=\xi_{\bm+\bn}$, respectively, $u_{\bm,\bn}^*\xi_{\bm+\bn}=\xi_\bm\otimes\xi_\bn$, it is clear that $x_\bn\mapsto\xi_\bn\otimes x_\bn$ defines an (adjointable!) embedding of the discrete super/subproduct system $F^\#$ into the discrete super/subproduct system $\bfam{\DG_\bt(F^\#)}_{\bt\in\N_0^d}$. Equations \eqref{TOcodef} give the hint. We choose:
\beqn{
\xi_\bn
~:=~
(\I_{\RO{n_1-1,n_1}}\otimes\ldots\otimes\I_{\RO{0,1}})\otimes\ldots\otimes(\I_{\RO{n_d-1,n_d}}\otimes\ldots\otimes\I_{\RO{0,1}}).
}\eeqn
We, thus, have proved:

\bthm \label{disexpthm}
Every discrete \nbd{d}parameter super/subproduct system $F^\#$ is contained in its discrete exponential super/subproduct system $\bfam{\DG_\bt(F^\#)}_{\bt\in\N_0^d}$ as a super/subproduct subsystem.
\ethm

\bcor\label{embedFcor}
$F^\#$ embeds into some $F'^{\#'}$ if and only if $\DG^\#(F^\#)$ embeds into some $E^{\#'}$.
\ecor

\proof
The forward implication follows by second quantizing the embedding of $F^\#$ into $F'^{\#'}$.

For the backwards operation, suppose $\DG^\#(F^\#)$ embeds into $E^{\#'}$. Then $F^\#$, isomorphic to some subsystem of $\bfam{\DG_\bt(F^\#)}_{\bt\in\N_0^d}$, embeds into the discrete system  $E_{\N_0^d}^\#=\bfam{E_\bt}_{\bt\in\N_0^d}$, too.\qed

% \lf
% Recall that, by Proposition \ref{submoprop}, this includes the answer to the question when a superproduct system embeds into a product system.

% \newpage

\bulletline
We conclude this section with two lengthy remarks. They both point to future work on exponential systems or work where exponential systems play a role. The second remark also addresses the problem that after almost three decades of product systems (five decades, if we take into account what started with Streater \cite{Str69}, Araki \cite{Ara70}, Parthasarathy and Schmidt \cite{PaSchm72}, and Guichardet \cite{Gui72}, in all of which Fock type product systems occur without really denominating the structure) the terminology could benefit from more consequent choices.

\brem
Adjointable (time continuous) morphisms between one-parameter time ordered product systems $\DG^\odot(F)\rightarrow\DG^\odot(F')$ have been characterized in Barreto, Bhat, Liebscher, and Skeide \cite[Section 5.2]{BBLS04}. It is most convenient to allow for unbounded morphisms%
% $a_t\in\sL^{a,bil}(\DG_t(F),\DG_t(F'))$
, because their parametrization is simpler. More specifically, every morphism sends units to units; so, necessarily a morphism leaves invariant the algebraic time ordered product system generated algebraically by the units with product maps restricted to the algebraic tensor product. (\hl{Time continuous} means just that the morphism (bounded or not) sends \it{continuous} units to continuous units; see Appendix \ref{FockAPP} for how the continuous units are made. It is noteworthy, that time continuous endomorphisms of $\DG^\odot(F)$ correspond exactly to the strongly continuous local cocycles for the CCR-flow of which $\DG^\odot(F)$ is the associated product system.) By \cite[Theorem 5.2.1]{BBLS04}, the adjointable, possibly unbounded, time continuous morphisms are characterized one-to-one by matrices
\beqn{
\sMatrix{\gamma&\eta^*\\\eta'&\alpha}
~\in~
\sB^{a,bil}\sMatrix{\sMatrix{\cB\\F},\sMatrix{\cB\\F'}}.
}\eeqn
(We think this remains true also for not necessarily adjointable morphisms, if we replace $\sB^{a,bil}$ with $\sB^{bil}$.) In the adjointable case, \cite[Corollary A.7]{BBLS04}, which characterizes the contractive positive morphisms on $\DG^\odot(F)$ (when applied to the positive morphism $a_t^*a_t$) also will provide a characterization of the contractive morphisms. (In the non-adjointable case, one should switch to von Neumann modules.) Contractive for morphisms is almost the same as bounded. (Only if $\snorm{\alpha^{\odot n}}$ should remain uniformly bounded in $n$, there are bounded morphisms that are not contractions. For instance, if $F$ (or $F'$) is \phantomsection\hl{nilpotent}\index{correspondence!nilpotent} to some degree, meaning $F^{\odot n}=\zero$ for some $n$, just all morphisms are bounded and there are morphisms that are not contractive.)

The morphisms that are second quantizations, correspond precisely to matrices where only $\alpha$ is different from $0$. This shows that we may not hope to get all morphisms as second quantizations.

The natural question to characterize morphisms between all sorts of exponential systems and to find out to what extent these can be characterized in a similar manner, raises a completely new area of questions. It is important to note that the ``parameters'' $F$ and $F'$ relate to a description of the exponential product systems in an ``infinitesimal way''. This possibility to put, in the one-parameter case, everything into a single correspondence, is strictly related to the simple structure of discrete one-parameter product systems; see Observation \ref{d1pob}. Already in the case of \nbd{d}parameter \bf{product} systems, the $F_i$ for each marginal system may but need not be the same $F$ for all $i=1,\ldots,d$. In order to get hold of a suitable $F$ (even in the one-parameter case), the question whether or not a given exponential super/subproduct system embeds into a product system, in the affirmative case, will be paired with the natural question, if this product system can be chosen exponential. There are more related questions that are formulated and under investigation by several authors for Hilbert spaces, and only partial answers to them exist; for correspondences, difficulties will still increase even more. A very much related task would be to characterize exponential systems intrinsically. For product systems this is done by type-classification. (Arveson systems are exponential, if they are \hl{type I}, that is, if they are generated by their units. For \nbd{C^*}correspondences, they have to be also \hl{spatial}=existence of a central unital unit; see Skeide \cite{Ske06d}. For von Neumann correspondences, by \cite{BBLS04} we get back that type I is enough, while in the \nbd{C^*}case, by \cite{BLS10} this may fail.) We will need to understand what this means for super/subproduct systems. Some aspects, those related directly to units, will be addressed in the following remark.
\erem

\brem \label{copurem}
In this conclusive remark we address the problem that the -- now well-estab\-lished and, therefore, unchangeable -- terminology \it{product system} together with already two ramifications, \it{superproduct system} and \it{subproduct system}, can cause sometimes a certain head\-ache.

We are used to understand terminologies such as \it{super-} and \it{submartingale} as, usually complementary, generalizations of the common part \it{martingale}. Here, \it{complementary} is to be understood in the sense that having both super and sub, what we get is a martingale. So far, this is fine also with product systems and its super and sub versions. However, this is almost all where the ramification of product systems into super and sub parallels the corresponding ramification for martingales.

It is easy to transform a super-martingale into a sub-martingale, or \it{vice versa}, simply by multiplying it with $-1$. Correspondingly, the theory of super-martingales and the theory of sub-martingales coincide. (Well, it gets more complicated if additional properties, like positive super- and submartingales, are taken into consideration. But, let us consider these as second order effect.) There is no such way to transform a superproduct system into a subproduct system; and \it{vice versa}. However, an adjointable superproduct system can be transformed into something with a coproduct structure. Still, the adjoints of the product maps of a superproduct system do not turn the superproduct system into a subproduct system. (For that, the adjoints should be isometries, too, so that we actually have a product system to begin with.) However, the adjoints \bf{are} coproduct maps that fulfill everything required in the definition of subproduct system, but isometricity. The same goes with the adjoints of the coproduct maps of an adjointable subproduct system, which fulfill everything required from a product of a superproduct system, but isometricity. Therefore, if we drop the condition of isometricity, at least for adjointable product/coproduct maps, we get a sort of duality.

Since the name \it{product system} is already occupied, we propose the following definition -- for future work, only; it will be used only in the remainder of this section and not anywhere else in these notes. (Gerhold and Skeide \cite{GeSk20b} adopted it.)

\bdefin \phantomsection
Let $E^\circ=\bfam{E_t}_{t\in\bS}$ be a family over $\bS$ of \nbd{\cB}correspondences $E_t$.

If there is a family of maps $\pi_{s,t}\in\sB^{bil}(E_s\odot E_t,E_{st})$ fulfilling everything but, possibly, isometricity from the definition of superproduct system, we say $E^\circ$ is a \hl{productive}\index{systems!productive and coproductive} system and denote it by $E^<$.

If there is a family of maps $\gamma_{s,t}\in\sB^{bil}(E_{st},E_s\odot E_t)$ fulfilling everything but, possibly, isometricity from the definition of subproduct system, we say $E^\circ$ is a \hl{coproductive} system and denote it by $E^>$.

A (co)productive system with contractive (co)product maps is called a \hl{contractive}\index{systems!productive and coproductive!contractive} (co)prod\-uct system.
\edefin

Note that we canceled the ``-ive'' from the contractive version. Since \it{product system} does already exist, it should be clear that a \it{contractive product system} must be something different. And we accept that the difference actually consists not in being more special but in being more general (like super and sub).

It is noteworthy that (without doing the effort of an explicit formulation) Theorems \ref{dpexpthm} and \ref{dbexpthm} remain true for contractive (co)product systems, having contractive exponential systems. (Apart from the obvious if-and-only-if statements regarding preservation if isometricity and coisometricity of the (co)product maps, we may add a statement that density of the images under the (co)product is preserved.) Contractivity is not strictly necessary. However, like for second quantized morphisms, we have to control the norms of the discrete (co)products uniformly, and contractivity is the easiest way to do this; and frequently, boundedness is equivalent to contractivity.

For adjointable systems we had to prove only one version of Theorems \ref{dpexpthm} and \ref{dbexpthm}, because the other version follows by taking adjoints. So, we get the impression that for adjointable systems, the two theories fall together. But also this is not entirely correct. The point is that the notions of morphisms and units are different. If we have a morphism between adjointable productive systems, then it is the adjoint of that morphism (provided it exists) that is a morphism between the adjoint coproductive systems, obviously in the opposite direction. There is, in general, no way to show that the morphism itself is a morphism also for adjoint structure. We dispense with a discussion of the mixed versions like morphism from a productive system into a coproductive system. We also dispense with a further discussion the differences, and refer the reader to Bhat and Mukherjee \cite{BhMu10}, where several types of morphisms between (continuous time one-parameter) subproduct systems of Hilbert spaces are discussed. (Up to four types have to be distinguished, corresponding to the $2\times 2$ possibilities to look at the productive structure or at the coproductive structure of domain and codomain of the potential morphism.) We only have a slightly closer look at multiplicative and and comultiplicative units.

We defined units for a superproduct system as multiplicative units and units for subproduct system as comultiplicative units. (It would be more consequent to call comultiplicative units \it{counits}. But in these notes we are interested only in the two described cases, so we prefer not to bother the reader, except for the discussion in this remark.) This is in both cases the strongest possibility, because in both cases the isometric map, that is, the map that does not making loose information, has been used. Additionally, for the non-adjointable systems in both cases, this is the only map available. When we switch to general adjointable (co)productive systems, no such property as isometry singles out what is the better notion. Again, when we drop adjointability, there is only one possibility left.

Anyway, why do we worry about multiplicative units and comultiplicative units in one fixed adjointable system, when it is clear what would be the natural notion, namely, (co)multiplicative units for a (co)productive system?  The reason is that there is already an example where both sorts of units occur and where it is actually more important to look at the ``wrong'' sort of units. In fact, the comultiplicative units for subproduct system are units also in any containing product system; but even if the containing product system is the \it{generated} one (Theorem \ref{indlimthm}), we do not get all units in that way. On the other hand, when we project units of a containing product system down to the subproduct system, what we get are multiplicative units. Bhat and Mukherjee \cite[Theorem 10]{BhMu10} have shown that for (time continuous one-parameter) subproduct systems of Hilbert spaces, the multiplicative units of the subproduct system (modulo a growth condition) are in one-to-one correspondence with the units of the generated Arveson system. (In \cite{BhMu10}, comultiplicative and multiplicative units are called strong and weak, respectively, units.) It is not known if this generalizes (appropriately) to (one-parameter) subproduct systems of correspondences. (``Appropriately'', would surely include that the subproduct system has to sit complementedly in the containing product system, while in that case, the requirement that the subproduct system be adjointable (so that, \it{multiplicative unit} is meaningful), is automatic by Corollary \ref{asubscor}.) We may also ask if the results from \cite{BhMu10} generalize to subproduct systems sitting in superproduct systems, or how everything generalizes beyond the one-parameter case. Another question would be, if there are analogue results for superproduct systems sitting in (sub)product systems. The latter question is part of a whole set of question where we know (at least, partial) answers for subproduct systems, while analogue questions for superproduct systems have never been asked. We see a lot of potential, here.

Last but not least, the notion of (co)productive systems (where they arise) will lead to the question whether they can be \it{dilated} to product systems or, possibly, dilated in stages to an intermediate super/subproduct system. An embedding of a (co)productive system into a (co)product system (in the sense that the (co)product of the embedded systems is recovered by restriction of the unitary (co)product of the containing system) can never be obtained, if the to be embedded system does not have an isometric (co)product. So, how does embedding of an adjointable subproduct system into a a product system read, when we consider the subproduct system as a contractive product system? The answer is, that the embedded family does not take its (coisometric) product  by restriction of the product of the product system, but that we have to project the product down by a family of projections onto the subsystem (as is done in \cite[Definition 5.1]{ShaSo09}). It is easy to capture the properties of this family of projections, and to present the outcome as  definition of what we would call a dilation of a contractive product system to a  product system. We do not go any further, but, we wish to mention that we do have at least two (quite different) frameworks where proper contractive coproduct systems do occur quite naturally.

A similar discussion can be done for the suggested morphisms. (Except for the two structure preserving sorts of morphisms, one for the multiplicative was to to look at the system and one for comultiplicative way, there is another way, called weak morphism in \cite{BhMu10}, that is weaker than both.) We dispense with giving details. But the whole discussion shows how, as frequently happens when new mathematics is produced, that terminology created when the first and most special instance of a structure occurs, is adopted to mean only the special case. Some generalizations (attaining the same name for a more general structure) work well and are intuitive (like the generalization of product system from Hilbert spaces to correspondences), and others would feel awkward (like if we would call a productive system just a product system). Of course, the terminology \it{productive system} is provocative; this is intentional. We use it (in these notes, at least) only in this section.
\erem

\newpage

\section{CP-Semigroups and subproduct systems} \label{CPspsSEC}

As we mentioned already, the notion of subproduct system (or inclusion system) arises from CP-semigroups, by simply looking a the GNS-construction of the individual members of the CP-semigroup; 
%%%% BO
this was observed in \cite[Section 3]{ShaSo09}. Actually, Shalit and Solel mostly considered 
%%%% EO 
a construction that is \it{dual} to the GNS-construction and leads to the so-called \it{Arveson-Stinespring correspondence}. This approach is limited to von Neumann algebras. We discuss the relation in Appendix \ref{vNAPP}. While the subproduct system of GNS-correspondences of a CP-semigroup does not require a construction -- accepting that there is the GNS-construction, the subproduct system is simply there --, we cannot derive every subproduct system as a GNS-subproduct system. To get an analogue of the result \cite[Corollary 2.10]{ShaSo09} that every subproduct system of von Neumann correspondences is the Arveson-Stinespring subproduct system of a (normal) CP-semigroup, we have to seek in the (strict) Morita equivalence class of subproduct systems. This is what we do later in this section. We start with discussing the GNS-subproduct system.

Recall that if $E^\odot$ is a (sub)product system of \nbd{\cB}correspondences over a monoid $\bS$ with a (contractive) unit $\xi^\odot$ (so that, by definition of unit, $\cB$ is unital), then $T_t:=\AB{\xi_t,\bullet\xi_t}$ defines a (contractive) CP-semigroup on $\cB$ over the monoid $\bS^{op}$. (The computation for units in subproduct systems is the same as the computation in \eqref{unitcomp} for units in product systems we did after Definition \ref{unitdef}, and also holds for superproduct systems.)

Now suppose we start with a CP-semigroup $T=\bfam{T_t}_{t\in\bS^{op}}$ over $\bS^{op}$. We do the GNS-construction $(\sE_t,\xi_t)$ for each $T_t$. Then,
\beqn{
w_{s,t}
\colon
\xi_{st}
~\longmapsto~
\xi_s\odot\xi_t
}\eeqn
defines a bilinear isometry $\sE_{st}=\cls\cB\xi_{st}\cB\rightarrow\sE_s\odot\sE_t$. The family $\sE^\bodot=\bfam{\sE_t}_{t\in\bS}$ with the structure maps $w_{s,t}$ is, clearly, a subproduct system over $\bS$, and the elements $\xi_t\in\sE_t$ form a unit $\xi^\odot$ for $\sE^\bodot$. Obviously, the pair $(\sE^\bodot,\xi^\odot)$ is determined uniquely up to unit-intertwining isomorphism by the CP-semigroup $T$. We call $\sE^\bodot$ the \phantomsection\hl{GNS-subproduct system}\index{GNS-subproduct system}\index{systems!subproduct!GNS-subproduct system of a CP-semi\-group}\index{subproduct system!GNS-subproduct system of a CP-semi\-group} \index{CP-semigroup!GNS-subproduct system of}of $T$ with \hl{cyclic unit}\index{GNS-subproduct system!cyclic unit}\index{unit!for a subproduct system!cyclic} $\xi^\odot$.

\bob \label{uGNSsupob}
If $T$ is the CP-semigroup arising from a unit $\xi^\odot$ in a (super)product system, then the subproduct subsystem generated, by the unit as in Example \ref{unitgenex}, is (isomorphic to) the GNS-subproduct system. We come back to this in Section \ref{DilSPSpSEC}.
\eob

\bob \label{elemob}
It is clear that subproduct system with cyclic unit and CP-semigroup determine each other. (We exploited this already in Example \ref{unitgenex}.) The more important it is to note that the subproduct system alone tells something, yes, but not everything about the CP-semigroup. For instance, if $\bfam{c_t}_{t\in\bS}$ is a (contraction) semigroup in $\cB$, then the GNS-subproduct system of the elementary CP-semigroup\index{GNS-subproduct system!of an elementary CP-semigroup}\index{CP-semigroup!GNS-subproduct system of!elementary}\index{elementary!CP-semigroup!GNS-subproduct system of} $T_t:=c_t^*\bullet c_t$ (over $\bS^{op}$!) consists of the ideals $\cls\cB c_t\cB$ generated by $c_t$ (with unit $\xi_t=c_t$ and embeddings $c_{st}\mapsto c_s\odot c_t=c_sc_t\in \sE_s\odot\sE_t=\cls\cB c_s\cB c_t\cB$ also being an ideal in $\cB$; see also Examples \ref{nonadex} and \ref{cnonadex}). If $\cB$ is simple (and $c_t\ne0$ for $t\ne0$), then this GNS-subproduct system is just the trivial \bf{product} system. Conversely, if the GNS-subproduct systems is the trivial product system (or a subproduct subsystem of the trivial one), then the unit is necessarily formed by elements $c_t$ of a semigroup in $\cB$ giving, therefore, back $T_t$ as $c_t^*\bullet c_t$. So, the (nonzero) elementary CP-semigroups on a simple \nbd{C^*}algebra are exactly those that have the trivial product system as GNS-subproduct system. But only the unit tells us which is the CP-semigroup.
\eob

\brem
Existence of a unit depends on that $\cB$ is required unital. We like to mention, however, that the construction of the GNS-correspondences does not. In fact, the GNS-corre\-spondence $E$ of a CP-map $T$ from $\cA$ to $\cB$ is obtained by defining on the \nbd{\cA}\nbd{\cB}bimodule $\cA\otimes\cB$ the (semi)inner product $\AB{a\otimes b,a'\otimes b'}=b^*T(a^*a')b'$ and doing (Hausdorff) completion. (Only for getting the cyclic vector as the image of $\U\otimes\U$ we would need units in the algebras.) Note that for all $a$ and all approximate units $\bfam{u_\lambda}_{\lambda\in\Lambda}$ for $\cB$ the limit
\beqn{
a\otimes\U
~:=~
\lim_\lambda a\otimes u_\lambda
}\eeqn
exists in $E$ and does not depend on the choice of $\bfam{u_\lambda}_{\lambda\in\Lambda}$. (Attention: Nothing like this is true for $v_\mu\otimes b$, if $T$ is not strict!) Now if $S$ is a CP-map from $\cB$ to $\cC$ we get a \nbd{\cB}\nbd{\cC}correspondence $F$ in the same way, and we also get an \nbd{\cA}\nbd{\cC}correspondence $G$ from the CP-map $S\circ T$. We want to define a map from $G$ into $E\odot F$ by
\beqn{
a\otimes c
~\longmapsto~
\lim_\lambda(a\otimes\U)\odot(u_\lambda\otimes c).
}\eeqn
Again, one easily checks that the limit exists and does not depend on the choice of the approximate unit. (Essentially, because the statement is true for the elements $bu_\lambda\otimes c$ of $F$ for every $b$.) Clearly, the map is a bilinear isometry. Iterations behave associatively. In this case, of course, a CP-semigroup $\cB$ gives a subproduct system $\sE_t$ together with the family of cyclic linear maps $b\mapsto b\otimes\U\in\sE_t$ and is determined by them.
\erem

Let us return to unital \nbd{C^*}algebras. It is clear that not every subproduct system is the GNS-subproduct system of a CP-semigroup. For instance,
%%%% BO 
the GNS-subproduct systems consisting of (nonzero) Hilbert spaces are necessarily the ones coming from (nonzero) CP-semigroups on $\C$, and thus are all one-dimensional. 
%%%% EO
But there are product systems (\it{a fortiori}, there are subproduct systems) of Hilbert spaces that are not one-dimensional. Even worse, there are unit-less \it{Arveson systems} (that is, continuous one-parameter product systems of separable Hilbert spaces), which, therefore, do not contain a single GNS-subproduct subsystem.

On the other hand, Arveson systems and other subproduct systems of Hilbert spaces do arise from normal CP-semigroups on $\sB(H)$ where $H$ is a Hilbert space, that is, a Hilbert \nbd{\C}module. Several aspects of the following construction of subproduct systems of \nbd{\cB}correspondences from a \phantomsection\hl{strict}\index{CP-semigroup!strict on $\sB^a(E)$} CP-semigroup $T$ on $\sB^a(E)$ where $E$ is a full Hilbert \nbd{\cB}module (strict means that each $T_t$ is a strict map), have appeared over the years in several papers.

\brem \label{CPASrem}
The Arveson system of a normal (one-parameter) CP-semigroup on $\sB(H)$ has been defined by Bhat \cite{Bha96} through dilation, as that of the dilating \nbd{E}semigroup in the unique minimal dilation. Skeide \cite{Ske03c} has constructed directly from the GNS-subproduct system of that CP-semigroup a subproduct system of Hilbert spaces; the Arveson system generated by that subproduct system coincides with Bhat's Arveson system. Calculating this subproduct system for Powers' sum of spatial CP-semigroups and comparing it with the one used in Skeide \cite{Ske06} to construct the \it{product} of spatial product systems,  the Arveson system of the Powers sum is identified with the spatial product of the Arveson systems of the summands. Bhat and Mukherjee \cite{BhMu10} recovered the construction of the spatial product as a special case of their \it{amalgamated} product of Arveson systems. The construction of the subproduct system of Hilbert spaces for a CP-semigroup in $\sB(H)$, has been generalized to strict (still one-parameter) CP-semigroups on $\sB^a(E)$ by Bhat, Liebscher, and Skeide \cite{BLS08,Ske10}. The key operation of this construction is a strict version of Morita equivalence for correspondences, which is also at the heart of the construction of product systems out of \nbd{E}semigroups (Skeide \cite{Ske09}) and the representation theory of $\sB^a(E)$ (Muhly, Skeide, and Solel \cite{MSS06}). \it{Morita equivalence} for correspondences is from Muhly and Solel \cite{MuSo00} and for product systems from Skeide \cite{Ske09}. But, actually, Morita equivalence of correspondences is a functor that respects tensor products and inclusions, so it makes sense also for super- and subproduct systems; that is what we explain in the sequel.
\erem

For the balance of this section, we do \bf{not} require \nbd{C^*}algebras to be unital.

A correspondence $M$ from $\cB$ to $\cC$ is a \phantomsection\label{MEdef}\hl{Morita equivalence}\index{Morita equivalence}\index{correspondence!Morita equivalence} if there is an \hl{inverse}\index{correspondence!inverse} correspondence $N$ from $\cC$ to $\cB$ such that $M\odot N\cong\cB$ and $N\odot M\cong\cC$. (This definition is proposed in Skeide \cite[Section 2]{Ske16} and we follow it for the whole discussion. \cite{Ske16} also shows that it is equivalent to one of the usual definitions (Lance \cite[Section 7]{Lan95}), and that we may (and, of course, will) choose the isomorphisms such that they are compatible with $(M\odot N)\odot M=M=M\odot(N\odot M)$ and with $(N\odot M)\odot N=N=N\odot(M\odot N)$.) If $E$ is a \nbd{\cB}correspondence, then we may define the \phantomsection\hl{Morita equivalent}\index{Morita equivalence!of correspondences over a \nbd{C^*}algebra}\index{correspondence!Morita equivalence of} \nbd{\cC}correspondence $F:=N\odot E\odot M$.

A full Hilbert \nbd{\cB}module $E$ is a Morita equivalence from the compacts $\sK(E)$ to $\cB$. In fact, the inverse of $E$ is given by the \phantomsection\hl{dual}\index{correspondence!dual} correspondence $E^*:=\CB{x^*\colon x\in E}=\sK(E,\cB)\subset\sB^a(E,\cB)$ with inner product $\AB{x^*,y^*}:=xy^*$ and bimodule operations 
%%%% BO
%$bx^*a:=(b^*xa^*)^*$. 
$bx^*a:=(a^*xb^*)^*$. 
%%%% EO
Moreover, we shall always use the canonical identifications $E^*\odot E=\cB$ via
%%%% BO
% Usually we write $=$ and not $\mapsto$ for identifications, including again in this sentence. 
$x^*\odot y=\AB{x,y}$ and $E\odot E^*=\sK(E)$ via $x\odot y^*=xy^*$. 
%%%% EO
(Note that $x\odot y^*\odot z=x\AB{y,z}$ no matter which tensor product is calculated first, and that $x^*\odot y\odot z^*=(z\odot y^*\odot x)^*$.) Clearly, if we have \nbd{\sK(E)}correspondences $F_i$, then $E_i:=E^*\odot F_i\odot E$ are \nbd{\cB}correspondences. And, clearly, this operation respects tensor products in the sense that if $F_1\odot F_2\cong F_3$ via an isomorphism $u$ then $E_1\odot E_2\cong E_3$ via the induced isomorphism $\id_{E^*}\odot u\odot\id_E$. So, a (super)(sub)product system of \nbd{\sK(E)}correspondences gives rise to a (super)(sub)product system of \nbd{\cB}correspondences.

Now suppose we have \nbd{\sB^a(E)}correspondences $F_i$. In order that
\beq{ \label{Metpc}
E_1\odot E_2
~=~
(E^*\odot F_1\odot E)\odot(E^*\odot F_2\odot E)
~=~
E^*\odot F_1\odot\sK(E)\odot F_2\odot E
}\eeq
gives $E_3=E^*\odot F_3\odot E$, it is necessary to ``make disappear'' the $\sK(E)$ in between $F_1$ and $F_2$. Our \nbd{\sB^a(E)}correspondences will frequently be full (even with unit vectors) over $\sB^a(E)$, so that in this case there is no way that $F_1\odot\,\sK(E)=\cls F_1\sK(E)$ could coincide with $F_1$. Likewise, we cannot guarantee, in general, that $\sK(E)\odot F_2=\cls\sK(E)F_2$ is equal to $F_2$. (For instance, the identity correspondence $\sB^a(E)$ does not fulfill that.) However, we frequently can guarantee, for instance, when speaking about strict CP-semigroups on $\sB^a(E)$, that $\sK(E)\odot F_2\odot E=\cls\sK(E)(F_2\odot E)=F_2\odot E$. We, then, say the correspondence $F_2\odot E$ from $\sB^a(E)$ to $\cB$ is \phantomsection\label{SCdef}\hl{strict}. We borrow the following result from Skeide and Sumesh \cite[Lemma 3.1]{SkSU14}, and repeat also the simple proof. (This lemma proves \cite[Theorem 3.2]{SkSU14}, a unified version of the \it{KSGNS-construction} (see Lance \cite[Section 5]{Lan95}) and the (strict) representation theory of $\sB^a(E)$ (see \cite{MSS06}). The latter case we illustrate in Example \ref{EPSex}, while the former reoccurs more or less in Observation \ref{opGNSob}.)

\bprop\label{striGNSprop}
Let $E$ be a Hilbert \nbd{\cB}module and let $F$ be a Hilbert \nbd{\cC}module (where $\cB$ and $\cC$ are not necessarily unital \nbd{C^*}algebras). Let $T\colon\sB^a(E)\rightarrow\sB^a(F)$ be a CP-map and denote by $\sF$ its GNS-correspondence. Then $T$ is strict if and only if the correspondence $\sF\odot F$ from $\sB^a(E)$ to $\cC$ is strict.
\eprop

\proof
Suppose $\sF\odot F$ is strict. Denote by $\zeta\in\sF$ the cyclic vector, and define $\zeta\odot\id_F\in\sB^a(F,\sF\odot F)$ as $y\mapsto\zeta\odot y$. Then, clearly, $T=\AB{\zeta,\bullet\zeta}=(\zeta\odot\id_F)^*\bullet(\zeta\odot\id_F)$ is strict.

Conversely, suppose $T$ is strict. This means, in particular, that for every bounded approximate unit $\bfam{u_\lambda}_{\lambda\in\Lambda}$ for $\sK(E)$ and for every $a\in\sB^a(E)$, the bounded net $T(a^*(\id_E-u_\lambda)^*(\id_E-u_\lambda)a)$ converges strictly in $\sB^a(F)$ to $0$. For $a\in\sB^a(E)$, $a'\in\sB^a(F))$, and $y\in F$, it follows that the net $u_\lambda(a\zeta a'\odot y)=u_\lambda a\zeta\odot a'y$ converges in norm to $a\zeta\odot a'y=a\zeta a'\odot y$.\qed

\bcor\label{stri1cor}
If $T$ is strict, then
\beqn{
\sF\odot F
~=~
\sK(E)\odot\,\sF\odot F
~=~
E\odot E^*\odot\,\sF\odot F
~=~
E\odot\sE
}\eeqn
as \nbd{\sB^a(E)}\nbd{\cC}correspondences, where we define the \nbd{\cB}\nbd{\cC}correspondence $\sE:=E^*\odot\,\sF\odot F$.
\ecor

\bcor\label{stri2cor}
Let $T$ be a strict CP-semigroup over $\bS^{op}$ on $\sB^a(E)$ and denote by $\bfam{\sF_t}_{t\in\bS}$ its GNS-subproduct system. Then the \nbd{\cB}correspondences $\sE_t:=E^*\odot\,\sF_t\odot E$ $(t\ne0)$ and $\sE_0:=\cB$ form a subproduct system, the \phantomsection\hl{subproduct system of \nbd{\cB}correspondences}\index{CP-semigroup!strict on $\sB^a(E)$!subproduct system of \nbd{\cB}cor\-re\-spond\-ences associated with}\index{subproduct system!of \nbd{\cB}cor\-re\-spond\-ences associated with a strict CP-semigroup} associated with $T$.
\ecor

Usually, we shall require that $E$ is full, so that $E^*\odot\,\sF_0\odot E=E^*\odot E=\cB=\sE_0$, automatically. Since all $\sE_t$ may be viewed also as correspondences over $\cB_E$, we may simply reduce to full $E$ by making $\cB$ smaller.

By \eqref{Metpc}, we always get $\sE_s\odot\sE_t\subset E^*\odot\sF_s\odot\sF_t\odot E$. We emphasize that the crucial contribution of Proposition \ref{striGNSprop} to Corollary \ref{stri2cor}, is that strictness of the $T_t$ leads to equality $\sE_s\odot\sE_t=E^*\odot\sF_s\odot\sF_t\odot E$, and only then by $\sF_s\odot\sF_t\supset\sF_{st}$, to the subproduct system property of $\bfam{\sE_t}_{t\in\bS}$. Note, however, that strictness is not necessary. 

\bex \label{nonstriPSex}
Let $E=H$ be an infinite-dimensional Hilbert space, and let $T$ be the CP-map on $\sB(H)$ defined by $T(a):=\id_H\vp(a)k^*k$ where $\vp$ is a state that annihilates $\sK(H)$ and where $k\in\sK(H)$. Note that $T$ is not strict. If $(G,\eta)=$ GNS-$\vp$, then GNS-$T\!=$ $(\sF:=G\otimes\sK(H),\zeta:=\eta\otimes k)$ with $\sB(H)$ acting from the left on the factor $G$. We find $\sF\odot\sF=(G\otimes\sK(H))\odot(G\otimes\sK(H))=G\otimes(\sK(H)\odot G)\otimes\sK(H)=\zero$, because $\sK(H)\odot G=\zero$. Also, $\sE:=E^*\odot\sF\odot E=(H^*\odot G)\otimes(\sK(H)\odot H)=\zero$, because $H^*\odot G=\zero$. So, despite the semigroup $T^n$ is not strict, we obtain $\sE_m\odot\sE_n=E^*\odot\sF_m\odot\sF_n\odot E$ for all $m,n\in\N_0$. Note that both the $\sF_n$ and the $\sE_n$ form even a product system.
\eex

On the other hand, if  $\bfam{\sF_t}_{t\in\bS}$ happens to be a product system, then by \eqref{Metpc} we always get that $\bfam{\sE_t}_{t\in\bS}$ is a superproduct system; by the discussion preceding Example \ref{nonstriPSex}, $\bfam{\sE_t}_{t\in\bS}$ is a product system 
%%%% BO 
% After struggling with this, I reached the conclusion that all we have to do to make this correct is to replace: "if and only if $T$ is strict." with: 
if $T$ is strict.
% but I am not sure this is what you wanted. 
%%%% EO
The following example captures the most prominent situation when $\bfam{\sF_t}_{t\in\bS}$ happens to be a product system.	

\bex\label{EPSex}
If $\vt$ is an \nbd{E}semigroup on $\sB^a(E)$, then $\sF_t={_{\vt_t}}\vt_t(\id_E)\sB^a(E)$ and the $\sF_t$ form a product system (with $a_s\odot a'_t\mapsto \vt_t(a_s)a'_t$). Consequently, the $E_t:=E^*\odot{_{\vt_t}}\vt_t(\id_E)\sB^a(E)\odot E$ $=E^*\odot{_{\vt_t}}\vt_t(\id_E)E$ form a superproduct system of \nbd{\cB}correspondences, which is a product system if $\vt$ is strict. (The product induced by \eqref{Metpc} simplifies to a form we reproduce in \eqref{absPS}, where we need it.) This is the product system associated with a strict \nbd{E}semigroup as in Skeide \cite{Ske09} (preprint 2004).

Actually, in \cite{Ske09} we discuss only \nbd{E_0}semigroups, and we use the notation $E_t=E^*\odot_t E$. We shall use the same notation for \nbd{E}semigroups. 
%%%% BO new
% Nothing inserted, just this comment: If I recall correctly, I somehow got the idea that what is written is "we shall *not* use the same notation here", and then I searched to see if we really do not use it. So it was just a misunderstanding, sorry. 
%%%% EO 
Note, however, that $_{\vt_t}E$ is not a correspondence unless $\vt_t$ is unital. 
\eex

Note that the situation of endomorphisms in Example \ref{EPSex} is not necessary for that $\bfam{\sF_t}_{t\in\bS}$ is a product system. In fact, every non-Markov CP-semigroup on $\C$ is an example.

\bob\label{opGNSob}
Note that not only $\zeta$, but every element $Y\in\sF$ gives rise to an operator
%%%% Bo
% $L_Y:=Y\odot\id_E\in\sB^a(F,\sF\odot F)=\sB^a(F,E\odot\sE)$;
$L_Y:=Y\odot\id_F\in\sB^a(F,\sF\odot F)=\sB^a(F,E\odot\sE)$; 
%%%% EO
see Skeide \cite[Section 7]{Ske09} for details. Clearly, $L_Y^*L_{Y'}=\AB{Y,Y'}$ and $L_{Ya}=L_Ya$, so we may identify $\sF$ as a concrete operator \nbd{\sB^a(F)}submodule of $\sB^a(F,E\odot\sE)$. This submodule is $\cls\sB^a(E)L_\zeta\sB^a(F)$; it is \it{strictly dense} 
%%%% BO
% We haven't defined strictly on corners. --Orr. 
(equivalently, \nbd{*}strongly dense)
%%%% EO
in $\sB^a(F,E\odot\sE)$. 

Conversely, suppose we have a correspondence $\sE$ from $\cB$ to $\cC$ and a map $L\in\sB^a(F,E\odot\sE)$, and suppose we define the (strict) CP-map $T(a)=L^*aL$. Then the GNS-correspondence of $T$ is $\sF:=\cls\sB^a(E)L\sB^a(F)\subset\sB^a(F,E\odot\sE)$ with cyclic vector $L$. It is \nbd{*}strongly and, therefore, strictly dense in $\sB^a(F,E\odot\sE)$ if and only if $\cls\sB^a(E)LF=E\odot\sE$. If $E$ is full (or if the \phantomsection\hl{range ideal}\index{range ideal} $\cB_E:=\cls\AB{E,E}$ acts nondegenerately on $\sE$, so that $\sE$ is a also a correspondence from $\cB_E$ to $\cC$), then we get
\beqn{
E^*\odot\,\sF\odot F
~=~
E^*\odot E\odot\sE
~=~
\sE.
}\eeqn
(Otherwise, we get $\cls\cB_E\sE$, with which we may replace $\sE$.) Comparing this identification with Corollary \ref{stri1cor}, we see that we recover $\sE$ as $\cls\bCB{(x^*\odot\id_\sE)Ly\colon x\in E,y\in F}$.
\eob

GNS-subproduct systems need not be adjointable. We discuss several examples in Section \ref{EXnonadSEC}. However, if a subproduct system is adjointable, then, as we now show, it comes from a strict CP-semigroup on a suitable $\sB^a(E)$. In these notes, we limit ourselves to CP-semigroups and subproduct systems without continuity conditions, even if the monoid is topological. Consequently, we may not expect that our constructions give ``small'' objects, for instance, preserving separability.

\bthm\label{adSPS-CPthm}
Let $\sE^\bodot=\bfam{\sE_t}_{t\in\bS}$ be an adjointable subproduct system of \nbd{\cB}correspondences over a cancellative monoid $\bS$. Then there are a (full) Hilbert \nbd{\cB}module $E$ and a strict CP-semigroup\index{CP-semigroup!strict on $\sB^a(E)$!constructed from subproduct system of \nbd{\cB}cor\-re\-spond\-ences} $T$ over $\bS^{op}$ on $\sB^a(E)$ such that $T$ has $\sE^\bodot$ as associated subproduct system of \nbd{\cB}corre\-spondences.
\ethm

\proof
Since the $w_{s,t}$ are adjointable, we may define the \it{associative product} $x_sy_t:=w_{s,t}^*(x_s\odot y_t)\in\sE_{st}$. Define $E:=\bigoplus_{s\in\bS}\sE_s$. Then $v_t\colon E\odot\sE_t\rightarrow E$ defined by
\beqn{
v_t
\colon
(x_s\odot y_t)
~\longmapsto~
w_{s,t}^*(x_s\odot y_t)
}\eeqn
is an adjointable(!) partial isometry in $\sB^a(E)$ with range $\bigoplus_{s\in\bS}\sE_{st}$. (Here is, where we need that $\bS$ is cancellative.) Moreover, the \it{action} $xy_t:=v_t(x\odot y_t)$ iterates associatively with the \it{product}, respectively, the $v_t^*$ iterate coassociatively with the \it{coproduct}, that is,
\beqn{
(\id_E\odot w_{s,t})v_{st}^*
~=~
(v_s^*\odot\id_t)v_t^*.
}\eeqn
So, the strict CP-maps $T_t$ defined by
\beqn{
T_t(a)
~:=~
v_t(a\odot\id_t)v_t^*
}\eeqn
form a semigroup $T$ over $\bS^{op}$ on $\sB^a(E)$. Indeed,
\bmun{
T_t(T_s(a))
~=~
v_t\bfam{(v_s(a\odot\id_s)v_s^*)\odot\id_t}v_t^*
~=~
v_t(v_s\odot\id_s)(a\odot\id_s\odot\id_t)(v_s^*\odot\id_s)v_t^*
\\
~=~
v_{st}(\id_E\odot w_{s,t}^*)(a\odot\id_s\odot\id_t)(\id_E\odot w_{s,t})v_{st}^*
~=~
v_{st}(a\odot\id_{st})v_{st}^*
~=~
T_{st}(a).
}\emun
(Note that in the step $(\id_E\odot w_{s,t}^*)(a\odot\id_s\odot\id_t)(\id_E\odot w_{s,t})=(a\odot\id_{st})$ we used that $w_{s,t}$ is an isometry.)

The subproduct system of \nbd{\cB}correspondences $E^\bodot$ associated with $T$ is given by the correspondences $E_t:=E^*\odot\,\sF_t\odot E$ where $\sF_t$ is the GNS-correspondence of $T_t$. From Observation \ref{opGNSob}, we know that $\sF_t=\cls\sB^a(E)v_t^*\sB^a(E)$ and that $\sE_t=\cls(E^*\odot\id_t)v_t^*E\cong E_t$ (via the canonical isomorphism $(x^*\odot\id_t)v_t^*x'\mapsto x^*\odot\zeta_t\odot x'$). The subproduct system coproduct of $\sF^\bodot=\bfam{\sF_t}_{t\in\bS}$ maps $v_{st}^*\in\sF_{st}$ to $v_s^*\odot v_t^*\in\sF_s\odot\,\sF_t$ (and this determines the coproduct). The coproduct of the subproduct system $E^\bodot$ is defined by
\beqn{
E_{st}
~=~
E^*\odot\sF_{st}\odot E
~\longrightarrow~
E^*\odot\sF_s\odot\sF_t\odot E
~=~
E^*\odot\sF_s\odot E\odot E^*\odot\sF_t\odot E
~=~
E_s\odot E_t.
}\eeqn
Since we stick to $E^*\odot\,\sF_{st}\odot E=\cls E^*\odot\zeta_{st}\odot E$ and since $\zeta_{st}$ is mapped to $\zeta_s\odot\zeta_t$, this means we have to write the latter as an element of $\sF_s\odot\,\sF_t=\sF_s\odot E\odot E^*\odot\,\sF_t$. So, let $\bfam{\sum_{i=1}^{n_\lambda}{z'_i}^\lambda{z^\lambda_i}^*}_{\lambda\in\Lambda}$ be a bounded approximate unit for $\sK(E)$ consisting of finite-rank operators. It follows that the coproduct map is
\beqn{
E_{st}
~\ni~
x^*\odot\zeta_{st}\odot y
~\longmapsto~
\lim_\lambda\sum_{i=1}^{n_\lambda}(x^*\odot\zeta_s\odot{z'_i}^\lambda)\odot({z^\lambda_i}^*\odot\zeta_t\odot y)
~\in~
E_s\odot E_t.
}\eeqn
We, therefore, have to check that
\beqn{
\sE_{st}
~\ni~
(x^*\odot\id_{st})v_{st}^*y
~\longmapsto~
\lim_\lambda\sum_{i=1}^{n_\lambda}((x^*\odot\id_s)v_s^*{z'_i}^\lambda)\odot(({z^\lambda_i}^*\odot\id_t)v_t^*y)
~\in~
\sE_s\odot\sE_t.
}\eeqn
gives back $w_{s,t}$. Now,
\beqn{
((x^*\odot\id_s)v_s^*z')\odot((z^*\odot\id_t)v_t^*y)
~=~
(x^*\odot\id_s\odot\id_t)(v_s^*\odot\id_t)(z'z^*\odot\id_t)v_t^*y.
}\eeqn
Inserting for $z'z^*$ the approximate unit and taking into account that the action of $\sB^a(E)$ at that place is strict, we get
\bmun{
\lim_\lambda\sum_{i=1}^{n_\lambda}((x^*\odot\id_s)v_s^*{z'_i}^\lambda)\odot(({z^\lambda_i}^*\odot\id_t)v_t^*y)
~=~
(x^*\odot\id_s\odot\id_t)(v_s^*\odot\id_t)v_t^*y
\\~=~
(x^*\odot\id_s\odot\id_t)(\id_E\odot w_{s,t})v_{st}^*y
~=~
w_{s,t}(x^*\odot\id_{st})v_{st}^*y.\qedsymbol
}\emun
\noqed\vspace{-2ex}

\brem
In showing that the associated subproduct system $E^\bodot$ coincides with the original one $\sE^\bodot$, we did every effort to produce a proof that does not depend on the specific nature of the $v_t$. The only important point is that $\sB^a(E)v_t^*E$ is total in $E\odot\sE_t$. This easily follows, by applying the specific $v_t^*$ to $x_t\in\sE_t\subset E$ which is sent to $\om\odot x_t$ where we defined the \hl{vacuum} $\om:=\U\in\sE_0\subset E$. Letting act $y\om^*$ gives any $y\odot x_t$. (Anyway, the concrete $v_t$ do not make the proof any simpler.) If, in the case of topological monoids with right invariant measures, we take direct integrals instead of direct sums, the solution of this problem, though possible, is more cumbersome.
\erem

\brem
The $v_t$ are isometries if and only if $\sE^\bodot$ is a product system and $T$ is an \nbd{E}semi\-group. (This possibility to construct an \nbd{E}semigroup for a product system is well-known \it{folklore}.) Since the $v_t$ are definitely not unitaries, $T$ is not an \nbd{E_0}semigroup. For well over 15 years, even for Arveson systems there were only hard proofs of the fact that every Arveson system comes from an \nbd{E_0}semigroup; see Arveson \cite{Arv89,Arv90,Arv90a,Arv89a} and Liebscher \cite{Lie09} (preprint 2003). The elementary proof for Arveson systems in Skeide \cite{Ske06}, made it possible to resolve the same problem also for one-parameter product systems of correspondences in Skeide \cite{Ske07,Ske16}. The adaptation of this proof would work only for very special monoids; see Murugan and Sundar \cite{MurSu18}.

We also ask the question if every subproduct system is the subproduct system of a strict Markov semigroup on some $\sB^a(E)$. (Here, we would need $v_t$ that are coisometries. This is as difficult as the construction of \nbd{E_0}semigroups for product systems.)
\erem

We close this section, answering the natural question how the CP-semigroup constructed in Theorem \ref{adSPS-CPthm} behaves under isomorphism of subproduct systems.

\bob \label{Tembedob}
Suppose we have two adjointable subproduct systems $\sE^\bodot$ and $\sE'^\bodot$ over a cancellative monoid $\bS$, and construct $E^{(')}$, $v^{(')}_t$, and $T^{(')}$ as in the proof of Theorem \ref{adSPS-CPthm}. Let $i^\bodot=\bfam{i_t}_{t\in\bS}$ be an adjointable embedding of $\sE^\bodot$ into $\sE'^\bodot$. (That is, recall, the $i_t$ are adjointable isometries $\sE_t\rightarrow\sE'_t$ satisfying $i_t=\id_\cB$ and $(i_s\odot i_t)w_{s,t}=w'_{s,t}i_{st}$.) Then we may define an adjointable isometry $i\colon E\rightarrow E'$ acting fiberwise as $i_t$. Clearly, $iv_t=v'_t(i\odot i_t)$. It follows that
\beqn{
i_tT_t(a)i_t^*
~=~
(i_tv_t)(a\odot\id_t)(i_tv_t)^*
~=~
v'_t(iai^*\odot i_ti_t^*)v'^*_t.
}\eeqn
If the embedding is an isomorphism, then the latter is precisely $T'_t(iai^*)$, so that the CP-semi\-groups $T$ and $T'$ are conjugate.

(Note: If $i_t$ is not surjective, then there is no possibility to interpret this in terms of $T'_t$, because $E'$ is full (so that $a'\odot\id_t$ determines $a'$ whatever $a'$ is, while $a'\odot i_ti_t^*$ does not) and $v'^*_t$ is surjective, while $v'_t$ is injective. In fact, if $\sE'^\bodot$ is a product system, then $T'$ is an \nbd{E}semigroup, and if $i_tT_t(a)i_t^*$ would equal $T'_t(iai^*)$, then also $T$ would have to be an \nbd{E}semigroup, meaning that also $\sE^\bodot$ had to be a product system. What we can say, is that the CP-map $T_t$ is dominated by the CP-map  $i_t^*T'_t(i\bullet i^*)i_t$.)
\eob

\newpage

\section[\sc{Examples:} CP-semigroups with nonadjointable GNS-sub\-product systems]{Examples: CP-semigroups with nonadjointable GNS-sub\-product systems}\label{EXnonadSEC}

In Theorem \ref{adSPS-CPthm} we showed for cancellative monoids that every subproduct system comes in some sense from a CP-semigroup -- provided it is adjointable, leading to that also the GNS-subproduct system of that CP-semigroup is adjointable. The more important it is to answer the question, whether there are GNS-subproduct systems that are not adjointable. In this section, we find examples of such CP-semigroups.

In the existent literature on subproduct systems, starting with Shalit and Solel \cite{ShaSo09}, with each discrete one-parameter subproduct system, there are associated operator algebras generalizing the well-known \it{Cuntz-Pimsner(-Toeplitz) algebras} (Pimsner \cite{Pim97} and (non-selfadjoint) \it{tensor algebras} (Muhly and Solel \cite{MuSo98}), which, in our language, would correspond to discrete one-parameter product systems. While in the von Neumann case considered in \cite{ShaSo09}, adjointability is automatic, Viselter \cite{Vis10,Vis12} started to consider the \nbd{C^*}case; here, adjointability is necessary to work and included in the definition of subproduct system. There are several 
%%%% BO new
%forthcoming 
subsequent
%%%% EO
papers by Dor-On and Markiewicz \cite{D-OMa14,D-OMa17} adopting this setting. 
%%%% BO new 
% We report a non-adjointable example
In this section we provide several examples of CP-semigroups with nonadjointable GNS-subproduct systems. Example \ref{doronex} below is 
%%%% EO 
from an unpublished note \cite{D-O14} by Dor-On that motivated him to propose in \cite{D-O18} a different construction  without adjointability.

\bex \label{nonadex}
Choose $c\in\cB$. \index{CP-semigroup!GNS-subproduct system of!non-adjointable}Then the GNS-correspondence of the elementary CP-map $T:=c^*\bullet c$ is the ideal generated by $c$ with cyclic vector $c$. (This is true even if $\cB$ is not necessarily unital.) For $c,d\in\cB$ denote by $(E,c)$, $(F,d)$, and $(G,cd)$ the GNS-constructions for $T:=c^*\bullet c$, $S:=d^*\bullet d$, and $S\circ T=(cd)^*\bullet(cd)$, respectively. Note that $E\odot F=\cls\cB c\cB d\cB$. (See also Observation \ref{elemob}.)

Suppose we have $c$ and $d$ such that $G$ is an essential ideal in $\cB$, and such that there exists $b\in\cB$ with $cbd\notin G$ (so that, in particular, $G$ is a nontrivial ideal). It follows that $G$ is a nontrivial essential ideal also in $E\odot F$, and as such a noncomplemented Hilbert submodule of $E\odot F$. Consequently, the embedding of the GNS-correspondence of $S\circ T$ into the tensor product of those for $T$ and for $S$ is not adjointable. If we find an elementary CP-semigroup $T_t=c_t^*\bullet c_t$ for some semigroup $c_t$ in $\cB$ such that $c_t=c$ and $c_s=d$, then its GNS-subproduct system is not adjointable.

Here is an explicit choice: Let $\cC$ be a unital \nbd{C^*}algebra (for instance, $\sB(H)$ for an infinite-dimensional separable Hilbert space $H$) with a nontrivial essential ideal $I$ (for instance $\sK(H)$) with a self-adjoint element $k\in I$ that generates $I$ as an ideal. Put $\cB=M_2(\cC)$ and $c=\sMatrix{k&\U\\0&-k}$. Then $c^2=\sMatrix{k^2&0\\0&k^2}$ generates the essential ideal $M_2(I)$, and $b=\sMatrix{0&0\\\U&0}$ is in $\cB$ such that $cbc=\sMatrix{k&\U\\-k^2&-k}$ is not in $M_2(I)$. So the CP-semigroup $T_n=(c^n)^*\bullet(c^n)$ over $\N_0$ does not have an adjointable GNS-subproduct system.
\eex

\lf
Of course, exponentiating Example \ref{nonadex} as in the Example-Section \ref{EXexpSEC}, will produce a continuous time example. We prefer to add here another construction.

\bex \label{cnonadex}
Bhat and Skeide \cite[Lemma 2.4(1)]{BhSk15} show how to interpolate a discrete semigroup of of isometries $\breve{\s}^n$ in $\sB^r(\breve{E})$ to a strongly continuous semigroups of isometries $\s_t$ in $\sB^r(L^2\SB{0,1}\otimes\breve{E})$ in the sense that $\s_n=\id_{L^2}\otimes\,\breve{\s}^n$ for all $n\in\N_0$. 
%%%% BO 
The semigroup $\s_t$ is simply defined by the formula
\vspace{-2ex}
\beqn{
\s_t
~:=~
(u_t\otimes\id_{\breve{F}})
\bfam{\I_{\RO{0,1-(t-n_t)}}\otimes\breve{\s}^{n_t}+\I_{\RO{1-(t-n_t),1}}\otimes\breve{\s}^{n_t+1}},
\vspace{-1ex}
}\eeqn
where $n_t$ is the largest integer $\le t$, $\I_S$ is the indicator function of the set $S$, and $u_t$ is the unitary right shift modulo $1$ on $L^2\SB{0,1}$, that is, $u_t f(x) = f((x-t) \operatorname{mod} 1)$. It is clear that the above formula for interpolation also works if $\breve{\s}$ is not necessarily an isometry. Replacing $\breve{s}$ with the contraction $c$ from Example \ref{nonadex} we get a contraction semigroup $\s_t$ (adjointable, because $c$ is). 
Clearly, whatever is true for the GNS-correspondences of the CP-maps $(c^n)^*\bullet(c^n)$ and their tensor products, is also true for the GNS-correspondences of the CP-maps $(\s_n)^*\bullet(\s_n)=\id_{\sB(L^2)}\otimes(c^n)^*\bullet(c^n)$.
\eex

\lf
The preceding examples are essentially noncommutative, and it is natural to ask whether all GNS-subproduct systems over a commutative \nbd{C^*}algebra $\cB$ are adjointable. The answer is negative, but it is much harder to exhibit a commutative example with all details (let alone find one). Here is the promised example from Dor-On \cite{D-O14}.

\bex  \label{doronex}
Let $\cB = C[0,1]$ and consider the Markov map
\beqn{
\textstyle
T 
\colon
f
 ~\longmapsto~
 \BSB{
x ~\mapsto~\frac{1}{2}f(\frac{1}{2}) + \frac{1}{2}f(x^2)
}
}\eeqn
on $\cB$. It can be shown that the GNS-subproduct system of the Markov semigroup $\bfam{T^n}_{n\in\N_0}$ is non-adjointable; we omit the details. 
\eex

% \OW{REREAD: It might be of interest to describe the context in which Dor-On reached the above counter example. 
% Every Markov map gives rise to a Markov semigroup over $\N_0$, and thus to a GNS-subproduct system. 
% In turn, there is a way to construct certain operator algebras -- the \bf{tensor algebra} and the \bf{Cuntz-Pimsner algebra}  -- from an adjointable subproduct system \cite{Vis10,Vis12}. So, one obtains a way to attach certain operator algerbas to certain Markov maps, and this can serve as an invariant of the map, by which one can study it. This approach was taken by Dor-On and Markiewicz in \cite{D-OMa14,D-OMa17}, in the case of Markov maps on commutative von Neumann algebras. Dor-On was interested in generaling those works to the case of commutative C*-algebras. The fact that GNS-subproduct systems arising from Markov maps on commutative \nbd{C^*}algebras may be non-adjointable, caused that approach to break down, and led Dor-On to encode such maps using tensor algebras of \nbd{C^*}correspondences in \cite{D-O18}. }
% \erem
%%%% EO

\newpage

\section{Dilations and superproduct systems}\label{DilSPSpSEC}

In Example \ref{EPSex}, we have recovered how to associate with a strict \nbd{E}semigroup $\vt$ over the monoid $\bS^{op}$ on $\sB^a(E)$ (where $E$ is a full Hilbert \nbd{\cB}module) a product system $E^\odot=\bfam{E_t}_{t\in\bS}$ over $\bS$ of \nbd{\cB}correspondences $E_t:=E^*\odot_t E$. In this section, we examine to what extent we can do the same for the \nbd{E}semigroup $\theta$ of a dilation $(\cA,\theta,p)$. It turns out that we are led only to \it{superproduct systems} -- superproduct systems containing, for strong dilations, the GNS-subproduct system of the dilated semigroup -- but not always a product system. Still, this allows us to show existence of Markov (CP-)semigroups that do not admit any (strong) dilation; see the Example-Section \ref{EXsubnsupSEC}. The superproduct systems of weak dilations may contain or not contain the GNS-subproduct system. The latter exhibit such \it{bad} behaviour that we call the other dilations \it{good}. We give an illustration in the Example-Section \ref{EXBexSEC}; note that Bhat's Example\index{Bhat's example} \ref{Bex}, is a discrete one-parameter example.

The obstacles to obtain a product system are two-fold. Already for the $E_t$ constructed from an \nbd{E}semigroup on $\sB^a(E)$, for that the (always associative) \it{product}
\beq{ \label{absPS}
(x^*\odot_s x')\odot(y^*\odot_t y')
~\longmapsto~
(x^*\odot_s x')(y^*\odot_t y')
~:=~
x^*\odot_{st}\vt_t(x'y^*)y'
}\eeq
among the $E_t$ generates all $E_{st}$ on the right-hand side, we need that $\vt_t$ is strict; otherwise, we only get a superproduct system. This raises the question, what does `strict' mean for the endomorphisms $\theta_t$ of a strong dilation $(\cA,\theta,p)$. (It would be clear, if the dilation is full so that $\cA=\sB^a(\cA p)$; see Proposition \ref{convprop} and the discussion preceding it. We come back to this question 
%%%% BO new
%then
when
%%%% EO 
we discuss topological questions in Section \ref{topSEC} closely related to notions of minimality in Section \ref{minSEC}.) Fortunately, the question about the right topology in which each $\theta_t$ should be continuous, goes away in the von Neumann-case; they should be normal and normality does not depend on how we represent $\cA$. Unfortunately, in the following construction of a superproduct system from $(\cA,\theta,p)$, there are more serious, algebraic, obstacles why in general we only get a superproduct system (see Examples \ref{hypexex} and \ref{discex}). Also these obstacles will be discussed in more detail in Section \ref{topSEC}. For now, let us concentrate on getting the construction.

Recall from Skeide \cite{Ske02}, that if $E$ has a unit vector $\xi$, then the product system of the strict \nbd{E}semigroup $\vt$ on $\sB^a(E)$ may also be obtained as follows: Put $\sE_t=\vt_t(\xi\xi^*)E$ and define a left action of $\cB$ on $\sE_t$ by $bx_t:=\vt_t(\xi b\xi^*)x_t$. Then
\beqn{
x_s\odot y_t
~\longmapsto~
\vt_t(x_s\xi^*)y_t
}\eeqn
turns $\bfam{\sE_t}_{t\in\bS}$ into a product system $\sE^\odot$ (a superproduct system, if $\vt$ is not strict). This imitates Bhat's construction \cite{Bha96} for Hilbert spaces, where unit vectors exist in abundance; \cite{Ske02} just added the ``correct'' left action.  As mentioned in Skeide \cite[Section 1]{Ske04p} (version 4), the maps $x^*\odot_t y\mapsto\vt_t(\xi x^*)y$ establish an isomorphism $E^\odot\rightarrow\sE^\odot$ with inverse given by $\vt_t(\xi\xi^*)z\mapsto\xi^*\odot_t z$.

By Proposition \ref{convprop}, if we have a full dilation $(\cA,\theta,p)$, then it can be written as a module dilation $(E,\vt,\xi)$, where $E=\cA p$ and $\xi=p$ and, by fullness, $\sB^a(E)=\cA$, identifying also $\vt$ with $\theta$. Comparing the roles of $\xi$ and $p$, we are ready to imitate the construction from \cite{Ske02} directly in terms of $(\cA,\theta,p)$.

\bthm \label{E-supPSthm}
Let $\cA$ be a \nbd{C^*}algebra with unit\index{dilation, weak!superproduct system of}\index{systems!superproduct!of a triple (dilation)}\index{superproduct system!of a triple (dilation)} $\U$, let $\theta$ be an \nbd{E}semigroup over $\bS^{op}$ on $\cA$, and let $p$ be a projection in $\cA$. Put $\cB:=p\cA p$, and define the Hilbert \nbd{\cB}submodules $\sE_t:=\theta_t(p)E$ of the Hilbert \nbd{\cB}module $E:=\cA p$. Define the left action $b.x_t:=\theta_t(b)x_t$ of $\cB$ on $\sE_t$. Then the maps $v_{s,t}$ defined by
\beqn{
v_{s,t}
\colon
x_s\odot y_t
~\longmapsto~
\theta_t(x_s)y_t
}\eeqn
$(x_s\in\theta_s(p)\cA p$, $y_t\in\theta_t(p)\cA p)$ turn $\sE^\podot=\bfam{\sE_t}_{t\in\bS}$ into a superproduct system  over $\bS$.
\ethm

\proof
Note that the action of $\theta_t(a)$ defines a ($*$!)representation of the \nbd{C^*}algebra $\cA$ on $E=\cA p$, which, therefore, is a contraction into $\sB^a(E)$. Since the unit of $\cB$ is $p$, the left action of $\cB$ on $\sE_t$ is unital, so the $\sE_t$ are, indeed, correspondences over $\cB$. Also,
\beq{\label{spPeq}
\theta_t(\theta_s(p)ap)\theta_t(p)a'p
~=~
\theta_{st}(p)\theta_t(ap)a'p
~\in~
\sE_{st}.
}\eeq
(Recall that $\theta$ is a semigroup over $\bS^{op}$ so that $\theta_t\circ\theta_s=\theta_{st}$.) For $x_s=\theta_s(p)a_1p$, $y_t=\theta_t(p)a_2p$, $x'_s=\theta_s(p)a_3p$, $y'_t=\theta_t(p)a_4p$, we see that
\bmun{
\AB{x_s\odot y_t,x'_s\odot y'_t}
~=~
\BAB{(\theta_s(p)a_1p)\odot(\theta_t(p)a_2p),(\theta_s(p)a_3p)\odot(\theta_t(p)a_4p)}
\\
~=~
\BAB{\theta_t(p)a_2p,\bAB{\theta_s(p)a_1p,\theta_s(p)a_3p}.\theta_t(p)a_4p}
~=~
\BAB{\theta_t(p)a_2p,\theta_t((\theta_s(p)a_1p)^*\theta_s(p)a_3p)\theta_t(p)a_4p}
\\
~=~
\BAB{\theta_t((\theta_s(p)a_1p)\theta_t(p)a_2p,\theta_t((\theta_s(p)a_3p)\theta_t(p)a_4p}
~=~
\bAB{\theta_t(x_s)y_t,\theta_t(x'_s)y'_t},
}\emun
so the $v_{s,t}$ are isometries. Of course, $\sE_0=\theta_0(p)\cA p=\cB$ and $v_{0,t},v_{t,0}$ fulfill the marginal conditions. Finally,
\baln{
(x_r\odot y_s)\odot z_t
&
~\longmapsto~
\theta_t(\theta_s(x_r)y_s)z_t
~=~
\theta_{st}(x_r)\theta_t(y_s)z_t,
\\
x_r\odot(y_s\odot z_t)
&
~\longmapsto~
\theta_{st}(x_r)\theta_t(y_s)z_t.
}\ealn
Thus, $\sE^\podot=\bfam{\sE_t}_{t\in\bS}$ is a superproduct system over $\bS$.\qed

\bob \label{fu-mod-ob}
%%%% BO 
Note that $\sE_t$ is, clearly, isomorphic to 
%%%% BO new
% $E^*\odot_t E$, 
$E^*\odot_t E:=E^*\odot_{\theta_t}\theta_t(\U_\cA)E$, 
% I added something to make sure it is clear how this should be defined (if I got it right this time...)
%%%% EO (new)
where the tensor product over $\sK(E)\subset\cA$ is for the left action on $E$ is via $\theta_t$ like the one of $\vt_t$ in Example \ref{EPSex}. Also the superproduct system structure is the right one, given by \eqref{absPS}.
%%%% EO
\eob

In Theorem \ref{E-supPSthm} we did \bf{not} require that $(\cA,\theta,p)$ is a dilation.

\bthm \label{sdilunithm}
Let $(\cA,\theta,p)$ and $\sE^\podot$ be as in Theorem \ref{E-supPSthm}, and put $\xi_t:=\theta_t(p)p$.

\begin{enumerate}
\item \label{DU1}
If $(\cA,\theta,p)$ is a strong dilation, then the $\xi_t$ form a unit $\xi^\odot$ for $\sE^\podot$.

\item \label{DU2}
If the $\xi_t$ form a unit $\xi^\odot$ for $\sE^\podot$, then $(\cA,\theta,p)$ is a weak dilation.

\item \label{DU3}
In either case, the superproduct system $\sE^\podot$ of $\theta$ contains the GNS-sub\-product system $\Bfam{\bfam{\cls\cB\xi_t\cB}_{t\in\bS}\,,\,\xi^\odot}$ as a subproduct subsystem.
\end{enumerate}
\ethm

\proof
\ref{DU1}. 
If we have a strong dilation, then $\theta_t(p)p=\theta_t(\U)p$, so
\beqn{
\xi_s\xi_t
~=~
v_{s,t}(\theta_s(p)p\odot\theta_t(p)p)
~=~
\theta_t(\theta_s(p))\theta_t(p)p
~=~
\theta_t(\theta_s(p))\theta_t(\U)p
~=~
\theta_t(\theta_s(p))p
~=~
\theta_{st}(p)p.
}\eeqn

\ref{DU2}. 
Whatever $(\cA,\theta,p)$ is, we always have $\AB{\xi_t,(pap)\xi_t}=p\theta_t(pap)p$. So, $(\cA,\theta,p)$ is a weak dilation if and only if the $T_t:=\AB{\xi_t,\bullet\xi_t}$ form a semigroup. The latter is the case, if the $\xi_t$ form a unit.

\ref{DU3}. 
This follows from Example \ref{unitgenex}.\qed

\brem
Note that the proof of Theorem \ref{sdilunithm}\eqref{DU1} is slightly more elaborate than the proof of the analogue statement in \cite[Proposition 3.1]{Ske00}, where $p$ is increasing.
\erem

\bob \label{puniob}
From the formula in the proof, we also find $\xi_s\xi_t=\theta_{st}(p)\theta_t(p)p$ in the generality of Theorem \ref{E-supPSthm}. So, the question if the $\xi_t$ form a unit amounts to the question if $\theta_{st}(p)\theta_t(p)p=\theta_{st}(p)p$. Written in the form
\beq{ \label{gcond}
\theta_{st}(p)\theta_t(\U-p)p
~=~
0
}\eeq
(taking also into account that $\theta_{st}(p)=\theta_{st}(p)\theta_t(\U)$), we see explicitly that this condition is weaker than the condition
\beq{ \label{scond}
\theta_t(\U-p)p
~=~
0
}\eeq
to be strong in Theorem \ref{strequivprop}\eqref{sp1}. It cannot be overemphasized that both conditions, \eqref{scond} and \eqref{gcond}, also imply the statement that the triple $(\cA,\theta,p)$ actually is a dilation. Frequently, they are the only applicable criteria that allow to check that the $T_t$ form a semigroup. See, in particular, the discussion of compressions in Subsections \ref{minSEC}\ref{primSSEC} and \ref{minSEC}\ref{compSSEC}.
\eob

\bex \label{pPSelemex}
Suppose $\theta$ is an elementary \nbd{E}semigroup (over $\bS^{op}$) on $\cA$, so that there is a semigroup $w$ (over $\bS$) of coisometries $w_t\in\cA$ such that $\theta_t=w_t^*\bullet w_t$. For a projection $p\in\cA$, we find $\sE_t:=w_t^*pw_t\cA p=w_t^*p\cA p=w_t^*\cB$. We easily check that $b\mapsto w_t^*b$ establishes an isomorphism of superproduct systems from the trivial product system $\cB^\odot$ onto the superproduct system $\sE^\podot$, which, therefore, is a product system. (Of course, the image $\bfam{w_t^*p}_{t\in\bS}$ of the unit $\U_\cB^\odot$ for $\cB^\odot$ has no chance but to be a unit for $\sE^\podot$. But this is not the unit that interests us in this context.)

Now suppose that $(\cA,\theta,p)$ is a solidly elementary dilation, so that the $c_t:=pw_tp$ form a contraction semigroup. Then $\xi_t:=w_t^*pw_tp=w_t^*c_t$ is the image of $c_t$ under the isomorphism, and since the $c_t$ form a unit for $\cB^\odot$, so do the elements $\xi_t$ in $\sE_t$.
\eex

\lf
Theorem \ref{sdilunithm} is about the question whether for the triple $(\cA,\theta,p)$ the $\xi_t:=\theta_t(p)p$ form a unit. By the importance of the consequences, we call such triples \phantomsection\hl{good}\index{good triples (dilations)}. Part \ref{sdilunithm}\eqref{DU2} tells us that good triples are weak dilations, which we, therefore, call \hl{good dilations}\index{dilation, weak!good}. Example \ref{dwnsex}, taking into account also Example \ref{pPSelemex}, furnishes good solidly elementary dilations (with product system) in the discrete one-parameter case, which may be strong or not. By the procedure in Example \ref{cwnsex}, this lifts to the continuous time one-parameter case. So, being good is not sufficient for being strong.

Bhat's Example\index{Bhat's example} \ref{EXBexSEC}\ref{BexSSEC} gives a weak dilation of an elementary (even scalar!) CP-semigroup which is not good. (Here, we do not know if we can get a continuous time example.) So, the condition in \ref{sdilunithm}\eqref{DU2} is not necessary. Tensor products and direct sums of dilations are, clearly, good if and only if both constituents are good. Therefore, while all dilations of Markov semigroups are strong, hence, good, discrete one-parameter non-Markov CP-semigroup that decay \it{sufficiently fast} admit dilations that are not good. (Namely, if the CP-semigroup can be written as product, hence, tensor product of another CP-semigroup with the scalar CP-semigroup from Example \ref{EXBexSEC}\ref{BexSSEC}.)

Theorem \ref{sdilunithm}\eqref{DU1} does have a partial converse in Theorem \ref{minstriunithm}. Namely, a \it{minimally strict dilation} is good if and \bf{only} if the dilation is strong. Strictness of an algebraically minimal dilation is a much clearer (and natural) condition, which in the von Neumann case is, anyway, replaced by the omnipresent condition to be normal. Since algebraic minimality can always be achieved, and since the appropriate minimalization procedure preserves goodness, we get, essentially, that existence of a good dilation implies existence of a strong dilation. (``Essentially'', refers to that we get the statement in the von Neumann case.) The topological issues are explained in Section \ref{topSEC}, while questions of minimality are dealt with in Section \ref{minSEC}.

We have now reached the basic connection between CP-semigroups and subproduct systems in Section \ref{CPspsSEC} and between dilations and superproduct systems in the present section. The discussion in Theorem \ref{sdilunithm} shows the importance of the situation where the subproduct system sits in the superproduct system, that is, of good dilations. And in Section \ref{leftdilSEC} (in particular, Theorem \ref{Markmodthm}) we will see the importance of the question whether the superproduct system of a good dilation embeds into a product system. This appears to be a good moment to collect the questions that are suggested by these and other considerations. We put them down in the next section. There, we also recall the answers we can already give, anticipate the answers we will give later on in these notes, and single out the questions to which we will not be able to give an answer in what follows.

\newpage

\section[\sc{Questions:} With and without answers]{Questions: With and without answers} \label{QSEC}

Let us rest for a moment to see what we have achieved in the preceding sections, and to which natural questions this leads.
\begin{itemize}
\item
With every weak dilation -- actually, with every triple -- $(\cA,\theta,p)$ we have associated a superproduct system. If the dilation (or the triple) is full, that is if $\cA=\sB^a(\cA p)$ in the sense of Proposition \ref{convprop}, and if it is \hl{strict}\phantomsection\index{dilation, weak!strict}\index{strict!dilation} (or \hl{normal}\index{dilation, weak!normal}\index{normal!dilation} on the von Neumann case) in the sense that each $\theta_t$ is strict (or normal in the von Neumann case), then, by Example \ref{EPSex}, the superproduct system is even a product system.

\item
If $(\cA,\theta,p)$ is a strong dilation, then the superproduct system has a unit giving back the dilated CP-semigroup and, therefore, contains the GNS-subproduct system of the dilated semigroup. If the dilation is full and strict (or normal in the von Neumann case), then the GNS-subproduct system is even contained in a product system.

\end{itemize}
{Conclusion:}
\begin{itemize}
\item
A subproduct system that does not embed into a superproduct system cannot come in any way from a CP-semigroup that admits a strong dilation.

\item
A subproduct system that does not embed into a product system cannot come in any way from a CP-semigroup that admits a strict (or normal) strong module (or full) dilation.
\end{itemize}
%%%% BO new
%By ``in any way'', we mean the three known possibilities to construct subproduct systems from CP-semigroups. We mentioned that we may not expect to obtain every subproduct system as GNS-subproduct system of a CP-semigroup: The GNS-subproduct system; the associated subproduct system of \nbd{\cB}correspondences; the Arveson-Stinespring subproduct system. 
%
% PLEASE CHECK (flipped around):
%
By ``in any way'', we mean the three known possibilities to construct subproduct systems from CP-semigroups: The GNS-subproduct system; the subproduct system of \nbd{\cB}correspondences associated with a CP-semigroup as in Corollary \ref{stri2cor}; and the Arveson-Stinespring subproduct system. We mentioned that we may not expect to obtain every subproduct system as GNS-subproduct system of a CP-semigroup. 
%%%% EO
However, we have learned in Theorem \ref{adSPS-CPthm} that at least every adjointable subproduct system over a cancellative monoid is the subproduct system of \nbd{\cB}correspondences of a strict CP-semigroup on some $\sB^a(E)$ (relating to the GNS-subproduct systems via (strict) Morita equivalence), and we have learned from \cite[Corollary 2.10]{ShaSo09} that every subproduct system of von Neumann \nbd{\cB'}correspondences (automatically adjointable) over submonoids of $\R^d$ is the Arveson-Stinespring subproduct system of a normal CP-semigroup on some $\sB^a(E)$. Note that in the von Neumann case, \cite[Corollary 2.10]{ShaSo09} (where applicable) and Theorem \ref{adSPS-CPthm} are \it{equivalent via commutant}. This is so, essentially, because Theorem \ref{adSPS-CPthm} may be rephrased saying that every subproduct system of von Neumann \nbd{\cB}correspondences is Morita equivalent to a GNS-subproduct system (namely, of \nbd{\sB^a(E)}correspondences) while the commutant $\cB=\cB''$ of $\cB'$ in \cite[Corollary 2.10]{ShaSo09} exactly depends up to Morita equivalence on the representation of $\cB'$; see also Appendix \ref{vNAPP}\ref{vNcomm}.

This raises several questions. (In particular, the answer to Question \ref{Q2}, improving on results from \cite{ShaSo09, ShaSk11}, illustrates how the notion of the superproduct system of a dilation allows to show that there exist CP-semigroups with no strong dilation.)

% \newpage
\bemp[Question:~] \label{Q1}
For that the superproduct system of a dilation (triple) $(\cA,\theta,p)$ be a product system, is it enough that the dilation be normal or `nice' in some other reasonable topology?

\bf{Answer:} No. See Example \ref{hypexex}. (Note that this example, a strong dilation, is not a full or module dilation. If it was full, then we know we get a product system.)
\eemp

\bemp[Question:~] \label{Q2}
Does there exist an (adjointable) subproduct system that does not embed into a superproduct system?

\bf{Answer:} Yes. In Section \ref{EXsubnsupSEC} we analyze the known example of a(n adjointable) subproduct system that does not embed into a product system from Shalit and Solel \cite[Proposition 5.15]{ShaSo09}, and we show that it actually does not even embed into a superproduct system. Therefore: There exists a CP-semigroup that does not admit whatsoever strong dilations, reinforcing the result from \cite{ShaSo09} that this CP-semigroup does not admit a full dilation and no \it{minimal} dilation to a von Neumann algebra.

By unitalization, the preceding example extends to Markov semigroups: There exists a Markov semigroup that admits no weak$=$strong dilation.
That is, we get an example of a CP-semigroup (namely, the stated unitalizations) that do not admit any dilation, no matter whether weak or strong. However:
\eemp

\bemp[Question:~]\label{Q3}
We know, there do exist weak dilations that are not strong; see Section \ref{EXwnsSEC}. But, does there exist a CP-semigroup with no strong but a weak dilation?

\bf{Answer:} Unknown.

However, if the answer is no, that is, if the answer is that every CP-semigroup that admits a dilation also admits a strong dilation, then we would obtain that embeddability of the GNS-sub\-product system into a superproduct system is necessary for existence of any dilation.

Of course, the examples for 
%%%% BO
proper (=non-Markov) 
%%%% EO
CP-semigroups that admit no strong dilations from Question \ref{Q2} are candidates for being counter examples also here. We have to examine whether they admit any dilation, necessarily a weak one.
% \OW[OPEN]{(We leave this open. --Orr)}

Note that by the discussion following Theorem \ref{sdilunithm} and based on Theorem \ref{minstriunithm}, we have to search among the weak dilations for which the $\xi_t:=\theta_t(p)p$ do not form a unit, that is, which are not \it{good}, and which are \it{algebraically minimal}.
% \OW[OPEN: Which weak, not strong, are algebraically minimal, and and how does the minimal version look like? Anyway, algebraic minimality is a question that is relevant for ALL our ``strange'' dilations.]{(I think these two todos can be erased, they are left as explicit open questions in the text. --Orr)}

% \bthmn
% Suppose the CP-semigroup $T$ on the unital \nbd{C^*}algebra $\cB$ is is fundamentally good. Then:
% \begin{enumerate}
% \item
% All algebraically minimal strict (or normal) weak dilations of $T$ are strong.

% \item
% Therefore, existence of a weak dilation implies existence of a strong dilation.

% \item
% All superproduct systems of weak dilations contain the GNS-subproduct system.

% \item
% $T$ admits a weak dilation if and only if its unitalization admits a weak, hence, strong dilation. (See also Questions \ref{Q4} and \ref{Q5}.)
% \end{enumerate}
% \ethmn
% Strictly, speaking Points 2-4 are true only if we speak about \it{weak dilations whose minimalization is strict}. It is true in the von Neumann case for normal dilations.

% We remain with the (unanswered) question if all CP-semigroups are fundamentally good.
\eemp

% \newpage
\bemp[Question:~]\label{Q4}
For existence of a strong dilation, unitalization is important, because by Theorem \ref{uninonunithm} the question is equivalent to the question if the unitalization has a dilation. (By the first Corollary \ref{unicor1}, even existence of an \nbd{E_0}dilation of the unitalization.) This equivalence is relevant also for several other questions. But here we are particularly interested in the observation that in the one-parameter case we have the well-known inductive limit construction from Bhat and Skeide \cite[Section 5]{BhSk00} that, starting from a product system with unital unit, provides a module dilation of the Markov semigroup determined by that unit. In \cite[Section 8]{BhSk00} and in Skeide \cite{Ske08}, several unitalization procedures have been applied to get the same result for arbitrary one-parameter CP-semigroups. We ask:

Can we promote the inductive limit construction from \cite[Section 5]{BhSk00} to more general monoids?

\bf{Answer:} Yes, for Ore monoids. (See Theorem \ref{Markmodthm}.) It is 
%%%% BO
one of 
%%%% EO
the major results of Section \ref{leftdilSEC} that by Theorem \ref{Oreindthm}, from a product system over an Ore monoid and a unital unit for that product system, we can construct a strict module dilation of the Markov semigroup determined by the unit.

Consequently, for Markov semigroups over the opposite of an Ore monoid, existence of a strict module (\nbd{E_0})dilation is equivalent to that the GNS-subproduct system embeds into a product system.

For existence of strict strong module dilations of CP-semigroups, the condition of embeddability is, of course, a necessary condition. However, even for CP-semigroups over the opposite of an Ore monoid, to deduce from embeddability existence of strong dilations, one has to look for embeddability of the GNS-subproduct system of the unitalization. It is not enough to have embeddability for GNS-subproduct system of the CP-semigroup itself:
\eemp

\bemp[Question:~] \label{Q5}
Suppose the GNS-subproduct system of a CP-semigroup embeds into a product system. Does this mean that the GNS-subproduct system of its unitalization embeds into a product system?

\bf{Answer:} No. In Observation \ref{elemob} we have illustrated that the GNS-subproduct system of an elementary CP-semigroup does not only embed into a product system, namely, into the trivial product system, but that it is actually isomorphic to the trivial product system if the elementary CP-semigroup acts on a simple \nbd{C^*}algebra. By \cite[Theorem 5.14]{ShaSo09}, which is based on Parrot's construction Example \ref{Parex}, we get an elementary CP-semigroup on $\sB(H)$ ($H$ even finite-dimensional so that $\sB(H)$ is a simple \nbd{C^*}algebra) with no strong dilation, despite the fact that the GNS-subproduct system does not only embed but, actually, is already a product system.

In Lemma \ref{pswtpslem}, we provide a precise condition for when the answer is positive, at least when there is a strong dilation.
\eemp

% \newpage
\bemp[Question:~] \label{Q6}
We know that embeddability of the GNS-subproduct system of a CP-semi\-group into a superproduct system is a necessary condition for existence of a strong dilation. For existence of a strict (or normal) module dilation, the embedding has to be even into a product system. On the other hand, at least for a large class of monoids, embeddability of the GNS-subproduct system of a Markov semigroup (for instance, of that of the unitalization of a CP-semigroup) into a product system guarantees existence of a strict module dilation. Apart from Question \ref{Q3} about existence of weak but no strong dilations, this leaves us with the one big question:

Does embeddability of the GNS-subproduct system of a Markov semigroup into a \bf{super}product system guarantee existence of a dilation?

\bf{Answer:} Unknown.

However, if we could show that a superproduct system containing the GNS-subproduct system does embed into a product system, we would show even existence of a module dilation.
\eemp

\bemp[Question:~] \label{Q7}
This point consists actually of several questions that all regard the problem of, given the GNS-subproduct system sitting in a superproduct system (like, for instance, the one deduced from a strong or a good dilation), how to find a product system containing both (and, thus, strictly related to answer the preceding question). We phrase the situation as a superproduct system with unit $\sE^\podot\ni\xi^\odot$, meaning that $\cls\cB\xi_t\cB$ are the members of the GNS-subproduct system of the CP-semigroup determined by $\xi^\odot$. The ``multi-example'' constructed in the Example-Section \ref{EXN02SEC} is of great help to sort out at least some things.

Does there, possibly, sit a product system $E^\odot$ in between, that is, $\sE^\podot\supset E^\odot\ni\xi^\odot$?

\bf{Answer:} No, not always. In Section \ref{EXN02SEC}, we construct a strong module dilation whose product system $E^\odot$ does not contain any other product subsystem $E'^\odot\ni\xi^\odot$ containing the GNS-subproduct system, while the superproduct system $\sE^\podot\ni\xi^\odot$ of the minimalized dilation is a proper superproduct subsystem of $E^\odot$, therefore, not containing any $E'^\odot\ni\xi^\odot$.

However, in this example we see the situation $E^\odot\supset\sE^\podot\ni\xi^\odot$. Is this situation, maybe, typical? Is it always possible to embed $\sE^\podot\ni\xi^\odot$ into a product system $E^\odot$?

\bf{Answer:} Unknown -- as unknown, of course, as the answer to Question \ref{Q6}.

In the example, $E^\odot$ is generated by $\xi^\odot$ (and, therefore, by the GNS-subproduct system) as a product system. Is, possibly, the pair $(E^\odot,\xi^\odot)$ determined uniquely by this property? Equivalently, does the semigroup have a unique GNS-\bf{product} system (as in the one-parameter case)?

\bf{Answer:} No. For the pair $(E^\odot,\xi^\odot)$ from Section \ref{EXN02SEC} from a pair of commuting CP-maps, the superproduct system generated by $\xi^\odot$ is a proper superproduct subsystem of $E^\odot$. However, if the two CP-maps commute strongly, then by Section \ref{strcomSEC} we may construct a different (strong) module dilation, and the product system of this dilation \bf{is} generated by $\xi^\odot$ as a \bf{superproduct} system. (So in one case case the superproduct system generated by $\xi^\odot$ is proper, in the other case it is not.) We see, if it should be possible to construct a containing product system, then we may not expect that the construction be \it{universal}.
\eemp

\newpage

\section{Superproduct systems and unitalizations} \label{unisupSEC}

After this account, in the present and the following two sections, we continue our analysis of superproduct systems of dilations, providing also some of the answers anticipated in the answers to the questions in the preceding section. In this section, we analyze the relation between the superproduct system of a strong dilation of a CP-semigroup and that of its unitalization. Then, in the next section, we examine so-called left dilations that relate (super)product systems to \nbd{E}semigroups, providing also the promised inductive limit construction for product systems over Ore monoids with unital units, leading to a strict module \nbd{E_0}dilation of the related Markov semigroup. Finally, in Section \ref{EXpropsupSEC}, we present examples of dilations that have associated a proper superproduct system (illustrating, that this phenomenon is not of a topological nature but algebraic).

\lf
In Theorems \ref{wstrunithm} and \ref{uninonunithm}, we have seen how to relate strong dilations of CP-semigroups with dilations of the Markov semigroups obtained by unitalization. For a CP-semigroup $T$, let $(\wh{\cA},\wh{\theta},\wh{p})$ be a strong dilation of $\wt{T}$, and denote by $(\cA,\theta,p)$ the strong dilation of $T$ as constructed in the proof of Theorem \ref{uninonunithm}. (Recall that, in the notation of the proof of that result, we have $\wh{\cA}\supset\wh{p}\wh{\cA}\wh{p}=\wt{\cB}\supset\cB\ni p=\U_\cB$ so that $q:=\wh{p}-p=\U_{\wt{\cB}}-\U_\cB\in\wh{\cA}$ and we have $\cA:=\U\wh{\cA}\U$ with $\U:=\wh{\U}-q$ is invariant under $\wh{\theta}$ so that we may define the (co)restriction to $\cA$, $\theta$.) Since $p$ and $\wh{\theta}_t(p)$ are both in $\cA$, it follows that $\U p=p$ and that $\wh{\theta}_t(p)\U=\wh{\theta}_t(p)$. Consequently, the superproduct system $\wh{\sE}^\podot=\bfam{\wh{\sE}_t}_{t\in\bS}$ with $\wh{\sE}_t:=\wh{\theta}_t(\wh{p})\wh{\cA}\wh{p}$ of $(\wh{\cA},\wh{\theta},\wh{p})$ contains the superproduct system $\sE^\podot=\bfam{\sE_t}_{t\in\bS}$ with $\sE_t:=\theta_t(p)\cA p$ in the sense that
\beqn{
\sE_t
~=~
\wh{\theta}_t(p)\U\wh{\cA}\U p
~=~
\wh{\theta}_t(p)\wh{\sE}_tp,
}\eeqn
and also $\xi_t=\wh{\theta}_t(p)p=\wh{\theta}_t(p)\wh{\theta}_t(\wh{p})\wh{p}p=\wh{\theta}_t(p)\wh{\xi}_tp$. In other words, we are in the situation discussed in the last part of Section \ref{SPSpbSEC}, with the ideal $\cI=\cB$ in $\wt{\cB}$. Since $p$ is a central projection in $\wt{\cB}$, we see that $\wh{x}_t\mapsto\wh{\theta}_t(p)\wh{x}_tp$ defines a bilinear projection in $\sB^{a,bil}(\wh{\sE}_t)$ onto $\sE_t$. The following result is an immediate consequence of Corollary \ref{asubscor}.

\bprop
$\sE^\podot$ is adjointable, if $\wh{\sE}^\podot$ is adjointable.
\eprop

\bcor
If $\wh{\sE}^\podot$ is a product system, then $\sE^\podot$ is adjointable.
\ecor

To say what $\wh{\sE}^\podot$ or $\sE^\podot$ being a product system implies for the other one, is not possible in general.
% \OW[OPEN]{Check: $\wh{\sE}\leftrightarrow\wt{\sE}$! (I don't understand, this is addressed below. --Orr)}
The situation improves considerably, when we switch from a general dilation $(\wh{\cA},\wh{\theta},\wh{p})$ of $\wt{T}$ to $(\wh{\cA},\wh{\theta},\wh{p})=(\wt{\cA},\wt{\theta},\wt{p})$, the unitalization of $(\cA,\theta,p)$ (as explained in Theorem \ref{wstrunithm}). In that case, we denote the superproduct system by $\wt{\sE}^\podot$. The big difference is that, here, $q$ is a central projection not only for $\wt{\cB}$, but for $\wt{\cA}$, too.

The crucial property that made work Corollary \ref{PSIPScor} has been referred to as \it{triangular} in Remark \ref{PSIPSrem}.

\bcor
$\sE^\podot$ is \hl{triangular}\index{product system!triangular} in $\wt{\sE}^\podot$ with respect to $p$, that is, $p\wt{\sE}_tq=\zero$.
\ecor

\proof
$aq=0$ for all $a\in\cA$, and $\theta_t(p)\wt{a}\in\cA$ for all $\wt{a}\in\wt{\cA}$.\qed

\lf
And by Corollary \ref{PSIPScor}:

\bprop \label{sppspsprop}
$\sE^\podot$ is a product system, if $\wt{\sE}^\podot$ is a product system.
\eprop

As for the converse, we have to be more specific.

\blem \label{pswtpslem}
$\wt{\sE}^\podot$ is a product system if and only if $\sE^\podot$ is a product system and
\beq{\label{PSuniPSeq}
\ol{\theta_t(\U-\theta_s(\U))\cA p}
~=~
\cls\theta_t(\U-\theta_s(\U))\theta_t(\cA p)\cA p
}\eeq
for all $t,s\in\bS$.
\elem

\proof
Recall that $\cls\wt{\sE}_s\wt{\sE}_t=\wt{\sE}_{st}$ if and only if $\cls\wt{\theta}_{st}(\wt{p})\wt{\theta}_t(\wt{\cA}\wt{p})\wt{\cA}\wt{p}\ni\wt{\theta}_{st}(\wt{p})\wt{a}\wt{p}$ for all $\wt{a}\in\wt{\cA}$. Since $\wt{\sE}_tq=\wt{\theta}_t(\wt{p})\wt{\cA}\wt{p}q=\C q$, we always have $\cls\wt{\sE}_s\wt{\sE}_tq=\C q=\wt{\sE}_{st}q$. Also, since $\wt{\sE}^\podot$ is triangular,
\beqn{
\cls p\wt{\sE}_s\wt{\sE}_tp
~=~
\cls p\wt{\sE}_sp\wt{\sE}_tp
~=~
\cls\sE_s\sE_t.
}\eeqn
So, $p\wt{\sE}_{st}p=\sE_{st}$ and $\cls p\wt{\sE}_s\wt{\sE}_tp=\cls\sE_s\sE_t$ coincide if and only if $\sE^\podot$ is a product system. The only components still to be compared are $q\wt{\sE}_{st}p$ with $\cls q\wt{\sE}_s\wt{\sE}_tp$.

First of all, $\wt{\theta}_{st}(q)=\wt{\theta}_{st}(\wt{\U}-\U)=\wt{\U}-\theta_{st}(\U)$. Next, since $\wt{\cA}p=\cA p\subset\cA$, in the expressions to be compared we may write $\U-\theta_{st}(\U)$ for $\wt{\theta}_{st}(q)$. Finally, we observe that
\beqn{
\U-\theta_{st}(\U)
~=~
(\U-\theta_t(\U))+\theta_t(\U-\theta_s(\U))
}\eeqn
(recall that $\theta$ is a semigroup over $\bS^{op}$), where, as $\theta_t(\U)$ is decreasing, $\U-\theta_t(\U)$ and $\theta_t(\U-\theta_s(\U))$ are a pair of orthogonal projections. In a typical element $\wt{\theta}_{st}(q)\wt{\theta}_t(\wt{a'}\wt{p})a'' p\in\wt{\theta}_{st}(q)\wt{\theta}_t(\wt{\cA}\wt{p})\cA p=\wt{\theta}_{st}(q)\wt{\theta}_t(\wt{\cA}\wt{p})\wt{\cA}p$ we have two cases, namely, $\wt{a'}=q$ (so that $\wt{a'}\wt{p}=q$) and $\wt{a'}=a'\in\cA$ (so that $\wt{a'}\wt{p}=a'p$). For the first case, we have $\wt{\theta}_{st}(q)\wt{\theta}_t(q)=\wt{\theta}_t(q)$. So, with $(\U-\theta_t(\U))a''p$ we get every element of the form $(\U-\theta_t(\U))ap\in q\wt{\sE}_{st}p$. For the second case, we note that $(\U-\theta_t(\U))\theta_t(a'p)=(\U-\theta_t(\U))\theta_t(\U)\theta_t(a'p)=0$. So, in the range of $\theta_t(\U-\theta_s(\U))$ we get all elements $\theta_t(\U-\theta_s(\U))\theta_t(a'p)a''p$ and what is in their closed linear span, but not more. Therefore, given that $\sE^\podot$ is a product system, $\wt{\sE}^\podot$ is a product system if and only if \eqref{PSuniPSeq} holds.\qed

\lf
Explicit positive examples, where $\wt{\sE}^\podot$ is a product system, can be found in Bhat and Skeide \cite[Section 8]{BhSk00} and, more generally, in Skeide \cite{Ske08}. An explicit negative example, where $\sE^\podot$ is a product system but $\wt{\sE}^\podot$ is not, is missing.

\newpage

\section{Superproduct systems and left dilations} \label{leftdilSEC}

We leave the discussion of superproduct systems of unitalizations, and imitate now the inductive limit construction of a dilation from a product system and a unital unit from Bhat and Skeide \cite{BhSk00}. The generalization is twofold. Firstly, we start from a superproduct system (leading to an \nbd{E}semigroup if and only if the superproduct system is a product system), and secondly, we replace $\R_+$ by an arbitrary Ore monoid. The reasons why we discuss this construction for superproduct systems, are twofold. Firstly, doing the construction for superproduct systems does not create any problems as compared with doing it for product systems; secondly, only actually doing the construction for superproduct systems, explains why in the end the output necessarily fails to give what we want if the superproduct system is proper.

To understand better what we wish to get, let us return to the situation in the beginning of Section \ref{DilSPSpSEC}, culminating in Theorem \ref{E-supPSthm}. Recall (Example \ref{EPSex}) that the representation theory of $\sB^a(E)$ applied to a strict \nbd{E_0}semigroup $\vt$ does not only provide us with a product system $E_t=E^*\odot_tE$, but that it allows to recover $\vt_t$ as $v_t(\bullet\odot\id_t)v_t^*$,  for unitaries $v_t\colon E\odot E_t\rightarrow E$. (After all, the scope of the representation theory \bf{is} to recover a homomorphism $\vt_t$ as amplification with a multiplicity correspondence $E_t$. The fact that for a semigroup the multiplicity correspondences form a product system, is a ``second order effect''.) The $v_t$ are defined by $x\odot(y^*\odot _tz)\mapsto\vt_t(xy^*)z$, and they satisfy the associativity condition that the \hl{product}
\beqn{
(x,y_t)
~\longmapsto~
xy_t
~:=~
v_t(x\odot y_t)
}\eeqn
iterates associatively with the product system structure: $(xy_s)z_t=x(y_sz_t)$. Generally, if $E$ is a full(\bf{!}) Hilbert module and a family of unitaries $v_t\colon E\odot E_t\rightarrow E$ fulfills the associativity condition, since Skeide \cite{Ske06,Ske07} we say the $v_t$ form a \phantomsection\hl{left dilation}\index{left dilation!of a product system}\index{product system!left dilation of} of the product system $E^\odot$ to the (full!) Hilbert module $E$. If we have a left dilation of $E^\odot$ to $E$, then, thanks to associativity, the maps $\vt_t:=v_t(\bullet\odot\id_t)v_t^*$ define a strict \nbd{E_0}semigroup on $\sB^a(E)$ (over $\bS^{op}$) and, thanks to fullness of $E$, the product system of this \nbd{E_0}semigroup is (isomorphic to) $E^\odot$. Note that $v_0$ is the canonical identification, automatically. (All these statement are practically \it{verbatim} generalization of  statements in Skeide \cite{Ske16} to arbitrary monoids $\bS$. We get them by restricting the following proposition to product systems.)

Formally, we extend the definition of left dilation also to superproduct systems. However:

\bprop\label{lsdilPSprop}
Let $\sE^\podot$ be a superproduct system over the monoid $\bS$ with a (unitary!) left dilation $\bfam{v_t}_{t\in\bS}$ to a (full!) Hilbert \nbd{\cB}module $\sE$. Then
\beqn{
\vt_t(a)
~:=~
v_t(a\odot\id_t)v_t^*
}\eeqn
defines a strict \nbd{E_0}semigroup $\vt$ over $\bS^{op}$ on $\sB^a(\sE)$. Moreover, the \index{left dilation!of a superproduct system implies product system}\index{superproduct system!left dilation of, implies product system}product system associated with $\vt$ is isomorphic to $\sE^\podot$.

Consequently, a superproduct system that admits a left dilation, is a product system.
\eprop

\proof
From associativity we get
\bmu{\label{sgass}
\vt_t\circ\vt_s(a)
~=~
v_t((v_s(a\odot\id_s)v_s^*)\odot\id_t)v_t^*
~=~
v_t(v_s\odot\id_t)(a\odot\id_s\odot\id_t)(v_s^*\odot\id_t)v_t^*
\\
~=~
v_{st}(\id_\sE\odot v_{s,t})(a\odot\id_s\odot\id_t)(\id_\sE\odot\,v_{s,t}^*)v_{st}^*
~=~
v_{st}(a\odot\id_{st})(\id_\sE\odot\,v_{s,t})(\id_\sE\odot\,v_{s,t}^*)v_{st}^*
\\
~=~
v_{st}(a\odot\id_{st})(\id_\sE\odot\,v_{s,t}v_{s,t}^*)v_{st}^*.
}\emu
Since $\vt_t$ and $\vt_s$ are unital and since $v_{st}$ is unitary, inserting $a=\id_\sE$, we see that $\id_\sE\odot\,v_{s,t}v_{s,t}^*$ is the identity of $\sE\odot\sE_{st}$. (Since $\sE$ is full, this alone is already enough to show, by tensoring with $\id_{\sE^*}$, that $v_{s,t}v_{s,t}^*=\id_{st}$ so that $\sE^\podot$ is a product system. Here, we recover that in a different way.) In conclusion of \eqref{sgass}, we get
\beqn{
\vt_t\circ\vt_s(a)
~=~
v_{st}(a\odot\id_{st})v_{st}^*
~=~
\vt_{st}(a),
}\eeqn
that is, the $\vt_t$ form a semigroup over $\bS^{op}$.

Clearly, $\vt_t$ is a strict unital endomorphism of $\sB^a(\sE)$. It is immediate that the multiplicity correspondence $E_t:=\sE^*\odot_t\sE$ is isomorphic to $\sE_t$ via $x^*\odot_t(yz_t)\mapsto\AB{x,y}z_t$. The computation
\bmun{
(x^*\odot_s(yz_s))\odot(x'^*\odot_t(y'z'_t))
~\longmapsto~
x^*\odot_{st}\vt_t((yz_s)x'^*)(y'z'_t)
~=~
x^*\odot_{st}v_t((yz_s)x'^*\odot\id_t)v_t^*(y'z'_t)
\\
~=~
x^*\odot_{st}v_t((yz_s)x'^*\odot\id_t)(y'\odot z'_t)
~=~
x^*\odot_{st}v_t((yz_s)\odot\AB{x',y'}z'_t)
~=~
x^*\odot_{st}yz_s\AB{x',y'}z'_t
}\emun
shows that the product of the product system $E^\odot$ of the \nbd{E_0}semigroup coincides with the (super)product $\AB{x,y}z_s\odot\AB{x',y'}z'_t\mapsto\AB{x,y}z_s\AB{x',y'}z'_t$ of $\sE^\podot$. Therefore, $E^\odot$ and $\sE^\odot$ are isomorphic superproduct systems. Since $E^\odot$ is a product system, so is $\sE^\podot$.\qed

\lf
The construction of a superproduct system for a triple $(\cA,\theta,p)$ in Theorem \ref{E-supPSthm} was motivated by modifying the construction of a product system from a strict \nbd{E}semigroup on $\sB^a(E)$; and we anticipated already that there are triples where we get only a superproduct system and that this does not disappear under continuity conditions. Let us now see what we get, if we imitate also the construction of the left dilation obtained from a strict \nbd{E_0}semigroup on $\sB^a(E)$.

First of all, already the fact that we get a left dilation (unitaries!) depends on that we start with an \nbd{E_0}semigroup. Already, if we start with a strict \nbd{E}semigroup on $\sB^a(E)$, then what we get are isometries $v_t\colon E\odot E_t\rightarrow E$. Here, we have to add that $v_0$ is unitary. As long as $E$ is full (guaranteeing suitable uniqueness) and $v_0$ is unitary, we say the $v_t$ form a \phantomsection\hl{left semidilation}\index{left dilation!semi}\index{product system!left semidilation of}. The isometries are onto the complemented submodule $\vt_t(\id_E)E$ of $E$ and, therefore, adjointable, that is, the left semidilation is \hl{adjointable}\index{left dilation!semi!adjointable}\index{product system!left semidilation of!adjointable}\index{adjointable!left semidilation}.

Also here, we relax the notion of (adjointable) left semidilations to superproduct systems. We get a bit less than Proposition \ref{lsdilPSprop}:

\bob \label{lsdilPSob}
Proposition \ref{lsdilPSprop} remains true for adjointable left semidilations and \nbd{E}semi\-groups provided $\sE^\podot$ \bf{is} a product system. (Equation \eqref{sgass} depends only on associativity, and the conclusion to get the next equation goes through directly, because by assuming that $\sE^\podot$ \bf{is} a product system, we do not have to show first that the $v_{s,t}$ are coisometries.) However, the proof also shows that we always have $\vt_t\circ\vt_s(\id_\sE)\le\vt_{st}(\id_\sE)$.  And in order to verify if a superproduct system admitting an adjointable left semidilation is a product system, it is enough to check if in this inequality we have equality, that is, to check the semigroup property of $\vt$ at the single element $\id_\sE$.

It is clear that also for left semidilations the $\vt_t$ cannot form a semigroup, unless $\sE^\podot$ is a product system. (If $\vt$ is a semigroup, then, as before for left dilations, the structure of $\sE^\podot$ is isomorphic to the structure of the product system associated with $\vt$.) We see in a minute that there exist proper superproduct systems that admit adjointable left semidilations (so that the associated $\vt_t$ cannot form a semigroup).

(It is interesting to compare the situation in \eqref{sgass} with the corresponding step in the proof of Theorem \ref{adSPS-CPthm}, where the $v_t$ also look like a sort of left dilation. We emphasize however, that there we were speaking about (adjointable) subproduct systems (guaranteeing in the verification of the semigroup property the $v_{s,t}=w_{s,t}^*$ to be coisometries), while here we are speaking about superproduct systems.)
\eob

\lf
It is clear that without adjointability of the left semidilation, there is no such analogue of Proposition \ref{lsdilPSprop} discussed in the preceding observation, because $\vt_t$ cannot be defined. However:

\bprop \label{semiadjprop}
A superproduct system admitting an adjointable left semidilation, is adjointable.
\eprop

\proof
Associativity of the $v_t$ with the $v_{s,t}$ is $v_{st}(\id_\sE\odot\,v_{s,t})=v_t(v_s\odot\id_t)$. If all (isometric!) $v_t$ have adjoints, we get $\id_\sE\odot\,v_{s,t}=v_{st}^*v_t(v_s\odot\id_t)$. Since the right-hand side is adjointable, so is $\id_\sE\odot\,v_{s,t}$. Since $\sE$ is full, also $v_{s,t}$ is adjointable. (In full generality, if there is $a\in\sB^{bil}(F_1,F_2)$ and full $E$ such that $\id_E\odot\,a$ is adjointable, then tensoring with $\id_{E^*}$, taking also into account that $E^*\odot E\odot F_i=F_i$, shows that $a$ is adjointable. In our situation, where $\sE$ usually has a unit vector $\xi$, we get that $v_{s,t}=(\xi^*\odot\id_{st})(\id_\sE\odot\,v_{s,t})(\xi\odot\id_{st})$ is adjointable.)\qed

\lf
What we get from a triple $(\cA,\theta,p)$ in the notations of Theorem \ref{E-supPSthm}, is the following.

\bthm \label{E-lsemdilthm}
The maps $v_t$ defined by
\beqn{
v_t
\colon
x\odot y_t
~\longmapsto~
\theta_t(x)y_t
}\eeqn
$(x\in\cA p$, $y_t\in\theta_t(p)\cA p)$ define a left semidilation of $\sE^\podot$ to $\sE:=\cA p$.
\ethm

The proof of isometricity (and associativity) goes exactly as in the proof of Theorem \ref{E-supPSthm}, replacing $x_s\in\sE_s$ (and $x_r\in\sE_r$) with general $x\in\sE$.

% \bob
% For $a\in\cA$, denote $\sB^a(\sE)\ni a_\sE\colon x\mapsto ax$ (see Proposition \ref{convprop}\eqref{cp1}).
% \eob

\bemp[Question.] \label{pstriQ}
Are the $v_t$ adjointable (implying that also $\sE^\podot$ would be adjointable)? Equivalently, is $v_t(\sE\odot\sE_t)=\cls\theta_t(\cA p)\cA p$ complemented in $\sE$?

One can show that the answer is affirmative, if $\theta_t$ is \phantomsection\hl{\nbd{p}relatively strict}\index{prelatively@\nbd{p}relatively strict|bf} \index{strict!prelatively@\nbd{p}relatively}\index{dilation, weak!prelatively@\nbd{p}relatively strict}in the sense that for a bounded approximate unit of $\ls\cA p\cA$ the action of the image of this approximate unit under $\theta_t$ on $\sE$ is strictly Cauchy in $\sB^a(\sE)$. (In that case, the limit is a projection in $\sB^a(\sE)$ onto $\cls\theta_t(\cA p\cA)\cA p$ and this equals $\cls\theta_t(\cA p)\theta_t(\cA)\cA p=\cls\theta_t(\cA p)\cA p$. Indeed, the latter is, clearly, not smaller than the former, and by applying the projection we see that it is not bigger, either.) So, for reasonably good triples, superproduct system and left dilation are adjointable. In the Example-Section \ref{EXpropsupSEC}, we provide even \nbd{p}relatively strict dilations of Markov semigroups that lead to adjointable superproduct systems that are not product systems and, consequently, to adjointable left semidilations that do not determine \nbd{E}semigroups. We come back to these discussions of topological character in Section \ref{topSEC}.
\eemp

\lf
We now perform the promised inductive limit construction. Let $\sE^\podot$ be a superproduct system over the Ore monoid $\bS$ with structure maps $v_{s,t}\colon\sE_s\odot\sE_t\rightarrow\sE_{st}$, and let $\xi^\odot$ be a unital unit for $\sE^\podot$. As usual, we define the product $x_sy_t:=v_{s,t}(x_s\odot y_t)$. Since $\xi_s$ is a unit vector, by $x_t\mapsto\xi_sx_t$ we define isometric (right linear but, in general, not bilinear) embeddings $\sE_t\rightarrow\sE_{st}$. (In the verification
\beqn{
\AB{\xi_sx_t,\xi_sx'_t}
~=~
\AB{v_{s,t}(\xi_s\odot x_t),v_{s,t}(\xi_s\odot x'_t)}
~=~
\AB{\xi_s\odot x_t,\xi_s\odot x'_t}
~=~
\AB{x_t,\AB{\xi_s,\xi_s}x'_t}
~=~
\AB{x_t,x'_t}
}\eeqn
that this map is an isometry, we used that we may leave out an isometry if it occurs on the left in both arguments of an inner product. So, nothing like this would work for unital units of subproduct systems.) Since an Ore monoid is right-reversible, it is directed. Since it is cancellative, the element $r$ that illustrates that $s\le t=rs$, is unique. We, therefore, may define unique isometries $\gamma_{t,s}\colon\sE_s\rightarrow\sE_t$ as $x_s\mapsto\xi_rx_s$, whenever $s\le t$. By marginal conditions and associativity of the product of $\sE^\podot$, we get
\baln{
\gamma_{t,t}
&
~=~
\id_t,
&
\gamma_{t,s}\gamma_{s,r}
&
~=~
\gamma_{t,r}
}\ealn
for all $r\le s\le t\in\bS$, that is,  the $\gamma_{t,s}$ form an \phantomsection\hl{inductive system}\index{inductive!system}. Forming the \hl{inductive limit}\index{inductive!limit} over that system, we get a Hilbert \nbd{\cB}module $\sE:=\limind_t\sE_t$ with isometries $k_t\colon\sE_t\rightarrow\sE$ fulfilling
\baln{
\sE
&
~=~
\ol{\bigcup_{t\in\bS}k_t\sE_t},
&
k_t\gamma_{t,s}
&
~=~
k_s
}\ealn
for all $s\le t\in\bS$. (See \cite[Appendix]{BhSk00} for details about the construction. The brackets around ``(with $s\ge T$ ...)'' in\cite[p.569, l.1]{BhSk00} are misleading and should be canceled.) Note that $\sE$ and the family $k_t$ are determined by the stated properties up to suitable unitary equivalence. The same is true for the following standard \hl{universal property}\index{inductive!limit!universal property}: Given a uniformly bounded family $a_t\in\sB^r(\sE_t,F)$ such that $a_s=a_t\gamma_{t,s}$ for all $s\le t\in\bS$, there is a unique $a\in\sB^r(\sE,F)$ such that $ak_t=a_t$ for all $t\in\bS$.

\bob\label{indlimob} ~

\begin{enumerate}
\item\label{il1}
The inner product of $k_tx_t$ and $k_sy_s$ has no choice but being determined by finding $t',s'$ such that $t't=s's$ and computing
\beqn{
\AB{k_tx_t,k_sy_s}
~=~
\AB{k_{t't}\gamma_{t't,t}x_t,k_{s's}\gamma_{s's,s}y_s}
~=~
\AB{\gamma_{t't,t}x_t,\gamma_{s's,s}y_s},
}\eeqn
both elements in $\sE_{t't}=\sE_{s's}$. (Actually, proving this is well-defined, proves that the inductive limit has an inner product.)

\item\label{il2}
If $s\le t$, then $k_s\sE_s=k_t\gamma_{t,s}\sE_s\subset k_t\sE_t$.

\item\label{il3}
For each $t\in\bS$ the set $\bS t$ is \phantomsection\hl{cofinal}\index{cofinal} in $\bS$, that is, for each $r\in\bS$ exists $st\in\bS t$ such that $r\le st$. Consequently, $\sE=\ol{\bigcup_{s\in\bS}\sE_{st}}$.

\item\label{il4}
If all $\gamma_{t,s}$ are adjointable, then so are all $k_t$. (Indeed, for $x_s\in\sE_s$ find $r$ such that $t\le r,s\le r$ and put $k_t^*(k_sx_s):=\gamma_{r,t}^*\gamma_{r,s}x_s$. This does not only define an adjoint $k_t^*$ of $k_t$, but it also shows that that $\gamma_{r,t}^*\gamma_{r,s}\in\sB^a(\sE_s,\sE_t)$ does not depend on the choice of $r\ge t,s$.)

Note that adjointability of our concrete $\gamma_{st,t}\colon x_t\mapsto\xi_sx_t=v_{s,t}(\xi_s\odot x_t)$ depends on adjointability of the product maps $v_{s,t}$ of the superproduct system. The maps $\xi_s\odot\id_t\colon\sE_t\rightarrow\sE_s\odot\sE_t$ have adjoints, namely, $\xi_s^*\odot\id_t$.

\item\label{il5}
But even if all $\gamma_{t,s}$ are adjointable and even if all $a_t$ in the universal property are adjointable, the operator $a$ need not be adjointable. (Indeed, let $\cB:=\wt{C_0(\N)}$ denote the convergent sequences. Let $\sE_n:=\I_{\SB{1,n}}\cB$ denote the sequences that are $0$ starting from the $n+1$ term. The inductive system is the canonical embedding $\gamma_{n+m,n}\colon\sE_n\rightarrow\sE_{n+m}$. The inductive limit is $C_0(\N)$ with $k_n$ again the canonical embeddings into $C_0(\N)$. Let $a_n$ denote the canonical embeddings into $\wt{C_0(\N)}=\cB$. They determine $a$ to be the canonical embedding of $C_0(\N)$ into $\cB$. Although all $\gamma_{n,m}$ and all $a_n$ are adjointable, $a$ is not.)
%%%% BO
% I removed the parenthesis around the example above, since the following sentence refers to it; it seems wrong to refer to a parenthetical remark from outside parenthesis. --Orr
%%%% EO

An increasing net (even a sequence) of projections need not converge strongly. (This follows from the same example.) If it converges strongly, however, then the limit is a projection, too.
\end{enumerate}
\eob

\noindent
Like in \cite{BhSk00}, the elements $\xi_t$ of the unit, under the embeddings $k_t$ are mapped to the same vector $\xi\in\sE$. (Indeed, for $t,s\in\bS$ choose $t',s'\in\bS$ such that $t't=s's$. Then $k_t\xi_t=k_{t't}\gamma_{t't,t}\xi_t=k_{t't}(\xi_{t'}\xi_t)=k_{s's}(\xi_{s'}\xi_s)=k_{s's}\gamma_{s's,s}\xi_s=k_s\xi_s$.) Obviously, $\xi$ is a unit vector, so that, in particular, $\sE$ is full.

Like in \cite{BhSk00}, it is our scope to extend the embedding $v_{s,t}\colon\sE_s\odot\sE_t\rightarrow\sE_{st}$ to an embedding $v_t\colon\sE\odot\sE_t\rightarrow\sE$, roughly speaking, by sending ``$s\to\infty$''. More precisely, we define $v_t$ on the dense subspace $\bigcup_{s\in\bS}k_s\sE_s\odot\sE_t$ as
\beqn{
v_t(k_sx_s\odot y_t)
~:=~
k_{st}(x_sy_t).
}\eeqn
For $r,s\in\bS$, choose $r',s'\in\bS$ such that $s's=r'r$. From
\bmun{
\AB{k_{st}(x_sy_t),k_{rt}(x'_ry'_t)}
~=~
\AB{k_{s'st}(\xi_{s'}x_sy_t),k_{r'rt}(\xi_{r'}x'_ry'_t)}
~=~
\AB{\xi_{s'}x_sy_t,\xi_{r'}x'_ry'_t}
\\
~=~
\AB{v_{s's,t}(\xi_{s'}x_s\odot y_t),v_{r'r,t}(\xi_{r'}x'_r\odot y'_t)}
~=~
\AB{\xi_{s'}x_s\odot y_t,\xi_{r'}x'_r\odot y'_t}
\\
~=~
\AB{y_t,\AB{\xi_{s'}x_s,\xi_{r'}x'_r}y'_t}
~=~
\AB{y_t,\AB{k_sx_s,k_rx'_r}y'_t}
~=~
\AB{k_sx_s\odot y_t,k_rx'_r\odot y'_t},
}\emun
we see that this, indeed, (well-)defines an isometry. With the usual product notation $xy_t:=v_t(x\odot y_t)$, from
\beqn{
(k_rx_r)(y_sz_t)
~=~
k_{rst}(x_r(y_sz_t))
~=~
k_{rst}((x_ry_s)z_t)
~=~
(k_{rs}(x_ry_s))z_t
~=~
((k_rx_r)y_s)z_t
}\eeqn
we see that the isometries form a left semidilation. Like in \cite{BhSk00}, we get
\beq{\label{ktxi}
\xi x_t
~=~
v_t(\xi\odot x_t)
~=~
v_t(k_s\xi_s\odot x_t)
~=~
k_{st}(\xi_sx_t)
~=~
k_tx_t.
}\eeq

\lf
Unlike in \cite{BhSk00}, the  $v_t$ need not be unitary. Effectively:

\bprop \label{lduniprop}
The $v_t$ are unitary (that is, they form a left dilation) if and only if $\sE^\podot$ is a product system.
\eprop

\proof
We know from Proposition \ref{lsdilPSprop} that for having a left dilation, $\sE^\podot$ has to be a product system. Conversely, if $\sE^\podot$ is a product system, then $v_t(k_s\sE_s\odot\sE_t)=k_{st}\cls\sE_s\sE_t=k_{st}\sE_{st}$. By Observation \ref{indlimob}\eqref{il3}, the latter increase to a dense subspace of $\sE$.\qed

\bthm\label{Oreindthm}
Let $\sE^\odot$ be a product system over an Ore monoid $\bS$ with a unital unit $\xi^\odot$, and construct inductive limit $\sE\ni\xi$ and left dilation $v_t$ as above. Then $(\sE,\vt,\xi)$ with $\vt_t(a):=v_t(a\odot\id_t)v_t^*$ is an \nbd{E_0}dilation of the Markov semigroup $T$ over $\bS^{op}$ with $T_t:=\AB{\xi_t,\bullet\xi_t}$. 
\ethm

\proof
By $\xi=k_t\xi_t=\xi\xi_t$, we find
\vspace{-1ex}
\beqn{
\AB{\xi,\vt_t(\xi b\xi^*)\xi}
~=~
\AB{\xi\odot\xi_t,((\xi b\xi^*)\odot\id_t)(\xi\odot\xi_t)}
~=~
\AB{\xi\odot\xi_t,\xi\odot b\xi_t}
~=~
\AB{\xi_t,b\xi_t}.\qedsymbol
}\eeqn
\noqed

\vspace{-3ex}
\noindent
Since, as we know, every strong strict full dilation comes along with a product system (Theorem \ref{E-supPSthm} and Example \ref{EPSex}) containing the GNS-subproduct system (Theorem \ref{sdilunithm}\eqref{DU1} and Observation \ref{uGNSsupob}), we get even more than what we asked in Question \ref{Q4}:

\bthm \label{Markmodthm}
Let $T$ be a Markov semigroup over the opposite of an Ore monoid. Then $T$ admits a strict full (or module) dilation if and only if the GNS-subproduct system of $T$ embeds into a product system.
\ethm

\brem
Let us repeat that the embeddability of subproduct systems in product systems was \bf{the} motivation for Shalit and Solel \cite{ShaSo09}. The pair \cite[Theorem 5.12 and Corollary 5.10]{ShaSo09} states the preceding result in the \nbd{d}parameter case for von Neumann algebras in terms of the \it{Arveson-Stinespring subproduct system}, the \it{commutant} 
%%%% BO postarxiv
%if
of 
%%%% EO postarxiv
the GNS-subproduct system; see Appendix \ref{vNAPP}\ref{vNcomm} for details.
\erem
We close with a few general results with decreasing degree of proved usefulness. The following observation about the left semidilation in Theorem \ref{E-lsemdilthm} and its \it{left subdilations} turns out to be useful in Subsection \ref{minSEC}\eqref{1-p-SSEC}.%, while the last two results about \hl{left quasi-(semi)dilations} (that is, a left (semi)\-dilation except for that it need not be to a full module) still have to find a useful application.

\bob \label{lsubdilob}
Let $(\cA,\theta,p)$ be a triple as in Theorem \ref{E-supPSthm} and denote be $v_t\colon\sE\odot\sE_t\rightarrow\sE$ its left semidilation as in Theorem \ref{E-lsemdilthm}. Furthermore, suppose $(\sF,\sF^\odot)$ is a \hl{left subdilation}\index{left dilation!semi!left subdilation of}\index{product system!left semidilation of!left subdilation of} of $v_t$, that is, $\sF\subset\sE$ and $\sF^\odot$ is a (by Proposition \ref{lsdilPSprop}, necessarily) product subsystem of $\sE^\podot$ such that $\sF\sF_t\subset\sF$ and such that the (co)restrictions of the $v_t$ define a left dilation of $\sF^\odot$ to $\sF$. Denote by $\vt$ the \nbd{E_0}semigroup on $\sB^a(\sF)$ induced by that left dilation. We examine how $\vt$ sits in $\theta$.
\begin{enumerate}
\item \label{lsd1}
For $a\in\cA$, denote $\sB^a(\sE)\ni a_\sE\colon x\mapsto ax$ (see Proposition \ref{convprop}\eqref{cp1}). Then
\beqn{
\theta_t(a)_\sE v_t(x\odot y_t)
~=~
\theta_t(a)\theta_t(x)y_t
~=~
\theta_t(a_\sE x)y_t
~=~
v_t((a_\sE x)\odot y_t)
~=~
\bfam{v_t(a_\sE\odot\id_t)v_t^*}\,v_t(x\odot y_t).
}\eeqn
In particular, if $(\sF,\sF^\odot)=(\sE,\sE^\podot)$ (so that $v_t$ is a left dilation), then $\theta_t(a)_\sE=\vt_t(a_\sE)$. Moreover, if $\sB^a(\sE)\subset\cA$ (in an obvious way; see the discussion in Section \ref{topSEC}, in particular, around Definition \ref{contdefi} and Theorem \ref{contthm}; but also the discussion about \nbd{p}relative strictness in Question \ref{pstriQ} is not irrelevant, here), then $\theta_t(a)=\vt_t(a)$ for all $a\in\sB^a(\sE)\subset\cA$ so that $\vt$ is the (co)restriction of $\theta$ to 
%%%% BO new
% $\sB^a(E)$.
$\sB^a(\sE)$.
%%%% EO

\item \label{lsd2}
In the general case $(\sF,\sF^\odot)\subset(\sE,\sE^\podot)$ we limit ourselves to the situation where $\sB^a(\sF)\subset\cA$. By this we mean that there is a projection $P\in\cA$ such that $P\sE=\sF$ and such that $P\cA P$ and $\sB^a(\sF)$ are isomorphic via $PaP\mapsto\bSB{(PaP)_\sF\colon y\mapsto PaPy}$.

We easily verify that $P\theta_t(a)P=\vt_t(a)$ for all $a\in\sB^a(\sF)\subset\cA$. (Clearly, $\theta_t(a)Py=\vt_t(a)y\in\sF$ for $y\in\sF$, so that also $P\theta_t(a)Py=\vt_t(a)y$.) Except for that we do not require that $(\cA,\theta,p)$ is a dilation nor that $P\ge p$, this is what we phrase in Subsection \ref{minSEC}\ref{primSSEC} as $P$ \it{compresses} $\theta$ on $\cA$ to $\vt$ on $P\cA P=\sB^a(\sF)$.

\item \label{lsd3}
All questions about containment of $\sB^a(\sE)$ or of $\sB^a(\sF)$ in $\cA$ disappear when $\theta$ is a normal \nbd{E}semigroup on the a von Neumann algebra $\cA\ni p$. (Simply, an approximate unit for $\sF(\sF)\subset\cA$ will converge strongly in $\cA$ to $P$.) See, again, the already mentioned discussion in Section \ref{topSEC}.
\end{enumerate}
\eob

\lf\noindent
After this useful observation, here are two results about \phantomsection\hl{left quasi-(semi)dilations}\index{left dilation!quasi-}\index{product system!left quasi-dilation of} (that is, a left (semi)\-dilation to a module that need not be full) still have to find a useful application.% Our hope was that they allowed to show that the superproduct system of a Markov dilation always contained a product system plus unit. But, this is true only in the one-parameter case; see Subsection \ref{minSEC}\ref{1-p-SSEC}. In Subsection \ref{EXN02SEC}\ref{spNpsSSEC}, we obtain discrete two-parameter counter examples.

\blem\label{semisublem}
Let $\sE^\podot=\bfam{\sE_t}_{t\in\bS}$ be a superproduct system over the monoid $\bS$ and let $v_t\colon\sE\odot\sE_t\rightarrow\sE$ be a left quasi-semidilation of $\sE^\podot$ to $\sE$.

Then
\vspace{-2ex}
\beqn{
E
~:=~
\cls\bigcap_{\substack{n\in\N\\t_n,\ldots,t_1\in\bS}}\sE\sE_{t_n}\ldots\sE_{t_1}
\vspace{-1.5ex}
}\eeqn
is a Hilbert submodule of $\sE$, for each $t\in\bS$
\vspace{-1ex}
\beqn{
E_t
~:=~
\bCB{x_t\in\sE_t\colon Ex_t\subset E}
\vspace{-1ex}
}\eeqn
is a subcorrespondence of $\sE_t$, the $E_t$ form a superproduct subsystem $E^\podot$ of $\sE^\podot$, and the (co) restrictions of the $v_t$ to maps $E\odot E_t\rightarrow E$ form a left quasi-semidilation of $E^\podot$ to $E$.

Moreover, if the left quasi-semidilation $v_t$ is \hl{properly semi} (that is, if $v_t$ is nonunitary for some $t\in\bS$), then $E\ne\sE$.
\elem

\proof
By definition, $E$ is a closed linear subspace of $\sE$, which is invariant under right multiplication by $b\in\cB$ because each $\sE_{t_1}$ is.

$E_t$ is closed under addition because $E$ is, $E_t$ is closed under right multiplication because $E$ is, and $E_t$ is closed under limits because $E$ is. It is closed under left multiplication because $E$ is closed under right multiplication. (If $x_t\in E_t$, then $E(bx_t)=(Eb)x_t\subset Ex_t\subset E$, so $bx_t\in E_t$.)

If $x_s\in E_s$ and $y_t\in E_t$, then $E(x_sy_t)=(Ex_s)y_t\subset Ey_t\subset E$, so $x_sy_t\in E_{st}$. In other words, the correspondences $E_t$ form a superproduct subsystem of $\sE^\podot$. Of course, the (co)restrictions of $v_t$ form a left quasi-semidilation.

If $v_t$ is nonunitary, then $E\subset\cls\sE\sE_t=v_t(\sE\odot\sE_t)\subsetneq\sE$.
\qed

\lf
We say a pair $(E^\podot,E)$ of a superproduct subsystem $E^\podot$ of $\sE^\podot$ and a Hilbert submodule $E$ of $\sE$ fulfilling $EE_t\subset E$ is a \hl{left quasi-semisubdilation} of $v_t$.

\bprop\label{semisubprop}
Under the same hypotheses as in Lemma \ref{semisublem}: If $\bfam{({E^i}^\podot,E^i)}_{i\in I}$ is a family of left quasi-semisubdilations of $v_t$, then
\vspace{-1ex}
\beqn{
(E^\podot,E)
~:=~
\Bfam{\bigcap_{i\in I}{E^i}^\podot,\bigcap_{i\in I}E^i}
\vspace{-2.5ex}
}\eeqn
is a left quasi-semisubdilation of $v_t$, too.
\eprop

\proof
Since $x_s\in\bigcap_{i\in I}E^i_s$ and $y_t\in\bigcap_{i\in I}E^i_t$ implies $x_sy_t\in E^i_{st}$ for all $i$, we see that $E^\podot$ is superproduct subsystem of $\sE^\podot$. Since $x\in\bigcap_{i\in I}E^i$ and $y_t\in\bigcap_{i\in I}E^i_t$ implies $xy_t\in E^i$ for all $i$, we see that that $(E^\podot,E)$ is a left quasi-semisubdilation of $v_t$.\qed

\lf
The problem is that in both propositions $E$ and $E_t$ $\zero$. In fact, attempting to apply the results to the superproduct system of a dilation, contrary to what our hopes were, we were not able to shows that $E$ and $E_t$ contain $\xi$ and $\xi_t$, respectively. (In cases where this is true, both results together would imply existence of product subsystem $E^\odot\subset\sE^\podot$ containing the GNS-subproduct system.) It is true in the one-parameter case; see Subsection \ref{minSEC}\ref{1-p-SSEC}. In Subsection \ref{EXN02SEC}\ref{spNpsSSEC}, we obtain discrete two-parameter counter examples.

\newpage

\section[\sc{Examples:} Dilations with proper superproduct systems]{Examples: Dilations with proper superproduct systems} \label{EXpropsupSEC}

% \OW[POSTPONED checking! (Assuming mathematics is correct, checking and polishing to be done later.) (I checked, and corrected minor things -- Orr)]{I checked. --Orr}
We discuss two examples of a triple $(\cA,\theta,p)$, actually a dilation of a Markov semigroup, where $\theta$ is a one-parameter \nbd{E_0}semigroup that is strict and strongly time continuous in any reasonable topology on $\cA$ (different from the norm topology), and where the superproduct system according to Theorem \ref{E-supPSthm} is proper.

\bex \label{hypexex}
We assume known the notions from Appendix \ref{FockAPP} on the time ordered Fock module $\DG(F)$ over a \nbd{\cB}correspondence $F$ (the `time ordered sector' of the full Fock module $\sF(L^2(\R_+,F))$) including the definition of $L^2(\R_+,F)$, the time ordered product system $\DG^\odot(F)=\bfam{\DG_t(F)}_{t\in\T_+}$ (with $\DG_t(F)$ the `time ordered sector' of the full Fock module $\sF(L^2(\SB{0,t},F))$) and its product $u_{s,t}$, the left dilation $v_t\colon\DG(F)\odot\DG_t(F)\rightarrow\DG(F)$, and the CCR-flow over $F$ (the \nbd{E_0}semigroup determined by that left dilation). Recall, too, that $\DG(F)$ and left dilation can be identified with the inductive limit of $\DG^\odot(F)$ over the vacuum unit $\om^\odot$ according to Theorem \ref{Oreindthm}. Recall also the definition of exponential vectors
\beqn{
\ee(x)
~:=~
\sum_{n\in\N_0}\Delta_n x^{\odot n},
}\eeqn
whenever $x\in L^2(\R_+,F)$ is such that the sum exists (for instance, for a step function), and how exponential vectors to step function can be composed by taking products of pieces from exponential units $\ee_t(y):=\ee(\I_{\RO{0,t}}y)\in\DG_t(F)$ ($y\in F$).

For each $x$ such that $\ee(x)$ (and, consequently, also $\ee(-x)$) exists, we define the \phantomsection\hl{hyperbolic vectors}\index{hyperbolic vector}
% \baln{
% \ee_+(x)
% &
% ~:=~
% \frac{\ee(x)+\ee(-x)}{\sqrt{2}},
% &
% \ee_-(x)
% &
% ~:=~
% \frac{\ee(x)-\ee(-x)}{\sqrt{2}}.
% }\ealn
\baln{
\ee_+(x)
&
~:=~
\frac{\ee(x)+\ee(-x)}{2},
&
\ee_-(x)
&
~:=~
\frac{\ee(x)-\ee(-x)}{2i}.
}\ealn
Clearly, $\ee_+(x)$ and $\ee_-(x)$ are the components of $\ee(x)$ in the even part $\DG^+(F)$ and (modulo the factor $\frac{1}{i}$) in the odd part $\DG^-(F)$, respectively, of $\DG(F)$. Since the $\ee(x)$ are total in $\DG(F)$, the hyperbolic vectors $\ee_\pm(x)$ are total in the subspaces $\DG^\pm(F)$.

(One checks that  $\AB{\ee_+(x),\ee_+(y)}=\cosh\AB{x,y}$ and $\AB{\ee_-(x),\ee_-(y)}=\frac{1}{i}\sinh\AB{x,y}$, while, of course, $\AB{\ee_\pm(x),e_\mp\ee(y)}=0$; we do not need that.)

Let $\C^2$ be the diagonal subalgebra of $M_2$, and define the flip automorphism $\f\colon\rtMatrix{z\\w}\mapsto\rtMatrix{w\\z}$. When we consider $\C^2$ as  the identity \nbd{\C^2}correspondence we denote it by $\C^2_+$, while $\C^2_-$ denotes the \nbd{\C^2}correspondence obtained by equipping the Hilbert \nbd{\C^2}module $\C^2$ with the left action $b.b':=\f(b)b'$. (See also Skeide \cite{Ske01a}.) Obviously, $\C^2_-\odot\C^2_-$ can be identified with $\C^2_+$ via $\rtMatrix{z\\w}\odot\rtMatrix{z'\\w'}=\rtMatrix{wz'\\zw'}$, so that
\baln{
{\C^2_-}^{\odot(2n)}
&
~=~
\C^2_+,
&
{\C^2_-}^{\odot(2n+1)}
&
~=~
\C^2_-,
}\ealn
via $\rtMatrix{z_{2n}\\w_{2n}}\odot\ldots\odot\rtMatrix{z_1\\w_1}=\rtMatrix{w_{2n}z_{2n-1}\ldots w_2z_1\\z_{2n}w_{2n-1}\ldots z_2w_1}$, and analogously for the odd part. We get that
\beqn{
\DG(\C^2_-)
~=~
\DG(\C)\otimes\C^2
}\eeqn
as Hilbert \nbd{\C^2}module (and analogously for the Hilbert \nbd{\C^2}module structure of $\DG_t(\C^2_-)$), while as correspondence we get
\beqn{
\DG(\C^2_-)
~=~
(\DG^+(\C)\otimes\C^2_+)
~\oplus~
(\DG^-(\C)\otimes\C^2_-).
}\eeqn

Since $(xb)y_t=x(by_t)$, since $\DG(F)$ is \nbd{\C^2}spanned by $\ee(\vp)\otimes\U$, since $\DG_t(F)$ is \nbd{\C^2}spanned by $\ee_+(\vp)\otimes\U$ and $\ee_-(\vp)\otimes\U$, and since $\DG^\pm_t(\C^2_-)$ are also left \nbd{\C^2}spanned by $\ee_\pm(\vp)\otimes\U$, in the identification of $\DG(F)\odot\DG_t(F)$ with $\DG(F)$ it is enough to understand $(\ee(\vp)\otimes\U)(\ee_+(\psi_t)\otimes\U)=\ee(\vp)\ee_+(\psi_t)\otimes\U$ and $(\ee(\vp)\otimes\U)(\ee_-(\psi_t)\otimes\U)=\ee(\vp)\ee_-(\psi_t)\otimes\U$, where the products in the left factor are those of $\DG(\C)$ with $\DG_t(\C)$.

Note that $\sB^a(\DG(\C^2_-))=\rtMatrix{\sB(\DG(\C))\\\sB(\DG(\C))}$. Denoting by $\sS$ the CCR-flow on $\DG(\C)$ (while $\vt$ denotes the CCR-flow on $\DG(\C^2_-)$), we find
\baln{
\vt_t\sMatrix{a\\a'}(\ee(\vp)\otimes\U)(\ee_+(\psi_t)\otimes\U)
&
~=~
\sMatrix{\sS_t(a)\\\sS_t(a')}(\ee(\vp)\otimes\U)(\ee_+(\psi_t)\otimes\U),
\intertext{and}
\vt_t\sMatrix{a\\a'}(\ee(\vp)\otimes\U)(\ee_-(\psi_t)\otimes\U)
&
~=~
\sMatrix{\sS_t(a')\\\sS_t(a)}(\ee(\vp)\otimes\U)(\ee_-(\psi_t)\otimes\U).
}\ealn
Let $\Om$ denote the vacuum of $\DG(\C)$ (while $\om=\Om\otimes\U$ is the vacuum of $\DG(\C^2_-)$). Put
\beqn{
p
~:=~
\om\rtMatrix{1\\0}\om^*
~=~
%%%% BO
\rtMatrix{\Om\Om^*\\0}.
%Replaced \rtMatrix{\Om\Om^*\\0}\otimes\U, check if correct.
%%%% EO
}\eeqn
So,
\baln{
\vt_t(p)(\ee(\vp)\otimes\U)(\ee_+(\psi_t)\otimes\U)
&
~=~
\Om\ee_+(\psi_t)\otimes\sMatrix{1\\0},
\intertext{and}
\vt_t(p)(\ee(\vp)\otimes\U)(\ee_-(\psi_t)\otimes\U)
&
~=~
\Om\ee_-(\psi_t)\otimes\sMatrix{0\\1}.
}\ealn
Applying $p$ to both, taking also into account that $(\Om\Om^*)\ee_+(\vp)=\Om$ and $(\Om\Om^*)\ee_-(\vp)=0$, it follows that $\vt_t(p)\ge p$.

The projection $q:=\rtMatrix{\sid_{\DG(\C)}\\0}$ is the central projection in $\sB^a(\DG(\C^2_-))$ such that
\beqn{
q\sB^a(\DG(\C^2_-))
~=~
\cls^s\sB^a(\DG(\C^2_-))p\sB^a(\DG(\C^2_-))
~=~
\rtMatrix{\sB^a(\DG(\C))\\0}.
}\eeqn
(It coincides with multiplication from the right with the (central) element $\rtMatrix{1\\0}\in\C^2$.) We have $q(\ee(\vp)\ee_\pm(\psi_t)\otimes\U)=\ee(\vp)\ee_\pm(\psi_t)\otimes\rtMatrix{1\\0}$, but $\vt_t(q)(\ee(\vp)\ee_-(\psi_t)\otimes\U)=\ee(\vp)\ee_-(\psi_t)\otimes\rtMatrix{0\\1}$, so that $q$ is not increasing.

Denote the members of the superproduct system associated with $(\sB^a(\DG(\C^2)),\vt,p)$ by $E_t$. We find
\beqn{
E_t
~:=~
\vt_t(p)\sB^a(\DG(\C^2))p
~=~
\vt_t(p)\DG(\C^2)\rtMatrix{1\\0}\om^*
~=~
\sMatrix{\Om\DG_t^+(\C)\\0}\om^*.
}\eeqn
The left action of $\cB:=p\sB^a(\DG(\C^2))p=\C p\ni\lambda p$ is simply multiplication with $\lambda$. Effectively, this correspondence over $\cB\cong\C$ is just the Hilbert space $\DG_t^+(\C)$, the even part of $\DG_t(\C)$. The 
%%%% BO 
% tensor product 
product 
%%%% EO
$x_s\odot y_t\mapsto x_sy_t:=\vt_t(x_s)y_t$ is recovered from the above formulae as $\DG_s^+(\C)\otimes\DG_t^+(\C)\rightarrow\DG_{(s+t)}^+(\C)$ as restriction of that of $\DG_t(\C)$. Clearly, restricted to the even part, this is not onto the even part of $\DG_{s+t}(\C)$ as products of two odd factors are missing. So the $E_t$ form a superproduct system that is not a product system.

Summing up, the triple $(\sB^a(\DG(\C^2)),\vt,p)$ is a strong \nbd{E_0}dilation of the trivial Markov semigroup $T_t=\id_\C$ on $\C$ with an associated superproduct system that is proper.
\eex

Of course, restricting the preceding example to integer times, we get a discrete example. We prefer, however, to look at the general inherent structure, making it an independent example.

Before discussing the example we wish to clarify in more generality some things about discrete one-parameter \nbd{E_0}semigroups and their product systems.

\bob \label{corrE0ob}
Suppose we have a (full) Hilbert \nbd{\cB}module $E$, a \nbd{\cB}correspondence $F$, and a unitary $v\colon E\odot F\rightarrow E$. These ingredients give rise to several things. Firstly, and independently of $E$ and $v$, every correspondence $F$ gives rise to a discrete product system $F^\odot$ with  $F_n:=F^{\odot n}$ and the tensor product as product operation; by Observation \ref{d1pob}, up to isomorphism, all discrete one-parameter product systems arise in that way, so $F$ stands, in a sense, for the product system it generates. Secondly, $v$ gives rise to a unital (strict) endomorphism $\vt:=v(\bullet\odot\id_F)v^*$ of $\sB^a(E)$ and, therefore, to a whole \nbd{E_0}semigroup $\vt_n:=\vt^n$.

Of course, it appears to be clear that the product system of the semigroups $\vt_n$ is $F^\odot$. The easiest way to see this, is to establish a left dilation $v_n$ of $F^\odot$ to $E$ giving back $\vt^n$ as $\vt^n=v_n(\bullet\odot\id_{F^{\odot n}})v_n^*$. (Uniqueness for full $E$ identifies the product system up to isomorphism.) A left dilation that does this job, is
\beqn{
v_n
~:=~
v(v\odot\id_F)\ldots(v\odot\id_{F^{n-1}}),
}\eeqn
or, recursively, $v_0\colon E\odot\cB\rightarrow E$ the canonical map and $v_{n+1}=v(v_n\odot\id_F)$. (Since the product of $F^\odot$ is rebracketting, that is, \it{cum grano salis} the identity, and since by definition $v_n$ can be computed iteratively in one and only one way, $v_n$ is associative. And from the recursion it follows $v_0(\bullet\odot\id_{F^{\odot 0}})v_0^*=\id_E$ and
\beqn{
v_{n+1}(\bullet\odot\id_{F^{\odot {n+1}}})v_{n+1}^*
~=~
v(v_n\odot\id_F)(\bullet\odot\id_{F^{\odot n}}\odot\id_F)(v_n\odot\id_F)^*v^*
~=~
\vt (v_n (\bullet\odot\id_{F^{\odot n}})v_n^*),
}\eeqn
so that $v_n (\bullet\odot\id_{F^{\odot n}})v_n^*=\vt_n$ implies the same statement for $n+1$.)

Finding for a given correspondence $F$ suitable (full) $E$ and $v$, is more tricky. We make some comments in Remark \ref{corrE0rem}. There, we also discuss the case $\cB=\C$ which, plus some extra properties, is needed to make sure the following example is not nonempty.
\eob

\bex \label{discex}
For Hilbert spaces $H$ and $G$, suppose we have a unitary $V\colon H\otimes G\rightarrow H$. (Note: If $H\ne\zero$, so that $H$ is a full Hilbert \nbd{\C}module, this means $\dim G\ge1$. If, additionally, $\dim G\ge2$, this means $\dim H\ge\infty=\aleph_0$.) We assume the usual product notation $hg:=$ $V(h\otimes g)$. As in Observation \ref{corrE0ob}, we define the unital endomorphism $\Theta:=V(\bullet\otimes\id_G)V^*$ of $\sB(H)$. Then the \nbd{E_0}semigroup $\bfam{\Theta^n}$ has the product system $\bfam{G^{\otimes n}}$.

Define the full Hilbert \nbd{\C^2}module $E:=H\otimes\C^2$. Decompose $G$ (non-trivially) as $G=G_+\oplus\,G_-$ (so that $\dim G\ge2$) and put $F:=(G_+\otimes\C^2_+)\oplus(G_-\otimes\C^2_-)$. Define the unitary $v\colon E\odot F\rightarrow E$ as
\beqn{
\Bfam{h\otimes\sMatrix{z\\w}}\odot(g_+\otimes\U+g_-\otimes\U)
~\longmapsto~
hg_+\otimes\sMatrix{z\\w}+hg_-\otimes\sMatrix{w\\z}
}\eeqn
Again, as in Observation \ref{corrE0ob}, we define the unital endomorphism $\vt:=v(\bullet\otimes\id_F)v^*$ of $\sB^a(E)$. Then the \nbd{E_0}semigroup $\bfam{\vt^n}$ has product system $\bfam{F^{\odot n}}$. Note that
\beqn{
F^{\odot n}
~=~
\bigoplus_{\ve\in\CB{+,-}^n}\bfam{G_{\ve_1}\otimes\ldots\otimes G_{\ve_n}}\otimes\C^2_{\ve_1\ldots\ve_n},
}\eeqn
where we think of `$+=+1$' and `$-=-1$' so that the `product' $\ve_1\ldots\ve_n$ is ``$+$'' if there is an even number of ``$-$'' and ``$-$'' if there is an odd number of ``$-$'' in $\ve$.

Now assume that there are unit vectors $\Om\in H$ and $\Om_+\in G_+$ such that $\Om\Om_+=\Om$. (This means, the projection $\Om\Om^*$ is increasing for $\Theta$ and the vector state $\AB{\Om,\bullet\Om}$ on $\sB(H)$ is invariant for $\Theta$. If $\Om\Om^*$ increases to the identity on $H$, then the \nbd{E_0}semigroup $\Theta^n$ is said to be in \phantomsection\hl{standard form}\index{semigroup!E0@\nbd{E_0}!standard form}\index{E0-semigroup@\nbd{E_0}semigroup!standard form}.) Also $\om:=\Om\otimes\U$ and $\om_1:=\Om_+\otimes\U$ satisfy, $\om\om_1=\om$. So, $\om\om^*$ is increasing for $\vt$ and the vector expectation $\AB{\om,\bullet\om}$ on $\sB^a(E)$ is invariant for $\vt$. The \nbd{E_0}semigroup $\vt^n$ is in standard form, if $\Theta^n$ is.

As in Example \ref{hypexex}, define
\beqn{
p
~:=~
\om\sMatrix{1\\0}\om^*
~\in~
\sB^a(E)
}\eeqn
so that $pE=\Om\otimes\rtMatrix{\C\\0}=\om\rtMatrix{\C\\0}$. Then
\beqn{
\vt(p)pE
~=~
v(p\odot\id_F)v^*v(\om\odot\om_1)\sMatrix{\C\\0}
~=~
v((p\om)\odot\om_1)\sMatrix{\C\\0}
~=~
\Om\otimes\sMatrix{\C\\0},
}\eeqn
so $\vt(p)\ge p$. In other words, $(\sB^a(E),(\vt^n),p)$ is a (strong \nbd{E_0})dilation of the trivial Markov semigroup $T_n=\id_\C$ on $p\sB^a(E)p=p\C\cong \C$.

The superproduct system of $(\sB^a(E),(\vt^n),p)$ is
\bmun{
\sE_n
~:=~
\vt_n(p)\sB^a(E)p
~=~
v_n(p\odot\id_{F^{\odot n}})v_n^*E\sMatrix{1\\0}
~=~
v_n(p\odot\id_F)(E\odot F^{\odot n})\sMatrix{1\\0}
~=~
v_n(pE\odot F^{\odot n})\sMatrix{1\\0}
\\
~=~
\om\sMatrix{1\\0}F^{\odot n}\sMatrix{1\\0}
~\subset~
\om F^{\odot n}.
}\emun
The latter is the (super)product system of $(\sB^a(E),(\vt^n),\om\om^*)$, and we know it is a product system isomorphic to $F^\odot$. Effectively, is is easy to see that the isomorphism is just $\om y_n\mapsto y_n$. Since the product among the $\sE_n$ is (by definition!) the restriction of the product among the $\om F^{\odot n}$, the superproduct system is isomorphic to the subsystem $\sMatrix{1\\0}F^{\odot n}\sMatrix{1\\0}$ with tensor product $\odot$ as product operation. We find
\beqn{
\sMatrix{1\\0}F^{\odot n}\sMatrix{1\\0}
~=~
\bigoplus_{\substack{\ve\in\CB{+,-}^n\\\ve_1\ldots\ve_n=+}}\bfam{G_{\ve_1}\otimes\ldots\otimes G_{\ve_n}}\otimes\sMatrix{1\\0}
~\cong~
\bigoplus_{\substack{\ve\in\CB{+,-}^n\\\ve_1\ldots\ve_n=+}}\bfam{G_{\ve_1}\otimes\ldots\otimes G_{\ve_n}}
}\eeqn
as Hilbert subspace of $G^{\otimes n}$, and also the product is just the tensor product of subspaces of $G^{\otimes n}$. Clearly, this product is not surjective because it misses all elements in the \nbd{(m+n)}term that split into an \nbd{m}term and an \nbd{n}term each of which with an odd number of ``$-$'' in them.

Summarizing: $(\sB^a(E),(\vt^n),p)$ is a strong dilation of a Markov semigroup with associated a proper superproduct system. This dilation is strict; also, all modules are von Neumann modules, so that, being strict, this dilation is also normal. We see that, unlike the case that an \nbd{E_0}semigroup on $\sB^a(E)$ ($E$ a full Hilbert \nbd{\cB}module) has associated a product system (of \nbd{\cB}correspondences) if and only if the \nbd{E_0}semigroup is strict, the problem of a dilation  having associated  only a superproduct system is not (only) a topological  question. (Recall: The \nbd{E_0}semigroup on $\sB^a(E)$ has a product system (of \nbd{\C^2}correspondences), which is isomorphic to the (super)product system of the strong dilation $(\sB^a(E),(\vt^n),\om\om^*)$. It is the superproduct system (of \nbd{p\C}correspondences, that is, of Hilbert spaces) of the strong dilation $(\sB^a(E),(\vt^n),p)$ that is proper.)
\eex

\brem \label{corrE0rem}
To answer the question if, given $F$, there exist full $E$ and unitary $v\colon E\odot F\rightarrow E$, is easy to answer for Hilbert spaces and more tricky to answer for correspondences that are not Hilbert spaces. Since $E$ is (strongly) full, obviously, also $F$ has to be (strongly) full (in the von Neumann case). In Skeide \cite{Ske09} it is shown that for strongly full $F$, the answer is affirmative in the von Neumann case; it is affirmative in the \nbd{C^*}case, provided $F$ is full is over a \nbd{\sigma}unital \nbd{C^*}algebra. (In either case, separability of $F$ or its pre-dual can be preserved.)

For Hilbert spaces $G$, the answer is easy: Just take any unit vector $\Om_1\in G$ and input the pair $((G^{\otimes n}),(\Om_1^{\otimes n}))$ into the inductive limit construction in Theorem \ref{Oreindthm}. This provides us not only with a (nonzero, that is, full) Hilbert space $H$ and a unitary $V\colon H\otimes G\rightarrow H$, but also with a unit vector $\Om$ such that $\Om\Om_1=\Om$. This shows that Example \ref{discex} is not working on the empty set, and concludes this example. The construction is well-known; just that it is known since long, rather with inverted order in tensor products (see also Appendix \ref{vNAPP}\ref{vNcomm}).

For Hilbert spaces, there is not really a difference between finding a unitary $V\colon H\otimes G\rightarrow H$ or a unitary $W\colon G\otimes H\rightarrow H$. The difference is that the latter defines a \it{representation} of the correspondence $G$, while the former defines an anti-representation. For separable $G$, identifying the elements of an ONB of $G$ with the generators of the Cuntz algebra $\cO_{\mathsf{dim}\,G}$, we get a(n anti-)representation of  $\cO_{\mathsf{dim}\,G}$. This well-known relation between representations of Cuntz algebras and unital endomorphisms of $\sB(H)$ generalizes to more general correspondences.

Faithful (\it{nondegenerate} or \it{essential}) representations of a \nbd{\cB}correspondence $F$ on Hilbert space $H$ correspond (one-to-one!) with pairs consisting of a faithful \nbd{\cB}\nbd{\C}correspondence $H$ (that is, a Hilbert space with a faithful nondegenerate representation of $\cB$) and a left linear unitary $w\colon F\odot H\rightarrow H$. For finding faithful $H$, necessarily $F$ has to be faithful, too. But this is all that is needed: Hirshberg \cite{Hir05a} has shown existence of faithful nondegenerate representations for faithful \nbd{C^*}correspondences that are also full. \cite{Ske09} has shown the von Neumann case and the \nbd{C^*}case in full generality, that is, for all faithful von Neumann and \nbd{C^*}correspondences.
\erem

\newpage

\section{Product systems over products}\label{compSEC}

The scope of this section is to understand how product systems over a product $\bS=\bS^1\times\ldots\times \bS^d$ of monoids $\bS^k$ are made up out of their $d$ \it{marginal} product systems over $\bS^k$ and how, conversely, product systems over $\bS^k$ can be put together to form a product system over the product $\bS$. (Our special interest is to understand \nbd{d}parameter product systems out of one-parameter product systems, $\bS=\N_0^d$ or $\bS=\R_+^d$. Since for the construction of one-parameter product systems, we have the powerful Theorem \ref{indlimthm}, it is a promising strategy to understand how to put together $d$ one-parameter product systems to a \nbd{d}parameter product system.) It is easy to find a subset of the set of product maps of the product system over $\bS$ that determine the whole structure and necessary conditions that this subset has to satisfy. Showing that these are sufficient for the converse direction, is, however, tricky and involves a detailed analysis of a subset of permutations, the \it{order improving, partially order preserving} permutations (to be dealt with in Appendix \ref{popAPP}).

We start by discussing the case $d=2$ for Hilbert spaces passing, then, to modules; not so much to see what \bf{is} possible for Hilbert spaces, but rather to explain what is \bf{not} possible for modules (and, therefore, has to be avoided) and how it can be replaced by a generalizable version.

Let ${E^1}^\otimes$ and ${E^2}^\otimes$ be product systems of Hilbert spaces over $\bS^1$ and $\bS^2$, respectively. We easily verify that the family $\bfam{E^1_{t_1}\otimes E^2_{t_2}}_{(t_1,t_2)\in\bS^1\times\bS^2}$ with product $(x^1_{s_1}\otimes x^2_{s_2})(y^1_{t_1}\otimes y^2_{t_2})=(x^1_{s_1}y^1_{t_1})\otimes(x^2_{s_2}y^2_{t_2})$ is a product system over $\bS^1\times\bS^2$. The product involves the natural \phantomsection\hl{flip}\index{flip} isomorphism $E^2_{s_2}\otimes E^1_{t_1}\rightarrow E^1_{t_1}\otimes E^2_{s_2}$. Since there is no such flip for tensor products of correspondences, there is no such \it{external tensor product} of product systems of correspondences.%
\footnote{ \label{exttenFN}
Well, we may form a \phantomsection\it{truly external tensor product}\index{product system!external tensor product of} of product systems of correspondences\index{tensor product!external}\index{correspondence!external tensor product} even over different algebras. But then, the result would be a product system of correspondences over the (or, better, some) tensor product of the algebras. We wish, however, to stay inside the category of \nbd{\cB}correspondences for a fixed $\cB$, so, we do not want such a \it{truly external tensor product}.
}

% \lf
\brem \label{HSTPrem}
Note that there is also the tensor product within the category of product systems of Hilbert spaces over a fixed monoid $\bS$ with product $(x^1_s\otimes x^2_s)(y^1_t\otimes y^2_t)=(x^1_sy^1_t)\otimes(x^2_sy^2_t)$. Also here, for correspondences there is no such tensor product. (See, however, again Footnote \ref{exttenFN}.) For the category of \phantomsection\hl{spatial}\index{product system!spatial}\index{spatial!product system} tensor product systems over $\R_+$ (that is, there is a central unital unit), Skeide \cite{Ske06d} has replaced it by the spatial product. (As pointed out by Bhat and Mukherjee \cite{BhMu10}, the construction requires to embed a subproduct system into a product system as in Theorem \ref{indlimthm}. This, therefore, does not work for arbitrary monoids.) This is not what we are interested in. Our interest is in constructing multi-parameter product systems out of one-parameter product systems.
\erem

On the other hand, suppose $E^\odot$ is a product system over $\bS^1\times\bS^2$. Then the \phantomsection\hl{marginal}\index{marginal product systems of a product system over a product} families ${E^1}^\odot:=\bfam{E_{(t_1,0)}}_{t_1\in\bS^1}$ and ${E^2}^\odot:=\bfam{E_{(0,t_2)}}_{t_2\in\bS^2}$ are product systems over $\bS^1$ and $\bS^2$, respectively. Moreover, $E^1_{t_1}\odot E^2_{t_2}$ and $E_{(t_1,t_2)}$ are isomorphic via $u_{(t_1,0),(0,t_2)}$. Therefore, the product system structure of $E^\odot$ may be rewritten in terms of the family $\bfam{E^1_{t_1}\odot E^2_{t_2}}_{(t_1,t_2)\in\bS^1\times\bS^2}$. The induced product maps for this family, made to make the family $u_{(t_1,0),(0,t_2)}$ into an isomorphism of product systems, are
\vspace{-1ex}
\beq{ \label{1o2prod}
(u^1_{s_1,t_1}\odot u^2_{s_2,t_2})(\id_{E^1_{s_1}}\odot(u_{(t_1,0),(0,s_2)}^*u_{(0,s_2),(t_1,0)})\odot\id_{E^2_{t_2}}).
}\eeq
(Effectively, applying to this map $u_{(s_1t_1,0),(0,s_2t_2)}$, taking also into account that
\vspace{-1ex}
\beqn{
(u_{(s_1t_1,0),(0,s_2t_2)})(u^1_{s_1,t_1}\odot u^2_{s_2,t_2})
~=~
u_{(s_1,0),(t_1,0),(0,s_2),(0,t_2)}
~=~
u_{(s_1,0),(t_1,s_2),(0,t_2)}(\id_{E^1_{s_1}}\odot u_{(t_1,0),(0,s_2)}\odot\id_{E^2_{t_2}})
}\eeqn
(where we use the $n$th iterated product notation from Observation \ref{d1pob}), we get
\vspace{-1ex}
\bmun{
(u_{(s_1t_1,0),(0,s_2t_2)})(u^1_{s_1,t_1}\odot u^2_{s_2,t_2})(\id_{E^1_{s_1}}\odot(u_{(t_1,0),(0,s_2)}^*u_{(0,s_2),(t_1,0)})\odot\id_{E^2_{t_2}})
\\
~=~
u_{(s_1,0),(t_1,s_2),(0,t_2)}(\id_{E^1_{s_1}}\odot u_{(0,s_2),(t_1,0)}\odot\id_{E^2_{t_2}})
~=~
u_{(s_1,0),(0,s_2),(t_1,0),(0,t_2)}
\\
~=~
u_{(s_1,s_2),(t_1,t_2)}(u_{(s_1,0),(0,s_2)}\odot u_{(t_1,0),(0,t_2)}),
}\emun
that is, the family $u_{(t_1,0),(0,t_2)}$, indeed, intertwines the product of the families $\bfam{E^1_{t_1}\odot E^2_{t_2}}_{(t_1,t_2)\in\bS^1\times\bS^2}$ and $E^\odot$, also showing that \eqref{1o2prod} is, indeed, a product.) We see that in the product in \eqref{1o2prod}, apart from the structure maps $u^k_{s_k,t_k}$ of the marginals ${E^k}^\odot$\!, there occurs a family of isomorphisms $u^{1,2}_{t_1,s_2}:=u_{(t_1,0),(0,s_2)}^*u_{(0,s_2),(t_1,0)}\colon E^2_{s_2}\odot E^1_{t_1}\rightarrow E^1_{t_1}\odot E^2_{s_2}$.

\bex \label{HSprodex}
For instance, if $E^\otimes$ is just the external tensor product of two product systems of Hilbert spaces as described above, then $u^{1,2}_{t_1,s_2}$ is nothing but the flip. But even for Hilbert spaces the structure of a product system over $\bS^1\times\bS^2$ is not determined by the structure of the marginals. It particular, it need not be isomorphic to an external tensor product. In fact, in Corollary \ref{N0disocor} we determine the structure of all discrete \nbd{d}parameter product systems $E^\odot$. In particular, for a correspondence $E$, every bilinear unitary $\sF_{1,2}\in\sB^{bil,a}(E\odot E)$ (thought of as an operator $E_2\odot E_1\rightarrow E_1\odot E_2$) gives rise to a product system structure on $\bfam{E^{\odot n_1}\odot E^{\odot n_2}}_{(n_1,n_2)\in\N_0^2}$. (Every product system $E^\odot$ over $\N_0^2$ with (pairwise) isomorphic marginals is isomorphic to one of that form, but we do not need that here.) Two such operators $\sF_{1,2}$ and $\sF'_{1,2}$ give rise to isomorphic product systems if and only if there exist unitaries $a_k\in\sB^{bil,a}(E)$ (thought of as operators $E_k\rightarrow E_k$)) such that $(a_1\odot a_2)\sF_{1,2}=\sF'_{1,2}(a_2\odot a_1)$. It is clear that, for instance, $\sF_{1,2}=\id$ and $\sF'_{1,2}=\f$ (the flip for some Hilbert space $H$ with $\dim H\ge2$) violate this condition.
\eex

To see more clearly what we aim at for general $d$, let us repeat, also slightly reformulating, what we achieved in the case $d=2$, and arrange it into the following three steps.
\begin{enumerate}
\item \label{step1}
We started with a product system $E^\odot=\bfam{E_\bt}_{\bt\in\bS}$ over $\bS=\bS^1\times\bS^2$, and replaced it by a family written in terms of the marginals ${E^k}^\odot$ of the form $\bfam{E_{(t_1,0)}\odot E_{(0,t_2)}}_{\bt\in\bS}=\bfam{E^1_{t_1}\odot E^2_{t_2}}_{\bt\in\bS}$, which is ``pointwise'' isomorphic to the family $E^\odot$ via $u_{(t_1,0),(0,t_2)}$.

\item \label{step2}
We noted that the product system structure of $E^\odot$ can also be expressed in terms of the family $\bfam{E_{t_1}\odot E_{t_2}}_{\bt\in\bS}$ equipping it with the product
\vspace{-1ex}
\beqn{
u_{(s_1t_1,0),(0,s_2t_2)}^*u_{\bs,\bt}(u_{(s_1,0),(0,s_2)}\odot u_{(t_1,0),(0,t_2)})
}\eeqn
obtained by unitary equivalence from the product $u_{\bs,\bt}$.

\item \label{step3}
We convinced ourselves that this product can be written as
\vspace{-1ex}
\beqn{
(u^1_{s_1,t_1}\odot u^2_{s_2,t_2})(\id_{E^1_{s_1}}\odot u^{1,2}_{t_1,s_2}\odot\id_{E^2_{t_2}})
}\eeqn
in a way where only the operators $u^{1,2}_{t_1,s_2}:=u_{(t_1,0),(0,s_2)}^*u_{(0,s_2),(t_1,0)}\colon E^2_{s_2}\odot E^1_{t_1}\rightarrow E^1_{t_1}\odot E^2_{s_2}$ and the products $u^k_{s_k,t_k}$ of the marginals occur.
\end{enumerate}
Carrying out Steps \ref{step1} and \ref{step2} for general $d$ obviously generalizes without any question: Starting from a product system $E^\odot=\bfam{E_\bt}_{\bt\in\bS}$ over $\bS=\bS^1\times\ldots\times\bS^d$ we may pass to the family $\bfam{E^1_{t_1}\odot\ldots\odot E^d_{t_d}}_{\bt\in\bS}$ consisting of the tensor products of the \phantomsection\hl{marginal systems}\index{marginal product systems of a product system over a product|bf} ${E^k}^\odot:=\bfam{E^k_{t_k}:=E_{(0,\ldots,0,t_k,0,\ldots,0)}}_{t_k\in\bS^k}$, which is ``pointwise'' isomorphic to the family $E^\odot$ via $u_{(t_1,0,\ldots,0),\,\ldots\,,(0,\ldots,0,t_d)}$. And the product $u_{\bs,\bt}$ of $E^\odot$ lifts to the unitarily equivalent product
\vspace{-1ex}
\beqn{
u_{(s_1t_1,0,\ldots,0),\,\ldots\,,(0,\ldots,0,s_dt_d)}^*u_{\bs,\bt}(u_{(s_1,0,\ldots,0),\,\ldots\,,(0,\ldots,0,s_d)}\odot u_{(t_1,0,\ldots,0),\,\ldots\,,(0,\ldots,0,t_d)})
\vspace{-1ex}
}\eeqn
on the other family.

What about Step \ref{step3}? Or better, what is it we gained in the case $d=2$ by expressing the product according to Step \ref{step2} in the form in Step \ref{step3}, that makes it desirable to do the same for general $d$? Answer: The product as in Step \ref{step3} expresses in a very concise form what, apart from the structure of the marginals it needs to compute the product, namely: A flip operation that, in $(E^1_{s_1}\odot E^2_{s_2})\odot(E^1_{t_1}\odot E^2_{t_2})$ brings things into ``the right'' order $(E^1_{s_1}\odot E^1_{t_1})\odot(E^2_{s_2}\odot E^2_{t_2})$ so that, then, the product of the marginals can be applied. In fact, this is very much how we defined the external tensor product in the case of Hilbert spaces with the usual flip; just that, for correspondences, the flip does, in general, not exist; and even if, for Hilbert spaces, the flip does exist, Example \ref{HSprodex} shows that the flip is by far not the only possible way. Note, too, the gradual change from the scope of conveniently describing the structure of a given product system over $\bS^1\times\bS^2$ in terms of its marginals (and some extras, the flips $u^{1,2}_{t_1,s_2}$), to rather starting from product systems ${E^k}^\odot$ over $\bS^k$ and to define on the family of tensor products a product extending the marginal products (with the help of some flips). Let us  transfer this to general $d$.

So, suppose we have $d$ product systems ${E^k}^\odot$ $(k=1\ldots,d)$, each over its own $\bS^k$, with product maps $u^k_{s_k,t_k}$. For each $\bt=(t_1,\ldots,t_d)\in\bS^1\times\ldots\times\bS^d=:\bS$, we define
\vspace{-1ex}
\beqn{
E_\bt
~:=~
E^1_{t_1}\odot\ldots\odot E^d_{t_d}.
\vspace{-1ex}
}\eeqn
Henceforth, to lighten notation, we will identify $\bS^k$ with the submonoid $(0,\ldots,0,\bS^k,0,\ldots,0)$ of $\bS$ and, therefore, $t_k\in\bS^k$ with $(0,\ldots,0,t_k,0,\ldots,0)\in\bS$. We also identify $E_{t_k}=E^1_0\odot\ldots\odot E^{k-1}_0\odot E^k_{t_k}\odot E^{k+1}_0\odot\ldots\odot E^d_0$ with $E^k_{t_k}$, so that the marginal families are identified with ${E^k}^\odot$. Our intention is to define a product system structure $u_{\bs,\bt}$ on the family $E^{1,\ldots,d}:=\bfam{E_\bt}_{\bt\in\bS}$ that extends the product of the marginals. To that goal, we wish to send
\begin{subequations} \label{b-a}
\vspace{-.5ex}
\beq{ \phantomsection\label{before}
E_\bs\odot E_\bt
~=~
(E^1_{s_1}\odot\ldots\odot E^d_{s_d})\odot(E^1_{t_1}\odot\ldots\odot E^d_{t_d}),
\vspace{-1ex}
}\eeq
from which the arguments of $u_{\bs,\bt}$ stem, by applying successively suitable next-neighbour flips into
\vspace{-1.5ex}
\beq{ \phantomsection\label{after}
(E^1_{s_1}\odot E^1_{t_1})\odot\ldots\odot(E^d_{s_d}\odot E^d_{t_d}),
}\eeq
\end{subequations}
to which, then, we can apply the marginal products as $u^1_{s_1,t_1}\odot\ldots\odot u^d_{s_d,t_d}$. At first sight, the flips that might allow us to succeed with our purpose, might even depend not only on the pair $E^i_{s_i}\odot E^j_{t_j}$ they have to flip, but on the entire context (that is, on all indices and ``times'' of the surrounding tensor factors). However, if we assume for a moment that a product $u_{\bs,\bt}$ does exist, we know that we can achieve the flip by the operators $u_{t_j,s_i}^*u_{s_i,t_j}$. This special flip is also compatible with the product map in the sense of the following little lemma, which, to free it from too many indices, we formulate generally.

\blem\label{switchlem}
Let $E^\odot$ be a product system over a monoid $\bS$. Choose $n\in\N$ and $t_k\in\bS$ $(1\le k\le n)$. Suppose that for some $i$ we have $t_it_{i+1}=t_{i+1}t_i$. Then the $n$th iterated product is \hl{invariant} under the flip $u_{t_{i+1},t_i}^*u_{t_i,t_{i+1}}$, that is,
\vspace{-1ex}
\beqn{
u_{t_1,\ldots,t_n}
~=~
u_{t_1,\ldots,t_{i-1},t_{i+1},t_i,t_{i+2},\ldots,t_n}(\id_{t_1}\odot\ldots\odot\id_{t_{i-1}}\odot~u_{t_{i+1},t_i}^*u_{t_i,t_{i+1}}\odot\id_{t_{i+2}}\odot\ldots\odot\id_{t_n}).
}\eeqn
\elem

\vspace{.5ex}
\noindent
(The proof is very similar to the computation after \eqref{1o2prod}, and we omit it.) So, in any product system, assuming also that $\bS$ is abelian, we may compute the $n$th iterated product also after ``flipping around'' the spaces of the $n$ factors (though, not their elements, as the flip does, in general, not exist) with the help of the operators $u_{t_{i+1},t_i}^*u_{t_i,t_{i+1}}$ before without changing the result. Returning to our family $E^{1,\ldots,d}$ this means that, if there does exist a product, then we may compute it by, first, bringing the form in \eqref{before} into the form in \eqref{after} and, then, applying the marginal products. To achieve this reordering, we may iterate any selection of flips $u_{t_j,s_i}^*u_{s_i,t_j}$ that compose to give a suitable permutation to end up in \eqref{after}. Well, actually, we need only to bring spaces of higher index $i$ on the left of a space with lower index $j$ into the right order. We might flip also if $j=i$ (Lemma \ref{switchlem} shows, it does not change the result), but it is not necessary to do so. So, we need the flips only for $j<i$. (It will turn out that the suitable permutation that does the job is, then, unique.) We summarize:

\bob \label{prodob}
Suppose we have a family of bilinear unitaries $u^{j,i}_{t_j,s_i}\colon E^i_{s_i}\odot E^j_{t_j}\rightarrow E^j_{t_j}\odot E^i_{s_i}$ ($1\le j<i\le d$, $s_i\in\bS^i$, $t_j\in\bS^j$). Then:
\begin{enumerate}
\item \label{prod1}
The prescription indicated around Equations \eqref{b-a} over-determines maps $u_{\bs,\bt}$ on the family $E^{1,\ldots,d}$; `over-determines' in the sense that we will definitely be able to compose a suitable permutation out of flips in the prescribed way (to which we, then, apply the marginal products), but it is not \it{a priori} clear if this is well-defined.

\item \label{prod2}
If so, and if the $u_{\bs,\bt}$ do define a product (associativity!) turning the family $E^{1,\ldots,d}$ into a product system ${E^{1,\ldots,d}}^\odot$, then necessarily $u_{t_j,s_i}^*u_{s_i,t_j}=u^{j,i}_{t_j,s_i}$.

\item \label{prod3}
By the Steps \eqref{step1} -- \eqref{step3} (following Example \ref{HSprodex}), every product system over $\bS$ is isomorphic to some ${E^{1,\ldots,d}}^\odot$.

\item \label{prod4}
If the $u^{j,i}_{t_j,s_i}$ define a product, then necessarily
\vspace{-.5ex}
\beqn{
(\id^k_{t_k}\odot~u^{j,i}_{s_j,r_i})(u^{k,i}_{t_k,r_i}\odot\id^j_{s_j})(\id^i_{r_i}\odot~u^{k,j}_{t_k,s_j})
~=~
(u^{k,j}_{t_k,s_j}\odot\id^i_{r_i})(\id^j_{s_j}\odot~u^{k,i}_{t_k,r_i})(u^{j,i}_{s_j,r_i}\odot\id^k_{t_k})
}\eeqn
for each $k<j<i$ and $r_i\in\bS^i,s_j\in\bS^j,t_k\in\bS^k$. (This is so, because, by Lemma \ref{switchlem}, it does not matter which of the two possibilities to realize the permutation $E^i_{r_i}\odot E^j_{s_j}\odot E^k_{t_k}\rightarrow E^k_{t_k}\odot E^j_{s_j}\odot E^i_{r_i}$ we choose, before applying the triple product $u_{t_k,s_j,r_i}$.)

\item \label{prod5}
If the $u^{j,i}_{t_j,s_i}$ define a product, then necessarily
\vspace{-.5ex}
\baln{
(u^j_{s_j,t_j}\odot\id^i_{r_i})(\id^j_{s_j}\odot~u^{j,i}_{t_j,r_i})(u^{j,i}_{s_j,r_i}\odot\id^j_{t_j})
&
~=~
u^{j,i}_{s_jt_j,r_i}(\id^i_{r_i}\odot~u^j_{s_j,t_j}),
\\
(\id^j_{t_j}\odot~u^i_{r_i,s_i})(u^{j,i}_{t_j,r_i}\odot\id^i_{s_i})(\id^i_{r_i}\odot~u^{j,i}_{t_j,s_i})
&
~=~
u^{j,i}_{t_j,r_is_i}(u^i_{r_i,s_i}\odot\id^j_{t_j})
}\ealn
for $j<i$ and $s_j,t_j\in\bS^j;r_i,s_i\in\bS^i$. (This means it does
%%%% BO
not
%%%% EO 
matter if we, first, bring the three sites into increasing order and, then, multiply the two sites from the same marginal, or if we, first, multiply  the two sites from the same marginal and, then, bring the (remaining two) sites into increasing order. Again, the statement follows by suitable application of Lemma \ref{switchlem}.)
\end{enumerate}
\eob
This observation means: We capture any product system structure on $E^{1,\ldots,d}$ that extends the marginal products by appropriately using the flips $u^{j,i}_{t_j,s_i}$ ($j<i$). We capture any product system $E^\odot$ over the product $\bS$ up to isomorphism by its marginals and the flips. In either case, necessarily, the flips are recovered from the product as $u^{j,i}_{t_j,s_i}=u_{t_j,s_i}^*u_{s_i,t_j}$. Moreover, they fulfill the two sets of necessary conditions stated in \ref{prod4} and \ref{prod5}. We now show that these two sets of conditions are also sufficient.

The key ingredients for the proof are to be found in the following definition and lemma. The proof of the lemma is postponed to Appendix \ref{popAPP}. But first let us fix some notation. For $n\in\N$ we shall denote \phantomsection\index{Nn@$\N_n$}$\N_n:=\CB{1,\ldots,n}$. For each $n$, the permutation group $S_n$\index{permutation} is the set of all bijections $\sigma$ on $\N_n$. For $j\ne i$, we have the usual \hl{transposition}\index{transposition} of $j$ and $i$ as
\vspace{-.5ex}
\beqn{
\tau_{j,i}(\ell)
~=~
\begin{cases}
i&\ell=j,
\\
j&\ell=i,
\\
\ell&\text{otherwise}.
\end{cases}
}\eeqn
Obviously, $\tau_{j,i}=\tau_{i,j}=\tau_{j,i}^{-1}$. It is well known that the \hl{next-neighbour transpositions}\index{transposition!next-neighbour} $\tau_k:=\tau_{k,k+1}$ $(k\in\N_{n-1})$ generate $S_n$ (as a semigroup), that they fulfill the relations
\bal{\label{Snrel}
\tau_k^2
&
=
e,
&
~~~~~~
\tau_k\tau_{k'}
&
=
\tau_{k'}\tau_k
~~~
(\abs{k-k'}\ge2),
~~~~~~
&
\tau_k\tau_{k+1}\tau_k
&
=
\tau_{k+1}\tau_k\tau_{k+1}
~~~
(k\le n-2),~~~~~~
}\eal
and that every set of $n-1$ elements $\tau'_k$ in a monoid fulfilling these relation generates a group isomorphic to $S_n$ (via $\tau'_k\mapsto\tau_k$).

\bdefi \label{admisdef}
Let $E_1,\ldots,E_p$ be correspondences over $\cB$ and for each $1\le j<i\le p$ let $\sF_{j,i}$ be a bilinear unitary $E_i\odot E_j\rightarrow E_j\odot E_i$. Let $f\colon\N_q\rightarrow\N_p$ be a function and put
\beqn{
E_f
~:=~
E_{f(1)}\odot\ldots\odot E_{f(q)}.
}\eeqn

For $j<i$ and $\vk\in\N_{q-1}$, we say the pair $(\sF_{j,i},\vk)$ is \phantomsection\hl{admissible}\index{admissible!pair} with respect to $f$ (via $\vk$) if $f(\vk)=i$ and $f(\vk+1)=j$, that is, if in  the two neighbouring sites $E_{f(\vk)}$ and $E_{f(\vk+1)}$ in $E_f$ we have the tensor product $E_{f(\vk)}\odot E_{f(\vk+1)}=E_i\odot E_j$. If $(\sF_{j,i},\vk)$ is admissible we denote by $\sF_{j,i;\vk}^{\odot\ssid}$ the amplification of $\sF_{i,j}$ acting on the sites $E_{f(\vk)}\odot E_{f(\vk+1)}$ of $E_{f(1)}\odot\ldots\odot E_{f(q)}$, namely,
\beqn{
\sF_{j,i;\vk}^{\odot\ssid}
~:=~
\id_{f(1)}\odot\ldots\odot\id_{f(\vk-1)}\odot\,\,\sF_{j,i}\odot\id_{f(\vk+2)}\odot\ldots\odot\id_{f(q)}.
}\eeqn
(Note that the codomain of $\sF_{j,i;\vk}^{\odot\ssid}$ is $E_{f\circ\,\tau_\vk}$.) A chain $(\sF_{j_1,i_1},\vk_1),\ldots,(\sF_{j_m,i_m},\vk_m)$ is \phantomsection\hl{admissible}\index{admissible!chain (of pairs)} with respect to $f$ (via $\vk_1,\ldots,\vk_m$) if  each $(\sF_{j_k,i_k},\vk_k)$ is admissible with respect to $f\circ\tau_{\vk_1}\circ\ldots\circ\tau_{\vk_{k-1}}$. (This means that domains and codomains of the operators in $\sF_{j_m,i_m;\vk_m}^{\odot\ssid}\circ\ldots\circ\sF_{j_1,i_1;\vk_1}^{\odot\ssid}$ match. Note, too, that the codomain of the whole thing is $E_{f\circ\,\tau_{\vk_1}\circ\ldots\circ\,\tau_{\vk_m}}$.)

An admissible chain $(\sF_{j_1,i_1},\vk_1),\ldots,(\sF_{j_m,i_m},\vk_m)$ is \phantomsection\hl{maximal}\index{admissible!chain (of pairs)!maximal} if there is no $(\sF_{j_{m+1},i_{m+1}},\vk_{m+1})$ such that the chain $(\sF_{j_1,i_1},\vk_1),\ldots,(\sF_{j_{m+1},i_{m+1}},\vk_{m+1})$ is still admissible.
\edefi

\brem
It is clear that admissibility of $(\sF_{j,i},\vk)$ or of a chain $(\sF_{j_1,i_1},\vk_1),\ldots,(\sF_{j_m,i_m},\vk_m)$ is a question that depends exclusively on the function $f$ and the indices in question, and not on the concrete nature of the isomorphisms $\sF_{j,i}$. Nevertheless, for the formulations that follow, we found it more convenient to speak of admissible pairs of operators $\sF_{j,i}$ and indices $\vk$ rather than speaking of admissible triples of indices $(j,i;\vk)$. See, however, Appendix \ref{popAPP}.
\erem

\blem \label{pi_flem}
In the situation of Definition \ref{admisdef}:
\begin{enumerate}
\item
For each $f$, there do exist maximal admissible chains $(\sF_{j_1,i_1},\vk_1),\ldots,(\sF_{j_m,i_m},\vk_m)$.

\item
The permutation $\tau_{\vk_1}\circ\ldots\circ\,\tau_{\vk_m}$ for a maximal admissible chain does not depend on the choice of the chain. We denote that permutation by $\sigma_f$.

\item
If the $\sF_{j,i}$ fulfill the \phantomsection\hl{detailed exchange conditions}\index{detailed exchange conditions|bf}\index{exchange operator!detailed exchange conditions}
\beq{\label{Tijcond}
(\id_k\odot \sF_{j,i})(\sF_{k,i}\odot\id_j)(\id_i\odot \sF_{k,j})
~=~
(\sF_{k,j}\odot\id_i)(\id_j\odot \sF_{k,i})(\sF_{j,i}\odot\id_k)
}\eeq
$(1\le k<j<i\le p)$, then the operator $\sF_{j_m,i_m}^{\odot\ssid}\circ\ldots\circ\sF_{j_1,i_1}^{\odot\ssid}$ for a maximal admissible chain does not depend on the choice of the chain. We denote this operator by $\pi_f$.
\end{enumerate}
\elem

\noindent
Note that the function $f\colon\N_q\rightarrow\N_p$ can also be viewed as a \nbd{q}tuple $(f(1),\ldots,f(q))$ with values in $\N_p$. The tuple that belongs to $f\circ\sigma_f$ is increasing. We insisted that the permutation $\sigma_f$ is constructed from next-neighbour transpositions, but using only those transpositions that are really necessary in order to bring a pair of neighbouring elements in the tuple into increasing order. Not only transpositions that would increase disorder are not allowed, but also transpositions that would flip positions that have the same value of $f$ are forbidden. This means we change two places if the indices are not in order, but the order of places where the indices coincide is never changed. The former property means that we use only flips that are \phantomsection\hl{order improving}\index{order improving}\index{permutation!order improving}; the latter property is what we will call in Appendix \ref{popAPP} \hl{partially order preserving}\index{partially order preserving}\index{permutation!partially order preserving}; and these properties are responsible for that the permutation $\sigma_f$ is unique. All this will be made precise in Appendix \ref{popAPP}, to which we delegate the proof of Lemma \ref{pi_flem}.

We are now ready to complete the proof of the following theorem, showing also a reverse of Observation \ref{prodob}\eqref{prod4} and\eqref{prod5}.

\bthm \label{prodthm}
In the situation of Observation \ref{prodob}: \index{product system!over a product, the structure of}\index{structure!of a product system over a product}The formula $u^{j,i}_{t_j,s_i}=u_{t_j,s_i}^*u_{s_i,t_j}$ establishes a one-to-one correspondence between
\begin{enumerate}
\item
product systems structures $u_{\bs,\bt}$ on $E^{1,\ldots,d}$ that extend the marginal products, and 

\item
families of operators $u^{j,i}_{t_j,s_i}\colon E^i_{s_i}\odot E^j_{t_j}\rightarrow E^j_{t_j}\odot  E^i_{s_i}$ $(1\le j<i\le i,t_j\in\bS^j,s_i\in\bS^i)$ that satisfy the conditions in Observation \ref{prodob}\eqref{prod4} and \eqref{prod5}.
\end{enumerate}
\ethm

\proof
The forward implication is contained in Observation \ref{prodob}. So, let us assume we have a family of operators $u^{j,i}_{t_j,s_i}\colon E^i_{s_i}\odot E^j_{t_j}\rightarrow E^j_{t_j}\odot  E^i_{s_i}$ $(1\le j<i\le i,t_j\in\bS^j,s_i\in\bS^i)$ satisfying the conditions in \eqref{prod4} and\eqref{prod5}.

For showing that $u_{\bs,\bt}$ is well defined, we put $E_k:=E^k_{s_k}\oplus E^k_{t_k}$ so that \eqref{before} is contained in
\begin{subequations} \label{2b-a}
\beq{ \phantomsection\label{2before}
E_1\odot\ldots\odot E_d\odot E_1\odot\ldots\odot E_d,
}\eeq
while \eqref{after} is contained in
\beq{ \phantomsection\label{2after}
E_1\odot E_1\odot\ldots\odot E_d\odot E_d.
}\eeq
\end{subequations}
Further, for $j<i$ we define $\sF_{j,i}\colon E_i\odot E_j\rightarrow E_j\odot E_i$ to act direct-summand-wise as the appropriate $u^{j,i}_{\beta_j,\alpha_i}$ ($\beta_j\in\CB{t_j,s_j}$, $\alpha_i\in\CB{s_i,t_i}$). It is clear that verifying that the $\sF_{j,i}$ fulfill \eqref{Tijcond}, amounts to just verifying the (by hypotheses, satisfied) conditions of the $u^{j,i}_{\bullet,\bullet}$ for many different ``time'' arguments in Observation \ref{prodob}\eqref{prod4}. So, according to Lemma \ref{pi_flem}, for the function $f\colon\N_{2d}\rightarrow\N_d$ such that $(f(1),\ldots,f(2d))=(1,\ldots,d,1,\ldots,d)$ there is a unique permutation $\sigma_f$ putting the tuple into $(f\circ\sigma(1),\ldots,f\circ\sigma(2d))=(1,1,\ldots,d,d)$ that can be composed by flipping only next-neighbour-sites that are not in order. Moreover, there is a unique isomorphism $\pi_f$ from \eqref{2before} to \eqref{2after} that can be obtained by composing suitable amplifications of the $\sF_{j,i}$ corresponding to these flips. It is also clear that the restriction of $\pi_f$ to \eqref{before} maps onto \eqref{after}. So, we get an isomorphism
\beqn{
\pi_2
\colon
(E^1_{s_1}\odot\ldots\odot E^d_{s_d})\odot(E^1_{t_1}\odot\ldots\odot E^d_{t_d}),
~\longrightarrow~
(E^1_{s_1}\odot E^1_{t_1})\odot\ldots\odot(E^d_{s_d}\odot E^d_{t_d}),
}\eeqn
and putting $u_{\bs,\bt}:=(u^1_{s_1,t_1}\odot\ldots\odot u^d_{s_d,t_d})\pi_2$ we define a product that has been constructed precisely according to our prescription. That this definition does not depend on how we put that prescription into practise, is precisely the statement that $\pi_2$ is unique.

To show associativity of that product, we start by observing that, exactly as in the construction of $\pi_2$ (but, now, starting from $E_k:=E^k_{r_k}\oplus E^k_{s_k}\oplus E^k_{t_k}$), there is a unique isomorphism
\beqn{
\pi_3
\colon
(E^1_{r_1}\odot\ldots\odot E^d_{r_d})\odot(E^1_{s_1}\odot\ldots\odot E^d_{s_d})\odot(E^1_{t_1}\odot\ldots\odot E^d_{t_d})
~\longrightarrow~
(E^1_{r_1}\odot E^1_{s_1}\odot E^1_{t_1})\odot\ldots\odot(E^d_{r_d}\odot E^d_{s_d}\odot E^d_{t_d})
}\eeqn
that can be composed out of order improving next-neighbour-flips obtained by meaningful amplifications of suitable $u^{j,i}_{\bullet,\bullet}$. We claim that both $u_{\br\bs,\bt}(u_{\br,\bs}\odot\id_\bt)$ and $u_{\br,\bs\bt}(\id_\br\odot u_{\bs,\bt})$ coincide with $u_{\br,\bs,\bt}:=(u^1_{r_1,s_1,t_1}\odot\ldots\odot u^d_{r_d,s_d,t_d})\pi_3$ (proving associativity of our product). And we show it only for one, because the other case is analogue. (For now, $u_{\br,\bs,\bt}$ is just an abbreviation; the notation will be justified, once we show associativity.)

Let us compute $u_{\br,\bs\bt}(\id_\br\odot u_{\bs,\bt})$. We have
\bmun{
u_{\br,\bs\bt}(\id_\br\odot u_{\bs,\bt})
~=~
(u^1_{r_1,s_1t_1}\odot\ldots\odot u^d_{r_d,s_dt_d})\pi'_2(\id_\br\odot(u^1_{s_1,t_1}\odot\ldots\odot u^d_{s_d,t_d})\pi_2)
\\
~=~
(u^1_{r_1,s_1t_1}\odot\ldots\odot u^d_{r_d,s_dt_d})\pi'_2(\id^1_{r_1}\odot\ldots\odot\id^d_{r_d}\odot u^1_{s_1,t_1}\odot\ldots\odot u^d_{s_d,t_d})(\id_\br\odot\pi_2),
}\emun
where $\pi'_2$ is the version of $\pi_2$ adapted to the ``times'' $(r_1,\ldots,r_d)$ and $(s_1t_1,\ldots,s_dt_d)$, in the same way as $\pi_2$ is the one adapted to the ``times'' $(s_1,\ldots,s_d)$ and $(t_1,\ldots,t_d)$ we discussed above. It is our scope to ``commute'' $\pi'_2$ through $\id_\br\odot u^1_{s_1,t_1}\odot\ldots\odot u^d_{s_d,t_d}$ with the help of the conditions in Observation \ref{prodob}\eqref{prod5}. $\pi'_2$ is composed out of (amplifications of) next-neighbour flips $u^{j,i}_{s_jt_j,r_i}$ that come to act on $E^i_{r_i}\odot u^j_{s_j,t_j}(E^j_{s_j}\odot E^j_{t_j})$ ($j<i$). By the first half of these conditions in \ref{prodob}\eqref{prod5}, the flip passes through from the left of $\id^i_{r_i}\odot u^j_{s_j,t_j}$, by which the latter transforms into $u^j_{s_j,t_j}\odot\id^i_{r_i}$, to the right, and by passing through the form is now a `mini'-permutation composed out of two flips that does to $E^j_{s_j}\odot E^j_{t_j}$ in $E^i_{r_i}\odot(E^j_{s_j}\odot E^j_{t_j})$ exactly the same as $u^{j,i}_{s_jt_j,r_i}$ does to $u^j_{s_j,t_j}(E^j_{s_j}\odot E^j_{t_j})$ in $E^i_{r_i}\odot u^j_{s_j,t_j}(E^j_{s_j}\odot E^j_{t_j})$. Doing this for every single transposition in $\pi'_2$, we get
\beqn{
\pi'_2(\id^1_{r_1}\odot\ldots\odot\id^d_{r_d}\odot u^1_{s_1,t_1}\odot\ldots\odot u^d_{s_d,t_d})
~=~
((\id^1_{r_1}\odot u^1_{s_1,t_1})\odot\ldots\odot(\id^d_{r_d}\odot u^d_{s_d,t_d}))\pi_{2'},
}\eeqn
where $\pi_{2'}$ reorders
\beqn{
E^1_{r_1}\odot\ldots\odot E^d_{r_d}\odot(E^j_{s_1}\odot E^j_{t_1})\odot\ldots\odot(E^j_{s_d}\odot E^j_{t_d})
}\eeqn
(the range of $(\id_\br\odot\pi_2)$) into
\beqn{
(E^1_{r_1}\odot(E^j_{s_1}\odot E^j_{t_1}))\odot\ldots\odot(E^d_{r_d}\odot(E^j_{s_d}\odot E^j_{t_d})).
}\eeqn
Since $\pi_{2'}(\id_\br\odot\pi_2)$ is made up out of order improving flips only, by the uniqueness statement in Lemma \ref{pi_flem} (applied appropriately to the (co)restriction $\pi_3$!), we have $\pi_{2'}(\id_\br\odot\pi_2)=\pi_3$. So, $u_{\br,\bs\bt}(\id_\br\odot u_{\bs,\bt})=u_{\br,\bs,\bt}$. The other equality, $u_{\br\bs,\bt}(u_{\br,\bs}\odot\id_\bt)=u_{\br,\bs,\bt}$ follows in the same way, now, making use of the other half of the conditions in  \ref{prodob}\eqref{prod5}.\qed

\lf
We  briefly address the question when two sets of flips, $u^{j,i}_{t_j,s_i}$ and $u'^{j,i}_{t_j,s_i}$, define isomorphic product system structures ${E^{1,\ldots,d}}^\odot$ and ${E'^{1,\ldots,d}}^\odot$ on $E^{1,\ldots,d}$. Necessarily, automorphisms $a_\bt$ of the correspondences $E_\bt$, to form an isomorphism intertwining the two products $u_{\bs,\bt}$ and $u'_{\bs,\bt}$, have to restrict to automorphisms $a^k_{t_k}$ of the marginals ${E^k}^\odot$. And, necessarily, $a_\bt=a^1_{t_1}\odot\ldots\odot a^d_{t_d}$. One easily checks that a family of automorphism ${E^k}^\odot$ defines an isomorphism from ${E^{1,\ldots,d}}^\odot$ to ${E'^{1,\ldots,d}}^\odot$ if and only if
\beqn{
(a^1_{s_1}\odot a^1_{t_1}\odot\ldots\odot a^d_{s_d}\odot a^d_{t_d})\pi_2
~=~
\pi'_2(a^1_{s_1}\odot\ldots\odot a^d_{s_d}\odot a^1_{t_1}\odot\ldots\odot a^d_{t_d}),
}\eeqn
where $\pi_2$ is the permutation from the proof of Theorem \ref{prodthm} and where $\pi'_2$ is the same permutation expressed, however, using $u'^{j,i}_{t_j,s_i}$ instead of $u^{j,i}_{t_j,s_i}$. Clearly, this happens if and only if the $a^k_{t_k}$ intertwine appropriately the next neighbour flips composing the permutation, that is, if and only if
\beq{ \label{isoprodcond}
(a^j_{t_j}\odot a^i_{s_i})u^{j,i}_{t_j,s_i}
~=~
u'^{j,i}_{t_j,s_i}(a^i_{s_i}\odot a^j_{t_j}).
}\eeq
For instance, if all $a^k_{t_k}$ are identities, then necessarily $u^{j,i}_{t_j,s_i}=u'^{j,i}_{t_j,s_i}$: An isomorphism between different product system structures on $E^{1,\ldots,d}$ is necessarily nontrivial on at least one of the marginals ${E^k}^\odot$. On the other hand, if $u^{j,i}_{t_j,s_i}=u'^{j,i}_{t_j,s_i}$, then \eqref{isoprodcond} tells which of the elements in $\euf{Aut}({E^1}^\odot)\times\ldots\times\euf{Aut}({E^d}^\odot)$ lift to elements in $\euf{Aut}({E^{1,\ldots,d}}^\odot)$. Finally, we may ask, given ${E^{1,\ldots,d}}^\odot$, for which elements in $\euf{Aut}({E^1}^\odot)\times\ldots\times\euf{Aut}({E^d}^\odot)$ the definition
\beqn{
u'^{j,i}_{t_j,s_i}
~:=~
(a^j_{t_j}\odot a^i_{s_i})u^{j,i}_{t_j,s_i}(a^i_{s_i}\odot a^j_{t_j})^*
}\eeqn
defines a product system ${E'^{1,\ldots,d}}^\odot$. The answer is for \bf{all}! We omit the straight-forward proof.

% \newpage
\bob
Let us recall that in this section we did two things: Firstly, we showed that every product system $E^\odot$ over a product $\bS=\bS^1\times\ldots\times\bS^d$ decomposes, in a very specific sense, into the tensor product $E^{1,\ldots,d}$ of its marginals ${E^k}^\odot$.
% (As explained in Remark \ref{HSTPrem}, for Hilbert spaces there is the \it{external-type} tensor product; but this captures by far not all product systems structures. It is important that our description in terms of $E^{1,\ldots,d}$ and the flips, really captures all product systems over a product $\bS$.) 
Secondly, we examined how to define a product system structure over $\bS$ on the tensor product $E^{1,\ldots,d}$ of given product systems ${E^k}^\odot$ over $\bS^k$.
%
%It is noteworthy that the second, constructive, part also works for $d$ superproduct systems ${E^k}^\odot$ over $\bS^k$.
%\OW[New! See commented pdf.]{???? TRUE ???? Erase if doubts}
\eob

\newpage

\section{An application: Strongly commuting CP-semigroups} \label{strcomSEC}

% \OW[FOR ORR]{FOR ORR: I thought, it'd be a good idea to take out the relation with semigroups over products, giving you also the occasion and space to to discuss strongly commuting stuff (possibly in this generality, or at least for $\R^d$ with general $d$) and explain their connection with your work on $\R^2$. (See also the hint in the beginning of Section \ref{N0dSEC}!)

% I thought, here you can explore everything about strong commutativity (also in relation with the preceding section) in terms of (sub)PS and units, while the relation with your original formulation with Baruch would be in the vN-appendix.

% I commented several ToDos in the end out (around the Remark and regarding spatial things; simply see a preceding version if you want to recall them). Before forgetting them completely, I still would like to think a bit about them a little ....)}

In the preceding section, we entirely captured the structure of a product system $E^\odot$ over a product $\bS:=\bS^1\times\ldots\times\bS^d$ of monoids $\bS^k$, in terms of the marginal product systems ${E^k}^\odot$ over $\bS^k$ and of the flip operators $u^{j,i}_{t_j,s_i}$ ($j<i$). Actually, in these notes we are interested in CP-semigroups, and product systems are only a mean to an end to understand better the CP-semigroups, in particular, to construct dilations for them. The relation between CP-semigroups $T$ and product systems $E^\odot$ is made in terms of units $\xi^\odot$ for $E^\odot$ by the formula $T_\bt=\AB{\xi_\bt,\bullet\xi_\bt}$. The following corollary of Theorem \ref{prodthm} settles entirely the structure of units of $E^\odot$ in terms of units for the marginal product systems and the flips.

\bcor \label{uniprodcor}
Let ${\xi^k}^\odot$ denote units for the marginal subsystems ${E^k}^\odot$ of a product system $E^\odot$ over $\bS^1\times\ldots\times\bS^d$. Then there is a unit $\xi^\odot$ for $E^\odot$ such that $\xi_{t_k}=\xi^k_{t_k}$ if and only if
\begin{subequations}
\beq{\label{ud1}
\xi^i_{s_i}\xi^j_{t_j}
~=~
\xi^j_{t_j}\xi^i_{s_i}
}\eeq
for all $j<i$, $s_i\in\bS^i$, $t_j\in\bS^j$. If the product system $E^\odot$ has been obtained from ${E^k}^\odot$ and $u^{j,i}_{t_j,s_i}$ as in Theorem \ref{prodthm}, then this equation reads
\beq{\label{ud2}
u^{j,i}_{t_j,s_i}(\xi^i_{s_i}\odot\xi^j_{t_j})
~=~
\xi^j_{t_j}\odot\xi^i_{s_i}.
}\eeq
\end{subequations}
\ecor

\proof
Even if $E^\odot$ has not been constructed as in Theorem \ref{prodthm} but is given, we may define $u^{j,i}_{t_j,s_i}:=u_{t_j,s_i}^*u_{s_i,t_j}$, and \eqref{ud1} is equivalent to \eqref{ud2} after applying $u_{t_j,s_i}^*$ to the former. So, it is sufficient to assume that $E^\odot$ is given as in Theorem \ref{prodthm}. But here, \eqref{ud2} just means that
\beqn{
\pi_2(\xi^1_{s_1}\odot\ldots\odot\xi^d_{s_d}\odot\xi^1_{t_1}\odot\ldots\odot\xi^d_{t_d})
~=~
\xi^1_{s_1}\odot\xi^1_{t_1}\odot\ldots\odot\xi^d_{s_d}\odot\xi^d_{t_d},
}\eeqn
and the statement follows because all ${\xi^k}^\odot$ are units.\qed

\lf
While the form in \eqref{ud1} is surely more intuitive and perfectly symmetric, the form in \eqref{ud2} reminds very much a condition that occurs in (one of) the definition(s) of \it{strongly commuting} CP-maps:

\bdefi \label{scCPdefi}
Let $S$ and $T$ be CP-maps on a unital \nbd{C^*}algebra and put $(\sE,\xi):=$ GNS-$T$ and $(\sF,\zeta):=$ GNS-$S$. We say, $T$ and $S$ commute \phantomsection\hl{strongly}\index{strongly commuting!CP-maps}\index{CP-map!pair!strongly commuting} if there is an isomorphism
\beqn{
u\colon
\sE\odot\sF
~\longrightarrow~
\sF\odot\sE
\vspace{-2ex}
}\eeqn
such that
\vspace{-1ex}
\beqn{
u(\xi\odot\zeta)
~=~
\zeta\odot\xi.
}\eeqn
\edefi
Of course, strongly commuting CP-maps commute. (See \eqref{unitcomp}.)

\brem
The, somewhat mysterious, notion of strongly commuting CP-maps has been first introduced by Solel \cite{Sol06} in the context of dilations of discrete time two-parameter CP-semigroups; the case considered there is normal CP-maps on a von Neumann algebra $\cB$ and formulated in terms of their Arveson-Stinespring correspondences (over the commutant $\cB'$; see Appendix \ref{vNAPP}). In \cite{Sha11}, Shalit shows that the von Neumann case of Definition \ref{scCPdefi} is equivalent to the original definition in \cite{Sol06}. Only Definition \ref{scCPdefi} makes sense also in the \nbd{C^*}case. And Definition \ref{scCPdefi} gives at least some hint to intuition, which, possibly, makes appear the notion a bit less mysterious.
\erem

In order to not interrupt the flow, we give a collection of examples in the end of this section in Subsection \ref{SCex}; but we will refer to some of them in the discussion. Our scope is to reflect on how strong commutativity has occurred so far, how strong commutativity relates to \eqref{ud2}, and how Corollary \ref{uniprodcor} can help to improve the nice consequences of strong commutativity.

\lf
Let us repeat what our ultimate aim is: For a CP-semigroup $T$ (over $\bS^{op}$ on a unital \nbd{C^*}al\-gebra $\cB$) we wish to construct a product system $E^\odot$ (over $\bS$) and a unit $\xi^\odot$ for that product system such that the subproduct system generated by the unit is the GNS-subproduct system of the CP-semigroup, that is, such that $\AB{\xi_t,\bullet\xi_t}=T_t$; Observation \ref{uGNSsupob}. We, then, wish to use the pair $(E^\odot,\xi^\odot)$ to construct a (strong  and strict module) dilation $(E,\vt,\xi)$ of $T$. This has been the strategy since Bhat and Skeide \cite{BhSk00} (see also the very first part of the introduction in Section \ref{intro}), also in Muhly and Solel \cite{MuSo02} (dealing with same question as \cite{BhSk00}), Solel \cite{Sol06} (dealing with the discrete two-parameter case), and Shalit \cite{Sha08} (dealing with the continuous two-parameter case). (The latter three are for the von Neumann case only, and the ingredients first have to be translated in terms of a duality called \it{commutant} as explained in Skeide \cite[Section 2]{Ske03c}, before we see that we are speaking about the same strategy; see Appendix \ref{vNAPP}.)

The step $(E^\odot,\xi^\odot)\leadsto(E,\vt,\xi)$, we can do for Ore monoids, whenever the unit is unital, that is, whenever $T$ is a Markov semigroup; Theorems \ref{Oreindthm} and \ref{Markmodthm}. The nonunital case, remains problematic.
% \OW[OPEN: What to do if the unit of a PS is nonunital? How to at least recover Baruch's nonunital two-parameter dilation? (The only method we have for constructing the dilation of a nonunital semigroup is by passing to the unitalization; anything else will be a new construction. --Orr)]{OPEN! Important!}
(See Question \ref{Q5}. Unitalization helps in so far, as we may pass from the GNS-subproduct system of the nonunital CP-semigroup to that of its unitalization. But, it is the latter subproduct system which we have to embed into a product system obtaining, then, a dilation of the unitalization which `contains' also a dilation of the original CP-semigroup. Knowing that the GNS-subproduct system of the original CP-semigroup embeds into a product system, does not help here; for an (indirectly related) idea about the obstacles, see Lemma \ref{pswtpslem}.)
So, in this section we will mainly concentrate on how to construct product systems for GNS-subproduct systems; the final step to (module) dilation will be limited to the unital case.

Now if $\bS=\bS^1\times\ldots\times\bS^d$ is a product, Theorem \ref{prodthm} tells us exactly how a product system over $\bS$ is determined by its marginal product systems over $\bS^k$ and, \it{vice versa}, how, given product systems over $\bS^k$, we may (or not) construct out of them a product system over $\bS$ having these as marginal systems. And Corollary \ref{uniprodcor} tells us how, in this situation, a unit has to look like. Theorem \ref{prodthm} and Corollary \ref{uniprodcor} do not tell us, how we get product systems over $\bS^k$ and units for them; but if we have them, these results tell us the only way how they can be composed to a product system over $\bS$ and a unit for that product system.

In other words, if we have $d$ commuting CP-semigroups $T^k$ (so that they determine a CP-semigroup $T$ over $\bS$ with marginals $T^k$) and if, for some reason, we know (as, for instance, in the one-parameter case $\bS^k=\N_0$ or $\bS^k=\R_+$, or in other situations where Theorem \ref{indlimthm} applies) how to embed the GNS-subproduct system of each $T^k$ into a product system, then we are in the situation where we wish to apply Corollary \ref{uniprodcor}. And it is here, where strong commutativity enters:

Roughly speaking, the product systems and units constructed from the marginal CP-semi\-groups $T^k$ by applying Theorem \ref{indlimthm} to the GNS-subproduct system of $T^k$, may be composed to a product system over $\bS$ and unit if and only if they fulfill the ``right'' condition of commuting strongly. In other words, we restrict the situation of Corollary \ref{uniprodcor} to product systems obtained as in Theorem \ref{indlimthm} (or a better way to capture its properties) and make its hypotheses a definition -- and hope that at least sometimes a less complicated formulation is possible.

The product systems ${E^k}^\odot$ of the marginals $T^k$ constructed as in Theorem \ref{indlimthm}, are not only \hl{generated} by their unit (there is no smaller \bf{product} subsystem containing the unit or, equivalently, the GNS-subproduct system of the corresponding CP-semigroup), but they are \hl{spanned} by the unit (there is no smaller \bf{superproduct} subsystem containing the unit, or, equivalently, the product system is spanned by the unit as in Theorem \ref{supintthm}). Clearly, the same is true for the family $E^{1,\ldots,d}$, once we are able to turn it into a product system with unit.

So, the definition we aim at, presumes that the marginal objects are in some sense, a sense more restrictive than just being \it{generated}, spanned by the cyclic elements that form the unit. But, this is also exactly the difference between the situation in \eqref{ud2} for fixed $j<i$, $s_i\in\bS^i$, and $t_j\in\bS^j$, and the statement that the two (commuting!) CP-maps $T^i_{s_i}:=\AB{\xi^i_{s_i},\bullet\xi^i_{s_i}}$ and $T^j_{t_j}:=\AB{\xi^j_{t_j},\bullet\xi^j_{t_j}}$ commute strongly according to Definition \ref{scCPdefi}. The key requirement in Definition \ref{scCPdefi} is that the two correspondences $\sE$ and $\sF$ are the GNS-correspondences of the occurring CP-maps, that is, they are generated in a very strong sense, namely, as correspondences, by their cyclic vectors $\xi$ and $\zeta$, respectively.

In the situation of \eqref{ud2} (not fixing on only two CP-maps, but taking into account that they sit in semigroups and that the correspondences are part of product systems), this will happen rather rarely. Indeed, if, for instance, $T^j_{t_j}$ is composed, that is, if $t_j=t''_jt'_j$ so that $T^j_{t_j}=T^j_{t'_j}\circ T^j_{t''_j}$ and $E^j_{t_j}\cong E^j_{t''_j}\odot E^j_{t'_j}$, the vector $\xi^j_{t''_j}\xi^j_{t'_j}$ will be cyclic for $E^j_{t_j}$ only in very particular circumstances (for instance if $T^j_{t'_j}$ is an endomorphism). So, rarely will the correspondences $E^j_{t_j}$ occurring in \eqref{ud2} be the GNS-correspondences $\sE^j_{t_j}$ of the CP-maps $T^j_{t_j}$. Adding to the conditions in \eqref{ud2} that the $u^{j,i}_{t_j,s_i}$ (co)restrict to isomorphisms $\sE^i_{s_i}\odot\sE^j_{t_j}\rightarrow\sE^j_{t_j}\odot \sE^i_{s_i}$ is an extra condition (meaning that all members from different marginals commute strongly); an extra condition that frequently turns even out to be unnecessary.

We see, from the beginning, the notion of strong commutativity in the context of semigroups cannot be as simple as in Definition \ref{scCPdefi} for just two CP-maps. In the sequel, we present the discrete and continuous time two-parameter case, deriving some of the known results (sometimes also improving on them), and lift them to the \nbd{d}parameter case in view of Corollary \ref{uniprodcor}. Then, we discuss the general definition prepared in the preceding paragraphs.

\lf
\subsection{\normalsize The discrete $d$--parameter case} \label{SCdd}
Well, actually, everything that can be said about the \bf{discrete} \nbd{d}parameter case, is said in Section \ref{N0dSEC}. In particular, Observation \ref{N0dob} settles the general structure of product systems over $\N_0^d$ that contain the GNS-subproduct system of a CP-semigroup over $\N_0^d$: Product systems over $\N_0^d$ are (necessarily, up to isomorphism!) obtained as in Theorem \ref{N0dthm} from $d$ correspondences $E_k$ with isomorphisms $\sF_{j,i}\colon E_i\odot E_j\rightarrow E_j\odot E_i$ $(j<i)$, satisfying  the detailed exchange conditions in \eqref{Tijcond} and vectors $\xi_k\in E_k$ fulfilling $\sF_{j,i}(\xi_i\odot\xi_j)=\xi_j\odot\xi_i$ $(j<i)$. In view of Corollary \ref{uniprodcor}, the detailed exchange relations are just the relations in Observation \ref{prodob}\eqref{prod4} restricted to the correspondences $E_k$ that generate the whole product system $\bfam{E_k^{\odot n_k}}$, and the conditions on  the $\xi_k$ are the restrictions of \eqref{ud2}. All other relations in Observation \ref{prodob}\eqref{prod4} and \eqref{ud2} follow by construction; and the relations in Observation \ref{prodob}\eqref{prod5} (also fulfilled automatically) have disappeared entirely.

In the two-parameter case $d=2$, also the conditions in Observation \ref{prodob}\eqref{prod4} are vacuous. We remain with two correspondences and vectors $E_1\ni\xi_1$ and $E_2\ni\xi_2$ and an isomorphism $\sF_{1,2}\colon E_2\odot E_1\rightarrow E_1\odot E_2$ subject to the one and only one condition $\sF_{1,2}(\xi_2\odot\xi_1)=\xi_1\odot\xi_2$. In Theorem \ref{N02Mdilthm}, we exploit that
%%%% BO 
% to recover the Markov case of Solel's 
to obtain the Markov case of Solel's result
% Solel's result dilates a (possibly Markov) CP-semigroup to a (possibly non-unital) E-semigroup. It is certain that Solel's dilation is *not proved* to give a Markov dilation if the CP maps are Markov. I do not know (and I don't remember whether I knew for sure) whether Solel's dilation is *definitely not* Markov (in some case) when the CP maps are Markov. 
%%%% EO
\cite[Theorem 5.13]{Sol06} that every normal CP-semigroup over $\N_0^2$ on a von Neumann algebra has a (strong, normal module) dilation. (Of course, \bf{here }unitalization helps to recover Solel's result in full generality; see Corollary \ref{N02dilcor}.) 

But all this has nothing to do with \it{strongly commuting}. Now, we return to arbitrary $d$ as described above, but, to make the connection with \it{strongly commuting}, we assume that $E_k$ is $\sE_k\ni\xi_k$, the GNS-correspondence of $T_k:=\AB{\xi_k,\bullet\xi_k}$, so that $T_j$ and $T_i$ commute strongly for $j<i$. 
%%%% BO
% This section is referenced below (in the section on exponentiating, Example 15.7 in old numbering)), and when checking one sees that there was really no definition of strongly commuting here. So I added it. 
% I also commented out the following sentence which is not clear:
% What goes for general $E_k$, goes also here for $E_k=\sE_k$
% We say that the $d$ maps $T_1, \ldots, T_k$ \hl{commute strongly} if there are isomorphisms $\sF_{j,i}\colon E_i\odot E_j\rightarrow E_j\odot E_i$ $(j<i)$, satisfying  the detailed exchange conditions in \eqref{Tijcond} and vectors $\xi_k\in E_k$ fulfilling $\sF_{j,i}(\xi_i\odot\xi_j)=\xi_j\odot\xi_i$ $(j<i)$. As we explained above, this determines a product system over $\N_0^d$ with a unit that recovers $T_\bn$ as $\AB{\xi_\bn, \bullet \xi_\bn}$. 
%%%% EO
Additionally, we get that the marginals of the product system constructed in Theorem \ref{N0dthm}, are the GNS-product systems of the marginal CP-semigroups $\bfam{T_k^n}$. If the CP-maps $T_k$ are unital, we also get a dilation. (In the case $d=2$, this recovers \cite[Proposition 5.15]{Sol06} restricted to unital CP-maps.) We would get a dilation, if strongly commuting CP-maps had strongly commuting unitalizations; but this may fail by Example \ref{SCex5}. And, now, for general $d$, the $\sF_{j,i}$ do have do have to satisfy the detailed exchange conditions in \eqref{Tijcond}, which, for $d=2$, are vacuous. 
% \OW[OPEN: Same as last ToDo. See why Baruch gets this result and check if it translates. (The Muhly Solel approach breaks dilation into two parts: embedding of the subproduct system into a product system and then dilation of the representation to an isometric representation. Solel has an Ando-type result that every contractive representation of a PS over $\N^2$ has an isometric dilation. In Proposition 5.15 Solel assumes strongly commuting, so the embedding of subPS into PS is easy. Note that Baruch also has also Theorem 5.13 which doesn't require strongly commuting; there he uses a construction like ours to essentially embed a subPS into a PS. I added some words about this in the text (marked by BO, EO) .--Orr)]{OPEN! Important!}

\bob \label{SCnonpob}
We know that for a single CP-map $T$, the members $T^n$ of the (discrete one-parameter) semigroup generated by $T$, of course commute, but need not commute strongly among themselves. Indeed, Example \ref{SCex3} presents a CP-map $T$ that does not commute strongly with $T^2$. Of course, $T$ commutes strongly with itself. So, even if $T$ and $S:=T$ commute strongly, the marginals of the CP-semigroup $(n_1,n_2)\mapsto S^{n_2}\circ T^{n_1}$ over $\N_0^2$ need not commute strongly, pointwise. Indeed, with $T$ from the example, we get that $T$ does not commute strongly with $S^2$, in fact, not with any $S^n$ $(n\ge2)$.
\eob

We see: Neither may we expect that the correspondences that occur in Corollary \ref{uniprodcor}, are the GNS-correspondences; also in the discrete \nbd{d}parameter case this may be true, in general, at most for the generating correspondences $\sE_k$. Nor may we expect that the marginal CP-semigroups commute strongly, pointwise; in the discrete two-parameter case (for which the notion of strongly commuting CP-maps has been invented in \cite{Sol06}), requiring strong commutativity pointwise for the whole marginal semigroups, may even restrict applicability where no such restriction is necessary. In the strongly commuting continuous time two-parameter case \cite{Sha08}, this caused no little headache.

\lf
\subsection{\normalsize The continuous time $d$--parameter case ...} \label{SCdc}
While in the discrete \nbd{d}parameter case the simplification of Theorem \ref{prodthm} and Corollary \ref{uniprodcor} to Theorem \ref{N0dthm} and Observation \ref{N0dob}, is enormous, in the continuous time case, we cannot actually say much more than what Theorem \ref{prodthm} and Corollary \ref{uniprodcor} say for $\bS^k=\R_+$. What we can say is, as always, that for $d=2$ the conditions in Observation \ref{prodob}\eqref{prod4} are vacuous. But the conditions in Observation \ref{prodob}\eqref{prod4} remain. (The second half of Section \ref{N0dSEC} is dedicated to explain that they remain there in a very substantial way.)

Now, every nonzero $t_k\in\R_+$ can be decomposed. So, except in trivial cases, none of the correspondences $E^k_{t_k}\ni\xi^k_{t_k}$ will coincide with the GNS-correspondence $\sE^k_{t_k}:=\cls\cB\xi^k_{t_k}\cB\subset E^k_{t_k}$ of $T^k_{t_k}=\AB{\xi^k_{t_k},\bullet\xi^k_{t_k}}$. The best we may hope for, is that (like in the discrete case) the ${E^k}^\odot$ are the GNS-product systems of the CP-semigroups $T^k$. But, honestly, we do not see an advantage to just restrict the situation in Corollary \ref{uniprodcor} to that situation, as long as there do not occur useful consequences in the construction of the product system over $\R_+^d$ from the fact that each ${E^k}^\odot$ is spanned by ${\xi^k}^\odot$.

We obtain useful consequences, if we require, now really, that the marginal CP-semigroups commute strongly pointwise, that is, for all $j<i$ and $t_j,s_i\in\R_+$ there is an isomorphism $\upsilon^{j,i}_{t_j,s_i}\colon\sE^i_{s_i}\odot\sE^j_{t_j}\rightarrow\sE^j_{t_j}\odot\sE^i_{s_i}$ such that $\upsilon^{j,i}_{t_j,s_i}(\xi^i_{s_i}\odot\xi^j_{t_j})=\xi^j_{t_j}\odot\xi^i_{s_i}$. (We should not forget that this condition is not required and need not be fulfilled in the discrete \nbd{d}parameter case.) Since we assume that the ${E^k}^\odot$ are the GNS-product systems, they are spanned by tensor products of pieces of units. So, it is natural to try to define $u^{j,i}_{t_j,s_i}$ (acting on products of pieces of units) by iterating $\upsilon^{j,i}_{t_j,s_i}$ (acting on the pieces) suitably. This means, in particular, that the $u^{j,i}_{t_j,s_i}$, if they exist, have to (co)restrict to the  $\upsilon^{j,i}_{t_j,s_i}$ (and, then, are determined uniquely by them). So, if the $u^{j,i}_{t_j,s_i}$ exist, then the $\upsilon^{j,i}_{t_j,s_i}$ necessarily fulfill the conditions as the $u^{j,i}_{t_j,s_i}$ Observations \ref{prodob}\eqref{prod4} (being vacuous for $d=2$) and \eqref{prod5}. Switching to the coproducts $u^{k^{\scriptstyle*}}_{s_k,t_k}$ of ${E^k}^\odot$ (which (co)restrict to the isometric coproduct maps of the GNS-subproduct system of $T^k$), the conditions in Observation \ref{prodob}\eqref{prod5} read:
%%%% BO 
% 
% I rewrote to have this in perfect alignment with 14.4(5), by taking adjoints (Please check!)
%
%
\begin{subequations} \phantomsection\label{SCSG}
\beq{ \phantomsection\label{SCSGa}
\xymatrix{
\sE^i_{r_i+s_i}  \odot \sE^j_{t_j}		\ar[d]_{u^{i^{\scriptstyle*}}_{r_i,s_i}\odot\,\ssid^j_{t_j}}\ar[rrrr]^{\upsilon^{j,i}_{t_j,r_i+s_i}}
&&&&
\sE^j_{t_j} \odot \sE^i_{r_i+s_i}  \ar[d]^{\ssid^j_{t_j}\odot\,u^{i^{\scriptstyle*}}_{r_i,s_i}}
\\
\sE^i_{r_i} \odot \sE^i_{s_i} \odot \sE^j_{t_j}		\ar[rrrr]_{(\upsilon^{j,i}_{t_j,r_i}  \odot\,\ssid^i_{s_i})(\ssid^i_{r_i}\odot\, \upsilon^{j,i}_{t_j,s_i})}
&&&&
\sE^j_{t_j} \odot \sE^i_{r_i} \odot \sE^i_{s_i} 
}
}\eeq
\beq{ \phantomsection\label{SCSGb}
\xymatrix{
\sE^i_{r_i} \odot \sE^j_{s_j+t_j}  	\ar[d]_{\ssid^i_{r_i}\odot\,u^{j^{\scriptstyle*}}_{s_j,t_j}}		\ar[rrrr]^{\upsilon^{j,i}_{s_j+t_j,r_i}}
&&&&
\sE^j_{s_j+t_j}  \odot \sE^i_{r_i}\ar[d]^{u^{j^{\scriptstyle*}}_{s_j,t_j}\odot\,\ssid^i_{r_i}}
\\
\sE^i_{r_i}\odot \sE^j_{s_j} \odot \sE^j_{t_j}												\ar[rrrr]_{(\ssid^j_{s_j}\odot\, \upsilon^{j,i}_{t_j,r_i})(\upsilon^{j,i}_{s_j,r_i}  \odot\,\ssid^j_{t_j})}
&&&&
\sE^j_{s_j} \odot \sE^j_{t_j} \odot \sE^i_{r_i} 
}
}\eeq
\end{subequations}
%
%
% I rewrote to have this in perfect alignment with 14.4(5), by taking adjoints (please check!)
%
%
% \begin{subequations} \label{SCSG}
% \beq{ \label{SCSGa}
% \parbox{10cm}{\xymatrix{
% \sE^j_{t_j} \odot \sE^i_{r_i+s_i}  		\ar[d]_{\ssid^j_{t_j}\odot\,u^{i^{\scriptstyle*}}_{r_i,s_i}} \ar[rrrr]^{\upsilon^{i,j}_{r_i+s_i,t_j}}
% &&&&
% \sE^i_{r_i+s_i} \odot \sE^j_{t_j}  \ar[d]^{u^{i^{\scriptstyle*}}_{r_i,s_i}\odot\,\ssid^j_{t_j}}
% \\
% \sE^j_{t_j}  \odot  \sE^i_{r_i} \odot \sE^i_{s_i} 	\ar[rrrr]_{(\ssid^i_{r_i}\odot\, \upsilon^{i,j}_{s_i,t_j})(\upsilon^{i,j}_{r_i,t_j}  \odot\,\ssid^i_{s_i})}
% &&&&
% \sE^i_{r_i} \odot \sE^i_{s_i} \odot \sE^j_{t_j}  
% }}
% }\eeq
% \beq{ \label{SCSGb}
% \parbox{10cm}{\xymatrix{
% \sE^j_{s_j+t_j} \odot \sE^i_{r_i}    \ar[d]_{u^{j^{\scriptstyle*}}_{s_j,t_j}\odot\,\ssid^i_{r_i}}		\ar[rrrr]^{\upsilon^{i,j}_{r_i,s_j+t_j}}
% &&&&
% \sE^i_{r_i} \odot \sE^j_{s_j+t_j}  \ar[d]^{\ssid^i_{r_i}\odot\,u^{j^{\scriptstyle*}}_{s_j,t_j}}
% \\
% \sE^j_{s_j} \odot \sE^j_{t_j} \odot \sE^i_{r_i}												\ar[rrrr]_{(\upsilon^{i,j}_{r_i,s_j}  \odot\,\ssid^j_{t_j})(\ssid^j_{s_j}\odot\, \upsilon^{i,j}_{r_i,t_j})}
% &&&&
% \sE^i_{r_i}  \odot \sE^j_{s_j} \odot \sE^j_{t_j}
% }}
% }\eeq
% \end{subequations}
%%%% EO
A careful analysis of the inductive limit leading to Theorem \ref{indlimthm}, shows that Equations \eqref{SCSG} are enough to well-define isomorphisms $u^{j,i}_{t_j,s_i}$. (Actually, this all goes through, whenever all $\bS^k$ fulfill the hypotheses of Theorem \ref{indlimthm}.) And a careful application of our results on permutations from Section \ref{compSEC} and Appendix \ref{popAPP}, shows that the conditions in Observation \ref{prodob}\eqref{prod4} turn over to the $u^{j,i}_{t_j,s_i}$ when fulfilled for the $\upsilon^{j,i}_{t_j,s_i}$. We, therefore, obtain a product system and unit for $d$ one-parameter semigroups (or over $\bS^k$ as in Theorem \ref{indlimthm}) that commute strongly in this very strong sense.

While the conditions in Observations \ref{prodob}\eqref{prod4} are verbatim for the isomorphisms $\upsilon^{j,i}_{t_j,s_i}$ (and vacuous for $d=2$), we preferred to write down the condition in Observations \ref{prodob}\eqref{prod4} explicitly and in the form \eqref{SCSG}. One reason is that they do no longer involve only isomorphisms; and it is (we found) usually more desirable to choose a form where isometries occur, not coisometries, which is also adapted to the situation that we start with the marginal GNS-subproduct systems. Another reason is that the extra conditions in the form \eqref{SCSG} appeared for the continuous time two-parameter case in the \it{corrigendum} \cite{Sha08CORR} to \cite{Sha08}, and that these conditions made work the construction of product systems and, hence, a dilation in \cite{Sha08} in the unital von Neumann case. 
%%%% BO 
Later the conditions were used again in \cite{Sha11} to construct a dilation for nonunital CP-semigroups on $\sB(H)$. 
\cite{Sha08} and \cite{Sha11} also show that the constructed dilation is strongly continuous, if the dilated CP-semigroups are. We believe that also for the construction here, the methods from Skeide \cite[Appendix A]{Ske16} to show continuity generalize; we do not address these questions here.
%%%% EO

\lf
\subsection{\normalsize ... and beyond?}

Let us rest for a moment and see what we have. Apart from Definition \ref{scCPdefi} for a pair of CP-maps and $d$ mutually strongly commuting CP-maps in \ref{SCdd}, we have the following two definitions of strongly commuting CP-semigroups:

\bdefi \label{SCSGdefi}
Let $d\ge2$ and for each $k=1,\ldots,d$ let $\bS^k$ be a monoid fulfilling the hypotheses of Theorem \ref{indlimthm} and let $T^k$ be a CP-semigroup over ${\bS^k}^{op}$. Denote by $({\sE^k}^\bodot,{\xi^k}^\odot)$ the GNS-subproduct system of $T^k$, and denote by $({E^k}^\odot,{\xi^k}^\odot)$ the GNS-product system of $T^k$ according to Theorem \ref{indlimthm}.\phantomsection\index{GNS-PS-strongly}\index{GNS-subPS-strongly}\index{CP-semigroup!strongly commuting!GNS-PS-strongly}\index{CP-semigroup!strongly commuting!GNS-subPS-strongly}\index{strongly commuting!CP-semigroups!GNS-PS-strongly}\index{strongly commuting!CP-semigroups!GNS-subPS-strongly}
\begin{enumerate}
\item
We say $T^1,\ldots,T^d$ \hl{commute GNS-PS-strongly} if there exist isomorphisms $u^{j,i}_{t_j,s_i}\colon E^i_{s_i}\odot E^j_{t_j}\rightarrow E^j_{t_j}\odot E^i_{s_i}$ for $j<i$, $s_i\in\bS^i$, and $t_j\in\bS^j$ fulfilling the conditions in Observations \ref{prodob}\eqref{prod4} and \eqref{prod5}, and in Equations \eqref{ud2}.

\item
We say $T^1,\ldots,T^d$ \hl{commute GNS-subPS-strongly} if there exist isomorphisms $\upsilon^{j,i}_{t_j,s_i}\colon\sE^i_{s_i}\odot\sE^j_{t_j}\rightarrow\sE^j_{t_j}\odot\sE^i_{s_i}$ for $j<i$, $s_i\in\bS^i$, and $t_j\in\bS^j$ fulfilling the conditions in Observation \ref{prodob}\eqref{prod4}, and in Equations \eqref{SCSG} and \eqref{ud2}.
\end{enumerate}
\edefi

\bob
For $d$ mutually strongly commuting CP-maps $T_k$, the two definitions capture what we said in Subsection \ref{SCdd}. Namely, the discrete one-parameter semigroups generated by the $T_k$, $T^k$, commute GNS-PS-strongly. By Observation \ref{SCnonpob}, they need not commute strongly, pointwise; in particular, they need not commute GNS-subPS-strongly.
\eob

\bthm \label{SCSGthm} ~

\begin{enumerate}
\item
If the $T^1,\ldots,T^d$ commute GNS-subPS-strongly, then they commute GNS-PS-strongly. (The $u^{j,i}_{t_j,s_i}$ fulfill and are determined by the condition that they (co)restrict to the $\upsilon^{j,i}_{t_j,s_i}$.)

\item
If the $T^1,\ldots,T^d$ commute GNS-PS-strongly, then, by Corollary \ref{uniprodcor}, we construct a product system over $\bS:=\bS^1\times\ldots\times\bS^d$ and a unit for the CP-semigroup $T$ over $\bS^{op}$ generated by them.

\item
If the $T^1,\ldots,T^d$ are even Markov and $\bS$ is Ore, then by Theorem \ref{Oreindthm} we construct a (strict strong module) dilation for $T$.
\end{enumerate}
\ethm

\noindent
We just mentioned Observation \ref{SCnonpob}, which is an example of two discrete one-parameter CP-semigroups that commute GNS-PS-strongly, but not GNS-subPS-strongly. It appears that in the continuous time two- and \nbd{d}parameter case we have no such example. In fact, despite the condition to commute GNS-subPS-strongly appears restrictive as compared with GNS-PS-strongly commuting, we do have many examples for GNS-subPS-strongly commuting continuous time one-parameter semigroups. A trivial example are tensor products of (two or $d$) CP-semigroups that always commute GNS-subPS-strongly. (This situation is similar to the external-type tensor product of Arveson systems as in Remark \ref{HSTPrem}.) \nbd{E}Semigroups that commute, commute GNS-subPS-strongly. (This goes exactly as the two-map case in Example \ref{SCex1}.) Elementary CP-semigroups that commute, commute GNS-subPS-strongly. (See Observation \ref{elemob}. Actually, whatever $\bS$ and the semigroup $\bfam{c_t}_{t\in\bS}$ is, the GNS-correspondences $\sE_t=\cls\cB c_t\cB$ of $T_t=c_t^*\bullet c_t$ and their tensor products $\sE_s\odot\sE_t=\cls\cB c_s\cB c_t\cB$ sit in the trivial product system over $\bS$. The two ideals $\sE_s\odot\sE_t$ and $\sE_t\odot\sE_s$ both coincide with the intersection of the ideals $\sE_s$ and $\sE_t$. Applying this to the \nbd{d}parameter case, we see that suitably $\upsilon^{j,i}_{t_j,s_i}$ are given by the identity.)

Here is a class of examples arising from exponentiating discrete examples as in Section \ref{EXexpSEC}, and for which we are not able to say in general, how strong any sort of commutativity for discrete semigroups is reflected by their exponentiated versions:

\bex
Let $S_k$ be $d$ mutually strongly commuting CP-maps as in \ref{SCdd}. We construct $(F_k,\zeta_k)=$ GNS-$S_k$, and input this in Theorem \ref{N0dthm} and Observation \ref{N0dob} to obtain a product system $F^\odot$ over $\N_0^d$ and a unit $\zeta^\odot$ for $F^\odot$ such that $\AB{\zeta_\bn,\bullet\zeta_\bn}=:S\!_\bn$ is the discrete \nbd{d}parameter semigroup generated by the $S\!_k$.

We input $F^\odot$ in the construction of the exponential product systems $\DG^\odot(F^\odot)$ from Section \ref{EXexpSEC}. Necessarily (by the \it{only-if} direction of Theorem \ref{prodthm} spoken out in Observation \ref{prodob}), we may consider $\DG^\odot(F^\odot)$ as obtained from the marginal product systems $\DG^\odot(F_k)$ via unitary flips $u^{j,i}_{t_j,s_i}:=u_{t_j,s_i}^*u_{s_i,t_j}$. As in the beginning of this section, which started by adding to the generalities of Theorem \ref{prodthm} the statements about units in Corollary \ref{uniprodcor}, also here we have the obvious statement that the exponential units $\ee_t(\zeta_k)$ for the marginals $\DG^\odot(F_k)$ of $\DG^\odot(F^\odot)$ (see Appendix \ref{FockAPP}), obviously, satisfy the conditions in Corollary \ref{uniprodcor}. The corresponding unit $\xi^\odot$ for $\DG^\odot(F^\odot)$ generates a semigroup $T_\bt:=\AB{\xi_\bt\bullet\xi_\bt}$ with marginal semigroups $T^k_{t_k}:=\AB{\ee_{t_k}(\zeta_k),\bullet\ee_{t_k}(\zeta_k)}=e^{t_k\AB{\zeta_k,\bullet\zeta_k}}$; see also Footnote \ref{expFN}. (Obviously, the maps $T_\bt$ and $T^k_{t_k}$ are completely positive; but they are not contractive, unless $\zeta_k=0$. Since the $T^k$ are uniformly continuous, they may be normalized by scalar semigroups to be contractive. We avoid that, and relax, in this example, CP-semigroup to mean not necessarily contractive CP-semigroup.)

By construction, the CP-semigroups $T^k$ commute. Unfortunately, this is as far as the construction of continuous time \nbd{d}parameter semigroups by \it{exponentiating} discrete \nbd{d}parameter semigroups goes unproblematically. Of course, we can ask whether and in which sense the semigroups $T^k$ commute possibly strongly. But, let us mention that the product subsystem of $\DG^\odot(F_k)$ generated by the unit $\ee_{t_k}(\zeta_k)$ need not even be isomorphic to some $\DG^\odot(\sF_k)$. (In Bhat, Skeide, and Liebscher \cite{BLS10} there is a counter example where the generated product subsystem does not contain any vectors that commute with the algebra; in particular, there is no vacuum. In the von Neumann case , Barreto, Bhat, Liebscher, and Skeide \cite{BBLS04} show, that the generated subsystem is at least isomorphic to (the strong closure of) some $\DG^\odot(\sF_k)$; but, it need not contain the original vacuum of $\DG^\odot(F_k)$.)

\eex

%%%% BO 
% I changed things regarding the example (15.7 in old numbering), and moved it into the next subsection. --Orr.
%%%% EO

In Section \ref{qconvSEC}, we present a rich class of nontrivial positive examples, \it{quantized convolutions semigroups} (Arveson \cite{Arv02} and Markiewicz \cite{Mar03}), which illustrates that our machinery works in the continuous time \nbd{d}parameter case, generalizing (and proving it in a totally different way) Shalit \cite[Theorem 3.4.12]{Sha09} from $d=2$ to arbitrary $d$.

\lf
Here are the promised examples of commuting CP-maps with the analysis whether or not they commute strongly:

\lf
\subsection{\normalsize Examples of commuting -- strongly or not -- pairs of CP-maps} \label{SCex}
We give some (known) examples (possibly, with new reasoning sketched), because we feel we owe to give a first idea of how that notion applies -- or not -- to the reader who sees this notion for the first time. But, not to interrupt the flow, we recommend to skip the examples at first reading.

\bex \label{SCex1}
Two endomorphisms that commute, commute strongly.

%%%% BO
Indeed, the GNS-correspondence of an endomorphism $\vt$ of $\cB$ is $_\vt\cB$, that is, the Hilbert \nbd{\cB}module $\vt(\U)\cB$ with left action via $\vt$, and the cyclic vector is $\xi_\vt=\vt(\U)\U$. If we have two of them, $\vt$ and $\vt'$, then $_\vt\cB\odot{_{\vt'}}\cB\cong\cB_{\vt'\circ\vt}$ (via $b\odot b'\mapsto\vt'(b)b'$) and $\xi_\vt\odot\xi_{\vt'}\mapsto\xi_{\vt'\circ\vt}$. Form this the statement follows, once $\vt$ and $\vt'$ commute.
%%%% EO

By a similar but slightly more involved argument, a CP-map $T$ and an automorphism $\alpha$ that commute, commute strongly.

Indeed, put $(\sE,\xi):=$ GNS-$T$. Denote by $_\alpha\sE$ the Hilbert right \nbd{\cB}module $\sE$, but with left action $b.x:=\alpha(b)x$. Then $_\alpha\cB\odot\sE\cong{_\alpha}\sE$ via $b\odot x\mapsto bx$. Denote by $\sE^\alpha$ the left \nbd{\cB}module $\sE$ but with the inner product $\AB{x,y}_\alpha:=\alpha(\AB{x,y})$, turning $\sE^\alpha$ into a correspondence over $\cB$ with (the only possible) right multiplication $x.b=x\alpha^{-1}(b)$. Then $\sE\odot{_\alpha}\cB\cong\sE^\alpha$ via $x\odot b\mapsto x\alpha^{-1}(b)$. Using also that $\alpha\circ T=T\circ\alpha$, it follows that
\beqn{
_\alpha\sE
~\ni~
b\xi b'
~\longmapsto~
\alpha^{-1}(b)\xi\alpha^{-1}(b')
~\in~
\sE^\alpha
}\eeqn
defines an isomorphism which does the job.

These were first observed (for von Neumann algebras) in Solel \cite[Lemma 5.4]{Sol06}.
\eex

\bex \label{SCex2}
All CP-maps on $M_n$ that commute, commute strongly.

Indeed, if we put $(\sE,\xi):=$ GNS-$T$ and $(\sF,\zeta):=$ GNS-$S$, by appropriate application of Appendix \ref{vNAPP}\ref{vNBGmod} we get $\sE\odot\sF=M_n\otimes\eH\otimes\eG$ and $\sF\odot\sE=M_n\otimes\eG\otimes\eH$ (where $\sE=M_n\otimes\eH$ and $\sF=M_n\otimes\eG$), and projections $p\in\sB(\eH\otimes\eG)$ and $q\in\sB(\eG\otimes\eH)$ such that the subcorrespondences generated by $\xi\odot\zeta$ and $\zeta\odot\xi$, are $(\id_{\C^n}\otimes p)(\sE\odot\sF)$ and $(\id_{\C^n}\otimes q)(\sF\odot\sE)$, respectively. Since $T$ and $S$ commute, both subcorrespondences are isomorphic to the GNS-correspondence of the same CP-map $S\circ T=T\circ S$ via the unique isomorphism obtained by extending $\xi\odot\zeta\mapsto\zeta\odot\xi$. Therefore, $p$ and $q$ have the same (finite) rank, and since $\eH\otimes\eG\cong\eG\otimes\eH$ are finite-dimensional, also $\id_{\eH\otimes\eG}-p$ and  $\id_{\eG\otimes\eH}-q$ have the same (finite) rank. Therefore the isomorphism between the subcorrespondences extends as an isomorphism $u\colon\sE\odot\sF\rightarrow\sF\odot\sE$, the isomorphism we seek.

This was first observed in \cite[Proposition A.1]{Sha08}.
\eex

\bex \label{SCex3}
Actually, ``inventing'' a commuting pair of CP-maps (to examine, then, if they commute strongly), is not so easy a task. (From the view point of linear algebra, CP-maps behave just like generic linear maps; thus, the problem is not easier than getting commuting linear maps on a vector space.) An efficient way, is to look at different members in a one-parameter CP-semigroup, such as $T$ and $T^2$. Maybe a bit surprisingly, we do get this way even examples of commuting CP-maps that do not commute strongly. The following example is from Shalit \cite[Section A.5]{Sha08}.

%%%% BO

% I think we should do it simpler, and self contained. I changed the exposition to the following, but then I changed back when re-reading found that this is referred to in the last example of the section and changed that too, which shortens it somewhat. --Orr

We take the commutative algebra $\C^n$, the diagonal subalgebra of $M_n$. Obviously (see also Skeide \cite[Corollary 1.7.9]{Ske01}), the structure of a correspondence $\sE$ over $\C^n$ is determined (up to isomorphism) by a matrix of Hilbert spaces $\bfam{H_{i,j}}_{i,j=1,\ldots,n}$ (up to matrix-entry-wise isomorphism). The algebra acts from the left and from the right in the ``diagonal way''. The \nbd{i}component of the inner product of matrices $X,Y\in\sE$  is $\AB{X,Y}_i=\sum_j\AB{x_{j,i},y_{j,i}}$. If $\sF=\bfam{G_{i,j}}_{i,j=1,\ldots,n}$ is another \nbd{\C^n}correspondence, then the \nbd{i,j}entry of the tensor product $\sE\odot\sF$ is $\bigoplus_k H_{i,k}\otimes G_{k,j}$. It follows that $\dim(\sE\odot\sF)_{i,j}=\sum_k\dim H_{i,k}\cdot\dim G_{k,j}$.

A CP-map $T$ on $\C^n$ is given by a matrix $\bfam{T_{i,j}}$ with entries in $\R_+$. The Hilbert space $H_{i,j}$ of the GNS-correspondence is one-dimensional if $T_{i,j}\ne0$ and zero-dimensional if $T_{i,j}=0$. Suppose that $n\ge3$ and suppose that $T_{i,j}\ne0$ for all but one pair $(i_0,j_0)$ with $i_0\ne j_0$. Then for $S:=T^2$, we have $S_{i,j}\ne0$ for \bf{all} $i,j$. (The latter definitely fails for $n=2$, as Markov maps of the form $\rtMatrix{1-b~&b\\&1}$ show.) Let $\sF$ be the GNS-correspondence of $S$. Then for $j\ne j_0$,
\beqn{
\dim(\sE\odot\sF)_{i_0,j}
~=~
\sum_k1-\delta_{k,j_0}
~=~
n-1,
\text{~~~~~~while~~~~~~}
\dim(\sF\odot\sE)_{i_0,j}
~=~
\sum_k1
~=~
n.
}\eeqn
So, $\sE\odot\sF\ncong\sF\odot\sE$ and $T$ and $T^2$ commute but do not commute strongly.

The structure is taken from the concrete example
$T=
\rtMatrix{
\frac{1}{2}&0&\frac{1}{2}
\\
\frac{1}{4}&\frac{1}{2}&\frac{1}{4}
\\
\frac{1}{4}&\frac{1}{2}&\frac{1}{4}
}
$
from Shalit \cite[Section A.5]{Sha08}.

% We take the commutative algebra $\C^n$, and write $e_i$ for the standard basis vector in $\cB$. A CP-map $T$ on $\C^n$ is given by a matrix $\bfam{T_{i,j}}$ with entries in $\R_+$. If we consider the standard construction of the GNS-correspondence $\sE$ as a quotient of $\cB \otimes \cB$, then it is then easy to check that $\bCB{e_i \otimes e_j \colon T_{j,i} \neq 0}$ is a basis for $\sE$. Likewise, if $\sF$ is the GNS-correspondence of a CP map $S$ on $\cB$, then one may check that $\bCB{(e_i \otimes e_j) \odot (1 \otimes e_l) \colon S_{lj} T_{j,i} \neq 0}$ forms a basis for $\sE \odot \sF$. 

% Suppose now that $n\ge3$ and suppose that $T_{p,q}\ne0$ for all but one pair $(p_0,q_0)$ with $p_0\ne q_0$; for example, take
% $T=
% \rtMatrix{
% \frac{1}{2}&0&\frac{1}{2}
% \\
% \frac{1}{4}&\frac{1}{2}&\frac{1}{4}
% \\
% \frac{1}{4}&\frac{1}{2}&\frac{1}{4}
% }
% $. Then for $S:=T^2$, we have $S_{i,j}\ne0$ for \bf{all} $i,j$. (The latter definitely fails for $n=2$, as Markov maps of the form $\rtMatrix{1-b~&b\\&1}$ show.) Let $\sF$ be the GNS-correspondence of $S$. Then for $k\ne p_0$,
% \beqn{
% \dim(e_{q_0} \sE\odot\sF e_k)
% ~=~
% \# \bCB{j \colon S_{kj} T_{j,q_0} \neq 0}
% ~=~
% n-1,
% }\eeqn
% while
% \beqn{
% \dim(e_{q_0} \sF\odot\sE e_k)
% ~=~
% \# \bCB{j \colon T_{kj} S_{j,q_0} \neq 0}
% ~=~
% n.
% }\eeqn
% So, $\sE\odot\sF\ncong\sF\odot\sE$ and $T$ and $T^2$ commute but do not commute strongly.
%%%% EO
\eex

\bex \label{SCex4}
In the preceding example, we proved that the commuting CP-maps do not commute strongly by showing that the GNS-correspondences have non-isomorphic tensor products. It is natural to ask if there is an example where $\sE\odot\sF$ and $\sF\odot\sE$ are isomorphic, but no isomorphism can intertwine $\xi\odot\zeta$ and $\zeta\odot\xi$. Such an example exists:

Let $H$ be a Hilbert space, let $v\in\sB(H)$ a proper isometry (so that $H$ is infinite-di\-mensional), and suppose $v$ has a (unital) eigenvector $\vk$. Put $T:=v^*\bullet v$ (so that GNS-$T=(\sB(H),v)$) and put $S:=\id_H\AB{\vk,\bullet\vk}$ (so that GNS-$S=(H\otimes\sB(H),\vk\otimes\id_H)$, where $\sB(H)$ acts on the left factor of $H\otimes\sB(H)$). It follows that $\sB(H)\odot(H\otimes\sB(H))\cong H\otimes\sB(H)$ via $b\odot(h\otimes b')\mapsto bh\otimes b'$, and $(H\otimes\sB(H))\odot\sB(H)\cong H\otimes\sB(H)$ via $(h\otimes b)\odot b'\mapsto h\otimes bb'$. But $v\odot(\vk\otimes\id_H)=v\vk\otimes\id_H$ generates the whole thing $H\otimes\sB(H)$ as a correspondence, while the correspondence generated by $(\vk\otimes\id_H)\odot v=\vk\otimes v$ leaves the nonzero complement $H\otimes(\id_H-vv^*)\sB(H)$. No automorphism of $H\otimes\sB(H)$ can intertwine these to vectors, so $T$ and $S$, which clearly commute, do not commute strongly, despite having isomorphic tensor products of their GNS-correspondences.

This example is Solel \cite[Example 5.5]{Sol06}.

It is noteworthy that it does not matter, if we take the tensor product $H\otimes\sB(H)$, that is, in the \nbd{C^*}category, as we did, or if we pass to strong closures, that is, the von Neumann category; only the latter is really equivalent to \cite[Example 5.5]{Sol06}. It might be tempting to try to get a separable example in the \nbd{C^*}case, by passing to $\cB:=\wt{\sK(H)}$ for a separable, infinite-dimensional Hilbert space. However, it turns out (as a matter of Fredholm index and its stability under compact perturbations; we do not give details) that in this case the two tensor products of the GNS-correspondences are no longer isomorphic.
\eex

\bex \label{SCex5}
Note that both concrete examples in \ref{SCex3} and \ref{SCex4}, actually, are Markov. On the other hand, the Markov matrix $T:=\rtMatrix{1-b~&b\\&1}$ that illustrated how the key ingredient to the argument in \ref{SCex3} breaks down, is even a unitalization of the CP-map $z\mapsto(1-b)z$ on $\C$, which commutes strongly with its square. This raises the question how strong commutation of CP-maps is related to strong commutation of their unitalizations, in general.

We start with the usual situation $(\sE,\xi)=$ GNS-$T$ and $(\sF,\zeta)=$ GNS-$S$, and construct $(\wt{\sE},\wt{\xi})=$ GNS-$\wt{T}$ and $(\wt{\sF},\wt{\zeta})=$ GNS-$\wt{S}$ as in \ref{GNSuni}. So, $\wt{\sE}=\sE\oplus\ol{\wh{\xi}\wt{\cB}}$ and $\wt{\xi}=\xi\oplus\wh{\xi}$ where $\wh{\xi}=\sqrt{\rule[1.3ex]{0pt}{1ex}\,\smash{\wt{\U}-T(\U)}}=\wt{\U}-\U+\sqrt{\U-T(\U)}$ and where the left action of $\wt{\cB}$ on $\wh{\xi}$ is the unitalization of the \nbd{0}action of $\cB$, and similarly for $(\wt{\sF},\wt{\zeta})$. We have $\U\wt{\sE}=\sE=\sE\U$ and $(\wt{\U}-\U)\wt{\sE}=\ol{\wh{\xi}\wt{\cB}}$.
Therefore,
\baln{
\wt{\sE}\odot\wt{\sF}
&
~=~
(\sE\odot\sF)\oplus(\sE\odot\ol{\wh{\zeta}\wt{\cB}})\oplus(\ol{\wh{\xi}\wt{\cB}}\odot\wt{\sF)}
\\
&
~=~
(\sE\odot\sF)\oplus~~~~~\zero~~~~~\oplus\,\ol{\wh{\xi}\odot\wt{\sF}}
\\
&
~=~
(\sE\odot\sF)~~~~~~~~~~~~~~~~~~~\oplus~{^{\wh{~}}}\Bfam{\ol{(\U-T(\U))\sF}}\,\oplus~{^{\wh{~}}}\Bfam{\ol{(\wt{\U}-S(\U))\wt{\cB}}},
}\ealn
where by the pre-superscript ${^{\wh{~}}}\Bfam{}$ we indicate (where necessary) that \nbd{\wt{\cB}}correspondence obtained by equipping the Hilbert \nbd{\wt{\cB}}module inside the brackets with the left action obtained by unitalization of the \nbd{0}action of $\cB$. Moreover,
\beqn{
\wt{\xi}\odot\wt{\zeta}
~=~
(\xi\odot\zeta)\oplus\sqrt{\U-T(\U)}\zeta\oplus\sqrt{\rule[1.1ex]{0pt}{1ex}\,\smash{\wt{\U}-S(\U)}}
}\eeqn
Necessarily, if $\wt{T}$ and $\wt{S}$ commute strongly, then so do $T$ and $S$. (A suitable isomorphism $\wt{u}$ for the unitalization restricts, after multiplying with $\U$ from the left to a suitable $u$ for the original CP-maps.)

The conclusion in the opposite direction need not be true. We return to the CP-map $\wt{T}:=\rtMatrix{1-b~&b\\&1}$ on $\C^2$ from \ref{SCex3}, which describes the unitalization of the CP-map $T$ on $\C$ given by multiplication with $1-b$. We assume $0<b<1$, so that the $T$ is neither $0$ nor Markov. We take $\wt{S}:=\rtMatrix{0~&1\\&1}$, the unitalization of the zero-map $S$ on $\C$. Clearly, $T$ and $S$ commute strongly. (Actually, this is true for all CP-maps on $\C=M_1$, but for $S=0$ it is particularly obvious.) Returning to expressing the GNS-correspondences as square matrices as explained in Example \ref{SCex3},
%%%% BO
we find
\baln{
\sMatrix{\C&\C\\&\C}
~=~
\wt{\sE}
&
~\ni~
\wt{\xi}
~=~
\sMatrix{\sqrt{1-b}&\sqrt{b}\\&1},
&
\sMatrix{0&\C\\&\C}
~=~
\wt{\sF}
&
~\ni~
\wt{\zeta}
~=~
\sMatrix{0&1\\&1}.
}\ealn
Therefore,
\baln{
\wt{\sE}\odot\wt{\sF}
&
~=~
\sMatrix{\C&\C\\&\C}\odot\sMatrix{0&\C\\&\C}
~=~
\sMatrix{0&\C^2\\&\C},
&
\wt{\sF}\odot\wt{\sE}
&
~=~
\sMatrix{0&\C\\&\C}\odot\sMatrix{\C&\C\\&\C}
~=~
\sMatrix{0&\C\\&\C}.
}\ealn
Since the two tensor products are not isomorphic, the two CP-maps $\wt{T}$ and $\wt{S}$ commute, but not strongly.
% By counting the number of $i,j,k$ such that $\wt{T}_{lj}\wt{S}_{ji} \neq 0$, and then counting those for which $\wt{S}_{lj}\wt{T}_{ji} \neq 0$ (as in Example \ref{SCex3}), we see that the two CP-maps $\wt{T}$ and $\wt{S}$ commute, but not strongly.
%%%% EO

% Suppose we do have $u$ illustrating the $T$ and $S$ commute strongly. Can we find a unitary
% \beqn{
% ^{\wh{~}}u
% \colon
% {^{\wh{~}}}\Bfam{\ol{(\U-T(\U))\sF}}\,\oplus~\ol{(\wt{\U}-S(\U))\wt{\cB}}
% ~\longrightarrow~
% {^{\wh{~}}}\Bfam{\ol{(\U-S(\U))\sE}}\,\oplus~\ol{(\wt{\U}-T(\U))\wt{\cB}}
% }\eeqn
% (left linearity is automatic due to the trivial form of the left action) that acts on $(\wt{\U}-\U)\wt{\xi}\odot\wt{\zeta}$ as
% \beqn{
% \sqrt{\U-T(\U)}\zeta\oplus\sqrt{\rule[1.1ex]{0pt}{1ex}\,\smash{\wt{\U}-S(\U)}}
% ~\longmapsto~
% \sqrt{\U-S(\U)}\xi\oplus\sqrt{\rule[1.1ex]{0pt}{1ex}\,\smash{\wt{\U}-T(\U)}}?
% }\eeqn
% First of all, multiplying with $\U$ from the right, we get rid of the twiddles. We get
% \beqn{
% \sqrt{\U-T(\U)}\zeta\oplus\sqrt{\U-S(\U)}
% ~\longmapsto~
% \sqrt{\U-S(\U)}\xi\oplus\sqrt{\U-T(\U)}.
% }\eeqn
% It is easy to get examples, where it is not possible to find such $^{\wh{~}}u$.

\eex

% \OW{I am still not giving up the possibility to use two (or more) spatial one-parameter product systems and their spatial product to construct two- (or multi-)parameter product systems! But there is something weird there, I don't understand yet. We can try and discuss this in the case of time-ordered product systems to understand better the structure.}

% \OW{Now, after the remark, I was dreaming of formulating a theorem answering when the family $E^1_{t_1}\odot\ldots\odot E^d_{t_d}$ can be turned into product system over $\bS^1\times\ldots\times\bS^d$. That is actually, why I put a separate section here that should be introductory to all the following sections. But I realize that in order to formulate that we will need isomorphisms $u^{i,j}_{t_i,s_j}\colon E^j_{s_j}\odot E^i_{t_i}\rightarrow E^i_{t_i}\odot E^j_{s_j}$ for $i<j$, and all that permutation stuff.

% I simply try and see what comes out. So in case I don't cancel this OW (that is, why you can read it), to keep the message about my intentions in mind. But for the flux of reading ignore it; and, of course, ignore also the remark. (All my remarks may be ignored safely at the place where they stand, otherwise I would not call them remarks. A remark either provides further information, or might be referred to from a later place. But it is NEVER necessary for the argument being discussed around the remark.)
% }

%%%% BO 
% The following Example is old Example 15.7, moved here and changed. 
% I added some observations which are partly written out in my old file "ExamplesOfTwoParamCpSG.pdf". 
%%%% EO

\newpage

\section[\sc{Examples:} Quantized convolution semigroups]{Examples: Quantized convolution semigroups} \label{qconvSEC}

Let $\om\colon H\times H\rightarrow\R$ be a \hl{symplectic form}\index{symplectic form} on the real vector space $H$. (This means, $\om$ is bilinear and anti-symmetric.) Then the \hl{\nbd{C^*}algebra of canonical commutation relations}\index{symplectic form!CCR-algebra over}\index{CCR-algebra! $CCR(H,\om)$ over a sympectic form} $CCR(H,\om)$ is defined%
\footnote{
For $CCR(H,\om)$ we refer to the very useful monograph \cite{Pet90} by Petz. The only pitfall is that, despite being motivated very much by physics, his Hilbert spaces have inner products that are anti-linear in their second argument.
\vspace{1ex}
}
as the universal \nbd{C^*}algebra generated by symbols $w_h$ subject to the relations
\baln{
w_h^*
&
~=~
w_{-h},
&
w_hw_{h'}
&
~=~
w_{h+h'}e^{i\om(h,h')}.
}\ealn
(The latter relation implies the commutation relation $w_hw_{h'}=e^{i2\om(h,h')}w_{h'}w_h$.) Not only does the universal \nbd{C^*}algebra exist. (Indeed, since $\om(h,h)=0$, the $*$-algebra generated by $w_h$ is unital with unit $\U=w_0$, and the $w_h$ are all unitary; hence, for all $a$ in this $*$-algebras, the supremum of $\norm{\pi(a)}$ over all representations $\pi$ is finite and defines a \nbd{C^*}(semi)norm.) But by \cite[Theorem 2.1]{Pet90}, every nonzero homomorphism is, actually, isometric, hence, an isomorphism onto its range. By this uniqueness result, we have
\beqn{
CCR(H_1\oplus H_2,\om_1\oplus\om_2)
~\cong~
CCR(H_1,\om_1)\otimes CCR(H_2,\om_2),
}\eeqn
and the \nbd{C^*}norm of the tensor product is unique.

For every Hilbert space $H$, we define a symplectic form $\om(h,h'):=-\frac{1}{2}\Im\AB{h,h'}$, which is also \phantomsection\hl{nondegenerate}\index{symplectic form!nondegenerate} (that is, $\om(h,h')=0$ for all $h'$ implies $h=0$).%%
\footnote{
The factor $-\frac{1}{2}$ is a convention to be compatible with the works \cite{Arv02,Mar03}. The minus has it origin in the choice of the argument in which inner products linear. Anyway, if $\om$ is a symplectic form, then so is $\lambda\om$ for real $\lambda$.
\vspace{1ex}
}
Every nondegenerate symplectic form on a finite-dimensional space (necessarily \nbd{2n}dimensional) arises in this way by establishing an isomorphism of real vector spaces between $H$ and $\C^n$. The symplectic space $H$ decomposes as $\C^n=\C\oplus\ldots\oplus\C$. So for the finite-dimensional case, it is enough to understand $CCR(\C,\om)$ with $\om(z,z')=\frac{yx'-xy'}{2}$.

Let $H$ be a finite-dimensional Hilbert space. A representation $\rho\colon w_h\mapsto W_h$ of $CCR(H,\om)$ is \phantomsection\hl{continuous}\index{CCR-algebra!continuous representation} if the function $h\mapsto W_h$  is strongly (or, equivalently, weakly) continuous. The Stone-von Neumann theorem asserts that every unital continuous representation is an amplification of the (unique, up to unitary equivalence) irreducible one.%
\footnote{
This is proved by reduction to the (complex) one-dimensional case $CCR(\C,\om)$. There is a bit of confusion, in what exactly the (Mackey-)Stone-von Neumann theorem means in that case. The version in \cite[Theorem 1.2]{Pet90} states just uniqueness of the irreducible continuous representation. (To show what we stated, one has to add to this that any cyclic representation is irreducible.) Other versions put emphasis on the canonical commutation relations satisfied by the generators of the one-parameter unitary groups $W_x$ and $W_{iy}$. Also here, in the identification of $W_x$ with $e^{iqx}$ and of $W_{iy}$ with $e^{ipy}$, where the generator $q$ and $p$ should satisfy the \it{Heisenberg commutation relation} $\SB{q,p}=i$, unfortunate sign conventions (from physics!) have to be sorted out. For suitable sign conventions, the group $e^{iqx}$ acts on $L^2(\R)\ni f\colon \lambda\mapsto f(\lambda)$ as multiplication by $e^{i\lambda x}$ (usual convention) and $e^{ipy}$ acts as right shift by $y$ (usual would be the left shift by $y$, so our generator $p$ is minus the usual one). The recent (self-contained) proof in Bhat and Skeide \cite[Observation 2.7]{BhSk15}, reduces the statement to the structure result about pure isometry semigroups (namely, being multiples of the right shift).
}

Henceforth, we put $H:=\C$. We do not need to know the concrete form of the (unital) representation $\rho\colon w_z\mapsto W_z\in\sB(L)$ for some Hilbert space $L$ (the letter $L$ to remind that it is $L^2(\R,K)$ for some $K$), but only that it is continuous. We shall identify $CCR(\C,\om)$ with $\cB:=\rho(CCR(\C,\om))\subset\sB(L)$. Note that $\ol{\cB}^s=\sB(L)$, if (and only if) we choose the unique irreducible continuous representation.

Following Arveson \cite{Arv02} and Markiewicz \cite[Part II]{Mar03} 
%%%% BO new
% there was an overfull hbox here
(see Davies' book \cite[Chapter 8]{Dav76}, too), 
% (see also Chapter 8 in Davies' book \cite{Dav76}), 
%%%% EO
for every  (Borel) probability measure $\mu$ on $\C$, we define a Markov map $T_\mu:=\int W_z^*\bullet W_z\,\mu(dz)$ on $\cB$. (Well, actually, for compatibility with the rest of our notes, as compared with \cite{Arv02,Mar03}, we switched $W_z^*$ and $W_z$, which amounts to exchanging $z$ and and $-z$.) While it is clear that this  defines a (normal) Markov map on $\sB(L)$ (that is, what \cite{Arv02,Mar03} actually looked at), it follows from
\beq{ \label{hatmuT}
T_\mu(W_w)
~=~
W_w\,\hat{\mu}(w),
}\eeq
where
\beqn{
\hat{\mu}(w)
~:=~
\int e^{i2\om(w,z)}\,\mu(dz),
}\eeqn
that $T_\mu$ maps $\cB$ into $\cB$. Now that we know that $T_\mu$ is a Markov map into $\cB$, we may even forget the continuous representation $\rho$, from which we derived that piece of information, and work directly with \eqref{hatmuT} as definition of $T_\mu$.

Putting $z=x+iy$, $w=u+iv$, and so forth, we see that $w\mapsto\hat{\mu}(iw)=\int e^{i(ux+vy)}\,\mu(dz)$ is nothing but the \hl{characteristic function} (or \hl{Fourier transform}) of the probability measure $\mu$. Of course, the (injective!) map $\mu\mapsto\hat{\mu}$ shares the property of characteristic functions that $\wh{\mu*\nu}=\hat{\mu}\hat{\nu}$. Therefore, $T_{\mu*\nu}=T_\mu\circ T_\nu$ (recovering the restriction of \cite[Proposition 1.7]{Arv02} to probability measures), so that $\mu\mapsto T_\mu$  is a homomorphism of monoids (actually, an isomorphism) from  the monoid of probability measures on $\C$ under convolution (with the point measure $\delta_0$ at $0$ as neutral element) $\bS_*$ onto the monoid of Markov maps $T_\mu$ (with $T_{\delta_0}=\id_\cB$), the so-called \index{}\hl{quantized convolution semigroup}\index{quantized convolution semigroup}\index{CP-semigroup!quantized convolution semigroup}.

Clearly, convolution of probability measures on $\C$ is commutative (as manifest also from $\hat{\mu}\hat{\nu}=\hat{\nu}\hat{\mu}$). So, $T_*:=\bfam{T_\mu}_{\mu\in\bS_*}$ is an abelian Markov semigroup over $\bS_*$. Definition \ref{SCSGdefi} gives two versions of strongly commuting for \nbd{d}parameter semigroups. It, therefore, makes sense to ask to what extent the \nbd{d}parameter subsemigroups of $T_*$  are strongly commuting. Actually, we now shall show that the whole semigroup $T_*$ is, in a sense, \it{GNS-subPS-strongly commuting}. This implies that that every \nbd{d}parameter subsemigroup of $T_*$ is GNS-subPS-strongly commuting (in a sense even stronger than Definition \ref{SCSGdefi}) and, therefore, by Theorem \ref{SCSGthm}, admits a (module \nbd{E_0})dilation.

Recall from Appendix \ref{FockAPP} the definition of $L^2(\Om,F)$ for a measure space $\Om$ and a correspondence $F$ as external tensor product $L^2(\Om)\otimes F$. Here, where the measurable space $\C$ (or, sometimes, $\C^n$) is fixed, we rather write $L^2(\mu)$ and the correspondence will be $\cB$. Hence, we will be concerned with $L^2(\mu,\cB)=L^2(\mu)\otimes\cB$ or its higher dimensional analogues $L^2(\mu_1\times\ldots\times\mu_n,\cB)=L^2(\mu_1\times\ldots\times\mu_n)\otimes\cB=L^2(\mu_1)\otimes\ldots\otimes L^2(\mu_n)\otimes\cB$. (Probability measures are finite, hence, \nbd{\sigma}finite.) We may consider $\cB$ as the identity correspondence over $\cB$, turning $L^2(\mu,\cB)$ into a correspondence over $\cB$ with pointwise left action on the ``functions'' in $L^2(\mu,\cB)$. Or we may consider $\cB$ as a Hilbert \nbd{\cB}module, that is, a \nbd{\C}\nbd{\cB}correspondence, turning $L^2(\mu,\cB)$ into a Hilbert \nbd{\cB}module, on which we may define \it{ad hoc} a different left action of $\cB$. Since the elements $f\otimes\U$ generate $L^2(\mu,\cB)$ as Hilbert module, it is enough to define a left action of $\cB$ by indicating the action of all $W_w$ on these elements. (Well-definedness follows by standard arguments, if the defining relations of $CCR(\C,\om)$ are fulfilled under taking matrix elements with elements $f\otimes\U$ from the generating set.)

From convolution semigroups of \it{L\'{e}vy processes} (quantum or not), we may expect a GNS-correspondence of the form $L^2(\mu)\otimes\cB$ as a Hilbert \nbd{\cB}module and an \it{ad hoc} left action of $\cB$. (See the conclusive Remark \ref{LPrem}.) We, indeed, confirm this, by defining a left action via
\beqn{
W_w(f\otimes\U)
~:=~
e^{i2\om(w,\bullet)}f\otimes W_w,
}\eeqn
where $e^{i2\om(w,\bullet)}$ is the multiplication operator determined by the function $z\mapsto e^{i2\om(w,z)}$ (acting on the function $f\colon z\mapsto f(z)$ in $L^2(\mu)$ by sending it to the function $e^{i2\om(w,\bullet)}f\colon z\mapsto  e^{i2\om(w,z)}f(z)$ in $L^2(\mu)$). We claim,
\beqn{
\text{GNS-}T_\mu
~=~
\Bfam{
\sE_\mu:=L^2(\mu)\otimes\cB
,
\xi_\mu:=1\otimes\U
},
}\eeqn
where $1\in L^2(\mu)$ is the constant function $z\mapsto 1$. Clearly, $\AB{\xi_\mu,W_w\xi_\mu}=T_\mu(W_w)$. So, the only thing that remains to be shown is that the elements $W_w\xi_\mu$ generate $L^2(\mu)\otimes\cB$ as Hilbert module. But this follows from $W_w\xi_\mu W_{-w}=e^{i2\om(w,\bullet)}1\otimes\U$ and because the functions $e^{i2\om(w,\bullet)}1=e^{i2\om(w,\bullet)}$ (now, interpreted as elements of $L^2(\mu)$) are total in $L^2(\mu)$.

\brem
One might be tempted to define $\sE_\mu$ as $L^2(\mu)\otimes\cB$ with canonical left action of $\cB$ on itself, and $\xi_\mu$ as the function $\breve{\xi}_\mu\colon z\mapsto W_z$. The problem is that this function is not an element of $L^2(\mu)\otimes\cB$ as soon as the support of $\mu$ is uncountable. (In the approximation by a sequence of finite sums of simple tensors, only countably many elements of $\cB$ occur. But the span of the range of the function $\breve{\xi}_\mu$ is not separable.) It is, in a canonical fashion, an element of the strong closure $L^2(\mu)\sbars{\otimes}\ol{\cB}^s=\ol{L^2(\mu)\otimes\ol{\cB}^s}^s=\ol{L^2(\mu)\otimes\cB}^s$. (In the case of irreducible $\rho$, this is equal to $\sB(L,L^2(\mu)\otimes L$.) The closed \nbd{\cB}subbimodule $\breve{\sE}_\mu:=\cls\cB\breve{\xi}_\mu\cB$ is even a correspondence over $\cB$ (isomorphic to $\sE_\mu$ via bilinear extension of $\breve{\xi}_\mu\mapsto\xi_\mu$). But it does, in general (if $\mu$ is not atomic), not coincide with the subbimodule $L^2(\mu)\otimes\cB$ of $\ol{L^2(\mu)\otimes\cB}^s$.
\erem

We return to the correspondences $\sE_\mu$. Since $(f\otimes W_w)\odot(g\otimes\U)=(f\otimes\U)\odot W_w(g\otimes\U)=(f\otimes\U)\odot (e^{i2\om(w,\bullet)}g\otimes W_w)$, we see that the elements $(f\otimes\U)\odot(g\otimes\U)$ generate $\sE_\mu\odot\sE_\nu$ as a Hilbert \nbd{\cB}module. Obviously, the left action reads
\beqn{
W_w(f\otimes\U)\odot(g\otimes\U)
~=~
(e^{i2\om(w,\bullet)}f\otimes W_w)\odot(g\otimes\U)
~=~
(e^{i2\om(w,\bullet)}f\otimes\U)\odot(e^{i2\om(w,\bullet)}g\otimes W_w)
}\eeqn
Clearly,
\beqn{
u_{\nu,\mu}
\colon
(f\otimes\U)\odot(g\otimes\U)
~\longmapsto~
(g\otimes\U)\odot(f\otimes\U)
}\eeqn
defines an isomorphism $\sE_\mu\odot\sE_\nu\rightarrow\sE_\nu\odot\sE_\mu$ of correspondences. And, clearly,
\beq{ \label{QCprod4}
(\id_\pi\odot~u_{\nu,\mu})(u_{\pi,\mu}\odot\id_\nu)(\id_\mu\odot~u_{\pi,\nu})
~=~
(u_{\pi,\nu}\odot\id_\mu)(\id_\nu\odot~u_{\pi,\mu})(u_{\nu,\mu}\odot\id_\pi)
}\eeq
for all $\mu,\nu,\pi\in\bS_*$. (Both sides act as
\beqn{
(f\otimes\U)\odot(g\otimes\U)\odot(h\otimes\U)
~\longmapsto~
(h\otimes\U)\odot(g\otimes\U)\odot(f\otimes\U)
}\eeqn
from $\sE_\mu\odot\sE_\nu\odot\sE_\pi$ to $\sE_\pi\odot\sE_\nu\odot\sE_\mu)$.) Obviously,
\beq{ \label{QCxisflip}
u_{\nu,\mu}(\xi_\mu\odot\xi_\nu)=\xi_\nu\odot\xi_\mu.
}\eeq

Now the structure $v_{\mu,\nu}\colon\sE_{\mu*\nu}\rightarrow\sE_\mu\odot\sE_\nu$ of the GNS-subproduct system $\sE^\bodot=\bfam{\sE_\mu}_{\mu\in\bS_*}$ (over $\bS_*^{op}=\bS_*$) of $T_*$ is defined as
\beqn{
e^{i2\om(w,\bullet)}1\otimes\U
~=~
W_w\xi_{\mu*\nu}W_{-w}
~\longmapsto~
W_w\xi_\mu\odot\xi_\nu W_{-w}
~=~
(e^{i2\om(w,\bullet)}1\otimes\U)\odot(e^{i2\om(w,\bullet)}1\otimes\U).
}\eeqn
It follows that both paths
\beqn{
\sE_{\mu*\nu}\odot\sE_\pi
~\longrightarrow~
\sE_\mu\odot\sE_\nu\odot\sE_\pi
~\longrightarrow~
\sE_\pi\odot\sE_\mu\odot\sE_\nu
}\eeqn
and
\beqn{
\sE_{\mu*\nu}\odot\sE_\pi
~\longrightarrow~
~~~\sE_\pi\odot\sE_{\mu*\nu}~~~
~\longrightarrow~
\sE_\pi\odot\sE_\mu\odot\sE_\nu
}\eeqn
act as
\bmun{
W_w\xi_{\mu*\nu}W_{-w}\odot W_{w'}\xi_\pi W_{-w'}
~=~
(e^{i2\om(w,\bullet)}1\otimes\U)\odot(e^{i2\om(w',\bullet)}1\otimes\U)
\\
~\longmapsto~
(e^{i2\om(w',\bullet)}1\otimes\U)\odot(e^{i2\om(w,\bullet)}1\otimes\U)\odot(e^{i2\om(w,\bullet)}1\otimes\U)
\\
~=~
W_{w'}\xi_\pi W_{-w'}\odot W_w\xi_\mu W_{-w}\odot W_w\xi_\nu W_{-w}.
}\emun
\begin{subequations} \label{QCSCSG}
This shows the written out version
\beq{ \label{QCSCSGa}
(u_{\pi,\mu}\odot\id_\nu)(\id_\mu\odot u_{\pi,\nu})(v_{\mu,\nu}\odot\id_\pi)
~=~
(\id_\pi\odot v_{\mu,\nu})u_{\pi,\mu*\nu}
}\eeq
of Diagram \eqref{SCSGa}; the proof for the written out version
\beq{ \label{QCSCSGb}
(\id_\nu\odot u_{\pi,\mu})(u_{\nu,\mu}\odot\id_\pi)(\id_\mu\odot v_{\nu,\pi})
~=~
(v_{\nu,\pi}\odot\id_\mu)u_{\nu*\pi,\mu}
}\eeq
of Diagram \eqref{SCSGb} is analogue.
\end{subequations}

\bdefi
A CP-semigroup over an abelian monoid is \hl{fully strongly commuting}\phantomsection\index{fully strongly commuting}\index{CP-semigroup!strongly commuting!fully, over an abelian monoid}\index{strongly commuting!CP-semigroups!fully, over an abelian monoid} if for its GNS-subproduct system there exist isomorphisms $u_{\nu,\mu}\colon\sE_\mu\odot\sE_\nu\rightarrow\sE_\nu\odot\sE_\mu$ fulfilling Equations \eqref{QCprod4}, \eqref{QCxisflip}, and \eqref{QCSCSG}.
\edefi

We, thus, have proved:

\bthm
The Markov semigroup $T_*$ is fully strongly commuting.
\ethm

\bcor
For $k=1,\ldots,d$ let $\bS^k$ be $\N_0$ or $\R_+$. Let $\mu$ be a \nbd{d}parameter convolution semigroup in over $\bS:=\bS^1\times\ldots\times \bS^d$. Then
\begin{enumerate}
\item
Each marginal one-parameter Markov semigroup $T^k=\bfam{T_{\mu_{t_k}}}_{t_k\in\bS^k}$ is fully strongly commuting.

\item
The $T^1,\ldots,T^d$ commute GNS-subPS-strongly.

\item
The \nbd{d}parameter Markov semigroup $\bfam{T_{\mu_{\bt}}}_{\bt\in\bS}$ is fully strongly commuting and admits a strict module \nbd{E_0}dilation.
\end{enumerate}
\ecor

\brem
In \cite[Part II]{Mar03}, Markiewicz examines the Arveson system of the normal continuous time one-parameter Markov semigroup $\bfam{T_{\mu_t}}_{t\in\R_+}$ on $\sB(L)=\ol{\cB}^s$ obtained from a one-parameter convolution semigroup $\mu$ and the unique irreducible continuous representation of $CCR(\C,\om)$. Shalit studied the \nbd{2}parameter case, and proved in \cite[Theorem 3.4.12]{Sha09} that the marginal semigroups (in our words) commute GNS-subPS-strongly (and, therefore, by Shalit \cite{Sha08,Sha08CORR}, admit an \nbd{E_0}dilation).
\erem

\brem \label{LPrem}
It is well known that a weakly continuous convolution semigroup $\bfam{\mu_t}_{t\in\R_+}$ of probability measures on (the locally compact abelian Lie group) $\C$ characterizes (up to stochastic equivalence) a L\'{e}vy process with values in $\C$. By $S_t:=\id*\mu_t$, where
\beqn{
\SB{\id*\mu_t(f)}(z)
~:=~
\int f(z-w)\,d\mu_t(w),
}\eeqn
we define a Markov semigroup $S=\bfam{S_t}_{t\in\R_+}$ on $\cC:=C_b(\C)$ or $\cC:=\wt{C_0(\C)}$. A similar statement is true for all quantum L\'{e}vy processes, and Skeide \cite{Ske05b} points out that the GNS-correspondences always have the form $L^2(\mu_t)\otimes\cB$ with nontrivial left action, while the $L^2(\mu_t)$ form a subproduct system of Hilbert spaces. One may check that also the Markov semigroup $S_*\colon\mu\mapsto\id*\mu$ is fully strongly commuting.

Of course, being indexed by the same monoid, it would be interesting to examine the relation between the semigroups $S_*$ and $T_*$.
\erem

\newpage

\section{Product systems over $\N_0^d$} \label{N0dSEC}

Conditions \eqref{prod4} and \eqref{prod5} in Observation \ref{prodob}, with the many ``time'' indices form quite a ``salad''\!, whose verification in general means quite an effort. In Section \ref{strcomSEC} we have seen that even in the continuous two-parameter case $\R_+^2$, where \eqref{prod4} are vacuous, \eqref{prod5} remain tedious. In the discrete \nbd{d}parameter case $\bS^k=\N_0$, so $\bS=\N_0^d$, the marginals ${E^k}^\odot$ are (up to isomorphism) determined by indicating a single correspondence $E_k:=E^k_1$ (see Observation \ref{d1pob}). This allows to substitute the conditions in Observation \ref{prodob}\eqref{prod4} and \eqref{prod5} with a much simpler set of conditions involving only the $E_k$ and flips among them. It turns out to be convenient, not to reduce the problem to the general Theorem \ref{prodthm}, but to give an independent treatment based directly on Lemma \ref{pi_flem}, that is, based on the essence of Appendix \ref{popAPP}.

We denote by $\be_k:=(0,\ldots,0,1,0,\ldots,0)$ (the $1$ in the $k$th place) the elements of the canonical basis of $\N_0^d$.

\bthm \label{N0dthm}
Let $E_1,\ldots,E_d$ be correspondences over $\cB$ and for $1\le j<i\le d$ let $\sF_{j,i}\colon E_i\odot E_j\rightarrow E_j\odot E_i$ be isomorphisms. Then there exists a product on the family $\bfam{E_1^{\odot n_1}\odot\ldots\odot E_d^{\odot n_d}}_{\bn\in\N_0^d}$ fulfilling $u_{\be_j,\be_i}^*u_{\be_i,\be_j}=\sF_{j,i}$ $(j<i)$ if and only if the $\sF_{j,i}$ fulfill the detailed exchange conditions in \eqref{Tijcond}. \index{product system!over $\N_0^d$}\index{structure!of a product system over a product!the case $\N_0^d$}Moreover, such a product is unique.

\ethm

\proof
Necessity and uniqueness being plain (compare with the discussion leading to Observation \ref{prodob}), we only show sufficiency of \eqref{Tijcond}.

Both well-definedness and associativity of the product, simply amount to the observation that (by Lemma \ref{pi_flem} and Proposition \ref{detexprop}) there are permutations in $S_{m_1+\ldots+m_d+n_1+\ldots+n_d}$ and in $S_{\ell_1+\ldots+\ell_d+m_1+\ldots+m_d+n_1+\ldots+n_d}$ obtainable by admissible flips that, when acting on the following tensor products, put
\beqn{
(E_1^{\odot m_1}\odot\ldots\odot E_d^{\odot m_d})\odot(E_1^{\odot n_1}\odot\ldots\odot E_d^{\odot n_d})
}\eeqn
(for well-definedness) and
\beqn{
(E_1^{\odot\ell_1}\odot\ldots\odot E_d^{\odot\ell_d})\odot(E_1^{\odot m_1}\odot\ldots\odot E_d^{\odot m_d})\odot(E_1^{\odot n_1}\odot\ldots\odot E_d^{\odot n_d})
}\eeqn
(for associativity), into non-decreasing order
\beqn{
E_1^{\odot(m_1+n_1)}\odot\ldots\odot E_d^{\odot(m_d+n_d)}
}\eeqn
and
\beqn{
E_1^{\odot(\ell_1+m_1+n_1)}\odot\ldots\odot E_d^{\odot(\ell_d+m_d+n_d)}
}\eeqn
 respectively, and that (in either case) the permutation is unique.\qed

\lf
The conditions in Observation \ref{prodob}\eqref{prod5} have disappeared; they are automatic. Like in the case $\R_+^2$, in the case $\N_0^2$ the conditions in \ref{prodob}\eqref{prod4} are vacuous. Therefore:

\bcor \label{N02cor}
Let $E_1$ and $E_2$ be correspondences over $\cB$. Then any isomorphism $\sF_{1,2}\colon E_2\odot E_1\rightarrow E_1\odot E_2$ defines a unique product on the family $\bfam{E_1^{\odot n_1}\odot E_2^{\odot n_2}}_{\bn\in\N_0^2}$ fulfilling $u_{\be_1,\be_2}^*u_{\be_2,\be_1}=\sF_{1,2}$. 
\ecor

%%%% BO
\brem
The fact that product system structures on $\bfam{E_1^{\odot n_1}\odot E_2^{\odot n_2}}_{\bn\in\N_0^2}$ correspond to isomorphisms $E_2\odot E_1\rightarrow E_1\odot E_2$ was already observed and used in Solel \cite[Section 4]{Sol06}. Without explicit reference to product systems and related techniques, Power and Solel \cite{PowSo11} examined in depth the non-selfadjoint operator algebras associated to product systems over $\N_0^2$ in the case when $E_1$ and $E_2$ are finite dimensional Hilbert spaces; see \cite[Theorem 5.10]{PowSo11}. One may dig out from \cite{PowSo11} implications regarding the classification of product systems of finite dimensional Hilbert spaces over $\N_0^2$ as easy corollaries. We do not spell this out; we just mention that the necessary and sufficient condition for product systems being isomorphic that follows from their results (what they call \it{product unitary equivalent}) is the same (in their context) as Equation \eqref{flipmorph} below. 
\erem
%%%% EO 

\bob \label{N0dob}
Obviously, if we have vectors $\xi_k\in E_k$ such that $\sF_{j,i}(\xi_i\odot\xi_j)=\xi_j\odot\xi_i$, then the units $\xi_k^{\odot n}$ for the marginal product systems $E_k^{\odot n}$ fulfill the hypotheses of Corollary \ref{uniprodcor}, whether the $\xi_k$ generate the $E_k$ (that is, the CP-maps are mutually strongly commuting), or not.
\eob

Suppose we have two product systems $E^\odot$ and $E'^\odot$ over $\N_0^d$ obtained from correspondences $E^{(')}_1,\ldots,E^{(')}_d$ and flips $\sF^{(')}_{j,i}$. Clearly, if the maps $a_\bn\in\sB^{bil}(E_\bn,E'_\bn)$ form a morphism $E^\odot\rightarrow E'^\odot$, then they are determined by the $a_k:=a_{\be_k}\in\sB^{bil}(E_k,E'_k)$ ($k=1,\ldots,d$). A necessary condition these $a_k$ have to fulfill, is that for all $j<i$
\beq{ \label{flipmorph}
(a_j\odot a_i)\sF_{j,i}
~=~
(a_j\odot a_i)u^*_{\be_j,\be_i}u_{\be_i,\be_j}
~=~
u'^*_{\be_j,\be_i}a_{\be_j+\be_i}u_{\be_i,\be_j}
~=~
u'^*_{\be_j,\be_i}u'_{\be_i,\be_j}(a_i\odot a_j)
~=~
\sF'_{j,i}(a_i\odot a_j).
}\eeq
It requires only a little moments thought to convince ourselves that, basically because the products $u$ and $u'$ are composed out of (amplifications of) flips $\sF_{j,i}$ and $\sF'_{j,i}$, respectively, these conditions are also sufficient. (Compare this with the discussion in the end of Section \ref{compSEC}.)

\bcor \label{N0disocor}
The two product systems $E^\odot$ and $E'^\odot$ are isomorphic if and only if we can find bilinear unitaries $a_k\colon E_k\rightarrow E'_k$ fulfilling \eqref{flipmorph}. This applies in particular, to automorphisms.
\ecor

The considerations maybe reversed, in a sense.

\bcor \label{N0dindcor}
Suppose we have bilinear unitaries $a_k\colon E_k\rightarrow E'_k$. Then bilinear unitaries $\sF_{j,i}$ $(j<i)$ fulfill \eqref{Tijcond} if and only if $\sF'_{j,i}:=(a_j\odot a_i)\sF_{j,i}(a_i^*\odot a_j^*)$ fulfill \eqref{Tijcond}. If one, hence, both are fulfilled, then the product systems $E^\odot$ and $E'^\odot$ are isomorphic via the isomorphism determined by the maps $a_k$.
\ecor

In the following section we use necessity of the conditions in Theorem \ref{N0dthm} to construct a subproduct system over $\N_0^3$ that does not embed into a superproduct system. (Corollaries \ref{N0disocor} and \ref{N0dindcor} will reoccur in the more special situation $E_k=E$ for all $k$, but in the more general context of subproduct systems.) The positive statement in Corollary \ref{N02cor}, instead, we use in Section \ref{EXN02SEC} to construct a product system for every Markov semigroup over $\N_0^2$ on a von Neumann algebra $\cB$, showing this way that all Markov semigroups (and, therefore, all CP-semigroups) over $\N_0^2$ possess a dilation. The dilation obtained, will also provide a whole bunch of counter examples to many desirable properties.

Comparing the special Theorem \ref{N0dthm} with the general Theorem \ref{prodthm}, it is perfectly legitimate to ask, what is responsible for the disappearance in the former of the conditions in \ref{prodob}\eqref{prod5} in the latter. The answer becomes clear if we realize that in Theorem \ref{N0dthm} we also do require by far not all conditions in \ref{prodob}\eqref{prod4} that would occur when we did apply Theorem \ref{prodthm} directly, but actually only a very small part of them. The conditions in \ref{prodob}\eqref{prod4} involve all members of the marginal product systems. But, in the discrete one-parameter case, the marginal product systems are \it{generated} by the 1-elements $E_k$, and in Theorem \ref{N0dthm} we write down only those conditions from \ref{prodob}\eqref{prod4} that involve only their \it{generating} correspondences $E_k$. In fact, all higher $u^{j.i}_{n_j,m_i}$ are to be computed by recursion from $u^{j,i}_{1,1}=\sF_{j,i}$. And both the missing relations in \ref{prodob}\eqref{prod4} and all relations in \ref{prodob}\eqref{prod5} turn out to be fulfilled -- simply, because by Theorem \ref{N0dthm} they actually come from a product system. (It is fair to remark, that the proof of Theorem \ref{N0dthm} makes much fuller use of Appendix \ref{popAPP} (in the form of Lemma \ref{pi_flem}) than the proof of Theorem \ref{prodthm}, where only a small part of the non decreasing partially order preserving permutations occur.)

Of course, having noticed that for Theorem \ref{N0dthm} (the analogue of) the conditions in \ref{prodob}\eqref{prod4} are enough, we may ask if the conditions in \ref{prodob}\eqref{prod5} are really to be required explicitly for Theorem \ref{prodthm} to hold. Or, in case we find $u^{j,i}_{t_j,s_i}$ satisfying \ref{prodob}\eqref{prod4} but not \ref{prodob}\eqref{prod5}, we may ask if we might be able to modify $u^{j,i}_{t_j,s_i}$, preserving ``essential parts of their structure'', to fulfill also \ref{prodob}\eqref{prod5}. Example \ref{4not5ex} below, shows that the answer to the first question is a clear no. Example \ref{4not5ex} also shows that the answer to the second question, is no if we fill the wishy-washy term ``essential parts of their structure'' that have to be preserved with a concrete meaning in the following sense:

\bob
In the situation of Observation \ref{prodob} (that is, product systems ${E^k}^\odot$ $(k=1\ldots,d)$, each over its own $\bS^k$, and bilinear unitaries $u^{j,i}_{t_j,s_i}\colon E^i_{s_i}\odot E^j_{t_j}\rightarrow E^j_{t_j}\odot E^i_{s_i}$ ($1\le j<i\le d$, $s_i\in\bS^i$, $t_j\in\bS^j$), and the question whether or not the $u^{j,i}_{t_j,s_i}$ determine a product system structure on the family $E^{1,\ldots,d}:=\bfam{E_\bt}_{\bt\in\bS:=\bS^1\times\ldots\times\bS^d}$), in which Theorem \ref{prodthm} is formulated: Assume the $u^{j,i}_{t_j,s_i}$ fulfill \ref{prodob}\eqref{prod4}.

Then, by Theorem \ref{N0dthm}, for each $\tau=(\tau_1,\ldots,\tau_d)\in\bS$, the maps  $u^{j,i}_{\tau_j,\tau_i}$ induce the structure of a product system over $\N_0^d$ on the family $E^{1,\ldots,d}_\tau:=\bfam{{E^1_{\tau_1}}^{\!\!\odot n_1}\odot\ldots\odot{E^d_{\tau_d}}^{\!\!\odot n_d}}_{\bn\in\N_0^d}$. We say a change from the $u^{j,i}_{t_j,s_i}$ to $u'^{j,i}_{t_j,s_i}$ (aiming at having the new $u'^{j,i}_{t_j,s_i}$ satisfying not only \ref{prodob}\eqref{prod4} but also \ref{prodob}\eqref{prod5}) is \hl{essentially structure preserving} if for each $\tau$ the product system structures induced on $E^{1,\ldots,d}_\tau$ by $u^{j,i}_{t_j,s_i}$ and by $u'^{j,i}_{t_j,s_i}$ are isomorphic.
\eob

Clearly, with the isomorphisms
\beqn{
u^1_{\tau_1,\ldots,\tau_1}\odot\ldots\odot u^d_{\tau_d,\ldots,\tau_d}
\colon
{E^1_{\tau_1}}^{\!\!\odot n_1}\odot\ldots\odot{E^d_{\tau_d}}^{\!\!\odot n_d}
~\longrightarrow~
E_{\bfam{\tau_1^{n_1},\ldots,\tau_d^{n_d}}}
}\eeqn
(with the \nbd{n_k}fold product $u^k_{\tau_k,\ldots,\tau_k}$ of the marginal system ${E^k}^\odot$) we induce an isomorphic product system structure on the family $\Bfam{E_{\bfam{\tau_1^{n_1},\ldots,\tau_d^{n_d}}}}_{\bn\in\N_0^d}$. This subfamily of the family $E^{1,\ldots,d}$ may be considered as obtained by restricting the indexing set  $\bS$ of the latter to its submonoid $\tau^{\N_0^d}:=\tau_1^{\N_0}\times\ldots\times\tau_d^{\N_0}$ (for each $\tau$ isomorphic to $\N_0^d$). As $\Bfam{E_{\bfam{\tau_1^{n_1},\ldots,\tau_d^{n_d}}}}_{\bn\in\N_0^d}$ is a product system, it has its flips $u^{j,i}_{\tau_{\smash j}^{\,n_{\smash j}},\tau_{\smash i}^{\,m_{\smash i}}}$. If the whole family should be a product system, then, necessarily for $\tau$ and $\tau'$ the $u^{j,i}_{\tau_{\smash j}^{\,n_{\smash j}},\tau_{\smash i}^{\,m_{\smash i}}}$ obtained from $\tau$ and the $u^{j,i}_{{\tau'}_{\smash j}^{\,n'_{\smash j}},{\tau'}_{\smash i}^{\,m'_{\smash i}}}$ obtained from $\tau'$ have to coincide whenever $(\tau^\bm,\tau^\bn)=(\tau'^{\bm'},\tau'^{\bn'})$.

If, starting only from the conditions in \ref{prodob}\eqref{prod4}, the ideas based on this fact should turn out to be enough to prove also \ref{prodob}\eqref{prod5}, then it would be worth to detail these ideas. (For instance, to check if, for fixed $\tau$, the elements $\tau^\mu$ and $\tau^\nu$ ($\mu,\nu\in\N_d$) are compatible in that way, it is sufficient to check if $\tau^\mu$ and $\tau^\nu$ both are compatible with $\tau$, individually. Thinking this over carefully, would allow to settle the rational \nbd{d}parameter case $\Q_+^d$ considering $\Q_+^d$ as inductive limit over $\bfam{\nu\cdot\N_0^d}_{\nu\in\N_0^d}$.) Unfortunately, we will now see that \ref{prodob}\eqref{prod4} is not enough in a very special case.

\bex \label{4not5ex}
Choose a Hilbert space $H$ with $\dim H\ge2$, put $d=2$, and put $E_1:=H=:E_2$, so that $E^{1,2}=\bfam{H^{\otimes n_1}\otimes H^{\otimes n_2}}_{(n_1,n_2)\in\N_0^2}=\bfam{H^{\otimes(n_1+n_2)}}_{(n_1,n_2)\in\N_0^2}$. Define $u^{1,2}_{1,1}:=\f$, the flip on $H^{\otimes 2}$. Define $u^{1,2}_{1,2}:=\id_{H^{\otimes 3}}$. For any other $(n,m)\in\N^2$ choose an arbitrary unitary $u^{2,1}_{n,m}$ in $\sB(H^{\otimes(m+n)})$, and if one of the indices is $0$ then choose the identity of the corresponding tensor power of $H$ (as required from the marginal conditions). Then these $u^{1,2}_{n,m}$ fulfill the conditions in \ref{prodob}\eqref{prod4}, because $d=2$ and for $d=2$ these conditions are vacuous.

By Theorem \ref{N0dthm}, if there is a product system structure on $E^{1,2}$, then it is determined uniquely by $\sF_{1,2}:=u^{1,2}_{1,1}$. In particular, necessarily the flip $u^{1,2}_{2,1}$ would have to equal (by \ref{prodob}\eqref{prod5}!) $(u^{1,2}_{1,1}\otimes\id_H)(\id_H\otimes u^{1,2}_{1,1})=(\f\otimes\id_H)(\id_H\otimes\f)$. This is just the flip $H\otimes H^{\otimes 2}\rightarrow H^{\otimes 2}\otimes H$. Obviously, it does not coincide with our choice $u^{1,2}_{1,2}:=\id_{H^{\otimes 3}}$, showing that flips fulfilling \ref{prodob}\eqref{prod4} need not fulfill \ref{prodob}\eqref{prod5}. Moreover, the identity and the flip $(231)$ on $H^{\otimes 3}$ are, of course, not unitarily equivalent. So no unitary transform of $u^{1,2}_{2,1}$ would ever erase this incompatibility. This means (by Corollaries \ref{N0disocor} and \ref{N0dindcor}), the product system structure on $\bfam{H^{\otimes n_1}\otimes H^{\otimes 2n_2}}_{(n_1,n_2)\in\N_0^2}$ induced by restriction of the product system structure of $E^{1,2}$ coming from $u^{1,2}_{1,1}$ by restriction to the subfamily over $\N_0\times 2\N_0\subset\N_0^2$, is not isomorphic to the product system structure induced by $u^{1,2}_{2,1}$.
\eex

% \OW[ORR (I did not try this. --Orr.)]{P.S.: I think the example may be embedded in $\R_+^d$, showing that without continuity your erratum is substantial. On the other hand I would not be surprised that under continuity conditions it might be possible to find s sufficient condition in terms of the generators of the one parameter semigroups and that this condition has similarity with Theorem \ref{N0dthm}. (I did not try this. --Orr.)}

\newpage

\section[\sc{Examples:} Subproduct systems that do not embed into superproduct systems]{Examples: Subproduct systems that do not embed into superproduct systems} \label{EXsubnsupSEC}

The scope of this section is to exhibit subproduct systems $F^\bodot =\bfam{F_\bn}_{\bn\in\N_0^3}$ over $\bb{N}_0^3$ that cannot be embedded into a superproduct system, in order to obtain examples of Markov semigroups that do not admit any dilation. By exponentiation (Theorem \ref{dpexpthm}), \it{en passant}, we also obtain subproduct systems $\DG^\bodot(F^\bodot)=\bfam{\DG_\bt(F^\bodot)}_{\bt\in\R_+^d}$ over $\R_+^d$ which, by Corollary \ref{embedFcor}, do not embed into a superproduct system. We start with the following example of a subproduct system $G^\botimes = \bfam{G_\bn}_{\bn\in\N_0^3}$ of Hilbert spaces from Shalit and Solel \cite{ShaSo09}. In \cite[Proposition 5.15]{ShaSo09}, it was shown that $G^\botimes$ cannot be embedded into a product system. We, first, show that $G^\botimes$ cannot even be embedded into a superproduct system; this includes a new proof of \cite[Proposition 5.15]{ShaSo09}. Then, analyzing the structure of this example, we provide a far-reaching generalization in Theorem \ref{truncFijthm}. In the end, we input these non-embeddable subproduct systems (and their exponentials) in Theorem \ref{adSPS-CPthm} to obtain 
%%%% BO 
CP-semigroups with GNS subproduct systems that do not embed into a superproduct system. Combining this with Theorem \ref{sdilunithm}, we obtain CP-semigroups with no strong dilation, and in fact with no good dilation (see the discussion following Example \ref{pPSelemex}). Invoking Proposition \ref{Mcharprop} and Theorem \ref{uninonunithm}, we obtain examples of Markov semigroups with no (weak) dilation whatsoever.
%%%% EO

\bex \label{subPSnembex}
Let $\CB{\be_1,\be_2,\be_3}$ denote the standard basis of $\bb{N}_0^3$.  
We set $G_0 = \C$, $G_{\be_i} = \C^2$ and $G_{\be_i +\be_j} = \C^2 \otimes \C^2$ for all $i,j=1,2,3$. 
For all other values of $\bn \in \bb{N}_0^3$, we define $G_\bn = \zero$. 
The coproduct maps $\,w_{\be_i,\be_j} \colon G_{\be_i + \be_j} \rightarrow G_{\be_i} \otimes G_{\be_j}\,$ are defined to be equal to the identity for all $i,j$ except for $i=3, j=2$, while
$w_{\be_3,\be_2} \colon G_{\be_3 + \be_2} \rightarrow G_{\be_3} \otimes G_{\be_2}$ is defined to be the flip $\f\colon x \otimes y\mapsto y \otimes x$, for $x,y \in \C^2$. 
All other coproduct maps $w_{\bm,\bn}$ are defined as they must be in order to (successfully) turn the family $\bfam{G_\bn}_{\bn\in\N_0^3}$ into a subproduct system $G^\botimes$.
\eex

\blem \label{sup->PSlem}
Let $E^\podot$ be a superproduct system over $\N_0^d$. If the product maps $v_{\be_i,\be_j}$ are unitary for all $i,j$, then the superproduct subsystem of $E^\podot$ generated by the family $S_\bn$ with $S_{\be_i}:=E_{\be_i}$ $(i=1,\ldots,d)$ and $S_\bn:=\zero$ $(\bn\notin\CB{\be_1,\ldots,\be_d})$ is a product system.
\elem

\proof
Recall the notation for the $n$th iterated product $v_{\bn_1,\ldots,\bn_n}\colon E_{\bn_1}\odot\ldots\odot E_{\bn_n}\rightarrow E_{\bn_1+\ldots+\bn_n}$ (see Observation \ref{d1pob}). With this, in the big union in Theorem \ref{supintthm} for the computation of the \nbd{\bn}member of the generated superproduct system, the only terms that survive are
\beqn{
E_{\bn_1,\ldots,\bn_n}
~:=~
v_{\bn_1,\ldots,\bn_n}(E_{\bn_1}\odot\ldots\odot E_{\bn_n})
}\eeqn
where $\bn_i\in\CB{\be_1,\ldots,\be_d}$ such that $\bn_1+\ldots+\bn_n=\bn$.

We claim that the $E_{\bn_1,\ldots,\bn_n}$, for fixed $\bn$, do not depend on the decomposition of $\bn$. Note that two compositions of the same $\bn$ can differ at most by a permutation of the occurring $\bn_1,\ldots,\bn_n\in\CB{\be_1,\ldots,\be_d}$.  Indeed, from
\beqn{
v_{\bn_1,\ldots,\bn_n}
~=~
v_{\bn_1,\ldots,\bn_{i-1}+\bn_i,\ldots,\bn_n}(\id^{\odot(i-2)}\odot v_{\bn_{i-1},\bn_i}\odot\id^{\odot(n-i)}),
}\eeqn
with the same equation with $\bn_i$ and $\bn_{i-1}$ exchanged, and with
\beqn{
v_{\bn_{i-1},\bn_i}(E_{\bn_{i-1}}\odot E_{\bn_i})
~=~
E_{\bn_{i-1}+\bn_i}
~=~
v_{\bn_i,\bn_{i-1}}(E_{\bn_i}\odot E_{\bn_{i-1}}),
}\eeqn
we see that
\baln{
E_{\bn_1,\ldots,\bn_{i-1},\bn_i,\ldots,\bn_n}
~=~
&
v_{\bn_1,\ldots,\bn_{i-1},\bn_i,\ldots,\bn_n}(E_{\bn_1}\odot\ldots\odot E_{\bn_{i-1}}\odot E_{\bn_i}\odot\ldots\odot E_{\bn_n})
\\
~=~
&
v_{\bn_1,\ldots,\bn_i,\bn_{i-1},\ldots,\bn_n}(E_{\bn_1}\odot\ldots\odot E_{\bn_i}\odot E_{\bn_{i-1}}\odot\ldots\odot E_{\bn_n})
~=~
E_{\bn_1,\ldots,\bn_i,\bn_{i-1},\ldots,\bn_n}.
}\ealn
So, $E_{\bn_1,\ldots,\bn_n}$ is invariant under next-neighbour transpositions of its indices, hence, under all permutations. Therefore, the superproduct subsystem of $E^\podot$ generated by $S_\bn$ is given by the family $F_\bn:=E_{\bn_1,\ldots,\bn_n}$ for some $\bn_i\in\CB{\be_1,\ldots,\be_d}$ such that $\bn_1+\ldots+\bn_n=\bn$.

We claim, the $F_\bn$ form a product system. Indeed,
%%%% BO
\bmun{
v_{\bm,\bn}(F_\bm\odot F_\bn)
~=~
v_{\bm,\bn}(v_{\bm_1,\ldots,\bm_m}(E_{\bm_1}\odot\ldots\odot E_{\bm_m})\odot v_{\bn_1,\ldots,\bn_n}(E_{\bn_1}\odot\ldots\odot E_{\bn_n}))
 \\
~=~
 v_{\bm_1,\ldots,\bm_m,\bn_1,\ldots,\bn_n}(E_{\bm_1}\odot\ldots\odot E_{\bm_m}\odot E_{\bn_1}\odot\ldots\odot E_{\bn_n})
~=~
F_{\bn+\bm},
}\emun
%%%% EO
so $v_{\bm,\bn}$ is also surjective, hence, unitary.\qed

\bcor
If the subproduct system from Example \ref{subPSnembex} embeds into a superproduct system, then it embeds into a product system.
\ecor

\proof
The only thing to be observed is that the product system constructed in Lemma \ref{sup->PSlem} also contains the members to $\bn=\bn_1+\bn_2$ of the contained subproduct system.\qed

\lf
With this, we are done showing that the subproduct system in Example \ref{subPSnembex} does not embed into a superproduct system, reducing the statement to \cite[Proposition 5.15]{ShaSo09}. But we wish to recover the latter with a new proof, applying our results about the structure of product systems over products of monoids.

Suppose the subproduct system $G^\botimes$ does embed into the product system $H^\otimes$. We may assume that $H^\otimes$ is generated (as a superproduct system) by $G^\botimes$, so $H_{\be_i}=G_{\be_i}=\C^2$ and $H_{\be_i,\be_j}=G_{\be_i,\be_j}=\C^2\otimes\C^2$. The exchange maps of the product system over $\N_0^3$ occurring in Theorem \ref{N0dthm} are
\beqn{
\sF_{j,i}
~:=~
u_{\be_j,\be_i}^*u_{\be_i,\be_j}
~=~
w_{\be_j,\be_i}w_{\be_i,\be_j}^*,
}\eeqn
for $1\le j<i\le 3$. This is,
\beqn{
\sF_{j,i}
~=~
\begin{cases}
\,\f&j=2,i=3,
\\
\,\id&\text{otherwise.}
\end{cases}
}\eeqn
Inserting this in \eqref{Tijcond}, for the only nontrivial possibility $k=1,j=2,i=3$ we get the necessary condition
\beqn{
\id\otimes \f
~=~
\f\otimes\id,
}\eeqn
a condition that, obviously, is not satisfied. This proves \cite[Proposition 5.15]{ShaSo09}. Altogether:

\bitemp[{\protect\cite[Proposition 5.15]{ShaSo09}$'$}.] \label{notEmbedprop}
The subproduct system in Example \ref{subPSnembex} does not embed into a superproduct system.
\eitemp

Example \ref{subPSnembex} can be widely generalized:

\bthm \label{truncFijthm}
Let $E$ be a correspondence over $\cB$ and for $i,j=1,\ldots,d$ let $u_{i,j}$ be unitaries in $\sB^{a,bil}(E\odot E)$. Put $\sF_{j,i}:=u_{j,i}^*u_{i,j}$.
\begin{enumerate}
\item \label{Fij1}
The family $F^\times=\bfam{F_\bn}_{\bn\in\N_0^d}$ with $F_0=\cB$, $F_{\be_i}=E$, $F_{\be_i+\be_j}=E\odot E$, and all remaining $F_\bn=\zero$, is turned into a subproduct system $F^\bodot(U)$ over $\N_0^d$ by defining $w_{\be_i,\be_j}:=u_{i,j}^*$ and all other $w_{\bm,\bn}$ in the only possible way. Moreover:
\begin{enumerate}
\item \label{Fij1a}
Two such subproduct systems $F^\bodot(U)$ and $F^\bodot(U')$ are isomorphic if and only if we can find bilinear unitaries $a_i\colon E\rightarrow E$ such that
\beqn{
(a_j\odot a_i)\sF_{j,i}
~=~
\sF'_{j,i}(a_i\odot a_j)
}\eeqn
for $j<i$.

\item \label{Fij1b}
$F^\bodot(U)$ is isomorphic to $F^\bodot(T)$ where $T$ has the \hl{upper triangular form} $t_{i,j}=\id_{E\odot E}$ for $1\le j\le i\le d$ and $t_{j,i}=\sF_{j,i}$ for $1\le j<i\le d$.

\end{enumerate}

\item \label{Fij2}
$F^\bodot(U)$ embeds into a superproduct system ($E^\podot$, say) if and only the $\sF_{j,i}$ satisfy Equation \eqref{Tijcond}. If so, then $F^\bodot(U)$ also embeds into a product system (namely, into the superproduct subsystem of $E^\podot$ generated by $F^\bodot(U)$).
\end{enumerate}
\ethm

\proof
They main part of \eqref{Fij1} is clear, and \eqref{Fij2} follows from the preceding discussion. The parts that require work are \eqref{Fij1a} and \eqref{Fij1b}.

The only freedom for an isomorphism from $F^\bodot(U)$ to $F^\bodot(U')$ consists in bilinear unitary maps $a_i\colon E\rightarrow E$ (thought of as $ F_{\be_i}\rightarrow F'_{\be_i}$) and $a_{i,j}\colon E\odot E\rightarrow E\odot E$ (thought of as $ F_{\be_i+\be_j}\rightarrow F'_{\be_i+\be_j}$, so that, necessarily, $a_{i,j}=a_{j,i}$). For being an isomorphism, these maps have to satisfy
\beqn{
(a_i\odot a_j)w_{\be_i,\be_j}
~=~
w'_{\be_i,\be_j}a_{i,j},
}\eeqn
and nothing else. So, given $F^\bodot(U)$ and $F^\bodot(U')$, there exists an isomorphism between them if and only if we can find $a_i$ such that the matrix $w'^*_{\be_i,\be_j}(a_i\odot a_j)w_{\be_i,\be_j}$ is symmetric. Putting $a_{i,j}:=w'^*_{\be_i,\be_j}(a_i\odot a_j)w_{\be_i,\be_j}=w'^*_{\be_j,\be_i}(a_j\odot a_i)w_{\be_j,\be_i}=:a_{j,i}$, this is the same as the condition in \eqref{Fij1a}.

We can also change the point of view and say that, given $F^\bodot(U)$ and choosing arbitrary $a_i$ and $a_{i,j}=a_{j,i}$, then these define an isomorphism from $F^\bodot(U)$ to $F^\bodot(U')$, if we choose $U'$ to satisfy $u'^*_{i,j}=w'_{\be_i,\be_j}=(a_i\odot a_j)w_{\be_i,\be_j}a_{i,j}^*$. To prove \eqref{Fij1b}, we put $a_i:=\id$ and choose $a_{i,j}=a_{j,i}$ to satisfy $w'_{\be_i,\be_j}=(a_i\odot a_j)w_{\be_i,\be_j}a_{i,j}^*=w_{\be_i,\be_j}a_{i,j}^*=\id_{E\odot E}$ for $j\le i$, that is, $a_{i,j}=a_{j,i}:=w_{\be_i,\be_j}$. Then, for $j<i$, we have $w'_{\be_j,\be_i}=w_{\be_j,\be_i}a_{i,j}^*=w_{\be_j,\be_i}w_{\be_i,\be_j}^*=\sF_{j,i}$, as claimed.\qed

\brem
Of course, embeddability of a subproduct system into a (super)product system does not change under isomorphism. Nevertheless, as in Corollary \ref{N0dindcor}, one shows by hand that the condition in \eqref{Tijcond} is invariant under the transformation $\sF_{j,i}\to\sF'_{j,i}:=(a_j\odot a_i)\sF_{j,i}(a_i^*\odot a_j^*)$.
\erem

\bex \label{subPSnotembedc}
By Theorem \ref{dbexpthm} and Corollary \ref{embedFcor},  the non-embeddable subproduct systems $F^\bodot$ in Theorem \ref{truncFijthm} (including Example \ref{subPSnembex}) gives rise to a continuous time ordered subproduct system $\DG^\bodot(F^\bodot) = \bfam{\DG_\bt(F^\bodot)}_{\bt\in\R_+^3}$ over $\R_+^3$ which cannot be embedded into a superproduct system. 
\eex

% \bob \label{trunPSisoob}
% \OW{CHECK!}
% The results in Theorem \ref{truncFijthm} can be used to understand the isomorphism problem for product systems over $\N_0^d$ with isomorphic marginals. In fact:
% \begin{enumerate}
% \item
% Every product system $E^\odot$ over $\N_0^d$ with (pairwise) isomorphic marginals gives rise to a subproduct system as in Theorem \ref{truncFijthm} by truncating $E^\odot$ to elements $E_\bn$ with degree of $\bn$ not bigger than $2$. Moreover, $E^\odot$ is generated by this subproduct subsystem.

% \item
% The structure of the product system is determined by the exchange operators. In other words, the structure of the product system is determined by the structure of the contained subproduct system.
% \end{enumerate}
% Therefore, two product system over $\N_0^d$ are isomorphic if and only if their exchange operators satisfy the conditions in Theorem \ref{truncFijthm}\eqref{Fij1a}.
% \eob

\lf
We now pass to the CP-semigroups constructed in Theorem \ref{adSPS-CPthm}.

Choose any \nbd{d\times d}matrix $U$ with unitary entries $u_{i,j}$, and consider the
%%%% BO 
(adjointable) 
%%%% EO 
subproduct system $F^\bodot:=F^\bodot(U)$ as in Theorem \ref{truncFijthm}\eqref{Fij1}. Following the construction of a CP-semigroup for $F^\bodot$ from Theorem \ref{adSPS-CPthm}, we put
\beqn{
F
~:=~
\bigoplus_{\bn\in\N_0^d}F_\bn
~=~
\cB\oplus\Bfam{\underbrace{\bigoplus_{1\le i\le d}F_{\be_i}}_{=:F_\delta}}\oplus\Bfam{\underbrace{\bigoplus_{1\le j\le i\le d}F_{\be_i+\be_j}}_{=:F_{\delta\times\delta}}},
}\eeqn
and we define $v_\bn\colon F\odot F_\bn\rightarrow F$ by $v_\bn(x_\bm\odot y_\bn):=w_{\bm,\bn}^*(x_\bm\odot y_\bn)$ $(=:x_\bm y_\bn)$. In order to determine $T_\bn(a):=v_\bn(a\odot\id_\bn)v_\bn^*$, it is convenient to decompose $a\in\sB^a(F)$ as
\beqn{
a
~=~
\Matrix{
a_{0,0}&a_{0,\delta}&a_{0,\delta\times \delta}
\\
a_{\delta,0}&a_{\delta,\delta}&a_{\delta,\delta\times \delta}
\\
a_{\delta\times \delta,0}&a_{\delta\times \delta,\delta}&a_{\delta\times \delta,\delta\times \delta}
}
\in
\sB^a
\Matrix{
F_0
\\
F_\delta
\\
F_{\delta\times \delta}
}
.
}\eeqn
Of course, $T_0=\id_{\sB^a(F)}$ and $v_\bn=0$ for $n\ge3$ and $\bn=\bn_1+\ldots+\bn_n$ ($\bn_i\in\CB{ \be_1,\ldots,\be_d}$), so $T_\bn=0$ for such $\bn$. It remains to calculate $T_{\be_i}$ and $T_{\be_i+\be_j}$.

In general, $v_\bn^*$ is $0$ on $(v_\bn F)^\perp$, so it is sufficient to compute $T_\bn(a)$ on elements of the form $v_\bn(x\odot y_\bn)$. We find
\beqn{
T_\bn(a)(v_\bn(x\odot y_\bn))
~=~
v_\bn(ax\odot y_\bn).
}\eeqn

If $\bn=\be_i+\be_j$, then only $x_0\in F_0$ gives a non-zero contribution to $x\odot y_\bn$, so, from $a$ only the first column survives. Likewise from $ax$ only the component of $ax_0$ in $F_0$ survives in $(ax_0)\odot y_\bn$. Therefore, in this case $T_\bn(a)=T_\bn(a_{0,0})$, so,
\beqn{
T_\bn(a_{0,0})(v_\bn(x_0\odot y_\bn))
~=~
v_\bn(a_{0,0}x_0\odot y_\bn)
~=~
a_{0,0}v_\bn(x_0\odot y_\bn),
}\eeqn
so, $T_\bn(a)=a_{0,0}p_\bn$, where $p_\bn=v_\bn v_\bn^*$ is the projection onto $v_\bn (F\odot F_\bn)=F_\bn$.

Now suppose $\bn=\be_i$. Then from $a$ only the block $\rtMatrix{a_{0,0}&a_{0,\delta}\\a_{\delta,0}&a_{\delta,\delta}}$ contributes to $T_\bn(a)$, while the result $T_\bn(a)$ has the form
\beqn{
\Matrix{
0&&
\\
&T_\bn(a_{0,0})&T_\bn(a_{0,\delta})
\\
&T_\bn(a_{\delta,0})&T_\bn(a_{\delta,\delta})
}
.
}\eeqn
The range of $v_\bn$ is $v_\bn((F_0\oplus F_\delta)\otimes F_\bn)=F_\bn\oplus\bfam{\bigoplus_{1\le j\le d}F_{\be_j+\bn}}$. For computing $T_\bn(a)$ it is convenient decompose $a$ further as
\beqn{
a
~=~
\Matrix{
a_{0,0}&a_{0,1}&\ldots&a_{0,d}&
\\
a_{1,0}&a_{1,1}&\ldots&a_{1,d}&
\\
\vdots&\vdots&&\vdots&
\\
a_{d,0}&a_{d,1}&\ldots&a_{d,d}&
\\
&&&&0
}
.
}\eeqn
Then
\beqn{
T_\bn(a)
~=~
\Matrix{
0&&&&
\\
&a_{0,0}&a_{0,1}\odot\id&\ldots&a_{0,d}\odot\id
\\
&a_{1,0}\odot\id&u_{1,i}(a_{1,1}\odot\id)u_{1,i}^*&\ldots&u_{1,i}(a_{1,d}\odot\id)u_{d,i}^*
\\
&\vdots&\vdots&&\vdots
\\
&a_{d,0}\odot\id&u_{d,i}(a_{d,1}\odot\id)u_{1,i}^*&\ldots&u_{d,i}(a_{d,d}\odot\id)u_{d,i}^*
}
.
}\eeqn
(Recall that $\bn=\be_i$!) Without making special choices for $u_{i,j}$, there is not much more that can be said about $T$, but that $T$ does not admit a dilation if (and only if) the $u_{i,j}$ violate the condition in Theorem \ref{truncFijthm}\eqref{Fij2}. We may pass to the upper triangular form of $U$. (Note that, by Observation \ref{Tembedob}, this just means passing to a conjugate CP-semigroup.) Then only the $u_{j,i}$ in the matrix for which $j<i$ remain (possibly) non-trivial. (Note that the size of this block depends on $i$ of $\bn=\be_i$.) If all $u_{j,i}$ are trivial (so that $T$ admits a dilation), then $T_\bn(a)$ is just $(pap)\odot\id_\bn$, where $p$ is the projection onto $\cB\oplus F_\delta$.

\bex \label{nondilCPex}
Let us return to our Example \ref{subPSnembex}. So, $d=3$, $E=\C^2$, and only $u_{2,3}=\f$ is non-trivial. For $i=1,2$ we are in the situation where all $u_{j,i}$ occurring in the matrix are trivial. For $i=3$ we have
\beqn{
T_\bn(a)
~=~
\Matrix{
0&&&&
\\
&a_{0,0}&a_{0,1}\otimes\id&a_{0,2}\otimes\id&a_{0,3}\otimes\id
\\
&a_{1,0}\otimes\id&a_{1,1}\otimes\id&(a_{1,2}\otimes\id)\f&a_{1,3}\otimes\id
\\
&a_{2,0}\otimes\id&\f(a_{2,1}\otimes\id)&\f(a_{2,2}\otimes\id)\f&\f(a_{2,3}\otimes\id)
\\
&a_{3,0}\otimes\id&a_{3,1}\otimes\id&(a_{3,2}\otimes\id)\f&a_{3,3}\otimes\id
}
.
}\eeqn
Despite this example looks to some extent like the simplest possible, it appears that the action on concrete matrices of operators on a \nbd{1+3\times 2+6\times 4=31}dimensional Hilbert space looks everything but easy to understand.
\eex

%%%% BO
The main result of \cite{ShaSk11} asserts that there are three commuting Markov maps such that the Markov semigroup over $\N_0^3$ generated by them has no dilation acting on a type I factor. For the above example we may argue as indicated in the opening paragraph of this section, and, using \ref{notEmbedprop}, we find the following strengthening of \cite[Theorem 1.5]{ShaSk11}. 

\bcor \label{noDilcor}
There are three commuting CP maps (namely, $T_{\be_i}$ from Example \ref{nondilCPex}) on $M_{31}$ (the algebra of \,$31 \times 31$ complex matrices) such that the CP-semigroup over $\N_0^3$ generated by them has no good dilation, and in particular no strong dilation. Moreover, there are three commuting Markov maps on $\wt{M_{31}}=M_{31}\oplus\,\C\wt{\U}$ such that the Markov semigroup over $\N_0^3$ generated by them has no (weak) dilation whatsoever. 
%Let $\cB = M_{31}$ (the algebra of $31 \times 31$ complex matrices). There are three commuting CP maps on $\cB$ (namely, $T_{\be_i}$ from Example \ref{nondilCPex}) such that the CP-semigroup over $\N_0^3$ generated by them has no good dilation, and in particular no strong dilation. Moreover, there are three commuting Markov maps on $\wt{\cB}=\cB\oplus\C\wt{\U}$ such that the Markov semigroup over $\N_0^3$ generated by them has no (weak) dilation whatsoever. 
\ecor

\proof
If a (strict or normal) CP-semigroup on $\sB^a(E)$ has good dilation, then its GNS-sub\-product system $\sF^\bodot$ embeds into a superproduct system $F^\podot$ (of course, of \nbd{\sB^a(E)}corre\-spond\-ences). By the consideration anticipating Example \ref{nonstriPSex}, (independently of strictness or normality) the passage from $F_t$ to $E_t:=E^*\odot F_t\odot E$ gives a superproduct systems of \nbd{\cB}correspondences, and since this passage respects inclusions, this superproduct system contains the subproduct system (now, strictness or normality really is required) of \nbd{\cB}correspondences associated with $T$.

Since the (clearly both normal and strict) CP-semigroup emerging, as indicated in the corollary, from Example \ref{nondilCPex} (where $M_{31}=\sB^a(E)$ for the Hilbert \nbd{\C}module $E=\C^{31}$) has a subproduct system of \nbd{\cB}correspondences (with $\cB=\C$) which by Proposition \ref{notEmbedprop} does not embed into a superproduct system, this CP-semigroup cannot have a good dilation. Its unitalization cannot have any dilation, because otherwise by Theorem \ref{uninonunithm} the original CP-semigroup admits a strong (hence, good) dilation.\qed

\lf
As mentioned in Example \ref{subPSnotembedc}, by exponentiating a non-embeddable subproduct system over $\N_0^3$ as in Theorem \ref{dbexpthm}, we obtain by Corollary \ref{embedFcor}  a subproduct system over $\R_+^3$ that does not embed into a superproduct system. Using the same reasoning as above, we obtain the following continuous time version. 

\bcor \label{noDilContcor}
There is a CP-semigroup over $\R_+^3$ that has no good dilation, and in particular no strong dilation. Moreover, there exists a Markov semigroup over $\R_+^3$ that has no (weak) dilation whatsoever.
\ecor

\newpage

\section{Topology enters}\label{topSEC}

Of course, all our objects that we looked at so far, were \nbd{C^*}objects and the maps on or between them were bounded (mostly by automatic continuity, because we were looking at \nbd{C^*}objects). But, we successfully managed to avoid almost every reference to a topology different from norm topology. A prominent exception is the discussion on Section \ref{CPspsSEC} about \it{strict} CP-semigroups on $\sB^a(E)$. (This discussion includes the fact that strict \nbd{E}semigroups on $\sB^a(E)$ have a product system, not just a superproduct system, and leads to the statement of Theorem \ref{Markmodthm} that Markov semigroups over the opposite of an Ore monoid have a \it{strict} module dilation if and only if the GNS-subproduct system embeds into a product system.) However, also strictness appeared only in form of (equivalent) nondegeneracy conditions. The crucial Proposition \ref{convprop} (explaining that full dilations and module dilations are essentially the same thing) is phrased in entirely algebraic terms by characterizing $\sB^a(E)$ as double centralizers of the pre-\nbd{C^*}algebra $\sF(E)$ (another automatic continuity, provided that $E$ is a Hilbert module).

We have now reached our limits in proceeding with the discussion of general dilations $(\cA,\theta,p)$ without introducing a topology in which the individual maps $\theta_t$ are continuous (for instance, identifying the algebra $\cA$ of a full dilation as a multiplier algebra, induces on $\cA$ a natural strict topology).

The scope of this Section is, mainly, two-fold. The first basic scope is to examine to what extent the \nbd{C^*}framework can be pushed further. For instance, we attempt to look how much of $\sB^a(E)$ (and its natural strict topology) in the full case, survives to the general case. If we have a pair $(\cA,p)$ of a \nbd{C^*}algebra $\cA$ and a projection $p$, then, by Proposition \ref{convprop}\eqref{cp1}, we have  a Hilbert module $E:=\cA p$ over $\cB:=p\cA p$ and, always, $\cA$ contains $\sF(E)=\ls\cA p\cA$. But, does $\cA$ contain $\sB^a(E)=M(\sF(E))$? We give a precise (and, somehow, unique) meaning to this question in Definition \ref{contdefi}, and prove in Theorem \ref{contthm} that the answer is yes if and only if $\cA=\sB^a(E)\oplus\sF(E)^\perp$. The problem how to define a suitable topology on the direct summand $\sF(E)^\perp$ remains, however, unsolved -- in the \nbd{C^*}case, at least. And this is the second basic scope, namely, illustrating how passing to the von Neumann case (with its natural topologies) sorts out the problems that have a topological origin. For instance, in Theorem \ref{vNcontthm}\eqref{vNct1}, we show that if $\cA$ is a von Neumann algebra, then $\cA$ always contains $\sB^a(E)$. The passage to the von Neumann case removes problems of a topological nature and leaves us with the true structural problems. We examine some of them (in the \nbd{C^*}setting under the assumption that $\cA$ contains $\sB^a(E)$) in the end of the section.

But before starting properly, let us insert two considerations (in between $\bullet\bullet\bullet$s). One is another strong reason for why we need a weaker topology than the norm topology, based on continuity of $\theta_t$ with $t$. (We said, we do not tackle real problems about this sort of continuity in these notes. But the few statements in several remarks about what \bf{can} be done, get a precise meaning, using the notation we introduce here.) This first one can be entirely skipped. The other one is a very quick review about von Neumann modules, aiming at the reader who knows \nbd{W^*}modules, and to give the reader who knows none of them at least a bit of a hold \it{without} having to work through Appendix \ref{vNAPP} before.

\bulletline
There is another reason why we need a topology on $\cA$ more suitable than the norm topology. Although, in these notes, we do not address questions of continuity of $T$ or $\theta$ with respect to $t$, in the end, we surely wish that CP-semigroups $T$ that are continuous in $t$ in some topology, will have dilating \nbd{E}semigroups that are continuous in some topology, too. For continuity of operators $\theta_t$ on the Banach space $\cA$ with $t\in\bS$, without any further structure on $\cA$ (for instance, being a multiplier algebra like $\sB^a(E)=M(\sK(E))$), there are only \phantomsection\hl{strong Banach continuity}\index{semigroup!continuous!Banach (strongly$=$weakly)} ($t\mapsto\theta_t(a)$ continuous for all $a\in\cA$) and \hl{weak Banach continuity} ($t\mapsto\vp\circ\theta_t(a)$ continuous for all $a\in\cA,\vp\in\cA^*$). However:

\bob \label{pstrcontob}
For $\bS=\R_+$, suppose $p=\theta_t(p)$ for all $t$. Then
\beqn{
T_t(b)
~=~
p\theta_t(b)p
~=~
\theta_t(pbp)
~=~
\theta_t(b).
}\eeqn
In other words, $\theta$ restricts to an \nbd{E_0}semigroup $T$ on $p\cA p$, that is, $T$ dilates itself. Therefore, whenever $T$ is a nonendomorphic Markov semigroup, $\theta_t(p)$ cannot be constant. However, the nonconstant family of projections $\theta_t(p)$ is not norm continuous. In other words, $\theta$ is not strongly Banach continuous. In the one-parameter case, it cannot even be weakly Banach continuous, because the latter implies the former; see Engel and Nagel \cite[Theorem I.1.6]{EnNa06}.%
\footnote{
Regarding questions of measurability of one-parameter semigroups, still today there is no better reference than Hille and Phillips \cite{HiPhi57}. (For instance, we do not know results that did not go eventually back to \cite[Theorem 10.5.5]{HiPhi57}, about what beyond measurability for $t>0$ it needs to make a semigroup continuous at $t=0$.) But regarding continuity, we recommend the compact text Engel and Nagel \cite{EnNa06}, which contains everything one might like to need easily detectable within the first 37 pages.
}
\eob

% When we leave the situation of full (or module) dilations, we have to face the problem to find out how much of the good behaviour of module dilation is preserved. A question of particular interested is which topology can play the role of the strict topology. In the sequel, we go a bit into this direction in the \nbd{C^*}framework put up until now. Quickly, it gets nasty, or we are not able to find answers to the natural question we will be asking. The reader will appreciate to know that larger parts of these problems will disappear when we switch to \it{von Neumann algebras}, \it{von Neumann modules}, and \it{von Neumann correspondences}. Sometimes, the technical discussion will be a bit more involved, and generally the spirit of discussions in the von Neumann case is a bit more operator theoretic than the \nbd{C^*}case; a price we gladly will pay looking at what we gain.

Before going into discussing specific problems, let us review some general aspects of von Neumann modules, enough so -- we hope -- to allow the reader knowing von Neumann algebras to go through the remainder of the section without having to consult all the time Appendix \ref{vNAPP}. (To some extent this concerns also the just stated observation about continuity.) We choose to work with \it{von Neumann algebras}, not with \phantomsection\hl{\nbd{W^*}algebras}\index{von Neumann!W@\nbd{W^*}} (that is, \nbd{C^*}algebras with a \it{pre-dual}). A \phantomsection\hl{von Neumann algebra}\index{von Neumann!algebra|bf}, for us, is a \nbd{*}invariant subset $\cA\subset\sB(H)$ (where $H$ is some Hilbert space) such that $\cA$ coincides with its \it{double commutant} $\cA''$. (Equivalently, $\cA$ is a strongly closed \nbd{*}algebra acting nondegenerately on $H$.) In particular, a von Neumann algebra is always a concrete operator algebra and comes, therefore, shipped with the Hilbert space $H$ on which it acts. Also a \it{von Neumann module} $E$ over a von Neumann algebra $\cB\subset\sB(G)$ is always a module of operators; either it is given immediately as a subset of $E\subset\sB(G,H)$ (the same $G$, where $\cB$ acts), a so-called \phantomsection\hl{concrete von Neumann module}\index{von Neumann!module!concrete|bf}\index{concrete!von Neumann module|see {von Neumann}} (Skeide \cite{Ske06b}), or it is just a \phantomsection\hl{von Neumann module}\index{von Neumann!module}\index{module!von Neumann|see {von Neumann}} (Skeide \cite{Ske00b}) that can be turned into a concrete von Neumann module; it is just important to recall that, even for pre-Hilbert modules over pre-\nbd{C^*}algebras, once the representation of the algebra is chosen, there remains essentially (that is, up to suitable unitary equivalence) exactly one way to turn a module into an operator module, namely, the procedure in Footnote \ref{StineFN}. There are several (at least three) equivalent ways to assure that $E$ is a von Neumann module (one is strong closure). Likewise, a (\phantomsection\hl{concrete}) \hl{von Neumann correspondence}\index{von Neumann!correspondence}\index{correspondence!von Neumann|see {von Neumann}}\index{von Neumann!correspondence!concrete|bf}\index{concrete!von Neumann correspondence|see {von Neumann}} is a (concrete) von Neumann module with a nondegenerate left action of another von Neumann algebra that is normal; also here there are several equivalent possibilities to say when a left action is normal. We refer the reader to Appendix \ref{vNAPP} for a more detailed discussion of von Neumann modules.

The topologies we use on von Neumann algebras (and their analogues for von Neumann modules) are the \it{strong topology} (also known as \it{strong operator topology}) and \it{normality} of positive maps. (The latter is not actually continuity in a topology; it is order continuity, and as such equivalent to \nbd{\sigma}weak continuity. See also here Appendix \ref{vNAPP} for details about other topologies and some reasons for our choice.) The only thing going beyond Appendix \ref{vNAPP}, is that for the \it{weak topology} (also known as \it{weak operator topology}) it is sometimes useful to recall that a \phantomsection\hl{weakly operator continuous}\index{semigroup!continuous!von Neumann (strongly$=$weakly)} semigroup (that is, $t\mapsto\AB{x,T_t(b)y}$ is continuous for all $b\in\cB\subset\sB(G)$ and $x,y\in G$) of normal CP-maps is \hl{strongly operator continuous} (that is, $t\mapsto T_t(b)x$ is continuous for all $b\in\cB\subset\sB(G)$ and $x\in G$), too.%
\footnote{
In many papers it is claimed, possibly presenting the easy argument for strong \bf{right} continuity of $T$, that this statement is obvious. But left continuity is tricky and has a different more involved proof. The first formal (and direct) proof seems to be in Markiewicz and Shalit \cite{MaSha10}. A proof by first showing that the minimal dilation is weakly continuous, and then appealing to the fact that for \nbd{E}semigroups the implication is obvious and turns over easily to the dilated CP-semigroup, can be found in Skeide \cite[Appendix A.2]{Ske16}
}
In general, we will speak of \hl{strongly continuous} semigroups both for the \nbd{C^*}case and for the von Neumann case, meaning strongly Banach continuous in the former and strongly operator continuous in the latter. This is more or less unambiguous, because by the result \cite{Ell00} of Elliot a strongly Banach continuous normal CP-semigroup on a von Neumann algebra is uniformly continuous; so it makes not much sense to speak of strongly Banach continuous semigroups in the von Neumann case. However, if $E$ is a von Neumann module, then it \bf{does} make sense to distinguish between the strong topology of the space $\sB^a(E)$ of operators on the Banach space $E$ and the strong topology of the subspace $\sB^a(E)\subset\sB(H)$ of the space of operators $\sB(H)$ on the Hilbert space $H$. To minimize confusion, for semigroups on $\sB^a(E)$ we will always (also in the \nbd{C^*} case) speak of \phantomsection\hl{strictly continuous}\index{semigroup!continuous!BaE@$\sB^a(E)$ (strongly$=$strictly)} when referring to the strong operator topology of $\sB^a(E)$. (This is justified by the fact that on $\sB^a(E)$ the strong and the strict topology coincide on bounded subsets and that strong continuity of $t\mapsto T_t(a)$ \bf{is} a matter of bounded subsets.) We speak of \phantomsection\hl{strongly continuous}\index{semigroup!continuous!BaEvN@$\sB^a(E)\subset\sB(H)$ (strongly$=$weakly)} when referring to the strong operator topology of $\sB(H)\supset\sB^a(E)$ when $E$ is a von Neumann module.

Now in the situation of Observation \ref{pstrcontob} there is no longer a problem: Also non-constant $\theta_t(p)$ may be strictly continuous, or strongly continuous in the von Neumann case.

\bulletline
The only situation where we could pose a reasonable continuity condition (going beyond boundedness) on the individual maps $\theta_t$ of a dilation $(\cA,\theta,p)$, was strictness in the case when the dilation is full. And, as explained in Example \ref{EPSex} (taking into account also Observation \ref{fu-mod-ob}), in the case of full dilations, strictness implies that the superproduct system of the dilation according to Theorem \ref{E-supPSthm} is actually a product system. That we were able to speak about strictness depended very much on how $p$ sits inside $\cA$; it is a property relative to $p$ that may change if we dilate another CP-semigroup to the same $\cA$ but with different $p$, and $\theta$ may be strict for one $p$ but not for the other.

When $\cA\subset\sB(H)$ is a von Neumann algebra, this has the enormous advantage that the strong topology does not depend on the choice of other ingredients of the dilation; what is a normal map (unlike the strong topology) does not even depend on the concrete identification of $\cA$ as a subset of $\sB(H)$ but only on the \nbd{*}algebra structure (more precisely, on the order structure, which is determined by the algebraic structure) of $\cA$. So all subsets of $\cA$ (like $\cA p$ and its submodules, or $\ls\cA p\cA$, which occurred in several circumstances), inherit their own strong topologies (and all other topologies a von Neumann algebra has). And the requirement of the \nbd{E}semigroup $\theta$ to be normal, is independent of whatever $p$ might be.

In the \nbd{C^*}case we know that strictness of $\theta$ for a full dilation grants that the superproduct system of a dilation is a product system; for non-full dilations, we do not even know what strict is. Now, in the von Neumann case, we can raise for arbitrary dilations the question, whether normality grants that the super product system of the dilation is a product system. And, while we promised that, usually, passing to the von Neumann case will sort out problems caused by missing topological properties, we start right with an example of a problem (a normal dilation, having a proper superproduct system) which, before passing to von Neumann algebras, could not even be formulated when the notion of strictness was missing:

\bob \label{hypevNob}
All the modules that occur in Example \ref{hypexex} are von Neumann modules; all algebras are von Neumann algebras. (This always happens, when $\cB$, hence, also $\wt{\cB}$ are finite-dimensional like $\C$ and $\wt{\C}=\C^2$.) All \nbd{E}semigroups are normal. So, by Example \ref{hypexex} (and also Example \ref{discex}) we get an example also for the von Neumann case that the superproduct system of a dilation may be proper, even if the \nbd{E}semigroup is normal. The \nbd{E}semigroup in Example \ref{hypexex} is even strongly continuous, so that also continuity with time does not sort out the problem.
\eob

The definition of full dilation (part of Definition \ref{dildef}) and Proposition \ref{convprop} (that gives meaning to parts of the definition) are not affected by whether $\cA$ is a von Neumann algebra or not. Just, if $\cA\subset\sB(H)$ \bf{is} a von Neumann algebra and $G:=pH$ the Hilbert subspace of $H$ onto which $p$ projects, then $E:=\cA p$ is a von Neumann module over $\cB:=p\cA p\subset\sB(G)$. More precisely, $E\subset\sB(G,H_p)$ with $H_p:=\ol{\cA G}$ is a concrete von Neumann module. (The only question, strong closure, is proved by exactly the same line as the proof of Proposition \ref{convprop}\eqref{cp1}, just replacing everywhere norm closure with strong closure.) So, apart from the extra information that $E$ is actually a von Neumann module provided $\cA$ is a von Neumann algebra, Proposition \ref{convprop} and the definition of full dilation remain really untouched.

$\sB^a(E)$ is the multiplier algebra of $\sF(E)$; and this is independent of whether $E$ is a Hilbert module or even a von Neumann module. But it means that (the unit ball of) $\sF(E)$ is strictly dense in (the unit ball of) $\sB^a(E)$. If $E$ is a von Neumann module, then this implies the same statement for the strong operator topology inherited from the containing $\sB(H)$ (or $\sB(H_p)$ in the case $E=\cA p$). For any dilation, $\cA$ contains $\sF(E)=\ls\cA p\cA$; but does $\cA$ contain also $\sB^a(E)$? If $\cA$ is a von Neumann algebra, the answer is yes. ($\cA$ is strongly closed, and a strictly converging net, converges \it{a fortiori} strongly in the strong topology of the von Neumann algebra $\cA\subset\sB(H)$.) The general answer is based on the following simple facts:

\bob
~

\begin{enumerate}
\item
If $\cI$ is a (not necessarily closed) \nbd{*}ideal in a \nbd{C^*}algebra $\cA$, then with $\cI^\perp:=\CB{a\in\cA\colon ia=0~(i\in\cI)}$ we get an essential ideal $\cI+\cI^\perp\cong\cI\oplus\cI^\perp$ in $\cA$. (The only not fairly obvious statement about the norm of $i+j$ ($i\in\cI$ and $j\in\cI^\perp$) follows by analyzing the limit of $\bnorm{\frac{i+j}{\max\CB{\norm{i},\norm{j}}}}^{2n}$.)

\item
If $\cI$ is essential in $\cA$, then the multiplier algebra $M(\cI)$ contains $\cA$ in a unique way. (Why? Because $\cI$ being essential, the canonical homomorphism $\cA\ni a\mapsto(a\bullet,\bullet a)\in M(\ol{\cI})\subset M(\cI)$ is faithful. As for uniqueness, of course, we will require that the containment of $\cI$ via $M(\cI)\supset\cA\supset\cI$ is the same as the containment of $\cI$ in its multiplier algebra. So, if $\alpha\colon\cA\rightarrow M(\cI)$ is any other monomorphism satisfying $\alpha(i)=i$, by using a bounded approximate unit $u_\lambda$ for $\cI$, we get
\beqn{
 \alpha(a)i
~=~
\lim_\lambda\alpha(a)u_\lambda i
~=~
\lim_\lambda\alpha(a)\alpha(u_\lambda)i
~=~
\lim_\lambda\alpha(au_\lambda)i
~=~
\lim_\lambda au_\lambda i
~=~
ai,
}\eeqn
so the action of $\alpha(a)$ on $i$ is $a$.)
\end{enumerate}
\eob

\bdefi \label{contdefi}
For a \nbd{*}ideal $\cI$ in $\cA$, we say $\cA$ \hl{contains} $\sB^a(\cI)$ if the canonical homomorphism $\sB^a(\cI)\rightarrow\sB^a(\cI)\oplus0\subset M(\cI)\oplus M(\cI^\perp)=M(\cI\oplus\cI^\perp)$ is into $\cA\subset M(\cI\oplus\cI^\perp)$.
\edefi

\bthm \label{contthm}
Suppose $\cA$ is a \nbd{C^*}algebra, $p\in\cA$ is a projection, and put $E:=\cA p$. Then $\cA$ contains $\sB^a(E)$ $(=M(\sF(E)))$ if and only if $\cA=\sB^a(E)\oplus\sF(E)^\perp$.
\ethm

\proof
If $\cA=\sB^a(E)\oplus\sF(E)^\perp\supset\sF(E)\oplus0=\sF(E)$, then the canonical homomorphism to be constructed to check containment is, indeed, $\sB^a(E)\rightarrow\sB^a(E)\oplus0$.

Conversely, suppose $\cA$ contains $\sB^a(E)$. Denote by $P\in M(\sF(E))\oplus M(\sF(E)^\perp)$ the unit of $\sB^a(E)=M(\sF(E))\subset\cA$. Clearly, $P$ is central in $M(\sF(E))\oplus M(\sF(E)^\perp)$, so, \it{a fortiori} is it central in $\cA\subset M(\sF(E))\oplus M(\sF(E)^\perp)$. Since $\sB^a(E)\oplus 0\subset\cA\subset M(\sF(E))\oplus M(\sF(E)^\perp)$, we get
\beqn{
\sB^a(E)
~=~
P(\sB^a(E)\oplus 0)
~\subset~
P\cA
~\subset~
P(M(\sF(E))\oplus M(\sF(E)^\perp))
~=~
\sB^a(E),
}\eeqn
that is $P\cA=\sB^a(E)$. Therefore, $\cA=P\cA\oplus(\U-P)\cA$. (If $\cA$ is nonunital, then $\U\in M(\sF(E)\oplus\sF(E)^\perp)$ is not in $\cA$. However, in the decomposition $a=Pa+(a-Pa)$, since $Pa$ is in $\cA$, also $(\U-P)a=a-Pa$ is in $\cA$ for all $a\in\cA$.) Obviously, $(\U-P)\cA=\sF(E)^\perp$.\qed

\lf
\bdefi \label{sfulldefi}
We say the pair $(\cA,p)$ is \hl{semifull} or the corner $p\cA p$ is \hl{semifull}\index{semifull}\index{dilation, weak!full!semifull} in $\cA$ if $\cA$ contains $\sB^a(E)$. Likewise, a triple (in particular, a dilation) $(\cA,\theta,p)$ is \hl{semifull} if $(\cA,p)$ is.
\edefi

%%%% BO
% I think the following is an example: 
Not all dilations are semifull:

\bex \label{nsfex}
Let $\cA = \sK(H) + \id_H\C$, with $H$ infinite-dimensional, and let $p$ be a rank-one projection. Then $E = \cA p = H$ and $\ls\cA p\cA = \sF(E) = \sF(H)$ is just the algebra of all finite rank operators on $H$. So, $M(\sF(E)) = \sB(H)$ is not contained in $\cA$. Letting $\theta_t$ be the identity, we get a dilation $(\cA,\theta,p)$ of the identity semigroup (on $p\cA p=\id_H\C$). For a less trivial example, we may take $H$ to carry the minimal coisometric dilation $\bfam{w^n}$ ($w$ as in \ref{classmin}) of the the contraction semigroup $\bfam{c^n}$ for some $c\in\C$ with $\abs{c}<1$, so that $(w^*H)^\perp$ is finite-dimensional and, therefore, $w^*\bullet w$ leaves $\cA$ invariant. Then the 
%%%% BO new
%restriction  
(co)restriction 
% :-) 
%%%% EO 
of the corresponding (strong) solidly elementary dilation (Section \ref{EXwnsSEC}) to $\cA=\sK(H)+\id_H\C$, does the job.

% Finally, we construct an example of a strong dilation of a proper CP-semigroup which is not semifull. Let $t$ be a proper nonzero contraction generating an elementary CP-semigroup on the finite dimensional space $G$, and let $v \in \sB(H)$ (where $H \supset pH = G$) be its minimal coisometric dilation. Now let $\cA$ be the C*-algebra generated by $\sK(H)$ and by all the projections ${v^*}^k v^k$, $k \in \N_0$. Now let $\theta$ be the elementary \nbd{E}semigroup on $\cA$ determined by $v$. Then $(\cA,\theta,p)$ is a strong dilation; again we have that $E = \cA p = H$ and that $M(\sF(E)) = \sB(H)$ is not contained in $\cA$ (note that $\cA$ is separable.) 
\eex
% old TODO: 
%\OW[OPEN]{At first sight, it seems unrealistic that  we get that direct sum decomposition in the \nbd{C^*}case. On the other hand, I have not the slightest idea how a counter example, a dilation where $\cA$ does not contain $\sB^a(E)$, could exist.}
%%%% EO

\bob
We might think that we could use the embedding $\cA\subset M(\cI\oplus\cI^\perp)$ to introduce a strict topology on all of $\cA$. But in particular in the case when $\cA$ contains $M(\cI)$, which we would prefer, we get that $\cA=M(\cI)\oplus\cI^\perp$. Since in our applications $\cA$ is unital and, of course, $M(\cI)$ is unital, so is $\cI^\perp$. So, on $\cI^\perp$ we get nothing but the norm topology, and we know from Observation \ref{pstrcontob} that this is not suitable.

At this point, we give up our attempts to find a ``natural'' strict topology on $\cI^\perp$ that would allow us to define a ``reasonable'' notion of strictness for a general dilation in the \nbd{C^*}case, but content ourselves with the good thing that at least in the von Neumann case we have normality.
\eob

Fortunately, as we mentioned already, in the von Neumann case everything is considerably simpler. $\sB^a(E)$ is the strongly closed ideal in $\cA$ spanned by $\sF(E)$. The unit $P$ of $\sB^a(E)$ is a central projection, so that $\cA=P\cA\oplus(\U-P)\cA$. So this is compatible with our earlier definition of containment (thanks to the easy \it{if} direction of Theorem \ref{contthm}).

Only because in the literature some notions of minimality in the von Neumann case refer to the so-called central cover of a projection in a von Neumann algebra, we add the following piece of \it{folklore} that relates central covers to full dilations. The \phantomsection\hl{central cover}\index{central cover} of a projection $p$ in a von Neumann algebra $\cA$ is the smallest central projection $c(p)\in\cA$ such that $c(p)\ge p$. (Since the product of two central projections $Q_i\ge p$ is a central projection fulfilling $Q_i\ge Q_1Q_2\ge p$, the set of central projections $Q\ge p$ is directed decreasingly; since $\U\ge p$, it is non-void. So, since we are in a von Neumann algebra, the net converges to its unique minimum $c(p)\ge p$, which is also a central projection.) Now, since $p\in\sB^a(E)\subset\cA$, we have $P=\id_E\ge p$. If $Q\ge p$ is another central projection, then $Qapa'=aQpa'=apa'$, so that $Q$ acts as identity on $\sB^a(E)$. Therefore, $Q\ge P$, that is $P=c(p)$.

We summarize for case of dilations.

\bthm \label{vNcontthm}
Let $(\cA,\theta,p)$ be a normal dilation on the von Neumann algebra $\cA$. Then:
\begin{enumerate}
\item \label{vNct1}
$\cA$ contains $\sB^a(E)=P\cA=c(p)\cA$, that is, the dilation is semifull.

\item \label{vNct2}
The dilations is full if and only if $c(p)=\U$.
\end{enumerate}
\ethm

\noindent
The latter condition is also equivalent to $\cls\cA G=H$ (where $\cA\subset\sB(H)$). In the \nbd{C^*}frame\-work this condition may not even be phrased without first representing $\cA$ faithfully (and, of course, nondegenerately) on a Hilbert space. But, whether or not this condition is fulfilled may depend on the representation. (For instance, if $\cA=\sB(H)$ and $p$ a rank-one projection as in Example \ref{nsfex}, then the condition is fulfilled. However if we represent $\cA=\sB(H)$ faithfully not normally as $\cA\subset\sB(H')$, then $\cls\cA G'$ is the subspace of $H'$ on which $\cA$ acts normally.)

\lf
Now suppose $(\cA,\theta,p)$ is just a unital \nbd{C^*}algebra with an \nbd{E}semigroup and a projection as in Theorem \ref{E-supPSthm}, and assume the pair $(\cA,p)$ is semifull. (In the von Neumann case, this is automatic; in the \nbd{C^*}case, we require it.) Then $\cA$ decomposes as $\cA^0\oplus\cA^1=\rtMatrix{\cA^0\\\cA^1}$ with $\cA^1:=\sB^a(E)\ni p$ and $\cA^0:=\sF(E)^\perp$. Like for any direct sum decomposition, we denote $\theta_t=\rtMatrix{\theta^{00}_t&\theta^{01}_t\\\theta^{10}_t&\theta^{11}_t}$ (with homomorphisms $\theta^{ij}_t\colon\cA^j\rightarrow\cA^i$) and the semigroup property (recall that $\theta$ is a semigroup over $\bS^{op}$) reads
\beqn{
\SMatrix{\theta^{00}_{st}&\theta^{01}_{st}\\\theta^{10}_{st}&\theta^{11}_{st}}
~=~
\SMatrix{\theta^{00}_t&\theta^{01}_t\\\theta^{10}_t&\theta^{11}_t}
\circ
\SMatrix{\theta^{00}_s&\theta^{01}_s\\\theta^{10}_s&\theta^{11}_s}
~=~
\SMatrix{
\theta^{00}_t\circ\theta^{00}_s+\theta^{01}_t\circ\theta^{10}_s~&~\theta^{00}_t\circ\theta^{01}_s+\theta^{01}_t\circ\theta^{11}_s
\\
\theta^{10}_t\circ\theta^{00}_s+\theta^{11}_t\circ\theta^{10}_s~&~\theta^{10}_t\circ\theta^{01}_s+\theta^{11}_t\circ\theta^{11}_s~
}.
}\eeqn
In general, $p$ would decompose into $\rtMatrix{p^0\\p^1}$ for two projections $p^i\in\cA^i$. The superproduct system would look like
\beqn{
\sE_t
~:=~
\theta_t(p)\cA p
~=~
\SMatrix{
(\theta^{00}_t(p^0)+\theta^{01}_t(p^1))\cA^0p^0
\\
(\theta^{10}_t(p^0)+\theta^{11}_t(p^1))\cA^1p^1
}
~=~
\SMatrix{
\sE^{00}_t\oplus\sE^{01}_t
\\
\sE^{10}_t\oplus\sE^{11}_t
}
}\eeqn
with $\sE^{ij}_t:=\theta^{ij}_t(p^j)\cA^ip^i$. We omit writing down the intricate product of that superproduct system, as we do not need it in that generality. In our situation $\cA^1=\sB(E)\ni p$, we have $p=\rtMatrix{0\\p^1}$ with $p=p^1\in\cA^1$. The superproduct system looks
\beqn{
\sE_t~
:=~
\theta_t(p)\cA p
~=~
\theta_t(p^1)\cA p^1
~=~
(\theta^{01}_t(p^1)+\theta^{11}_t(p^1))\cA^1p^1
~=~
\theta^{11}_t(p^1)\cA^1p^1
~=~
\sE^{11}_t
}\eeqn
with product
\bmun{
v_{s,t}
\colon
(\theta^{11}_s(p^1)a^1_sp^1)\odot(\theta^{11}_t(p^1)a'^1_tp^1)
~\longmapsto~
\theta_t\bfam{\theta^{11}_s(p^1)a^1_sp^1}\theta^{11}_t(p^1)a'^1_tp^1
\\
~=~
\theta^{11}_t\circ\theta^{11}_s(p^1)\theta^{11}_t(a^1_s)\theta^{11}_t(p^1)a'^1_tp^1.
}\emun

\bprop
If $\sE^\podot$ is a product system, then, necessarily, $\theta^{11}_t\circ\theta^{11}_s(p)=\theta^{11}_{st}(p)$, that is, $\theta^{10}_t\circ\theta^{01}_s(p)=0$. This means: The part of $\theta_s(p)$ that is not in $\sB^a(E)$, never comes back.
\eprop

\proof
If $\theta^{10}_t\circ\theta^{01}_s(p)\in\cA^1=\sB^a(E)$ is not $0$, then the elements in $\cA^1p^1=\cA p=E$ ``see that''. Consequently, $\theta^{11}_t\circ\theta^{11}_s(p^1)\theta^{11}_t(\cA^1)\theta^{11}_t(p^1)\cA^1p^1\subset\theta^{11}_t\circ\theta^{11}_s(p^1)\cA^1 p^1$ cannot span all of $\theta^{11}_{st}(p^1)\cA p^1=\sE_t$.\qed

\lf
What can we say for the opposite direction? Well, not really much without topological assumptions; this is, why we discuss it in this section. Topological assumptions means, in the von Neumann case with normal endomorphisms, we can say quite a lot; and if we want to say something about the \nbd{C^*}case, we need strictness. Looking carefully at the proposition, we note that, for fixed $s$ and $t$, we are actually speaking about the homomorphisms $\theta^{11}_t$, $\theta^{11}_s$, $\theta^{11}_{st}$, and the composition of the first two, $\theta^{11}_t\circ\theta^{11}_s$, in 
%%%% BO new
the
%%%% EO 
situation where the difference of  $\theta^{11}_{st}$ and the latter, $\theta^{10}_t\circ\theta^{01}_s$, is also a homomorphism. All these homomorphisms act on $\cA^1=\sB^a(E)$; and here we know what strict is. Strictness of $\theta^{11}_t$ is exactly what has been defined as \phantomsection\it{\nbd{p}relatively strict}\index{prelatively@\nbd{p}relatively strict} for $\theta_t$ in Question \ref{pstriQ} (for the not necessarily semifull case).

\blem
Let $E$ be a Hilbert module with a unit vector $\xi\in E$, and let $\vp$ and $\psi$ be strict endomorphisms of $\sB^a(E)$ such that $\rho:=\vp-\psi$ is an endomorphism, too. Then $\vp=\psi$ if (and, of course, only if) $\rho(\xi\xi^*)=0$.
\elem

\proof
Recall (see the beginning of Section \ref{DilSPSpSEC}) that for every strict endomorphism $\vt$ on $\sB^a(E)$ we have a correspondence $\sE_\vt$ and an adjointable isometry $v_\vt\colon E\odot\sE_\vt\rightarrow E$ (onto $\vt(\id_E)E$) such that $\vt=v_\vt(\bullet\odot\id_{\sE_\vt})v_\vt^*$, and recall that $\sE_\vt$ and $v_\vt$ are unique if $E$ is full. Recall, further, that when $E$ has a unit vector (so that $E$ is full), we may write $\sE_\vt:=\vt(\xi\xi^*)E$,~ $b.x_\vt:=\vt(\xi b\xi^*)x_\vt$,~ and $v_\vt(x\odot y_\vt):=\vt(x\xi^*)y_\vt$.

Now, $\vp=\psi+\rho$. Therefore, $\sE_\vp=\sE_\psi+\sE_\rho$ and $v_\vp=v_\psi+v_\rho$. If $\rho(\xi\xi^*)=0$, then $\sE_\rho=\zero$, hence, $v_\vp=v_\psi$. Consequently, $\vp=v_\vp(\bullet\odot\id_E)v_\vp^*=v_\psi(\bullet\odot\id_E)v_\psi^*=\psi$.\qed

\lf
\bthm \label{p-PSthm}
Let the pair $(\cA,p)$ be semifull and let $\theta$ be a \nbd{p}relatively strict\index{prelatively@\nbd{p}relatively strict} \nbd{E}semigroup on $\cA$. Then the following conditions are equivalent:

\begin{enumerate}
\item \label{p-PS1}
The $\theta^{11}_t$ form a strict semigroup $\theta^{11}$ on $\cA^1:=\sB^a(E)$.

\item \label{p-PS2}
$\theta^{11}_t\circ\theta^{11}_s(p)=\theta^{11}_{st}(p)$ for all $s,t\in\bS$.

\item \label{p-PS3}
$\sE^\podot$ is a product system (the product system of $\theta^{11}$).
\end{enumerate}
Moreover, under any of these conditions, $(\cA,\theta,p)$ is a dilation if and only if $(\cA^1,\theta^{11},p)$ is.
\ethm

\proof
By the proposition, we have \eqref{p-PS3}$\Rightarrow$\eqref{p-PS2}. To see \eqref{p-PS2}$\Rightarrow$\eqref{p-PS1}, apply the lemma to $\psi:=\theta^{11}_t\circ\theta^{11}_s$ and $\vp:=\theta^{11}_{st}$. For \eqref{p-PS1}$\Rightarrow$\eqref{p-PS3}, we simply recall that $\sE_t=\sE^{11}_t$, including the correct superproduct system structure. So if one is a product system then so is the other.\qed

\bob
If $\cA$ is a von Neumann algebra, so that every pair $(\cA,p)$ is semifull, we may weaken the hypothesis to normal $\theta$. The only things that change are, then, that in \ref{p-PSthm}\eqref{p-PS3} and in the proof of the lemma we have to pass to the tensor product of von Neumann correspondences (strong closure). So, normal dilations on von Neumann algebra $\cA$ may be (co)restricted to a (normal) dilation on $\sB^a(E)$ (necessarily contained in $\cA$) if and only if their superproduct system is a product system of von Neumann correspondences.
\eob

Theorem \ref{p-PSthm} and (in particular) its von Neumann version tell to some extent that there is not much we can do in the direction of restriction of a given dilation to $\cA^1=\sB^a(E)$, when the superproduct system of this dilation is not a product system. They also provide a powerful (because easy-to-check) criterion for when this is the case. (We suggest as an exercise to see how Example \ref{hypexex} fits into the discussion.)

Trying (co)restriction to $\cA^1\subset\cA$ is actually a \it{compression} with the projection $P=\id_E$. Compression is already among the strategies that occur when we discuss minimalities; we postpone these to
%%%% BO 
% the next section.
Section \ref{minSEC}. 
%%%% EO

What we might try is also to see whether we can \it{extend} the \nbd{E}semigroup to an algebra $\wh{\cA}$ containing $\cA$. If we wish a module dilation, we may try to identify $\sK(E)^\perp$ as another $\sB^a(E^0)$ so that $\cA$ is the diagonal subalgebra of $\wh{\cA}=\sB^a\rtMatrix{E^0\\E}=\rtMatrix{\sB^a(E^0)&\sB^a(E,E^0)\\\sB^a(E^0,E)&\sB^a(E)}$. Even in the case when $\vt^{11}$ and $\vt^{00}$ are \nbd{E_0}semigroups (leaving no space for the original $\vt$ being non-diagonal, so that in this case there would be nothing to do for us), in order to allow an extension to $\wh{\cA}$, they necessarily need to have the same product system; see \cite[Theorem 6.7]{Ske16} for the one-parameter case. A first problem to look at, is the general form of the product system of any strict \nbd{E}semigroup on $\wh{\cA}$ and to see how the given $\sE^{ij}_t$ fit into that. (They must appear somehow as subspaces of that product system.) We do not tackle this question in these notes.

\newpage

\section[\sc{Examples:} Bhat's example; a ``multi-example'']{Examples: Bhat's example; a ``multi-example''}\label{EXBexSEC}

We report an example due to Bhat (starting after the proof of \cite[Theorem 3.2]{Bha03}). Bhat shows that his example is \it{incompressible} (see Definition \ref{incomprdefi}; Bhat calls it \it{atomic}) but not \it{minimal} (see Theorem \ref{1pamcnMcor}). Here, we show that it is a true ``multi-example'', which serves as a counter example for several more properties relevant in several places in these notes. Some places did occur already; others will occur in the next section on minimality. Therefore, the present section also serves to some extent as a motivation for the notions of minimality we shall introduce and discuss in the next section.

Bhat's example is the dilation of a discrete one-parameter CP-semigroup on $\sB(\C)=\C$; it illustrates that already in the simplest of all one-parameter cases possible, not all is as simple and understood as literature may try to make us believe. In fact, Bhat's example is non-Markov, while the statements it contradicts, in the one-parameter Markov case are true. 
%%%% BO 
%(For instance, every dilation of a one-parameter Markov semigroup can be compressed to the unique minimal one; see 
%\OW[ORR: Right now, I am not sure. If it is true, then we should include it (in way that is similarly stringent) in Section \ref{minSEC}\ref{1-p-SSEC}, and not a comma more complicated.]{where?}.) 
% (For instance, every dilation of a normal one-parameter Markov semigroup on a von Neumann algebra can be compressed to the unique minimal one; see\cite{Arv97} or \cite[Theorem 8.9.7]{Arv03}.)
%%%% EO
(After all, Bhat's example is not a good dilation -- a thing that does not happen in the Markov case.)
The more important it is to be aware that the other ``multi-example'', which we shall discuss in the conclusive Section \ref{EXN02SEC}, is Markov but (discrete) two-parameter; it is, therefore, a completely different type of example.

We have to recall some facts about (one-parameter, but not only) CP-semigroups (on $\sB(H)$, but not only) and their dilations; they are (almost) all well-known, but we recover them easily in our terminology.

\lf
\subsection{\normalsize CP-semigroups and their dilations under compression} \label{CPcompSSEC}
Let $(\cA,\theta,P)$ be a dilation and assume that the dilated CP-semigroup $T_t:=P\theta_t(P\bullet P)P$ on $P\cA P$ can be compressed further by $p\le P$ in $\cA$ to a CP-semigroup $S_t:=pT_t(p\bullet p)p=pP\theta_t(Pp\bullet pP)Pp=p\theta_t(p\bullet p)p$ on $p\cA p\subset P\cA P$. Then $(\cA,\theta,p)$ is a dilation of $S$.\index{compression!of a CP-semigroup}\index{CP-semigroup!compression of}

\lf
\subsection{\normalsize Kraus decomposition for CP-maps and endomorphisms}
If $T$ is a normal CP-map on $\sB(G)$, then there are $c_i\in\sB(G)$ such that $T=\sum_ic_i^*\bullet c_i$. (This is so, because for every such CP-map there is a von Neumann correspondence $E$ over $\sB(G)$, for instance, the strong closure of the GNS-correspondence $\sE$, and a vector $\xi\in E$ such that $T=\AB{\xi,\bullet\xi}$. \index{Kraus decomposition}\index{CP-map!Kraus decomposition}\index{Kraus decomposition!for endomorphisms}By Appendix \ref{vNAPP}\ref{vNBGmod}, we have $E=\sB(G,G\otimes\eH)=\sB(G)\sbars{\otimes}\eH$, and for any ONB $\bfam{e_i}$ of $\eH$, the (unique) $c_i$ such that $\xi=\sum_ic_i\otimes e_i$ do the job.)

Moreover, $T$ is an endomorphism if and only if the $c_i$ fulfill $c_ic_j^*=\id_G\delta_{i,j}$. These conditions show, that in this case $\xi$ is a coisometry. And $v:=\xi^*$ is the isometry $G\otimes\eH\rightarrow G$ such that $T=v(\bullet\otimes\id_\eH)v^*$. (See the discussion of \it{left semidilation} following Proposition \ref{lsdilPSprop}.)

\lf
%%%% BO postarxiv
\subsection{\normalsize Kraus decomposition and dilation} \label{KrowSSEC}
% \subsection{\normalsize Kraus decomposition and dilations} \label{KrowSSEC}
%%%% EO postarxiv
Kraus decompositions do not behave nicely under composition. In the continuous time case, this makes them rather useless. (The deeper reason is that the multiplicity spaces $\eH_t$ for each $T_t$ are closely related to the associated Arveson system, and that except in the one-dimensional case, the ONBs of $\eH_t$ that determine the Kraus decompositions, cannot be chosen compatible with the structure of the Arveson system.) At least in the discrete case, where all semigroups are generated by a single map, we can say a bit more. (Although, paralleling Section \ref{EXwnsSEC}, we could also discuss the non-strong case, we limit ourselves to the strong case, which we really need.)

Let $(\sB(H),\bfam{\theta^n}_{n\in\N_0},p)$ be a strong dilation of the semigroup $\bfam{T^n}_{n\in\N_0}$ on $\sB(G)=p\sB(H)p$, where $G:=pH$ and $T:=p\theta(p\bullet p)p$. Let $w_i$ be a Kraus decomposition for $\theta$, so that $w_iw_j^*=\id_H\delta_{i,j}$. By Observation \ref{semgenob}, the triple $(\sB(H),\bfam{\theta^n}_{n\in\N_0},p)$ is a strong dilation if and only if
\beqn{
\sum_i pw_i^*papw_ip
~=~
T(pap)
~=~
p\theta(a)p
~=~
\sum_i pw_i^*aw_ip
}\eeqn
for all $a\in\sB(H)$. It follows that $0=p\theta(\id_H-p)p=\sum_i((\id_H-p)w_ip)^*((\id_H-p)w_ip)$, so, necessarily $w_ip=pw_ip$. (Compare this with Proposition \ref{strsolprop}.) Clearly, the $c_i:=w_ip$ form a Kraus decomposition for $T$.

The family $\bfam{w_i^*}$ is what is known as an \it{isometric dilation} of the \it{row contraction} $\bfam{c_i^*}$ (satisfying $\sum_ic_i^*c_i\le\id_G$, which means exactly that $T$ is a contraction). Every discrete one-parameter semigroup admits a strong dilation. (This is the discrete version of Bhat and Skeide \cite[Section 8]{BhSk00}. It also follows putting together Theorems \ref{indlimthm} and \ref{Oreindthm} with Theorem \ref{uninonunithm}.) We, therefore, see that $T$ can be recovered from \bf{some} Kraus decomposition $\bfam{c_i}$ that admits a \phantomsection\hl{coisometric dilation}\index{dilation!coisometric, of a column contraction}\index{Kraus decomposition!coisometric dilation} $\bfam{w_i}$ in the sense that $\bfam{w_i^*}$ is an isometric dilation of $\bfam{c_i^*}$.\index{dilation!isometric, of a row contraction}

But what, if a (contractive) $T$ is given by a Kraus decomposition $\bfam{c_i}$ to begin with? Can we always arrange the strong dilation $(\sB(H),\bfam{\theta^n}_{n\in\N_0},p)$ for $\bfam{T^n}_{n\in\N_0}$, such that $\theta$ has a Kraus decomposition $\bfam{w_i}$ that allows to recover the given $c_i$ as $w_ip$? In other words, can we always find a coisometric dilation $\bfam{w_i}$ of $\bfam{c_i}$?

The answer is well known, and affirmative; see Bunce \cite{Bun84} and Popescu \cite{GPop89}, who solved the problem for countably many $w_i$ (earlier Frazho \cite{Fra82} treated the case of two). 
%%%% BO new
% Regarding the references: the way I see it, Popescu obtained the result independently (in 1987, in Romania), and it makes sense to cite them both. There are sufficient differences between the proofs, and one can argue that Popescu's approach gives more information, e.g. uniqueness as well as the structure of the dilation. (Popescu's approach also seemed to survive better with time). 
%%%% EO 
We can recover this result (for arbitrary dimension!) by passing to correspondences as in the preceding point, setting $\xi:=\sum_ic_i\otimes e_i\in\sB(G,G\otimes\eH)=:E$ ($\eH$ a Hilbert space with ONB $\bfam{e_i}$) and apply ideas from Skeide \cite[Theorem 1.2]{Ske08}. The version there, is for the continuous time one-parameter case; the much simpler restriction to the discrete case, is what we report here:

\begin{itemize}
\item
Do to $(E,\xi)$ exactly the same as in \ref{GNSuni} to obtain $(\wt{E},\wt{\xi})$ where $\AB{\wt{\xi},\bullet\wt{\xi}}=\wt{T}$ is the unitalization of $\AB{\xi,\bullet\xi}=T$; just that now $E$ is not necessarily the GNS-correspondence of $T$, but may be bigger.

\item
 Input product system and (unital) unit generated by $(\wt{E},\wt{\xi})$ in Theorem \ref{Oreindthm}, to obtain a module dilation $(\wt{\sE},\bfam{\wt{\vt}^n}_{n\in\N_0},\wt{\xi}_0)$ of $\bfam{\wt{T}^n}_{n\in\N_0}$.

\item
This dilation, analyzing (as in \cite{Ske08}) its structure, turns out to be the unitalization of a (necessarily strong) module dilation $(\sE,\bfam{\vt^n}_{n\in\N_0},\xi_0)$ of $\bfam{T^n}_{n\in\N_0}$. (See also Subsections \ref{minSEC}\ref{uniminSSEC} and \ref{1-p-SSEC}.)% (See also Corollaries \ref{fufcor} and \ref{fufdilcor}.)

\item
Let $v\colon\sE\sodots E\rightarrow\sE$ denote the isometry such that $\vt=v(\bullet\odot\id_E)v^*$. Recall that, putting $H:=\sE\odot G$, we have $\sE=\sB(G,H)$ and that $\sE\sodots E=\sE\sbars{\otimes}\eH=\sB(G,H\otimes\eH)$, so that $v$ may also be identified with the operator $v\odot\id_G\in\sB(H\otimes\eH,H)\odot\id_G$, and in this picture $\vt = v(\bullet\otimes\id_\eH)v^*$ on $\sB(H)$. Then, clearly, the $w_i:=(\id_H\otimes e_i^*)v^*\in\sB(H)$ form a Kraus decomposition for $\vt$.

\item
Since for every $a\in\sB(H)=\sB^a(\sE)$, we have that $a(\id_H\otimes e_i^*)=(\id_H\otimes e_i^*)(a\otimes\id_\eH)$ and $(a\otimes\id_\eH)v^*=v^*\vt(a)$, we get
\beqn{
\xi_0\xi_0^*w_i\xi_0\xi_0^*
~=~
\xi_0\xi_0^*(\id_H\otimes e_i^*)v^*\xi_0\xi_0^*
~=~
(\id_H\otimes e_i^*)v^*\vt(\xi_0\xi_0^*)\xi_0\xi_0^*.
}\eeqn
Note that $\id_H\otimes e_i^*=(\id_\sE\odot\id_G)\otimes e_i^*=\id_\sE\odot(\id_G\otimes e_i^*)$, the amplification of the operator $\id_G\otimes e_i^*\in\sB^{a,bil}(E,\sB(G))=\sB^{bil}(G\otimes\eH,G)$. Recall, too, that $E$ (the member of the product system for $n=1$) is recovered as $\xi_0E=\vt(\xi_0\xi_0^*)\sE$ and that the element corresponding to $\xi\in E$ is $\xi_0\xi=\vt(\xi_0\xi_0^*)\xi_0$. Therefore,
\bmun{
\xi_0\xi_0^*w_i\xi_0\xi_0^*
~=~
(\id_H\otimes e_i^*)v^*(\xi_0\xi)\xi_0^*
~=~
(\id_\sE\odot(\id_G\otimes e_i^*))(\xi_0\odot\xi)\xi_0^*
\\
~=~
\xi_0((\id_G\otimes e_i^*)\xi)\xi_0^*
~=~
\xi_0c_i\xi_0^*
~=~
c_i
~\in~
\sB(G)
~=~
\xi_0\xi_0\sB(H)\xi_0\xi_0^*.
}\emun
\end{itemize}

\noindent
It is noteworthy that the obtained coisometric dilation of $\bfam{c_i}$ is minimal in the sense used in that theory. (This corresponds to the Property (2) in \cite[Theorem 1.2]{Ske08}.) The dilation of $\bfam{T^n}_{n\in\N_0}$ is minimal in the sense of Subsection \ref{minSEC}\ref{1-p-SSEC} (Theorem \ref{1pamcnMcor}) if and only if $(E,\xi)=$ GNS-$T$, that is, if and only if the Kraus decomposition of $T$ is minimal; see Example \ref{Krausex}.

%%%% EO

\lf
\subsection{\normalsize The example}\label{BexSSEC}

\bex \label{Bex}\index{Bhat's example|bf}
Now, we are ready to discuss Bhat's example. Let $G:=\C^2$. Define the elements
\baln{
c_1
&
~:=~
\frac{1}{\sqrt{2C}}\SMatrix{2&1\\-1&0},
&
c_2
&
~:=~
\frac{1}{\sqrt{6C}}\SMatrix{0&1\\3&0},
&
c_3
&
~:=~
\frac{1}{\sqrt{3C}}\SMatrix{0&1\\0&0}
}\ealn
in $\sB(G)=M_2$. We readily verify that the CP-map $T:=\sum_ic_i^*\bullet c_i$ and the semigroup $T^n$ generated by it, act as
\beqn{
T\SMatrix{a&b\\c&d}
~=~
\frac{1}{C}\SMatrix{2(a+d)-(b+c)&a\\a&a},
~~~~~~~~~~~~
T^n\SMatrix{a&b\\c&d}
~=~
\frac{2(a+d)-(b+c)}{4}\family{\frac{2}{C}}^n\SMatrix{2&1\\1&1}
~~~(n\ge2). 
}\eeqn
We find $\norm{T}=\norm{T(\id_G)}=\frac{5+\sqrt{13}}{2C}$, so $T$ is a contraction if and only if $C\ge\frac{5+\sqrt{13}}{2}$. Let $(\sB(H),\bfam{\theta^n}_{n\in\N_0},P)$ denote the strong dilation obtained as above from the Kraus decomposition $c_1,c_2,c_3$. So, $\theta=\sum_iw_i^*\bullet w_i$, where the $w_1,w_2,w_3$ form a 
%%%% BO
% (minimal) 
% minimal
%%%% EO
coisometric dilation of $c_1,c_2,c_3$.
% (Since the Kraus decomposition is, clearly, minimal,
%%%% BO 
% the dilation generated by $\theta$ is minimal, too.)
% $\theta$ is the minimal dilation of $T$, too; we shall not require this.)
%%%% EO

We define $p$ to be the projection onto the first component in $\C^2=G$. Then, $p\rtMatrix{a&b\\c&d}p=\rtMatrix{a&0\\0&0}$ and
\beqn{
pT^n\SMatrix{a&0\\0&0}p
~=~
\family{\frac{2}{C}}^n\SMatrix{a&0\\0&0}
}\eeqn
for all $n$ (including $n=1$). So, $p$ compresses $\bfam{T^n}_{n\in\N_0}$ to a semigroup $\bfam{S^n}_{n\in\N_0}$ and, identifying $p\in\sB(G)=P\sB(H)P\subset\sB(H)$, we get a dilation $(\sB(H),\bfam{\theta^n}_{n\in\N_0},p)$ of $\bfam{S^n}_{n\in\N_0}$. 
\eex

\lf
So far, the exposition is from Bhat \cite{Bha03}. Bhat produced this example to obtain an instance for a dilation exhibiting the following bad behaviour in terms of certain notions of minimality, which we will discuss in the following section: An \it{incompressible} dilation that is not the (unique, in the one-parameter case) \it{algebraically minimal} strong or full dilation. His proof is by direct computations.

We state Bhat's result, and some more, in a moment with references to the following section for notions of minimality and the necessary results about them. But, our proof is different, and makes use of things we can state and prove now referring only to what we know already.

About the dilation $(\sB(H),\bfam{\theta^n}_{n\in\N_0},p)$ of $\bfam{S^n}_{n\in\N_0}$ we have the following facts:

\begin{itemize}
\item
The dilation is full and (not only normal, but, thanks to finite multiplicity $w_1,w_2,w_3$, even) strict. Therefore, the superproduct system $\sE_t:=\theta^n(p)\sB(H)p$, by Observation \ref{fu-mod-ob} and Example \ref{EPSex}, is a product system.

Since $p=ee^*$ is a rank-one projection, we can identify this product system with the Arveson system formed by the Hilbert spaces $H_n:=\theta^n(p)\sB(H)pe=\theta^n(p)H$. (Of course, Bhat, the first one associating in \cite{Bha96} Arveson systems to CP-semigroups on $\sB(G)$, knew this fact. But in his presentation he opted to not refer to product systems.) 
%%%% BO 
%%%% EO

\item
By very direct computation, the $\theta_n(p)p$ do not form a unit for that product system, that is, the dilation is not good. (Indeed, if they formed a unit then we must have $\theta(\theta(p)p)p=\theta^2(p)p$ or
\beqn{
\sum_{j,i}w_j^*w_i^*pw_ipw_jp
~=~
\sum_{j,i}w_j^*w_i^*pw_iw_jp.
}\eeqn
By the orthogonality relations, this means $pw_ipw_jp=pw_iw_jp$ for all $i,j$. By $p\le P$, we may insert lots of $P$s, and obtain, taking also into consideration that the dilation $(\sB(H),\bfam{\theta^n}_{n\in\N_0},P)$ is strong, the necessary condition
\beqn{
pc_ipc_jp
~=~
pc_ic_jp.
}\eeqn
Clearly, $i=1$ and $j=2$, or $i=2,3$ and $j=1,2$ violate that condition; see below.)
\end{itemize}
We preferred to let the reader see this direct computation. But, the result also follows from the following key lemma about the product system of the dilation.

\blem \label{Bkeylem}
There is no proper product subsystem $\sF^\odot$ of $\sE^\odot$ satisfying 
%%%% BO 
$\theta_n(p)p\in\sF_n$ for all $n\in\N_0$.
%%%% EO
\elem

\proof
%%%% BO
We will work directly with the Hilbert spaces $H_n$, and we shall write $\xi_n := \theta_n(p)pe$. We will show that there is no proper subsystem containing $\xi_1$ and $\xi_2$. 
%%%% EO

Already in the preceding computations, the expressions $pc_ip=e\AB{e,c_ie}e^*$ and $pc_ic_jp=e\AB{e,c_ic_je}e^*$ showed that the vector $d\in\C^3$ with coordinates $d_i:=\AB{e,c_ie}$ and the matrix $D\in M_3$ with entries $d_{i,j}:=\AB{e,c_ic_je}$ play a role. We find
\baln{
d
&
~=~
\sqrt{\frac{2}{C}}\,\SMatrix{1\\0\\0},
&
D
&
~=~
\frac{1}{6C}\SMatrix{9&3\sqrt{3}&0\\-\sqrt{3}&3&0\\-\sqrt{6}&~~3\sqrt{2}~~&0}.
}\ealn
Note that product subsystems are exactly the images under projection morphisms $p^\odot=\bfam{p_n}_{n\in\N_0}$, and that $p_n=u_{1,\ldots,1}(p_1\otimes\ldots\otimes p_1)u_{1,\ldots,1}^*$ is determined by $p_1$. So, we have to find projections $p_1$ with $p_1\xi_1=\xi_1$ and $(p_1\otimes p_1)u_{1,1}^*\xi_2=u_{1,1}^*\xi_2$. Passing to $\id_1-p_1$, this happens exactly, if we find a projection $q$, such that
\beqn{
q\xi_1
~=~
0,
\text{~~~~~~~~~~~~}
(q\otimes\id_1)u_{1,1}^*\xi_2
~=~
0,
\text{~~~~~~and~~~~~~}
(\id_1\otimes q)u_{1,1}^*\xi_2
~=~
0.
}\eeqn
And since validity of these equations for a projection $q$ implies validity for any subprojection of $q$, we may assume $q=xx^*$ where either $x$ is a unit vector or $x=0$. We shall show $x=0$ and, thus, $p_1=\id_1$, completing the proof.

By $\theta(p)h=\sum_iw_i^*e\AB{w_i^*e,h}$ we see that that the vectors $w_i^*e$ form an ONB of $H_1$. Put $x=\sum_iw_i^*e x_i$. We have $\xi_1=\sum_iw_i^*ed_i=w_1^*e\sqrt{\frac{2}{C}}$. By $q\xi_1=0$, we conclude $x_1=0$. Now,
\beqn{
u_{1,1}(w_i^*e\otimes w_j^*e)
~=~
\theta(w_i^*ee^*)w_j^*e
~=~
w_j^*w_i^*e.
}\eeqn
Similarly, we have $\xi_2=\sum_{j,i}w_j^*w_i^*ed_{i,j}$, so,
\vspace{-.5ex}
\beqn{
u_{1,1}^*\xi_2
~=~
\sum_{j,i}(w_i^*e\otimes w_j^*e)d_{i,j}.
\vspace{-1.5ex}
}\eeqn
From $(q\otimes\id_1)u_{1,1}^*\xi_2=0$, we conclude
\vspace{-.5ex}
\beqn{
\sum_i\bar{x}_id_{i,j}
~=~
0
\vspace{-1.5ex}
}\eeqn
for all $j$. With $x_1=0$, this fixes a possible solution to $x=\rtMatrix{0\\-\sqrt{2}\\1}x_3$. From $(\id_1\otimes q)u_{1,1}^*\xi_2=0$, we conclude
\vspace{-.5ex}
\beqn{
\sum_jd_{i,j}\bar{x}_j
~=~
0
% \vspace{-.5ex}
}\eeqn
for all $i$. Therefore, $x_3=0$, hence, $x=0$, hence, $q=0$, hence, $p_1=\id_1$.\qed

\newpage
Obviously, $\sE^\odot$ is  not one-dimensional. (In fact, since $w_1^*, w_2^*, w_3^*$ are isometries with pairwise orthogonal ranges, $\dim H_1=\dim\sE_1e=\dim\cls\CB{w_i^*e}=3$.)

\bcor
The $\theta(p_n)p$ do not form a unit.
\ecor

\proof
If the $\theta(p_n)p$ formed a unit, then $\sF_n:=\theta(p_n)p\C$ was a proper product subsystem.\qed

\lf
This has, taking in, step by step, also notions and results from the next section, several consequences

\begin{itemize}
\item
Immediate is that the dilation is full and not good, hence, by Theorem \ref{sdilunithm}\eqref{DU1}, not strong.

\item
Whatever \it{compression} (Definition \ref{comprdef}) of the dilation, is full (Observation \ref{compob}\eqref{cob2}) and not good (Corollary \ref{goodccor}), hence, not strong.

\item
Bhat shows (for $C>\frac{5+\sqrt{13}}{2}$)  that the dilation is not \it{primary} (Definition \ref{vNprimdef}). We are not interested in this result and refer to \cite{Bha03}. We are interested in primary dilations. And like Bhat, we pass to a primary dilation (which is, by the preceding observations, full and not good) by compressing Example \ref{Bex} with $Q:=\bigvee_n\theta_n(p)$ (Theorem \ref{vNprimthm}\eqref{vNp3}).

\item
By Observation \ref{primpPSob}, the primary version of the dilation has the same (super)product system. By Lemma \ref{Bkeylem} and Theorem \ref{fpincthm}, it is incompressible. We, therefore, obtain a full and incompressible (hence, primary) dilation that is not good (hence, not strong).

\item
The von Neumann subalgebra generated by $\theta_{\N_0}(p)$ is not all of $\sB(QH)$, that is, the dilation is not (normally) \it{algebraically minimal} (Definition \ref{namindef}). (Indeed, we show that the subspace $K:=Q(w_3^*H\cap e^\perp)$ of $QH$ is nonzero and perpendicular to the smallest subspace $L$ invariant under all $\theta_{\N_0}(p)$ and containing $e$. As for $K\ne\zero$: Since the $w^Q_i:=Qw_iQ$ form a Kraus decomposition for $\theta^Q$, they fulfill the same orthonormality relations as the $w_i$, but relative to $QH$; since $QH$ is infinite dimensional ($\dim QH\ge\dim Q\sE_n=3^n$ for all $n$), so is its image under the isometry ${w_3^Q}^*$; since $e\in QH$, the space $K$ is infinite-dimensional, too. As for $\AB{K,L}=\zero$: Observe that $Qe=e$ and $Q\theta^n(p)=\theta^n(p)Q$, so $QL=L$; it is enough, to show that $L\subset L':=\cls(w_3^*H)^\perp\cup\CB{e}$, because, clearly, $\AB{K,L'}=\zero$. Now, $L'\ni e$, $pL'=e\C\subset L'$, and for $n\ge1$ and $x\in L'$
\beqn{
\theta^n(p)x
~=~
\sum_{i_1,\ldots,i_n}w_{i_n}^*\ldots w_{i_1}^*e\AB{e,w_{i_1}\ldots w_{i_n}x}
~=~
%%%% BO
\sum_{i_n\ne3}\sum_{i_1,\ldots,i_{n-1}}w_{i_n}^*\ldots w_{i_1}^*e\AB{e,w_{i_1}\ldots w_{i_n}x},
%%%% EO
}\eeqn
because $w_3e=0$ and $w_3y=0$ for $y\in (w_3^*H)^\perp$, so, $\theta^n(p)x\in(w_3^*H)^\perp\subset L'$; in other words, $L'$ contains $e$ and is invariant under $\theta_{\N_0}(p)$, so that $L'\supset L$.
\end{itemize}

%%%% BO
% This way it makes more sense to me: 
\noindent
The (normal) \it{algebraic minimalization} is neither full nor good. 
We do not know if the algebraically minimal version is incompressible or not.
% \OW[OPEN (Help! No energy left. --Orr.)]{If $Q'\in\ol{\sB(QH)_\infty}^s$ is a projection that compresses ${\theta^Q}^\infty$ further (so that $Q'\le Q$), how far can the $Q'w^Q_iQ'=Q'w_3Q'$ be away from satisfying the orthonormality relations, so that $\theta{Q'}$ defines an endomorphism such that ${\theta{Q'}}^n={\theta^n}^{Q'}$ on the subalgebra $Q'\ol{\sB(QH)_\infty}^sQ'$?}
But, by Theorem \ref{vNincfthm}, if it is incompressible, then its superproduct system is not a product system.
%%%% EO

\newpage

\section{Minimality}\label{minSEC}

Better, this section should be called `minimalities'. In fact, already in the one-parameter case there are several notions of minimality, which in the end turn out to be equivalent. Once we leave the one-parameter case we no longer have this luxury.\index{minimal!dilation|see {dilation, weak}}\index{dilation, weak!minimality}
%\OW[Isn't this even a bit incorrect? Once we drop conditions that assure fullness of the minimal dilation, we know we have non-full dilations that are minimal in several other senses also in the one-parameter case.]{Question!}

Rather than saying what \bf{is} minimal, it is somehow easier to imagine situations that are agreeably \bf{not} minimal. Is the algebra $\cA$ too big (in the sense that we can easily pass to a smaller one)? Is the unit of $\cA$ too big (in the sense that we can compress the dilation with a projection, becoming the new unit)? Can we obtain minimality from non-minimality in this sense? In other words, if a dilation is, for some obvious reason, not minimal, can we pass to a dilation that is not non-minimal (at least not for the same reason)? If there exists a dilation that is (in what sense ever) minimal, is this dilation unique?

It is noteworthy that already in the one-parameter case, guaranteeing the possibility to \it{minimalize} a given dilation is limited to strong dilations. (Example \ref{EXBexSEC}\ref{BexSSEC} gives us a dilation that is not minimal and cannot be minimalized (by compression) further.) And in the two-parameter case it is well-known that the minimal dilation need not be unique. (We present some examples below, see \ref{minnuniex}.) It is remarkable, that the most common notion of minimality in the one-parameter case involves the condition to be, in our terminology, a \it{full} dilation. It is one of our biggest open questions, whether or not existence of a dilation implies existence of a full dilation. (If yes, then a Markov semigroup that admits a dilation, also admits a module dilation and, therefore, the GNS-subproduct system embeds into a product system. We, then, had the latter embedding property as an \it{iff}-criterion for existence of dilations for a Markov semigroup over the opposite of an Ore monoid.)

\subsection{\normalsize Algebraic minimality} \label{algminSSEC}

Suppose $(\cA,\theta,p)$ is a dilation (say, of a CP-semigroup $T$). One possibility of this dilation being not what merits to be called \it{minimal}, is when we find a (for the time being, not necessarily unital) \nbd{C^*}subalgebra $\breve{\cA}\ni p$ such that $\theta$ (co)restricts to an \nbd{E}semigroup $\breve{\theta}$ on $\breve{\cA}$. We may ask whether $(\breve{\cA},\breve{\theta},p)$ (apart from $\breve{\cA}$ being possibly nonunital) is again a dilation. Note that this is actually two questions. The first is, whether $(\breve{\cA},\breve{\theta},p)$ is a dilation, the second is whether it is a dilation of the same CP-semigroup $T$. The answer to the first question is, clearly, yes. As for the second question, the dilated CP-semigroup is simply the (co)restriction $\breve{T}$ of $T$ to $\breve{\cB}:=p\breve{\cA}p\subset p\cA p=:\cB$. (Indeed, $T_t(p\breve{a}p)=p\theta_t(p\breve{a}p)p=p\breve{\theta}_t(p\breve{a}p)p\in p\breve{\cA}p=\breve{\cB}$.) So, $(\breve{\cA},\breve{\theta},p)$ is a dilation of the same semigroup $T$ if and only if $\breve{\cB}=\cB$, that is, if and only if $\cB\subset\breve{\cA}$.

We see, to be a dilation of the same CP-semigroup $T$, we have $\cB\subset\breve{\cA}$, hence, $\theta_t(\cB)=\breve{\theta}_t(\cB)\subset\breve{\cA}$. So, in this case, $\breve{\cA}$ always contains the \nbd{C^*}algebra  $\cA_\infty:=C^*(\theta_\bS(\cB))$ generated by all $\theta_t(\cB)$. On the other hand, $\cA_\infty$, clearly, \bf{is} invariant under all $\theta_t$, so we always may pass to the minimal $\cA_\infty$.\phantomsection\index{minimalization (algebraic)|see {dilation, weak}}\index{algebraically minimal|see {dilation, weak}}

Unfortunately, $\cA_\infty$ is usually, indeed, nonunital. (Unless $\theta_t(p)$ is constant; see Observation \ref{pstrcontob}.) At least, if $\theta$ is an \nbd{E_0}semigroup, we may add the unit $\U_\cA$ to $\cA_\infty$. (Note that $p(\U_\cA\C+\cA_\infty)p$ still equals $\cB$.) But, this is somehow not really what we want. What we want is a unital algebra obtained as closure in $\cA$ or completion of $\cA_\infty$ in a topology in which each $\theta_t$ is continuous. Surprisingly, we will have to distinguish not between two cases ($C^*$ and von Neumann case), but there will be a third one.

Let us start with the easiest case.

\bdefi \label{namindef}
A dilation $(\cA,\theta,p)$ on the von Neumann algebra $\cA$ is \hl{normally algebraically minimal}\index{dilation, weak!algebraically minimal!normally}\index{algebraically minimal!normally} if $\cA=\ol{\cA_\infty}^s$. If the context is clear, then we will just say $(\cA,\theta,p)$ is algebraically minimal.
\edefi

If a dilation $(\cA,\theta,p)$ on the von Neumann algebra $\cA$ is normal, then, by (co)restriction, we obtain its \phantomsection\hl{normal algebraic minimalization}\index{dilation, weak!algebraic minimalization!normal}\index{algebraically minimal!normally!minimalization} $(\ol{\cA_\infty}^s,\theta^{\infty^s},p)$. (Note that if $(\cA,\theta,p)$ is not normal, then it is unclear, whether or not $\theta$ (co)restricts to $\ol{\cA_\infty}^s$.) If the context is clear, then we will write just $\theta^\infty$ instead of $\theta^{\infty^s}$.

\lf
The striking feature of the von Neumann case is that the definition of algebraically minimal works for all $(\cA,\theta,p)$ and the minimalization procedure works for all normal $(\cA,\theta,p)$. On the contrary, it is one of the insights of Section \ref{topSEC}, that in the \nbd{C^*}case, for saying what a strict dilation is we have to limit ourselves to full dilations.

\bdefi \label{samindef}
A full dilation $(\cA,\theta,p)$ is \hl{strictly algebraically minimal}\index{dilation, weak!algebraically minimal!strictly}\index{algebraically minimal!strictly} if $\cA=\ol{\cA_\infty}^{stri}$.
\edefi

Also here, \phantomsection\index{algebraically minimal!strictly!impossibility of minimalization}if $(\cA,\theta,p)$ is not strictly algebraically minimal (since $\cA=\cls^{stri}\cA p\cA$ and $\cA_\infty\supset\cA_\infty p\cA_\infty$ this means $\ol{\cA_\infty p}\ne\cA p$), provided $\theta$ is strict we may (co)restrict to a triple $(\ol{\cA_\infty}^{stri},\theta^{\infty^{stri}},p)$. Like $\cA_\infty$ itself, the triple fulfills everything we expect from a dilation (of the same semigroup $T$). But, here, it is unclear, whether or not the strict closure of $\ol{\cA_\infty}^{stri}$ in $\cA$ is unital. (So the main scope of closing $\cA_\infty$, getting a unital algebra, might fail.) Moreover, even if $(\ol{\cA_\infty}^{stri},\theta^{\infty^{stri}},p)$ is a dilation, it need not be full. (Bhat's Example\index{Bhat's example} \ref{Bex} is an instance.) We discuss full dilations (and their minimalizations) in a later paragraph of this section.

\lf
In the third case,  instead of complaining about (the almost always) nonunitality of $\cA_\infty$, we utilize
%%%% BO new
% whooops sorry!
%%%% EO 
it to introduce a different topology on $\cA$ at least when $\cA$ is suitably generated by $\cA_\infty$, and call this \it{\nbd{\infty}algebraically minimal}. This effort is justified by two things. Firstly, starting with an arbitrary dilation, \nbd{\infty}algebraic minimality is obtainable (though the minimalized version does, usually, not sit in the original dilation), at least, for dilations of Markov semigroups over sufficiently nice monoids; see Corollary \ref{istrMricor}. Secondly, for \nbd{\infty}algebraically minimal dilations, Theorem \ref{minstriunithm} relates strongness to goodness. (A warning: The hypotheses of the former result make the latter trivial. But, the results also holds for the von Neumann case, where algebraic minimality can be obtained, and for the full strict case.)

\bdefi \label{damindef}
A dilation $(\cA,\theta,p)$ is \hl{\nbd{\infty}algebraically minimal}\index{dilation, weak!algebraically minimal!\nbd{\infty}strictly}\index{algebraically minimal!\nbd{\infty}strictly} if $\cA=M(\cA_\infty)$, that is, if $\cA_\infty$ is an ideal in $\cA$ and if the canonical homomorphism $\cA\rightarrow M(\cA_\infty)$ is an isomorphism.

A dilation is \index{dilation, weak!fully minimal|bf}\hl{\nbd{\infty}strict}\index{dilation, weak!\nbd{\infty}strict} if for each $t\in\bS$ the (co)restriction of $\theta_t$ to $\cA_\infty$ extends to a (necessarily unique) strict endomorphism $\theta^{\infty^M}_t$ of $M(\cA_\infty)$. 

A dilation is \phantomsection\hl{minimally strict}\index{dilation, weak!minimally strict} if it is \nbd{\infty}algebraically minimal and \nbd{\infty}strict.
\edefi

If a dilation $(\cA,\theta,p)$ is \nbd{\infty}strict, we call the resulting $(M(\cA_\infty),\theta^{\infty^M},p)$ (necessarily a dilation) the \hl{\nbd{\infty}algebraic minimalization}\index{dilation, weak!algebraic minimalization!\nbd{\infty}strict|bf}\index{algebraically minimal!\nbd{\infty}strictly!minimalization|bf} of the \nbd{\infty}strict dilation $(\cA,\theta,p)$.

\bthm \label{minstriunithm}
A minimally strict dilation $(\cA,\theta,p)$ is strong if and only if it is good.
\ethm

\proof
``Only if'' is Theorem \ref{sdilunithm}\eqref{DU1} (for arbitrary dilations). For the ``if''-direction, suppose the dilation is good, that is, suppose the $\theta_t(p)p$ form a unit for the associated superproduct system. (Recall that we assume the semigroup indexed by $\bS^{op}$ and the superproduct system by $\bS$.) By Observation \ref{puniob}, this means suppose that (for all fixed $t$)
\beqn{
\theta_{st}(p)\theta_t(p)p
~=~
\theta_{st}(p)p
}\eeqn
for all $s$. We have 
\beqn{
0
~=~
\theta_{st}(p)p-\theta_{st}(p)\theta_t(p)p
~=~
\theta_t(\theta_s(p))(\theta_t(\U)-\theta_t(p))p
~=~
\theta_t(\theta_s(p))\theta_t(\U-p)p
}\eeqn
for all $s$, hence,
\beqn{
0
~=~
\theta_t(\theta_s(pap))\theta_t(\U-p)p.
}\eeqn
Since $\theta_t$ is \nbd{\infty}strict and since the ball of $\cA$ is strictly generated by the ball of $\alg^*\theta_\bS(p\cA p)$, it follows that 
\beqn{
0
~=~
\theta_t(a)\theta_t(\U-p)p
}\eeqn
for all $a\in\cA$. In particular, $0=\theta_t(\U)\theta_t(\U-p)p=\theta_t(\U-p)p$. By Proposition \ref{strequivprop}, $(\cA,\theta,p)$ is strong.\qed

\bcor
~~Let $(\cA,\theta,p)$ be an \nbd{\infty}strict dilation. Then its \nbd{\infty}algebraic minimalization $(M(\cA_\infty),\theta^\infty,p)$ is strong if and only if $(\cA,\theta,p)$ is good.
\ecor

\proof
This follows from $\theta_t(p)p=\theta^\infty_t(p)p\in\cA_\infty$ and $\cA\supset\cA_\infty\subset M(\cA_\infty)$.\qed

\lf
The proof of Theorem \ref{minstriunithm} is based on that the ball of $\cA_\infty$ is dense in a certain topology in the ball of $\cA$ and that $\theta_t$ respects limits in that topology. It, therefore, remains valid for dilations that are algebraically minimal in one of the other two senses. And the corollary remains true whenever we do get a minimalization. (Whether the latter is granted by a theorem or required as a hypothesis, does not matter.) We collect this in the following observation and remark.

\bemp[Observation. The von Neumann case.] \label{algminvNob}
A normal and normally algebraically minimal dilation is strong if and only if it is good. And an arbitrary normal dilation is good if and only if its normal algebraic minimalization is strong.
\eemp

\bemp[Remark. The strict case.] \label{algminstriob}
Every strict and strictly algebraically minimal (therefore, by definition, full) dilation is strong if and only if it is good. And, \it{cum grano salis}, an arbitrary strict and full dilation is good if and only if its strict algebraic minimalization is strong. \it{Cum grano salis} means, that either $\ol{\cA_\infty}^{stri}$ is unital (then the statement is properly true) or $\ol{\cA_\infty}^{stri}$ is nonunital (then the statement applies when we allow dilations to act also on nonunital \nbd{C^*}algebra as we did temporarily in the beginning of this paragraph on algebraic minimality).
\eemp

\lf
We note that the \nbd{\infty}strict topology can be quite strong. Just look at an \nbd{E}semigroup $\theta$ on $\cA$; then $(\cA,\theta,\U_\cA)$ is a dilation of $\theta$ itself and $\cA_\infty=\cA=M(\cA_\infty)$ so that the \nbd{\infty}strict topology is just the norm topology. Nevertheless, this dilation is clearly \nbd{\infty}strict. In general, we do not know if a dilation is \nbd{\infty}strict. However, if the dilated semigroup is Markov and the monoid is ``nice'', then we can say more. We prepare with a general tool.

\blem \label{istrstrlem}
Let $(\cA,\theta,p)$ be a dilation. A bounded net $a_\lambda\in M(\cA_\infty)$ converges \nbd{\infty}strictly to $a$ if and only if $a_\lambda\theta_t(p)$ and $a_\lambda^*\theta_t(p)$ converge to $a\theta_t(p)$ and $a^*\theta_t(p)$, respectively, for all $t$.
\elem

\proof
Clearly, the condition is necessary. So assume the condition holds. Since the net is bounded, it is sufficient to control strict convergence on the total subset $\cA_\infty$ of elements of the form
\beq{ \label{itot}
\theta_{t_1}(pa_1p)\ldots\theta_{t_n}(pa_np)
~=~
\theta_{t_1}(p)\Bfam{\theta_{t_1}(pa_1p)\ldots\theta_{t_n}(pa_np)},
}\eeq
on which the condition guarantees convergence.
\qed

\bthm \label{istrMrithm}
Let $(\cA,\theta,p)$ be an \nbd{\infty}algebraically minimal dilation of a Markov semigroup over the opposite of a right-reversible monoid $\bS$. Then the dilation is minimally strict.
\ethm

\proof
Suppose the net $a_\lambda\in M(\cA_\infty)$ converges \nbd{\infty}strictly to $a$. Then $\lim_\lambda a_\lambda\theta_t(p)=a\theta_t(p)$ for all $t$, and likewise for the adjoint net.

Recall from Example \ref{incrpex} that, under the stated hypotheses, the projections $\theta_t(p)$ are, indeed, increasing over the directed index set $\bS$. Therefore, for $s$ and $t$ we may choose $r$ such that $r=s's\ge s$ and $r\ge t$, so that
\beqn{
\theta_s(a_\lambda)\theta_t(p)
~=~
\theta_s(a_\lambda)\theta_r(p)\theta_t(p)
~=~
\theta_s(a_\lambda)\,\theta_s\circ\theta_{s'}(p)\theta_t(p)
~=~
\theta_s(a_\lambda\theta_{s'}(p))\theta_t(p).
}\eeqn
 If we choose $\lambda$ big enough, then $a_\lambda\theta_{s'}(p)$ is close (in norm) to $a\theta_{s'}(p)$, and likewise for the adjoint net. By Lemma \ref{istrstrlem}, $\theta_s$ is \nbd{\infty}strict.\qed

\bcor \label{istrMricor} \index{dilation, weak!algebraic minimalization!\nbd{\infty}strict}\index{algebraically minimal!\nbd{\infty}strictly!minimalization}
Every dilation of a Markov semigroup over the opposite of a right reversible monoid possesses a (unique) \nbd{\infty}strict minimalization.
\ecor

We do not know if the hypothesis that the dilated semigroup is Markov can be dropped and weakened to strong dilation. (Since $\bfam{\wt{\cA}}_\infty$ is generated by $\cA_\infty$ and all $\theta_t(\wt{\U}-\U)=\wt{\U}-\theta_t(\U)$, we see that the component in $\cA\subset\wt{\cA}$ of $\bfam{\wt{\cA}}_\infty$ is generated by $\cA_\infty$ and all $\theta_t(\U)-\theta_s(\U)$. If $\theta$ is nonunital, we may expect that the multiplier algebra of $M\bfam{\bfam{\wt{\cA}}_\infty}$ differs quite bit from $M(\cA_\infty)+(\wt{\U}-\U)\C$; and if $\theta$ is unital, then we do not need these considerations, because, by Observation \ref{E0strongob}, the dilated semigroup necessarily \bf{is} Markov.)

This discussion shows that the strict topology induced by $\cA_\infty$ is a tricky one. Minimally strict dilations exist, if a dilation exists -- well, at least, under the conditions stated in Corollary \ref{istrMricor}. And minimally strict dilations, by Theorem \ref{minstriunithm}, relate the question of existence of strong dilation to the existence of good dilations. But, recall!, the two results, unfortunately, do not play together well, because Corollary \ref{istrMricor} is limited to Markov semigroups, for which the statement of Theorem \ref{minstriunithm} is vacuous because all dilations of Markov semigroups are strong, hence, good. And, as just, explained, the usual unitalization procedure, in this case does not seem to help.

Only in the von Neumann case, where we can guarantee minimalization, we can say in general:

\bthm \label{strgothm}
A normal CP-semigroup admits a normal strong dilation if and only if it admits a normal good dilation.
\ethm

The former can be obtained as the normal algebraic minimalization of the latter.

\lf
\subsection{\normalsize Primary dilations and compressions} \label{primSSEC}

The following result can be taken to motivate another natural notion that, when missing, indicates that a dilation cannot be considered a minimal one.

\blem \label{ptaulem}
Let $(\cA,\theta,p)$ be a dilation of a Markov semigroup over the opposite of a right-reversible monoid $\bS$. Then the net $\bfam{\theta_t(p)}_{t\in\bS}$ is an increasing contractive selfadjoint approximate unit for $\cA_\infty$.
\elem

\proof
For each $a$ of the form in \eqref{itot}, $\theta_t(p)a$ is eventually constant and, then, equal to its limit $a$. By boundedness of the net, controlling on this total subset is enough.\qed

\bcor \label{primcor}
In this situation, the \nbd{\infty}algebraically minimalized dilation is \hl{primary}\index{dilation, weak!primary!under the standing hypotheses}\index{primary!under the standing hypotheses}, in the sense that $\theta_t(p)$ converges strictly in $\cA=M(\cA_\infty)$ to $\U$. Obviously, it is an \nbd{E_0}dilation.
\ecor 

Anticipating Lemma \ref{striminlem} from Subsection \ref{minSEC}\ref{fmsSSEC}, we obtain that Corollary \ref{primcor} remains valid for 
%%%% BO 
%strictly minimal 
strictly algebraically minimal
%%%% EO
(hence, full) dilations, in the sense that $\theta_t(p)$ increases strictly in $\cA=\sB^a(E)$ to $\U$.  (Note that we do not require that the dilation is strict.) Recall, however, that we have no statement yet that would allow strict minimalization; and there is not much we can do about it.

What we can do something about is, starting from an arbitrary strict full dilation, making it primary -- provided $\theta_t(p)$ converges strictly in $\cA=M(\cA_\infty)$ to some $P$. Clearly (under the standing hypotheses), if $P$ exists, then it a projection satisfying $P\ge p$.

%%%% BO
\bdefi \label{comprdef}
We say a projection $P\in\cA$ \hl{compresses}\index{strong compression|see {compression}}\index{dilation, weak!compression}\index{compression!of a dilation}\index{dilation, weak!compression!strong}\index{compression!of a dilation!strong} a dilation $(\cA,\theta,p)$ if $P\ge p$ and if $(\cA,\theta,P)$ is a dilation of an \nbd{E}semigroup. If the dilation $(\cA,\theta,P)$ is strong, we say $P$ compresses $(\cA,\theta,p)$ \hl{strongly}.
\edefi
%%%% EO

Recall that $(\cA,\theta,P)$ being a dilation, means $\theta^P_t:=P\theta_t(P\bullet P)P$ is a CP-semigroup (here, an \nbd{E}semigroup) on $\cA^P:=P\cA P$.

\bob \label{compob}
~

\begin{enumerate}
\item \label{cob1}
By $P\ge p$ it follows that the \hl{compression} $(\cA^P,\theta^P,p)$ is a dilation of the same semigroup.

\vspace{1ex}
(On the contrary, if $(\wt{\cA},\wt{\theta},\wt{p})$ is the unitalization of a strong dilation $(\cA,\theta,p)$ as in Theorem \ref{wstrunithm}, then $(\wt{\cA},\wt{\theta},\U)$ is a dilation of $\theta$. But $\U$ does not compress the dilation $(\wt{\cA},\wt{\theta},\wt{p})$ because $\U\ngeq\wt{p}$ so that the dilation $(\cA,\theta,p)$ obtained by compressing with $\U$ acts on a different algebra. This situation is rather like the compression of CP-semigroups discussed in Subsection \ref{EXBexSEC}\ref{CPcompSSEC} and applied in the construction of Bhat's Example\index{Bhat's example} \ref{Bex}.)

\item \label{cob2}
A compression of a full dilation is full. (This is most obvious, by writing the full dilation as module dilation: If $P\in\sB^a(E)$, then $P\sB^a(E)P=\sB^a(PE)$.)

\item \label{cob3}
A (not necessarily strong) compression of a strong dilation is a strong dilation. (Indeed, we have $\theta_t(\U-p)p=0$ and $\U-p\ge P-p$, so, $\theta_t(P-p)p=\theta_t(P-p)\theta_t(\U-p)p=0$. Therefore, $\theta^P_t(P-p)p=P\theta_t(P-p)Pp=0$.)

\end{enumerate}
\eob

\noindent
This raises the question, when does a projection compress an endomorphism to an endomorphism?

\bprop \label{endocprop}
Suppose $\theta\colon\cA\rightarrow\cA$ is an endomorphism of the \nbd{C^*}algebra $\cA$ and $p\in\cA$ a projection. Then the compression $T:=p\theta(p\bullet p)p$ is homomorphic on $\cB:=p\cA p\subset\cA$ if and only if $p$ commutes with all elements of $\theta(\cB)$.
\eprop

\proof
By definition, $T$, obviously being a \nbd{*}map, is homomorphic on $\cB$, if and only if
\bmun{
p\theta(pa'p)\theta(pap)p
~=~
p\theta((pa'p)(pap))p
~=~
T((pa'p)(pap))
\\
~=~
T(pa'p)T(pap)
~=~
p\theta(pa'p)p\theta(pap)p
}\emun
for all $a,a'\in\cA$. Of course, having $p$ commute with all $\theta(pap)$, this condition is fulfilled. Conversely, assuming the condition is fulfilled, putting $a'=a^*$, we get $(\theta(pap)p)^*(\theta(pap)p)=(\theta(pap)p)^*p(\theta(pap)p)$. Generally, $c^*c=c^*pc$ implies $(c-pc)^*(c-pc)=c^*c-c^*pc=0$, so $c=pc$. If also $c^*=pc^*$, we get $cp=pcp=(pc^*p)^*=(c^*p)^*=pc$, so $p$ commutes with $c$ (and with $c^*$). Consequently, $p$ commutes with all $\theta(pap)$.\qed

\lf
Consequences for strongness:

\bcor \label{strcompreqcor}
If $(\cA^P,\theta^P,p)$ is the compression of a dilation $(\cA,\theta,p)$ by $P$, then
\beqn{
\theta_t(\U-p)p
~=~
\theta_t(\U-P)Pp
+
\theta^P_t(P-p)p.
}\eeqn
\ecor

\proof
$\theta^P_t(P-p)p=P\theta_t(P-p)Pp=\theta_t(P-p)Pp=\theta_t(\U-p)p-\theta(\U-P)Pp$.\qed

\lf
With this, for strong compressions we also get the opposite of Observation \ref{compob}\eqref{cob3}.

\bcor \label{strcomprcor}
A strong compression of a dilation is strong if and only if the original dilation is strong.
\ecor

Consequences for superproduct systems and goodness:

\bcor \label{comPPScor}
The superproduct system ${\sE^P}^\podot$ of a compression by $P$ is a superproduct subsystem of $\sE^\podot$ containing the elements $\xi_t=\xi^P_t$.
\ecor

\proof
From $\sE^P_t=\theta^P_t(p)P\cA Pp=P\theta_t(PpP)P\cA p=P\theta_t(p)\cA p\subset\theta_t(p)\cA p=\sE_t$, including the correct left action of $\cB=p\cA p$, and from
\beqn{
\theta^P_t(\theta^P_s(p)Pap)
~=~
P\theta_t\Bfam{P\Bfam{P\theta_s(p)ap}P}P
~=~
\theta_t\Bfam{P\theta_s(p)ap}P
}\eeqn
so that $u_{s,t}(x_s\odot y_t)=u^P_{s,t}(x_s\odot y_t)$ for all $x_s\in\sE^P_s$ and $y_t\in\sE^P_t$, it follows that ${\sE^P}^\podot\subset\sE^\podot$. Moreover, $\xi^P_t=P\theta_t(PpP)Pp=\theta_t(p)p=\xi_t$.\qed

\bcor \label{goodccor}
A compression of a dilation is good if and only if the original dilation is good.
\ecor

Consequences for algebraic minimality:

\bcor \label{amcccor}
A projection compressing an algebraically minimal dilation is central.% Moreover, the compression is algebraically minimal, again.
\ecor

\proof
We have $\theta_t(pap)=\theta_t(PpapP)$, so $P\theta_t(pap)=\theta_t(pap)P$ for all generating elements $\theta_t(pap)$ of $\cA$.\qed

\lf
We postpone to Corollary \ref{compamincor}, that such a compression is algebraically minimal, again.

\lf
Recall from Observation \ref{E0strongob} and Example \ref{E0weakex} that only a strong compression of an \nbd{E_0}dilation is necessarily an \nbd{E_0}dilation and that there exist compressions of \nbd{E_0}dilations that are not \nbd{E_0}dilations. Of course, an \phantomsection\label{E0c}\hl{\nbd{E_0}compression} (that is, a compression to an \nbd{E_0}dilation) is, by Proposition \ref{Mcharprop}, a strong compression.

\lf
The preceding consequences of Proposition \ref{endocprop} were general. We now return to the situation of Lemma \ref{ptaulem} and Corollary \ref{primcor}:

\bcor \label{primcompcor}
If $(\cA,\theta,p)$ is a full Markov dilation over the opposite of an Ore monoid and if $\theta_t(p)$ increases strictly to $P$, then $P\ge p$ is a projection compressing $(\cA,\theta,p)$. The compressed dilation $(\cA^P,\theta^P,p)$, which is strict if $(\cA,\theta,p)$ is, is a full primary, hence, \nbd{E_0}dilation. In particular, $P$ compresses $(\cA,\theta,p)$ strongly.

Moreover, if $(\cA,\theta,p)$ is primary (so that $P=\U$), then every compression is full and primary, hence, a strong \nbd{E_0}compression. It is strict, if $(\cA,\theta,p)$ is.
\ecor

\bob \label{primob}
The strict dilation constructed in Theorem \ref{Oreindthm} from a product system over an Ore monoid with unital unit is, clearly, primary. More generally, a module dilation $(E,\vt,\xi)$ of a Markov semigroup over the opposite of an Ore monoid can be recovered as the inductive limit over the unit of its product system as in Theorem \ref{Oreindthm} if and only if the dilation is primary. (This is simply so, because the canonical injection $k_t$ into the inductive limit $\sE$ sends the member $\sE_t$ of the product system onto $\vt_t(\xi\xi^*)\sE$, and $\sE$ is spanned by all  these ranges. See also the proof of Theorem \ref{Ocomprthm} below.)
\eob

\lf
When, in the \nbd{C^*}case, we are neither in the \nbd{\infty}algebraic minimal case nor in the full case, then we have no topology around in which $\theta_t(p)$ can increase to something; actually, without the standing hypotheses we cannot even say that the $\theta_t(p)$ form an increasing net.

The situation is much nicer in the von Neumann case. Here, we need neither that the dilated semigroup is Markov nor that the monoid is right reversible. Recall that a family $\bfam{p_t}$ of projections in a von Neumann algebra $\cA$ has a \hl{supremum} $P:=\bigvee_tp_t$ in $\cA$ characterized uniquely by the properties that $P$ is a projection fulfilling $P\ge p_t$ for all $t$ and that $Q\ge P$ for any other projection fulfilling this condition. 

%%%% BO
% Hi Michael, this is standard material found in any von Neumann algebras text, and I think we don't have anything new to say here. So I am commenting out, and if you disagree you can uncomment. --Orr

% The test that previously appeared: 

\brem\label{infprem}
It is well-known that supremum $\vee$ and infimum $\wedge$ (defined in an analogue manner) turn the set of projections in a von Neumann algebra into a \it{complete lattice} relative to the partial order $\le$. However, standard proofs of existence (first, of $\bigwedge_tp_t$ and, then, putting $\bigvee_tp_t=\U-\bigwedge_t(\U-p_t)$), such as Takesaki \cite[Proposition V.1.2]{Tak02}, define $\bigwedge_tp_t$ to be the projection onto the intersection of all $p_tH$ and, then, show it is in $\cA''$. For later use (Theorem \ref{vNE0compthm}), we wish to have a concrete approximation. For two, $p\wedge q$ can be computed as strong limit $\lim_{n\to\infty}p(qp)^n$, which is in $\cA$; see Proposition \ref{minpqprop}. This lifts to finitely many; and the net of all the \it{minima} over all finite subchoices is decreasing over inclusion, so that it has a strong limit in $\cA$, too.
\erem
%%%% EO

\bdefi \label{vNprimdef}
A dilation $(\cA,\theta,p)$ on a von Neumann algebra $\cA$ is \hl{primary}\index{dilation, weak!primary}\index{primary} if $\bigvee_t\theta_t(p)=\U$.
\edefi

\bthm \label{vNprimthm}
Let $(\cA,\theta,p)$ be a dilation on a von Neumann algebra $\cA$.
\begin{enumerate}
\item \label{vNp1}
If $(\cA,\theta,p)$ is normally algebraically minimal, then it is a primary.

\item \label{vNp2}
If $(\cA,\theta,p)$ is a primary dilation of a Markov semigroup, then it is an \nbd{E_0}dilation.

\item \label{vNp3}
If $(\cA,\theta,p)$ is normal, then $\bigvee_t\theta_t(p)$ compresses $(\cA,\theta,p)$ to a (normal) primary dilation.\index{dilation, weak!primary!``primarization'' of a normal dilation}\index{primary!``primarization'' of a normal dilation}

\item \label{vNp4}
If $(\cA,\theta,p)$ is primary, then so is any of its compressions.
\end{enumerate}
\ethm

\proof
Put $P:=\bigvee_t\theta_t(p)$.

\eqref{vNp1}.~
$P\ge\theta_t(p)$ for all $t$. Therefore, $P$ acts as $\U$ on all elements of the form in \eqref{itot}, hence, on all elements in $\cA=\ol{\cA_\infty}^s$.

\eqref{vNp2}.~
Since $\U\ge\theta_t(p)$, we have $\theta_s(\U)\ge\theta_{ts}(p)\ge\theta_t(p)$ (see Example \ref{incrpex}!) for all $t$. So, $\theta_s(\U)\ge\U$.

\eqref{vNp3}.~
Since $\theta_t(P)\le P$ (supremum over a smaller set; here we use normality of $\theta_t$),we have $\theta_t(P\cA P)=P\theta_t(P)\theta_t(\cA)\theta_t(P)P\subset P\cA P$. So, $\theta$ leaves $\cA^P$ invariant, and the compression $\theta^P$ coincides with the (co)restriction of $\theta$.

\eqref{vNp4}.~
Let $Q$ denote a compressing projection. Then $\theta^Q_t(p)=Q\theta_t(p)$, so $\bigvee_t\theta^Q_t(p)=Q\bigvee_t\theta_t(p)=Q$.\qed

\bcor \label{amE0ccor}
Every compression of a normally algebraically minimal Markov dilation is an \nbd{E_0}com\-pression and, therefore, strong.
\ecor

\bob \label{primpPSob}
In the situation of \eqref{vNp3}:
\begin{enumerate}
\item \label{pPSob1}
We have $\sE_t^P:=\theta^P_t(p)\cA p=\theta_t(p)\cA p=\sE_t$. So, $(\cA,\theta,p)$ and its ``\it{primarization}'' have the same superproduct system.

\item \label{pPSob2}
In Example \ref{E0weakex}, it is easy to show that $\bigvee_t\vt_t(q)=p$. (The span of all $\vt_t(q)f$ contains all $\I_{\RO{s,t}}f$ ($0\le s<t$), which form a total subset of $L^2(\R_+)$.) On the other hand, $\vt_t(p)=\I_{\RO{t,\infty}}\ngeq\I_{\R_+}=p$. Therefore, $p$ does not compress $(\sB(H),\vt,q)$ to its ``\it{primarization}'' $(\sB(G),\vt^p,q)$ strongly: If the dilation $(\cA,\theta,p)$ is not Markov, then the compression by $P$ need not be strong.
\end{enumerate}
\eob

% \bob
% In general, it is unclear if a non-primary normal dilation can be compressed by $\bigvee_t\theta_t(p)$ to a primary dilation. Only in the one-parameter case, or, more generally, if the ``partitions'' $\bJ_t$ are directed by refinement, we can the non-Markov case to the Markov case via the first main theorem of Skeide \cite{Ske08}.)
% \OW[Check!]{Check!}
% \eob

\lf
\subsection{\normalsize Incompressible dilations} \label{compSSEC}

Non-primary dilations can (at least, under suitable assumptions) be compressed to primary dilations and algebraically minimal dilations are primary. But, sometimes, primary dilations may be compressed further. (Examples \ref{hypexex} and \ref{discex} are primary dilations of the identity semigroup $T_t=\id_\C$ on $\C$. But $T$, being an endomorphism semigroup, is also its own dilation -- and, thus, the result of compressing the given primary dilation (with the same $p$) to the identity dilation.) Surely, compressible dilations are not candidates for what what we would call \it{minimal}.

\bdefi \label{incomprdefi}
A dilation is (\hl{strongly}) \hl{incompressible}\index{dilation, weak!incompressible@(strongly) incompressible}\index{incompressible@(strongly) incompressible} if the only (strong) compression it admits is by $P=\U$. (See Definition \ref{comprdef}.)
\edefi

Among all compressions of a given dilation, is there a smallest incompressible one? Well, for arbitrary compressions, we do not know. But, for strong compressions, the answer, in the von Neumann case, is yes.

\bthm \label{vNE0compthm}
Let $(\cA,\theta,p)$ be a normal dilation on a von Neumann algebra $\cA$, and denote $P:=\bigwedge Q$ where $Q$ runs over all projections that compress $(\cA,\theta,p)$ strongly. Then $P$ is the unique smallest projection that compresses $(\cA,\theta,p)$ strongly to the unique strongly incompressible dilation among all strong compressions of $(\cA,\theta,p)$.
\ethm

\proof
Recall that a projection $Q\ge p$ defines a strong compression if and only if $Q$ commutes with all $\theta_t(Q\cA Q)$ (Proposition \ref{endocprop}) and $\theta_t(\U-Q)Q=0$ for all $t$ (Proposition \ref{strequivprop}).

If, for $i=1,2$, we have a projection $Q_i\ge p$ that commutes with all $\theta_t(Q_i\cA Q_i)$, then also the minimum $q:=Q_1\wedge Q_2=\lim_{n\to\infty}Q_1(Q_2Q_1)^n$ (see the proof of Proposition \ref{minpqprop}) majorizes $p$ and commutes with all $\theta_t(q\cA q)$. So, all compressed maps $\theta^q_t:=q\theta_t(q\bullet q)q$ are endomorphisms of $\cA^q:=q\cA q$. This turns over to $Q_1\wedge\ldots\wedge Q_n$ and passing through the limit over the decreasing net of all finite choices, taking also into account normality of $\theta_t$, we get the same statement for $P$.

Now, if $\theta_t(\U-Q_i)Q_i=0$ for all $t$, then from
\beqn{
\theta_t(Q_i)\theta_t(\U)q
~=~
\theta_t(Q_i)q
~=~
\theta_t(Q_i)Q_iq
~=~
\theta_t(\U)Q_iq
~=~
\theta_t(\U)q
}\eeqn
we get
\beqn{
\theta_t(\U-Q_1(Q_2Q_1)^n)q
~=~
\theta_t(\U)q-\theta_t(\U)q
~=~
0.
}\eeqn
By normality, $\theta_t(\U-q)q=0$. Now, by normality ($:=$order continuity!) of $\theta_t$ and joint strong continuity of the product on bounded subsets, we get $\theta_t(\U-P)P=0$.

Altogether, also $P$ defines a strong compression. Of course, by the usual standard arguments, $P$ is minimal and unique.\qed

\brem \label{sgrem}
If, in the definition of $Q$, we would have let run the infimum over \bf{all} compressions, we still would get a family of homomorphisms $\theta^P_t$. But, only the condition $\theta_t(\U-P)P=0$ guarantees that this family is a semigroup. In fact, we do not know whether $\theta^{Q_1\wedge Q_2}_t$ always defines a semigroup, or not. (Bhat's result \cite[Corollary 2.3]{Bha02} in the one-parameter case, does not look applicable easily to $Q_1\wedge Q_2$, not speaking about the question whether it could be generalized to other monoids.)

By Theorem \ref{sdilunithm}\eqref{DU2} and by Observation \ref{puniob}, also just checking the \it{goodness condition} $\theta_{st}(P)\theta_t(\U-P)P=0$ would be enough to get a dilation. But, we do not know if $\theta_{st}(Q_i)\theta_t(\U-Q_i)Q_i=0$ implies $\theta_{st}(q)\theta_t(\U-q)q=0$. So, we do not know if the infimum over \it{good compressions} $Q$ results in a good compression $P$.
\erem

\lf
Recall from Corollary \ref{amcccor} that a projection compressing an algebraically minimal dilation is necessarily central. In this case, we have $Q_1\wedge Q_2=Q_1Q_2$, so that the proof of Theorem \ref{vNE0compthm} simplifies drastically. Another consequence of Corollary \ref{amcccor}:

\bcor \label{amfinccor}
A normally algebraically minimal full dilation is incompressible.
\ecor

\proof
By Theorem \ref{vNcontthm}\eqref{vNct2}, the only central projection $Q\ge p$ is $\U$.\qed

\bthm \label{vNincfthm}
A normal incompressible dilation is full if and only if its superproduct system is a product system.
\ethm

\proof
Every (normal or strict) full 
%%%% BO 
%%%% EO
dilation has a product system as superproduct system. On the other hand, if the normal dilation is not full, then, by Theorem \ref{vNcontthm}, it is at least semifull, and if a semifull but not full dilation has a product system as superproduct system, then, by Theorem \ref{p-PSthm}, it is compressible.\qed

\lf
As we said, we do not know, if an arbitrary dilation can always be compressed to an incompressible one. But for full primary dilations we can give a sufficient criterion for incompressibility in terms of the superproduct system of the dilation (the latter, necessarily a product system, because the dilation is full).

\bthm \label{fpincthm}
A normal primary full dilation with a product system that is generated by its elements $\xi_t$, is incompressible.
\ethm

\proof
If the normal dilation $(\cA,\theta,p)$ is primary and full, then $E=\cls^s\sE_\bS$, where $\sE_\bS:=\bigcup_{t\in\bS}\sE_t$. Therefore, $\cA=\sB^a(E)=\cls^sEpE^*=\cls^s\sE_\bS\sE_\bS^*$.

If $P$ is a non-trivial projection, then $\cls^s\sE_\bS\sE_\bS^*=\cA\ne P\cA P=\cls^s(P\sE_\bS)(P\sE_\bS)^*$. Therefore, for at least one $t$ we have $\sE_t\ne P\sE_t$.

Now, if $P$ is an (arbitrary) projection compressing the dilation, then the $P\sE_t$ form the product system of the (full!) compression, by Corollary \ref{comPPScor}, a (super)product subsystem of $\sE^\odot$ containing all $\xi_t$. If $P$ was non-trivial, then the product subsystem was proper. Therefore, $P=\U$.\qed

\lf
Observe that the dilation is \bf{not} required strong. In the case of a (necessarily strong) full primary dilation of a Markov semigroup over the opposite of an Ore monoid, we can even show necessity of this criterion; and we can show most of it both in the \nbd{C^*}case and in the von Neumann case. So, we switch to the situation we sketched already in Observation \ref{primob}, and write a full dilation as module dilation $(E,\vt,\xi)$ and we will assume that we compute its product system as $\sE_t:=\vt_t(\xi\xi^*)E$ with unit $\xi_t:=\vt_t(\xi\xi^*)\xi$ (which, as an element of $E$, is simply $\xi$). The product is $u_{s,t}\colon x_s\odot y_t\mapsto\vt_t(x_s\xi^*)y_t$.

\bthm \label{Ocomprthm}
Let $(E,\vt,\xi)$ be a strict primary (hence, $E_0$) module dilation of a Markov semigroup over the opposite of an Ore monoid $\bS$. Then the formula
\beqn{
q_t
~=~
\vt_t(\xi\xi^*)q
}\eeqn
establishes an order isomorphism between
\begin{enumerate}
\item
the set of projections $q\in\sB^a(E)$ compressing the dilation $(E,\vt,\xi)$, and

\item
the set of projection morphisms $q^\odot=\bfam{q_t}_{t\in\bS}$ of $\sE^\odot$ fulfilling $q_t\xi_t=\xi_t$.
\end{enumerate}
\ethm

\noindent
Recall from Corollary \ref{primcompcor}, that every compression is an \nbd{E_0}compression (hence, strong) to a strict primary full \nbd{E_0}dilation.
 
\bcor
$(E,\vt,\xi)$ is incompressible if and only if $\sE^\odot$ has no proper complemented product subsystem containing the unit $\xi^\odot$.
% \OW[MICHAEL]{Don't we need to restrict only to complementable subsystems, or perhaps to the von Neumann case? See also the Remark 21.31 below. --Orr}
\ecor

In this case, we say $\sE^\odot$ is \hl{generated} by $\xi^\odot$, hence, by the GNS-subproduct system of the dilated Markov semigroup contained in $\sE^\odot$ as the subproduct subsystem generated ($=$spanned) by $\xi^\odot$; see Observation \ref{uGNSsupob}. The fact that $\sE^\odot$ need not be  \hl{spanned} by $\xi^\odot$ (that is, $\xi^\odot$ may span only a proper superproduct subsystem of $\sE^\odot$; see Theorem \ref{supintthm}), will still keep us busy in the remainder. It represents a key problem regarding minimality.

\lf
\proof[Proof of Theorem \ref{Ocomprthm}.~]
Let us first be more precise about the statements we made in Observation \ref{primob}. Suppose we have any strict module dilation $(E,\vt,\xi)$ of a Markov semigroup, and construct its product system $\sE^\odot$ and unit $\xi^\odot$ (as indicated before the theorem). Then the inductive limit construction proving Theorem \ref{Oreindthm} provides us with $\sE$ and a left dilation $v_t$ of $\sE^\odot$ determining a dilation $(\sE,\vt^\xi,\xi)$ of the same Markov semigroup. Since the canonical maps $k_t\colon\sE_t\rightarrow\sE$ of the inductive limit act as $k_tx_t:=\xi x_t$, where we defined the product $xy_t:=v_t(x\odot y_t)$, since $\vt^\xi_t(\xi\xi^*)\sE=k_t\sE_t$, and since the $k_t\sE_t$, by definition, span the inductive limit $\sE$, the dilation $(\sE,\vt^\xi,\xi)$, clearly, is primary. If also $(E,\vt,\xi)$ is primary, then $E\ni\sE_t\ni x_t:=\vt_t(\xi\xi^*)x\mapsto k_tx_t\in\sE$ extends to a unique unitary $E\rightarrow\sE$ establishing a unitary equivalence between the two dilations $(E,\vt,\xi)$ and $(\sE,\vt^\xi,\xi)$. We, henceforth, will identify $(E=\sE,\vt=\vt^\xi,\xi)$ and write only $(\sE,\vt,\xi)$, assuming our primary dilation has been constructed by means of Theorem \ref{Oreindthm}.

Let $q\ge p$ be a projection compressing  $(\sE,\vt^\xi,\xi)$ to an \nbd{E_0}dilation. We know, that happens if and only if $q$ is increasing and if it commutes with all $\vt_t(q\sB^a(E)q)$. By the latter property, $q$ maps $\sE_t=\vt_t(\xi\xi^*)\sE$ into $\sE_t$. So we may define a projection $q_t\in\sB^a(\sE_t)$ by (co)restricting $q$ to $\sE_t$. By the same commutation property, $q_t$ is actually in $\sB^{a,bil}(\sE_t)$. Put $\sF_t:=q_t\sE_t$. We know, $u_{s,t}(x_s\odot y_t)=\vt_t(x_s\xi^*)y_t$ ($x_s\in\sE_s$, $y_t\in\sE_t$) is total in $\sE_{st}$. So, (very pedestrian, to be clear in each single step, also taking into account that $q$ commutes with with expressions like $\vt_t(\xi\xi^*)=\vt_t(q\xi\xi^*q)$ and $\vt_t(qx\xi^*)=\vt_t(qx\xi^*q)$)
\bmun{
q_{st}u_{s,t}(x_s\odot y_t)
~=~
\vt_{st}(\xi\xi^*)q\vt_t(x_s\xi^*)y_t
~=~
q\vt_{st}(\xi\xi^*)\vt_t(x_s\xi^*)y_t
~=~
q\vt_t(x_s\xi^*)y_t
\\
~=~
q\vt_t(q)\vt_t(x_s\xi^*)y_t
~=~
q\vt_t(qx_s\xi^*)y_t
~=~
\vt_t(qx_s\xi^*)qy_t
~=~
\vt_t(q_sx_s\xi^*)q_ty_t
~=~
%%%% BO
u_{s,t}(q_sx_s\odot q_ty_t)
%%%% EO
}\emun
is total in $\sF_{st}$. It follows that the $\sF_t$ form a product subsystem $\sF^\odot$ of $\sE^\odot$ and that the $q_t$ form a projection morphism $q^\odot$ for $\sE^\odot$ projecting onto $\sF^\odot$.

Since $(E,\vt,\xi)$ is primary, we recover $q$ from $q^\odot$ as strong limit $q=\lim_tq_t$ (when considering $q_t$ as the maps on $\sE$ as which they, originally, have been defined). Clearly, if we have two such projections, $q$ and $q'$, then $q\ge q'$ if and only if $q_t\ge q'_t$ for all $t\in\bS$. So, the map $q\mapsto q^\odot$ is an order isomorphism 
%%%% BO
% I commented out the following, I don't understand it (what is affine in this context?) and I don't think it is used. 
% (unavoidably affine by the definition of order) 
%%%% EO
onto its range. The last remaining question, is this map surjective?, takes its positive answer just from the universal property of the inductive limit. (Indeed, if $q^\odot$ is a contractive morphism fulfilling $q_t\xi_t=\xi_t$, then the family $a_t:=k_tq_t\in\sB^a(\sE_t,\sE)$ is contractive and, for $t=rs\ge s$, satisfies
\beqn{
a_t\gamma_{t,s}x_s
~=~
\xi q_{rs}(\xi_rx_s )
~=~
\xi(q_r\xi_r)(q_sx_s)
~=~
\xi(q_sx_s)
~=~
a_sx_s.
}\eeqn
So, there is a unique contraction $q\in\sB^r(\sE)$ satisfying $qk_t=a_t=k_tq_t$. Obviously, if $q^\odot$ is a projection morphism, then $q$ is a projection. One easily verifies that $v_t(q\odot q_t)v_t^*=q$, by testing it on the total subset of vectors $(k_sx_s)y_t$. So, clearly, $q$ is increasing and
\bmun{
~~~~~~
q\vt_t(qaq)(xy_t)
~=~
q\vt_t(q)\vt_t(aq)(xy_t)
~=~
q(aqx)y_t
\\
~=~
(qaqx)(q_ty_t)
~=~
\vt_t(qaq)((qx)(q_ty_t))
~=~
\vt_t(qaq)q(xy_t).
~~~~~~
}\emun
It follows that $q$ defines a compression of $(E,\vt,\xi)$ to an \nbd{E_0}dilation and that $q\mapsto q^\odot$ gives back the $q^\odot$ we started with.)\qed

\lf
Needless to say that all this goes through also in the von Neumann case. (Existence of $q$ when given $q^\odot$ in the last part of the proof, is actually a bit simpler, because we know that the increasing net $v_t(\xi\xi^*\odot q_t)v_t^*$ of projections converges strongly to some projection $q$.)

\brem
In the \nbd{C^*}case, we think it might be interesting to replace the projection morphisms $q^\odot$ \bf{on} the product system of the to-be-compressed dilation, by (not necessarily adjointable) embedding morphisms $i^\odot$ from another product system $\sF^\odot$ (with unit $\zeta^\odot$) into $\sE^\odot$ (sending $\zeta^\odot$ to $\xi^\odot$). We might, then, ask how the two primary dilations constructed from $(\sF^\odot,\zeta^\odot)$ and from $(\sE^\odot,\xi^\odot)$ are related. We do not follow this, here.
\erem

\lf
Observe that if $P$ compresses a dilation $(\cA,\theta,p)$, then $P\cA_\infty P=(P\cA P)_\infty$. (See the proof of Corollary \ref{amcccor}.) Likewise, for the strong closures in the von Neumann case. In general, for a projection $P\in M(\cA_\infty)$, we have $M(P\cA_\infty P)=PM(\cA_\infty)P\subset M(\cA_\infty)$. Therefore, in all three cases (\nbd{\infty}, strictly, and normally), compression works together nicely with algebraic minimality.

\bcor \label{compamincor}
If a dilation $(\cA,\theta,p)$ is algebraically minimal, then so is any of its compressions.
\ecor

And together with Theorem \ref{vNE0compthm}:

\bcor \label{vNincmincor}
Every CP-semigroup on a von Neumann algebra that possesses a normal dilation, possesses an algebraically minimal strongly incompressible normal dilation.
\ecor

% \brem
% We just wish to mention that the proof of Theorem \ref{vNE0compthm} simplifies considerably, if the (normal) dilation we wish to maximally compress is also algebraically minimal. In this case, a projection $Q\ge p$ that compresses $(\cA,\theta,p)$, commuting with all $\theta_t(Q\cA Q)\subset\theta_t(\cB)$, has to commute with all $\ol{\cA_\infty}^s=\cA$. So, $Q$ is a central projection. Therefore, in this case, $Q_1\wedge Q_2=Q_1Q_2$ and all the verifications of the properties of $Q_1\wedge Q_2$ and $\bigwedge Q$ are much easier. Still, we do not know whether compressing with the infimum over \bf{all} (not only strong) compression, leads to a semigroup.
% \erem

\lf
\subsection{\normalsize Fully minimal dilations} \label{fmsSSEC}

In Examples \ref{hypexex} and \ref{discex}, $p$ is increasing and $\cB=p\C$, so, $\cA_\infty=C^*\CB{\theta_t(p)\colon t\ge0}$ is commutative; hence, the algebraic minimalization cannot be full. (As said in the beginning of Subsection \ref{minSEC}\ref{compSSEC}, it is compressible by $p$, because the dilated semigroup is an \nbd{E_0}semigroup.) These are examples for both the \nbd{C^*}case and the normal von Neumann case. It turns out that in the existing literature (mainly on the one-parameter case), the definitions of minimality in use refer, implicitly or explicitly, to dilations that are algebraically minimal \bf{and} full. We speak of \hl{fully minimal}\index{dilation, weak!fully minimal} dilations.%
\footnote{ \label{ArvminFN}
In Chapter 8 of the book \cite{Arv03}, Arveson examines algebraically minimal strong dilations of normal strongly continuous one-parameter CP-semigroups on von Neumann algebras.
%%%% BO 
% Several changes. 
In \cite[Definition 8.3.5]{Arv03} he calls such a dilation \it{minimal} if it is full, and  in \cite[Definition 8.9.3]{Arv03} he calls a dilation \it{minimal over} the dilated CP-semigroup if it is strongly incompressible. 
%%%% BO new
% Regarding insertion of "strongly": correct. In this context, Arveson's compressions are only via increasing projections (see Definition 8.9.1, which is used in Definition 8.9.3, in his book). 
%%%% EO (new)
He then proves in \cite[Theorem 8.9.5]{Arv03} that a (strong) dilation of a Markov semigroup is strongly incompressible if and only if it is algebraically minimal and full (that is, what we call \it{fully minimal}), and that these conditions are equivalent to Bhat's definition of minimality in \cite{Bha96} for $\sB(H)$. (We recover these results in Subsection \ref{1-p-SSEC}.)
%In \cite[Section 8.9]{Arv03}, he defines \it{minimal} dilation of a Markov semigroup as being incompressible, and shows in \cite[Theorem 8.9.5]{Arv06} that this is equivalent to being algebraically minimal and full (expressed using the central cover of $p$ as in our Theorem \ref{vNcontthm}\eqref{vNct2}), and that is is equivalent to Bhat's definition of minimality in \cite{Bha96} for $\sB(H)$; all this in the one-parameter case.
%%%% EO 
Our next section shows that already in the discrete two-parameter case, incompressible and fully minimal are not the same. Bhat's definition cannot even be formulated directly (making reference to the total order of the monoid); a weaker version of it could be formulated for any monoid, and it is clear that this formulation guarantees that an algebraically minimal dilation is fully minimal.
}

Depending on which version of algebraic minimality we choose, \it{a priori} there are three versions: Fully (normally)(strictly)(\nbd{\infty})minimal. Recall that the definition of \it{full} means the same, $\cA=\sB^a(E)$ with $E:=\cA p$, in the \nbd{C^*}case and in the von Neumann case, and that also in the von Neumann case $\sB^a(E)$ is spanned strictly by $\sF(E)$. The difference between Definitions \ref{namindef} and \ref{samindef} remains, however, in which topology $\cA_\infty$ generates $\sB^a(E)$.

\bthm \label{fminthm}
The full dilation $(\cA,\theta,p)$ (on the von Neumann algebra $\cA$) is strictly (normally) algebraically minimal if and only if (the strong closure of) $E_\infty:=\cA_\infty p$ coincides with $E:=\cA p$.
\ethm

\proof
In the \nbd{C^*}case, we have $\ol{\cA_\infty}^{stri}p=\ol{\cA_\infty p}=\cA_\infty p$. Therefore, if $E_\infty\subsetneq E$, then $\ol{\cA_\infty}^{stri}$ cannot be equal to $\cA=\sB^a(E)$, while if $E_\infty=E$, then $\ol{\cA_\infty}^{stri}\supset\sF(E_\infty)=\sF(E)$ coincides with $\sB^a(E)=\cA$.

In the von Neumann case, the Hilbert module $E_\infty$ will usually not be strongly closed. $\sF(E_\infty)$ will not coincide with $\sF(E)$ and, in this case, also the strict closure of $\sF(E_\infty)$ in $\sB^a(E)$ (!) is not all of $\sB^a(E)$, because $\ol{\sF(E_\infty)}^{stri}E\subset E_\infty$. But the strong closure of $\sF(E_\infty)$ contains $\sF(\ol{E_\infty}^s)$ and, therefore, coincides with $\sB^a(\ol{E_\infty}^s)$. From this, the statement follows.\qed

\lf
Clearly, if, in the normal von Neumann case, $(\ol{\cA_\infty}^s,\theta^{\infty^s},p)$ is full, then it is the strong compression of $(\cA,\theta,p)$ with $\id_{\ol{E_\infty}^s}$.

\bcor \label{incmnfcor}
If a normal module dilation $(E,\vt,\xi)$ is not (normally) fully minimal but strongly incompressible (\it{a fortiori}, if it is incompressible), then its algebraic minimalization is not full.
\ecor

And as a special case of Theorem \ref{vNincfthm}:

\bcor \label{algminincfullcor}
If a normal normally algebraically minimal dilation is incompressible, then it is full if and only if its superproduct system is a product system.
% \OW[MICHAEL]{I am confused please clarify how this follows. Theorem \ref{vNincfthm} is about incompressible, not strongly incompressible. --Orr}
\ecor

We refer to both situations in the theorem as just \phantomsection\hl{fully minimal}\index{fully minimal}\index{dilation, weak!fully minimal|bf}, while \hl{fully \nbd{\infty}minimal}\index{fully minimal!fully \nbd{\infty}minimal}\index{dilation, weak!fully minimal!fully \nbd{\infty}minimal} means that $M(\cA_\infty)=\cA=\sB^a(E)$. Recall that for a general dilation, in the \nbd{C^*}case $M(\cA_\infty)$ need not even be a subset of $\cA$. In the von Neumann case (like for every multiplier algebra of a faithfully represented \nbd{C^*}algebra), we have $\cA_\infty\subset M(\cA_\infty)\subset\ol{M(\cA_\infty)}^s=\ol{\cA_\infty}^s\subset\cA$. At least in the fully minimal case, we get $M(\cA_\infty)\subset\cA$ also in the \nbd{C^*}case.

\blem \label{striminlem}
Suppose the full dilation $(\cA,\theta,p)$ is strictly algebraically minimal. Then sending ~$\SB{a\colon a_\infty\mapsto aa_\infty}\in M(\cA_\infty)=\sB^a(\cA_\infty)$~ to ~$\SB{a_p\colon a_\infty p\mapsto (aa_\infty)p}\in\cA=\sB^a(\cA_\infty p)$,~ defines an injective unital homomorphism $M(\cA_\infty)\rightarrow\sB^a(E)$, which is, clearly, strictly continuous.
\elem

\proof
Since $E=\cA_\infty\subset\cA_\infty$, we may restrict $a\in\sB^a(\cA_\infty)$ to arguments from $E$. Since $a\in\sB^a(\cA_\infty)$ is right linear, it satisfies $a(a_\infty p)=(aa_\infty)p\in E$, so that we may define $a\mapsto a_p$ as stated. Since the inner product of $E=\cA_\infty p$ is just the restriction of the inner product of $\cA_\infty$, the map is a unital homomorphism. Since $EE^*\subset\cA_\infty$ generates $\sB^a(E)\supset\cA_\infty$ strictly and since $a(xy^*)=(a_px)y^*$, we that $a_p=0$ implies $a_p\sB^a(E)=\zero$ implies $a\cA_\infty=\zero$ implies $a=0$.\qed

\bcor \label{strimincor}
If $(\cA,\theta,p)$ is a fully minimal dilation, then  $M(\cA_\infty)\subset\cA$.
\ecor

\lf
Being fully \nbd{\infty}minimal (that is, $M(\cA_\infty)=\sB^a(E)$),  is a tough requirement for a dilation. (We have the somewhat trivial Example \ref{minnuniex}, and probably the same methods show that all dilations of semigroups on $M_n$ as in Example \ref{dwnsex} from minimal coisometric dilations will do.) In any case, if $M(\cA_\infty)=\sB^a(E)$, then under the standing hypotheses (Markov over the opposite of right reversible) we know the dilation is \nbd{\infty}strict; while we do not know if it is also strict.

\lf
\subsection{\normalsize (Non-)uniqueness of minimal dilations} \label{nuniminSSEC}

It is a striking feature of the one-parameter case (see also Footnote \ref{ArvminFN}) that the fully minimal dilation (always exists and) is unique (in our terminology, up to unitary \nbd{\xi}intertwining equivalence among module dilations $(E,\vt,\xi)$). Nothing like this survives to the general case.\index{minimal!dilation!non-uniqueness of}\index{fully minimal!non-uniqueness of fully minimal dilations}\index{dilation, weak!fully minimal!non-uniqueness of}

\bex \label{minnuniex}
Discrete one-parameter CP-semigroups $T=\bfam{T_n}_{\in\N_0}$ on $\C=\sB(\C)$ may be dilated in various ways; see Sections \ref{EXwnsSEC} and \ref{EXpropsupSEC} and Example \ref{Bex}. In Example \ref{dwnsex}, we have dilated them by considering them as elementary semigroups $T_n=\bar{c}_n\bullet c_n$, constructing coisometric dilations $w$ of $c=c_1$. In particular, the minimal coisometric dilation in \eqref{classmin}, in compact form is $w\in\sB\rtMatrix{\C\\\ell^2(\N)}$ given by the matrix
\beqn{
w
~:=~
\SMatrix{c&\delta\\&\ve},
}\eeqn
where $\delta=\sqrt{1-\abs{c}^2}e_1^*$ and $\ve$ is the adjoint of the isometric right shift on $\ell^2(\N)$, and the projection such that $pwp=wp=pc$ is $p=e_0e_0^*=\rtMatrix{1&0\\0&0}$. For the corresponding dilation by the elementary \nbd{E}semigroup $\theta$ defined by $\theta_n=w_n^*\bullet w_n$, it is easy to show by induction that the \nbd{*}algebra generated by $\theta_0(p),\ldots,\theta_n(p)$ is contained in $M_{n+1}$ (acting on $e_0,\ldots,e_n$), so that $\cA_\infty\subset\sK\rtMatrix{\C\\\ell^2(\N)}$. On the other hand, by induction, $\ls\CB{\theta_0(p)e_0,\ldots,\theta_n(p)e_0}=\rtMatrix{\C\\\C^n}$. So, $\cA_\infty$ contains all rank-one operators and, therefore, coincides with $\sK\rtMatrix{\C\\\ell^2(\N)}$. Note that this dilation is fully minimal both considered as \nbd{C^*}case or considered as von Neumann case.

Now if we have $d\ge2$ such semigroups $T^k$ and their minimal dilations $(\sB\rtMatrix{\C\\\ell^2(\N)},\theta^k,p^k)$, then we obtain a dilation of the discrete \nbd{d}parameter semigroup $\bfam{T^1_{n_1}\circ\ldots\circ T^d_{n_d}}_{\bn\in\N_0^d}$ by taking the tensor product $(\sB\rtMatrix{\rtMatrix{\C\\\ell^2(\N)}^{\text{\tiny$\otimes d$}}},\theta^1\otimes\ldots\otimes\theta^d,p^1\otimes\ldots\otimes p^d)$. Also this dilation is fully minimal both in the \nbd{C^*}case and in the von Neumann case.

However, if all $T^k$ are the same $T$ (with the dilation  from before), so that the \nbd{d}parameter semigroup is just $T_\bn:=T_{n_1+\ldots+n_d}$, obviously, $\theta_\bn:=\theta_{n_1+\ldots+n_d}$ also defines a fully minimal dilation. Clearly, this dilation cannot be the same as the dilation obtained via tensor product, because the former satisfies $\theta_{\be_i}=\theta_{\be_j}$ for all $i\ne j$, while the latter does not.
\eex

\lf
\subsection{\normalsize Minimality under unitalization} \label{uniminSSEC}

Unitalization of a dilation $(\cA,\theta,p)$ (see Section \ref{monoSEC}, after Proposition \ref{strequivprop}) is limited to strong dilations (Theorem \ref{wstrunithm}). Several of the results in the preceding subsections are limited to Markov semigroups. We now examine to what extent unitalization can help to promote results from Markov semigroups to CP-semigroups. That is, we look how conditions of minimality behave under unitalization.

For the balance, $(\cA,\theta,p)$ is a \bf{strong dilation} (of a CP-semigroup $T$) so that (by Theorem \ref{wstrunithm}) its unitalization $(\wt{\cA},\wt{\theta},\wt{p})$ is a strong \nbd{E_0}dilation (of the unitalization of $T$, $\wt{T}$). Recall that $\wt{p}:=p+(\wt{\U}-\U)$ and that $\wt{p}\wt{\cA}\wt{p}=(\wt{\U}-\U)\C+\cB$ is really isomorphic to $\C\oplus\cB\cong\wt{\cB}$ and will be used as $\wt{\cB}:=\wt{p}\wt{\cA}\wt{p}$. It is just important to be aware that the unit of $\wt{\cB}\subset\wt{\cA}$ is not the unit of $\wt{\cA}$, $\wt{\U}$, but $\U_{\wt{\cB}}=\wt{p}$, the same way as $\U_\cB=p$.

We find $\wt{E}:=\wt{\cA}\wt{p}=(\wt{\U}-\U)\C+E\cong(\wt{\U}-\U)\C\oplus E$, the Hilbert \nbd{\wt{\cB}}module direct sum of the ideal $(\wt{\U}-\U)\C$ in $\wt{\cB}$ (considered as Hilbert \nbd{\wt{\cB}}module) and the Hilbert \nbd{\cB}module $E$ (considered as Hilbert \nbd{\wt{\cB}}module).

\bcor \label{fufcor}
$(\cA,p)$ is full if and only if $(\wt{\cA},\wt{p})$ is full. Clearly, the strict topology of $\sB^a(E)$ coincides with the relative strict topology of $\sB^a(E)\subset\sB^a(\wt{E})$.
\ecor

% \proof
% $\sB^a(\wt{E})=\sB^a\rtMatrix{(\wt{\U}-\U)\C\\E}=\rtMatrix{(\wt{\U}-\U)\C\\\sB^a(E)}$ and $\sF(\wt{E})=\rtMatrix{(\wt{\U}-\U)\C\\\sF(E)}$, so, $M(\sF(\wt{E}))=M\rtMatrix{(\wt{\U}-\U)\C\\\sF(E)}=\rtMatrix{M((\wt{\U}-\U)\C)\\M(\sF(E))}=\rtMatrix{(\wt{\U}-\U)\C\\\sB^a(E)}$.\qed

\proof
$\sB^a(\wt{E})=\sB^a\rtMatrix{(\wt{\U}-\U)\C\\E}=\rtMatrix{(\wt{\U}-\U)\C\\\sB^a(E)}=\wt{\sB^a(E)}$, so $\cA=\sB^a(E)$ if and only $\wt{\cA}=\sB^a(\wt{E})$.\qed

\lf
More delicate is the discussion of algebraic minimality, full or not. We recall from the discussion after Corollary \ref{istrMricor} that $(\wt{\cA})_\infty$ is generated by $\wt{\U}-\U$, $\cA_\infty$ and all $\theta_t(\U)-\theta_s(\U)\in\cA=\sB^a(E)$. So, if $\cA_\infty$ generates $\cA$ in a topology and that topology is preserved under the passage from $\cA$ to $\wt{\cA}$, then $(\wt{\cA})_\infty$ generates $\wt{\cA}=\sB^a(\wt{E})$ in that topology. 

% This is for both the \nbd{C^*}case and the von Neumann case. The second half of the proof actually shows that, in the \nbd{C^*}case, strict algebraic minimality (recall that this is, by definition, about full dilations) is preserved.

\bcor \label{uniamincor}
~If $(\cA,\theta,p)$ is strictly (normally) algebraically minimal, then $(\wt{\cA},\wt{\theta},\wt{p})$ is strictly (normally) algebraically minimal.

% $(\cA,\theta,p)$ is strictly algebraically minimal (hence, full) if and only if $(\wt{\cA},\wt{\theta},\wt{p})$ is strictly algebraically minimal (hence, full).
\ecor

% \proof
% By the proof of Lemma \ref{striminlem}, $\cA_\infty$ is strictly dense in $\sB^a(E)$ (if and) only if the span of $\cA_\infty p\cA_\infty\subset\sF(E)\cap\cA_\infty$ is strictly dense. By the discussion after Corollary \ref{istrMricor}, we have that $(\wt{\cA})_\infty$ is generated by $\wt{\U}-\U$, $\cA_\infty$ and all $\theta_t(\U)-\theta_s(\U)\in\cA=\sB^a(E)$. So, if $\cA_\infty$ generates $\cA$ strictly, then $(\wt{\cA})_\infty$ does so for $\wt{\cA}=\sB^a(\wt{E})$. On the other hand, if $\cA_\infty$ does not generate $\cA$ strictly, then there is no chance that $(\wt{\cA})_\infty\wt{p}(\wt{\cA})_\infty=(\wt{\U}-\U)\C+\cA_\infty p\cA_\infty$ does so for $\wt{\cA}=\sB^a(\wt{E})=(\wt{\U}-\U)\C+\sB^a(E)$.\qed

Recall that in the strict case this includes, by definition, fullness. As mentioned after Corollary \ref{istrMricor}, the \nbd{\infty}strict topology under the passage from $\cA$ to $\wt{\cA}$ gets, usually, stronger.

For normal dilations, we get a partial inverse to Corollary \ref{uniamincor}.

\bprop \label{uniprimprop}
Suppose $(\cA,\theta,p)$ is normal and primary. If $(\wt{\cA},\wt{\theta},\wt{p})$ is algebraically minimal, then so is $(\cA,\theta,p)$. In particular, normal primary $(\cA,\theta,p)$ is fully minimal if and only if $(\wt{\cA},\wt{\theta},\wt{p})$ is fully minimal.
\eprop

\proof
If $(\cA,\theta,p)$ is primary, then $\U=\bigvee_{t\in\bS}\theta_t(p)$ is in $\ol{\cA_\infty}^s$. Since $\theta_t(\cA_\infty)\subset\cA_\infty$ and $\theta_t$ is normal, also $\theta_t(\U)\in\ol{\cA_\infty}^s$. So, if $\ol{\cA_\infty}^s\subsetneq\cA$, then $\ol{(\wt{\cA})_\infty}^s\subset\ol{\cA_\infty}^s+\wt{\U}\C\subsetneq\wt{\cA}=\cA+\wt{\U}\C$.\qed

\lf
As a preparation to apply the following lemma in the next subsection, we show:

\bprop \label{Pcsprop}
Let $(\cA,\theta,p)$ be a (strong as everywhere in this subsection, although here we do not unitalize!) normal dilation over the opposite of a totally directed monoid $\bS$. Then the compression with $P:=\bigvee_{t\in\bS}\theta_t(p)$ (to obtain the ``primarization'' $(\cA^P,\theta^P,p)$ of $(\cA,\theta,p)$ as in Theorem \ref{vNprimthm}\eqref{vNp3}) is strong.
\eprop

\proof
To show $\theta_t(\U-P)P=0$ it suffices to show $\theta_t(\U-P)\theta_s(p)=0$ for all $s\in\bS$. Since $\bS$ is totally directed we have $s=rt$ or $t=rs$ for suitable $r\in\bS$, leading to $\theta_t(\U-P)\theta_s(p)=\theta_t\bfam{(\U-P)\theta_r(p)}$ or $\theta_t(\U-P)\theta_s(p)=\theta_s\bfam{\theta_r(\U-P)p}$. Since $P\ge\theta_r(p)$, we have $(\U-P)\theta_r(p)=0$. Since $P\ge p$ and by strongness, we get $\theta_r(\U-P)\le\theta_r(\U-p)\le\U-p$, so $\theta_r(\U-P)p=0$.\qed

 \lf
As far as compression under unitalization is concerned, a projection $\wt{Q}\in\wt{\cA}$ has the form $q+Q$ where $Q\in\cA$ is a projection and where $q\in(\wt{\U}-\U)\C$ is either $0$ or $\wt{\U}-\U$. If $\wt{Q}\ge\wt{p}\ge\wt{\U}-\U$, then necessarily $q=\wt{\U}-\U$. So, we have a one-to-one correspondence between projections $\wt{\cA}\ni\wt{Q}\ge\wt{p}$ and $\cA\ni Q\ge p$ via $\wt{Q}=(\wt{\U}-\U)+Q$.

\blem \label{unicompelem}
Let $(\cA,\theta,p)$ be a strong dilation, choose a projection $Q\in\cA$, and put $\wt{Q}:=(\wt{\U}-\U)+Q\in\wt{\cA}$. Then the following conditions are equivalent:
\begin{enumerate}
\item \label{uni1}
$Q$ compresses $(\cA,\theta,p)$ strongly.

\item \label{uni2}
$\wt{Q}$ compresses $(\wt{\cA},\wt{\theta},\wt{p})$ strongly.

%%%% BO
% \item  \label{uni3}
% $\wt{Q}$ compresses $(\wt{\cA},\wt{\theta},\wt{p})$.
\end{enumerate}
\elem

\proof
Since $\wt{\U}-\wt{Q}=\U-Q$, we have
\beqn{
\wt{\theta}_t(\wt{\U}-\wt{Q})\wt{Q}
~=~
\wt{\theta}_t(\U-Q)\wt{Q}
~=~
\theta_t(\U-Q)Q.
}\eeqn
Therefore, $(\wt{\cA},\wt{\theta},\wt{Q})$ is a strong dilation (of the \it{a priori} CP-semigroup $\wt{\theta}^{\wt{Q}}$) if and only if $(\cA,\theta,Q)$ is a strong dilation (of the \it{a priori} CP-semigroup $\theta^Q$). 

If one, hence, both the projections $\wt{Q}$ and $Q$ satisfy this condition, then from $\wt{\theta}_t(\wt{Q}(\wt{\U}\lambda+a)\wt{Q})=\wt{\U}\lambda+\theta_t(Q-\U)\lambda+\theta_t(QaQ)$ we get
\beqn{
\wt{\theta}_t(\wt{Q}(\wt{\U}\lambda+a)\wt{Q})\wt{Q}
~=~
\wt{Q}\lambda+0+\theta_t(QaQ)Q,
}\eeqn
so that $\wt{Q}$ commutes with all $\wt{\theta}_t(\wt{Q}(\wt{\U}\lambda+a)\wt{Q})$ if and only if $Q$ commutes with all $\theta_t(QaQ)$. So, the CP-semigroup $\wt{\theta}^{\wt{Q}}_t$ is endomorphic if and only if the CP-semigroup $\theta^Q_t$ is.\qed

% \proof
% One easily verifies that $\wt{Q}$ commutes with all $\wt{\theta}_t(\wt{Q}(\wt{\U}\lambda+a)\wt{Q})=\wt{\U}\lambda+\theta_t(QaQ)$ if and only if $Q$ commutes with all $\theta_t(QaQ)$. So, $\wt{\theta}^{\wt{Q}}_t$ is endomorphic if and only if $\theta^Q_t$ is. Since $\wt{\U}-\wt{Q}=\U-Q$, we compute
% \beqn{
% \wt{\theta}_t(\wt{\U}-\wt{Q})\wt{Q}
% ~=~
% \wt{\theta}_t(\U-Q)\wt{Q}
% ~=~
% \theta_t(\U-Q)Q.
% }\eeqn
% So $(\cA,\theta,Q)$ is a strong dilation if and only if $(\wt{\cA},\wt{\theta},\wt{Q})$ is a strong dilation.\qed

\bcor \label{unicompecor}
$(\cA,\theta,p)$ is strongly incompressible if and only if $(\wt{\cA},\wt{\theta},\wt{p})$ is strongly incompressible.
%%%% BO 
% if and only if $(\wt{\cA},\wt{\theta},\wt{p})$ is incompressible.
%%%% EO
\ecor

\bob \label{unicompeob}
Under either condition in the lemma, $(\wt{\cA}^{\wt{Q}},\wt{\theta}^{\wt{Q}},\wt{p})$ is the unitalization of $(\cA^Q,\theta^Q,p)$. (Indeed, $\wt{Q}\wt{\cA}\wt{Q}=\wt{Q}((\wt{\U}-\U)\C+\cA)\wt{Q}=(\wt{\U}-\U)\C+\cA^Q=((\wt{\U}-\U)+Q)\C+\cA^Q=\wt{Q}\C+\cA^Q$ is, really, the unitalization of $\cA^Q$, once we recognize $\wt{Q}$ as the unit of $\wt{\cA^Q}$. And since the strong compression of an \nbd{E_0}dilation is an \nbd{E_0}dilation, we also recover $\wt{\theta}^{\wt{Q}}_t$ as the unitalization of $\theta^Q_t$.) 
\eob

Recovering $(\wt{\cA}^{\wt{Q}},\wt{\theta}^{\wt{Q}},\wt{p})$ as the unitalization of $(\cA^Q,\theta^Q,p)$, reconfirms, as we already know from Corollary \ref{strcomprcor}, that $(\cA^Q,\theta^Q,p)$ is strong.

\bob
The last dilation considered in Example \ref{E0weakex}, acting on $\sB\rtMatrix{\C\\H}$, does (co)re\-strict to a dilation acting on $\wt{\sB(H)}\cong\rtMatrix{\C\\\sB(H)}\subset\sB\rtMatrix{\C\\H}$, being the unitalization of  a dilation of the zero-semigroup (acting on $\sB(\zero)$) to $\sB(H)$. Clearly, $\wt{Q}:=p+q$ and $Q=p$ compress these dilations, and both compressions are not strong.
\eob

\lf
\subsection{\normalsize The one-parameter case} \label{1-p-SSEC}

In the literature, frequently one reads, in the one-parameter case, something like ``let $\theta$ be the unique minimal dilation of the CP-semigroup $T$''. However, already in this one-parameter case, the situation is actually a bit more involved than how these formulations make it appear. Let us look at what sort of bad behaviour
%%%% BO new
% the illustrated examples exhibited have so far:
the examples exhibited so far have illustrated: \index{CP-semigroup!one-parameter}\index{one-parameter case!dilation}
%%%% EO
\begin{itemize}
\item
In Subsection \ref{EXBexSEC}\ref{BexSSEC}, we explored Bhat's Example\index{Bhat's example} \ref{Bex} of a one-parameter dilation (even of a one-dimensional CP-semigroup on $\cB=\C$) which is full and incompressible, but not good (hence, not strong and necessarily non-Markov) and not algebraically minimal.  Its algebraic minimalization is not full and continues being not good. (We do not know whether the latter is compressible or not.)

\item
Every nontrivial algebraically minimal dilation of the Markov semigroup $T_t=\id_\cB$ on $\cB=\C$ (acting on the commutative algebra generated by $\theta_t(p)$) is necessarily nonfull. Examples arise by algebraic minimalization of the last paragraph of Example \ref{E0weakex} or the Examples \ref{hypexex} and \ref{discex}. But these dilations are compressible (to the trivial dilation of $T$ by itself, as explained in the beginning of Subsection \ref{minSEC}\ref{compSSEC}).
\end{itemize}
Effectively, the definitions of minimality the cited formulation ``unique minimal dilation'' is referring to, all include directly or in an equivalent way, that the ``\it{minimal} dilation'' is full and algebraically minimal -- and we just explained that even in the Markov case none of the two implies the other. (See also Footnote \ref{ArvminFN}.) In the situation where all semigroups act on some $\sB(H)$ only, fullness is automatic; but this excludes the powerful tool of unitalization, because $\wt{\sB(H)}$ is not another $\sB(H)$. In the case of general von Neumann algebras there is not really a reason to restrict the algebras $\cA$ on which the dilations act to the same Morita equivalence class as $\cB$. But is seems that the habit of expecting exactly that (namely, requiring fullness) has survived to the more general situation. Also practically all constructive results that exist so far in all dilation theory -- and the present notes are not an exception to that -- are full dilations or arise by restricting full dilations to not necessarily full ones.

\lf
In this subsection, we address, in the one-parameter (or, better, in the totally directed) von Neumann case, normal strong dilations that are fully minimal. We start with the Markov case, and recover uniqueness (Corollary \ref{1pamcMcor}) as well as Arveson's result (see again Footnote \ref{ArvminFN}) that fully minimal and strongly incompressible are equivalent, adding the new result that this implies incompressibility. Every Markov dilation can be compressed to the unique fully minimal one. Then, by applying our whole arsenal of results about properties under unitalization, we prove, as a new result, all these statements also in the strong non-Markov case (Theorems \ref{1pamcnMcor} and \ref{1pnMthm}).

Let us start with the Markov case. (Note that Parts \ref{1pM1} and \ref{1pM2} of the following Lemma also apply to the strict case, but recall that strict includes full in the definition and that fullness is the property we wish to leave apart for a moment.)

\blem \label{1pMlem}
Let $(\cA,\theta,p)$ be a normal dilation of a Markov semigroup over the opposite of a totally directed cancellative monoid $\bS$. Then:
\begin{enumerate}
\item \label{1pM1}
The left semidilation (Theorem \ref{E-lsemdilthm}) of the algebraic minimalization $(\ol{\cA_\infty}^s,\theta^\infty,p)$ 
%%%% BO new
% to $\sE^\infty$
to $\sE^\infty:=\ol{\cA_\infty}^s p$ 
%%%% EO
is a left dilation, its superproduct system ${\sE^\infty}^\podot$, therefore (Proposition \ref{lsdilPSprop}), a product system.

\item \label{1pM2}
The left dilation coincides with (is conjugate to) the left dilation constructed from the pair $({\sE^\infty}^\podot,\xi^\odot)$ leading to the module dilation $(\sE^\infty,\vt,\xi)$ in Theorem \ref{Oreindthm}.

\item \label{1pM3}
By Observation \ref{lsubdilob}, the projection $P^\infty=\id_{\sB^a(\sE^\infty)}\in\cA$ compresses $(\cA,\theta,p)$ (strongly) to the (full primary strict  \nbd{E_0})dilation $(\sE^\infty,\vt,\xi)$. Moreover, since $\sB^a(\sE^\infty)=\ol{\sF(\sE^\infty)}^s\subset\ol{\cA_\infty}^s$, the projection $P^\infty$ compresses also $(\ol{\cA_\infty}^s,\theta^\infty,p)$ (strongly) to $(\sE^\infty,\vt,\xi)$, which, therefore (Corollary \ref{compamincor}), is even fully minimal.

\item \label{1pM4}
The pair $({\sE^\infty}^\podot,\xi^\odot)$ is the unique minimal one, granted in Theorem \ref{indlimthm}, that is generated (even spanned) by its unit. The module (or full) primary, hence, \nbd{E_0}dilation $(\sE^\infty,\vt,\xi)$ is, therefore, unique and (Theorem \ref{fpincthm}) incompressible.
\end{enumerate}
\elem

\proof
\eqref{1pM1}~
Note that every generic element
\beq{ \label{sEtot}
\theta_{t_1}(pa_1p)\ldots\theta_{t_n}(pa_np)p
}\eeq
of $\sE^\infty:=\ol{\cA_\infty}^sp$ can be written in the form
\beq{ \label{dectot}
\theta_{s_m}(b_mp)\ldots\theta_{s_1}(b_1p)b_0
}\eeq
with $s_m\ge\ldots\ge s_1$ and $b_i\in\cB=p\cA p$. (This is the 
%%%% BO new 
%contents 
content 
%%%% EO 
of Bhat's reduction procedure\index{Bhat's reduction procedure}: If $t=s's\ge s\ge r$, then
\bmun{
~~~~~~
\theta_s(b)\theta_t(b')\theta_r(b'')
~=~
\theta_s(b)\theta_s(p)\theta_t(b')\theta_s(p)\theta_r(b'')
\\
~=~
\theta_s(b)\theta_s ( p\theta_{s'}(b')p ) \theta_r(b'')
~=~
\theta_s (bT_{s'}(b')) \theta_r(b''),
~~~~~~
}\emun
where in the step $\theta_r(b'')=\theta_s(p)\theta_r(b'')$ we used that $p$ is increasing. This, plus its adjoint to deal also with the case $t\ge r\ge s$, allows to bring, by induction from the right to the left, \eqref{sEtot} into the form \eqref{dectot}.) Note, too, that we get the analogue statement for $\sE^\infty_t:=\theta_t(p)\sE^\infty$ by applying $\theta_t(p)$ to \eqref{sEtot} and observing that, then, we may arrange \eqref{dectot} such that $s_m=t$. We know that the map $v_t\colon\sE^\infty\sodots\sE^\infty_t\rightarrow\sE^\infty$ defined by $x\odot y_t\mapsto\theta_t(x)y_t$ is isometric. By inserting in \eqref{dectot} (back to the generic element of $\sE^\infty$) $\theta_t(p)$ at that place that corresponds to the position where $t$ can be inserted in the chain $s_m\ge\ldots\ge s_1$ (possibly on the left when $t\ge s_m$ or on the right when $t\le s_1$; see the proof of Theorem \ref{totdirthm}), we see that $v_t$ is also surjective. So, the left semidilation given by $v_t$ is a left dilation.% (Obviously, instead of referring to Proposition \ref{lsdilPSprop} for recognizing ${\sE^\infty}^\podot$ as a product system, this can be inferred directly by restricting $v_t$ to $x_s\odot y_t\in\sE^\infty_s\sodots\sE^\infty_t\subset\sE^\infty\sodots\sE^\infty_t$.)

\eqref{1pM2}~
Clearly, the $\sE^\infty_t\subset\sE^\infty$ increase to $\sE^\infty$ and their mutual inner products are the same as those of the $k_t\sE^\infty_t$ in the inductive limit leading to Theorem \ref{Oreindthm}. Also the left dilations and the unit vectors $p=\xi$ are intertwined correctly.

\eqref{1pM3}~
There is nothing to add here.
% Let $\cA\ni a=P^\infty aP^\infty$ be an arbitrary element of $\sB^a(\sE^\infty) \subset\ol{\cA_\infty}^s$. Then, for the generic element $v_t(x\odot y_t)$ of $\sE^\infty$, we find
% \beqn{
% P^\infty\theta_t(a)P^\infty v_t(x\odot y_t)
% ~=~
% P^\infty\theta_t(a) \theta_t(x) y_t
% ~=~
% P^\infty\theta_t(ax) y_t
% ~=~
% P^\infty v_t(ax\odot y_t)
% ~=~
% \vt_t(a)v_t(x\odot y_t).
% }\eeqn
% It follows $\theta^{P^\infty}=\vt$, that is, $P^\infty$ compresses $(\cA,\theta,p)$ to a primary full \nbd{E_0}dilation. The compression to an \nbd{E_0}dilation is strong. Now, since $\sB^a(\sE^\infty) \subset\ol{\cA_\infty}^s$, the module dilation $(\sE^\infty,\vt,\xi)$ is a compression of $(\ol{\cA_\infty}^s,\theta^\infty,p)$, and by Corollary \ref{compamincor} it is algebraically minimal. 

\eqref{1pM4}~
This follows, because the typical elements of $\sE^\infty_t$ in \eqref{dectot} (recall that in $\sE^\infty_t$ means we choose $s_m=t$) transforms, if we choose $r_i$ such that $s_1=r_1$ and $s_{i+1}=r_{i+1}s_i$, into $b_m\xi_{r_m}\ldots b_1\xi_{r_1}b_0$.\qed

\bcor \label{1pamcMcor}
For a normal Markov dilation $(\cA,\theta,p)$ over the opposite of a totally directed cancellative monoid $\bS$, the following conditions are equivalent:
\begin{enumerate}
\item \label{amc1}
The dilation is fully minimal.

\item \label{amc2}
The dilation is incompressible.

\item \label{amc3}
The dilation is strongly incompressible.
\end{enumerate}
Moreover, any such dilation is determined uniquely by the Markov semigroup it dilates.
\ecor

\proof
By Corollary \ref{amfinccor}, \eqref{amc1}$\Rightarrow$\eqref{amc2}. Of course, \eqref{amc2}$\Rightarrow$\eqref{amc3}. And if the dilation is strongly incompressible, then, by the lemma, it has to be the unique one characterized there, hence, \eqref{amc3}$\Rightarrow$\eqref{amc1}.\qed
%since, by Corollary \ref{amE0ccor}, all compressions are strong, we have \eqref{amc3}$\Rightarrow$\eqref{amc2}, as well.\qed

\lf
We call this dilation the \phantomsection\hl{minimal}\index{one-parameter case!the minimal dilation}\index{dilation, weak!one-parameter!the minimal dilation}\index{minimal!dilation!in the one-parameter case} dilation, and reformulate:

\bthm \label{1pMthm}
Let $T$ be a normal Markov semigroup over the opposite of a totally directed cancellative monoid $\bS$. Then:
\begin{enumerate}
\item \label{1pMt1}
$T$ admits a normal Markov dilation, namely, the unique minimal one.

\item \label{1pMt2}
Every normal dilation of $T$ `contains'  the unique minimal dilation, in the sense that the former may always be compressed to the latter.
\end{enumerate}
\ethm

\lf\noindent
We now come to the case of CP-semigroups that are 
%%%% BO new
%non
not 
%%%% EO 
necessarily Markov. The strategy is to apply what the preceding theorem asserts to their unitalizations. But, this is limited to strong dilations. % -- and we know that there is no hope to get something this way, if the dilation is not at least a good one. 
Regarding good dilations, we know from the discussion leading to Theorem \ref{strgothm} that, passing to the algebraic minimalization of a good dilation, we get a strong one.

%Corollary \ref{1pamcMcor} remains true for strong dilations. More precisely:

\bthm \label{1pamcnMcor}
For a normal strong dilation $(\cA,\theta,p)$ over the opposite of a totally directed cancellative monoid $\bS$, the following conditions are equivalent:
\begin{enumerate}
\item \label{amcn1}
The dilation is fully minimal.

\item \label{amcn2}
The dilation is incompressible.

\item \label{amcn3}
The dilation is strongly incompressible.
\end{enumerate}
Moreover, any such dilation is determined uniquely by the CP-semigroup it dilates.
\ethm

\proof
We have \eqref{amcn1}$\Rightarrow$\eqref{amcn2} (by Corollary \ref{amfinccor}) and, obviously, \eqref{amcn2}$\Rightarrow$\eqref{amcn3}. Since $(\cA,\theta,p)$ is strong we may pass to the unitalization $(\wt{\cA},\wt{\theta},\wt{p})$. If $(\cA,\theta,p)$ is strongly incompressible, then, by Corollary \ref{unicompecor}, so is $(\wt{\cA},\wt{\theta},\wt{p})$. Therefore, by Corollary \ref{1pamcMcor}, $(\wt{\cA},\wt{\theta},\wt{p})$ is fully minimal, so, by Propositions \ref{uniprimprop} and \ref{Pcsprop}, $(\cA,\theta,p)$ is fully minimal, too. Thus, \eqref{amcn3}$\Rightarrow$\eqref{amcn1}. Since $(\wt{\cA},\wt{\theta},\wt{p})$ determines $(\cA,\theta,p)$ (see Observation \ref{preurecob}), also $(\cA,\theta,p)$ is determined uniquely by $T$.\qed

\lf
We call also this unique dilation the \phantomsection\hl{minimal}\index{one-parameter case!the minimal dilation}\index{dilation, weak!one-parameter!the minimal dilation}\index{minimal!dilation!in the one-parameter case} dilation. Here is the general analogue of Theorem \ref{1pMthm}:

\bthm \label{1pnMthm}
Let $T$ be a normal CP semigroup over the opposite of a totally directed cancellative monoid $\bS$. Then:
\begin{enumerate}
\item \label{1pnMt1}
$T$ admits a normal strong dilation, namely, the unique minimal one.

\item \label{1pnMt2}
Every normal strong dilation of $T$ `contains'  the unique minimal dilation, in the sense that the former may always be compressed to the latter.

\item \label{1pnMt3}
Every normal good dilation of $T$ `contains'  the unique minimal dilation, in the sense that the the algebraic minimalization of the former may always be compressed to the latter.
\end{enumerate}
\ethm

\lf\lf\lf

\newpage

\section[\sc{Examples:} The case $\N_0^2$; another ``multi-example'']{Examples: The case $\N_0^2$; another ``multi-example''}\label{EXN02SEC}

In this section, we examine a way to construct dilations of Markov semigroups over $\N_0^2$ on a von Neumann algebra. \index{two-parameter case!discrete}\index{discrete case!two-parameter}\index{dilation, weak!discrete  two-parameter case}The dilation is constructed by embedding the GNS-subproduct system into a product system with the aid of Corollary \ref{N02cor} and Observation \ref{N0dob} and, then, appealing to the inductive limit construction in Theorem \ref{Oreindthm}.  ($\N_0^2$ is an Ore monoid; and it is abelian so that the tedious reference to the opposite of an Ore monoid disappears.) The dilation is, therefore, full and primary (hence, an \nbd{E_0}dilation). It shares (sometimes under some extra conditions, but still for large classes of Markov semigroups) much of the bad behaviour we know from Bhat's Example \ref{EXBexSEC}\ref{BexSSEC}. The more important it is to emphasize that, here, we are speaking about a dilation of a Markov semigroup which, therefore, is strong, hence, good. (The main purpose of the discussion of Bhat's example\index{Bhat's example} was to find a dilation that is not good.)

Recall that a semigroup $T=\bfam{T_\bn}_{\bn\in\N_0^2}$ (no matter whether Markov, endomorphism, or just CP) is determined by the two commuting maps $T_i:=T_{\be_i}$ so that the marginal semigroups are given as $T^i_n=T_i^n$. Conversely, given two commuting maps $T_i$, by $T_\bn:=T_2^{n_2}\circ T_1^{n_1}$ we define a semigroup over $\N_0^2$.

Starting from the GNS-constructions $(\sE_i,\xi_i)$ of the commuting Markov maps $T_i$, we will construct $E_i\supset\sE_i\ni\xi_i$ and $\sF_{1,2}\colon E_2\odot E_1\rightarrow E_1\odot E_2$ satisfying $\sF_{1,2}(\xi_2\odot\xi_1)=\xi_1\odot\xi_2$ as ingredients for Corollary \ref{N02cor}. This, in turn, as ingredient to Theorem \ref{Oreindthm}, gives a module dilation. Actually, we will have $E_i=E$ and $\sF\colon E\odot E\rightarrow E\odot E$ playing the role of $\sF_{1,2}$. (This situation resembles most the situation of Theorem \ref{truncFijthm}; the more important it is to take into consideration Observation \ref{ntruncob} for avoiding a possible confusion. We leave to the reader to evaluate whether or not our decision to avoid many indices indicating if an $E $ in question works as $E_1$ or $E_2$, is helpful.) The construction we propose, requires, at a certain point, that a certain submodule possesses a complement; it, therefore, is limited to the von Neumann case. So, unless stated otherwise explicitly, $T_i$ are commuting normal Markov maps on a von Neumann algebra $\cB\subset\sB(G)$, and unless stated otherwise explicitly, the modules and correspondences (including GNS-$T_i$) are von Neumann (including tensor products of the latter). 

%%%% BO
It is fair to acknowledge that a first version of our proof of Theorem \ref{N02Mdilthm} was inspired by Solel's proof of \cite[Lemma 5.10]{Sol06}. 
%%%% EO

The trick is passing to the direct sum $E:=\sE_1\oplus\sE_2$. (Looking at how importantly this idea enters also crucial steps in the proof of Corollary \ref{N02cor}, such as Proposition \ref{E_nprop}, we are tempted to call this our \it{standard trick}.) Clearly, $E$ contains both $\xi_1$ and $\xi_2$, so the tensor product $E\odot E$, which we like to visualize as
\beqn{
E\odot E
~=~
\bigoplus_{i,j}\sE_i\odot\sE_j
~=~
\SMatrix{\sE_1\odot\sE_1&\sE_1\odot\sE_2\\\sE_2\odot\sE_1&\sE_2\odot\sE_2},
}\eeqn
contains both $\xi_1\odot\xi_2$ and $\xi_2\odot\xi_1$; but they live in distinct direct summands. However, since $T_1$ and $T_2$ commute, the map
\beqn{
\xi_2\odot\xi_1
~\longmapsto~
\xi_1\odot\xi_2
}\eeqn
extends to a unique isomorphism $\f\colon F_{21}\rightarrow F_{12}$ from the von Neumann subcorrespondence $F_{21}:=\cls^s\cB\xi_2\odot\xi_1\cB$ of $E\odot E$ generated by $\xi_2\odot\xi_1$ (isomorphic to the GNS-correspondence of $T_1\circ T_2$) to the von Neumann subcorrespondence $F_{12}:=\cls^s\cB\xi_1\odot\xi_2\cB$ of $E\odot E$ generated by $\xi_1\odot\xi_2$ (isomorphic to the GNS-correspondence of $T_2\circ T_1=T_1\circ T_2$). Since $F_{21}$ and $F_{12}$ are orthogonal and complemented in $E\odot E$, also $F=F_{21}+F_{12}$ is  complemented. We, therefore, may extend the automorphism $\f+\f^{-1}$ of $F$ by the identity on $F^\perp$ to an automorphism $\sF$ of $E\odot E$.

\brem
$\sF$ acts as identity on $\sE_i\odot\sE_i$. It, clearly (co)restricts to an automorphism of $(\sE_2\odot\sE_1)\oplus(\sE_1\odot\sE_2)$. It does, however, not (co)restrict to an isomorphism $\sE_2\odot\sE_1\rightarrow\sE_1\odot\sE_2$ (or \it{vice versa}). ($\sF$ only exchanges the parts that are generated by $\xi_2\odot\xi_1$ and $\xi_1\odot\xi_2$; the parts orthogonal to them are untouched.) In fact, we know from Examples \ref{SCex3} that $\sE_2\odot\sE_1$ and $\sE_1\odot\sE_2$ need not be isomorphic.
\erem

As announced, we input $E_1:=E\ni\xi_1$, $E_2:=E\ni\xi_2$, and $\sF_{1,2}:=\sF$ into Corollary \ref{N02cor} and obtain:
\begin{itemize}
\item
A unique product system structure on the family $E^\odot=\bfam{E^{\odot(n_1+n_2)}}_{\bn\in\N_0^2}=\bfam{E^{\odot n_1}\odot E^{\odot n_2}}_{\bn\in\N_0^2}$, satisfying $u_{\be_1,\be_2}^*u_{\be_2,\be_1}=\sF$.

\item
A unique unit $\xi^\odot=\bfam{\xi_1^{\odot n_1}\odot\xi_2^{\odot n_2}}_{\bn\in\N_0^2}$ satisfying $\xi_{\be_i}=\xi_i$.
\end{itemize}
By Theorem \ref{Oreindthm}, we obtain:

\bthm \label{N02Mdilthm}
Every (normal) Markov semigroup over $\N_0^2$ on a von Neumann algebra admits a 
%%%% BO
(normal)
%%%% EO 
strong primary module dilation $(\sE,\vt,\xi)$.
\ethm

\bob
We may ask to what extent this result may be generalized to other semigroups over 
%%%% BO 
% $\N_0^d$.
$\N_0^2$.
%%%% EO
\begin{itemize}

\item
The word ``normal'' may be left out, as long as $\cB$ is a von Neumann algebra. (This the next in a row of statements of that type in Bhat and Skeide \cite{BhSk00} and Barreto, Bhat, Liebscher and Skeide \cite{BBLS04}, which assert similar things.) In fact, the only thing we need is that we may pass to the strong closure of a right module $E$ over $\cB\subset\sB(G)$ in the canonically associated $\sB(G,E\odot G)$, maintaining the property that the obtained closure is still a Hilbert module over $\cB$. A possible left action survives that. Therefore, also all \nbd{C^*}tensor products may be closed strongly. So each single step in the construction of the dilation goes through also if the CP-maps are not normal.

\item
We do not know if the result can be lifted to Markov semigroups on \nbd{C^*}algebras. By passing to the bidual, we get a normal CP-semigroup $T^{**}$ on $\cB^{**}$, and a module dilation $(\sE,\vt,\xi)$ of the latter. What we do not know, is if $\sB^a(E)$ contains a subalgebra $\cA$ that is left invariant by $\vt$ such that $\xi\xi^*\cA\xi\xi^*=\xi\cB\xi^*$. (Note that $\xi\xi^*\sB^a(E)\xi\xi^*=\xi\cB^{**}\xi^*$.) In the one-parameter case, this problem can be shipped around; but in this case we do not need $\cB^{**}$, because we have a direct construction of the product system in the \nbd{C^*}case (which can easily be pushed forward to prove the von Neumann case). However, proving the one-parameter \nbd{C^*}case by reducing it (unnecessarily) to the von Neumann case via passage to the bidual semigroup, illustrates what goes wrong in the two-parameter case; we recommend this as an instructive exercise for everybody who wishes to prove the result for the two-parameter \nbd{C^*}case.
\end{itemize}
\eob

\brem
It is noteworthy that the construction of $E$ and $\sF$ can be generalized. We can take any $E$ containing suitable $\xi_1$ and $\xi_2$ as long as $\AB{\xi_1,\cB\xi_2}=\zero$. And in the definition of $\sF$ we can replace the identity on $F^\perp$ just with any automorphism. It would be interesting to compare dilations obtained this way with our ``most economic'' choice.
\erem

%%%% BO
By unitalization, Theorem \ref{N02Mdilthm} generalizes to (normal) CP-semigroups -- but only to some extent:

\bcor \label{N02dilcor}
Every (normal) CP-semigroup over $\N_0^2$ on a von Neumann algebra admits a (normal) strong full dilation.
\ecor

\proof
If $T$ is the CP-semigroup, by Theorem \ref{N02Mdilthm}, we get a full \nbd{E_0}dilation $(\wh{\cA},\wh{\theta},\wh{p})$ of $\wt{T}$ and, further, by Theorem \ref{uninonunithm}, a dilation of $T$. Since the left action of $q$ on $\wt{\sE}_i$ is right multiplication by $\wt{\U}_{\wt{\cB}}-\U_\cB\in\wt{\cB}$, we see that $q$ is central in $\wh{\cA}$. By Corollary \ref{qcentcor}, this dilation is full.\qed

\brem
Existence of a strong dilation was proved in Solel \cite[Theorem 5.13]{Sol06}, not addressing the question of fullness, nor the question of whether or not the constructed dilation is an \nbd{E_0}dilation when $T$ is Markov (we believe that this is not necessarily so). We do not know if the dilation is primary. But, of course, by Theorem \ref{vNprimthm}\eqref{vNp3}, a normal dilation may be compressed to primary dilation and, by Observation \ref{compob}\eqref{cob2}, the compression is full, again.
\erem

% \brem
% The dilation of a CP-semigroup obtained as in the Corollary need not be primary. Indeed, consider the dilation of the zero semigroup on $\C$; the dilation constructed by the above procedure is the semigroup itself, and it is not primary. 
% \erem
%%%% EO

\subsection{\normalsize The spanned superproduct system}

We have the product system $E^\odot$ and the unit $\xi^\odot$ giving back $T$ as $T_\bn=\AB{\xi_\bn,\bullet\xi_\bn}$. By Example \ref{unitgenex}, the subproduct system $S^\bodot$ of $T$ sits in $E^\odot$ as $S_\bn=\cls^s\cB\xi_\bn\cB$. We are interested in the superproduct subsystem $\sE^\podot$ of $E^\odot$ spanned by $S^\bodot$ in the sense of Theorem \ref{supintthm}. In particular, we shall show that under a suitable nontriviality condition on the Markov maps $T_i$, it is not a product system.

\bdefi
A pair of (normal) CP-maps $(T_1,T_2)$ is \hl{quasi-generic}\index{quasi-generic!CP-maps}\index{CP-map!pair!quasi-generic} if the (von Neumann) GNS-correspondences satisfy $F_{21}\subsetneq\sE_2\odot\sE_1$ or $F_{12}\subsetneq\sE_1\odot\sE_2$ under the canonical embeddings. The pair is \hl{generic}\index{generic!CP-maps}\index{CP-map!pair!generic} if both conditions are fulfilled.
\edefi

Most pairs of CP-maps are generic; this is why we call them that. For example, suppose that $(T_1,T_2)$ are CP maps on $\sB(G)$ with $\dim G < \infty$, and suppose that they have minimal Kraus decompositions $T_i = \sum_{k=1}^{n_i} {c_k^i}^* \bullet c_k^i$ ($i=1,2$), according to Appendix \ref{vNAPP}\ref{vNBGmod}. Then whenever the set $\bfam{c_k^1 c_\ell^2}_{k,\ell}$ is linearly dependent, we have a strict inclusion $F_{12}\subsetneq\sE_1\odot\sE_2$, because the former is the GNS representation of $T_1 \circ T_2$ while the latter has too big a dimension. An easy sufficient condition for the pair $(T_1, T_2)$ to be generic, is $n_1 n_2 > (\dim G)^2$. 
%On the other hand, we see that if one of the maps is elementary, then the pair is not quasi-generic. Likewise
CP maps on $\C^n$ are given by matrices with positive entries (see Example \ref{SCex3}), and as soon as in of both their products there occurs a matrix entry in the computation of which there occurs a nontrivial sum, the pair is generic. There are several more examples in Subsection \ref{strcomSEC}\ref{SCex}, which we will meet back later on in this section.

\bthm \label{nPSspthm}
If the pair $(T_1,T_2)$ is quasi-generic, then $\sE^\podot$ is not a product system.
\ethm

\proof
Clearly, $\sE_{\be_i}=S_{\be_i}\cong\sE_i$. If $\sE^\podot$ was a product system, then $\sE_{\be_1+\be_2}=\cls^s\sE_{\be_1}\sE_{\be_2}=\cls^s\sE_{\be_2}\sE_{\be_1}$. We have
\baln{
\cls^s\sE_{\be_1}\sE_{\be_2}
&
~=~
\SMatrix{~~~&\sE_1\odot\sE_2\\~&\big.},
&
\cls^s\sE_{\be_2}\sE_{\be_1}
&
~=~
\SMatrix{&F_{12}\\F_{21}^\perp\cap(\sE_2\odot\sE_1)&}.
}\ealn
The two coincide (if and) only if $(T_1,T_2)$ is not quasi-generic.\qed

\bob \label{SCncob}
So, for a quasi-generic commuting pair $(T_1,T_2)$ we have constructed a dilation whose contained GNS-subproduct system does not span a product system. On the other hand, if $T_1$ and $T_2$ commute even strongly, then we may construct another dilation following Example \ref{SCdd} from Section \ref{strcomSEC}; the product system of this dilation \bf{is} spanned by the contained GNS-subproduct system. We see that the property whether or not a subproduct system, when embedded into a product system, spans a product subsystem or only a superproduct subsystem, is not intrinsic to the subproduct system but depends on the embedding.

Clearly, the two corresponding primary strong module dilations cannot be conjugate.
\eob

\bob \label{ntruncob}
A warning about the following little pitfall is in place. We mentioned that the procedure to define a product system structure over $\N_0^2$ on $E_\bn=E^{\odot(n_1+n_2)}$, on a formal level resembles much the situation in Theorem \ref{truncFijthm}; but Theorem \ref{truncFijthm}\eqref{Fij2} asserts that a certain truncated subproduct system (looking carefully at the proof, always) spans a product system. This sounds like a potential contradiction to Theorem \ref{nPSspthm}. Note, however, the crucial condition in Theorem \ref{truncFijthm} that the embeddings $F_{\be_i+\be_j}\rightarrow F_{\be_i}\odot F_{\be_j}$ (in the notations of Theorem \ref{truncFijthm}) are unitary. In the situation of Theorem \ref{nPSspthm}, this happens if and only if the pair of Markov maps is not quasi-generic.
\eob

\lf
\subsection{\normalsize The generated product system} \label{spNpsSSEC}

One might think that the problem in Observation \ref{SCncob}, two non-conjugate dilations of the same Markov semigroup, would depend simply on that our product system is not generated by the contained GNS-subproduct system. While this might be true, it does not sort out the problem. In fact, we know from (the von Neumann version of) Theorem \ref{Ocomprthm} that the (von Neumann, hence, projectionable) product subsystems containing the unit are in one-to-one correspondence with the compressions of the given dilation. Therefore, if, by any chance, the unit (that is, the GNS-subproduct system) does not generate the whole product system, by Observation \ref{vNintob}, we may simply pass to the minimal product subsystem containing the unit, compress the dilation with the corresponding projection, and obtain a dilation that has a product system generated by the GNS-subproduct system.

We show now that if $\cB$ is a factor, we need not worry about compression: The product system is, then, generated by the contained GNS-subproduct system. We prepare this, formulating the following lemma first in the \nbd{C^*}case.

\blem
Let $E$ be an \nbd{\cA}\nbd{\cB}correspondence, let $F$ be a \nbd{\cB}\nbd{\cC}correspondence, and suppose $\cB$ is simple. Then $E\odot F=\zero$ implies $E=\zero$ or $F=\zero$.
\elem

\proof
Obvious. (If $0\ne x\in E$ and $F\ne\zero$, so that, by simplicity of $\cB$, $F$ is faithful, then there exists $y\in F$ such that $\abs{x}y\ne0$, so, $x\odot y\ne0$.)\qed

\bcor
Let $F_1,F_2$ be complemented subcorrespondences of a \nbd{\cB}correspondence $E$, and assume $\cB$ is simple. If $F_1 \odot E=E \odot F_2$ then $F_i = E$ $(i=1,2)$ or $F_i=\zero$ $(i=1,2)$. 
\ecor

\proof
We have
\baln{
E\odot E
&
~=~
(F_1\odot F_2)\oplus(F_1\odot F_2^\perp)\oplus(F_1^\perp\odot F_2)\oplus(F_1^\perp\odot F_2^\perp),
\\
F_1\odot E
&
~=~
(F_1\odot F_2)\oplus(F_1\odot F_2^\perp),
\\
E\odot F_2
&
~=~
(F_1\odot F_2)~~~~~~~~~~~~~~~~~~~~\,\oplus(F_1^\perp\odot F_2).
}\ealn
So, $F_1 \odot E=E \odot F_2~~\Longleftrightarrow~~F_1\odot F_2^\perp=\zero=F_1^\perp\odot F_2$.

If $E=\zero$, there is nothing to show.

If $E\ne\zero$, then $F_1\ne\zero$ or $F_1^\perp\ne\zero$.

If $F_1\ne\zero$, then $F_2^\perp=\zero$, so, $F_2=E\ne\zero$, hence, $F_1^\perp=\zero$, so, $F_1=E$.

If $F_1^\perp\ne\zero$, then $F_2=\zero$, so, $F_2^\perp=E\ne\zero$, hence, $F_1=\zero$.\qed

\bob
 We are interested in the case $F_1=F_2$, so, $F\odot E=E\odot F$ implies $F=E$ or $F=\zero$. Other special cases are $F_1=E$ or $F_2=E$. Obviously, the corollary generalizes to the situation $F_i\subset E_i$ where $E_1$ is an \nbd{\cA}\nbd{\cB}correspondence and $E_2$ is a \nbd{\cB}\nbd{\cC}correspondence. The example $E=\cB$ and $F$ an essential ideal in $\cB$ shows, that the condition that $F$ be complemented is indispensable.
\eob

\bob
Obviously, the lemma and its corollary remain true for von Neumann correspondences, with $\cB$ a factor. (The von Neumann tensor product $E\sodots F$ is $\zero$ if and only if the \nbd{C^*}tensor product is $\zero$, so if the former is $\zero$ we just apply the lemma to the latter. And the proof of the corollary works in every additive tensor category in which the lemma holds.)
\eob

\bthm \label{facpsgenthm}
If $\cB$ is a factor and if the pair $(T_1,T_2)$ is quasi-generic, then the smallest product subsystem of $E^\odot$ containing $S^\bodot$ is $E^\odot$.
\ethm

\proof
Suppose the family $F_\bn$ forms a product subsystem of $E^\odot$ containing $S^\bodot$. Since $F_{\be_i}$ must contain $\sE_i$, we find $F_{\be_1}=\sE_1\oplus G_2\subset\sE_1\oplus\sE_2=E$ and  $F_{\be_2}=G_1\oplus\sE_2\subset\sE_1\oplus\sE_2=E$ for suitable subcorrespondences $G_i\subset\sE_i$. Since the $F_\bn$ form a product system, we must have $\cls^s F_1F_2=\cls^s F_2F_1$. Therefore, the following two subsets of $\sE_{\be_1+\be_2}=E\odot E$,
\beqn{
F_1\odot F_2
~=~
\SMatrix{\sE_1\odot G_1&\sE_1\odot\sE_2\\G_2\odot G_1&G_2\odot\sE_2}
\text{~~~~~~and~~~~~~}
\sF(F_2\odot F_1)
~=~
\sF\SMatrix{G_1\odot\sE_1&G_1\odot G_2\\\sE_2\odot\sE_1&\sE_2\odot G_2},
}\eeqn
must coincide. $\sF$ does nothing to the diagonal entries, so we get $\sE_i\odot G_i=G_i\odot\sE_i $. Since $\cB$ is a factor, by (the von Neumann version of) the corollary of the lemma, we get, for each $i=1,2$, that either $G_i=\sE_i$($\ne\zero$, since $T_1$ is Markov) or $G_i=\zero$. Since the pair $(T_1,T_2)$ is quasi-generic, at least one of the two $F_{ij}^\perp\cap(\sE_i\odot\sE_j)$ $(i\ne j)$ is not $\zero$ and stays under the flip where it is. Therefore, the corresponding $G_i\odot G_j$ cannot be $\zero$, so $G_i\ne\zero$, hence $G_i=\sE_i$ for $i=1,2$. It follows that $F_{\be_i}=E_{\be_i}$ and, further, $F_\bn=E_\bn$.\qed

\lf
From the discussion in the beginning of the paragraph we conclude:

\bcor
The dilation $(\sE,\vt,\xi)$ of $T$ is incompressible.
\ecor

\lf
\subsection{\normalsize Algebraic minimality}

We have already Bhat's Example\index{Bhat's example} \ref{EXBexSEC}\ref{BexSSEC} for an incompressible (not good, therefore, not strong) module dilation of a discrete one-parameter semigroup on $\C$ that is not algebraically minimal. Now, we illustrate that, at least under an (admittedly, quite restrictive) extra condition, the strong primary dilation constructed to prove Theorem \ref{N02Mdilthm} is not algebraically minimal and can be incompressible.

The condition is that the pair $(T_1,T_2)$ is not generic. It still maybe quasi-generic, and in this case, by the preceding paragraphs, the dilation is also incompressible. Example \ref{SCex4} (that is, Solel \cite[Example 5.5]{Sol06}) shows that we are not working on the empty set. 

\lf
Let us start with just any normal Markov semigroup $T$ over $\N_0^2$ and construct the module dilation $(\sE,\vt,\xi)$ as for Theorem \ref{N02Mdilthm}. By Theorem \ref{fminthm}, to show that this full dilation is not algebraically minimal, we have to show that $\ol{\cA_\infty\xi}^s\ne\sE$.

Note that, for whatever fixed $\bm^1,\bm^2\in\N_0^2$ and strongly total subsets $\Sigma_\bn\subset E_\bn$, the set
\beqn{
\bigcup_{\bn\in\N_0^2}
\xi\Sigma_{\bm^1}\Sigma_\bn\Sigma_{\bm^2}
}\eeqn
is strongly total in $\sE$. (Indeed, $\Sigma_{\bm^1}\Sigma_\bn\Sigma_{\bm^2}$ is strongly total in $E_{\bm^1+\bn+\bm^2}$ and $\CB{\bm^1+\bn+\bm^2\colon\bn\in\N_0^2}$ is cofinal in $\N_0^2$.)

Recall that $E_\bn=E^{\odot(n_1+n_2)}$ so that a typical element of $E_\bn$ is $x_{n_1}\odot\ldots\odot x_1\odot y_{n_2}\odot\ldots\odot y_1$ with $x_i,y_i\in E$. Sometimes, it is important to remember that the $x_i$ are, actually, from $E_{\be_1}=E$, while the $y_i$ are, actually, from $E_{\be_2}=E$; in particular, it is possible that in the course of a computation this affiliation changes. We, therefore, will sometimes attach to an $x\in E$ a superscript $\be_i$ to indicate that $x^{\be_i}$ has to be viewed as an element of $E_{\be_i}=E$.

\bex
Suppose we have $x_2,y_2\in\sE_2\subset E$. Then $\sF(x_2\odot y_2)=x_2\odot y_2$. But, when we recall that $\sF$ defines $u_{\be_1,\be_2}u_{\be_2,\be_1}^*$ so that, actually, $x_2$ is an element of $E_{\be_2}$, while $y_2$ is an element of $E_{\be_1}$, the computation goes as
\vspace{-1ex}
\beqn{
u_{\be_1,\be_2}u_{\be_2,\be_1}^*
\colon
x_2^{\be_2}\odot y_2^{\be_1}
~=~
x_2\odot y_2
~\overset{\sF}{\longmapsto}~
x_2\odot y_2
~=~
x_2^{\be_1}\odot y_2^{\be_2}.
}\eeqn
\eex

\lf
To show that $\ol{\cA_\infty\xi}^s$ is not all of $\sE$, it is sufficient to find a nonzero submodule $F\subset\sE$, perpendicular to $\xi$ and invariant under all $\vt_\bn(\xi b\xi^*)$. (Indeed, under these hypothesis we have $\AB{F,\vt_\bn(\xi\cB\xi^*)\xi}=\AB{\vt_\bn(\xi\cB\xi^*)F,\xi}\subset\AB{F,\xi}=\zero$, so $F\ne\zero$ is perpendicular to $\cA_\infty\xi$.)

Following our notation, by $\sE_2^{\be_1}$ we indicate the direct summand $\sE_2$ of $E$ when considered as $E=E_{\be_1}$. We define $F:=\cls^s\sE\sE_2^{\be_1}\ne\zero$.

$F$ satisfies:
\begin{itemize}
\item
$\AB{\xi,F}=\zero$. (Indeed, $\AB{\xi,Xx_2^{\be_1}}=\AB{\xi\xi_{\be_1},Xx_2^{\be_1}}=\AB{\xi_{\be_1},\AB{\xi,X}x_2^{\be_1}}=0$ for $x_2^{\be_1}\in\sE_2^{\be_1}$.) In particular, $\vt_0(\xi\cB\xi^*)F=\zero\subset F$.

\item
$\vt_\bn(\xi\cB\xi^*)F\subset F$, whenever $\bn\ge\be_1$. (More generally, $\vt_{\be_1}(a)(Xx_2^{\be_1})=(aX)x_2^{\be_1}$, so $\vt_{\be_1}(a)F\subset F$ for every $a\in\sB^a(E)$. In particular, if $\bn=\bn'+\be_1\ge\be_1$, then $\vt_\bn(\xi\cB\xi^*)F=\vt_{\be_1}(\vt_{\bn'}(\xi\cB\xi^*))F\subset F$.)
\end{itemize}
So, the only question that remains is whether or not $\vt_{n\be_2}(\xi\cB\xi^*)F\subset F$ ($n\ge1$). Now, our special hypothesis that $(T_1,T_2)$ is not generic comes into the game.

\bthm \label{naminthm}
If $F_{12}=\sE_1\odot\sE_2$, then $(\sE,\vt,\xi)$ is not algebraically minimal.
\ethm

\proof
To get hold of how $\vt_{n\be_2}$ acts on $F$, we wish to write somehow elements $Xx_2^{\be_1}$ from $F$ in a form as (sum over) elements from $\sE E_{n\be_2}$, where we know how $\vt_{n\be_2}$ acts. We will actually show that $F\subset\cls^sFE_{n\be_2}$ so that $\vt_{n\be_2}(\xi\cB\xi^*)F\subset\cls^s\vt_{n\be_2}(\xi\cB\xi^*)FE_{n\be_2}=\zero\subset F$.

We prepare for a proof by induction. The set $\sE E_{\be_2}$ is strongly total in $\sE$, so the set $\sE E_{\be_2}\sE_2^{\be_1}$ is strongly total in $F$. In $Xy^{\be_2}z_2^{\be_1}$ (with $X\in\sE$, $y=y_1+y_2\in E$ ($y_i\in\sE_i\subset E$) and $z_2\in\sE_2\subset E$), we wish to write the element $y^{\be_2}z_2^{\be_1}$ of $E_{\be_2}\sE_2^{\be_1}\subset E_{\be_1+\be_2}$ as an element of $\cls^sE_{\be_1}E_{\be_2}=E_{\be_1+\be_2}$. Now, $E_{\be_1+\be_2}=E\odot E$ and $y^{\be_2}z_2^{\be_1}=u_{\be_2,\be_1}(y^{\be_2}\odot z_2^{\be_1})$ where
\beqn{
E_{\be_2}\odot E_{\be_1}
~\supset~
E_{\be_2}\odot\sE_2^{\be_1}
~\ni~
y^{\be_2}\odot z_2^{\be_1}
~=~
\SMatrix{~~~~~~&y_1\odot z_2\\&y_2\odot z_2}
~\in~
E\odot\sE_2
~\subset~
E\odot E.
}\eeqn
By definition, the product $u_{\be_2,\be_1}$ makes this
\beqn{
E_{\be_1+\be_2}
~\supset~
E_{\be_2}\sE_2^{\be_1}
~\ni~
y^{\be_2}z_2^{\be_1}
~=~
\sF\SMatrix{~~~~~~&y_1\odot z_2\\&y_2\odot z_2}
~\in~
\sF(E\odot\sE_2)
~\subset~
E\odot E.
}\eeqn
Here is where our hypothesis strikes:
\beqn{
E_{\be_1+\be_2}
~\supset~
E_{\be_2}\sE_2^{\be_1}
~\ni~
y^{\be_2}z_2^{\be_1}
~=~
\SMatrix{&\\\f^{-1}(y_1\odot z_2)~&~y_2\odot z_2}
~\in~
\sE_2\odot E
~\subset~
E\odot E.
}\eeqn
Consequently, $E_{\be_2}\sE_2^{\be_1}\subset\cls\sE_2^{\be_1}E_{\be_2}$ (because $E_{\be_1+\be_2}\ni x^{\be_1}y^{\be_2}=x\odot y$.) This allows to conclude by induction that
\bmun{
F
~=~
\cls^s\sE\sE_2^{\be_1}
~=~
\cls^s\sE E_{\be_2}\sE_2^{\be_1}
~\subset~
\cls^s\sE\sE_2^{\be_1}E_{\be_2}
~=~
\cls^sFE_{\be_2}
\\
~\subset~
\ldots
~\subset~
\cls^sFE_{n\be_2}.\qedsymbol
}\emun
\noqed

\vspace{-5ex}
\bob \label{amiMdilob}
Now, recall that Example \ref{SCex4} fulfills the hypotheses of all of the Theorems \ref{nPSspthm}, \ref{facpsgenthm}, and \ref{naminthm}. Therefore,
inputting Example \ref{SCex4} in the construction that proves Theorem \ref{N02Mdilthm}, we get a full Markov dilation that is incompressible (Theorem \ref{facpsgenthm} and its corollary) and not algebraically minimal (Theorem \ref{naminthm}). By algebraic minimalization, we get an algebraically minimal dilation that is not full (Corollary \ref{incmnfcor}). Since the former is primary and since $\ol{\cA_\infty\xi}^s\ne\sE$ (this is, how we proved Theorem  \ref{naminthm} by appealing to Theorem \ref{fminthm}), the superproduct system of the latter is a proper subsystem of the (super)product system of the former; by Theorem \ref{facpsgenthm}, it is a proper superproduct system. By maximal (strong, because \nbd{E_0})compression of the latter as in Theorem \ref{vNE0compthm}, we get an algebraically minimal dilation that is incompressible. Since, by Corollary \ref{comPPScor}, the superproduct system of the compression 
%%%% BO new
%is contained, 
is
a superproduct subsystem of the superproduct system of the
algebraically minimalized
dilation, 
%%%% EO
it is proper, too. By Corollary \ref{algminincfullcor}, the compression is not full.

For the punchline: We, therefore, got an algebraically minimal incompressible Markov dilation that is not full.
\eob

%%%% BO new
%WONDERFUL!
%%%% EO

% \OW[SOMEWHERE (See my remarks in blue. ---Orr.)]{
% We know, in the meantime, that a dilation of a normal CP-semigroup on some von Neumann algebra containing $\cB$ does not allow to conclude that we have a dilation of the original semigroup on $\cB$: We simply have, in general, no reason to assume that the von Neumann algebra $\cA$ on which the dilation of the normal semigroup lives, contains an invariant (unital) \nbd{C^*}subalgebra $\cA_0$ such that $p\cA_0p=\cB$. We should say this somewhere, possibly with an example.
% (This is addressed in the observation following the theorem on existence of dilation for two-parameter discrete. The fact that the problem can be shipped around in the case of one parameter (got it!) shows that it will be hard to give an example, I give up. --Orr)}

\newpage

\appendix
\addcontentsline{toc}{section}{Appendices:}

\section{A brief on von Neumann modules} \label{vNAPP}

Many results in these notes (basically, until Section \ref{EXsubnsupSEC}), we proved for (unital) \nbd{C^*}algebras and Hilbert modules; they (usually) promote easily to the situation of von Neumann algebras and von Neumann modules (frequently, adding also the word ``normal'' to the maps between such objects). Starting from Section \ref{topSEC}, more and more results (though not all) pop up which hold only in the von Neumann case (or, at least, we can prove them only in this setting). For both the easiness with which the results of the first 18 sections promote to the von Neumann setting as well as for having efficient proofs for the results that only hold in the latter, it is indispensable to have at hand the ``right language'' to deal with this setting; the ``right language'' for this is that of von Neumann modules and von Neumann correspondences.

A von Neumann algebra is, unlike an abstract \nbd{W^*}algebra, a concrete operator \nbd{*}algebra $\cB$ acting nondegenerately on a Hilbert space $G$ (distinguished by assuming $\cB$ is strongly closed or by assuming $\cB$ coincides with its double-commutant, $\cB''$); it is only natural that von Neumann modules and correspondences are modules of operators (distinguished by assuming they are strongly closed or, see \ref{vNvsrep}, a sort of double-commutant theorem). A bit surprisingly, despite the fact that it has been known how to (canonically!) transform a Hilbert module into an operator module practically as long as there exist Hilbert modules (no later than Rieffel \cite[Proposition 6.10]{Rie74a}), the first definition for the von Neumann objects among Hilbert modules and correspondences were rather \nbd{W^*}modules and \nbd{W^*}correspondences (Baillet, Denizeau, and Havet \cite{BDH88}). Only as late as Skeide \cite{Ske00b}, von Neumann modules have been defined (as as Hilbert modules over von Neumann algebras for which the canonically associated operator module is strongly closed), and, later in \cite{Ske06b}, concrete von Neumann modules (directly as strongly closed operator modules).

In these notes, we use consequently von Neumann modules and von Neumann correspondences that have to be ``\it{concreteified}'' by the canonical construction in Footnote \ref{StineFN} in Section \ref{intro}. The \it{concrete} categories only play a role in Subsection \ref{vNAPP}\ref{vNcomm} about the \it{commutant} of (concrete) von Neumann correspondences (see the crucial Observation \ref{noncob}), when we compare the results in these notes, based on the approach to dilation in Bhat and Skeide \cite{BhSk00} (GNS-correspondences) with results based on the approach in Muhly and Solel \cite{MuSo02} (Arveson-Stinespring correspondences).

% \newpage

\lf
\subsection{\normalsize Von Neumann modules and von Neumann correspondences} \label{vNmcSSEC}

Recall that a \hl{von Neumann algebra} is a concrete \nbd{*}algebra $\cB\subset\sB(G)$ of operators on a Hilbert space $G$ that is closed in the strong operator topology of $\sB(G)$ and that contains $\id_G$. (Equivalently, $\cB''=\cB$.) Von Neumann algebras are \nbd{W^*}algebras; to turn a \nbd{W^*}algebra into a von Neumann algebra, one, first, has to choose a faithful, normal, nondegenerate representation. 

Recall from Footnote \ref{StineFN}, that any (pre-)Hilbert module $E$ over a concrete operator algebra $\cB\subset\sB(G)$ can be turned into a concrete module of operators in $\sB(G,E\odot G)$ by amplification $E$ $\leadsto$ $E\odot\id_G$ and that $\sB^a(E)$ is acting (faithfully and nondegenerately) on the Hilbert space $E\odot G$ as $\sB^a(E)\odot\id_G$. Recall, too, that the whole construction is unique in the very specific sense explained in that footnote.

We follow Skeide \cite{Ske00b}, and say a Hilbert module $E$ over a von Neumann algebra $\cB\subset\sB(G)$ is a \phantomsection\hl{von Neumann \nbd{\cB}module}\index{von Neumann!module|bf} if it fulfills one of the two properties
\begin{itemize}
\item
$E\odot\id_G$ is strongly closed in $\sB(G,E\odot G)$, or

\item
$\sB^a\rtMatrix{\cB\\E}\odot\id_G\subset\sB\rtMatrix{G\\E\odot G}$ is a von Neumann algebra,
\end{itemize}
which, by \cite[Proposition 4.5]{Ske00b}, are equivalent. A third equivalent property is (\cite[Theorem 4.16]{Ske00b}) that $E$ is \phantomsection\hl{self-dual}\index{self-dual}\index{module!Hilbert!self-dual} (that is, every $\Phi\in\sB^r(E,\cB)$ has the form $\Phi x=\AB{y,x}$ for some $y\in E$).%
\footnote{
A Hilbert module over a \nbd{W^*}algebra is a \phantomsection\hl{\nbd{W^*}module}
\index{von Neumann!W@\nbd{W^*}}
 if it is self-dual; see Baillet, Denizeau, and Havet \cite{BDH88}. Therefore, a Hilbert module over a von Neumann algebra is a von Neumann module if and only if it is a \nbd{W^*}module. All this gives the intrinsic result about \nbd{W^*}modules that a Hilbert module $E$ over a \nbd{W^*}algebra $\cB$ is a \nbd{W^*}module if and only if its \phantomsection\hl{extended linking algebra}\index{linking algebra!extended} $\sB^a\rtMatrix{\cB\\E}$ is a \nbd{W^*}algebra; but it would be not so easy to prove this without, first, transforming the \nbd{W^*}module into a von Neumann module by choosing a faithful (normal) representation of $\cB$.
\vspace{1ex}
}
In particular, bounded right linear maps between von Neumann modules are adjointable.

Several features of von Neumann modules, immediately underline their simplicity. For instance, sitting as blocks in $\sB\rtMatrix{G\\E\odot G}$, all operations among the corners (these include the right and left actions on the module $E$ by elements in $\cB$ and in $\sB^a(E)$, respectively, and the inner products $\AB{x,y}=x^*y$ as well as the rank-one operators $xy^*$) share all the good (and, of course, also the bad) properties of multiplication in the von Neumann algebra $\sB^a\rtMatrix{\cB\\E}\odot\id_G$. In particular, they are jointly continuous on bounded subsets and separately continuous in general, for the strong operator topology. (Likewise, for the \nbd{*}strong, the \nbd{\sigma}strong, and the \nbd{\sigma}\nbd{*}strong operator topology.)%
\footnote{
Despite these properties can be formulated and proved directly for \nbd{W^*}modules, the direct approach is tedious -- more so, than necessary. The relevant topologies from \cite{BDH88}, the \it{\nbd{\sigma}topology} and the \it{\nbd{s}topology}, on the \nbd{W^*}module are exactly the \nbd{\sigma}weak and the \nbd{\sigma}strong topology of the corresponding von Neumann module (no matter from which representation the latter has been obtained).
}
By self-duality (that of the submodule, not that of the containing module), the strongly closed submodules of a von Neumann module are in one-to-one correspondence with the projections in $\sB^a(E)$. Every element $x$ of  a von Neumann \nbd{\cB}module $E$ admits a (unique) \phantomsection\hl{polar decomposition}\index{polar decomposition} $x=v\abs{\AB{x,x}}$ where $v\in E$ is a partial isometry (such that $\ker v=\ker x$). Consequently, a von Neumann module admits a \phantomsection\hl{quasi orthonormal system}\index{quasi orthonormal system} $(e_s,p_s)_{s\in S}$ (that is, the $e_s$ are elements in $E$ and the $p_s$ are projections in $\cB$ such that $\AB{e_s,e_{s'}}=p_s\delta_{s,s'}$) which is \hl{complete}\index{quasi orthonormal system!complete} (that is, $\sum_se_se_s^*=\id_E$ in the strong topology of $\sB^a(E)$); see \cite[Theorem 4.11]{Ske00b}.

One of the most striking (and useful) features of von Neumann modules is that they are easy to obtain.%
\footnote{
Better: Von Neumann modules are easy to obtain in a form that allows to easily compute the inner products of their elements. Frequently, in papers using \nbd{W^*}modules one reads ``let $E$ be a pre-Hilbert module over a \nbd{W^*}algebra, and let $F:=\sB^r(E,\cB)$ denote \phantomsection\hl{Paschke's self-dual completion}\index{Paschke's self-dual completion}\index{self-dual!Pascke's completion}'', referring to \cite[Corollary 4.3]{Pas73}. So, the problem is not, actually, to identify the right space $\sB^r(E,\cB)$, but to explain how the inner product of $E$ extends to elements from that space. This is done by approximating suitably elements of the latter by elements of the former, and proving that this makes sense. But it remains that to, actually, compute their inner products one does have to approximate them; on the contrary, for the (unique, minimal) self-dual extension we describe for von Neumann modules, we simply compute the inner product as multiplication of operators between Hilbert spaces.
\vspace{1ex}
}
More precisely: If $E$ is a pre-Hilbert module over a von Neumann algebra $\cB\subset\sB(G)$, then $\ol{E}^s$ (strong closure in $\sB(G,E\odot G)$) is the unique minimal von Neumann \nbd{\cB}module containing $E$. By \hl{unique minimal} we mean that whatever the von Neumann \nbd{\cB}module $F\supset E$, the canonical map $E\rightarrow F$ extends to a (unique bounded linear) map (necessarily an embedding) $\ol{E}^s\rightarrow F$.%
\footnote{
Just a little moment's thought (playing a bit around with different faithful representations of $\cB$) shows that we may replace ``unique minimal von Neumann \nbd{\cB}module containing $E$'' by ``unique minimal \nbd{W^*}module over $\cB$ containing $E$''.
}
(Later, we will see that $\ol{E}^s$  can also be obtained with the help of a \it{double commutant theorem}. In the latter form, the result is known since no later than Rieffel \cite[Proposition 6.10]{Rie74a}.) In this context, it is clear that we have a \phantomsection\hl{Kaplansky density theorem}:
\index{Kaplansky density theorem} The unit ball of $E$ is strongly dense in the unit ball of $\ol{E}^s$. (This is, only for instance, of outstanding importance in controlling properties of maps defined on tensor products, by controlling the properties only on simple tensors --, namely, simple tensors from bounded subsets.)

Here is an illustration of one aspect why the strong closure procedure is so effortless in applications. 

\bob \label{strextob}
If we have a (pre-)Hilbert module $E$ over a von Neumann algebra $\cB\subset\sB(G)$ and an element in $a\in\sB^a(E)\subset\sB(E\odot G)$, then the action of this operator $a$ on elements $x\in E\subset\sB(G,E\odot G)$ is simply by operator multiplication. Therefore, by separate strong continuity of this product, the action extends to $x=\lim_\lambda x_\lambda\in\ol{E}^s$ (strong limit) sending $x$ to $ax$ (operator multiplication), which by $ax=a\lim_\lambda x_\lambda=\lim_\lambda ax_\lambda$ is, clearly, an element of $\ol{E}^s$.

From such considerations it also follows that (co)isometries on (pre)Hilbert modules over von Neumann algebras extend as, automatically adjointable, (co)isometries to the strong closure. And for checking if the extension of an isometry is onto the closure, it is enough to check if the original range is strongly dense. (Of course, the whole observation remains true for maps not on but between modules.)
\eob

\lf
Of course, who wishes to work seriously with von Neumann algebras has to know all the usual locally convex topologies and their interplay very well; the same holds for von Neumann modules. But experience shows that for \it{organization of thoughts} it is very convenient to concentrate (in particular, in the formulation of definitions and theorems) on a single one, using one of the other topologies only if strictly necessary to conclude a proof; this topology, for us, is the strong topology. There is one exception, though: The strong topology is \bf{not} intrinsic to the algebraic structure, but depends on the faithful representation in question. In particular, the strong topology is, usually, not preserved even under ``good''($=$ normal) homomorphisms, while intrinsic topologies like all the \nbd{\sigma}versions are preserved under normal homomorphisms. So, in order to deal with such questions, we have to add one more ingredient. And we prefer, in fact, \hl{normality} as our ingredient of choice. Recall that normality is not a topological requirement, but means \hl{order continuity}: A positive linear map $T\colon\cA\rightarrow\cB$ is \phantomsection\label{normald}\hl{normal}\index{normal!positive linear map}\index{normal!order@$=$order continuity} if $\sup_\lambda T(a_\lambda)=T(\sup_\lambda a_\lambda)$ for all families $a_\lambda$. (It is sufficient to check $\lim_\lambda T(a_\lambda)=T(\lim_\lambda a_\lambda)$ for all increasing bounded nets.) Since order is a matter only of the algebraic structure, it is clear that normality is intrinsic to the algebraic structure. The crucial fact about normality is that, for positive maps, it is equivalent to \nbd{\sigma}weak continuity. Since the former is a question of bounded subsets only, so is the latter.

Now, what is \it{normality} in the context of von Neumann modules? We follow the philosophy to motivate terminology for modules in terms of their linking algebras. A linear map $T\colon E_\cB\rightarrow F_\cC$ is \it{normal} if it admits normal extension $\sT$ to a (necessarily positive!) map acting \phantomsection\hl{blockwise}\index{blockwise}\index{linking algebra!blockwise linear map} between the linking algebras $\sB^a\rtMatrix{\cB\\E}$ and $\sB^a\rtMatrix{\cC\\F}$ (meaning $\sT(\cB)\subset\cC$, $\sT(\sB^a(E))\subset\sB^a(F)$, $\sT(E)\subset F$, and, hence necessarily, $\sT(E^*)\subset F^*$).

\bex \label{normetaex}
A trivial example is the canonical injection $E\rightarrow E\odot\id_G\subset\sB(G,E\odot G)$. The extension to the canonical injection $\sB^a\rtMatrix{\cB\\E}\rightarrow\sB^a\rtMatrix{G\\E\odot G}$ is, by definition, normal.

This gets less trivial, if we replace $G$ carrying the identity representation, with a Hilbert space $K$ carrying another (for convenience, nondegenerate) representation $\pi\colon\cB\rightarrow\sB(K)$. The whole procedure from Footnote \ref{StineFN} goes through as before, just that now the resulting ``embeddings'' need no longer be faithful and earns a letter: $\eta\colon x\mapsto x\odot\id_K$ defines a linear map onto a Hilbert \nbd{\pi(\cB)}submodule of $\sB(K,E\odot K)$, satisfying $\AB{\eta(x),\eta(y)}=\pi(\AB{x,y})$ and $\eta(xb)=\eta(x)\pi(b)$. And there is a (nondegenerate) representation $\Pi\colon a\mapsto a\odot\id_K$ of $\sB^a(E)$ on $E\odot K$. We leave it as an instructive exercise to the reader (also making use of tools like polar decomposition for elements of $E$ and cyclic decomposition of $K$) to show that the representation $\prod:=\rtMatrix{\pi&\eta^*\\\eta&\Pi}$ of $\sB^a\rtMatrix{\cB\\E}$ on $\rtMatrix{K\\E\odot K}$ (and, therefore, also $\eta$) is normal if $\pi$ is normal. (One may add: If $E$ is \phantomsection\hl{strongly full}\index{strongly full}\index{von Neumann!module!strongly full} (that is, if $\cls^s\AB{E,E}=\cB$) and if (now really, and not only for convenience) $\pi$ is nondegenerate, then it suffices to show that $\eta$ is \nbd{\sigma}weakly continuous.)
\eex

With this in mind, dealing with von Neumann correspondences and their tensor products, is really easy. A correspondence $_\cA E_\cB$ from a von Neumann algebra $\cA$ to a von Neumann algebra $\cB$ is a \phantomsection\hl{von Neumann correspondence}\index{von Neumann!correspondence|bf}  if $E_\cB$ is von Neumann \nbd{\cB}module and if the \phantomsection\hl{Stinespring representation}\index{von Neumann!correspondence!Stinespring representation}\index{Stinespring representation!of a von Neumann correspondence} $\rho\colon\cA\rightarrow\sB^a(E)\rightarrow\sB(E\odot G)$ is normal. (Equivalent to the latter is the intrinsic criterion that for all $x\in E$, the map $\cA\ni a\mapsto\AB{x,ax}\in\cB$ is normal; see Skeide \cite[Lemma 3.3.2]{Ske01} for this and one more equivalence.)

\bex \label{nGNSex}
If the CP-map $T$ between von Neumann algebras $\cA$ and $\cB$ is normal, then the strong closure $\ol{\sE}^s$ of its GNS-correspondence $\sE$ is a von Neumann correspondence. Again, we leave the routine argument to the reader.
\eex

The (von Neumann) \phantomsection\hl{tensor product}\index{von Neumann!correspondence!tensor product}\index{tensor product!of von Neumann correspondences} of a von Neumann correspondence $_\cA E_\cB$ and a von Neumann correspondence $_\cB F_\cC$ is the von Neumann \nbd{\cA}\nbd{\cC}correspondence
\beqn{
E\sodots F
~:=~
\ol{E\odot F}^s.
}\eeqn
Let us invest the time to produce the thoughts of a moment, illustrating that this is meaningful. To compute $\ol{E\odot F}^s$ we have to tensor $E\odot F$ with the representation space $L$ of $\cC\subset\sB(L)$. Clearly (in the tensor category we are working in), $(E\odot F)\odot L=E\odot F\odot L=E\odot(F\odot L)$. Since $F$ is a von Neumann correspondence, the Stinespring representation of $\cB$ on $K:=F\odot L$ is normal. By Example \ref{normetaex}, so is the representation of $\cA$ on $E\odot(F\odot L)$. Hence, the von Neumann \nbd{\cC}module $\ol{E\odot F}^s$ with the left action of $\cA$ by this representation, is a von Neumann correspondence. (Compare this with Observation \ref{strextob}.)

\brem \label{etanrem}
Of course, again by Example \ref{normetaex}, also the map $\eta\colon x\mapsto x\odot\id_{F\odot L}$ is normal. Moreover, as an element of $\sB(L,E\odot F\odot L)\supset E\odot F$, the simple tensor $x\odot y$ is given simply by operator multiplication $\eta(x)y$ of the operator $\eta(x)\in\sB(F\odot L,E\odot F\odot L)$ and the operator $y\in\sB(L,F\odot L)$. By normality of $\eta$ and the Kaplansky density theorem, the operators $x\odot y=\eta(x)y$ with $x$ and $y$ from the unit balls of $E$ and $F$, respectively, form a strongly total subset of (the unit ball of) $E\sodots F$. (Even more: The unit ball of the space formed by their finite linear combinations is strongly dense in the unit ball of $E\sodots F$.) Observations like that can be extremely useful in arguments.
\erem

We are now ready to illustrate in a single, though sufficiently complex, example why the results in the first \ref{EXsubnsupSEC} sections, promote easily to the von Neumann case. This example is very typical; it is not only the only truly closed circle in these notes where existence of certain dilations for a CP-semigroup is equated in the form of an \it{iff}-theorem with a concise property of the GNS-subproduct system of the CP-semigroup; but, actually, also explaining in detail what exactly its ingredients and statement mean in the von Neumann case, does already cover quite a lot of the topics dealt with in the \nbd{C^*}case. The promotion of the remaining results should, then, be just a simple adaptation, and is left to the reader.

\lf
\subsection{\normalsize Theorem \ref{Markmodthm}; von Neumann version} \label{MmvNSSEC}

\bthm \label{vNMarkmodthm}
Let $T$ be a normal Markov semigroup over the opposite of an Ore monoid. Then $T$ admits a normal full dilation if and only if the GNS-subproduct system of von Neumann correspondences of $T$ embeds into a product system of von Neumann correspondences.
\ethm

The theorem has two directions. For the forward direction, we have to clarify how the construction of a (super)product system of a (full) dilation promotes to the von Neumann case. For the backwards direction, relying basically on Theorem \ref{Oreindthm}, we have to clarify the same for the inductive limit construction from Observation \ref{indlimob}.

First of all, let us repeat (see the paragraph following Observation \ref{hypevNob}) that the notion of fullness of a dilation is not affected by the passage to the von Neumann case. Just, if the algebra $\cA\subset\sB(H)$ on which the dilating \nbd{E}semigroup acts is a von Neumann algebra, then $E:=\cA p$ is a von Neumann module over $\cB=p\cA p\subset\sB(G)$, where $G=pH$, sitting already concretely (though, not necessarily nondegenerately) in $\sB(G,H)=\sB(H)p$. Then $\sB^a(E)$ is a von Neumann algebra and can be obtained both as strong closure or as strict closure of $\sF(E)=\ls EE^*$. And full means the von Neumann algebra $\sB^a(E)$ coincides with $\cA$. (Since $\cA$ acts nondegenerately on $H$, necessarily so does $E$ on $G$.)

The next to be clarified, is the notion of (sub)(super)product system of von Neumann correspondences ($E^{\bar{\bodot}^s}$) ($E^{\bar{\podot}^s}$) $E^{\sodots}$. Of course, this means very much the same as (sub)(super)product system of \nbd{C^*}correspondences; just that now the correspondences are (also) required to be von Neumann correspondences and that their tensor product $E_s\odot E_t$ of \nbd{C^*}correspondences is replaced by the tensor product $E_s\sodots E_t$ of von Neumann correspondences. (Recall the methods from Observation \ref{strextob} to check (co)isometricity.)

From Example \ref{nGNSex}, it follows that the strong closures $\ol{\sE}^s_t$ of the GNS-correspondences $\sE_t$ form a subproduct system $\sE^{\bar{\bodot}^s}$ of von Neumann correspondences and the elements $\xi_t\in\sE_t\subset\ol{\sE}^s_t$ form a unit, generating $\sE^{\bar{\bodot}^s}$ as a subproduct system of von Neumann correspondences.

As for the (super)product system of an \nbd{E}semigroup $\vt$ on $\sB^a(E)\ni p=\xi\xi^*$ we have learned two ways how to obtain it. The first, $E^*\odot_tE$ (see Example \ref{EPSex}), does already involve the tensor product over $\sB^a(E)$. This approach raises two questions, which we like to avoid. (See Footnote \ref{E*E=ptEFN}.) The second way, $E_t:=\vt_t(\xi\xi^*)E$ (see Theorem \ref{E-supPSthm}), which we follow here, has the advantage that $E_t$ defined in this way is manifestly a von Neumann \nbd{\cB}module; it is a von Neumann correspondence over $\cB$ if (and only if the restriction of) $\vt_t$ (to $\cB$) is normal.%
\footnote{ \label{E*E=ptEFN}
To the curious reader: If $E$ is a von Neumann \nbd{\cB}module and if $\vt$ is just any (not necessarily normal nor strict) endomorphism of the von Neumann algebra $\sB^a(E)$, then $E^*\odot_\vt E$ is, indeed, a von Neumann \nbd{\cB}module; and that even if instead of the tensor product of \nbd{C^*}correspondence, we take only the linear span of all $x\odot_\vt y$. The point is that if $E$ has a unit vector $\xi$, then
\beqn{
\sscls^s E^*\ul{\odot}_\vt E
~=~
\sscls^s\xi^*\sF(E)\ul{\odot}_\vt E
~=~
\sscls^s\xi^*\ul{\odot}_\vt\sF(E)E
~=~
\xi^*\ul{\odot}_\vt E
}\eeqn
(the latter, because for every Hilbert module we have $E\AB{E,E}=E$). Clearly, $\xi^*\odot_\vt x\mapsto\vt(\xi\xi^*)x$ establishes an isomorphism of the \nbd{C^*}correspondences $E^*\odot_tE$ and $\vt(\xi\xi^*)E$, which, therefore, are both von Neumann modules, because $\vt(\xi\xi^*)E$ is. Using the results of Skeide \cite[Section 4]{Ske09}, but taking now really into account that $\vt$ has to be normal,  it is easy to reduce the case without unit vector to the case with unit vector. With our choice to write $\vt(\xi\xi^*)E$, we need not worry about these questions.
\vspace{1ex}
}

Now recall that the $E_t$ form a product system of \nbd{C^*}correspondences if $\vt$ is strict, and recall that this is due to the fact that strictness of $\vt_t$ is equivalent to the nondegeneracy condition $\cls\vt_t(\sF(E))E=\vt_t(\id_E)E$. By exactly the same argument, we prove that the $E_t$ form a product system of von Neumann correspondences if $\vt$ is normal, the latter being equivalent to $\cls^s\vt_t(\sF(E))E=\vt_t(\id_E)E$. (See Muhly, Skeide, and Solel \cite[Section 1]{MSS06} for more details.) In particular, the superproduct system of a normal dilation of $T$ containing the GNS-subproduct system, passing to strong closures (transforming the superproduct system into a \bf{product} system of von Neumann correspondences, containing the unit, hence, the strong closure of the GNS-subproduct system), we have proved the forward direction of Theorem \ref{vNMarkmodthm}

To conclude also the backward direction, we do the inductive limit construction in Observation \ref{indlimob} (after all, if the $E^{\sodots}$ is a product system of von Neumann correspondences, then the $E_t$ \bf{do} form a superproduct system of \nbd{C^*}correspondences), but we replace the resulting Hilbert module over the von Neumann algebra $\cB$ by its strong closure $E$.%
\footnote{
It is noteworthy that an inductive limit of von Neumann correspondences obtained by the strong closure of the inductive limit considering them as \nbd{C^*}correspondences, is a von Neumann correspondence, too. (Very much along the lines of the proof of  \cite[Lemma 3.3.2]{Ske01}, we show that the Stinespring representation is normal.) This is relevant, when promoting also the inductive limit in Theorem \ref{indlimthm} to the von Neumann case.
}
It remains to show that the $v_t$ defined after Observation \ref{indlimob}, are unitary. This follows because the density arguments in the proof of Proposition \ref{lduniprop}, taking also into account Example \ref{normetaex} and Remark \ref{etanrem}, easily turn over to the von Neumann case.

\lf
\subsection{\normalsize Von Neumann $\sB(G)$--modules} \label{vNBGmod}

The simplest von Neumann algebra is, in a sense, $\sB(G)$. (In fact, $\sB(G)$ is a \hl{factor}$=$a von Neumann algebra with trivial center=a von Neumann algebra with no nontrivial (strongly closed) ideals, and decomposition theory shows that von Neumann algebras (with separable predual) may be decomposed into factors. And except for the trivial case $\dim G=1$, resulting in the only commutative factor, there is not really a big difference between different $\sB(G)$.)

A \nbd{\cB}(bi)module $E$ is \hl{trivial} if it is a tensor product $V\otimes\cB$ (or $\cB\otimes V$) with a vector space. A basis $\bfam{e_s}$ for $V$ promotes to a (bi)module basis $\bfam{e_s\otimes\U}$ for $E$; usual modules do not have bases, and even if a bimodule has a right module basis, the right module basis can usually not be chosen to be a bimodule basis. (Typical example: $E={_\vt}\cB$ for an endomorphism of $\cB$ that is not an inner automorphism.)

Already in Bhat and Skeide \cite{BhSk00} it has been pointed out that von Neumann \nbd{\sB(G)}mod\-ules have the form $E=\sB(G,H)$ with $\sB^a(E)=\sB(H)$ \cite[Proposition 13.9]{BhSk00}, and that von Neumann \nbd{\sB(G)}correspondences are trivial \cite[Theorem 13.11]{BhSk00}. (Accepting that quantum dynamics on $\cB\subset\sB(G)$ can be conveniently dealt with in terms of Hilbert modules and correspondences, the fact that these are so simple when $\cB=\sB(G)$ is responsible for that modules are not absolutely indispensable when the dynamics is limited to $\sB(G)$. But even then, modules are still useful to guide thought and simplify arguments. A good example is Fagnola and Skeide \cite{FaSk07} and the successive examples with Bhat \cite{BFS08}, which contain results that even in the \nbd{\sB(G)}case could have been proved only with enormous difficulty when not having modules in mind.) Let us see, why.

Let $E\subset\sB(G,E\odot G)$ be a von Neumann \nbd{\sB(G)}module. Since $\sB(G)$ contains all rank-one operators, $E$ contains all rank-one operators of the form $(x\odot\id_G)g'g^*=(x\odot g')g^*$ ($x\in E$ and $g,g'\in G$). Since the $x\odot g'$ are total in $E\odot G$, the closed linear space is $\sK(G,E\odot G)$, and since $E$ is strongly closed, we have $E=\sB(G,E\odot G)$. We know that $\sB^a(E)\subset\sB(E\odot G)$, in general. But, of course, every operator on $E\odot G$ leaves the module $\sB(G,E\odot G)$ invariant, so $\sB^a(E)=\sB(E\odot G)$.  We see, the category of von Neumann \nbd{\sB(G)}modules is naturally equivalent to the category of Hilbert spaces.

From now on we will write $H=E\odot G$.

Whenever $\dim H$ is infinite and not smaller than $\dim G$, it is easy to see that $E$ does not only admit a complete quasi orthonormal system but a complete orthonormal system (that is, $\AB{e_s,e_s}=\U$ for all $s$), that is, a (topological) module basis and, hence, is trivial. (Simply slice $H$ into copies of $G$.) Of course, if $\dim H<\dim G$, then $\sB(G,H)$ does not admit unit vectors (that is, isometries in $\sB(G,H)$).%
\footnote{
It is noteworthy that even for general von Neumann modules, missing dimension of $H$ (and being not strongly full, of course) is the only reason for failing to have a complete orthonormal system: For every strongly full von Neumann \nbd{\cB}module $E$ there is a cardinality $\en$ such that $\ol{E^\en}^s\cong\ol{\cB^\en}^s$. (See the discussion leading to \cite[Corollary 4.3]{Ske09}.) Or, better, for every von Neumann module, we can find $\en$ such that $\ol{E^\en}^s\cong\ol{\cB_E^\en}^s$.
\vspace{1ex}
}
But, right in the case $\dim H=\dim G\ge\aleph_0$, the cardinality of  complete orthonormal system is never unique. (We may slice $H=G$ into $\en$ copies of $G$ for any cardinality $\en\le\dim G$.)%
\footnote{
Be aware that for Hilbert modules the cardinality of a complete  orthonormal system (if it exists) is unique. (See, for instance, Landi and Pavlov \cite[Proposition 3.1]{LaPa12}.) Of course, \it{complete}, for Hilbert modules, means that the \nbd{\cB}linear combinations are norm dense, while, for von Neumann modules, they need be only strongly dense. A Hilbert \nbd{\cB}module over a von Neumann algebra $\cB$ with an infinite complete orthonormal system (as Hilbert module) can never be a von Neumann module, unless the von Neumann algebra $\cB$ is finite-dimensional.
}

Now, let $E=\sB(G,H)$ be a von Neumann \nbd{\sB(G)}correspondence. Taking into account that in the category of von Neumann correspondences, $G$ is really a Morita equivalence from $\sB(G)$ to $\C$, as in the discussion of Morita equivalence starting after Remark \ref{CPASrem}, we get
\beqn{
E
~=~
G\sodots G^*\sodots E\sodots G\sodots G^*
~=~
G\sodots\eH\sodots G^*
~=~
\ol{G\otimes\eH\otimes G^*}^s
~=~
\sB(G,G\otimes\eH),
}\eeqn
with the \nbd{\C}correspondence (that is, the Hilbert space) $\eH:=G^*\sodots E\sodots G=G^*\odot E\odot G$. We see that $H=E\odot G=G\otimes\eH$ with the canonical left action of $\sB(G)$ as $\sB(G)\otimes\id_\eH$. We see that, again, the category of von Neumann \nbd{\sB(G)}correspondences is naturally equivalent to the category of Hilbert spaces; just the Hilbert space ``belonging'' to $E$ is not $H$ but $\eH$. The passage between morphisms is established by the map $\id_{G^*}\odot\bullet\odot\id_G\colon\sB^a(E_1,E_2)\rightarrow\sB(\eH_1,\eH_2)$ and its inverse $\id_G\otimes\bullet\otimes\id_{G^*}\colon\sB(\eH_1,\eH_2)\rightarrow\sB^a(E_1,E_2)$.

They are even equivalent as tensor categories:

The tensor product of two von Neumann correspondences $E_1=\sB(G,G\otimes\eH_1)$ and $E_2=\sB(G,G\otimes\eH_2)$ is, obviously, $E_1\sodots E_2=\sB(G,G\otimes\eH_1\otimes\eH_2)$. Note that, in view of Remark \ref{etanrem}, we may identify $x_1\odot y_2$ with $(x_1\otimes\id_{\eH_2})y_2$ (operator multiplication). Applying this to the Fock modules we get $\sF(E)=\sB(G,G\otimes\sF(\eH))$ and $\DG(E)=\sB(G,G\otimes\DG(\eH))$. Recalling the operators on these Fock modules are $\sB(G\otimes\sF(\eH))$ and $\sB(G\otimes\DG(\eH))$, we recover the famous form of Fock space $\sF(\eH)$ or $\DG(\eH)$ tensor initial space $G$, on which in quantum probability for dynamics on $\sB(G)$ dilations are modelled. However:

\bob
It is noteworthy that the order of factors in $H=G\otimes\eH$ comes naturally from the operation of Morita equivalence of correspondences. The formula $x_1\odot y_2=(x_1\otimes\id_{\eH_2})y_2$, underlines that we are doing ``the right thing''. Indeed, if we opted to identify $E$ with $\sB(G,\eH\otimes G)$ (of course, possible), then we would meet a flip involved either in the analogue formula, or in that the tensor product of $E_1$ and $E_2$ would have to be $\sB(G,G\otimes\eH_2\otimes\eH_1)$. This shows that the habit in quantum probability of writing the ``initial space'' $G$ on the right, is likely to cause unnecessary flips in formulae. The same applies to the habit to write the factorization of the Fock space as \it{past} tensor \it{future}; not for nothing, in the time ordered Fock module, the \it{future} is on the left of the \it{past}. Theorems like the one whose von Neumann version we just proved in \ref{MmvNSSEC}, would not work, tensoring the unit (which stand for a piece of future) from the other side.
\eob

\lf
Note that $E=\sB(G,G\otimes\eH)$ may be identified with $\sB(G)\sbars{\otimes}\eH$ via $\sB(G,G\otimes\eH)\ni b\otimes x\colon g\mapsto bg\otimes x$. (Of course, while the formula $x_1\odot y_2=(x_1\otimes\id_{\eH_2})y_2$ makes sense and gives the right thing also here, the tensor product may also be recovered as $(b_1\otimes h_1)\odot(b_2\otimes h_2)=b_1b_2\otimes h_1\otimes h_2$; so also in multiple tensor products, the elements of the Hilbert spaces $\eH_i$ are tensored, while the algebra elements are multiplied together. It might be worthwhile to compare this with the tensor products of the GNS-constructions in Section \ref{qconvSEC}, which also have the form Hilbert space tensor algebra but in different order and with a nontrivial left action.) Note, too, that for every ONB $\bfam{e_s}$ of $\eH$, the elements $\U\otimes e_s$ form a complete orthonormal system for $E$.

\bex \label{Krausex}
Let $T$ be a normal CP-map on $\sB(G)$ and denote by $\sE\ni\xi$ its (strongly closed) GNS-construction. Then $\xi=\sum_sc_s\otimes e_s$ (with $c_s=\AB{\U\otimes e_s,\xi}\in\sB(G)$), so
\beqn{
T
~=~
\sum_sc_s^*\bullet c_s
}\eeqn
(\phantomsection\hl{Kraus decomposition}\index{Kraus decomposition|bf}). Moreover, if $F\ni\zeta$ is another von Neumann \nbd{\sB(G)}correspondence such that $T=\AB{\zeta,\bullet\zeta}=\sum_td_t^*\bullet d_t$ (where the $d_t$ form the Kraus decomposition arising from some ONB $\bfam{f_t}$ of $\eG:=G^*\sodots F\sodots G$), then, as we know, $\xi\mapsto\zeta$ extends as a bilinear isometry $v\colon E\rightarrow F$. It, therefore, has the form $v=\id_G\otimes\upsilon$, where $\upsilon$ is a (unique) isometry $\eH\rightarrow\eG$. We conclude:
\begin{itemize}
\item
$\dim\eG\ge\dim\eH$. Therefore, the minimal cardinality of a Kraus decomposition of $T$ is $\dim\eH$.

\item
$d_t=\sum_s\upsilon_{t,s}c_s$, where $\upsilon_{t,s}:=\AB{f_t,\upsilon e_s}$,
%%%% BO 
and $c_s=\sum_t\bar{\upsilon}_{t,s}d_t$. 
% $c_s=\sum_t\bar{\upsilon}_{t,s}d_t$, and these sums must converge strongly (being a normal projections of strongly convergent sums). 
%%%% EO
So, all Kraus decompositions of $T$ have the same strongly closed linear span.
\end{itemize}
Of course, just cardinality is not a meaningful notion of minimality for Kraus decomposition. A Kraus decomposition would merit to be called minimal, if no $d_t$ can be taken away still maintaining the same strongly closed linear space, that is, if the $d_t$ are strongly linearly independent in some sense. It is noteworthy that the numerical coefficients occurring the preceding sums are square summable. Therefore:

\bdefin
The Kraus decomposition $\bfam{d_t}$ is \phantomsection\hl{minimal}\index{Kraus decomposition!minimal}\index{minimal!Kraus decomposition}, if the $d_t$ are \hl{strongly \nbd{\ell^2}linearly independent}, that is, if $\sum_t\bar{\lambda}_td_t=0$ for some square summable coefficients $\bfam{\lambda_t}$ in $\C$ implies that $\lambda_t=0$ for all $t$.
\edefin

We show the following \it{folklore} with our methods.

\bthmn
The Kraus decomposition $\bfam{d_t}$ for $T$ is minimal if and only of the corresponding pair $(F,\zeta)$ is the GNS-construction of $T$.
\ethmn

\proof
First, observe that the family $\bfam{\lambda_t}$ being square summable, means exactly $g:=\sum_tf_t\lambda_t$ defines a vector in $\eG$. Therefore, the condition $\sum_t\bar{\lambda}_td_t=0$, means exactly that $\AB{\U\otimes g,\zeta}=0$. Since $\U\otimes g$ is central, this means $\AB{\U\otimes g,y}=0$ for all $y$ in the correspondence strongly generated by $\zeta$. Therefore, admitting nonzero solutions $g$ or not, is equivalent to the GNS-correspondence (sitting as what $\zeta$ generates in $(F,\zeta)$) admitting nontrivial complement or not.\qed

\lf
Note that $(F,\zeta)$ being the GNS-construction is equivalent to the more frequent notion of minimality for Kraus decomposition that the vectors $bd_t\otimes f_t$ form a total subset of $G\otimes\eG$.
\eex

\lf
\subsection{\normalsize The commutant of von Neumann correspondences and ``translations''} \label{vNcomm}

Apart from communicating a different point of view on von Neumann modules and correspondences, which are of independent interest and make part of a complete picture of the subject of these notes, the present subsection serves, in particular, to put our result into perspective with some earlier work which has that different point of view. The number of results  that fall under this duality of points of view, called \it{commutant}, is enormous, and we cannot even give an approximately complete account, here. While there do exist the surveys Skeide \cite{Ske05a,Ske08a} (to which we also will refer for some aspects), they do not cover all parts of the picture relevant for us. A result of this fact is that we cannot base the discussion on a single reference, but have to explore the context to some extent directly.

\brem
It all started with Bhat and Skeide \cite{BhSk00} and Muhly and Solel \cite{MuSo02}. After in \cite{BhSk00} the problem of finding a dilation of a one-parameter CP-semigroup on a (unital) \nbd{C^*}algebra $\cB$ has been solved with the help of a product system of correspondences over $\cB$, for a von Neumann algebra $\cB\subset\sB(G)$ \cite{MuSo02} presented the same dilation but constructed with the help of a product system of von Neumann correspondences over the commutant of $\cB$, $\cB'$. (The dilations are the same, simply because they both are minimal, and the minimal dilation is unique.) The natural question how these two constructions are related, has been answered in Skeide \cite{Ske03c}: By the \it{commutant} of von Neumann correspondences. The commutant of a von Neumann \nbd{\cB}correspondence (introduced in \cite{Ske03c} and, independently and for different purposes under the name of \it{\nbd{\sigma}dual} where $\sigma$ is a faithful normal representation of $\cB$ that has to be chosen, by Muhly and Solel \cite{MuSo04}) is a von Neumann \nbd{\cB'}correspondence. \cite{Ske03c} points out that the product system of \cite{MuSo02} is the commutant of the product system of\cite{BhSk00}. The commutant has been generalized to the commutant of von Neumann \nbd{\cA}\nbd{\cB}correspondences (first, in Muhly and Solel \cite{MuSo05}, under the name of \it{\nbd{\rho}\nbd{\sigma}dual}, and, then in Skeide \cite{Ske06b}), resulting in von Neumann \nbd{\cB'}\nbd{\cA'}correspondences. With this, Skeide \cite{Ske09} points out that also all the other ingredients of the constructions of \cite{BhSk00} and of \cite{MuSo02} translate into each other. (See also \cite[Section 6]{Ske08a}.)

Some more instances:
\begin{itemize}
\item
In the presence of cyclic separating vectors, the commutant of GNS-correspondences extends to a duality between normal Markov maps $T\colon\cA\rightarrow\cB$ and $T'\colon\cB'\rightarrow\cA'$, first observed by Albeverio and Hoegh-Krohn \cite{AlHK78}. Gohm and Skeide \cite{GoSk05} exploit this to translate existence of quasi orthonormal systems into existence of so-called \it{weak tensor dilations} of Markov maps. Applying the duality of Markov maps to the canonical embedding of a subalgebras $\cB$ of $\cA$ into $\cA$, taking into account also \it{Tomita conjugation}, we recover the \it{Accardi-Cecchini expectation} $\cA\rightarrow\cB$ \cite{AcCe82} (see also Longo \cite{Lon84}), which generalizes the \it{conditional expectation} to cases where the latter does not exist.

\item 
The paper \cite{Ske03c} also generalizes Arveson's construction of an Arveson system for an \nbd{E_0}semi\-group on $\sB(H)$, to the construction of a product system of von Neumann \nbd{\cB'}cor\-re\-spon\-dences from an \nbd{E_0}semigroup on $\sB^a(E)$ (also used by Alevras \cite{Ale04} for the special case of \nbd{E_0}semigroups on type II$_1$ factors), and pointed out that it is the commutant of all the other product systems of von Neumann \nbd{\cB}correspondences for this \nbd{E_0}semigroup (for instance, those constructed in these notes); this includes the relation between the product systems in \cite{BhSk00} and in \cite{MuSo02}. As a special case, this includes that the Arveson system of an \nbd{E_0}semigroup on $\sB(H)$ is the commutant of the Arveson system constructed by Bhat \cite{Bha96} (generalized in \cite{Ske02} to \nbd{E_0}semigroups on $\sB^a(E)$) and, therefore, anti-isomorphic to the latter. (The relation between that Arveson system and the Bhat system by anti-isomorphism and how this generalizes to commutants of product systems, is the motivating backbone of \cite{Ske08a}; its detailed explanation is distributed over the whole thing.)

\item
In the other direction, also the construction of \nbd{E_0}semigroups for product systems has two versions. These can be summarized as constructing for a product system a \it{left dilation} (see Section \ref{leftdilSEC}) or constructing a \it{right dilation}. Except for the terminology (which has been introduced later in \cite{Ske06}), the two problems have been phrased, compared, and resolved in the discrete one-parameter case in \cite{Ske09}; \cite{Ske09} is at the heart of every solution of the continuous time one-parameter case for general product system (\cite{Ske07,Ske09a,Ske11a,Ske16}) and motivated also the proofs for Arveson systems in \cite{Ske06,Arv06}, which, then, generalize to arbitrary product systems.

\item
The duality of Markov maps mentioned in the first issue, depends on the existence of a pair of preserved faithful (normal) states. If $(E,\xi)=$ $\ol{\text{GNS-$T$}}^s$, then the cyclic separating vectors serve to define the (cyclic!) element $\xi'$ in the commutant of $E$, $E'$, that determines $T'$ as $\AB{\xi',\bullet\xi'}$. Without that, we may construct the commutant of the GNS-correspondence, yes, but we have no distinguished cyclic vector $\xi'$. This raises the question to what the cyclic vector $\xi$ for $E$ corresponds for $E'$. In the case of \nbd{\cB}corre\-spondences (so $\cA=\cB$), the answer is, the elements of $E$ correspond to \it{covariant representations} of $E'$; see Muhly and Solel \cite{MuSo98}, and many of their forthcoming papers. This turns over to (sub)product systems and their units, respectively, covariant representations.
\end{itemize}
The last two issues of this remark will keep us busy in the remainder of this appendix. They have to be understood by those who wish to understand better the relation between several of those results that existed before (mainly the two- and  the discrete \nbd{d}parameter case) as they have been produced exclusively in the commutant picture; trying to understand them in the language of \cite{BhSk00}, was the starting point of this note.
\erem

\bemp[Von Neumann \nbd{\cB}modules \it{vs} representations of $\cB'$.~] \label{vNvsrep}
Let $E$ be a (pre-)Hilbert module over a von Neumann algebra $\cB\subset\sB(G)$. On the tensor product $E\odot G$, apart from the natural action of $\sB^a(E)$ as $\sB^a(E)\odot\id_G$, there is another natural action of those maps in $\sB(G)$ that are \nbd{\cB}\nbd{\C}bimodule maps; these maps form the set $\sB^{bil}(G)=\cB'$, the \hl{commutant} of $\cB$. With reference to the section entitled \it{lifting commutants} in Arveson's \cite{Arv69}, we call the representation $\rho'\colon b'\mapsto\id_E\odot b'$ the \phantomsection\hl{commutant lifting}\index{von Neumann!module!commutant lifting}\index{commutant lifting!of a von Neumann module}. A routine check shows that $\rho'$ is normal and, of course, $\rho'$ is unital (turning $E\odot G$ into a von Neumann \nbd{\cB'}\nbd{\C}correspondence).

For a \nbd{\cB}bimodule $M$, we denote by $C_\cB(M):=\CB{m\in M\colon bm=mb~(b\in\cB)}$ its \hl{\nbd{\cB}center}. For convenience, we repeat the following line (see Skeide \cite[Section 2]{Ske05c})
\beqn{
\SMatrix{\cB&E^*\\E&\sB^a(E)}''
~=~
\left\{\SMatrix{b'&\\&\rho'(b')}\colon b'\in\cB'\right\}'
~=~
\SMatrix{\cB&C_{\cB'}(\sB(E\odot G,G))\\C_{\cB'}(\sB(G,E\odot G))&\rho'(\cB')'},
}\eeqn
which tells us almost everything we have to know about $\rho'$. Immediately, we get the following (two-part) \phantomsection\hl{double commutant theorem}\index{von Neumann!module!double commutant theorem}\index{double commutant theorem!for von Neumann modules} for von Neumann modules: (1), $E$ is a von Neumann \nbd{\cB}module if and only if $E=C_{\cB'}(\sB(G,E\odot G))$; and (2), if $E$ is not a von Neumann module, then $C_{\cB'}(\sB(G,E\odot G))$ is its unique minimal self-dual extension.%
\footnote{
Provided that every bounded right linear map $\Phi$ from (the von Neumann \nbd{\cB}module!) $C_{\cB'}(\sB(G,E\odot G))$ to $\cB$ gives rise to a bounded map $\Phi\odot\id_G$ (necessarily in $C_{\cB'}(\sB(G,E\odot G))^*$!), we see that the intertwiner module $C_{\cB'}(\sB(G,E\odot G))$ is self-dual; see \cite[Proposition 6.10]{Rie74a} or \cite[Section 3]{Ske05c}.

The question if $\Phi\odot\id_G$ (well-defined on the algebraic tensor product $E~\ul{\odot}~G$) is bounded, seems to be crucial in all proofs of self-duality of von Neumann modules. Rieffel \cite{Rie74a} refers back to his older result \cite{Rie69} about Banach modules over \nbd{C^*}algebras. Skeide \cite[Remark 4.3]{Ske05c} explains a direct proof for von Neumann modules (recovered also in \cite{BMSS12}). We would like to mention that the proof for (pre-)Hilbert modules in Skeide \cite[Lemma 3.9]{Ske00b} (or \cite[Lemma 2.3.7]{Ske01}) is incomplete, and we do not see a possibility how to fix that directly. (The proof is correct if $E$ is a von Neumann module.) The statement is, of course, true and \it{folklore} in Hilbert module theory.
}

We see, a von Neumann \nbd{\cB}module $E$ can be recovered from a normal unital representation $\rho'$ of $\cB'$ on the Hilbert space $H:=E\odot G$. The piece that is missing to obtain a full duality between von Neumann \nbd{\cB}modules and representations of $\cB'$, is the following result; see, for instance, the proof of Muhly and Solel \cite[Lemma 2.10]{MuSo02}: If $\rho'$ is a unital representation of $\cB'$ on a Hilbert space $H$, then  $\rho'$ is normal (if and) only if  \,$\cls C_{\cB'}(\sB(G,H))G=H$.  (Roughly, if $\rho'$ is normal and unital, then the middle term in the above line (without the $'$) is a von Neumann algebra on $\rtMatrix{G\\H}$ and $\id_{\cB'}\oplus\rho'$ is an isomorphism from $\cB'$ onto that von Neumann algebra. The result follows by examining the central cover of $p_G$ in the commutant of that subalgebra of $\sB\rtMatrix{G\\H}$.) It follows that the von Neumann \nbd{\cB}module $C_{\cB'}(\sB(G,H))$ acts nondegenerately on $G$ so that (see Footnote \ref{StineFN}) $H\cong E\odot G$ and, under this identification, $\rho'(b')=\id\odot b'$.

One easily verifies that $E$ is strongly full if and only if $\rho'$ is faithful. More generally, if $p$ is the central projection such that $p\cB=\ol{\cB_E}^s$ and if $q$ is the central projection such that $q\cB=\ker\rho'$, then $p+q=\id_G$.
\eemp

\brem
Introducing, as in Skeide \cite{Ske06b}, the category of \hl{concrete von Neumann \nbd{\cB}modules}\index{von Neumann!module!concrete} as pairs $(E,H)$ of subspaces $E\subset\sB(G,H)$ fulfilling $E\cB\subset E$, $E^*E\subset\cB$, $\cls EG=H$, and $E=\ol{E}^s$, with morphisms $\sB^a(E_1,E_2)$, we get a true one-to-one correspondence with the objects $(\rho',H)$ and the morphisms $\sB^{bil}(H_1,H_2)$ of the category of normal unital representations of $\cB'$. The category of von Neumann \nbd{\cB}modules is naturally equivalent to both.
\erem

\bex
Suppose we have the von Neumann \nbd{\sB(G)}module $E=\sB(G,H)$. So, $\sB(G)'=\id_G\C$, $\rho'(\id_Gz)\colon x\odot g\mapsto x\odot gz$, $\sB^a(E)=\rho'(\id_G\C)'=\sB(H)$, and $E=C_{\id_G\C}(\sB(G,H))=\sB(G,H)$.
\eex

%%%% BO
% Here is a backwards example.
Here is an example at the opposite extreme. 
%%%% EO

\bex
Suppose $G:=L^2(\Om)$ and $\cB:=L^\infty(\Om)$, so that $\cB'=\cB$, and suppose $\rho'\colon\cB=\cB'\rightarrow\sB(H)$ is a normal unital representation. Then $E:=C_\cB(\sB(G,H))$ is a von Neumann \nbd{\cB}module. Moreover, since $\cB\cap\cB'=\cB\mapsto\id\odot\cB$ is a homomorphism onto the center of $\sB^a(E)=\rho'(\cB)'$ (faithful if and only of $\rho'$ is faithful), we see that $\sB^a(E)'=\rho'(\cB)=\sB^a(E)\cap\sB^a(E)'\subset\sB^a(E)$. Therefore, $\sB^a(E)$ is a type I von Neumann algebra and, up to (algebraic) isomorphism, every type I von Neumann algebra may be obtained that way.
\eex
%%%% BO new
% I had a question above but I figured it out, the commutant of a commutative algebra is type one .... (stupid me, I was probably working too tired). 
%%%% EO

\bemp[Von Neumann correspondences \it{vs} commuting pairs of representations.~] If we have a von Neumann \nbd{\cA}\nbd{\cB}correspondence $E$, then apart from the commutant lifting, which captures everything regarding the von Neumann \nbd{\cB}module structure of $E$, we have the Stinespring representation $\rho\colon\cA\rightarrow\sB^a(E)\rightarrow\sB(E\odot G)$ of $\cA$. Since $\rho(\cA)\subset\sB^a(E)=\rho'(\cB')'$, the ranges of $\rho'$ and $\rho$ mutually commute. Clearly, $\sB^{a,bil}(E)=\rho'(\cB')'\cap\rho(\cA)'$.

In general, by a \phantomsection\hl{commuting pair of representations}\index{commuting pair of representations}\index{von Neumann!correspondence!commutant lifting}\index{commutant lifting!of a von Neumann correspondence} $(\rho',\rho,H)$, we mean that $\rho'\colon\cB'\rightarrow\sB(H)$ and $\rho\colon\cA\rightarrow\sB(H)$ are normal unital representations with mutually commuting ranges. In the same way as representations $(\rho',H)$ correspond to (concrete) von Neumann \nbd{\cB}modules $E$ (and the morphisms $\sB^a(E_1,E_2)=\sB^{bil}({_{\cB'}}H_1,{_{\cB'}}H_2)$), the commuting pairs $(\rho',\rho,H)$ correspond to (\hl{concrete}) von Neumann \nbd{\cA}\nbd{\cB}correspondences $E$ (and the morphisms $\sB^{a,bil}(E_1,E_2)=\sB^{bil}({_{\cB'}}H_1,{_{\cB'}}H_2)\cap\sB^{bil}({_{\cA}}H_1,{_{\cA}}H_2)$); see again \cite{Ske06b}.

The functor called \hl{commutant}\index{von Neumann!correspondence!commutant lifting}\index{commutant!of a concrete von Neumann correspondence}, is best understood in the concrete categories. So, if we have a commuting pair $(\rho',\rho,H)$ we may define the von Neumann \nbd{\cA}\nbd{\cB}correspondence $E:=C_{\cB'}(\sB(G,H))$, that is, considering $\rho'$ as commutant lifting and turning the corresponding von Neumann \nbd{\cB}module into a correspondence via the Stinespring representation $\rho$. But, we also may look at the pair as $(\rho,\rho',H)$, considering $\rho$ as commutant lifting (recall that also $\cA$ is a von Neumann algebra, so $\cA=\cA''\subset\sB(K)$) resulting into a von Neumann \nbd{\cA'}module $E':=C_\cA(\sB(K,H))$ turned into a von Neumann \nbd{\cB'}\nbd{\cA'}correspondence via the Stinespring representation $\rho'$. Clearly, since also the morphisms are identical, we have a one-to-one functor between the category of concrete von Neumann \nbd{\cA}\nbd{\cB}correspondences and the category of concrete von Neumann \nbd{\cB'}\nbd{\cA'}correspondences, generalizing (clearly!) the commutant of von Neumann algebras (considering $\cB\subset\sB(G)$ the trivial correspondence over itself).
\eemp

\bex \label{vNBex}
For $\cA'=\cA=\C\subset\sB(\C)$, we get that the commutant of the von Neumann \nbd{\cB}module $E=C_{\cB'}(\sB(G,H))$ considered as \nbd{\C}\nbd{\cB}correspondence, is $E'=C_\C(\sB(\C,H))={_{\cB'}}H$.
\eex

\bex
If $E=\sB(G,G\otimes\eH)$ is a von Neumann \nbd{\sB(G)}correspondence (see Subsection \ref{vNAPP}\ref{vNBGmod}), then
\beqn{
E'
~=~
C_{\sB(G)}(\sB(G,G\otimes\eH))
~=~
\id_G\otimes\eH
~\cong~
\eH
}\eeqn
is the Hilbert space $\eH$.
\eex

\bex \label{AScex}
Let $T\colon\sB(K)\supset\cA\rightarrow\cB\subset\sB(G)$ be a normal CP-map and denote by $\sE$ the strong closure of its GNS-correspondence. Then $\sE':=C_{\cA}(\sB(K,\sE\odot G))$ is what is called the \hl{Arveson Stinespring correspondence}\index{correspondence!Arveson Stinespring}\index{Arveson Stinespring correspondence} of $T$ in Muhly and Solel \cite{MuSo02}. The map $\xi^*\in\sE^*=C_{\cB'}(\sB(\sE\odot G,G))$ gives rise to a map $\eta'\colon\sE'\rightarrow\sB(K,G)$ via $\eta'(x')=\xi^*x'$. Since $T$ is contractive, $\eta$ is completely contractive, it fulfills
\beqn{
\eta'(b'x'a')
~=~
b'\eta'(x')a',
}\eeqn
and it is \nbd{\sigma}weakly continuous. If $\cA=\cB$, then $\eta'(b'_1x'b'_2)=b'_1\eta'(x')b'_2$.

In general, a linear map $\eta$ from a \nbd{\cB}correspondence $E$ into $\sB(H)$ is called a \phantomsection\hl{covariant representation}\index{correspondence!covariant representation of} if (it is \nbd{\sigma}weakly continuous and if) there exists a nondegenerate (normal) representation $\pi$ such that $\eta(b_1xb_2)=\pi(b_1)\eta(x)\pi(b_2)$. One may show that every completely contractive \nbd{\sigma}weakly continuous covariant representation $\eta$ has the form $\eta(x)=\Xi'^*(x\odot\id_G)$ where $\Xi'$ is a unique element of $C_\cB(\sB(H,E\odot H))$. Obviously, $T'=\AB{\Xi',\bullet\Xi'}$ defines a CP-map on $\pi(\cB)'$, and doing to that CP-map as before, gives back $\eta$.

We see, a normal CP-map on $\cB\subset\sB(G)$ may be captured either by a von Neumann \nbd{\cB}correspondence $E$ with a vector $\xi\in E$ or by a von Neumann \nbd{\cB'}correspondence $E'$ with completely contractive \nbd{\sigma}weakly continuous covariant representation on $G$ where $\pi'=\id_{\cB'}$ is the defining representation of $\cB'$.
\eex

\lf
The following observation is necessary knowledge for working with commutants of von Neumann correspondences. For comparing structures that involve tensor products, to be discussed in the sequel, the observation is absolutely crucial.

\bob \label{noncob}
In the concrete categories (that is, when $E=C_{\cB'}(\sB(G,H))$ and when $E'=C_\cA(\sB(K,H))$), as explained in Footnote \ref{StineFN}, we have that $E\odot G$ is \it{canonically} isomorphic to $H$ via $x\odot g\mapsto xg$ and that $E'\odot K$ is \it{canonically} isomorphic to $H$ via $x'\odot k\mapsto x'k$. Of course, this means
\beqn{
E\odot G
~\cong~
E'\odot K
}\eeqn
``in some way''. But in writing $x\odot g$``$=$''$\sum x'_i\odot k_i$ (meaning that $xg=\sum x'_ik_i\in H$) or $x'\odot k$``$=$''$\sum x_i\odot g_i$ (meaning that $x'k=\sum x_ig_i\in H$), there is nothing canonical. It just means that there is the possibility of writing any element of $H$ as $h=\sum x_ig_i=\sum x'_ik_i$ and that the same element $h\in H$ corresponds to $\sum x_i\odot g_i\in E\odot G$ and to $\sum x'_i\odot k_i\in E'\odot K$. The more important it is to notice that these identifications, indeed, intertwine the correct actions of $\cB'$ and $\cA$ on $E\odot G\cong H\cong E'\odot K$. For $b'\in\cB'$ we find
\bmun{
\textstyle
(\id_E\odot b')\bfam{\sum x_i\odot g_i}
~=~
\sum x_i\odot b'g_i
~~\text{``=''}~~
\sum x_ib'g_i
~=~
\rho'(b')\bfam{\sum x_ig_i}
\\
\textstyle
~=~
\rho'(b')\bfam{\sum x'_ik_i}
~~\text{``=''}~~
(b'\odot\id_K)\bfam{\sum x'_i\odot k_i},
}\emun
and likewise for the action of $a\in\cA$.
\eob

\bemp[The tensor product under commutant.~]
Let $\cA\subset\sB(K)$, $\cB\subset\sB(G)$, $\cC\subset\sB(L)$ be von Neumann algebras, and let $E$ and $F$ be a von Neumann \nbd{\cA}\nbd{\cB}correspondence and a von Neumann \nbd{\cB}\nbd{\cC}correspondence, respectively. Then, in order to construct their von Neumann tensor product, we have to construct their \nbd{C^*}tensor product $E\odot F$, to embed this as $(E\odot F)\odot\id_L$ into $\sB(L,(E\odot F)\odot L)$, and to take the strong closure $E\sodots F:=\ol{(E\odot F)\odot\id_L}^s$ in $\sB(L,(E\odot F)\odot L)$. Instead of the last step, we also may compute
\beqn{
E\sodots F
~=~
C_{\cC'}(\sB(L,(E\odot F)\odot L)),
}\eeqn
where the Hilbert space $(E\odot F)\odot L$ carries the commutant lifting $\sigma'\colon c'\mapsto\id_{E\odot F}\odot\,c'$ of $\cC'$ and the Stinespring representation $\sigma\colon a\mapsto a\odot\id_{F\odot L}$. Actually, we have
\beqn{
(E\odot F)\odot L
~=~
E\odot F\odot L
~=~
E\odot(F\odot L),
}\eeqn
as in any tensor category, and we used that already in the definition of $\sigma$. We have
\beqn{
F\odot L
~~\text{``=''}~~
F'\odot G,
}\eeqn
as explained in Observation \ref{noncob}. We have a \bf{canonical}(!) isomorphism
\beqn{
E\odot(F'\odot G)
~\cong~
F'\odot(E\odot G)
}\eeqn
via the flip $x\odot(y'\odot g)\mapsto y'\odot(x\odot g)$. Plugging in $ax$ for $x$ and $c'y'$ for $y'$, we see that this flip intertwines the canonical actions of $\cA$ and of $\cC'$ on these spaces. (This is actually true for all $\sB^a(E)$ and all $\sB^a(F')$, respectively, which are sent to a commuting pair acting on $E\odot F'\odot G=F'\odot E\odot G$.) We have, again,
\beqn{
E\odot G
~~\text{``=''}~~
E'\odot K,
}\eeqn
as explained in Observation \ref{noncob}. And, again, we have
\beqn{
F'\odot(E'\odot K)
~=~
F'\odot E'\odot K
~=~
(F'\odot E')\odot K.
}\eeqn
Altogether, we have
\bmun{
(E\odot F)\odot L
~=~
E\odot(F\odot L)
~~\text{``=''}~~
E\odot(F'\odot G)
\\
~\cong~
F'\odot(E\odot G)
~~\text{``=''}~~
F'\odot(E'\odot K)
~=~
(F'\odot E')\odot K,
}\emun
including the transformation of the commutant lifting $\sigma'$ into the Stinespring representation and the transformation of the Stinespring representation $\sigma$ into the commutant lifting. We see:
\beqn{
(E\sodots F)'
~~\text{``=''}~~
F'\sodots E',
}\eeqn
that is, the tensor product is anti-multiplicative under commutant.

A once-for-all computation shows that in $(D\sodots E\sodots F)'~~\text{``=''}~~F'\sodots E'\odot D'$ it does not matter the commutant of which tensor product we compute first; under the identifications from Observation \ref{noncob} (basically going back to those in Footnote \ref{StineFN}), it is as associative as the tensor product in a tensor category. The result is a \hl{rule of thumb}: For von Neumann \nbd{\cB_i}\nbd{\cB_{i-1}}correspondences $E_i$ $(i=1,\ldots,n)$ over von Neumann algebras $\cB_i\subset\sB(G_i)$, compute
\bmun{
E_n\odot\ldots\odot E_2\odot E_1\odot G_0
~=~
E_n\odot\ldots\odot E_2\odot E'_1\odot G_1
~=~
E'_1\odot E_n\odot\ldots\odot E_2\odot G_1
\\
~=~
E'_1\odot E'_2\odot E_n\odot\ldots\odot E_3\odot G_2
~=~
\ldots
~=~
E'_1\odot E'_2\odot\ldots\odot E'_n\odot G_n,
}\emun
and it does not matter if we insert a computation of $\Bfam{(E_n\sodots\ldots\sodots E_{m+1})\sodots(E_m\sodots\ldots\sodots E_1)}'=(E_m\sodots\ldots\sodots E_1)'\sodots(E_n\sodots\ldots\sodots E_{m+1})'$ before.
\eemp

\bemp[Product systems under commutant.~] \label{PScomm}
Recall that (as always, in the sense of Observation \ref{noncob}), $\sB^{a,bil}(E)=\sB^{a,bil}(E')$. Therefore, if we have isomorphisms $u_{s,t}\colon E_s\sodots E_t\rightarrow E_{st}$ this gives rise to isomorphisms $u'_{t,s}:=u_{s,t}\colon E'_t\sodots E'_s\rightarrow E'_{st}$. A moment's thought shows that if $E^{\sodots}$ is a product system over $\bS$, then the $u'_{t,s}$ turn the $E'_t$ into a product system over $\bS^{op}$.

Likewise the commutant of the GNS-subproduct system (of von Neumann correspondences) of a normal CP-semigroup $T$ on $\cB\subset\sB(G)$ over $\bS^{op}$ is the subproduct system of Arveson-Stinespring correspondences in \cite{ShaSo09}. Since inductive limits are stable under commutant, the inductive limit in \cite{BhSk00} embedding the GNS-subproduct system into a product system, translates into the inductive limit from \cite{MuSo02} doing the same for the Arveson-Stinespring subproduct system. The covariant representations $\eta'_t$ of $E'_t$ on $\sB(G)$ constructed from the unit $\xi^\odot$ for $E^\odot$, fulfill $\eta'_{st}(x'_ty'_s)=\eta'_t(x'_t)\eta'_s(y'_s)$. (Recall that the $st$ in $\eta'_{st}$ is the product of $s$ and $t$ in $\bS$, hence, the product of $t$ and $s$ in $\bS^{op}$.)
\eemp

\bemp[Left dilations under commutant.~]  \label{LdilRdil}
Let the $v_t\colon E\sodots E_t\rightarrow E$ form a left dilation of the product system $E^{\sodots}$ to the (strongly full!) von Neumann \nbd{\cB}module $E$. Recall from Example \ref{vNBex} that $E$ has the commutant $E'={_{\cB'}}(E\odot G)=:{_{\cB'}}H$. The $v'_t:=v_t\in\sB^{bil}({_{\cB'}}(E'_t\odot H),{_{\cB'}}H)$ form what is called a \hl{right dilation} of the product system $E'^{\sodots}$ (over $\bS^{op}$!) to the (faithful!) von Neumann \nbd{\cB'}\nbd{\C}correspondence $H$, that is, the $v'_t$ are unitary and with the \hl{product} $x'_th:=v'_t(x'_t\odot h)$, we have $x'_t(y'_sh)=(x'_ty'_s)h$.

The right dilation gives rise to an \nbd{E_0}semigroup $v'_t(\id'_t\odot\bullet)v'^*_t$ on $\sB^{bil}(H)=\rho'(\cB')'=\sB^a(E)$ which coincides with the $\vt$ constructed from the left dilation $v_t$. It is noteworthy (and has already been described in \cite{Ske03c}) that the product system $E'^{\sodots}$ may be recovered as
\beqn{
E'_t
~=~
C_{\sB^a(E)}(\sB(H,{_{\vt_t}}H))
}\eeqn
(taking into account that this is a von Neumann correspondence over $\rho'(\cB')$ and that $\rho'$ is faithful) with product system structure $x'_ty'_s$ really given by operator multiplication in $\sB(H)$.
\eemp

\bemp[Dilations under commutant.~] \label{dilcomm}
Via $\wh{\eta}'_t(x'_t)h=x'_th$, the concept of right dilation is equivalent to the concept of \hl{faithful nondegenerate normal representation} of the product system $E'^{\sodots}$, namely, a family of injective maps $\wh{\eta}'_t\colon E'_t\rightarrow\sB(H)$ (not $\sB^{bil}(H)$!) that fulfills  $\cls\wh{\eta}'_t(E'_t)H=H$ (nondegeneracy), $\wh{\eta}'_{st}(x'_ty'_s)=\wh{\eta}'_t(x'_t)\wh{\eta}'_s(y'_s)$ and $\wh{\eta}'_t(x'_t)^*\wh{\eta}'_t(y'_t)=\rho'(\AB{x'_t,y'_t})$. (By the latter property, $\wh{\eta}'_t$ is necessarily completely contractive, and normal really means there is a normal extension of a \nbd{\sigma}weak map to the linking algebra.) It is noteworthy that $\wh{\eta}'_t$ together with the projection $p=\xi\xi^*\odot\id_G)$ is what is called an \hl{isometric fully coisometric dilation} of the completely contractive covariant representation $\eta'_t$ in \cite{MuSo02}. In fact, the approach to constructing dilations of CP-semigroups based on \cite{MuSo02}, is by, first, (where possible) embedding the Arveson-Stinespring subproduct system into a product system and, then, find an isometric dilation of the covariant representation with which the product system comes shipped. Namely: $p\wh{\eta}'_tp=\eta'_t$. In very much the same way, translating our proof of Theorem \ref{vNMarkmodthm} (that is, the von Neumann version of Theorem \ref{Markmodthm}) under commutant (and restricting it to submonoids of $\R^d$), we recover the proof of the same result as obtained earlier in \cite[Corollary 5.10 and Theorem 5.12]{ShaSo09}.

We see that finding a(n endomorphic) dilation of $T$ is the same as finding a(n isometric) dilation of the family $\eta'_t$. Note, however, that for the former we need not know what a product system is, while the latter without product systems cannot even be formulated. (Closely related is the notion of \hl{elementary dilation}\index{dilation!elementary}\index{elementary!dilation} of $T$ \cite{Ske11a,Ske16}, namely, an embedding $\vp\colon\cB\rightarrow\cA$ and a semigroup $\bfam{c_t}$ in $\cA$ such that $\vp\circ T_t(b)=c_t^*\vp(b)c_t$, mentioned already in Footnote \ref{elemFN}.)
\eemp

\vspace{1ex}
We stop our brief account on commutants and how this relates our results (based on \cite{BhSk00}) to other result (based on \cite{MuSo02}). The commutant of product systems \ref{PScomm} has been introduced in \cite{Ske03c} and discussed further in \cite{Ske09}. The latter also discusses the relation between left dilations and right dilations \ref{LdilRdil} and the relation of the latter with representation of product systems from \ref{dilcomm}. The relation with (isometric) dilations of covariant representation in \ref{dilcomm}, is explained carefully in \cite{Ske08}. Among other things addressed in  \cite{Ske05a,Ske08a}, these surveys also discuss carefully the Hilbert space case, that is, \nbd{E_0}semigroups $\vt$ on $\sB(H)$. In particular, it is explained that the Arveson system of $\vt$ (intertwiners) is the commutant of the Bhat system of $\vt$ (with a unit vector $\xi$, inputting $p=\xi\xi^*$ in Theorem \ref{E-supPSthm}) and that, therefore, they are anti-isomorphic and (by Tsirelson \cite{Tsi00p1}) the two need not be isomorphic.%
\footnote{Connes \cite{Con80p} looked at Hilbert spaces $H$ that are \nbd{\cA}\nbd{\cB}bimodules over \nbd{W^*}algebras $\cA$ and $\cB$. \phantomsection\index{von Neumann!correspondence!\it{versus} Connes correspondence}\index{Connes correspondence!\it{versus} von Neumann correspondence}In other words, we have a commuting pair $(\rho^\circ,\rho,H)$ of (unital normal) representations of the opposite \nbd{W^*}algebra of $\cB$, $\cB^\circ$, and of $\cA$. Of course, if we represent $\cA$ and $\cB$ as von Neumann algebras, then we are exactly in the representation picture of von Neumann \nbd{\cA}\nbd{{\cB^\circ}'}correspondences. In particular, if we opt to have $\cB$ (hence, $\cB^\circ$) in \it{standard representation}, then by \it{Tomita conjugation} we get $\cB^\circ\cong\cB'$. (Already mentioned in \cite[Remark 4.3]{Ske06b}.) Following the presentation in Takesaki \cite[Section IX.3]{Tak03a}, by how the action of the operator space $L:=L_\psi(\cD(H,\psi))$ is defined in\cite[Equation (11)]{Tak03a} and by  \cite[Lemma 3.3(iii)]{Tak03a},  $L$  is (a strongly closed submodule of) $C_{N^\circ}(\sB(L^2(N),H))$. By how the inner product is defined in \cite[Equation (23)]{Tak03a}, the \it{fusion product} of two Connes correspondences is exactly the commuting pair picture of the tensor product of the corresponding von Neumann correspondences. (Mentioned after \cite[Remark 6.3]{Ske08a} and recovered in \cite{BMSS12}. See also Thom \cite{Tho11}; see \cite{Ske12p} for a thorough discussion.)

The commuting pair picture of von Neumann correspondences is, therefore, a generalization of Connes correspondences and the tensor product of the former includes the \it{fusion product} of the latter. Unlike the latter, the former does not depend on tools like \it{Tomita-Takesaki theory} and \it{standard representation}, and the formulae do not depend manifestly on the choice of a faithful semi-finite normal weight.
}

\newpage

\section{Automatic adjointability of coisometries} \label{coisoAPP}

\blem
Let $E$ be a Hilbert \nbd{\cB}module and let $E_0$ be a closed submodule of $E$. If $a\in\sB^r(E)$ is a contractive idempotent with $aE=E_0$, then $a$ is a projection (that is, apart from being idempotent, $a$ is adjointable and $a^*=a$).
\elem

\proof
Let us first prove the statement for Hilbert spaces. In that case there exists a projection $p$ onto $E_0$. We write $a=p+b$, where $b:=a-p=a(\id_E-p)$. So, $ab=b=b(\id_E-p)$. For all $x\in(\id_E-p)E$ and $\lambda\in\C$, we get $a(x+\lambda bx)=a(\id_E-p)x+\lambda bx=(1+\lambda)bx$. Taking also into account that $x=(\id_E-p)x$ and $bx=pbx$ are orthogonal, we get
\beqn{
\norm{a}^2
~\ge~
\frac{(\lambda+1)^2\norm{bx}^2}{\lambda^2\norm{bx}^2+\norm{x}^2}
}\eeqn
whenever $x\ne0$. Assume $b\ne0$, so that there is $x=(\id_E-p)x$ such that $bx\ne0$ (and, \it{a fortiori}, $x\ne0$). We get
\beqn{
\norm{a}^2
~\ge~
\frac{(\lambda+1)^2}{\lambda^2+\frac{\norm{x}^2}{\norm{bx}^2}}
~=~
\frac{\lambda^2+2\lambda+1}{\lambda^2+\frac{\norm{x}^2}{\norm{bx}^2}}
~>~
1
}\eeqn
whenever $2\lambda+1>\frac{\norm{x}^2}{\norm{bx}^2}$. In other words, if $b\ne0$ then $a$ is not a contraction. Since $a$ \bf{is} a contraction, $b=0$ and $a=p=a^*$.

Now suppose we are in the general situation of the lemma. Choose a faithful representation of $\cB$ on a Hilbert space $G$ and put $H:=E\odot G$. Then the operator $a\odot\id_G$ is a contractive idempotent on $H$. Therefore, by the first part, it is self-adjoint. It follows
\beqn{
\AB{g,\AB{ax,x'}g'}
~=~
\AB{(a\odot\id_G)(x\odot g),x'\odot g'}
~=~
\AB{x\odot g,(a\odot\id_G)(x'\odot g')}
~=~
\AB{g,\AB{x,ax'}g'}
}\eeqn
for all $x,x'\in E$ and for all $g,g'\in G$, so, $\AB{ax,x'}=\AB{x,ax'}$ for all $x,x'\in E$.\qed

\bcor
Let $E$ and $F$ denote Hilbert \nbd{\cB}modules and let $E_0$ be a closed submodule of $E$. Then every contraction $w\in\sB^r(E,F)$ that restricts to a unitary $u\colon E_0\rightarrow F$ is adjointable.
\ecor

\proof
Denote by $i\colon E_0\rightarrow E$ the canonical injection and put $a:=iu^*w\colon E\rightarrow F\rightarrow E_0\rightarrow E$. Then $a$ is as in the lemma and, therefore, is the projection onto $E_0$. In particular, $i$ is adjointable and $i^*i=\id_{E_0}$. So, $w=ui^*a\colon E\rightarrow E\rightarrow E_0\rightarrow F$ is adjointable.\qed

\bob
Of course, $E_0$ is necessarily closed. It is also unique. (Indeed, $E_0$ is the complement of $\ker a=\ker w$.) By the corollary, a coisometry (see Remark \ref{adrem}) is always adjointable.
\eob

\brem
Adjointability of coisometries could also be proved, appealing to the statement that $\ker a^{\perp\perp}=\ker a$ for every $a\in\sB^r(E)$ (Frank \cite[Lemma 2.4]{FraM02}) -- if this statement was true. Unfortunately, as Kaad and Skeide \cite{KaSk21p} show, such statement is false. Therefore, here we produce a different proof.
\erem

% \bex \label{nonsurex}

% \eex

% \brem
% \OW[later (OPEN)]{A short comment here about Frank \cite{FraM02}? Imagine if the question also reoccurs in that intersection business for tensor products of \nbd{C^*}correspondences!}
% \erem

\newpage

\section{Intersections in tensor products} \label{laAPP}

Is the \hl{intersection property}\index{tensor product!intersection property|bf} $(E_1\odot F_1)\cap(E_2\odot F_2)=(E_1\cap E_2)\odot(F_1\cap F_2)$ true? If yes, this would guarantee that the intersection of (super)(sub)product subsystems of (super)(sub)product systems is again a (super)(sub)product subsystem. The inclusion of the R.H.S.\ in the L.H.S.\ is obvious, guaranteeing that the intersection of superproduct systems is indeed again a superproduct system. It is also clear that showing the intersection property for two, by well-ordering the index set and transfinite induction we get the intersection property for arbitrary families of subsystems. So let us concentrate on the intersection of two.

It is \it{folklore}, but probably not too well-known, that the intersection property is true for the tensor product of vector spaces. Let us, therefore, discuss this case first.

Let $V$ and $W$ be vector spaces with subspaces $V_0$ and $W_0$, respectively. $V_0\otimes W_0\ni x\otimes y\mapsto x\otimes y\in V\otimes W$ identifies $V_0\otimes W_0$ with the subspace of $V\otimes W$ spanned by simple tensors from $V_0$ and $W_0$. The injection is well-defined by the universal property of the tensor product. Injectivity is seen by extending a basis $B_0$ of $V_0$ to a basis $B$ of $V$ and a basis $C_0$ of $W_0$ to a basis of $W$, and by the fact that the simple tensors of basis vectors form a basis of the tensor product.

\bob \label{incob}
So far, this remains true for the tensor product of correspondences, with or without completion or strong closure, just because the embedding defined in the same way is, clearly, isometric.
\eob

Recall that for any basis $B$ of a vector space $V$, We have the \hl{dual basis} $\vp_b$ $(b\in B)$ of linear functionals on $V$ defined by $\vp_bb'=\delta_{b,b'}$.

\blem \label{iplem}
\hfill
$
(V_1\cap V_2)\otimes(W_1\cap W_2)
~=~
(V_1\otimes W_1)\cap(V_2\otimes W_2).
$
\hfill\hfill
{~}\elem

\proof
Only the inclusion $\supset$ requires thought. Choose bases $B_0$ and $C_0$ for $V_0:=V_1\cap V_2$ and $W_0:=W_1\cap W_2$, respectively. Extend the bases $B_0$ and $C_0$ to bases $B_i$ of $V_i$ and $C_i$ of $W_i$, respectively. Since both sides are, clearly, contained in $(\ls V_1\cup V_2)\otimes(\ls W_1\cup W_2)$, it is sufficient to show the statement when, $\ls V_1\cup V_2=V$ and $\ls W_1\cup W_2=W$. It follows that the three sets $B_i\backslash B_0$ and $B_0$ ($C_i\backslash C_0$ and $C_0$) are disjoint and that their union $B$ ($C$) is a basis for $V$ ($W$). Suppose $u=\sum_{b\in B,c\in C}b\otimes c\lambda_{b,c}$ is in $(V_1\otimes W_1)\cap(V_2\otimes W_2)$. Since $u\in V_1\otimes W$ we get
%%%% BO 
$\lambda_{b,c}=(\vp_b\otimes\vp_c)u=0$ whenever, $b\notin B_1$; likewise for $b\notin B_2$, so, $b\in B_0$ whenever $\lambda_{b,c} \neq 0$. Likewise, $c\in C_0$ for every such $\lambda_{b,c}$.
%%%% EO
So, $u\in\ls B_0\otimes C_0=(V_1\cap V_2)\otimes(W_1\cap W_2)$.\qed

\lf
It would be tempting to imitate that proof for Hilbert spaces replacing bases with ONBs. But in general it is not possible to to choose the sets $B_0,B_1,B_2$ such that $B_1\backslash B_0$ and $B_2\backslash B_0$ are mutually orthogonal. And for von Neumann modules the simple tensors of elements from QONBs are still complete, yes, but they need no longer be orthogonal.
% \OW[OPEN: Just a question of general interest. To be erased next round]{Does every tensor product of vN-correspondences admit a QONB consisting of product vectors?}

One also might assume that the intersection property for vector spaces turns over to the algebraic tensor product over an algebra, because the latter is a quotient of the former. But, the following example shows that quotienting and intersection need not be compatible.

\bex
Let $V=\C^3$ and define the subspaces $V_1=\C\rtMatrix{1\\1\\0}$, $V_2=\C\rtMatrix{0\\1\\1}$, and $\sN=\rtMatrix{\C\\0\\\C}$. Then $V_0:=V_1\cap V_2=\zero$. Define any seminorm on $V$ that has kernel $\sN$, and define $W:=V/\sN=\C$ and $W_i:=V_i + \sN = \CB{x+\sN\colon x\in V_i}\subset W$ $(i=0,1,2)$. Then $W_0=\zero$, but  $W_1=W_2=\C$, so, $W_1\cap W_2\ne W_0$.
%Let $V=\C^3$ and define the subspaces $V_1=\C\rtMatrix{1\\1\\0}$, $V_2=\C\rtMatrix{0\\1\\1}$, and $\sN=\rtMatrix{\C\\0\\\C}$. Then $V_0:=V_1\cap V_2=\zero$. Define any seminorm on $V$ that has kernel $\sN$, and define $W:=V/\sN=\C$ and $W_i:=V_i/(\sN\cap V_i)\subset W$ $(i=0,1,2)$.%
%\footnote{
%Note that $W_i$ coincides with the subset $\CB{x+\sN\colon x\in V_i}$ of $W$ via the canonical bijection $x+(\sN\cap V_i)\mapsto x+\sN$.
%}
%Then $W_0=\zero$, but  $W_1=W_2=\C$, so, $W_1\cap W_2\ne W_0$.
\eex

We do not know if an example with similar behaviour can be obtained for the (algebraic or topological) tensor products of correspondences. But, Observation \ref{vNintob} tells us that the intersection property holds for tensor products of von Neumann correspondences, and finite-dimensional correspondences as their (algebraic or topological) tensor products are von Neumann correspondences. So if a counter example for \nbd{C^*}correspondences can be obtained, it cannot be obtained that easily.

\bprop \label{minpqprop}
Let $V$ and $W$ be Hilbert spaces and suppose that the occurring subspaces are closed. Then Lemma \ref{iplem} remains true replacing the algebraic tensor product with that of Hilbert spaces.
\eprop

\proof
Denote by $p_i$ the projections onto the subspaces $V_i$. It is well known that the projection $p_0$ onto the intersection $V_0:=V_1\cap V_2$ can be obtained as strong limit
\beqn{
p_0
~=~
\lim_{n\to\infty}p_1(p_2p_1)^n.
}\eeqn
(Indeed, $p_1(p_2p_1)^n\ge p_1(p_2p_1)^{n+1}\ge0$, so the limit $p_0$ exists. Clearly, $p_0$ is a projection fulfilling $p_1p_0=p_0$, hence, $p_0\le p_1$. Clearly, $p_0(\U-p_2)p_0=p_0-p_0=0$, so $p_2p_0=p_0$, hence, $p_0\le p_2$. Clearly, $qp_0=q$ for each projection $q\le p_i$ $(i=1,2)$, so $q\le p_0$. Hence, $p_0=p_1\wedge p_2$.)
The same holds for the analogue projection $q_i$ onto $W_i$. Clearly, the projection onto $V_i\otimes W_i$ is $p_i\otimes q_i$. The projection onto $(V_1\otimes W_1)\cap(V_2\otimes W_2)$ can be obtained as
\beqn{
\lim_{n\to\infty}((p_1\otimes q_1))((p_2\otimes q_2)(p_1\otimes q_1))^n
~=~
\lim_{n\to\infty}(p_1(p_2p_1)^n)\otimes(q_1(q_2q_1)^n)
~=~
p_0\otimes q_0,
}\eeqn
so it coincides with the projection onto $V_0\otimes W_0$.\qed

\bob \label{vNintob}
Since also strongly closed submodules of a von Neumann module $E\subset\sB(G,H)$ are characterized by projections in $\sB(H)$, and since the last computation in the proof does not depend on that the projections labeled by $1$ and $2$ are in tensor position, but holds for arbitrary mutually commuting tuples of projections $(p_i)$ and $(q_i)$, the proof goes through also for the tensor product of von Neumann correspondences.
\eob

\newpage

\section{The time ordered Fock module} \label{FockAPP}

The\index{module!Fock|see {Fock, module}} \hl{symmetric Fock space}\index{Fock!space!symmetric} $\G(H)$ over a Hilbert space $H$ may be defined as the subspace of the \hl{full Fock space}\index{Fock!space!full} $\cF(H):=\bigoplus_{n\in\N_0}H^{\otimes n}$ ($H^{\otimes 0}=\om\C$ with the \hl{vacuum} unit vector $\om$) spanned by the symmetric tensors $x\otimes\ldots\otimes x$ and $\om$. See Skeide \cite[Proposition 8.1.4]{Ske01} for why this coincides with the usual definition where to each summand $H^{\otimes n}$ 
%%%% BO 
we apply 
%%%% EO
the symmetrization projection
\beqn{
p_n
\colon
x_1\otimes\ldots\otimes x_n
~\longmapsto~
\frac{1}{n!}\sum_{\sigma\in S_n}x_{\sigma(1)}\otimes\ldots\otimes x_{\sigma(n)}.
}\eeqn
It is well known that $\G(H\oplus G)\cong\G(H)\otimes\G(G)$ in a \it{canonical fashion} (in the first place, this means in an \it{associative way}). However, already writing down \it{the usual way} to do that
\beqn{
\sqrt{n!}p_n(x_1\otimes\ldots\otimes x_n)\otimes\sqrt{m!}p_m(y_1\otimes\ldots\otimes y_m)
~\longmapsto~
\sqrt{(n+m)!}p_{n+m}(x_1\otimes\ldots\otimes x_n\otimes y_1\otimes\ldots\otimes y_m)
}\eeqn
(see \cite[Theorem 8.1.6]{Ske01}), shows that doing the \it{natural} thing, is quite a relative notion. Effectively, it is better to consider that $\G(H)$ is (easily) spanned by the \phantomsection\hl{exponential vectors}\index{exponential vector}
\beqn{
\ee(x)
~:=~
\sum_{n\in\N_0}\frac{x^{\otimes n}}{\sqrt{n!}}
}\eeqn
($x\in H$). Verifying that $\AB{\ee(x),\ee(x')}=e^{\AB{x,x'}}$, shows that $\bfam{\G(H),i\colon x\mapsto\ee(x)}$ is the \it{Kolmogorov decomposition} for the \it{positive definite kernel} $(x,x')\mapsto e^{\AB{x,x'}}$ on the set $H$. It is easy to verify that the isomorphism $\G(H\oplus G)\rightarrow\G(H)\otimes\G(G)$ above with the nontrivial combinatorial factors, is nothing but $\ee(x)\otimes\ee(y)\mapsto\ee(x\oplus y)$.

One might be tempted to try the same for the \phantomsection\hl{full Fock module}\index{Fock!module!full} $\cF(E):=\bigoplus_{n\in\N_0}E^{\odot n}$ over a \nbd{\cB}correspondence $E$ ($\cB\ni\U$), where $E^{\odot 0}=\om\cB$ with the central \hl{vacuum} unit vector $\om$ (isomorphic to $\cB$ via $\om\mapsto\U$). However, (without further structure, like $E$ being \it{centered} (Skeide \cite{Ske98}) or generated by other elements fulfilling `\it{good}' commutation relations (Accardi and Skeide \cite{AcSk00a}; see \cite[Chapter 8]{Ske01}) the possibilities end with the analogue of the definition of $\G(H)$ and the (still spanning!) exponential vectors. In fact, the projections $p_n$ may be ill-defined (sending $0=x_1\odot x_2$ to $\frac{1}{2}x_2\odot x_1\ne0$), and $\ee(x)\odot\ee(y)\mapsto\ee(x\oplus y)$ is rarely an isometry (for instance, by the same $x=x_1$ and $y=x_2$ as in the preceding bracket). In general:  There is no `\it{good}' symmetric Fock module (yet). However, in quantum probability and quantum dynamics, the `\it{interesting}' symmetric Fock spaces are all of the type $\G(L^2(\Om,K))$ and the direct sum decompositions of $L^2(\Om,K)$, usually, emerge from partitioning the measure space $\Om$ into two (or more) parts. Here, the partial commutativity of $L^2(\Om,K)=L^2(\Om)\otimes K$ (regarding the function part) and presence of further structure on $\Om$ (like the order structure of $\Om=\R$), enable us to improve the situation quite a bit. The result is the \it{time ordered Fock module} as introduced in Bhat and Skeide \cite[Section 11]{BhSk00}. (In the case of Hilbert spaces it is known as the \it{Guichardet picture} \cite{Gui72} of the symmetric Fock space. Other names around are \it{chronological} or \it{time consecutive} Fock space.) Today we would probably call it \it{monotone Fock module} in reference to the outstanding role it plays in \it{monotone independence} introduced by Lu \cite{Lu97} and Muraki \cite{MurN97} in the scalar case, while Skeide \cite{Ske04} and M.\ Popa \cite{MPop08} provided the \it{amalgamated} version.

\lf
Recall that for each measure space $\Om$ and any correspondence $F$, the space $L^2(\Om,F)$ is defined as \phantomsection\hl{external tensor product}\index{tensor product!external} $L^2(\Om)\otimes F$, the completion of the pre-\nbd{C^*}correspondence $L^2(\Om)\,\ul{\otimes}\,F$ (over the same \nbd{C^*}algebras as $F$) with its obvious structure. (See \cite[Section 4.3]{Ske01} for details.) It is \bf{not} correct that the elements of $L^2(\Om,F)$ are representable as (equivalence classes of) \nbd{F}valued functions; see \cite[Example 4.3.13]{Ske01}. So it is not justified to think of them as functions. (On the other hand, while for Hilbert spaces it is true that every element of $L^2(\Om)\otimes K$ can be represented as a measurable, square-integrable \nbd{K}valued function, we do not know a reference that would guarantee that we get (up to measure $0$) all measurable, square-integrable \nbd{K}valued functions, which do form a Hilbert space containing $L^2(\Om)\otimes K$, when $K$ is nonseparable.) However it is safe define operators in terms of operations on functions $\Om\rightarrow F$, when restricting to $f\otimes y\colon\alpha\mapsto f(\alpha)y$ for $f$ from a subset of ``nice'' functions total in $L^2(\Om)$ and $y$ from a total subset of $F$. For any measurable subset $S\subset\Om$, point-wise multiplication with the indicator function $\I_S$ defines a projection in $\sB(L^2(\Om))$ (also denoted by $\I_S$). This projection extends (by amplification with $\id_F$) to a projection in $\sB^a(L^2(\Om,F))$ (also denoted by $\I_S$). Note that, clearly, $\sK(L^2(\Om,F))$ is $\sK(L^2(\Om))\otimes\sK(F)$ so that $\sB^a(L^2(\Om,F))$ is the \nbd{*}strong or strict completion of $\sK(L^2(\Om))\otimes\sK(F)$. (There is no really better characterization of this multiplier algebra of a (spatial) \nbd{C^*}tensor product.) If $F_1$ and $F_2$ are correspondences such that $F_1\odot F_2$ makes sense and if the measure spaces $\Om_1$ and $\Om_2$ are such that the Fubini theorem holds for $\Om_1\times\Om_2$ (for instance if both are \nbd{\sigma}finite), then
\bmun{
L^2(\Om_1,F_1)\odot L^2(\Om_2,F_2)
~=~
(L^2(\Om_1)\otimes L^2(\Om_2))\otimes(F_1\odot F_2)
\\
~=~
L^2(\Om_1\times \Om_2)\otimes(F_1\odot F_2)
~=~
L^2(\Om_1\times\Om_2,F_1\odot F_2).
}\emun
(Note that the hypothesis about the Fubini theorem, necessary in the step from the first to the second line, is missing in \cite[Observation 4.3.11]{Ske01}.)

\lf
For each $n\in\N$, denote the indicator function of $\CB{\bbm{a}=(\alpha_n,\ldots,\alpha_1)\in\R^n\colon 0\le\alpha_1<\ldots<\alpha_n}$ by $\Delta_n$. Let $F$ be a correspondence over a unital \nbd{C^*}algebra $\cB$. For each measurable $S\subset\R$, we denote by
\beqn{
\DG(S;F)
~:=~
\bigoplus_{n\in\N_0}\Delta_nL^2(S^n,F^{\odot n})
~\subset~
\sF(L^2(S,F))
~\subset~
\sF(L^2(\R,F))
}\eeqn
(with $\Delta_0=\id_{\om\cB}$) the \phantomsection\hl{time ordered Fock module}\index{Fock!module!time ordered} over $F$ relative to $S$. In particular, we put
\baln{
\DG(F)
&
~:=~
\DG(\R_+;F),
&
\DG_t(F)
&
~:=~
\DG(\RO{0,t};F).
}\ealn
For every $x\in L^2(S,F)$ ($\subset L^2(\R,F)$), we define the \phantomsection\hl{exponential vector}\index{exponential vector|bf}
\beqn{
\ee(x)
~:=~
\sum_{n\in\N_0}\Delta_nx^{\odot n}
~\in~
\DG(S;F)
~~~
(\subset~
\DG(\R;F)
)
}\eeqn
whenever the right-hand side exists. Caution is in place with the following notation to not confuse elements of $L^2(S,\F)$ with elements of $F$. For each $y\in F$ and $t\in\R_+$ we put
\beqn{
\ee_t(y)
~:=~
\ee(\I_{\RO{0,t}}y)
~\in~
\DG_t(F).
}\eeqn
Note that all $\ee_t(y)$ exist. In fact, taking into account that $\int_{\RO{0,t}^n}\Delta_n(\bbm{a})\,d\bbm{a}=\frac{t^n}{n!}$, we get%
\footnote{ \label{expFN}
Note that by $e^{t\AB{y,\bullet y}}$ we rather mean $e_\circ^{t\AB{y,\bullet y}}$, the exponential of the map $t\AB{y,\bullet y}$ in $\sB(\cB)$ with composition $\circ$ as product. We do not mean the map $b\mapsto e^{t\AB{y,by}}$.
}
\beqn{
\AB{\ee_t(y),\bullet\ee_t(y)}
~=~
e^{t\AB{y,\bullet y}}.
}\eeqn
Recall the following \it{folklore}: The set of indicator functions of intervals $\RO{a_n,b_n}\times\ldots\times\RO{a_1,b_1}$ to mutually disjoint intervals $\RO{ a_i,b_i}$, are total in $L^2(\R^n)$. (This is usually referred to as ``the diagonal has measure $0$'', which hits the point for $n=2$. A more accurate formulation would be ``any hyperplane in $\R^n$ has measure $0$''.) Applied to the time ordered sector, this means: The set of indicator functions of intervals
\beqn{
\RO{a_n,b_n}\times\ldots\times\RO{a_1,b_1}
\text{~~~with~~~}
b_n>a_n\ge b_{n-1}>a_{n-1}\ge\ldots\ge b_1>a_1,
}\eeqn
are total in $\Delta_nL^2(\R^n)$. With this, we easily check that the maps that send ``functions'' $G_m\in\Delta_mL^2(\RO{0,s}^m,F^{\odot m})\subset\DG_s(F)$ and $H_n\in\Delta_nL^2(\RO{0,t}^n,F^{\odot n})\subset\DG_t(F)$ to the ``function''
\beqn{
(\beta_m,\ldots,\beta_1,\alpha_n,\ldots,\alpha_1)
~\longmapsto~
G_m(\beta_m-t,\ldots,\beta_1-t)\odot H_n(\alpha_n,\ldots,\alpha_1)
}\eeqn
in $\Delta_{m+n}L^2(\RO{0,s+t}^{m+n},F^{\odot m})\subset\DG_{s+t}(F)$ defines a bilinear unitary $u_{s,t}\colon\DG_s(F)\odot\DG_t(F)\rightarrow\DG_{s+t}(F)$. Clearly, the $u_{s,t}$ turn the family $\DG^\odot(F)=\bfam{\DG_t(F)}_{t\in\R_+}$ into a product system, the \phantomsection\hl{time ordered product system}\index{product system!time orered|bf} over $F$. Replacing $G_m$ with a ``function'' in $\Delta_mL^2(\R_+^m,F^{\odot m})\subset\DG(F)$, we get unitaries (even bilinear) $v_t\colon\DG(F)\odot\DG_t(F)\rightarrow\DG(F)$, that form a left dilation of $\DG^\odot(F)$ to $\DG(F)$. The resulting strict \nbd{E_0}semigroup on $\sB^a(\DG(F))$ is called the (\hl{generalized}) \hl{CCR-flow}\index{CCR-flow@(generalized) CCR-flow} over $F$.

Note that for $t=t_n>\ldots>t_1>t_0=0$ we get
\beqn{
\ee_{t_n-t_{n-1}}(y_n)\ldots\ee_{t_1-t_0}(y_1)
~=~
\ee\bfam{\I_{\RO{t_0,t_1}}y_1+\ldots+\I_{\RO{t_{n-1},t_n}}y_n}
~\in~
\DG_t(F)
~\subset~
\DG(F).
}\eeqn
We see that we get all exponential vectors to step functions as products of $\ee_t(y)$ and, clearly, the exponential vectors to step functions are total. Actually, one may show that it is enough to vary $y$ in a total subset of $F$ having $0$ as accumulation point; see \cite{Ske00a,BBLS04} or \cite[Theorem 7.4.3]{Ske01}.

Clearly $\om^\odot=\bfam{\om_t}_{t\in\R_+}$ with $\om_t=\om\in\DG_t(F)\subset\DG(F)$, is a unital unit. Taking $\DG^\odot(F)$ and $\om^\odot$ as input for Theorem \ref{Oreindthm}, one easily verifies that the inductive limit is $\DG(F)$ (with $k_t\colon\DG_t(F)\rightarrow\DG(F)$ being the canonical injection $\DG(F)\supset\DG_t(F)\rightarrow\DG(F)$) and the left dilation from Theorem \ref{Oreindthm} coincides with $v_t$ given above. So, the \nbd{E_0}semigroup from Theorem \ref{Oreindthm} is just the CCR-flow over $F$. (The fact that the unit $\om^\odot$ is central so that the inductive system is over bimodule maps, explains better why the inductive limit is a correspondence in its own right. For the same reason, the CCR-flow is also what is called a \hl{noise}; see Skeide \cite{Ske06d} or the sections on \it{spatial} Markov semigroups in \cite{Ske16}.)

Clearly, the families $\ee^\odot(y)=\bfam{\ee_t(y)}_{t\in\R_+}$ form units -- as we know already, enough units to generate the whole product system $\DG^\odot(F)$ when, as already said, varying $y$ over a total subset $S\subset F$ with $0\in\ol{S}$. One can show that the units that are continuous when considered as function into $\DG_t(F)\subset\DG(F)$ can be parametrized one-to-one by pairs $(\beta,\zeta)\in\cB\times F$ as $\xi^\odot(\beta,\zeta)=\bfam{\xi_t(\beta,\zeta)}_{t\in\R_+}$ where the component of $\xi_t(\beta,\zeta)$ in the \nbd{n}particle sector of $\DG_t(F)$ is given by the ``function''
\beqn{
(\alpha_n,\ldots,\alpha_1)
~\longmapsto~
e^{(t-\alpha_n)\beta}\zeta e^{(\alpha_n-\alpha_{n-1})\beta}\odot\ldots\odot\zeta e^{(\alpha_2-\alpha_1)\beta}\odot\zeta e^{\alpha_1\beta}
}\eeqn
for $n\in\N$ and $\om e^{t\beta}$ for $n=0$. (Note that $\ee^\odot(y)=\xi^\odot(0,y)$ and, in particular, $\om^\odot=\xi^\odot(0,0)$.) This is the main result of Liebscher and Skeide \cite{LiSk01}. It is of outstanding importance in the classification theory of \hl{type I} product systems (those generated by their continuous units) to know that the explicit form of the generators of the semigroups $\AB{\xi_t(\beta,\zeta),\bullet\xi_t(\beta',\zeta')}$ of maps on $\cB$ is
\beqn{
\eL^{(\beta,\zeta),(\beta',\zeta')}
\colon
b
~\longmapsto~
\AB{\zeta,b\zeta'}+\beta^*b+b\beta',
}\eeqn
a so-called \phantomsection\hl{Christensen-Evans form}\index{Christensen-Evans form}. We do not need this in full generality, but only the case $(\beta',\zeta')=(\beta,\zeta)$ giving the generator of the CP-semigroup $T_t:=\AB{\xi_t(\beta,\zeta),\bullet\xi_t(\beta,\zeta)}$ as $L(b)=\AB{\zeta,b\zeta}+\beta^*b+b\beta$. (Christensen and Evans \cite{ChrEv79} showed that in the von Neumann case every generator has this form.)

\newpage

\section[Partially order preserving permutations and order improving transpositions]{Partially order preserving permutations and\\order improving transpositions} \label{popAPP}

The main scope of this appendix is to prove Lemma \ref{pi_flem}, but also to give a better feeling about this ``permutation stuff'' in Section \ref{compSEC} without distracting from that section's main scope. (Note that both properties in the title are not properties of a permutation \it{per se}, but properties of a permutation relative to a function $f$.) See the paragraph in front of Definition \ref{admisdef} for basic notation and facts about the symmetric group. We interpret the following lemma in the observation that separates it from its proof; but, see also the discussion following Lemma \ref{pi_flem}.

\blem\label{piplem}
Fix $n\in\N$ and $d\in\N$. Let $f$ denote a function $\N_n\rightarrow\N_d$. Then there exists a unique permutation $\sigma\in S_n$ fulfilling the following conditions for all $j<i$:\phantomsection\index{partially order preserving|bf}\index{permutation!partially order preserving|bf}\index{permutation!nondecreasing}
\baln{
f\circ\sigma(j)
&
~\le~
f\circ\sigma(i)
&&
\text{($\sigma$ is \hl{nondecreasing} for $f$)},
\\
f\circ\sigma(j)=f\circ\sigma(i)
~&
\Longrightarrow~
\sigma(j)<\sigma(i)
&&
\text{($\sigma$ is \hl{partially order preserving} for $f$)}.
}\ealn
\elem

\bob \label{tuplemob}
$f$ defines an \nbd{n}tuple $\bbf:=(f(1),\ldots,f(n))$ with entries in $\N_d$. We let act a permutation $\sigma\in S_n$ on that tuple by sending it to the tuple $\bbf\circ\sigma:=(f\circ\sigma(1),\ldots,f\circ\sigma(n))$. Let us see what the permutation of 
%%%% BO
% the Lemma 
the lemma
%%%% EO
does to the tuple $\bbf$, also explaining our terminology.

The first condition is clear; it simply means that the tuple $\bbf\circ\sigma$ is nondecreasing. The second condition means that whenever two positions $j<i$ in the new tuple have the same value, then the new order $j<i$ corresponds to the old order $\sigma(j)<\sigma(i)$ in which the sites were before applying the permutation. That is, the order of all entries in $f^{-1}(k)$ where $f$ assumes the same value $k$ is not changed by our permutation. That is, what we mean by a permutation $\sigma$ being partially order preserving for $f$.

We could have formulated all this more intuitively, by looking at the natural representation of $S_n$ on $(\N_d)^n\ni\bbf$ (faithful, if $d\ge2$), where $\bbf\mapsto\bbf\circ\sigma$ is actually the representing map for $\sigma^{-1}$. But the proof runs more smoothly with our choice and our choice is actually also what is needed in the applications.
\eob

\proof[Proof of Lemma \ref{piplem}.~]
If $n=1$, then there are no pairs with $j<i$. So the only permutation $e\in S_1$ is the unique permutation fulfilling all conditions. We shall proceed by induction. So, we assume the statement is true for a certain $n\in\N$ and and we choose a function $f\colon\N_{n+1}\rightarrow\N_d$.

Recall that $S_n$ can be identified with the subgroup $\CB{\sigma\in S_{n+1}\colon\sigma(n+1)=n+1}$ of $S_{n+1}$, consisting of all bijections on $\N_{n+1}$ that leave $n+1$ fixed. Every permutation $\sigma\in S_{n+1}$ can be written uniquely as $\sigma=\tau\circ\sigma'$ where $\sigma'\in S_n\subset S_{n+1}$ and where $\tau$ is either the identity (if $\sigma\in S_n$) or the flip exchanging $n+1$ and $\sigma(n+1)$ (if $\sigma\notin S_n$). Moreover, if we define $f':=(f\circ\tau)\upharpoonright\N_n$, then
\beq{\label{n+1->n}
f\circ\sigma(j)
~=~
\begin{cases}
f'\circ\sigma'(j)& j\le n,
\\
f\circ\sigma(n+1)&j=n+1.
\end{cases}
}\eeq
Of course, $f'\colon\N_n\rightarrow\N_d$ fulfilling this equation, is determined uniquely by $f$ and $\sigma$.

Now, $\sigma$ fulfills the conditions in the lemma relative to $f$ if and only if $\sigma'$ fulfills the conditions relative to $f'$ and $\sigma(n+1)=\max f^{-1}(\max f)$ (the biggest argument on which $f$ assumes its maximum). Indeed, if $\sigma$ fulfills the conditions, then from
\beqn{
f\circ\sigma(j)
~\le~
f\circ\sigma(n+1)
\text{~~~~~~and~~~~~~}
f\circ\sigma(j)=f\circ\sigma(n+1)
~\Longrightarrow~
\sigma(j)<\sigma(n+1)
}\eeqn
we conclude that  necessarily $\sigma(n+1)=\max f^{-1}(\max f)$. So, assuming that this necessary condition holds, the conditions in the lemma are true for $j<i=n+1$. For arguments $j\le n$ we have $f\circ\sigma=f'\circ\sigma'$. So, for $j<i\le n$ the conditions hold for $f\circ\sigma$ if and only if they hold for $f'\circ\sigma'$.

We are now ready to make the inductive step.

\item
Uniqueness: Suppose $\sigma\in S_{n+1}$ is a permutation fulfilling the conditions. Then $\sigma(n+1)=\max f^{-1}(\max f)$ is uniquely determined by $f$ and $\sigma'$ fulfills the conditions relative to $f'$. Now, $\sigma(n+1)$ is determined uniquely by $f$, so $\tau$ (depending only on $\sigma(n+1)$) is determined uniquely by $f$, so $f'$ is determined uniquely by $f$. Since, by hypothesis, $f'$ determines $\sigma'$, also $\sigma=\tau\circ\sigma'$ is determined uniquely by $f$.

\item
Existence: Put $\sigma(n+1):=\max f^{-1}(\max f)$. As before, define $\tau$ (depending on $\sigma(n+1)=n+1$ or $\sigma(n+1)\le n$) and put $f':=(f\circ\tau)\upharpoonright\N_n$. By the inductive assumption, there exists a (unique) permutation $\sigma'\in S_n$ fulfilling the conditions for $f'$. Then by \eqref{n+1->n} we see that $\sigma:=\tau\circ\sigma'$ fulfills the conditions for $f$.\qed

\lf
Following Observation \ref{tuplemob}, the permutation of this lemma is somewhat the permutation we want in Lemma \ref{pi_flem}. But its characterization in Lemma \ref{pi_flem}, which allowed us to work, is a different one. Now, let us see how to obtain the characterization in Lemma \ref{pi_flem}.

The permutation from Lemma \ref{pi_flem} is the unique partially order preserving permutation for $f$ that puts the tuple $\bbf$ into nondecreasing order. Without that condition, a permutation that is partially order preserving for $f$, may very well increase disorder. In fact, a transposition $\tau_\vk$ is partially order preserving for $f$ if and only if $f(\vk)\ne f(\vk+1)$; this may be $f(\vk)>f(\vk+1)$ (what we want), or $f(\vk)<f(\vk+1)$ (what we do not want). To sort that out we define the following:
\bdefi
Denote by $(j\lcom i)$ a pair $(j,i)$ with $j<i$. The \hl{number of inversions}\index{number of inversions!of a function $\N_n\rightarrow\N_d$} of $f$ is
\beqn{
\inv(f)
~:=~
\#\bCB{(j\lcom i)\colon f(j)>f(i)}.
}\eeqn
A permutation $\sigma$ is \phantomsection\hl{order improving}\index{order improving|bf}\index{permutation!order improving|bf} for $f$ if $\,\inv(f)-\inv(f\circ\sigma)>0$.
\edefi

\lf%\noindent
For fixed $1\le\vk<d$, writing in $\inv(f)$ all pairs $(j\lcom i)$ where $i$ or $j$ equals $\vk$ or $\vk+1$ separately, we get
\baln{
\inv(f)
~=~
&
\#\Bfam{\,\bCB{(j\lcom i)\colon f(j)>f(i)}\,\cap\,\bfam{\N_d\backslash\CB{\vk,\vk+1}}^2\,}
\\
&
\,+\,
\#\bCB{(j\lcom\vk)\colon f(j)>f(\vk)}
\,+\,
\#\bCB{(j\lcom\vk)\colon f(j)>f(\vk+1)}
\\
&
\,+\,
\#\bCB{(\vk+1 \lcom i)\colon f(\vk)>f(i)}
\,+\,
\#\bCB{(\vk+1\lcom i)\colon f(\vk+1)>f(i)}
\\
&
~~~~~~~~~~~\,+\,
\begin{cases}
\hfill0&f(\vk)\le f(\vk+1),
\\
1&f(\vk)>f(\vk+1).
\end{cases}
}\ealn
Since the first row remains unchanged under the transposition $\tau_\vk$, since the summands in the second row just interchange under the transposition $\tau_\vk$, and since the summands in the third row just interchange under the transposition $\tau_\vk$, the change of from $\inv(f)$ to $\inv(f\circ\tau_\vk)$ is determined just by the last summand. We see, a transposition $\tau_\vk$ is order improving for $f$ if and only if $f(\vk)>f(\vk+1)$, and an order improving transposition is partially order preserving automatically.

Following the terminology of Definition \ref{admisdef} (but without attaching, for the time being, operators to transpositions $\tau_\vk$), we say a chain $\tau_{\vk_1},\ldots,\tau_{\vk_m}$ of transpositions is \phantomsection\hl{admissible}\index{admissible!chain (of transpositions)} for $f$ if for every $k$ the transposition $\vk_k$ is order improving for $f\circ\tau_{\vk_1}\circ\ldots\circ\tau_{\vk_{k-1}}$.  Since in each step the inversion number is reduced by $1$, and since the inversion number is nonnegative, there exists chains $\tau_{\vk_1},\ldots,\tau_{\vk_m}$ that are \hl{maximal}\index{admissible!chain (of transpositions)!maximal} in the sense that no further $\tau_{\vk_{m+1}}$ can be found such that the chain $\tau_{\vk_1},\ldots,\tau_{\vk_{m+1}}$ still is admissible.

\bprop
Equivalent are:
\begin{enumerate}
\item \label{01}
$\inv(f)=0$.

\item \label{02}
$f$ is nondecreasing.

\item \label{03}
There exists no order improving $\tau_\vk$ for $f$.
\end{enumerate}
\eprop

\proof
Obviously, \ref{02} and \ref{03} are equivalent and \ref{01} implies either of them. The only thing that requires an argument is, why existence of $j<i$ with $f(j)>f(i)$ (that is, not \ref{01}) implies existence of $\vk$ with $f(\vk)>f(\vk+1)$ (that is, not \ref{03}). If $j+1=i$, we are done with $\vk=j$. Otherwise, either we have $f(j)>f(j+1)$ so that we are done with $\vk=j$, or we have $f(j+1)\ge f(j)>f(i)$ so that we are in the same situation with $j+1$ and $i$ instead of $j$ and $i$. After a finite number of reiterations we will find $\vk$ with $j\le\vk<i$ such that $f(\vk)>f(\vk+1)$.\qed

\lf
So, for each maximal chain, $f\circ\tau_{\vk_1}\circ\ldots\circ\tau_{\vk_m}$ is nondecreasing, that is, the permutation $\tau_{\vk_1}\circ\ldots\circ\tau_{\vk_m}$ is nondecreasing for $f$. Note, too, that if $\sigma'$ is partially order preserving for $f$ and if $\sigma''$ is partially order preserving for $f\circ\sigma'$, then $\sigma'\circ\sigma''$ is partially order preserving for $f$. Since each order improving $\tau_\vk$ is also partially order preserving, $\tau_{\vk_1}\circ\ldots\circ\tau_{\vk_m}$ is partially order preserving for $f$, too. In other words, $\tau_{\vk_1}\circ\ldots\circ\tau_{\vk_m}$ is the unique permutation $\sigma$ from Lemma \ref{piplem}; so, it does not depend on which maximal chain we chose.

We, thus, proved the first two parts of Lemma \ref{pi_flem}.

\brem
An interesting fact, which we do not need, is that every partially order preserving and order improving permutation $\sigma$ for $f$ may be decomposed into an admissible chain of transpositions for $f$. (One may reduce this to Lemma \ref{piplem} by replacing $f$ with a `finer' version (for bigger $d$) so that the permutation is the unique one for the new $f$.)
\erem

The statement of the preceding remark is, like everything we proved so far in this appendix, one of these statements that are intuitively clear, but where intuition is far away from a well formulated statement and where the proof is surprisingly tricky. We add now to these problems that we actually want to represent the transpositions by operators and that these operators really only make sense and can be composed only when attached to a partially order preserving transposition. Like in the proof of Theorem \ref{prodthm}, we appeal to the trick to define things on a bigger space. The bigger space has the advantage that it carries a representation of the full permutation group (whose theory is well understood). Then, we have to show that our particular partially order preserving permutations go between the subspaces we are actually interested in.

\blem\label{sTlem}
Let $E$ be a \nbd{\cB}correspondence and suppose that $\sF\colon E\odot E\rightarrow E\odot E$ is a self- inverse isomorphism of correspondences fulfilling
\beq{\label{sTcond}
(\sF\odot\id_E)(\id_E\odot\,\,\sF)(\sF\odot\id_E)
~=~
(\id_E\odot\,\,\sF)(\sF\odot\id_E)(\id_E\odot\,\,\sF).
}\eeq
For $n\ge2$ define $\sT_k:=\id_E^{\odot(k-1)}\odot\,\,\sF\odot\id^{\odot(n-k-1)}\colon E^{\odot n}\rightarrow E^{\odot n}$ $(k=1,\ldots,n-1)$. Then  setting $\pi(\tau_k):=\sT_k$ defines a representation $\pi$ of $S_n$ into the automorphism group of $E^{\odot n}$.
\elem

\proof
The $\sT_k$ fulfill the relations \eqref{Snrel}.\qed

\lf
Note that the condition on $\sF$ is also necessary for that $\pi$ defines a representation. We shall say an operator $\sF$ fulfilling the condition of the lemma is an \phantomsection\hl{exchange operator}\index{detailed exchange conditions}\index{operator!exchange}\index{exchange operator}\index{exchange operator!detailed exchange conditions}.

\bprop\label{E_nprop}
Let $E_1,\ldots,E_d$ denote \nbd{\cB}correspondences and suppose for all $1\le j<i\le d$ there is an isomorphism $\sF_{j,i}\colon E_i\odot E_j\rightarrow E_j\odot E_i$ of correspondences. Put $E:=E_1\oplus\ldots\oplus E_d$, and define $\sF_{i,i}:=\id_{E_i\odot E_i}$ and $\sF_{j,i}:=\sF_{i,j}^*$ for $j>i$. Then the operator $\sF\in\sB^{a,bil}(E\odot E)$ defined by
\beqn{
\sF\upharpoonright(E_i\odot E_j)
~:=~
\sF_{j,i}
}\eeqn
is a self-inverse isomorphism of correspondences.

$\sF$ is an exchange operator if and only if the $\sF_{j,i}$ fulfill the conditions in \eqref{Tijcond}, which we repeat here for convenience:
\beq{\label{Tijcond'}
(\id_k\odot \sF_{j,i})(\sF_{k,i}\odot\id_j)(\id_i\odot \sF_{k,j})
~=~
(\sF_{k,j}\odot\id_i)(\id_j\odot \sF_{k,i})(\sF_{j,i}\odot\id_k)
}\eeq
for all $1\le k<j<i\le d$.
\eprop

\proof
Of course, $\sF$ is left-linear, self-inverse and self-adjoint, so, unitary. In other words, $\sF$ is a self-inverse isomorphism. By applying \eqref{sTcond} to all subspaces $E_i\odot E_j\odot E_k$ $(i,j,k\in\N_d)$ individually, we see that it is necessary and sufficient, if \eqref{Tijcond'} holds for all $i,j,k\in\N_d$. So, the only thing that remains to be shown, is that validity of \eqref{Tijcond'} for $k<j<i$ implies validity for all $i,j,k$.

It is easy to check that \eqref{Tijcond'} holds, whenever at least two of the $i,j,k$ are equal without any condition on the $\sF_{i,j}$ $(i<j)$. So, we have to check \eqref{Tijcond'} for all $i\ne j\ne k\ne i$, using only validity for $k<j<i$. Apart from the case $k<j<i$, which is true by hypothesis, there are five more cases.\noqed\vspace{-2ex}

\begin{enumerate}
\item
$k<i<j$, so by \eqref{Tijcond'}
\beq{\label{T1}
(\id_k\odot \sF_{i,j})(\sF_{k,j}\odot\id_i)(\id_j\odot \sF_{k,i})
~=~
(\sF_{k,i}\odot\id_j)(\id_i\odot \sF_{k,j})(\sF_{i,j}\odot\id_k).
}\eeq
Multiplying from the left with $\id_k\odot \sF_{j,i}$ and from the right with $\sF_{j,i}\odot\id_k$, gives the result.

\item
$j<k<i$, so by \eqref{Tijcond'}
\beq{\label{T2}
(\id_j\odot \sF_{k,i})(\sF_{j,i}\odot\id_k)(\id_i\odot \sF_{j,k})
~=~
(\sF_{j,k}\odot\id_i)(\id_k\odot \sF_{j,i})(\sF_{k,i}\odot\id_j).
}\eeq
Multiplying from the left with $\sF_{k,j}\odot\id_i$ and from the right with $\id_i\odot \sF_{k,j}$, gives the result.

\item
$i<k<j$ follows from \eqref{T2} by inverting.

\item
$j<i<k$ follows from \eqref{T1} by inverting.

\item
$i<j<k$ follows from \eqref{Tijcond'} by inverting.\qedsymbol
\end{enumerate}

\bprop\label{detexprop}
For each permutation $\sigma\in S_n$ and all choices $k_1,\ldots,k_n\in\N_d$, the operator 
%%%% BO 
% $\pi(\sigma)$ 
$\pi(\sigma)$ (from Lemma \ref{sTlem} applied to $\sF$ from Proposition \ref{E_nprop}) 
%%%% EO
restricts to an isomorphism
\beqn{
E_{k_1}\odot\ldots\odot E_{k_n}
~\longrightarrow~
E_{k_{\sigma^{-1}(1)}}\odot\ldots\odot E_{k_{\sigma^{-1}(n)}}.
}\eeqn
\eprop

\proof
The formula behaves nicely under compositions $\pi(\sigma)\pi(\sigma')=\pi(\sigma\circ\sigma')$, so it is enough to understand it for the flips $\tau_k$. But this case is obvious.\qed

\lf
These two propositions prove also the third part of Lemma \ref{pi_flem}.

\newpage

\setcounter{section}{17}
\section[A road map]{A road map} \label{rmAPP}

Apart from an index in the end of these notes, also in this appendix we try to help the readers find their way through the jungle of definitions and results in the main text. While in the main text things are (with few exceptions) introduced where they are also discussed, here we try to organize the definitions in a more systematic and compact way. For instance, the reader will find together all sorts of dilation together in one or two points, instead of having their extra properties spread all over the whole work.

For things that are tight together in the main text, we will, to avoid too much mere repetition, references to the main text instead of giving a formulation here. Among these there are also many preliminaries. 

\begin{itemize}
\item
Recall from the conventions in the introduction (p.\pageref{IC}) our conventions how to write monoids $\bS$ (especially that we will write the neutral element \bf{always} as $0\in\bS$). And recall from there or from the part preceding Convention \ref{ucCONV} what a \it{semigroup over a monoid $\bS$ on a set $\cB$} is.

 \item
 Recall from the part preceding Convention \ref{ucCONV} what a CP-semigroup is, and from Convention \ref{ucCONV} that our CP-semigroups are assumed \it{contractive}. Our CP-semigroups act on unital \nbd{C^*}algebras $\cA,\cB,\cC,\ldots$, unless stated otherwise. (Important exceptions are $\sK(E)$ or the pre-\nbd{C^*}algebra $\sF(E)$ (p.\pageref{IC}) or $\cA_\infty$ (see ....).) Recall Definition \ref{spCPdefi} for special CP-semigroups (\nbd{E}, Markov, and \nbd{E_0}semigroups).
 
 \item
 Recall from the conventions in the introduction (p.\pageref{IC}) what a \it{CP-map} and its \it{GNS-construction} are. Recall from Footnote \ref{StineFN} the relation of the latter to the \it{Stinespring construction}. (Throughout these notes we reference many times to the way to transform an abstract Hilbert module into a concrete operator module explained in that footnote.)
 
 \item
 Recall from the conventions in the introduction (p.\pageref{IC}) basic preliminaries of the multiplier algebra of a \nbd{C^*}algebra and the notations like $\sB^a(E)$, $\sK(E)$, and $\sF(E)$. Digest in Proposition \ref{convprop} and its proof the multiplier algebra of a a pre-\nbd{C^*}algebra in general an in the special case $\sF(E)$.

\item
It is warmly recommended to \it{read} the unusual conventions in the introduction (p.\pageref{IUC}).
\end{itemize}
More basic notations and other things to not be forgotten throughout:

\begin{itemize}
\item
Additional structures of monoids are dealt with in Section \ref{PSmonoSEC}: \hl{right-reversible} (Definition \ref{rrefdefi}) $=$ \hl{directed} (Proposition \ref{r-i-dir-prop}) directed monoids; \hl{Ore} semigroups (Definition \ref{Oredefi}, their \hl{universal covering group}, and when they are partially ordered; (necessarily right-reversible) monoids that are \hl{totally directed} (Definition \ref{totdirdefi}; the \it{generalized interval partitions} $\bJ_t$ for a monoid $\bS$ following Lemma \ref{totdirlem}), which are partially ordered in general (Proposition \ref{Jpoprop}) and directed if $\bS$ is totally directed and cancellative (Theorem \ref{totdirthm}).

\item
Recall  that semigroups and their related product-system-like structures are indexed by \bf{opposite monoids}! See, again, Section \ref{PSmonoSEC}, especially the part that starts with Definition \ref{unitdef} and ends with Example \ref{incrpex}.

\end{itemize}

\lf
\subsection{\normalsize Dilations} \label{RdilSSEC}

These notes are entirely dedicated to so-called \it{weak dilations}. But, general dilations are also addressed very briefly -- firstly, to give a basis to the historic Remark \ref{dilhistrem} in which we explain what else is there and what the reader may not expect here, and, secondly, because we wish to underline that the \it{strongness condition} for a dilation is (unlike many other properties of dilation we use in these notes) not restricted to weak dilations.

\lf
A (general) dilation is always given as a quadruple $(\cA,\theta,i,\ep)$ consisting of a unital \nbd{C^*}al\-ge\-bra $\cA$, an \nbd{E}semigroup $\theta$ over $\bS$ on $\cA$, an embedding $i\colon\cB\rightarrow\cA$ of another unital \nbd{C^*}algebra $\cB$ into $\cA$, and an \hl{expectation} $\ep\colon\cA\rightarrow\cB$ (that is, the map $\E:=i\circ\ep$ is a conditional expectation onto $i(\cB)$).

It is important to not forget that $(\cA,\theta,i,\ep)$ can be a \it{dilation} in two ways:

\begin{itemize}
\item
Firstly, given a CP-semigroup $T=\bfam{T_t}_{t\in\bS}$ on $\cB$ (sometimes denoted as a pair $(\cB,T)$), the quadruple is a \hl{dilation of $T$} if Diagram \ref{basdil} commutes, that is, if $\ep\circ\theta_t\circ i=T_t$; see Definition \ref{dildef}.

\item
Secondly, given just $(\cA,\theta,i,\ep)$, the quadruple is a \hl{dilation} \it{per se} if the $T_t:=\ep\circ\theta_t\circ i$ defined by Diagram \ref{basdil} form a semigroup (necessarily CP); see Definition \ref{dildef}.
\end{itemize}
Of course, if $(\cA,\theta,i,\ep)$ is a dilation, then it is a dilation of the CP-semigroup formed by the maps $\ep\circ\theta_t\circ i$; but frequently it is convenient not to specify that CP-semigroup in advance.

\begin{itemize}
\item
The quadruple $(\cA,\theta,i,\ep)$ is \hl{strong} if $\ep\circ\theta_t=\ep\circ\theta_t\circ\E$, that is, if $\ep\circ\theta_t=T_t\circ\ep$. (See Diagram \eqref{strodil}.)

\item
If the quadruple is strong then it is a dilation, a \hl{strong} dilation. (Observation \ref{semgenob}; it is even sufficient to check the strongness condition only on a generating subset of $\bS$.)
\end{itemize}
Properties of general dilations (other than being strong):
\begin{itemize}
\item
A dilation is \hl{unital} if $i$ is unital, \hl{semireversible} it $\theta$ is an injective \nbd{E_0}semigroup, and \hl{reversible} if $\theta$ is an automorphism semigroup. (Definition \ref{dildef}; these adjectives occur only very occasionally. In any nontrivial situation, unital or reversible dilations cannot be also \it{weak}.)

\item
A dilation is \hl{Markov} if the dilated CP-semigroup is Markov. (Definition \ref{dildef}.) It is an \hl{\nbd{E_0}dilation} if $\theta$ is an \nbd{E_0}semigroup. (Observation \ref{E0strongob}.)

\item
A dilation is \hl{strict}, \hl{normal}, and so forth, if each map $\theta_t$ is strict (see the conventions p.\pageref{IC}), normal (see p.\pageref{normald} in Appendix \ref{vNAPP}), and so forth. A dilation  is \hl{strictly}, \hl{strongly}, and so forth, \hl{continuous} if the map $t\mapsto\theta_t$ is continuous in that topology; see the part of Section \ref{topSEC} that starts before Observation \ref{pstrcontob}. (Recall that we refer to continuity with $t$ only very occasionally in these notes.)
\end{itemize}
All terminology for general dilations is used (where applicable and not mutually excluding like, for instance, \it{unital dilation}) also for \it{weak dilations}. \it{Weak dilations} are \bf{the only dilations} considered in these notes. The basis for \it{weak dilations} (very near to classical dilation theory; see Example Section \ref{EXwnsSEC}) is the concept of \it{compression}:
\begin{itemize}
\item \phantomsection\label{pBT}
A projection $p\in\cB$ \hl{compresses} a CP-semigroup $(\cB,T)$ (to a CP-semigroup $(\cB^p,T^p)$) if the maps $T^p_t:=pT_t(p\bullet p)p$ on $\cB^p:=p\cB p$ form a semigroup $T^p$.

\item
A compression is \hl{strong} if $T^p_t(pbp)=pT_t(b)p$ for all $b\in\cB$. We say $p$ compresses $(\cB,T)$ strongly.
\end{itemize}
In this generality, compression is discussed only in the ultrashort Subsection \ref{EXBexSEC}\ref{CPcompSSEC} and the application in Example \ref{Bex}. (See also the parenthetical part of Observation \ref{compob}\eqref{cob1}. And if we would have recognized the trivial semigroup on $\rtMatrix{\C\\0}$ as compression of the trivial semigroup on $\C^2$, in Example \ref{hypexex} we would have avoided to have to show that the projection $p$ is increasing.) Everywhere else, we only compress \nbd{E}semigroups to CP-semigroups (\it{weak dilations}) or we compress dilations to other dilations (of the same CP-semigroup). We, therefore, rather meet pairs $(\cA,p)$ than pairs $(\cB,p)$ of a unital \nbd{C^*}algebra and a projection, and usually, $\cB=p\cA p$; this is one of the two settings in which \it{weak dilations} show up in these notes.
\begin{itemize}
\item
A dilation $(\cA,\theta,i,\ep)$ is \hl{weak} if $i(\cB)=p\cA p$, where we put $p:=i(\U_\cB)$, so that $\E=p\bullet p$. (Definition \ref{dildef}.)

This means that $p$ compresses the \nbd{E}semigroup $(\cA,\theta)$ to a CP-semigroup on $i(\cB)$, name\-ly, $(i(\cB)=p\cA p,i\circ T\circ i^{-1}=\theta^p)$, where $i^{-1}\colon i(\cB)\rightarrow\cB$ is just $\ep\upharpoonright i(\cB)$. 

Note that a weak dilation is strong if and only if the compression by $p$ is strong. We will, in these notes, not say a \it{strong weak dilation}, but just a \it{strong dilation}, as all our dilations are weak.

\item
By a \hl{weak dilation} $(\cA,\theta,p)$ we mean an \nbd{E}semigroup $(\cA,\theta)$ and a pair $(\cA,p)$ where $p$ is a projection in $\cA$ that compresses $(\cA,\theta)$ to a CP-semigroup on $\cB:=p\cA p\subset \cA$; in other words, $(\cA,\theta,i^{can},p\bullet p)$ (where $i^{can}$ is the canonical embedding of $\cB$ into $\cA$) is a weak dilation. (Convention \ref{subalgconv}.)

\item
By a \hl{module dilation} $(E,\vt,\xi)$ we mean a (weak with $p=\xi\xi^*$) dilation of the form $(\sB^a(E),\vt,\xi\bullet\xi^*,\AB{\xi,\bullet\xi})$, where $E$ is a Hilbert \nbd{\cB}module and $\xi$ is a \hl{unit vector} in $E$, that is, $\AB{\xi,\xi}=\U$. (Convention \ref{Ba(E)conv}.)
\end{itemize}
In these notes, there are no other dilations than weak dilations $(\cA,\theta,p)$ and module dilations $(E,\vt,\xi)$. Module dilations are a different way to look at a special class of weak dilations; this raises the question which weak dilations can be viewed as module dilation.

\begin{itemize}
\item
A pair $(\cA,p)$ is \hl{full} if the corner $p\cA p$ is \hl{strictly full} in $\cA$, that is, if $\cA$ is $M(\ls\cA p\cA)$, the \hl{multiplier algebra} of the pre-\nbd{C^*}algebra $\ls\cA p\cA$. (Definition \ref{dildef} and Proposition \ref{convprop}.)

(Noteworthy and recurrent notation explained in Proposition \ref{convprop}: The Hilbert \nbd{p\cA p}mod\-ule $E:=\cA p$ and its finite-rank operators $\sF(E)=\ls\cA p\cA$. Of course, the pair $(\sB^a(E),\xi\xi^*)$ is full.)

\item
A weak dilation $(\cA,\theta,p)$ is \hl{full} if $(\cA,p)$ is full. (Definition \ref{dildef}.)

\item
A weak dilation is full if and only if it is \it{conjugate} (under the translation explained in Proposition \ref{convprop}) to a module dilation.
\item
If $\cA$ is a von Neumann algebra, then each $(\cA,p)$ is full if and only if the \hl{central cover} $c(p)$ of $p$ (see the before Theorem \ref{vNcontthm}) is $\U$. (Theorem \ref{vNcontthm}\eqref{vNct2}.)
\end{itemize}
For completeness:
\begin{itemize}
\item
A pair $(\cA,p)$ is \hl{semifull} (Definition \ref{sfulldefi}) if $\cA$ \hl{contains} $\sB^a(E)$ (Definition \ref{contdefi}). A weak dilation is \hl{semifull} if $(\cA,p)$ is semifull

\item
$(\cA,p)$ is semifull if and only if $\cA=\sB^a(E)\oplus\sF(E)^\perp$. (Theorem \ref{contthm}.) Not all dilations are semifull. (Example \ref{nsfex}.)

\item
If $\cA$ is a von Neumann algebra, then each $(\cA,p)$ is semifull. (Theorem \ref{vNcontthm}\eqref{vNct1}.)
\end{itemize}
There are many examples for dilations that are not full. Typically, they are obtained by \it{algebraic minimalization} (see below and Subsection \ref{minSEC}\ref{algminSSEC}) of a full dilation. (This is natural, because, after all, most explicit constructions of dilations give full dilations.)

\begin{itemize}
\item
The algebraic minimalizations of the (full) dilations constructed in Example Section \ref{EXpropsupSEC} give dilations acting on a commutative algebra $\cA$ and, therefore, are not full. (See the beginning of Subsection \ref{minSEC}\ref{fmsSSEC}.)

\item
\it{Incompressible} (see below and Subsection \ref{minSEC}\ref{compSSEC}) full dilations that are not algebraically minimal have algebraic minimalizations that are not full. (Corollary \ref{incmnfcor}.) Instances are: The minimalization of Bhat's Example \ref{Bex} (one-parameter non-Markov) and the dilation discussed in Observation \ref{amiMdilob} (the two-parameter Markov dilation as constructed throughout the whole Example Section \ref{EXN02SEC}).
\end{itemize}
Fullness is a property of the pair $(\cA,p)$ alone. Now we \bf{fix} a triple $(\cA,\theta,p)$ and see which properties of $p$ inside this triple can make the triple a dilation and  which stronger properties $p$ can encode. We also consider \it{good dilations} (see Subsection \ref{spsSSEC} below).

\begin{itemize}
\item
Proposition \ref{Mcharprop}: Equivalent are
\begin{itemize}
\item
$p$ is \hl{increasing} for $\theta$, that is, $\theta_t(p)\ge p$ (or $(\U-\theta_t(p))p=0$);

\item
the maps $T_t$ are unital;

\item
the triple is a weak Markov dilation.
\end{itemize}
Moreover, each Markov dilation is strong. (Equation \ref{strdileq}.)

\item
Proposition \ref{strequivprop}: Equivalent are% \OW{This stronger than \ref{strequivprop}, but what we prove!}
\begin{itemize}
\item
$\theta_t(\U-p)p=0$;

\item
$p$ is \hl{coincreasing} (that is, $\theta_t(\U-p)\le\U-p$);

\item
the triple is a strong dilation.
\end{itemize}
Dilations with coincreasing $p$ are called \hl{regular} in \cite{Bha03}.

\item
Observation \ref{puniob}: Equivalent are
\begin{itemize}
\item
$\theta_{ts}(p)\theta_t(\U-p)p=0$; (attention: the semigroups in Observation \ref{puniob} are over $\bS^{op}$;)

\item
the triple is a good dilation. (See Subsection \ref{spsSSEC} below.)
\end{itemize}
\end{itemize}
Of course, Markov $\Rightarrow$ strong $\Rightarrow$ good. Important is to note that all statement in these three groups that do not say ``the triple is a dilation'' do not assume that the triple is a dilation; therefore, these statements provide, in particular, methods to show that a triple is a dilation. (See also Remark \ref{sgrem}.)

We list some more properties and examples. \phantomsection%\OW{Check if ex-list complete}
\begin{itemize}
\item
A strong \nbd{E_0}dilation is Markov. (Observation \ref{E0strongob}.) \nbd{E_0}dilations that are not strong need not be Markov. (Example \ref{E0weakex}.)

\item
A dilation that is not good. (Bhat's Example \ref{Bex}; see the whole Example Section \ref{EXBexSEC}.)

\item
Dilations that are good, but not strong. (Everything in the Example Section \ref{EXwnsSEC} on \it{elementary} CP-semigroups (see below) that is not strong, is at least good; Example \ref{pPSelemex}. (These are discrete one-parameter examples, but Example \ref{cwnsex} explains how they may be promoted to continuous time one-parameter examples.) An \nbd{E_0}dilation of a proper \nbd{E}semigroup; Example \ref{E0weakex}.)
\end{itemize}
We have just seen three (counter) examples or classes of such that are based on \it{elementary} CP-semigroups. See the (not necessarily complete) list following Example \ref{cwnsex} for more instances.

\begin{itemize}
\item
A CP-semigroup $T$ \bf{on} $\cB$ is \hl{elementary} if it has the form $T_t=c_t^*\bullet c_t$ for a (contraction) semigroup $c=\bfam{c_t}$ \bf{in} $\cB$.

\item
A weak dilation $(\cA,\theta,p)$ (necessarily of an elementary CP-semigroup) is \hl{solidly elementary} if $\theta$ has the form $\theta_t=w_t^*\bullet w_t$ for a (necessarily coisometric) semigroup $w$ in $\cA$ and the $c_t:=pw_tp$ form a semigroup. (That is, $w$ is a coisometric dilation of the contraction semigroup $c$. See the paragraphs in between Proposition \ref{strsolprop} and Remark \ref{FSzNrem}.)

\item
Proposition \ref{strsolprop} (reformulation): ~A strong dilation to an elementary \nbd{E}semigroup is solidly elementary. Example \ref{nonsolex}: Not all dilations to an elementary \nbd{E}semigroup are solidly elementary.

\item
Note that the definitions of elementary and solidly elementary fix the semigroup in $\cB$ that represents the semigroup on $\cB$ by conjugation. We did not address question about what happens when fix only the semigroup on $\cB$ but not the semigroup in $\cB$ that represents it. See, however, the discussion in Appendix \ref{vNAPP}\ref{vNBGmod} where uniqueness questions for the \hl{Kraus decomposition} of a single CP-map on $\sB(G)$ are addressed, and Subsection \ref{EXBexSEC}\ref{KrowSSEC} were these are applied to dilations of \hl{column contractions}.

\end{itemize}

\lf
\subsection{\normalsize ``Standing hypotheses''}

\begin{itemize}
\item
The \hl{standing hypotheses} many results in Section \ref{minSEC} refer to, mean a \it{Markov dilation over the opposite of an Ore monoid} as they also occur in Section \ref{leftdilSEC}. (In Section \ref{leftdilSEC}, they allow to construct an inductive limit (of right modules) for a \it{superproduct system} over a \it{unital unit} (Theorem \ref{E-lsemdilthm}; see also Subsection \ref{spsSSEC} below); in Section \ref{minSEC}, they make sure that the family of projections $\theta_t(p)$ is, indeed, an (increasing) net.)

\item
Another recurrent set of hypotheses we \bf{could} have called \it{standing}, is: \it{Let $\bS$ be a totally directed monoid (or such that the $\bJ_t$ are directed).} (See Theorem \ref{totdirthm}.) They enable the possibility to turn a subproduct system into a product system (Theorem \ref{indlimthm}; see also Subsection \ref{spsSSEC} below).

Only in Subsection \ref{minSEC}\ref{1-p-SSEC} occur the recurrent hypotheses \it{over the opposite of a totally directed cancellative monoid}, which stand for how far we can go with what holds in the usual one-parameter case. They unite the two preceding hypotheses on the monoid (recall that a totally directed monoid is right reversible), but they do occur without the hypothesis \it{Markov}; in fact, the generalization to not necessarily unital CP-semigroups is one of the main points in \ref{minSEC}\ref{1-p-SSEC}.
\end{itemize}

\noindent

\lf
\subsection{\normalsize Minimalities}

A triple $(\cA,\theta,p)$, apart from giving the \nbd{E}semigroup $\theta$, determines the subalgebra $\cB:=p\cA p$ of $\cA$. The smallest \nbd{C^*}subalgebra of $\cA$ containing $\cB$ and being invariant under $\theta$, is
\beqn{
\cA_\infty
~:=~
C^*(\theta_\bS(\cB)).
}\eeqn
This algebra is (almost always) nonunital. The triple is \hl{algebraically minimal} if $\cA$ is generated in the one or the other topology by $\cA_\infty$. In Subsection \ref{minSEC}\ref{algminSSEC} we discuss three variants:
\begin{itemize}
\item
If $\cA$ is a von Neumann algebra, then $(\cA,\theta,p)$ is \hl{normally algebraically minimal} if $\cA=\ol{\cA_\infty}^s$. (Definition \ref{namindef}.)

\item
If $(\cA,p)$ is full, then $(\cA,\theta,p)$ is \hl{strictly algebraically minimal} if $\cA=\ol{\cA_\infty}^{stri}$. (Definition \ref{samindef}.)

\item
$(\cA,\theta,p)$ is \hl{\nbd{\infty}algebraically minimal} if $\cA=M(\cA_\infty)$. (Definition \ref{damindef}.)
 
 \item
If it is clear about which variant we are speaking, then we say just $(\cA,\theta,p)$ is \hl{algebraically minimal}.
\end{itemize}
The definitions of \it{minimality} in the literature almost always amount to requiring that $(\cA,\theta,p)$ is algebraically minimal \bf{and} full. (See Subsection \ref{minSEC}\ref{fmsSSEC}, especially Footnote \ref{ArvminFN}.) In this case, we say (and distinguish) $(\cA,\theta,p)$ is \hl{fully} (normally)(strictly)(\nbd{\infty})\hl{minimal}.

\begin{itemize}
\item
If $(\cA,p)$ is full, then $\sF(\cA p)$ generates $\cA$ in the strict topology and, \it{a fortiori}, in the strong topology when $\cA$ is a von Neumann algebra. The only question is, hence, whether $\sF(\cA_\infty p)$ generates $\sF(\cA p)$. Theorem \ref{fminthm}:
\begin{itemize}
\item
The triple is fully strictly minimal if and only if $\cA_\infty p=\cA p$.

\item
The triple is fully normally minimal if and only if $\ol{\cA_\infty p}^s=\cA p$.
\end{itemize}
If the context is clear, we refer to either situation as \hl{fully minimal}.

\item
The case of \it{fully \nbd{\infty}minimal} is kept separate for one or more of the following reasons:
\begin{itemize}
\item
For  $(\cA,\theta,p)$, full or not, $M(\cA_\infty)$ need not be contained neither in $\cA$, nor in $\sB^a(\cA p)$, nor in $\sB^a(\cA_\infty p)$.

\item
If $(\cA,\theta,p)$ is fully minimal, then $M(\cA_\infty)\subset\cA$.  (Corollary \ref{strimincor}.) However, rarely, such a dilation is fully \nbd{\infty}minimal. (Exception: Example \ref{minnuniex}.)

\end{itemize}

\end{itemize}
If a triple (or dilation) is not algebraically minimal, the we might wish to \it{minimalize} it by passing to what $\cA_\infty$ does generate (in one of the mentioned topologies). For that it is necessary that each $\theta_t$ is continuous (in that topology).

\begin{itemize}
\item
Normally:
\begin{itemize}
\item
A normal dilation $(\cA,\theta,p)$ may be minimalized. (See after Definition \ref{namindef}.)

\item
If it is good, then its algebraic minimalization $(\ol{\cA_\infty}^s,{\theta^\infty}^s,p)$ is strong. Especially, a normal CP-semigroup admits a good dilation if and only if it admits a strong dilation. (Observation \ref{algminvNob}.)
\end{itemize}

\item
Strictly:
\begin{itemize}
\item
A strict (hence, full) dilation $(\cA,\theta,p)$ may be minimalized if and only if an approximate unit for $\cA_\infty$ converges strictly.

\item
The minimalization is strict (hence, full) if and only of the approximate unit can be chosen from $\ls\cA_\infty p\cA_\infty $.
\end{itemize}

\item
$\infty$:
\begin{itemize}
\item
The \nbd{\infty}strict topology can be quite strong. (See between Remark \ref{algminstriob} and Lemma \ref{istrstrlem}.)

\item
However: Under the standing hypotheses, every Markov dilation may be \nbd{\infty}mini\-ma\-lized (Corollary \ref{istrMricor}) and the minimalization is \nbd{\infty}strict (Theorem \ref{istrMrithm}).
\end{itemize}
\end{itemize}
We have already said (see p.\pageref{pBT}) what it means to say that $p\in\cB$ compresses the pair $(\cB,T)$. We now say what it means to compress triples (such as dilation).
\begin{itemize}
\item
A projection $P\in\cA$ \hl{compresses} $(\cA,\theta,p)$ (\hl{strongly}) if $P\ge p$ and if $P$ compresses $(\cA,\theta)$ (strongly) to an \nbd{E}semigroup $(\cA^P,\theta^P)$, that is (see Definition \ref{comprdef}), if $(\cA,\theta,P)$ is a (strong) dilation of $\theta^P$.

\item
Observation \ref{compob}:
\begin{itemize}
\item
If $(\cA,\theta,p)$ is a dilation (of $T$), then $(\cA^P,\theta^P,p)$ is a dilation (of the same $T$).

\item
A compression of a full dilation is full.

\item
A compression of a strong dilation is strong.
\end{itemize}

\item
A strong compression of a dilation $(\cA,\theta,p)$ is strong if and only if $(\cA,\theta,p)$ is strong. (Corollary \ref{strcomprcor}.)

\item
A compression of a dilation $(\cA,\theta,p)$ is good if and only if $(\cA,\theta,p)$ is good. (Corollary \ref{goodccor}.)

\item
A projection compressing an algebraically minimal dilation is central. (Corollary \ref{amcccor}.) The compression is algebraically minimal, too. (Corollary \ref{compamincor}.)

\item
An \hl{\nbd{E_0}compression} is a compression to an \nbd{E_0}semigroup. (See p.\pageref{E0c}.) An \nbd{E_0}com\-pression is strong. (Proposition \ref{Mcharprop}.) Only a strong compression of an \nbd{E_0}dilation is necessarily an \nbd{E_0}dilation. (Observation \ref{E0strongob} and Example \ref{E0weakex}.) Every compression of a normally algebraically minimal Markov dilation is an \nbd{E_0}compression, hence, strong. (Corollary \ref{amE0ccor}.)
\end{itemize}
A \it{compressible} dilation is certainly not what we would call \it{minimal}.

\begin{itemize}
\item
A dilation is (\hl{strongly}) \hl{incompressible} if the only (strong) compression it admits is by $P=\U$. (Definition \ref{incomprdefi}.)

\item
Among all strong compressions of a normal dilation, there is a unique minimal strongly incompressible one. (Theorem \ref{vNE0compthm}.) We do not know if the same is true for general compressions or \it{good} compressions. (Remark \ref{sgrem}.)

\item
A fully normally minimal dilation is incompressible. (Corollary \ref{amfinccor}.)

\item
We mentioned already Corollary \ref{incmnfcor}. (Incompressible and full but not fully minimal, in the von Neumann case implies that the minimalization is not full.)
\end{itemize}
There are more results involving the \it{superproduct system} of a dilation (product system if the dilation is full) to be discussed in the next subsection. These include frequently also the hypotheses that the dilation is \it{primary} -- another property that, when missing, excludes a dilation from be candidate for something to be called \it{minimal}.

\begin{itemize}
\item
Under the standing hypotheses: A (full)(\nbd{\infty}algebraically minimal) dilation $(\cA.\theta,p)$ is \hl{primary} if $\theta_t(p)$ increases (strictly)(\nbd{\infty}(strictly) to $\U$. (Corollary \ref{primcor} and the discussion following it.)

We save the von Neumann case (which makes perfect sense also under the standing hypothesis, where $\theta_t(p)$ increases strongly) to the following general definition:

\item
A dilation $(\cA.\theta,p)$ on a von Neumann algebra $\cA$ is \it{primary} if $\bigvee_{t\in\bS}\theta_t(p)=\U$. Definition \ref{vNprimdef}.)
\begin{itemize}
\item
A normally algebraically minimal dilation is primary. (Theorem \ref{vNprimthm}\eqref{vNp1}.)

\item
A  primary dilation has primary compressions. (Theorem \ref{vNprimthm}\eqref{vNp4}.)

\item
A  primary Markov dilation is an \nbd{E_0}dilation. (Theorem \ref{vNprimthm}\eqref{vNp2}.)
\end{itemize}
\item
Most importantly: A normal dilation can be compressed by $\bigvee_{t\in\bS}\theta_t(p)=\U$ to a primary dilation.
\end{itemize}
We could have provided also a \it{strict} definition of \it{primary} for full dilations:
\begin{itemize}
\item
A full dilation $(\cA,\theta,p)$ is \hl{strictly primary} if $E=\cls\theta_\bS(p)E$.

\item
Weakening that condition to $E=\cls^s\theta_\bS(p)E$, this is equivalent to the general definition restricted to the full case.

\item
In the minimal one-parameter case, this definition is exactly the minimality property stated in \cite[Theorem 1.2(2)]{Ske08a} to be true for the dilation constructed via unitalization of a given CP-semigroup when that theorem is applied to the GNS-product system. (This may be used to substitute in the step \eqref{amcn3}$\Rightarrow$\eqref{amcn1} of the proof of Theorem \ref{1pamcnMcor} the application of Proposition \ref{uniprimprop} to show that the dilation in question is primary.)

We did not discuss this in detail, because in the \nbd{C^*}case a dilation can in general not be \it{primarized}.
\end{itemize}
As a last point in this subsection, we just mention the entire short Subsection \ref{minSEC}\ref{nuniminSSEC} for questions about non-uniqueness of dilation that are in some sense minimal.

\lf
\subsection{\normalsize Product systems, superproduct systems, and subproduct systems} \label{spsSSEC}

\it{Superproduct systems} and \it{subproduct systems} are generalizations of \it{product systems}, which point into different directions; where they intersect, we have product systems. While in Section \ref{PSmonoSEC} we rather concentrate on the indexing monoid in the context of product systems, super- and subproduct systems and their generalities are introduced and discussed in Section \ref{SPSpbSEC}. We do not reproduce these two sections and the results in them in any detail in this appendix; it is up to the reader to decide, whether to read them entirely at once, or whether to follow the references throughout the main text to the definitions and results in these sections.

This appendix concentrates on their occurrence in dilation theory.

Just recall:

\begin{itemize}
\item
(Super)(sub)product systems over a monoid $\bS$ are families $\bfam{E_t}_{t\in\bS}$ of correspondences over $\cB$ with $E_0=\cB$. It is important, not to forget $E_0$ and the \hl{marginal} (or \hl{monoid}) \hl{conditions} relating to it.

\item
Definition \ref{SPSUdef}: In the case of \hl{superproduct systems} (denoted $\bfam{E_t}_{t\in\bS}=E^\podot$) we have isometric bilinear \hl{product maps} $v_{s,t}\colon E_s\odot E_t\rightarrow E_{st}$ subject to the same diagrams following Definition \ref{PSdefi} (for product systems), while in the case of \hl{subproduct systems} (denoted $\bfam{E_t}_{t\in\bS}=E^\bodot$) we have isometric bilinear \hl{coproduct maps} $v_{s,t}\colon E_{st}\rightarrow E_s\odot E_t$ subject to the diagrams following Definition \ref{PSdefi} with all arrows inverted.

(Super)(sub)product systems may be \hl{adjointable} in the sense that their structure maps have adjoints, or not.

(In Remark \ref{copurem}, which concludes the Example Section \ref{EXexpSEC}, we point out that it is sometimes convenient to relax the definition further to \hl{productive systems} and \hl{coproductive systems} by dropping isometricity from the hypotheses. This is strictly limited to where it occurs in Section \ref{EXexpSEC}.)

\item
For (super)product systems, we use the \hl{product} notation $x_sy_t:=v_{s,t}(x_s\odot y_t)$ -- a possibility that is painfully missing for subproduct systems.

\item
(Super)(sub)product systems may have \hl{units}, that is, families $\xi^\odot=\bfam{\xi_t}_{t\in\bS}$ of elements $\xi_t\in E_t$  satisfying $\xi_0=\U\in\cB=E_0$ and $\xi_s\odot\xi_t$``$=$''$\xi_{st}$ in the direction of the structure maps, that is, $\xi_s\xi_t=\xi_{st}$ for superproduct systems, and $w_{s,t}\xi_{st}=\xi_s\odot\xi_t$ for subproduct systems. (Definition \ref{SPSUdef}.)

(This is not be confused with weaker notions of units that are around in the literature!)

\item
Likewise, \hl{morphisms} $a_t\colon E_t\rightarrow F_t$ between (super)(sub)product systems (\bf{not necessarily of the same kind}!) are defined as $a_s\odot a_t$``$=$'$a_{st}$ using in each place only the isometric identification. (Anyway, in the non-adjointable case, this is the only possible choice. See Definition \ref{moemisdef}.)

(This is not be confused with weaker notions of morphisms that are around in the literature!)
\end{itemize}
The part ``not necessarily of the same kind'', is referring to the possibility to mix, when there are two families involved (like for morphisms or subfamilies) their ``kind''. For us it is particularly important to consider subproduct subsystems of superproduct systems. But the discussion in Section \ref{SPSpbSEC} is general.

\begin{itemize}
\item
Intuitively, a subfamily $\bfam{F_t}_{t\in\bS}$ of a (super)(sub)product system $\bfam{E_t}_{t\in\bS}$ is a (super)(sub) prod\-uct \hl{subsystem} (not necessarily of the same kind) if
\beqn{
F_{st}
~~~(\text{``$\supset$''})
~
(\text{``$\subset$''})~~~
F_s\odot F_t.
}\eeqn
Here, the structure map used to identify isometrically one side as a super-/subset of the other is the one of the family $\bfam{E_t}_{t\in\bS}$ -- there is no other! --, while the direction of the inclusion $\supset$ or $\subset$ is chosen depending on whether the family $\bfam{F_t}_{t\in\bS}$ ought to be a super- or a subproduct subsystem. (This is detailed carefully in Definition \ref{subdefi}.)

\item
It is important that, following Lemma \ref{isubcoilem}, every subsystem inherits a structure (of the indicated kind, that is super- or sub-) from the structure (not necessarily of the same kind) of the containing family.

After this, an inclusion means that the structure of the subfamily is the inherited one. (Convention \ref{subconv}.)

\item
Inclusion work nicely when iterated (Proposition \ref{bicontprop}), and there are no clashes when one the subfamily or the containing family happens to be a product system (Proposition \ref{PSsubprop}).
\end{itemize}
There are a number of results about the generation of subfamilies. We mention only the most relevant observations.

\begin{itemize}
\item
Foremost, is certainly that fact that we do not know whether the tensor product of \nbd{C^*}cor\-re\-spond\-ences possesses the \it{intersection property}, while the tensor product of von Neumann correspondences does so. (See Appendix \ref{laAPP}.) Consequently, from the easy inclusion we get the unique minimal superproduct subsystem generated by a subfamily (it is even spanned in an obvious way by the subfamily; see Theorem \ref{supintthm}) of a super product system of \nbd{C^*}correspondences, while we do not know if the same is possible for subproduct systems, \it{a fortiori}, for product subsystems. For von Neumann correspondences this problem does not occur.

\item
A set of units of a superproduct system (in particular, a single unit or all units) spans a subproduct subsystem. (Example \ref{unitgenex}.)
\end{itemize}
Several generation results (which we do not list here individually) are particularly powerful (especially regarding the question when the generated subfamily is actually a product system) when the monoid is totally directed and cancellative. The most important results in this situation is, however, not about generation of subsystem, but construction of the containing system that is not yet there to begin with:
\begin{itemize}
\item
Under the stated hypothesis on the monoid: Every subproduct system over such a monoid is contained (uniquely minimally) in a product system (spanned by it as a superproduct system). (Theorem \ref{indlimthm}.)	
\end{itemize}
CP-semigroups and subproduct systems:

\begin{itemize}
\item
The basic equation that shows that a unit $\xi^\odot$ for a (super)(sub)product system determines a CP-semigroup $T^\xi$ via $T^\xi_t:=\AB{\xi_t,\bullet\xi_t}$, is
\beqn{
\AB{\xi_{st},\bullet\xi_{st}}
~=~
\AB{\xi_s\odot\xi_t,\bullet\xi_s\odot\xi_t}
~=~
\AB{\xi_t,\AB{\xi_s,\bullet\xi_s}\xi_t}
}\eeqn
(see Equation \eqref{unitcomp} and the discussion following it). (Note that in the isometric identification of $\xi_{st}$ and $\xi_s\odot\xi_t$ it does not matter if this is done via the product map $v_{s,t}$ of a superproduct system over $\bS$ or via the coproduct map $w_{s,t}$ of a subproduct system over $\bS$.) It is important to note that the CP-semigroup is over the \bf{opposite} monoid of $\bS$, $\bS^{op}$.

\item
If $T$ is a CP-semigroup on $\cB$ over $\bS^{op}$ and $(\sE_t,\xi_t)=$ GNS-$T$ (see Equation \eqref{GNSeq} in the conventions in the introduction on p.\pageref{GNSeq}), then the $\sE_t$ form a subproduct system, the \hl{GNS-subproduct system} of $T$, with with coproduct defined by $w_{s,t}\xi_{s,t}:=\xi_s\odot\xi_t$, and the $\xi_t$ form a unit, the \hl{cyclic unit}.

\item
If there is a unit $\xi^\odot$ in a (super)(sub)product system such that $T=T^\xi$, then this (super)(sub)product system contains the GNS-subproduct system as subproduct subsystem via $\sE_t:=\cls\cB\xi_t\cB$. (See Example \ref{unitgenex} and Observation \ref{uGNSsupob} for a unit in a superproduct system; the statement for a unit in a subproduct system is not more than the uniqueness statement for GNS-constructions (p.\pageref{GNSeq}).)

A subproduct system can have many different units each of which generating the subproduct system and each of which determining a different CP-semigroup. (Observation \ref{elemob}.) But not all subproduct systems arise as GNS-subproduct systems. (Arveson systems that are not one-dimensional, do not.) This changes, when we look at subproduct systems up to (``strict'') Morita equivalence. (See below)

\item
GNS-Subproduct systems need not be adjointable. (Example Section \ref{EXnonadSEC}.)
\end{itemize}
Dilations and superproduct system:

\begin{itemize}
\item
Strict module dilations (or even just strict triples) come along with product systems. The first construction (\cite{Ske02}) is based on the existence of the unit vector $\xi\in E$. (See, again, below for another constructions based on (``strict'') Morita equivalence.) Immitating this construction, dilations (or triples) give rise to superproduct systems.

\item
For a triple $(\cA,\theta,p)$ ($\theta$ being an \nbd{E}semigroup over $\bS^{op}$!) we put $E_t:=\theta_t(p)\cA p$ ($\subset E$) and turn it into a correspondence over $\cB:=p\cA p$ by defining the left action of $b=pap\in\cB\subset\cA$ by $b.x_t:=\theta(b)x_t$. Then the $E_t$ form a superproduct system over $\bS$, the \hl{superproduct system} of the triple, via $v_{s,t}\colon x_s\odot y_t\mapsto\theta_t(x_s)y_t$. (Theorem \ref{E-supPSthm}.)

\item
If the $\xi_t:=\theta_t(p)p$ form a unit (that is, if $(\cA,\theta,p)$ is \hl{good}), then $(\cA,\theta,p)$ is a (\hl{good}) dilation and its superproduct system contains the GNS-subproduct system of the dilated CP-semigroup. (Theorem \ref{sdilunithm}.) We mentioned already that strong implies good and several other result that relate strong and good dilations in Subsection \ref{RdilSSEC}.

\item
The most important consequence: A CP-semigroup whose GNS-subproduct system does not embed into a superproduct system does not possess a good dilation. \it{A fortiori}, a Markov semigroup whose GNS-subproduct system does not embed into a superproduct system does not possess any dilation. (See the Example Section \ref{EXsubnsupSEC} for (even adjointable) subproduct systems that do not embed into superproduct systems; see the part about (``strict'') Morita equivalence below for why this is enough to construct examples for CP- and Markov semigroups of the stated type.

\item
If $(\cA,\theta,p)$ is full, we recover the construction for $\cA=\sB^a(E)$: The superproduct systems of strict (or normal) full triples are product systems. (Example \ref{EPSex}.) For an incompressible normal dilation, also the opposite is true: It is full if and only if its superproduct system is a product system. (Theorem \ref{vNincfthm}.) The superproduct system of a dilation (good or not) contains the superproduct system of any of its compressions including the elements $\xi_t:=\theta_t(p)p$. (Corollary \ref{comPPScor}.) And a normal full dilation is incompressible if its (super)product system has no proper product subsystem containing all $\xi_t$. (Theorem \ref{fpincthm}.) A normal dilation and its ``\it{primarization}'' (Theorem \ref{vNprimthm}\eqref{vNp3}) have the same superproduct system. (Observation \ref{primpPSob}\eqref{pPSob1}.) Under the standing hypotheses: For a strict full primary (hence, \nbd{E_0})dilation (of a Markov semigroup over the opposite of an Ore monoid), the formula $q_t=\theta_t(p)q$ establishes an order isomorphism between the set of all projections compressing the dilation and all projection morphisms $q_t$ satisfying $q_t\xi_t=\xi_t$. (Theorem \ref{Ocomprthm}.)

\item
The Example Section \ref{EXpropsupSEC} presents dilations that have proper superproduct systems. Other examples are the same as for non-full dilations arising from minimalization of full dilations as in Subsection \ref{RdilSSEC}.

\end{itemize}
A basic question is, whether the construction of a superproduct system from a dilation can be reversed. For this and other purposes, notions around the idea of \it{left dilation} of a (super)product system are important. (Note that the corresponding notions of \it{right dilations} exist, but occur rather under the ``commutant point of view'' (as \it{covariant representations}), which we ignore entirely in the main part of these notes, and report only briefly in Appendix \ref{vNAPP}\ref{vNcomm}.)

\begin{itemize}
\item
The concept of left dilation of a superproduct system $E^\podot$ to a Hilbert module $E$ is about families of isometries $v_t\colon E\odot E_t\rightarrow E$ that iterate associatively with superproduct system structure, that is, $(xy_s)z_t=x(y_sz_t)$ where we extend the product notation to $xy_t:=v_t(x\odot y_t)$. The strongest concept, \it{left dilation}, comes from the relation with strict \nbd{E_0}semigroups on $\sB^a(E)$, while there are mollifiers semi and quasi, which also may be combined.
\begin{itemize}
\item
The $v_t$ form a \hl{left dilation} if $E$ is full and the $v_t$ are unitary. (See before Proposition \ref{lsdilPSprop}.) For a left dilation, the maps $\vt_t:=v_t(\bullet\odot\id_t)v_t^*$ define a strict \nbd{E_0}semigroup $\vt$ and the product system associated with $\vt$ (as with any strict \nbd{E_0}semigroup on some $\sB^a(E)$) is $E^\podot$. In particular, a superproduct system admitting a left dilation is a product system. (Proposition \ref{lsdilPSprop}.) A product system admitting a left dilation is necessarily \hl{full}, that is, each $E_t$ is full.

\item
The $v_t$ form a \hl{left semidilation} if $E$ is full. A left semidilation is \hl{adjointable} if each isometry $v_t$ is adjointable. (See before Observation \ref{lsdilPSob}.) For an adjointable left semidilation we may define $\vt_t$ as before, and we obtain a strict \nbd{E}semigroup $\vt$ if and only if $E^\podot$ \bf{is} a product system; also here, calculating the product system of $\vt$ (as associated with any strict \nbd{E}semigroup on $\sB ^a(E)$) we get back $E^\podot$. (Observation \ref{lsdilPSob}.) The superproduct system of a triple $(\cA,\theta,p)$ (Theorem \ref{E-supPSthm}) comes along with a left semidilation $v_t\colon x\odot y_t\mapsto\theta_t(x)y_t$. (Theorem \ref{E-lsemdilthm}.) We do not know it that left semidilation is adjointable. But if a superproduct system admits an adjointable left semidilation, then the superproduct system is adjointable. (Proposition \ref{semiadjprop}.)

\item
The $v_t$ form a \hl{left quasi-dilation} if the $v_t$ are unitary. (The discussion starts after Observation \ref{lsubdilob}, but as we do not need any of the stated results we do not give more details. Dropping fullness of $E$, we still may define strict a \nbd{E_0}semigroup (or an \nbd{E}semigroup if we ad semi and adjointable) but we lose any uniqueness statement and existence is not a problem. (For instance, $E$ may be $\zero$.)

\item
The notion of \hl{left subdilation} of a left semidilation, which is discussed in Observation \ref{lsubdilob}, has some relevance in Subsection \ref{minSEC}\ref{1-p-SSEC}.
\end{itemize}

\item
A reverse of the construction of a superproduct system works under \it{the standing hypotheses} and for product systems: Suppose $E^\odot$ is a product system over an Ore monoid and $\xi^\odot$ is unital unit for $E^\odot$. (So, the $T_t:=\AB{\xi_t,\bullet\xi_t}$ form a Markov semigroup $T$ over $\bS^{op}$ and $E^\odot$ contains the GNS-subproduct system of $T$.) Then, the inductive limit construction in Observation \ref{indlimob} yields a Hilbert \nbd{\cB}module $E$ with a unit vector $\xi$ and a left dilation $v_t\colon E\odot E_t\rightarrow E$ such that with the \nbd{E_0}semigroup $\vt$ defined by $\vt_t:=v_t(\bullet\odot\id_t)v_t^*$ the triple $(E,\vt,p)$ is a strict primary module dilation of $T$. (Theorem \ref{Oreindthm}.)

In fact, $T$, under the standing hypotheses, admits a strict full dilation if and only if its GNS-subproduct system embeds into a product system. (Theorem \ref{Markmodthm}.)

\end{itemize}
Especially, for understanding why subproduct systems that do not embed into product systems, give rise to Markov (CP-)semigroups that have no (good) dilations, we need the following (making up the major part of Section \ref{CPspsSEC}, starting after Remark \ref{CPASrem}). Be aware that for this part $\cB,\cC,\ldots$ are not required to be unital.
\begin{itemize}
\item
If $E$ is full, then it may be viewed as a \hl{Morita equivalence} (see p.\pageref{MEdef}) from $\sK(E)$ to $\cB$. The left action of a correspondence $F$ from $\sK(E)$ to $\cC$ extends to a unique and (automatically) strict action of $\sB^a(E)$. A correspondence from $\sB^a(E)$ to $\cC$ is \hl{strict} (see p.\pageref{SCdef}) if it arises that way from a correspondence from $\sK(E)$ to $\cC$, that is, if and only if the action of $\sB^a(E)$ is strict, that is, if and only if the restriction of this action to $\sK(E)$ is nondegenerate. If $T\colon\sB^a(E)\rightarrow\sB^a(F)$ is a CP-map and $\sF$ its GNS-correspondence, then $T$ is strict if and only if $\sF\odot F$ is strict. (Proposition \ref{striGNSprop}.)

\item
If $T$ is a strict CP-semigroup on $\sB^a(E)$ and $\sF^\bodot$ its GNS-subproduct system, then the $\sE_t:=E^\odot\sF_t\odot E$ ($t\ne0$) and $\sE_0:=\cB$ form a subproduct system via of \nbd{\cB}correspondences, the \hl{subproduct system of \nbd{\cB}correspondences} associated with $T$, via
\beqn{
\sE_{st}
~=~
E^*\odot\sF_{st}\odot E
~\longrightarrow~
E^*\odot\sF_s\odot\sF_t\odot E
~\cong~
E^*\odot\sF_s\odot E\odot E^*\odot \sF_t\odot E
~=~
\sE_s\odot\sE_t,
}\eeqn
where at the step ``$\cong$'' we used that, by strictness, $\cls\sK(E)\sF_t\odot E=\sF_t\odot E$ and that $\sK(E)\cong E\odot E^*$. (Corollary \ref{stri2cor}

\item
In the case of a strict \nbd{E}semigroup $\vt$ on $\sB^a(E)$, this gives the construction of the product system of $\vt$. (Example \ref{EPSex}.) See Observation \ref{primob} and the beginning of the proof of Theorem \ref{Ocomprthm} for other constructions.

\item
Every adjointable subproduct system over a cancellative monoid $\bS$ arises as the subproduct system of \nbd{\cB}correspondences associated with a strict (or normal) CP-semigroup on $\sB^a(E)$ for some full Hilbert \nbd{\cB}module $E$. (Theorem \ref{adSPS-CPthm}.)

\item
This, together with existence of an adjointable subproduct system that does not embed into a superproduct system (Example Section \ref{EXsubnsupSEC}), shows existence of a CP-semigroup with no good dilation and (invoking Theorem \ref{uninonunithm} on unitalisations; see Subsection \ref{Runi1pSSEC} below) and existence of a Markov semigroup with no dilation. (Corollary \ref{noDilcor}.)

\end{itemize}
For Ore monoids, we get dilations of Markov semigroups provided we can embed their GNS-subproduct systems into a product system. (Theorem \ref{Oreindthm}.) Only in the one-parameter case (Theorem \ref{indlimthm}) or, for von Neumann algebras, in the discrete two-parameter case we can guarantee that such an embedding is possible in general. The latter is actually already a consequence of the following:
\begin{itemize}
\item
If $\bS$ is a product of monoids $\bS^k$, then a product system over $\bS$ decomposes, roughly, into a tensor product of the \hl{marginal} product systems over $\bS^k\subset\bS	$. (Section \ref{compSEC}.) The idea is to apply these result in the \nbd{d}parameter case(s) considering them a \nbd{d}fold product of the (solvable!) one-parameter case(s).

\item
Theorem \ref{prodthm} gives a complete description of product systems over a product in terms of their marginal product systems and certain flips among them.

\item
As explained in Section \ref{N0dSEC}, the situation simplifies drastically for $\bS=\N_0^d$. (Theorem \ref{N0dthm}.) Section \ref{N0dSEC} is crucial for both the construction of product systems over $\N_0^2$ (where the set of conditions (here, considered as sufficient) is void; Example Section \ref{EXN02SEC}) and the proof of non-embeddability of a certain subproduct system over $\N_0^3$ (where the conditions (here, considered as necessary) cannot be fulfilled; Example Section \ref{EXsubnsupSEC}).
\end{itemize}
Example Section \ref{EXexpSEC} shows how (super)(sub)product systems over $\N_0^d$ may be \hl{exponentiated} to \hl{time ordered} (super)(sub)product systems over $\R_+^d$ (generalizing the construction of type I Arveson systems, which is the case $d=1$ and $\cB=\C$). This machinery allows to promote discrete \nbd{d}parameter examples to continuous time \nbd{d}parameter examples.

\lf
\subsection{\normalsize Unitalizations and the one-parameter case} \label{Runi1pSSEC}

Many results about Markov semigroups $T$ on $\cB$ and their dilations may be put forward to CP-semigroups $T$ on $\cB$ and their strong(\bf{!}) dilations, by unitalizing CP-semigroups and dilations.

\begin{itemize}
\item
Unitalization is discussed in the last part of Section \ref{monoSEC} starting on p.\pageref{ulize}.
\begin{itemize}
\item
The unitalization of a general \nbd{C^*}algebra $\cB$ is $\wt{\cB}:=\cB\oplus\wt{\U}\C$. Recall that our \nbd{C^*}algebras are unital, so $\wt{\cB}$ is isomorphic to the \nbd{C^*}algebraic direct sum $\C\oplus\cB$, where the $1$ in $\C$ is the central projection $\wt{\U}-\U$ in $\wt{\cB}$.

\item
Maps $T$ on $\cB$ are unitalized to maps $\wt{T}$ on $\wt{\cB}$ as $\wt{T}(b+\wt{\U}\lambda):=T(b)+\wt{\U}\lambda$. The unitalization of a (contractive!) CP-semigroup $T$ on $\cB$ is a Markov semigroup $\wt{T}$ on $\wt{\cB}$. (See Point \ref{GNSuni}.)

\item
We unitalize a triple $(\cA,\theta,p)$ to the triple $(\wt{\cA},\wt{\theta},\wt{p})$, where $\wt{p}:=p+(\wt{\U}-\U)$. Then $(\wt{\cA},\wt{\theta},\wt{p})$ is a dilation, necessarily Markov and $E_0$, (of $\wt{T})$) if and only if $(\cA,\theta,p)$ is a \bf{strong} dilation (of $T$). (Theorem \ref{wstrunithm}.) Moreover, if $\wt{T}$ possesses any dilation, then $T$ possess a strong dilation. (Theorem \ref{uninonunithm}.)
\end{itemize}

\item
How special properties of a strong dilation $(\cA,\theta,p)$ and its unitalization $(\wt{\cA},\wt{\theta},\wt{p})$ imply each other, is subject of Subsection \ref{minSEC}\ref{uniminSSEC}. Full $\Longleftrightarrow$ full and the strict topology is preserved. (Corollary \ref{fufcor}.) Strong compression $Q$ $\Longleftrightarrow$ strong compression $\wt{Q}=Q+(\wt{\U}-\U)$. (Lemma \ref{unicompelem}.) Especially, strongly incompressible $\Longleftrightarrow$ strongly incompressible. (Corollary \ref{unicompecor}.) The strong compression of the unitalization is the unitalization of the strong compression. (Observation \ref{unicompeob}.) Less satisfactory: $(\cA,\theta,p)$ is strictly (normally) algebraically minimal $\Longrightarrow$ $(\wt{\cA},\wt{\theta},\wt{p})$ is strictly (normally) algebraically minimal. (Corollary \ref{uniamincor}.) However, if $(\cA,\theta,p)$ is normal and primary, then normally algebraically minimal $\Longleftrightarrow$ normally algebraically minimal.  (Proposition \ref{uniprimprop}.)

\item
Section \ref{unisupSEC} examines the relation between the superproduct system of a dilation and that of its unitalization. Most noteworthy is that, here, we have that a good property for the superproduct system of the unitalization of a dilation implies the same good property for that of the original dilation.

\end{itemize}
One of the most striking application of \it{unitalization} (but surely not the only one) is Subsection \ref{minSEC}\ref{1-p-SSEC}, where we discuss the one-parameter case. After (re)proving (mainly known, especially in the case $\cB=\sB(G)$) results about the \it{unique minimal} dilation of one-parameter Markov semigroups with our methods, we we use \it{unitalization} to push forward most of them to the not necessarily Markov case. We think that one reason for that these results on the non-Markov case are new even for $\cB=\sB(G)$, is that the unitalization of $\sB(G)$ is not another $\sB(G)$, but requires to know the results in the Markov case for a $\cB$ that is more general than $\sB(G)$.

\lf
\subsection{\normalsize The ``von Neumann versions''}

We mentioned a number of results that we can prove/are true only for von Neumann algebras. (For some of them, the reason is that they are on not necessarily full dilation, for which in the \nbd{C^*}case the strict topology is missing; some really use results that are available only for von Neumann algebras or modules.) On the contrary, most results in the \nbd{C^*}case do easily promote to the von Neumann case. To see this, some basic knowledge on von Neumann modules and correspondences is necessary. A quick version of these prerequisites (basically for readers who know them already or who know very well \nbd{W^*}modules) is in the little more than half a page in Section \ref{topSEC} that starts after Observation \ref{pstrcontob}. A more detailed introduction is in Appendix \ref{vNAPP}. The ``promotion'' of results in the \nbd{C^*}case to the von Neumann case, is explained for Theorem \ref{Markmodthm} in Subsection \ref{vNAPP}\ref{MmvNSSEC}; this covers practically all arguments. The last part about the \it{commutant} of von Neumann correspondences, Subsection \ref{vNAPP}\ref{vNcomm}, only serves the purpose of comparison with results that have already been obtained earlier starting from the \hl{Arveson-Stinespring correspondence} of a (normal) CP-map instead of its \it{GNS-correspondence}, which are commutants of each other.

\lf
\subsection{\normalsize Open questions}

A fairly complete list of direct open questions are those in Section \ref{QSEC} marked with ``Answer: Unknown." Some less exposed questions:
\begin{itemize}
\item
Can every normal dilation be compressed to a (unique) incompressible one? (This is the analogue of Theorem \ref{vNE0compthm}, where $Q$ runs over all compressions. What if $Q$ runs only over all good compressions (see Remark \ref{sgrem})? Does the answer change, if the normal dilation is assumed strong or good?

Especially: Does every CP-semigroup that possesses a normal dilation, posses an algebraically minimal incompressible normal dilation? (We know it is true for \it{strongly incompressible}; see Corollary \ref{vNincmincor}.)

\item
Is there a strong dilation which is not algebraically minimal, but has an algebraically minimal unitalization? (See Corollary \ref{uniamincor}.) This question splits into `full' and `not necessarily full'.

By Proposition \ref{uniprimprop}, such a dilation cannot be primary. (This is for normal dilations; but by the discussion above, also a full and strict example cannot be primary.)

\item
Do tensor products of \nbd{C^*}correspondences possess the \hl{intersection property}? (Appendix \ref{laAPP}.) Especially, are intersections of (sub)product subsystems of \nbd{C^*}correspondences again (sub)product subsystems?

\end{itemize}

\newpage

\setlength{\baselineskip}{2.5ex}

\phantomsection\addcontentsline{toc}{section}{\refname}

% \bibliography{mybib}

\begin{thebibliography}{{\Sort{Kummerer}}K{\"um}85}

\bibitem[AC82]{AcCe82}
L.~Accardi and C.~Cecchini, \emph{{Conditional expectations in von Neumann
  algebras and a theorem of Takesaki}}, J.\ Funct.\ Anal. \textbf{45} (1982),
  245--273.

\bibitem[AFQ92]{AFQ92}
L.~Accardi, F.~Fagnola, and J.~Quaegebeur, \emph{{A representation-free quantum
  stochastic calculus}}, J.\ Funct.\ Anal. \textbf{104} (1992), 149--197.

\bibitem[AHK78]{AlHK78}
S.~Albeverio and R.~Hoegh-Krohn, \emph{{Frobenius theory for positive maps of
  von Neumann algebras}}, Commun.\ Math.\ Phys. \textbf{64} (1978), 83--94.

\bibitem[AK92]{ArKi92}
W.~Arveson and A.~Kishimoto, \emph{{A note on extensions of semigroups of
  $*$--en\-do\-mor\-phisms}}, Proc.\ Amer.\ Math.\ Soc. \textbf{116} (1992),
  169--774.

\bibitem[Ale04]{Ale04}
A.~Alevras, \emph{{One parameter semigroups of endomorphisms of factors of type
  II$_1$}}, J.\ Operator Theory \textbf{51} (2004), 161--179.

\bibitem[And63]{And63}
T.~And\^{o}, \emph{{On a pair of commutative contractions}}, Acta Sci.
  Math.(Szeged) \textbf{24} (1963), 88--90.

\bibitem[Ara70]{Ara70}
H.~Araki, \emph{{Factorizable representations of current algebra}}, Publ.\
  Res.\ Inst.\ Math.\ Sci. \textbf{5} (1970), 361--422.

\bibitem[Arv69]{Arv69}
W.~Arveson, \emph{{Subalgebras of $C^*$--algebras}}, Acta Math. \textbf{123}
  (1969), 141--224.

\bibitem[Arv89a]{Arv89}
\bysame, \emph{{Continuous analogues of Fock space}}, Mem.\ Amer.\ Math.\ Soc.,
  no. 409, American Mathematical Society, 1989.

\bibitem[Arv89b]{Arv89a}
\bysame, \emph{{Continuous analogues of Fock space III: Singular states}}, J.\
  Operator Theory \textbf{22} (1989), 165--205.

\bibitem[Arv90a]{Arv90a}
\bysame, \emph{{Continuous analogues of Fock space II: The spectral
  $C^*$--algebra}}, J.\ Funct.\ Anal. \textbf{90} (1990), 138--205.

\bibitem[Arv90b]{Arv90}
\bysame, \emph{{Continuous analogues of Fock space IV: essential states}}, Acta
  Math. \textbf{164} (1990), 265--300.

\bibitem[Arv97a]{Arv97}
\bysame, \emph{{Minimal $E_0$--semigroups}}, Operator algebras and their
  applications (P.~Fillmore and J.~Mingo, eds.), Fields Inst.\ Commun., no.~13,
  American Mathematical Society, 1997, pp.~1--12.

\bibitem[Arv97b]{Arv97a}
\bysame, \emph{{The index of a quantum dynamical semigroup}}, J.\ Funct.\ Anal.
  \textbf{146} (1997), 557--588.

\bibitem[Arv02]{Arv02}
\bysame, \emph{{The heat flow of the CCR algebra}}, Bull.\ London Math.\ Soc.
  \textbf{34} (2002), 73--83.

\bibitem[Arv03]{Arv03}
\bysame, \emph{{Noncommutative dynamics and $E$--semigroups}}, Monographs in
  Mathematics, Springer, 2003.

\bibitem[Arv06]{Arv06}
\bysame, \emph{{On the existence of $E_0$--semigroups}}, Infin.\ Dimens.\
  Anal.\ Quantum Probab.\ Relat.\ Top. \textbf{9} (2006), 315--320.

\bibitem[AS00]{AcSk00a}
L.~Accardi and M.~Skeide, \emph{{Hilbert module realization of the square of
  white noise and the finite difference algebra}}, Math.\ Notes \textbf{86}
  (2000), 803--818, (Rome, Volterra-Pre\-print 1999/0384).

\bibitem[AS07]{AbSk07}
G.~Abbaspour and M.~Skeide, \emph{{~~~Generators of dynamical systems on
  Hilbert modules}}, Commun.\ Stoch.\ Anal. \textbf{1} (2007), 193--207,
  (ar\-Xiv: math.OA/0611097).

\bibitem[BBLS04]{BBLS04}
S.D. Barreto, B.V.R. Bhat, V.~Liebscher, and M.~Skeide, \emph{{Type I product
  systems of Hilbert modules}}, J.\ Funct.\ Anal. \textbf{212} (2004),
  121--181, (Preprint, Cottbus 2001).

\bibitem[BDH88]{BDH88}
M.~Baillet, Y.~Denizeau, and J.-F. Havet, \emph{{Indice d'une esperance
  conditionnelle}}, Compositio Math. \textbf{66} (1988), 199--236.

\bibitem[BFS08]{BFS08}
B.V.R. Bhat, F.~Fagnola, and M.~Skeide, \emph{{Maximal commutative subalgebras
  invariant for CP-maps: (Counter-)examples}}, Infin.\ Dimens.\ Anal.\ Quantum
  Probab.\ Relat.\ Top. \textbf{11} (2008), 523--539, (ar\-Xiv: 0804.1864v3).

\bibitem[Bha96]{Bha96}
B.V.R. Bhat, \emph{{An index theory for quantum dynamical semigroups}}, Trans.\
  Amer.\ Math.\ Soc. \textbf{348} (1996), 561--583.

\bibitem[Bha02]{Bha02}
\bysame, \emph{{Minimal isometric dilations of operator cocylces}}, Integr.\
  Equat.\ Oper.\ Th. \textbf{42} (2002), 125--141.

\bibitem[Bha03]{Bha03}
\bysame, \emph{{Atomic dilations}}, Advances in quantum dynamics (G.L. Price,
  B~.M. Baker, P.E.T. Jorgensen, and P.S. Muhly, eds.), Contemporary
  Mathematics, no. 335, American Mathematical Society, 2003, pp.~99--107.

\bibitem[BL05]{BhLi05}
B.V.R. Bhat and J.M. Lindsay, \emph{{Regular quantum stochastic cocycles have
  exponential product systems}}, Quantum Probability and Infinite Dimensional
  Analysis --- From Foundations to Applications (M.~Sch\"urmann and U.~Franz,
  eds.), Quantum Probability and White Noise Analysis, no. XVIII, World
  Scientific, 2005, pp.~126--140.

\bibitem[BLS08]{BLS08}
B.V.R. Bhat, V.~Liebscher, and M.~Skeide, \emph{{A problem of Powers and the
  product of spatial product systems}}, Quantum Probability and Infinite
  Dimensional Analysis --- Proceedings of the 28th Conference (J.C. Garcia,
  R.~Quezada, and S.B. Sontz, eds.), Quantum Probability and White Noise
  Analysis, no. XXIII, World Scientific, 2008, (ar\-Xiv: 0801.0042v1),
  pp.~93--106.

\bibitem[BLS10]{BLS10}
\bysame, \emph{{Subsystems of Fock need not be Fock: Spatial CP-semigroups}},
  Proc.\ Amer.\ Math.\ Soc. \textbf{138} (2010), 2443--2456, electronically Feb
  2010, (ar\-Xiv: 0804.2169v2).

\bibitem[BM10]{BhMu10}
B.V.R. Bhat and M.~Mukherjee, \emph{{Inclusion systems and amalgamated products
  of product systems}}, Infin.\ Dimens.\ Anal.\ Quantum Probab.\ Relat.\ Top.
  \textbf{13} (2010), 1--26, (ar\-Xiv: 0907.0095v1).

\bibitem[BMSS12]{BMSS12}
P.~Bikram, K.~Mukherjee, R.~Srinivasan, and V.S. Sunder, \emph{{Hilbert von
  Neumann Modules}}, Commun.\ Stoch.\ Anal. \textbf{6} (2012), 49--64,
  (ar\-Xiv: 1102.4663).

\bibitem[BS00]{BhSk00}
B.V.R. Bhat and M.~Skeide, \emph{{Tensor product systems of Hilbert modules and
  dilations of completely positive semigroups}}, Infin.\ Dimens.\ Anal.\
  Quantum Probab.\ Relat.\ Top. \textbf{3} (2000), 519--575, (Rome,
  Volterra-Pre\-print 1999/0370).

\bibitem[BS15]{BhSk15}
\bysame, \emph{{Pure semigroups of isometries on Hilbert $C^*$--modules}}, J.\
  Funct.\ Anal. \textbf{269} (2015), 1539--1562, electronically Jun 2015.
  Pre\-print, ar\-Xiv: 1408.2631.

\bibitem[Bun84]{Bun84}
J.W. Bunce, \emph{{Models for $n$--tuples of noncommuting operators}}, J.\
  Funct.\ Anal. \textbf{57} (1984), 21--30.

\bibitem[CE79]{ChrEv79}
E.~Christensen and D.E. Evans, \emph{{Cohomology of operator algebras and
  quantum dynamical semigroups}}, J.\ London Math.\ Soc. \textbf{20} (1979),
  358--368.

\bibitem[Con80]{Con80p}
A.~Connes, \emph{{Correspondences}}, His hand-written notes, unpublished, 1980.

\bibitem[CP61]{ClPr61}
A.H. Clifford and G.B. Preston, \emph{{The algebraic theory of semigroups I}},
  Mathematical Surveys, no.~7, American Mathematical Society, 1961.

\bibitem[Dav76]{Dav76}
E.B. Davies, \emph{{Quantum theory of open systems}}, Academic Press, 1976.

\bibitem[DO14]{D-O14}
A.~Dor-On, \emph{{Tensor algebras and subproduct systems arising from Feller
  chains}}, unpublished, 2014.

\bibitem[DO18]{D-O18}
\bysame, \emph{{Isomorphisms of tensor algebras arising from weighted partial
  systems}}, Trans.\ Amer.\ Math.\ Soc. \textbf{370} (2018), 3507--3549.

\bibitem[DOM14]{D-OMa14}
A.~Dor-On and D.~Markiewicz, \emph{{Operator algebras and subproduct systems
  arising from stochastic matrices}}, J.\ Funct.\ Anal. \textbf{267} (2014),
  1057--1120.

\bibitem[DOM17]{D-OMa17}
\bysame, \emph{{\nbd{C^*}envelopes of tensor algebras arising from stochastic
  matrices}}, Integr.\ Equat.\ Oper.\ Th. \textbf{88} (2017), 185--227.

\bibitem[DRS11]{DRS11}
K.R. Davidson, C.~Ramsey, and O.M. Shalit, \emph{{The isomorphism problem for
  some universal operator algebras}}, Adv.\ Math. \textbf{228} (2011),
  167--218, (ar\-Xiv: 1010.0729v2).

\bibitem[Ell00]{Ell00}
G.A. Elliott, \emph{{On the convergence of a sequence of completely positive
  maps to the identity}}, J.\ Austral.\ Math.\ Soc.\ Ser.\ A \textbf{68}
  (2000), 340--348.

\bibitem[EN06]{EnNa06}
K.-J. Engel and R.~Nagel, \emph{{A short course on operator semigroups}},
  Universitext, Springer, 2006.

\bibitem[Fow02]{Fow02}
N.J. Fowler, \emph{{Discrete product systems of Hilbert bimodules}}, Pac.\ J.\
  Math. \textbf{204} (2002), 335--375.

\bibitem[Fra82]{Fra82}
A.E. Frazho, \emph{{Models for noncommuting operators}}, J.\ Funct.\ Anal.
  \textbf{48} (1982), 1--11.

\bibitem[Fra02]{FraM02}
M.~Frank, \emph{{On Hahn-Banach type theorems for Hilbert $C^*$--modules}},
  Int.\ J.\ Math. \textbf{13} (2002), 675--693.

\bibitem[FS07]{FaSk07}
F.~Fagnola and M.~Skeide, \emph{{~~Restrictions of CP-semigroups to maximal
  commutative subalgebras}}, Banach Center Publications \textbf{78} (2007),
  121--132, (ar\-Xiv: math.OA/0703001).

\bibitem[GS99]{GoSi99}
D.~Goswami and K.B. Sinha, \emph{{Hilbert modules and stochastic dilation of a
  quantum dynamical semigroup on a von Neumann algebra}}, Commun.\ Math.\ Phys.
  \textbf{205} (1999), 377--403.

\bibitem[GS05]{GoSk05}
R.~Gohm and M.~Skeide, \emph{{Constructing extensions of CP-maps via tensor
  dilations with the help of von Neumann modules}}, Infin.\ Dimens.\ Anal.\
  Quantum Probab.\ Relat.\ Top. \textbf{8} (2005), 291--305, (ar\-Xiv:
  math.OA/0311110).

\bibitem[GS20]{GeSk20b}
M.~Gerhold and M.~Skeide, \emph{{Interacting Fock spaces and subproduct
  systems}}, Infin.\ Dimens.\ Anal.\ Quantum Probab.\ Relat.\ Top. \textbf{23}
  (2020), n.3 , Article 1 (53 pp), Pre\-print, ar\-Xiv: 1808.07037v3.

\bibitem[Gui72]{Gui72}
A.~Guichardet, \emph{{Symmetric Hilbert spaces and related topics}}, Lect.\
  Notes Math., no. 261, Springer, 1972.

\bibitem[Hal70]{Hal70}
P.R. Halmos, \emph{{Ten problems in Hilbert space}}, Bull.\ Amer.\ Math.\ Soc.
  \textbf{76} (1970), 887--933.

\bibitem[Hir04]{Hir04}
I.~Hirshberg, \emph{{$C^*$--Algebras of Hilbert module product systems}}, J.\
  Reine Angew.\ Math. \textbf{570} (2004), 131--142.

\bibitem[Hir05]{Hir05a}
\bysame, \emph{{Essential representations of $C^*$--correspondences}}, Int.\
  J.\ Math. \textbf{16} (2005), 765--775.

\bibitem[HKK98]{HKK98}
J.~Hellmich, C.~K\"ostler, and B.~K\"ummerer, \emph{{Stationary quantum Markov
  processes as solutions of stochastic differential equations}}, Quantum
  probability (R.~Alicki, M.~Bozejko, and W.A. Majewski, eds.), Banach Center
  Publications, vol.~43, Polish Academy of Sciences --- Institute of
  Mathematics, 1998, pp.~217--229.

\bibitem[HKK04]{HKK04p}
\bysame, \emph{{Noncommutative continuous Bernoulli sifts}}, Pre\-print,
  ar\-Xiv: \newline math.OA/0411565, 2004.

\bibitem[HP57]{HiPhi57}
E.~Hille and R.S. Phillips, \emph{{Functional analysis and semi-groups}},
  American Mathematical Society, 1957.

\bibitem[HP84]{HuPa84}
R.L. Hudson and K.R. Parthasarathy, \emph{{Quantum Ito's formula and stochastic
  evolutions}}, Commun.\ Math.\ Phys. \textbf{93} (1984), 301--323.

\bibitem[Kas80]{Kas80}
G.G. Kasparov, \emph{{Hilbert $C^*$--modules, theorems of Stinespring {\&}
  Voi\-cu\-les\-cu}}, J.\ Operator Theory \textbf{4} (1980), 133--150.

\bibitem[{\Sort{Kostler}}K{\"os}00]{Koes00}
C.~{\Sort{Kostler}}K{\"os}tler, \emph{{Quanten-Markoff-Prozesse und
  Quanten-Brownsche Bewegungen}}, Ph.D. thesis, Stuttgart, 2000.

\bibitem[KS]{KaSk21p}
J.~Kaad and M.~Skeide, \emph{{Kernels of Hilbert module maps: A
  counterexample}}, Pre\-print, ar\-Xiv: 2101.03030v1 (to appear in J.\
  Operator Theory).

\bibitem[KS92]{KueSp92}
B.~K\"ummerer and R.~Speicher, \emph{{Stochastic integration on the Cuntz
  algebra $O_\infty$}}, J.\ Funct.\ Anal. \textbf{103} (1992), 372--408.

\bibitem[{\Sort{Kummerer}}K{\"um}85]{Kuem85}
B.~{\Sort{Kummerer}}K{\"um}merer, \emph{{Markov dilations on $W^*$--algebras}},
  J.\ Funct.\ Anal. \textbf{63} (1985), 139--177.

\bibitem[Lac00]{Lac00}
M.~Laca, \emph{{From endomorphisms to automorphisms and back: Dilations and
  full corners}}, J.\ London Math.\ Soc. \textbf{61} (2000), 893--904.

\bibitem[Lan95]{Lan95}
E.C. Lance, \emph{{Hilbert $C^*$--modules}}, Cambridge University Press, 1995.

\bibitem[Lev20]{Lev20p}
E.~Levy, \emph{{(Positive) totally ordered noncommutative monoids -- how
  noncommutative can they be?}}, Pre\-print, ar\-Xiv: 2006.00886, 2020.

\bibitem[Lie09]{Lie09}
V.~Liebscher, \emph{{Random sets and invariants for (type II) continuous tensor
  product systems of Hilbert spaces}}, Mem.\ Amer.\ Math.\ Soc., no. 930,
  American Mathematical Society, 2009, (ar\-Xiv: math.PR/0306365).

\bibitem[Lon84]{Lon84}
R.~Longo, \emph{{Solution of the factorial Stone-Weierstrass conjecture}},
  Invent.\ Math. \textbf{76} (1984), 145--155.

\bibitem[LP12]{LaPa12}
G.~Landi and A.~Pavlov, \emph{{On orthogonal systems in Hilbert
  $C^*$--modules}}, J.\ Operator Theory \textbf{68} (2012), 487--500.

\bibitem[LS01]{LiSk01}
V.~Liebscher and M.~Skeide, \emph{{Units for the time ordered Fock module}},
  Infin.\ Dimens.\ Anal.\ Quantum Probab.\ Relat.\ Top. \textbf{4} (2001),
  545--551, (Rome, Volterra-Pre\-print 2000/0411).

\bibitem[Lu97]{Lu97}
Y.G. Lu, \emph{{An interacting free Fock space and the arcsine law}}, Prob.\
  Math.\ Statist. \textbf{17} (1997), 149--166.

\bibitem[Mar03]{Mar03}
D.~Markiewicz, \emph{{On the product system of a completely positive
  semigroup}}, J.\ Funct.\ Anal. \textbf{200} (2003), 237--280.

\bibitem[MS98]{MuSo98}
P.S. Muhly and B.~Solel, \emph{{Tensor algebras over $C^*$--correspondences:
  representations, dilations, and $C^*$--envelopes}}, J.\ Funct.\ Anal.
  \textbf{158} (1998), 389--457.

\bibitem[MS00]{MuSo00}
\bysame, \emph{{On the Morita equivalence of tensor algebras}}, Proc.\ London
  Math.\ Soc. \textbf{81} (2000), 113--168.

\bibitem[MS02]{MuSo02}
\bysame, \emph{{Quantum Markov processes (correspondences and dilations)}},
  Int.\ J.\ Math. \textbf{51} (2002), 863--906, (ar\-Xiv: math.OA/0203193).

\bibitem[MS04]{MuSo04}
\bysame, \emph{{Hardy algebras, $W^*$--correspondences and interpolation
  theory}}, Math.\ Ann. \textbf{330} (2004), 353--415, (ar\-Xiv:
  math.OA/0308088).

\bibitem[MS05]{MuSo05}
\bysame, \emph{{~~~Duality of $W^*$-correspondences and applications}}, Quantum
  Probability and Infinite Dimensional Analysis --- From Foundations to
  Applications (M.~Sch\"urmann and U.~Franz, eds.), Quantum Probability and
  White Noise Analysis, no. XVIII, World Scientific, 2005, pp.~396--414.

\bibitem[MS07]{MuSo07}
\bysame, \emph{{Quantum Markov semigroups (product systems and
  subordination)}}, Int.\ J.\ Math. \textbf{18} (2007), 633--669, (ar\-Xiv:
  math.OA/0510653).

\bibitem[MS10]{MaSha10}
D.~Markiewicz and O.M. Shalit, \emph{{Continuity of CP-semigroups in the
  point-strong operator topology}}, J.\ Operator Theory \textbf{64} (2010),
  149--154, (ar\-Xiv: 0711.0111v1).

\bibitem[MS13]{MaSr13}
O.T. Margetts and R.~Srinivasan, \emph{{Invariants for $E_0$--semigroups on
  II$_1$ factors}}, Commun.\ Math.\ Phys. \textbf{323} (2013), 21155--1184.

\bibitem[MS14]{MaSr14p}
\bysame, \emph{{\hfill Non-cocycle-conjugate $E_0$--semigroups on factors}},
  Pre\-print, ar\-Xiv: 1404.5934v2, 2014.

\bibitem[MS17]{MurSu17p}
S.P. Murugan and S.~Sundar, \emph{{\nbd{E_0^P}semigroups and product systems}},
  Pre\-print, ar\-Xiv: math.OA/1706.03928, 2017.

\bibitem[MS18]{MurSu18}
\bysame, \emph{{On the existence of \nbd{E_0}semigroups -- the multiparameter
  case}}, Infin.\ Dimens.\ Anal.\ Quantum Probab.\ Relat.\ Top. \textbf{21}
  (2018), 1850007.

\bibitem[MS19]{MurSu19}
\bysame, \emph{{An essential representation for a product system over a
  finitely generated subsemigroup of ${\mathbb Z}^d$}}, Proc. Indian Acad. Sci.
  Math. Sci. \textbf{129} (2019), 17.

\bibitem[MSS06]{MSS06}
P.S. Muhly, M.~Skeide, and B.~Solel, \emph{{Representations of $\sB^a(E)$}},
  Infin.\ Dimens.\ Anal.\ Quantum Probab.\ Relat.\ Top. \textbf{9} (2006),
  47--66, (ar\-Xiv: math.OA/0410607).

\bibitem[Mur97]{MurN97}
N.~Muraki, \emph{{Noncommutative Brownian motion in monotone Fock space}},
  Commun.\ Math.\ Phys. \textbf{183} (1997), 557--570.

\bibitem[Par70]{Par70}
S.~Parrot, \emph{{Unitary dilations for commuting contractions}}, Pacific J.\
  Math. \textbf{34} (1970), 481--490.

\bibitem[Pas73]{Pas73}
W.L. Paschke, \emph{{Inner product modules over $B^*$--algebras}}, Trans.\
  Amer.\ Math.\ Soc. \textbf{182} (1973), 443--468.

\bibitem[Pet90]{Pet90}
D.~Petz, \emph{{An invitation to the algebra of canonical commutation
  relations}}, Leuven University Press, 1990.

\bibitem[Pim97]{Pim97}
M.V. Pimsner, \emph{{A class of $C^*$--algebras generalizing both Cuntz-Krieger
  algebras and crossed products by $\Z$}}, Free probability theory (D.V.
  Voiculescu, ed.), Fields Inst.\ Commun., no.~12, 1997, pp.~189--212.

\bibitem[Pop89]{GPop89}
G.~Popescu, \emph{{Isometric dilations for infinite sequences of noncommuting
  operators}}, Trans.\ Amer.\ Math.\ Soc. \textbf{316} (1989), 523--536.

\bibitem[Pop08]{MPop08}
M.~Popa, \emph{{A combinatorial approach to monotonic independence over a
  $C^*$--algebra}}, Pac.\ J.\ Math. \textbf{237} (2008), 299--325, (ar\-Xiv:
  math.OA/0612570v3).

\bibitem[Pow04]{Pow04}
R.T. Powers, \emph{{Addition of spatial $E_0$--semigroups}}, Operator algebras,
  quantization, and noncommutative geometry, Contemporary Mathematics, no. 365,
  American Mathematical Society, 2004, pp.~281--298.

\bibitem[PS72]{PaSchm72}
K.R. Parthasarathy and K.~Schmidt, \emph{{Positive definite kernels, continuous
  tensor products, and central limit theorems of probability theory}}, Lect.\
  Notes Math., no. 272, Springer, 1972.

\bibitem[PS11]{PowSo11}
S.C. Power and B.~Solel, \emph{{Operator algebras associated with unitary
  commutation relations}}, J.\ Funct.\ Anal. \textbf{260} (2011), 1583--1614.

\bibitem[Rie69]{Rie69}
M.A. Rieffel, \emph{{Multipliers and tensor products of $L^p$--spaces of
  locally compact groups}}, Studia Math. \textbf{33} (1969), 71--82.

\bibitem[Rie74]{Rie74a}
\bysame, \emph{{Morita equivalence for $C^*$--algebras and $W^*$--algebras}},
  J.\ Pure Appl.\ Algebra \textbf{5} (1974), 51--96.

\bibitem[Sar65]{Sar65}
D.~Sarason, \emph{{On spectral sets having connected complement}}, Acta Sci.
  Math.(Szeged) \textbf{26} (1965), 289--299.

\bibitem[Sch93]{MSchue93}
M.~Sch\"urmann, \emph{{White noise on bialgebras}}, Lect.\ Notes Math., no.
  1544, Springer, 1993.

\bibitem[Sha08a]{Sha08}
O.M. Shalit, \emph{{$E_0$--dilation of strongly commuting CP$_0$--semigroups}},
  J.\ Funct.\ Anal. \textbf{255} (2008), 46--89, (ar\-Xiv: 0707.1760v2).

\bibitem[Sha08b]{Sha08b}
\bysame, \emph{{What type of dynamics arise in \nbd{E_0}dilations of commuting
  quantum Markov semigroups?}}, Infin.\ Dimens.\ Anal.\ Quantum Probab.\
  Relat.\ Top. \textbf{11} (2008), 393--403.

\bibitem[Sha09]{Sha09}
\bysame, \emph{{Product systems, subproduct systems, and the dilation theory of
  completely positive semigroups}}, Ph.D. thesis, Technion, Haifa, 2009,
  (ar\-Xiv: 1002.4920).

\bibitem[Sha10a]{Sha08CORR}
\bysame, \emph{{\it{Corrigendum} to \cite{Sha08}}}, J.\ Funct.\ Anal.
  \textbf{258} (2010), 1068--1069.

\bibitem[Sha10b]{Sha10}
\bysame, \emph{{Representing a product system representation as a contractive
  semigroup and applications to regular isometric dilations}}, Canad. Math.
  Bull. \textbf{53} (2010), 550--563.

\bibitem[Sha11]{Sha11}
\bysame, \emph{{E-dilation of strongly commuting CP-semigroups (the nonunital
  case)}}, Houston J.\ Math. \textbf{37} (2011), 203--232, (ar\-Xiv:
  0711.2885).

\bibitem[Ske98]{Ske98}
M.~Skeide, \emph{{Hilbert modules in quantum electro dynamics and quantum
  probability}}, Commun.\ Math.\ Phys. \textbf{192} (1998), 569--604, (Rome,
  Volterra-Pre\-print 1996/0257).

\bibitem[Ske00a]{Ske00b}
\bysame, \emph{{Generalized matrix $C^*$--algebras and representations of
  Hilbert modules}}, Mathematical Proceedings of the Royal Irish Academy
  \textbf{100A} (2000), 11--38, (Cott\-bus, Rei\-he Mathe\-ma\-tik 1997/M-13).

\bibitem[Ske00b]{Ske00a}
\bysame, \emph{{Indicator functions of intervals are totalizing in the
  symmetric Fock space $\Gamma(L^2(\R_+))$}}, Trends in contemporary infinite
  dimensional analysis and quantum probability (L.~Accardi, H.-H. Kuo,
  N.~Obata, K.~Saito, {Si Si}, and L.~Streit, eds.), Natural and Mathematical
  Sciences Series, vol.~3, Istituto Italiano di Cultura (ISEAS), Kyoto, 2000,
  Volume in honour of Takeyuki Hida, (Rome, Volterra-Pre\-print 1999/0395),
  pp.~421--424.

\bibitem[Ske00c]{Ske00}
\bysame, \emph{{Quantum stochastic calculus on full Fock modules}}, J.\ Funct.\
  Anal. \textbf{173} (2000), 401--452, (Rome, Volterra-Pre\-print 1999/0374).

\bibitem[Ske01a]{Ske01}
\bysame, \emph{{Hilbert modules and applications in quantum probability}},
  Ha\-bi\-li\-ta\-tions\-schrift, Cottbus, 2001, Available at
  {\footnotesize\url{ http://web.unimol.it/skeide/}}.

\bibitem[Ske01b]{Ske01a}
\bysame, \emph{{Tensor product systems of CP-semigroups on $\C^2$}}, J.\ Math.\
  Sci. \textbf{106} (2001), 2890--2895, (Rome, Volterra-Pre\-print 1999/0387).

\bibitem[Ske02]{Ske02}
\bysame, \emph{{Dilations, product systems and weak dilations}}, Math.\ Notes
  \textbf{71} (2002), 914--923.

\bibitem[Ske03]{Ske03c}
\bysame, \emph{{Commutants of von Neumann modules, representations of
  $\sB^a(E)$ and other topics related to product systems of Hilbert modules}},
  Advances in quantum dynamics (G.L. Price, B.M. Baker, P.E.T. Jorgensen, and
  P.S. Muhly, eds.), Contemporary Mathematics, no. 335, American Mathematical
  Society, 2003, (Preprint, Cottbus 2002, ar\-Xiv: math.OA/0308231),
  pp.~253--262.

\bibitem[Ske04a]{Ske04}
\bysame, \emph{{Independence and product systems}}, Recent developments in
  stochastic analysis and related topics (S.~Albeverio, Z.-M. Ma, and
  M.~R\"ockner, eds.), World Scientific, 2004, (ar\-Xiv: math.OA/0308245),
  pp.~420--438.

\bibitem[Ske04b]{Ske04p}
\bysame, \emph{{Unit vectors, Morita equivalence and endomorphisms}},
  Pre\-print, ar\-Xiv: math.OA/0412231v4 (Version 4), 2004.

\bibitem[Ske05a]{Ske05b}
\bysame, \emph{{Lévy processes and tensor product systems of Hilbert modules}},
  Quantum Probability and Infinite Dimensional Analysis --- From Foundations to
  Applications (M.~Sch\"urmann and U.~Franz, eds.), Quantum Probability and
  White Noise Analysis, no. XVIII, World Scientific, 2005, pp.~492--503.

\bibitem[Ske05b]{Ske05a}
\bysame, \emph{{Three ways to representations of $\sB^a(E)$}}, Quantum
  Probability and Infinite Dimensional Analysis --- From Foundations to
  Applications (M.~Sch\"urmann and U.~Franz, eds.), Quantum Probability and
  White Noise Analysis, no. XVIII, World Scientific, 2005, (ar\-Xiv:
  math.OA/0404557), pp.~504--517.

\bibitem[Ske05c]{Ske05c}
\bysame, \emph{{Von Neumann modules, intertwiners and self-duality}}, J.\
  Operator Theory \textbf{54} (2005), 119--124, (ar\-Xiv: math.OA/0308230).

\bibitem[Ske06a]{Ske06}
\bysame, \emph{{~A simple proof of the fundamental theorem about Arveson
  systems}}, Infin.\ Dimens.\ Anal.\ Quantum Probab.\ Relat.\ Top. \textbf{9}
  (2006), 305--314, (ar\-Xiv: math.OA/0602014).

\bibitem[Ske06b]{Ske06b}
\bysame, \emph{{Commutants of von Neumann correspondences and duality of
  Eilen\-berg-Watts theorems by Rieffel and by Blecher}}, Banach Center
  Publications \textbf{73} (2006), 391--408, (ar\-Xiv: math.OA/0502241).

\bibitem[Ske06c]{Ske06c}
\bysame, \emph{{Generalized unitaries and the Picard group}}, Proc.\ Ind.\ Ac.\
  Sc.\ (Math Sc.) \textbf{116} (2006), 429--442, (ar\-Xiv: math.OA/0511661).

\bibitem[Ske06d]{Ske06d}
\bysame, \emph{{The index of (white) noises and their product systems}},
  Infin.\ Dimens.\ Anal.\ Quantum Probab.\ Relat.\ Top. \textbf{9} (2006),
  617--655, (Rome, Volterra-Pre\-print 2001/0458, ar\-Xiv: math.OA/0601228).

\bibitem[Ske07]{Ske07}
\bysame, \emph{{$E_0$--semigroups for continuous product systems}}, Infin.\
  Dimens.\ Anal.\ Quantum Probab.\ Relat.\ Top. \textbf{10} (2007), 381--395,
  (ar\-Xiv: math.OA/0607132).

\bibitem[Ske08a]{Ske08}
\bysame, \emph{{Isometric dilations of representations of product systems via
  commutants}}, Int.\ J.\ Math. \textbf{19} (2008), 521--539, (ar\-Xiv:
  math.OA/0602459).

\bibitem[Ske08b]{Ske08a}
\bysame, \emph{{Product systems; a survey with commutants in view}}, Quantum
  Stochastics and Information (V.P. Belavkin and M.~Guta, eds.), World
  Scientific, 2008, (Preferable version: ar\-Xiv: 0709.0915v1!), pp.~47--86.

\bibitem[Ske09a]{Ske09a}
\bysame, \emph{{$E_0$--Semigroups for continuous product systems: The nonunital
  case}}, Banach J.\ Math.\ Anal. \textbf{3} (2009), 16--27, (ar\-Xiv:
  0901.1754v1).

\bibitem[Ske09b]{Ske09r2}
\bysame, \emph{{Free product systems}}, Mini-workshop: Product systems and
  independence in quantum dynamics (B.V.R. Bhat, U.~Franz, and M.~Skeide,
  eds.), Oberwolfach Reports, no. 09/2009, Mathematisches Forschungsinstitut
  Oberwolfach, 2009, available
  at:\\{\footnotesize\url{https://www.mfo.de/occasion/0908b}}, pp.~528--530.

\bibitem[Ske09c]{Ske09}
\bysame, \emph{{Unit vectors, Morita equivalence and endomorphisms}}, Publ.\
  Res.\ Inst.\ Math.\ Sci. \textbf{45} (2009), 475--518, (ar\-Xiv:
  math.OA/0412231v5 (Version 5)).

\bibitem[Ske10]{Ske10}
\bysame, \emph{{The Powers sum of spatial CPD-semigroups and CP-semigroups}},
  Banach Center Publications \textbf{89} (2010), 247--263, (ar\-Xiv:
  0812.0077).

\bibitem[Ske11]{Ske11a}
\bysame, \emph{{Nondegenerate representations of continuous product systems}},
  J.\ Operator Theory \textbf{65} (2011), 71--85, (ar\-Xiv: math.OA/0607362).

\bibitem[Ske12]{Ske12p}
\bysame, \emph{{Hilbert von Neumann modules versus concrete von Neumann
  modules}}, Pre\-print, ar\-Xiv: 1205.6413v2 (revised 2021, to appear), 2012.

\bibitem[Ske16]{Ske16}
\bysame, \emph{{Classification of $E_0$--semigroups by product systems}}, Mem.\
  Amer.\ Math.\ Soc., no. 1137, American Mathematical Society, 2016,
  electronically Oct 2015. Pre\-print, ar\-Xiv: 0901.1798v4.

\bibitem[Ske22a]{Ske08p}
\bysame, \emph{{Boolean couplings of spatial product systems}}, Pre\-print, in
  preparation, 2022.

\bibitem[Ske22b]{Ske08p2}
\bysame, \emph{{Free product systems generated by spatial tensor product
  systems}}, Pre\-print, in preparation, 2022.

\bibitem[SNFBK10]{SzNF2010}
B.~Sz.-Nagy, C.~Foias, H.~Bercovici, and K\'{e}rchy, \emph{{Harmonic analysis
  of operators on Hilbert space}}, Springer Science \& Business Media, 2010.

\bibitem[Sol06]{Sol06}
B.~Solel, \emph{{Representations of product systems over semigroups and
  dilations of commuting CP maps}}, J.\ Funct.\ Anal. \textbf{235} (2006),
  593--618, (ar\-Xiv: math.OA/0502423).

\bibitem[Sol08]{Sol08}
\bysame, \emph{{Regular dilations of representations of product systems}},
  Mathematical Proceedings of the Royal Irish Academy \textbf{108} (2008),
  89--110, (ar\-Xiv: math.OA/0504129).

\bibitem[SS09]{ShaSo09}
O.M. Shalit and B.~Solel, \emph{{Subproduct systems}}, Documenta Math.
  \textbf{14} (2009), 801--868, (ar\-Xiv: 0901.1422v2).

\bibitem[SS11]{ShaSk11}
O.M. Shalit and M.~Skeide, \emph{{Three commuting, unital, completely positive
  maps that have no minimal dilation}}, Integr.\ Equat.\ Oper.\ Th. \textbf{51}
  (2011), 55--63, (ar\-Xiv: 1012.2111v2).

\bibitem[SS14]{SkSU14}
M.~Skeide and K.~Sumesh, \emph{{CP-H-Extendable maps between Hilbert modules
  and CPH-semigroups}}, J.\ Math.\ Anal.\ Appl. \textbf{414} (2014), 886--913,
  electronically Jan 2014. Pre\-print, ar\-Xiv: 1210.7491v2.

\bibitem[Sti55]{Sti55}
W.F. Stinespring, \emph{{Positive functions on $C^*$--algebras}}, Proc.\ Amer.\
  Math.\ Soc. \textbf{6} (1955), 211--216.

\bibitem[Str69]{Str69}
R.F. Streater, \emph{{Current commutation relations, continuous tensor products
  and infinitely divisible group representations}}, Local quantum theory
  (R.~Jost, ed.), Academic Press, 1969.

\bibitem[Tak02]{Tak02}
M.~Takesaki, \emph{{Theory of operator algebras I}}, Encyclopaedia of
  Mathematical Sciences, no. 124 (number V in the subseries Operator Algebras
  and Non-Commutative Geometry), Springer, 2002, 2nd printing of the First
  Edition 1979.

\bibitem[Tak03]{Tak03a}
\bysame, \emph{{Theory of operator algebras II}}, Encyclopaedia of Mathematical
  Sciences, no. 125 (number VI in the subseries Operator Algebras and
  Non-Commutative Geometry), Springer, 2003.

\bibitem[Tho11]{Tho11}
A.~Thom, \emph{{A remark about the Connes fusion tensor product}}, Theory Appl.
  Categ. \textbf{25} (2011), 38--50, (ar\-Xiv: math0601045).

\bibitem[Tsi00]{Tsi00p1}
B.~Tsirelson, \emph{{From random sets to continuous tensor products: answers to
  three questions of W.\ Arveson}}, Pre\-print, ar\-Xiv: math.FA/0001070, 2000.

\bibitem[Ver16]{Ver16}
A.~Vernik, \emph{{Dilations of CP-maps commuting according to a graph}},
  Houston J.\ Math. \textbf{42} (2016), 1291--1329.

\bibitem[Vis10]{Vis10}
A.~Viselter, \emph{{Covariant representations of subproduct systems}}, {Proc.\
  Lond.\ Math.\ Soc.} \textbf{102} (2010), 767--800.

\bibitem[Vis12]{Vis12}
\bysame, \emph{{Cuntz-Pimsner algebras for subproduct systems}}, Int.\ J.\
  Math. \textbf{23} (2012), 1250081.

\bibitem[Wig39]{Wig39}
E.P. Wigner, \emph{{On unitary representations of the inhomogeneous Lorenz
  group}}, Adv.\ Math. \textbf{40} (1939), 149--204.

\end{thebibliography}

\newcommand{\Swap}[2]{#2#1}\newcommand{\Sort}[1]{}
\providecommand{\bysame}{\leavevmode\hbox to3em{\hrulefill}\thinspace}
\providecommand{\MR}{\relax\ifhmode\unskip\space\fi MR }
% \MRhref is called by the amsart/book/proc definition of \MR.
\providecommand{\MRhref}[2]{%
  \href{http://www.ams.org/mathscinet-getitem?mr=#1}{#2}
}
\providecommand{\href}[2]{#2}

% \newpage

% \lf\lf\noindent
% \bf{Acknowledgments.~}
% After some loose e-mail discussions starting with the Fields Workshop ``Advances in Quantum Dynamics'' 2007, we started working seriously on this pro\-ject when we met, again in Canada, in November 2009 in Waterloo, where OS was spending his first year as postdoc, and in January 2010 in Kingston, where MS was spending the first four months of his sabbatical. We are deeply grateful to our hosts, Ken Davidson in Waterloo and Roland Speicher in Kingston, who have made these encounters possible. MS is grateful to the math departments of Ben-Gurion University and the Technion for having him hosted numerous times, as well as for (partial) economic support from his department of economy. MS is also grateful to OS for having made this possible, and for having taken on once the adventure of a visit in Campobasso. OS is grateful to MS for the wonderful hospitality granted during that adventure. OS is also grateful for partial support granted by ISF Grants no. 195/16 and 431/20.

% Last but not least, it is our wish to thank the referees who took on the task of giving this manuscript a careful reading.

\lf\lf\noindent
%%%% BO postarxiv
%Orr Moshe Shalit: \it{Technion, Haifa},
Orr Moshe Shalit: \it{Technion - Israel Institute of Technology, Haifa 3200003, Israel},
%%%% EO postarxiv
E-mail: \href{mailto:oshalit@technion.ac.il}{\tt{oshalit@technion.ac.il}},\\
Homepage: \url{https://oshalit.net.technion.ac.il/}

\lf\noindent
Michael Skeide: \it{Dipartimento di Economia, Universit\`{a} degli Studi del Molise, Via de Sanctis, 86100 Campobasso, Italy},
E-mail: \href{mailto:skeide@unimol.it}{\tt{skeide@unimol.it}},\\
Homepage: \url{http://web.unimol.it/skeide/}

% \newpage
\newpage

\clearpage
\phantomsection\addcontentsline{toc}{section}{Index}
\printindex

\setlength{\parskip}{1.5ex}

% \listofOWs
% \listoftodos

\end{document}